\documentclass[10pt]{book}
\usepackage[english]{babel}
\usepackage{latexsym}
\usepackage{amssymb,amsbsy,amsmath,amsfonts,amssymb,amscd}
\usepackage{amsmath, amsthm}
\usepackage[utf8]{inputenc}
\usepackage{epsfig, graphicx}
\usepackage[dvipsnames]{xcolor}
\usepackage{color}
\usepackage{graphics}
\usepackage{wrapfig}
\usepackage{enumitem}
\usepackage{mathrsfs}
\usepackage{cite}
\usepackage{environ}
\usepackage{lipsum}
\usepackage{enumitem} 
\usepackage{sectsty}
\usepackage{tcolorbox}
\usepackage{dsfont}
\usepackage[titletoc]{appendix}
\usepackage{mathtools}
\usepackage{bm}
\usepackage{epigraph}
\definecolor{wine-stain}{rgb}{0.7,0,0}
\usepackage[]{hyperref}
\hypersetup{
    colorlinks=true,       
    linkcolor=wine-stain,          
    citecolor=black,        
    filecolor=magenta,      
    urlcolor=red           
}
\usepackage{imakeidx}
\usepackage{tikz}

\title{Multiphasic formulation of Vlasov equations and applications}
\author{Aymeric Baradat\footnote{CNRS, Institut Camille Jordan, Université Claude Bernard Lyon 1, UMR CNRS 5208, Villeurbanne,
France. (\href{mailto:aymeric.baradat@cnrs.fr}{aymeric.baradat@cnrs.fr})}, \, Lucas Ertzbischoff\footnote{CEREMADE, CNRS, Université Paris-Dauphine, PSL Research University, 75016 Paris, France (\href{mailto:ertzbischoff@ceremade.dauphine.fr}{ertzbischoff@ceremade.dauphine.fr})},
\, and Daniel Han-Kwan\footnote{CNRS, Laboratoire de Mathématiques Jean Leray, Nantes Université, UMR CNRS 6629, Nantes, France (\href{mailto:daniel.han-kwan@univ-nantes.fr}{daniel.han-kwan@univ-nantes.fr})}}
\date{}

\newcommand{\R}{\mathbb{R}}
\newcommand{\Z}{\mathbb{Z}}
\newcommand{\T}{\mathbb{T}}

\newcommand{\N}{\mathbb{N}}
\newcommand{\C}{\mathbb{C}}

\newcommand{\F}{\mathcal{F}}
\newcommand{\G}{\mathcal{G}}
\renewcommand{\H}{\mathcal{H}}

\newcommand{\QQ}{\boldsymbol{Q}}
\newcommand{\PHI}{\boldsymbol{\Phi}}

\newcommand{\cg}{\langle}
\newcommand{\cd}{\rangle}
\newcommand{\1}{\mathds{1}}

\newcommand{\pf}{_\#}
\newcommand{\eps}{\varepsilon}
\newcommand{\I}{\mathcal{I}}

\newcommand{\rhorho}{\boldsymbol{\rho}}
\newcommand{\0}{\mathtt{o}}
\newcommand{\vv}{\boldsymbol{v}}
\newcommand{\uu}{\boldsymbol{u}}
\newcommand{\rr}{\boldsymbol{r}}
\newcommand{\ww}{\boldsymbol{w}}
\newcommand{\FF}{\boldsymbol{F}}
\newcommand{\GG}{\boldsymbol{G}}

\newcommand{\DDD}{\mathsf{D}}

\newcommand{\ZZ}{\boldsymbol{Z}}
\newcommand{\YY}{\boldsymbol{Y}}

\newcommand{\lVe}{\left\Vert}
\newcommand{\rVe}{\right\Vert}

\makeatletter

\renewcommand{\chaptermark}[1]{%
  \markboth{%
    \ifnum\c@secnumdepth>\m@ne
      \chaptername\ \thechapter.\ %
    \fi
    #1%
  }{}%
}

\renewcommand{\sectionmark}[1]{%
  \markright{%
    \ifnum\c@secnumdepth>\z@
      \thesection.\ %
    \fi
    #1%
  }%
}

\makeatother

\newcommand{\vertiii}[1]{{\left\vert\kern-0.25ex\left\vert\kern-0.25ex\left\vert #1 
		\right\vert\kern-0.25ex\right\vert\kern-0.25ex\right\vert}}

\newcommand{\enstq}[2]{\left\{#1~\middle|~#2\right\}}
\newcommand{\LRVert}[1]{\left\Vert #1 \right\Vert}

\newcommand{\upinfty}{%
  \raisebox{0.2ex}{\scalebox{0.8}{$\scriptstyle\infty$}}%
}

\theoremstyle{plain}
\newtheorem{Thm}{Theorem}[section]
\newtheorem{Cor}[Thm]{Corollary}
\newtheorem{Prop}[Thm]{Proposition}
\newtheorem{Lem}[Thm]{Lemma}

\theoremstyle{definition}
\newtheorem{Def}[Thm]{Definition}
\newtheorem{Rem}[Thm]{Remark}
\newtheorem{Ex}[Thm]{Example}

\newtheorem{Assx}{Assumption}
\newenvironment{Ass}[1]
 {\renewcommand\theAssx{#1}\Assx}
 {\endAssx}

\numberwithin{equation}{section}

\DeclareMathOperator{\Lip}{Lip}
\DeclareMathOperator{\Div}{div}

\DeclareMathOperator{\Leb}{Leb}

\DeclareMathOperator{\D}{d\!}
\DeclareMathOperator{\Id}{Id}

\DeclareMathOperator{\tr}{tr}
\DeclareMathOperator*{\esssup}{ess\,sup}

\newcommand{\sectionhead}[2][]{%
  \let\savedsectionmark\sectionmark
  \def\sectionmark##1{\savedsectionmark{#1}}%
  \section{#2}%
  \let\sectionmark\savedsectionmark
}

\newcommand{\subsectionhead}[2][]{%
  \let\savedsubsectionmark\subsectionmark
  \def\subsectionmark##1{\savedsubsectionmark{#1}}%
  \subsection{#2}%
  \let\subsectionmark\savedsubsectionmark
}

\makeindex[columns=1, title=Index of Notation, intoc]

\begin{document} 

\frontmatter

	\maketitle

\chapter*{Abstract}
\addcontentsline{toc}{chapter}{Abstract}

This work is a mathematical study of the {\it multiphasic formulation of the Vlasov equation}, which consists in recasting this collisionless kinetic equation as a system of coupled pressureless Euler equations.

This framework allows to consider solutions that are only measure-valued in the velocity variable and is thus relevant to tackle physical problems  where rough velocity distributions, such as Dirac masses, naturally arise.

We specifically address the case of nonlinear Vlasov equations where the force field is one derivative more regular than some moments in velocity of the solution. Under this key assumption, a unified theory of local well-posedness at finite regularity (typically in Sobolev spaces) is developed from scratch for the multiphasic formulation, yielding corresponding results for the Cauchy problem of the associated Vlasov equation at low regularity. 
 We thoroughly apply this abstract theory to two classes of equations, namely Vlasov--Poisson type systems, and Vlasov--Navier--Stokes type systems. In addition to the justification of the monokinetic limit, it also leads to specific applications for each class of equations, allowing possibly rough velocity distributions. For Vlasov--Poisson type systems, we justify the semiclassical limit from Hartree equations, and we describe the nonlinear instability of homogeneous equilibria. For Vlasov--Navier--Stokes type systems,  we establish nonlinear asymptotic stability near monokinetic profiles.

 New results allowing a separation between the regularity in space
 and in velocity are proven, and along the way, we also provide new proofs of known results and generalize them to the case of rough solutions. 
 Finally, the flexibility of the framework is exploited to obtain extensions to more sophisticated systems such as the Vlasov--Poisson equation for ions, or the Vlasov equation coupled with the compressible Navier--Stokes system.

\chapter*{Preface}

\addcontentsline{toc}{chapter}{Preface}

\epigraph{I Contain Multitudes}{Walt Whitman $\&$ Bob Dylan}

This project was undertaken in 2019 to address a gap in the mathematical literature related to Vlasov equations, namely around their multiphasic formulation. Multiphasic formulations have already been well identified as relevant in various physical contexts; most notably (i) in general relativity, the so-called dust solutions have been often used as important special solutions to the Einstein field equations, since the early days of the theory, see e.g. \cite{Friedmann,Lemaitre,Tolman,Oppenheimer,Godel};  (ii) in fluid mechanics and plasma physics, the multiphasic formulation allows to study certain singular limits, see in particular the influential works of Zakharov \cite{Zak},  Brenier \cite{Br-var1,Br97}, and Grenier \cite{Gr96}.
However, a thorough mathematical development of this theory was missing.
Recently, the paper \cite{Bar} by the first author contributed to the emergence of the simple yet compelling idea that it can also help establishing properties of Vlasov equations for data with only measure regularity in velocity. Our initial goal was to elaborate on this idea, for classical equations such as the Vlasov--Poisson system.

As we have tried to put forward an abstract multiphasic framework that allows to treat many Vlasov models at once, this work got longer\footnote{both in length and in time of completion} than usual. Developing such an abstract framework required extra exposition, but we believe that in the end, it brings clarity and unifies results that would otherwise stand separately. In this way, the analysis of the concrete equations also appears more systematic. 

The abstract framework sometimes involves heavy definitions and notation, but the approach is elementary; the analysis requires little prior knowledge (besides basic notions in functional analysis and partial differential equations), and is developed from scratch. This makes the whole work self-contained, and hopefully accessible to many.

 We have chosen to thoroughly address a few aspects of the multiphasic theory: local well-posedness, and stability or instability of certain homogeneous equilibria. Of course, these choices only reflect our personal tastes and limitations; we hope other directions of research will be developed in the future.

This project has benefited from inspiring discussions with several colleagues over many years, among which Claude Bardos, C\'ecile Huneau, Mikaela Iacobelli, Ayman Moussa, Sung-Jin Oh, Tony Salvi, Arthur Touati, and particularly Yann Brenier. It is a pleasure to express our deep gratitude.

\bigskip

Claude Bardos sadly passed away before the completion of this work. He was a guiding light with exceptional dedication, generosity, and benevolence; he was the best model we could have hoped for.
We humbly dedicate this monograph to his memory.

\begin{flushright}
 Lyon, Nantes, Paris \\ 
 \emph{11th September 2026}
 \end{flushright}

	{ \hypersetup{linktoc=page, linkcolor=wine-stain}
		\tableofcontents
	}

    \mainmatter

    \chapter{Introduction}\label{Part-INTRODUCTION}

We start by introducing the basic material for the multiphasic formulation of the Vlasov equation. This introductory chapter mainly serves three purposes. 

\begin{itemize}
    \item First of all, in Section~\ref{Section-motivations}, we  provide some physical context around the nonlinear Vlasov equations we study in this monograph, namely Vlasov--Poisson (for plasmas) and Vlasov--Navier--Stokes (for sprays).
    Then, we explain four kinds of incentives to study solutions to Vlasov equations which are measure-valued in velocity. To conclude this section, we introduce the multiphasic formulation. 
    \item  Second, Section~\ref{Section-bird} presents a rather detailed overview of the main results of the manuscript in order to give the reader a flavour of what is to come.  
    \item Finally, in Section~\ref{sec:equivalence_models}, we justify the equivalence between the kinetic and multiphasic formulations, setting forth a sort of bilingual dictionary allowing to pass from the multiphasic language to the kinetic one, and {\it vice versa}; the precise notions of solution that we consider are also specified.
\end{itemize}

\section{Motivations and framework}\label{Section-motivations}
In this work, we study Vlasov equations which describe the dynamics in the phase space (position $x \in \Omega = \R^d$ or $\T^d$ -- the periodic torus $\R^d/(2\pi\Z)^d$ --  and momentum $v \in \R^d$, with $d \in \N {\setminus \lbrace 0 \rbrace}$) of a system of interacting particles. Broadly speaking, given a force field $F(t,x,v) \in \R^d$, the statistical description of a large system of particles evolving according to the dynamics
\begin{align*}
    \dfrac{\mathrm{d}}{\mathrm{d}t}\mathrm{X} =\mathrm{V}, \ \ \dfrac{\mathrm{d}}{\mathrm{d}t} \mathrm{V}=F(t,X,V), \ \ t>0,
\end{align*}
that corresponds to Newton's laws, can be encoded, in suitable regimes, in the evolution of a single distribution function $f(t,x,v) \in \R^+$. In the so-called mean-field regime (see e.g. \cite{jabin2014review, golse2016dynamics} and references therein), the former should be regarded as a probability measure on the phase space and will formally be a solution to the Vlasov equation 
\begin{equation*}
	\partial_t f + v \cdot \nabla_x f + \mathrm{div}_v \big( F f \big) = 0.
	\end{equation*}
As a matter of fact, the force field $F$ appearing in such equations often depends on the solution $f$ itself, making the whole nonlinear coupling particularly rich and challenging from the mathematical analysis point of view.

In this work, we study the multiphasic formulation of Vlasov equations, with the goal of developing a unified theory of nonlinear Vlasov equations for distribution functions that are only measures with respect to the velocity variable $v \in \R^d$. An important example to have in mind, and that will be motivated in this introduction, is the case where the distribution function is a Dirac mass in velocity -- or a sum of Dirac masses. This program will be achieved for a large class of Vlasov equations for which the force field is related to moments of the solution and, loosely speaking, \emph{gains} one derivative with respect to $x$.

Before going further, let us explain the interest of considering rough (in velocity) solutions to Vlasov equations.
The model systems we have in mind, and that will be discussed in detail later, are of two kinds and will actually guide the development of our abstract framework:

	\begin{itemize}
		
		\item \textbf{The Vlasov--Poisson equation (VP)}.  This is a fundamental model for {\it plasmas} \cite{vlasov1938vibrational,Landau} which is widely studied in the physics literature (see e.g. the reference books \cite{KrallT,Nicholson,chen2015introduction}), and which displays rich mathematical properties (see e.g. \cite{Glassey-kinetic} and references therein). A plasma is loosely speaking an ionized gas, made of ions (particles with positive charge, possibly several families of them) and electrons (particles with negative charge), in which collective effects prevail.
        In this model, the interactions between the particles are described by the Coulomb potential\footnote{This corresponds to the electrostatic approximation; in principle Maxwell equations should prevail for describing electromagnetic effects, but in the non-relativistic limit, they are well-approximated by the Poisson equation.} and result in a scalar potential $U=U(t,x)$ solving the Poisson equation with a source given by the sum of all charge densities. When only electrons are considered while ions are assumed to be immobile, the Vlasov--Poisson equation reads as 
		\begin{equation}
		\label{eq:VP_real}
		\left\{
		\begin{gathered}
		\partial_t f + v \cdot \nabla_x f - \nabla_x U\cdot \nabla_v f = 0,\\
		-\Delta_x U = \int f\D v - \iint f \D v \D x,\\
		f|_{t=0} = f_0.
		\end{gathered}
		\right.
		\end{equation}
        The potential $U$ is usually called the electric potential, while $E\vcentcolon= -\nabla_x U$ stands for the electric field.
The distribution function $f(t,x,v)$ describes the evolution in phase space of electrons in a neutralizing background of ions.
We refer to Chapter \ref{Part2-VP} for more mathematical references.

    Finally, we note that the gravitational analogue of this equation \cite{jeans1915theory}, in which case the interaction is attractive instead of repulsive, is also of great importance.

		\item \textbf{The Vlasov--Navier--Stokes equation (VNS)}. This is a model which is relevant for the description of {\it sprays} (or aerosols), that is to say  suspensions of fine particles or liquid droplets in a surrounding fluid, such as air or another gas. We assume that the fluid is  homogeneous and incompressible, and is described by its divergence-free velocity field $U=U(t,x) \in \R^d$ (and pressure $P=P(t,x) \in \R$ that can be viewed as an associated Lagrange multiplier) solving an incompressible Navier--Stokes equation.
         In this model, the coupling results from the interaction of the particles with the fluid. The Vlasov--Navier--Stokes system is written as follows:
		\begin{equation}
		\label{eq:VNS_real}
		\left\{
		\begin{gathered}
		\partial_t f + v \cdot \nabla_x f +\Div_v\big( (U-v)f\big) = 0,\\
		\partial_t U + (U \cdot \nabla_x) U + \nabla_x P - \Delta_x U = \int (v-U) f \D v ,\\
		\Div_x (U) = 0,\\
		f|_{t=0} = f_0, \quad U|_{t=0} = U_0.
		\end{gathered}
		\right.
		\end{equation}
	The force field  $U-v$ in the Vlasov equation is often referred to as the drag force, while the forcing $\int (v-U) f \D v$ in Navier--Stokes is often called the Brinkman force.
This system can be viewed as an important instance of {\it fluid-kinetic} equations, which are central in the theory of sprays \cite{oro} and \cite{Williams,ReitzBook}. We refer for instance to \cite{Jabin,desvillettes2010-model,BGLM} for introductions to  general classes of fluid-kinetic models describing suspensions of particles in a fluid, the (VNS) system being one of the main prototypes to build on.  We send the reader to Chapter \ref{Part3-VNS} for a more thorough review of the mathematical literature on these models.
	\end{itemize}
For both (VP) and (VNS), it is important to think of the force field in the Vlasov equation as the solution of another equation involving only moments in velocity of $f$: that is, if one introduces the local density and momentum
    \begin{align}\label{def:intro-moment01}
        \rho_f(t,x)\vcentcolon=\int_{\R^d} f(t,x,v) \, \mathrm{d}v, \ \   j_f(t,x)\vcentcolon=\int_{\R^d} v f(t,x,v) \, \mathrm{d}v,
    \end{align}
 then the (VP) system has a force field involving the solution $U=U[\rho_f]$ to the Poisson equation, and the (VNS) system has a force field $U-v$ where $U=U[\rho_f, j_f]$ is the solution to the Navier--Stokes equations. 
The common feature between these two systems that will be crucial in our analysis is that \emph{the force field is one derivative more regular than the moments at play}. Broadly speaking, this is related to the fact that these ``auxiliary'' equations display some regularizing effect, of either of elliptic or parabolic type.

\medskip

Before going further, let us first introduce some notation and a standard functional framework  that will be constantly used throughout this work (several other functional spaces will be introduced later). This will allow us to set some terminology.

\paragraph{Notation and functional spaces.}
The letter $d \in \N {\setminus \lbrace 0 \rbrace}$ refers to the dimension of the physical space. We denote the flat torus by $\T^d=\R /(2\pi \Z)^d$, which we endow with the normalized Lebesgue measure. If not specified by a subscript, all the subsequent functional spaces regarding the $x$ variable are considered on $\T^d$ or $\R^d$ if there is no ambiguity. Similarly, differential operators like $\nabla$, $\mathrm{D}$ or $\Delta$ refer to the $x$ variable. 
If $g=g(x)$ is a function defined on the torus $\T^d$, we often write its average\index{m@$\langle \cdot \rangle$: mean on the torus} as $$\langle g \rangle\vcentcolon=\int_{\T^d}g(x) \, \mathrm{d}x.$$ 
For quantities $A$ and $B$, we write $A\lesssim B$ if $A\le cB$ holds for some universal constant $c>0$. 
If the dependence of $c$ on a parameter $\kappa$ needs to be made explicit, we write $A \lesssim_\kappa B$.

\medskip

The main functional spaces used in this work are the following.

\begin{itemize}

\item Given a measure space $(Y,\nu)$, and writing the variable as $y \in Y$ then, for $p \in [1,\infty)$, the Lebesgue space $L^p_y$ stands for the space of measurable functions such that $\| \varphi \|_{L^p_y} \vcentcolon= \left(\int_Y |\varphi(y)|^p \, \D \nu(y) \right)^{1/p}<+\infty$, while $\| \varphi \|_{L^\infty_{y}} \vcentcolon= \esssup \vert  \varphi \vert$. Given a measurable map $\Theta: (Y,\nu) \to \R_+ $ and $1\leq p<\infty$, we also define the weighted $L^p_{y,\Theta}$ norm as 
\begin{equation*}
    \| \varphi \|_{L^p_{y,\Theta}} \vcentcolon= \left(\int_Y |\varphi(y)|^p  \Theta(y) \, \D \nu(y) \right)^{1/p}<+\infty.
\end{equation*}
If $p=\infty$, we define it as $\| \varphi \|_{L^\infty_{y,\Theta}} \vcentcolon= \esssup \vert \Theta \varphi \vert$. Correspondingly, given a Banach space $(\mathcal B, \| \cdot \|_{\mathcal B})$, we call $L^p_{y,\Theta}( \mathcal B)$ the set of Bochner measurable maps $\Phi: Y \to \mathcal B$ such that
\begin{equation}
\label{eq:def_weighted_Lp}
    \| \Phi \|_{L^p_{y,\Theta}(\mathcal B)} \vcentcolon= \left( \int_Y \| \Phi \|_{\mathcal B}^p  \, \Theta \,  \D \nu\right)^{1/p} < + \infty.
\end{equation}
\index{L@$L^{p}_{y,\Theta}$: $L^p_y$ space with a weight $\Theta$}

\item If $(I, \mu)$ is a measure space carrying the forthcoming relevant multiphasic parameter $\alpha \in I$ (see later in Section \ref{Section-introMultiphaseFormulation}), we will write $L^p_\alpha$ (with $p \in [1, +\infty]$) for the associated Lebesgue space. Functions $\mathrm{g}: I \rightarrow X$ (with values in some space $X$) will be seen as families denoted in bold, that is
\begin{align*}
    \mathbf{g}\vcentcolon=(\mathrm{g}^\alpha)_{\alpha \in I}.
\end{align*}
\index{G@$\mathbf{g}$: multiphasic family $(\mathrm{g}^\alpha)_{\alpha \in I}$}
To ease readability, and depending on the context, we will allow ourselves to alternatively write  $\mathbf{g} \in L^p_\alpha X$ or $g^\alpha \in L^p_\alpha X$, if there is no ambiguity.

\item For $n \in \N$, the Sobolev space $H^n$ on $\Omega= \T^d$ or $\R^d$ is defined as
$$
H^n \vcentcolon= \left\{ \varphi \in L^2(\Omega), \, \| \varphi \|_{H^n} \vcentcolon= \sum_{\ell=0}^n \left(\int_\Omega | \mathrm{D}^\ell \varphi|^2 \, \D x \right)^{1/2}<+\infty \right\}.
$$
\index{H@$H^k$: Sobolev space}
The space $H^{-n}$ is the dual of $H^n$ with respect to the $L^2$ topology. When needed, we will also consider such spaces (or some variants) on $\Omega \times \R^d$.
\item For $k \in \N$ and $\theta \in (0,1]$, the space of $\theta$-Hölder functions on $\Omega= \T^d$ or $\R^d$ is defined as 

$$
\mathscr{C}^{k,\theta}
\vcentcolon=
\left\{
f\in \mathscr{C}^k
, \ \forall |\sigma|=k, \ \ 
\sup_{x\neq y}
\frac{|\mathrm{D}^\sigma f(x)-\mathrm{D}^\sigma f(y)|}{|x-y|^\theta}
<\infty
\right\},
$$
\index{C@$\mathscr{C}^{k, \theta}$: Hölder space}
equipped with the norm
$$
\|f\|_{\mathscr{C}^{k,\theta}}
\vcentcolon=
\sum_{|\sigma|\le k}\|\mathrm{D}^\sigma f\|_{L^\infty}
+
\sum_{|\sigma|=k}
\sup_{x\neq y}
\frac{|\mathrm{D}^\sigma f(x)-\mathrm{D}^\sigma f(y)|}{|x-y|^\theta}.
$$

\item  The space of finite Radon measures on $\R^d$, denoted by $\mathcal{M}(\R^d)$ is defined as the dual of $\mathscr{C}_0(\R^d)$, the space of continuous functions tending to $0$ at infinity. In fact, all our equations on densities preserve mass, and so we can obtain controls in the stronger \emph{narrow} topology, that is, in duality with the space $\mathscr{C}_b(\R^d)$ of continuous and bounded functions. If a sequence $(m_n)$ converges to $m$ narrowly, we write that the convergence holds in $w^*-\mathcal{M}$. In this work, we will often refer to {\bf measure-valued solutions (in velocity)}, meaning that we consider a distribution function $f$ such that $(t,x) \mapsto f(t,x,\cdot) \in \mathcal{M}(\R^d)$.

\item We call $\mathcal P_p(\R^d)$ the set of nonnegative Radon measures  $G$ of mass $1$ having a finite $p$-moment, that is
\begin{equation*}
    \int |v|^p \D G(v) < + \infty.
\end{equation*}
\end{itemize}

Since we will constantly use functions depending both on spatial variable $x$ and on the so-called multiphasic parameter $\alpha \in \I$, we will often highlight the dependency in the variable in the associated norm thanks to a subscript, using for instance $H^n_x$, or $L^p_\alpha H^n_x$.

\bigskip

Let us now explain why a theory of measure-valued solutions (in velocity) can be of interest in the study of systems like (VP) and (VNS). The motivations and results described below will allow us to highlight the main topics of our work.

\subsection{Motivation I: monokinetic and sum of monokinetic data for Vlasov equations.} \label{motivation1}
It is common to consider particular solutions, referred to as {\it monokinetic solutions} to (VP) and (VNS), which correspond to the following ansatz for the distribution function:
    \begin{align}\label{def:monokinsol}
        f(t,x,v) = \rho(t,x) \otimes  \delta_{v=u(t,x)},
    \end{align}
    where $\delta$ denotes the Dirac mass. In the context of Vlasov--Poisson for plasmas, this is often coined as the “cold plasma” approximation since a Dirac mass in velocity can be viewed as a limiting case of a Maxwellian with zero temperature. For the Vlasov--Poisson system \eqref{eq:VP_real}, the ansatz \eqref{def:monokinsol} gives rise to the pressureless Euler--Poisson system
	\begin{equation}
	\label{eq:Euler--Poisson}
	\left\{
	\begin{aligned}
	\partial_t \rho + \Div( \rho u) &=0, \\
	\partial_t u + ( u \cdot \nabla ) u &= -\nabla U,\\
	-\Delta U &= \rho - \int \rho \D x, \\
	\rho |_{t = 0} = \rho_0, \ u |_{t = 0} &= u_0,
	\end{aligned}
	\right.
	\end{equation}
which is a classical system whose smooth dynamics and/or potential blow-up has been the object of several recent works (see e.g. \cite{ELT,loeper2005quasi,LLS,CCTT} or \cite{BCK, song2026h-highMach, EJK} for the related version with massless ions); the attractive version, corresponding to the gravitational case, has also been considered, see e.g. \cite{ERS}. Two-fluid versions of~\eqref{eq:Euler--Poisson} are also relevant, for instance when considering plasmas containing two species, which corresponds to bikinetic data of the form
    \begin{align}\label{def:bikinsol}
        f(t,x,v) = \rho_1(t,x) \otimes \delta_{v=u_1(t,x)} + \rho_2(t,x)\otimes   \delta_{v=u_2(t,x)},
    \end{align}
    and gives rise to the system
	\begin{equation}
	\label{eq:Euler--Poisson-2}
	\left\{
	\begin{aligned}
	\partial_t \rho^\alpha + \Div( \rho^\alpha u^\alpha) &=0, \qquad \qquad \alpha=1,2, \\
	\partial_t u^\alpha + ( u^\alpha \cdot \nabla ) u^\alpha &= - \nabla U,\\
	-\Delta U &= \rho^1 + \rho^2- \int  (\rho^1 + \rho^2) \D x, \\
	\rho^\alpha |_{t = 0} = \rho^\alpha_0, \ u^\alpha |_{t = 0} &= u^\alpha_0,
	\end{aligned}
	\right.
	\end{equation}
	see e.g. \cite{CGG}. Two-fluid models with particles of opposite charge -- electrons and ions -- may also be considered:
		\begin{equation}
	\label{eq:Euler--Poisson-3}
	\left\{
	\begin{aligned}
	\partial_t \rho^\alpha + \Div( \rho^\alpha u^\alpha) &=0, \qquad \qquad \alpha=1,2, \\
	\partial_t u^\alpha + ( u^\alpha \cdot \nabla ) u^\alpha &= (-1)^\alpha \nabla U,\\
	-\Delta U &= \rho^1 - \rho^2, \\
	\rho^\alpha |_{t = 0} = \rho^\alpha_0, \ u^\alpha |_{t = 0} &= u^\alpha_0, \\
	 \int \rho^1_0 \D x &= \int \rho^2_0 \D x,
	\end{aligned}
	\right.
	\end{equation}
	(or possibly more complicated data). 
For the Vlasov--Navier--Stokes system \eqref{eq:VNS_real}, the monokinetic ansatz \eqref{def:monokinsol} gives rise to the pressureless Euler--Navier--Stokes system
	\begin{equation}
	\label{eq:Euler-NS}
	\left\{
	\begin{aligned}
	\partial_t \rho + \Div( \rho u) &=0, \\
	\partial_t u + ( u \cdot \nabla ) u &= U-u,\\
		\partial_t U + (U \cdot \nabla) U + \nabla P - \Delta U &= \rho (u - U),\\
		\Div(U) &= 0, \\
	\rho |_{t = 0} = \rho_0, \ u |_{t = 0} = u_0, \ U |_{t=0} &= U_0,
	\end{aligned}
	\right.
	\end{equation}
	which has also been the focus of several recent works (see e.g. \cite{ChoiJungR3,HuangTangZouR3,li2025global, lemarié2025, Danchin-ENS2026}).
	
	Finally, one may consider $N$-kinetic distributions, combined with a regular kinetic part, which corresponds to mixed data of the form
  \begin{align}\label{def:Nkinsol}
        f(t,x,v) = \sum_{k=1}^N \rho_k(t,x)  \otimes \delta_{v=u_k(t,x)} + g(t,x,v),
    \end{align}
    possibly with $N = +\infty$.
    We aim at developing a general Cauchy theory associated with such data (and more general ones), both for (VP) and (VNS).  In particular, this allows one to properly embed the pressureless Euler dynamics into the genuine Vlasov dynamics. A closely related question is the so-called monokinetic limit, which arises when one considers a sequence of (regular) initial data converging to a Dirac mass in velocity: is the limiting dynamics governed by pressureless Euler?

  \subsection{Motivation II: instability of rough homogeneous equilibria for (VP)--type equations.} One can observe that any spatially homogeneous profile $\mu=\mu(v)$, that depends only on $v$, normalized so that $\int \D \mu = 1$, is a stationary solution to (VP). For an equilibrium $\mu(v)$ localized around two different velocities $v_1 \in \R^d$ and $v_2 \in \R^d$, there exists a well-known instability mechanism for (VP), often referred to as \textit{two-stream instability} (see \cite{BohmGross,PierceH,roberts1967nonlinear,KrallT,Davidson}). It appears when two beams of electrons in plasmas have very different velocities and is actually a manifestation of the celebrated \textit{Penrose criterion} for (spectral) instability \cite{Pen60}, which will be precisely stated later on in the monograph.

    In the mathematical literature, this instability has been, in particular, studied in \cite{GS, HKH, HKNsima} for smooth profiles $\mu(v)$. Let us state two results about the nonlinear instability in the sense of Lyapunov of such profiles.
    \begin{Thm}[from \cite{HKH,HKNsima}, periodic case $\Omega=\T^d$]
        \label{eq:insta-kin-1}
Assume that $\mu(v)$ is smooth and rapidly decaying. If $\mu$ satisfies the Penrose instability criterion, then the  Vlasov--Poisson system \eqref{eq:VP_real} is nonlinearly unstable about $\mu$. That is, for all $s >0$, there exists $\varepsilon_0>0$ such that, for all $\delta>0$, there exists a solution $f$ to~\eqref{eq:VP_real} such that $\| f_0- \mu \|_{H^s_{x,v}} \leq \delta$ but for $T_\delta =O (|\log \delta|)$, 
    $$
    \| f(T_\delta)- \mu \|_{H^{-s}_{x,v}} \geq \varepsilon_0. 
    $$
    \end{Thm}
   Typical equilibria satisfying the required assumptions consist of a sum of two separated Maxwellians centred at $v_1 \in \R^d$ and $v_2 \in \R^d$, with $|v_1- v_2|$ sufficiently large.
    
Instability has also been investigated for the one-dimensional two-fluid pressureless Euler--Poisson equation~\eqref{eq:Euler--Poisson-2} in \cite{CGG}. In this case, the equilibrium is given by the constant $(a_1,v_1,a_2,v_2)$.
 \begin{Thm}[from \cite{CGG}, periodic case $\Omega=\T^d$]
    \label{eq:insta-kin-2} Let $d=1$.
    If the equilibrium $(a_1,v_1,a_2,v_2)$ is linearly unstable\footnote{We shall not specify here the exact meaning of this notion, but this is an explicit condition which should be viewed as the Penrose instability criterion at the level of two fluids. See the upcoming Section~\ref{sec:insta-appli} of  Chapter~\ref{Part2-VP} for details.}, then the system~\eqref{eq:Euler--Poisson-2} is also nonlinearly unstable about it. 
    That is, for all $s >0$, there exists $\varepsilon_0>0$ such that, for all $\delta>0$, there exists a solution $(\rho^1,u^1,\rho^2,u^2)$ to~\eqref{eq:Euler--Poisson-2} such that $\| (\rho^1_0-a_1,u^1_0-v_1,\rho^2_0-a_2,u^2_0-v_2)\|_{H^s_x} \leq \delta$ but for $T_\delta =O (|\log \delta|)$, 
    $$
    \| (\rho^1(T_\delta)-a_1,u^1(T_\delta)-v_1,\rho^2(T_\delta)-a_2,u^2(T_\delta)-v_2)\|_{L^1_x\cap L^{\upinfty}_x} \geq \varepsilon_0. 
    $$
    \end{Thm}

In this work, we shall tackle the analysis of instabilities of \textit{rough profiles in velocity} for systems behaving like (VP). Roughly speaking, the aim is to extend Theorem~\eqref{eq:insta-kin-1} for equilibria that are only measures,  for instance 
\begin{align*}
    \mu(v)=a_1 \otimes \delta_{v=v_1}+  a_2  \otimes \delta_{v=v_2}, \quad a_1,a_2 \geq 0,
\end{align*}
which can be interpreted as the zero-temperature limit of the sum of two Maxwellians centred at $v_1$ and $v_2$.  This idea originates from the work \cite{Bar}, which studied the \emph{ill-posedness of singular Vlasov equations} around rough equilibria and also obtain \emph{almost-Lyapunov} instability results for the Vlasov--Poisson system.

We shall recover Theorem~\ref{eq:insta-kin-2} as part of a general statement, for a class of equations including~\eqref{eq:Euler--Poisson-2}.  
 We will therefore unify the results of  \cite{GS, HKH, HKNsima} for Vlasov--Poisson and  of \cite{CGG} for two-fluid pressureless Euler--Poisson.

For physical relevance, it is desirable that the unstable distribution function $f$, as obtained in Theorem~\ref{eq:insta-kin-1}, is nonnegative. This property is ensured in \cite{HKH} by imposing a technical condition, the so-called $\delta$ (or $\delta'$) condition, that is satisfied for many natural examples. We also aim at removing this condition.

Finally, in another direction, some interesting variants of the Vlasov--Poisson equation~\eqref{eq:VP_real} involve additional nonlinearities at the level of the Poisson equation: we have in particular in mind the Vlasov--Poisson equation for ions and the Vlasov--Monge--Amp\`ere system (which we will soon enough introduce and discuss in greater detail). We also aim at obtaining general instability statements taking into account these nonlinearities, which shall ultimately lead to a general unified treatment of the aforementioned Vlasov--Poisson type systems.

   \subsection{Motivation III: concentration to monokinetic profiles for (VNS)--type equations.}
   For equations like (VNS), there is no regular homogeneous equilibrium, because of the friction term in the drag force $U-v$.  As a matter of fact, concentration in velocity is expected as $t\to +\infty$, that is to say, the distribution function $f$ should converge to a Dirac mass in velocity.
   Interestingly, such Dirac masses almost cancel out the Brinkman force, that is the forcing in Navier--Stokes, provided that they are centred at a velocity close to the fluid velocity $U$;
    therefore, thanks to the dissipation in (NS), one may also expect the velocity field $U$ to reach an equilibrium.
    Concretely, one thus expects
    \begin{align}\label{def:alignmentVNS}
        f(t) - \rho^\infty(t,x) \otimes \delta_{v=v^\infty}\underset{t \rightarrow + \infty}{\rightharpoonup} 0, \ \ U(t) \underset{t \rightarrow + \infty}{\rightarrow} v^\infty, \ \ \text{for some} \ \ v^\infty \in \R^d,
    \end{align}
    and for some asymptotic profile $\rho^\infty(t,x)$. This reflects the alignment of all phases to a common velocity $v^\infty$.

    After the works \cite{Jab} that treated a different Vlasov--Stokes equation and \cite{ChKw} which put forward a conditional result of convergence to equilibrium for (VNS), this has been fully justified in a close to equilibrium regime in \cite{HKMM} on the torus and in \cite{HK-VNSR3} in the whole space (see also \cite{danchin2024fujita,danchin2025-TORUS2d}). The main result obtained in \cite{HKMM} can be stated as follows.
    
    \begin{Thm}[from \cite{HKMM} -- $\T^3$ case]
    \label{thm:VNSlargetime}
    Let $(U_0,f_0)$ be smooth initial data of finite kinetic energy such that $f_0$ is decaying sufficiently fast in velocity (uniformly in $x$) and is small enough in the sense that 
\begin{align}
\label{eq:closetoeq-VNS}
\mathscr{E}(U_0,f_0)+\Vert U_0 \Vert_{\dot{H}^{1/2}_x} \ll 1,
\end{align}
where $\mathscr{E}$ is a modulated energy functional defined as
\begin{align}\label{def:modulated-energyVNS}
\begin{aligned}
    \mathscr{E}(U(t), f(t))&\vcentcolon=\dfrac{1}{2}\int_{\T^3 \times \R^3} f(t,x,v) \vert v-\langle j_f(t) \rangle\vert^2 \, \mathrm{d} v \, \mathrm{d} x \\
&\quad+  \dfrac{1}{2}\int_{\T^3} \vert U(t,x)- \langle  U(t) \rangle \vert^2 \, \mathrm{d} x
+ \dfrac{1}{4} \vert \langle j_f(t) \rangle - \langle U(t) \rangle \vert ^2.
\end{aligned}
\end{align}
Let $(U,f)$ be a global weak solution  of \eqref{eq:VNS_real} in the sense that
\begin{equation}
    \label{def:energyspaceVNS}
    \begin{aligned}
    &U \in L^{\infty}(\R^+; L^2_x) \cap L^2(\R^+;H^1_x), \\ 
    &f \in  L^\infty(\R^+; L^1 \cap L^\infty_{x,v}), \ \  \vert v \vert^2 f \in L^{\infty}(\R^+; L^1_{x,v}),
     \end{aligned}
\end{equation}
and with initial data $(U_0,f_0)$. Then there exists an asymptotic profile $\rho^{\infty} \in L^{\infty}_x$ and $\lambda>0$  such that for all $t \geq 0$
\begin{align*}
\mathrm{W}_1 \Big( f(t) , \rho^{\infty}\left(x-t\mathcal{W}_\infty\right) \otimes \delta_{v=\mathcal{W}_\infty} \Big) + \left\Vert u(t) - \mathcal{W}_\infty \right\Vert_{L^2_x}  \lesssim \mathscr{E}(U_0,f_0)^{1/2}e^{-\lambda t},
\end{align*}
where $\mathcal{W}_\infty=\frac{\langle j_{f_0}+U_0 \rangle}{2}$ and where $\mathrm{W}_1$ is the $1$-Wasserstein distance on $\T^3 \times \R^3$.
    \end{Thm}

\index{W1@$\mathrm{W}_1$: $1$-Wasserstein distance}Let us recall that a 
definition of the $1-$ Wasserstein distance on a complete metric space $X$ (with the so-called Monge--Kantorovitch duality) is the following (see e.g. \cite{Santambrogio}) : if $\mu$ and $\nu$ are two measures on $X$ of same finite mass and with finite first-order moment, we have
\begin{align}\label{def:Wasserstein1}
\mathrm{W}_1(\mu,\nu)=\sup \enstq{ \left\vert \int_X \psi(x)\mathrm{d}\mu(x) - \int_X \psi(x) \mathrm{d}\nu(x) \right\vert} { \Vert \nabla \psi \Vert_{L^\infty(X)} \leq 1}.
\end{align}
This distance metrizes the weak convergence of measures.

\medskip

      Global in time weak (more precisely, Leray type) solutions of \eqref{eq:VNS_real} were for instance shown to exist in \cite{BDGM}.
   All aforementioned works \cite{HKMM, HK-VNSR3, danchin2024fujita, danchin2025-TORUS2d} 
   concerning the long-time asymptotics of \eqref{eq:VNS_real} rely on a remarkable modulated energy--dissipation identity, that was uncovered in \cite{ChKw}. Interestingly, the monokinetic behaviour is then encoded in the convergence to $0$ of the modulated energy $\mathscr{E}$ defined in \eqref{def:modulated-energyVNS}, and Assumption~\eqref{eq:closetoeq-VNS} corresponds to a close to equilibrium regime. 
    One may wonder if a more direct perturbative analysis (around Dirac masses) could be possible. This seems out of reach for  solutions belonging to the natural energy space \eqref{def:energyspaceVNS}.
    As a matter of fact, it even seems unclear how one can linearize \eqref{eq:VNS_real} around a Dirac mass in velocity.
    
   \begin{Rem}
       We can also refer to \cite{GHKM, EHKM, E1} for other asymptotics regarding different geometries or forcing terms that can, in some instances, prevent  concentration in velocity.
   \end{Rem} 
   
    One of the goals pursued in this work will be to offer a new perspective on asymptotic results for (VNS): in particular, we will develop a framework in which a perturbative analysis of monokinetic profiles is allowed.

\bigskip

To summarize Motivations I--II--III, nonlinear Vlasov equations such as (VP) and (VNS) raise challenging questions related to the 
dynamics of profiles that are rough in velocity.   From a broader point of view, we aim at decoupling the required regularity in space and velocity, namely allowing measure regularity in velocity while considering smoother dynamics in space. To study such singular behaviours, we will reformulate the Vlasov equations in the so-called multiphasic framework.

\subsection{Motivation IV: semiclassical limit of the Schr\"odinger--Poisson equation for mixed states}
\label{sec:semi-intro}

The last motivation comes from the derivation of Vlasov--Poisson, starting from quantum mechanics.
Consider the semiclassical Schr\"odinger--Poisson equation set in $\R^d$:
\begin{equation}
\label{eq:SPintro}
	\left\{
	\begin{gathered}
i \varepsilon \partial_t \Psi^\eps + \frac{\varepsilon^2}{2} \Delta_x \Psi^\eps =   V^\eps \Psi^\eps, \\
-\Delta_x V^\eps =  |\Psi^\eps|^2, \\
\Psi^\eps|_{t=0} =  \Psi^\eps_0,
\end{gathered}
\right.
\end{equation}
where the small parameter $\eps>0$ is the scaled Planck constant. In the semiclassical limit $\eps \to 0$, one expects to recover classical mechanics.
The case of (monophasic) WKB data, that is to say initial data of the form
\begin{equation}
\label{eq:monophasicquant0intro}
\Psi^\eps_0 (x)= \sqrt{\rho_0^{\eps} (x)} \exp \left( \frac{i}{\eps} S^{\eps}_0(x) \right),
\end{equation}
for this equation and for other general nonlinear Schr\"odinger has been abundantly studied in the mathematical literature, see in particular \cite{Ger93,Gr98,ZZM,Zhang03,AC,CarlesM,Mas}. The leading behaviour in the semiclassical limit $\eps \to 0$ of \eqref {eq:SPintro} is given by a solution to the pressureless Euler--Poisson system
\begin{equation*}
	\left\{
	\begin{gathered}
\partial_t \rho + \Div (\rho u ) =0, \\
\partial_t u + u \cdot \nabla u = - \nabla V, \\
-\Delta V =  \rho, \\
\rho|_{t=0}= \lim_{\eps \to 0}\rho_0^{\eps}, \quad u|_{t=0}= \lim_{\eps \to 0}\nabla S_0^{\eps}.
\end{gathered}
\right.
\end{equation*}
Let us study the mixed states version of the Schr\"odinger--Poisson equation~\eqref{eq:SPintro}, which corresponds to the Hartree--Poisson equation for operators:
\begin{equation}
\label{eq:SPintro-hartree}
	\left\{
	\begin{aligned}
i \varepsilon \partial_t \gamma^\eps &=\left[ -\frac{\varepsilon^2}{2} \Delta_x + V_\eps,\gamma^\eps\right], \\
-\Delta_x V^\eps &=  \rho^\eps, \\
\gamma^\eps|_{t=0} &=  \gamma^\eps_0,
\end{aligned}
\right.
\end{equation}
where $\gamma^\eps \in \mathscr{L}(L^2(\R^d))$ is a nonnegative trace-class operator and where $[A, B]=AB-BA$. In \eqref{eq:SPintro-hartree}, $\rho^\eps$ stands for the density associated with $\gamma^\eps$; denoting by $\gamma^\eps(t,x,y)$ its Schwartz kernel, there holds 
$\rho^\eps(t,x)= \gamma^\eps(t,x,x)$.

Consider that the initial condition writes of the form
\begin{equation}
\label{eq:init-hartree}
\gamma^\eps_0(x,y) = \int_I  \Psi_0^{\eps,\alpha} (x) \overline{\Psi_0^{\eps,\alpha}}(y) \, \D \mu(\alpha),
\end{equation}
where $(I, \mu)$ is a given measure space (see for instance \cite{BB2,deSuzzoni}). This generalizes the decomposition given by the spectral theorem, in which case $I= \N$. The work \cite{BB2} proposed to assume that  each $\Psi^{\eps,\alpha}_0$ is a WKB datum, that is to say 
\begin{equation}
\label{eq:multiphasicquant0intro-hartree}
\Psi^{\eps,\alpha}_0 (x)= \sqrt{\rho_0^{\eps,\alpha} (x)} \exp \left( \frac{i}{\eps} S^{\eps,\alpha}_0(x) \right).
\end{equation}
In this framework, for the cubic nonlinearity (i.e. when $V^\eps = \rho^\eps$), \cite{BB2} justified the semiclassical limit towards the so-called Vlasov--Benney equation, for data with analytic regularity in space. This strong regularity assumption appears necessary as the Vlasov--Benney is a singular equation that is in general ill-posed in Sobolev regularity \cite{HKN,Bar}.

On the other hand, as far as we know, this problem has not been properly explored with a Poisson nonlinearity as in \eqref{eq:SPintro-hartree}.  In this case, results with a finite regularity assumption are expected.
This would lead to a derivation of Vlasov--Poisson for a class of initial data different from those of \cite{LP,LS}.
Another motivation of this work will be to provide a robust framework for investigating this question and therefore justifying this semiclassical limit.

\subsection{The multiphasic formulation}\label{Section-introMultiphaseFormulation}

	Our purpose is to give a systematic and unified description of nonlinear Vlasov equations, seen as \emph{multiphasic fluid systems}. 
	The underlying starting point consists in interpreting the distribution function $f$ as a superposition of (possibly infinitely many) different monokinetic profiles.
	Concretely, let us consider an abstract Vlasov equation with a given force field $F=F(t,x,v)$, that is
	\begin{equation}
	\label{eq:Vlasov_general_intro}
	\partial_t f + v \cdot \nabla_x f + \mathrm{div}_v \big( F f \big) = 0.
	\end{equation}
	The idea is that the equation \eqref{eq:Vlasov_general_intro} is formally equivalent to a possibly infinite system of \textit{coupled pressureless Euler equations}. More precisely, let us introduce an abstract parameter $\alpha$ in some measure space $(I,\mu)$ and consider the following system describing the evolution of a family of \textbf{densities} $(\rho^\alpha)_{\alpha \in I} = (\rho^\alpha(t,x))_{\alpha \in I} \in \R_+^I$ and \textbf{velocity fields} $(v^\alpha)_{\alpha \in I} = (v^\alpha(t,x))_{\alpha \in I}\in (\R^d)^I$ indexed by $\alpha$:
	\begin{equation}
	\label{eq:general_system_intro}
	\left\{
	\begin{gathered}
	\partial_t \rho^\alpha + \Div( \rho^\alpha v^\alpha) =0, \\
	\partial_t v^\alpha + ( v^\alpha \cdot \nabla ) v^\alpha = F^\alpha,\\
	\rho^\alpha |_{t = 0} = \rho^\alpha_0, \ v^\alpha |_{t = 0} = v^\alpha_0.
	\end{gathered}
	\right.
	\end{equation}
	Above, since the original force field  $F$ is fixed, $F^\alpha$ is given by $F^\alpha(t,x)\vcentcolon = F(t,x, v^\alpha(t,x))$,
 which appears as a forcing term in the momentum equation. However, in view of application to genuine kinetic equations, we will consider more general forces $F^\alpha$, see \eqref{eq:link_force_Vlasov_multiphase} below.
 
	The key observation is that given a \emph{smooth} solution $(\rho^\alpha,v^\alpha)_{\alpha \in I}$ of \eqref{eq:general_system_intro}, the distribution function
	\begin{equation}
	\label{eq:f_intro}
	f(t,x,v) \vcentcolon= \int_I \rho^\alpha (t,x) \otimes  \delta_{v= v^\alpha(t,x)} \, \D \mu(\alpha),
	\end{equation}
	satisfies the Vlasov equation~\eqref{eq:Vlasov_general_intro} in the sense of distributions, with the initial condition 
	\begin{equation*}
	f|_{t=0}(x,v) \vcentcolon= \int_I \rho^\alpha_0(x) \otimes \delta_{v= v^\alpha_0(x)} \, \D \mu(\alpha).
	\end{equation*}
    Here, the former decomposition means that for every nice test function $\psi=\psi(v)$, one has
    \begin{align*}
        \langle f(t,x, \cdot), \psi \rangle=\int_{I} \rho^\alpha(t,x) \psi(v^\alpha(t,x)) \,  \D \mu(\alpha),
    \end{align*}
    and we can formally recover the local density and momentum associated with $f$ (see \eqref{def:intro-moment01}) as
    \begin{align*}
        \rho_f(t,x)=\int_{I} \rho^\alpha(t,x) \,  \D \mu(\alpha), \ \ j_f(t,x)=\int_{I} \rho^\alpha(t,x) v^\alpha(t,x) \,  \D \mu(\alpha).
    \end{align*}
 We refer to~\eqref{eq:general_system_intro}--\eqref{eq:f_intro} as the {\bf multiphasic formulation} of the Vlasov equation~\eqref{eq:Vlasov_general_intro}; in short, this can be interpreted as a parametrization of the Vlasov equation, which allows us to recast it as a system of pressureless Euler equations coupled through the forcing. As seen previously for $\mathrm{(VP)}$ and $\mathrm{(VNS)}$, we would actually like the force $F$ in the Vlasov equation to depend on the solution $f$ itself, through moments in velocity of $f$: therefore we will consider general forces of the form
	\begin{equation}
	\label{eq:link_force_Vlasov_multiphase}
	F^\alpha\vcentcolon = \mathrm{F}^\alpha \left( (\rho^\alpha)_{\alpha \in I},(v^\alpha)_{\alpha \in I} \right)
	\end{equation}
	for families of mapping $(\mathrm{F}^\alpha)_{\alpha \in I}$ that will satisfy suitable structural assumptions.

 Let us also observe that the case of monokinetic solutions \eqref{def:monokinsol} corresponds to $\rho^\alpha=\rho$ and $v^\alpha=v$ for $\mu$-almost any $\alpha \in I$ (assuming that $(I, \mu)$ is a probability space).  The multiphasic formulation allows to consider several other natural types of solutions; before further generalization, let us mention 
 \begin{itemize}
     \item the case of a sum of $n$ 
     Dirac masses in velocity, that is 
     $$f=\sum_{j=1}^n m_j \otimes  \delta_{v=v_j}$$ 
     for some given set of positive masses $(m_j)_{1 \leq j \leq n}$ and velocities $(v_j)_{1 \leq j \leq n}$. Here, one can take $$I=[\![1, n]\!], \quad \mathrm{d}\mu(\alpha)=\frac{1}{n}\sum_{j=1}^n \delta_{\alpha=j},$$ 
     with $\rho^\alpha(x)=n m_\alpha(x)$ and $v^\alpha(x)=v_\alpha(x)$. 
     \item the case of a smooth function $f=f(x,v)$, where one can take $$I=\R^d, \quad \mathrm{d}\mu(\alpha)= \Psi(\alpha) \mathrm{d}\alpha,$$ 
     for some smooth integrable distribution $\Psi(\alpha)$, with $\rho^\alpha(x)=f(x,\alpha)/\Psi(\alpha)$ and $v^\alpha(x)=\alpha$.
 \end{itemize}
 A more general study of this representation will be carried out in the subsequent Section \ref{sec:equivalence_models}.

 \medskip

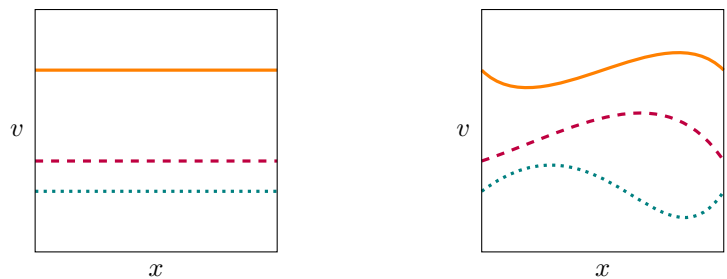
\begin{figure}
\centering
\begin{tikzpicture}[scale=0.8]
\draw[very thick, dotted, teal] (0,1) -- (4,1);
\draw[very thick, dashed, purple] (0,1.5) -- (4,1.5);
\draw[very thick, orange] (0,3) -- (4,3);
\draw (0,0) rectangle (4,4);
\draw (-0.3,2) node{$v$};
\draw (2,-0.3) node{$x$};
\end{tikzpicture}
\hspace{2cm}
\begin{tikzpicture}[scale=0.8]
\draw[very thick, dotted, teal] (0,1) .. controls (2,2.5) and (3,-0.5) .. (4,1);
\draw[very thick, dashed, purple] (0,1.5) .. controls (1.5,2) and (3,3) .. (4,1.5);
\draw[very thick, orange] (0,3) .. controls (1,2) and (3,4) .. (4,3);
\draw (0,0) rectangle (4,4);
\draw (-0.3,2) node{$v$};
\draw (2,-0.3) node{$x$};
\end{tikzpicture}
\caption{\label{figmultiphasic} Two instances of representations of multiphasic distributions in the case when $\mu$ is the sum of three Dirac masses, for instance on $I = \{1,2,3\}$. In fact, in these representations, we draw the graphs of the velocities $v^1$, $v^2$ and $v^3$. Each graph is associated a density, respectively $\rho^1$, $\rho^2$ and $\rho^3$, but they are not represented here. On the left, the velocities are constant in space. Along the evolution, they can evolve and give rise to distributions such as the one on the right. With this kind of representations, it appears clear that the multiphase formulation blows-up as soon as one slope on one graph becomes infinite. Also, let us mention that superposing infinitely many such graphs can perfectly give rise through~\eqref{eq:f_intro} to a distribution $f$ that is smooth in space \emph{and} velocity.}
\end{figure}

Writing this ansatz is a common approach within the Russian school of plasma physics or fluid mechanics (see e.g. \cite{ChesnoPavlov, Pavlov, FeraPavlov}\footnote{These references also refer to the Russian-language book \cite{silin} which we have not been able to access.}). The multiphasic formulation had notably remarkable applications in the early works of Zakharov about the Benney equation \cite{Zak, zakharov1981benney}\footnote{Interestingly,Zakharov's was first to rephrase  the Benney equation, i.e. the free-surface hydrostatic Euler equation, as a superposition of fluid layers as in \eqref{eq:general_system_intro}, and then to connect it to a Vlasov equation, namely the so-called Vlasov--Benney equation. A similar fruitful procedure to study stratified fluids comes from using the so-called \textit{isopycnal coordinates} (see e.g. \cite{bianchini2024hydrostatic, fradin2024well, bianchini-ertz2024review} for more in this direction).}. We may also refer to \cite{BB2,BB3} for recent work in this direction.

From the mathematical point view, the multiphasic formulation turns out well adapted to handle solutions which only have the regularity of a measure with respect to velocity.
In \cite{Gr96}, Grenier used this formulation to provide the first justification of the quasineutral limit of the Vlasov--Poisson system for data with analytic regularity with respect to $x$. On the other hand, it allows for distribution functions that are only measures (and possibly compactly supported) in the velocity variable. This approach was recently extended in \cite{GIRR} to the relativistic Vlasov--Maxwell system.  Alongside \cite{Zak, zakharov1981benney}, the multiphasic formulation is also useful in the theory of relaxed incompressible Euler equations, \emph{\`a la Brenier} (aka incompressible optimal transport), see \cite{Br-var1,Br-var2,Br97,Br-var3,Brehyqua,BreMoy}. Finally, it appeared to be fruitful to prove ill-posedness for singular Vlasov equations (such as the Vlasov--Benney equation) with no strong regularity assumption with respect to velocity, in the work \cite{Bar} by the first author, thus generalizing \cite{HKN}.

A common feature of \cite{Gr96,Bar,GIRR} is that they solve the corresponding initial-value problems~\eqref{eq:general_system_intro} for data with analytic regularity, invoking a Cauchy--Kowalevskaya theorem.

\medskip

We adopt in this work a different perspective, focusing on finite regularity -- typically Sobolev -- solutions to the multiphasic systems.
The price to pay is that, contrary to \cite{Gr96,Bar,GIRR}, we will not be able to study systems displaying a potential loss of derivatives; as a matter of fact, as already mentioned, we even require a gain of one derivative at the level of the force.
Moreover, in order to relate back with Vlasov equations, we shall only consider smooth solutions to \eqref{eq:general_system_intro}; as a result, shocks may occur, yielding a potential blow-up of the solutions.
However, shocks have in general to be understood as a failure of the multiphasic parametrization, and not as a true blow-up of the solution to the original Vlasov equation. Because of this aspect, the formulation~\eqref{eq:general_system_intro}--\eqref{eq:f_intro} may seem to be, at first sight, relevant only for local in time results. Another motivation of our work is to uncover situations where long-time results can be nevertheless obtained.

\begin{Rem}
	At this stage of the introduction, it is worth mentioning another popular parametrization of the Vlasov equation (in dimension $d=1$) that seems unrelated to the multiphasic representation: the so-called \textbf{water-bag formulation}  (see e.g. \cite{BF, NBwater, BDBF, Besse2009multi, besse2011waterbag, Bardos2013cauchy} and references therein), which is well-known in plasma physics and is convenient for numerical purposes, consists in writing the distribution function as a sum of indicator functions over level sets of velocity profiles. One often expects that they satisfy a system of hyperbolic type. More concretely, the ansatz for the kinetic distribution is the following:
    \begin{align*}
        f(t,x,v)=\sum_{i=1}^N \mathcal{A}_j \mathbf{1}_{u_j^-(t,x) < v < u_j^+(t,x)},
    \end{align*}
    for a fixed set of constant weights $\mathcal{A}_j\in \R$, and unknown velocities $u_j(t,x) \in \R$.
    For simplicity, let us illustrate this approach on the 1D Vlasov--Poisson system \eqref{eq:VP_real} set on $\T \times \R$ (in that case, the former ansatz can also be refered to as a superposition of electron layers, by analogy with vortex patches for the 2D incompressible Euler equations, see e.g. \cite{dziurzynski1987patches,MB,Roulley}). We obtain the following system for the velocities
    \begin{align*}
        \partial_t u_j^\pm + u_j^\pm \partial_x u_j^\pm&=-\partial_x U, \\
        -\partial_{xx} U&= \sum_{j=1}^N \mathcal{A}_j(u_j^+-u_j^-)-1.
    \end{align*}
    By setting $\rho_j\vcentcolon=\mathcal{A}_j(u_j^+-u_j^-)$ and $v_j\vcentcolon=(u_j^+ + u_j^-)/2$, we actually end up with the following coupled Euler--Poisson systems
    \begin{align*}
    \partial_t \rho_j + \partial_x (u_j \rho_j)&=0, \\
        \partial_t v_j + v_j \partial_x v_j+\rho_j^{-1} \partial_x P(\rho_j)&=-\partial_x U, \\
        -\partial_{xx} U= \sum_{j=1}^N \rho_j-1, \ \ P(\rho_j)&=\frac{\rho_j^3}{12 \mathcal{A}_j^2}.
    \end{align*}
    Interestingly, in contrast to the pressureless Euler systems \eqref{eq:general_system_intro} from the multiphase formulation set forward previously, the water-bag approach yields an additional cubic pressure (and is in particular hyperbolic).

    Let us mention that the water-bag parametrization seems to be, in full generality, restricted to Vlasov equations written in non-conservative form (in $(x,v)$).
\end{Rem}

\begin{Rem} Let us conclude by briefly mentioning two other contexts where the multiphasic formulation (or an avatar thereof) turns out relevant.

\begin{itemize}
    \item In general relativity, the exact multiphasic ansatz~\eqref{eq:f_intro} applied to solutions of the massless Einstein---Vlasov system gives rise to the so-called Einstein-null dust system \cite{Friedmann,Lemaitre,Tolman,Oppenheimer,Godel}.
  As a matter of fact, the multiphasic formulation was recently used for the resolution of the following two celebrated conjectures:
  the AdS instability conjecture in \cite{Moschidis,Moschidis2}, and  the Burnett (or reverse Burnett) conjecture about high-frequency limits of the Einstein vacuum equations, in \cite{HuneauL0,luk2020high,HuneauL1, HuneauL2,Arthur2,huneau2025}.

    \item  In quantum dynamics,  several avatars of the multiphasic representation are actually relevant.
       In the case of mixed states, the analogue of the Vlasov equation~\eqref{eq:Vlasov_general_intro} is the Hartree equation, and the multiphasic formulation is very reminiscent of what is obtained by the spectral theorem: that is, when the Hartree equation is recast as a system of coupled Schrödinger equations (see also the random-field formulation from \cite{deSuzzoni,CDSM}). 
The analogy becomes even stronger in the semiclassical limit, when considering WKB data, see Motivation IV in Subsection~\ref{sec:semi-intro}. 
In the case of pure states, that is when considering a single Schr\"odinger equation, justifying a \emph{multiphase} WKB approximation in the semiclassical limit involves another system of equations that can also be viewed as quantum analogues of the multiphasic systems considered in this work.

\end{itemize}
    \end{Rem}

\section{Overview of the results and structure of the monograph}
\label{Section-bird}

In this section, we provide a bird's-eye view of the monograph by highlighting some results and statements (but the description will not be comprehensive).
We then shall conclude this chapter with Section~\ref{sec:equivalence_models} where we explain the equivalence between the multiphasic and kinetic formulation.

The main body of this work is then divided into three chapters.

\begin{itemize}

\item In Chapter \ref{Part1-LWP}, we provide  abstract local in time well-posedness results for a broad class of multiphasic systems, for finite regularity data. 
  In Chapters \ref{Part2-VP} and \ref{Part3-VNS}, 
we thoroughly apply these general results to Vlasov--Poisson and Vlasov--Navier--Stokes type systems. 

\item In Chapter  \ref{Part2-VP}, dedicated to Vlasov--Poisson systems, we also study two other problems, namely the semiclassical limit of Schr\"odinger--Poisson systems and the nonlinear instability of homogeneous equilibria, both considered through the prism of the multiphasic formulation.
 
\item Chapter  \ref{Part3-VNS} is dedicated to the application for Vlasov--Navier--Stokes systems. We also study the problem of long-time dynamics around constant monokinetic equilibria, again with a multiphasic point of view.

\end{itemize}

It is precisely the goal of this section to provide an introduction to the content of each of these chapters.

\subsection{Overview of Chapter~\ref{Part1-LWP}: abstract local well-posedness theory}

 Let us now introduce the basic framework in which we will study the general multiphasic system of equations~\eqref{eq:general_system_intro}--\eqref{eq:f_intro}. Once the abstract measure space $(I, \mu)$ is fixed, we need to specify the structural assumptions on the force fields $\FF = (F^\alpha)_{\alpha \in I}$. These will depend on $(\rhorho, \vv)$ in the following way.
\medskip

\noindent \textbf{General assumptions on the force field, informal statement.}  There exist a regularity index $k_0 \in \N$ and a time $T>0$ such that for all $t\in [0, T]$ and $k \geq k_0$, 
 \begin{itemize}
     \item[\textbf{(H1)}] there exists a map
		\begin{equation*}
		\begin{aligned}
		\mathfrak{F}_t: \quad L^1\Big([0,t];  L^\infty_\alpha H^{-1}_x \times L^\infty_\alpha L^2_x\Big)&\longrightarrow L^1\Big([0,t];   L^\infty_\alpha L^2_x\Big) \\
		(\rhorho, \vv) &\longmapsto \FF,
		\end{aligned}
		\end{equation*}
		that is locally uniformly continuous;
     \item[\textbf{(H2)}] $\mathfrak{F}_t$  maps
		 $L^1\Big([0,t];  L^\infty_\alpha H^{k-1}_x \times L^\infty_\alpha H^k_x\Big)$ into $ L^1\Big([0,t];   L^\infty_\alpha H^k_x\Big)$ in a locally bounded way.
 \end{itemize}

\medskip
Roughly speaking, the former assumption means that the force field associated with  $(\rhorho, \vv)$  gains a derivative compared to the densities $\rhorho$ and has the same regularity as the velocities $\vv$.
The main local well-posedness result concerning the general multiphasic system~\eqref{eq:general_system_intro} of Chapter \ref{Part1-LWP} can be stated as follows.

\begin{Thm}\label{thm-loose-LWP}
\noindent \textbf{Well-posedness result for multiphasic systems, informal statement.}	
		Let $(I, \mu)$ be a probability space. Under the assumptions \textbf{(H1)} and \textbf{(H2)}, for any $k \geq k_0$ there exists $T^\star>0$ such that the system~\eqref{eq:general_system_intro} admits a unique local in time solution $(\rhorho, \vv) \in	L^\infty([0,T^\star];  L^\infty_\alpha H^{k-1}_x \times L^\infty_\alpha H^k_x )$ associated with $(\rhorho_0, \vv_0) \in L^\infty_\alpha H^{k-1}_x \times L^\infty_\alpha H^k_x$. 
		\end{Thm} 
As a consequence, we obtain a result at the level of the Vlasov equation. To ease readability, we state it only when the force field in the Vlasov equation involves a coupling with the first two moments of $f$ in velocity, that is $\rho_f$ and $j_f$ (as in \eqref{eq:VP_real} and \eqref{eq:VNS_real}), although more general couplings could be allowed. 
\begin{Thm}\label{thm-loose-LWPkin}
\noindent \textbf{Well-posedness result for kinetic systems, informal statement.}
Let $(I, \mu)$ be a probability space. Consider the following nonlinear Vlasov equation
\begin{equation}\label{eq:Vlasov-toymoments}
	\left\{
	\begin{aligned}
    \partial_t f + v \cdot\nabla_x f + \mathrm{div}_v(f F_f)=0, \\
    F_f=\mathrm{F}\left[v, \int_{\R^d} f(t,x,\mathrm{d}w), \int_{\R^d} w f(t,x,\mathrm{d}w) \right], \\
    f_{\mid t=0}=f_0,
\end{aligned}
	\right.
	\end{equation}
where the map
$$(\rhorho, \vv) \mapsto \left\lbrace \mathrm{F}\left[v^\alpha, \int_{I}\rho^\alpha \, \mathrm{d}\mu(\alpha), \int_{I} v^\alpha \rho^\alpha \, \mathrm{d}\mu(\alpha) \right] \right\rbrace_{\alpha \in I} $$
satisfies Assumptions \textbf{(H1)} and \textbf{(H2)}. Assume that 
$$f_0(x, \cdot)=\int_{I} \rho^\alpha_0(x) \otimes  \delta_{v=v^\alpha_0(x)} \, \mathrm{d}\mu(\alpha)$$
for some $(\rhorho_0, \vv_0) \in L^\infty_\alpha (H^{k-1}_x \cap \dot{H}^{-1}_x) \times L^\infty_\alpha H^k_x$ with $k \geq k_0$. Then there exists $T>0$ and a unique measure-valued solution $f$ to the Vlasov equation \eqref{eq:Vlasov-toymoments} on $[0,T]$ starting at $f_0$: it is given by
\begin{align}\label{eq:f_intro-THMloose}
  f(t,x, \cdot)=\int_{I} \rho^\alpha(t,x) \otimes  \delta_{v=v^\alpha(t,x)} \, \mathrm{d}\mu(\alpha),
\end{align}
where $(\rhorho, \vv)$ is the unique local in time solution from Theorem \ref{thm-loose-LWP}.
Furthermore, all the moments in velocity of $f$ belong to $L^\infty([0,T]; H^{k-1}_x)$.
\end{Thm}

In the torus case, it is standard that the Vlasov equation~\eqref{eq:Vlasov-toymoments} is locally well-posed in Sobolev spaces in $x$ and $v$, with high enough regularity and with additional weights in $v$ (see e.g. \cite[Prop 3.2]{HKapde}). Theorem~\ref{thm-loose-LWPkin} can thus be viewed as a generalization, allowing for initial data with neither regularity nor weights in velocity; this improvement is precisely enabled thanks to the multiphasic formulation. Nevertheless, it is also fair to point out that in the standard Sobolev case, the full regularity (in space and velocity) is propagated in time, whereas Theorem~\ref{thm-loose-LWPkin} propagates only full regularity in $x$ and retains merely measure regularity in velocity.

In the whole space case, this theorem allows to deal with data that are not necessarily $L^1$ integrable with respect to the space variable. In practice, the $\dot{H}^{-1}_x$ boundedness assumption for the densities is ensured by boundedness in $L^r$ for $r\in [1,2d/(d+2)]$ in dimension $d\geq 3$.

	\begin{Rem}
		\label{rem:bounded_velocity}
		The statements of Theorem \ref{thm-loose-LWP} and \ref{thm-loose-LWPkin} are only loose because we need to be more precise in \textbf{(H1)} and \textbf{(H2)} about the topologies for which our bounds are locally uniform and about the dependence of these bounds with respect to $t$. This will be done in the upcoming Section~\ref{sec:abstract_existence}. 
		Also, the dependence with respect to $\alpha$ required for $(\rhorho, \vv)$ is up to now a bit too restrictive (indeed, for $k$ large enough with respect to the dimension: Sobolev embedding implies that the associated distribution function $f$ as defined by \eqref{eq:f_intro-THMloose} is then necessarily compactly supported in velocity); this will  be relaxed in Section~\ref{sec:abstract_existence}.

	\end{Rem}

\begin{Rem}\label{Rem:notVMaxwell}
The  assumptions \textbf{(H1)} and \textbf{(H2)} on the force field, meaning that the coupling force gains at least one derivative compared with some moments of the distribution function (as defined by \eqref{eq:f_intro}), is  in general required to ensure well-posedness. 
Therefore it does not directly apply to the case where the force field {\it a priori} has the same regularity as the moments; then  analytic regularity may be required, see \cite{BHK} in the context of the so-called relativistic Vlasov--Maxwell system. As already alluded to, in analytic regularity, the force field is even allowed to loose one derivative with respect to moments, see e.g. \cite{Gr96} or the very recent  \cite{GIRR}. 
However, there are interesting situations where, due to fine ``hidden'' regularity properties for the {\bf full} system at hand, such a gain of derivative is actually possible (in the spirit of \cite{MJ}), and this is actually the case for the multiphasic formulation of the relativistic Vlasov--Maxwell system (see \cite{Tony}); in a slightly different but related context, this is also the case for the multiphasic formulation of the Einstein--Vlasov system\footnote{We thank C\'ecile Huneau, Tony Salvi and Arthur Touati for explaining these two facts to us.}.

\end{Rem}

We will show how this general framework 
can be applied to several classes of physically relevant systems,
mainly focusing on:
\begin{itemize}
\item \textbf{Vlasov--Poisson type systems} (see Chapter \ref{Part2-VP} and in particular Section \ref{Section-PresentationVPoissonandco} for an introduction to these models), with the main example of System~\eqref{eq:VP_real} for electrons. We also consider the related system for ions, as well as the Vlasov--Monge--Amp\`ere system for Newtonian gravitation. 
As already briefly alluded to,  we shall also justify a semiclassical limit from Schr\"odinger--Poisson for  mixed states to Vlasov--Poisson, as a kind of by-product of the abstract framework.

\item \textbf{Vlasov--Navier--Stokes type systems} (see Chapter \ref{Part3-VNS} and in particular Section \ref{Section-PresentationVNSandco} for an introduction to these models), the main prototype being the system~\eqref{eq:VNS_real}, set on the torus. In addition, we allow for a variant of this system where the fluid is assumed to be barotropic and compressible, while including the treatment of the Vlasov-(steady)Stokes system in the whole space. 
\end{itemize} 
Broadly speaking, we will prove that these different types of systems all enter into the abstract set-up developed in Chapter \ref{Part1-LWP}. The inferred local well-posedness result will then prepare the ground  for the subsequent long-time analysis performed in Chapters \ref{Part2-VP}--\ref{Part3-VNS}.

\begin{Rem}
Our theory can actually cover other examples which are not treated in this manuscript. These include: (i) other kinetic models for plasmas, such as the Vlasov--Darwin equations \cite{benachour2003global, pallard2006initial} or models related to the propagation of hydromagnetic waves \cite{gardner1960similarity, alonso2024well};
    (ii) alignment systems of Cucker--Smale type, with or without chemotaxis \cite{natalini2020mean,natalini2023mean};
    (iii) other models for sprays, such as the non-homogeneous incompressible Vlasov--Navier--Stokes system \cite{ChKw} or the gyrokinetic Vlasov--Euler system \cite{MoussaSueur}.
    \end{Rem}

\subsection{Overview of Chapter~\ref{Part2-VP}: application to Vlasov--Poisson type systems}

In Chapter~\ref{Part2-VP}, we start by applying the abstract local well-posedness results of Chapter~\ref{Part1-LWP} to Vlasov--Poisson type systems. For the sake of readability, we focus in this introduction on the standard Vlasov--Poisson equation for electrons already introduced in \eqref{eq:VP_real}, that is
\begin{equation}
		\label{eq:VPintro}
		\left\{
		\begin{gathered}
		\partial_t f + v \cdot \nabla_x f - \nabla_x U\cdot \nabla_v f = 0,\\
		-\Delta_x U = \left\{
    \begin{array}{ll}
         \displaystyle \int f\D v - \displaystyle \iint f\D v \D x, \qquad & \mbox{if }  \ \ \Omega=\T^d, \\
        \displaystyle \int f\D v, \qquad & \mbox{if} \ \  \Omega=\R^d,
    \end{array}
\right.
\\
            f|_{t=0}= f_0,
		\end{gathered} 
		\right.
		\end{equation}
though other variants are treated in Chapter~\ref{Part2-VP}.
In order to put forward some concrete results, we present here only statements pertaining to distributions which are measurable functions. 
We send the reader to Theorem~\ref{thm:vp-classical} for classical well-posedness results for the Vlasov--Poisson system.

\begin{Thm}
\label{thm:LWP-VP-intro}
Let $k\in \R$ such that $k>d/2$. 
Let $f_0 \in  L^1_v (H^k_x\cap \dot H^{-1}_x)  $. There exist $T>0$ and a unique weak solution to the Vlasov--Poisson system~\eqref{eq:VPintro} on $[0,T]$, with initial condition $f_0$.

The same holds when replacing $H^k_x$ by $\mathscr{C}^{0,\theta}_x$, the set of $\theta-$H\"olderian continuous functions, for any $\theta \in (0,1]$.
\end{Thm}

 We refer to Corollaries \ref{coro-cin-VPe}--\ref{coro-cin-VPe-BESOV} for a more complete version of this statement. In particular, we will also obtain some propagation of spatial regularity for the velocity moments and some time-continuity properties that loosely speaking corresponds to strong in space and weak in velocity regularity.
   
\begin{Rem}
A few early comments are in order (more will follow in Chapter~\ref{eq:vlasov-chapter2}). First of all, we note that no regularity is required in the velocity variable. We also remark that in the whole space case (that is when $\Omega=\R^d$), neither the total mass $\int f \, \mathrm{d} v \, \mathrm{d} x $  nor  the total kinetic energy  $\int f |v|^2 \, \mathrm{d} v \, \mathrm{d} x $ of the system is necessarily finite.
Furthermore, we can make the following other comments.
\begin{itemize}
    \item We only ask for $L^1$ integrability in velocity, and no higher integrability with respect to the velocity variable is needed. In particular no control of high moments in velocity is required to obtain this local well-posedness result. We are not asking either about weights in velocity, which is in sharp contrast with standard local well-posedness results (see e.g. \cite{HKapde}). 

    \item 
    In the version of this result with $\mathscr{C}^{0,\theta}_x$, this means that little assumption on the integrability in $x$ is required (see the upcoming Remark~\ref{rem:infinitemass} for a comparison with the existing literature).

    \item The control in the homogeneous Sobolev space $\dot{H}^{-1}_x$ is useful to ensure low frequency control in the Poisson equation and is only relevant in the whole space case. In practice, this control is ensured by a further integrability assumption. For instance, in dimension $d=3$, from Sobolev embedding, we know that $L^{6/5}_x \hookrightarrow \dot{H}^{-1}_x $, so that it is sufficient to ask that $f_0 \in L^1_v (H^k_x \cap L^q_x)$ (or for  $L^1_v (\mathscr{C}^{0,\theta}_x\cap L^q_x$) for $q \in [1,6/5]$).

\end{itemize}

\end{Rem}

The second class of results proved in Chapter~\ref{Part2-VP} is a derivation of the Vlasov--Poisson system as the semiclassical limit of the Schrödinger--Poisson system for mixed states, for a class of \textit{multiphasic WKB} data, a problem introduced in Motivation IV. This will come from an application of the stability estimates required for the local well-posedness theory in the abstract multiphasic framework. We send the reader to Theorem~\ref{thm:class-semiclass} for a glimpse of known results concerning the derivation of Vlasov--Poisson in the semiclassical limit.

\begin{Thm}
Let $\Omega=\R^d$ and $(I, \mu)$ be a probability space. Consider an initial family of functions of WKB type
$$\Psi_0^{\eps,\alpha}= a_0^{\eps,\alpha} (x) \exp \left( \frac{i}{\eps} S_0^{\eps,\alpha}(x) \right), \ \ \eps \in (0,1), \ \ \alpha \in I,$$
where the amplitude $a_0^{\eps,\alpha}$ is nonnegative and the phase $S_0^{\eps,\alpha}(x)$ is real-valued, and $(\sqrt{\boldsymbol{a}^\eps_0}, \nabla \boldsymbol{S}_0^\eps)$ are uniformly bounded in a Sobolev space with sufficiently large index. Then there exists a time $T>0$ independent of $\eps$ and a solution $(\gamma^\eps)_{\eps \in (0,1)}$ to the Hartree--Poisson system 
\begin{equation}
\label{eq:SPintro-hartree-re}
	\left\{
	\begin{aligned}
i \varepsilon \partial_t \gamma^\eps &=\left[ -\frac{\varepsilon^2}{2} \Delta_x + V_\eps,\gamma^\eps\right], \\
-\Delta_x V^\eps &=  \rho^\eps, \\
\gamma^\eps|_{t=0} &=  \gamma^\eps_0,
\end{aligned}
\right.
\end{equation}
on $[0,T]$ with initial condition
$$\gamma_0^\eps(x,y) = \int_I  \Psi_0^{\eps,\alpha} (x) \overline{\Psi_0^{\eps,\alpha}}(y) \, \D \mu(\alpha).$$ 
Moreover, if $(\sqrt{\boldsymbol{a}^\eps_0}, \nabla \boldsymbol{S}_0^\eps)$ converges strongly to $(\rhorho_0, \vv_0)$ as $\eps  \rightarrow 0$, the Wigner distribution associated with $\gamma_\eps$ satisfies 
$$W[{\gamma^\eps}] \rightharpoonup f$$ 
as $\eps \rightarrow 0$, where  $f$ is the unique weak solution to the Vlasov--Poisson system \eqref{eq:VPintro} on $[0,T]$ associated with the initial condition
$$
f_0 = \int_I \rho^\alpha_0 \otimes \delta_{v=v^\alpha_0} \, \D \mu(\alpha).
$$
\end{Thm}
We refer to the upcoming Definition~\ref{def:wigner} for the definition of the Wigner distribution. More concretely, we have
\begin{Cor}
    Let  $f_0 \in L^\infty_{v,\langle v\rangle^p} (H^k_x \cap \dot{H}^{-1}_x)$ with $k>d/2$ and $p>d$, such that $f_0 \geq 0$, 
    and such that $\sqrt{f_0} \in L^\infty_{v,\langle v\rangle^p} (H^k_x \cap \dot{H}^{-1}_x)$.
    Then the unique local solution $f$ to Vlasov--Poisson associated to $f_0$ (obtained from Theorem~\ref{thm:LWP-VP-intro}) can be obtained as the semiclassical limit of the Hartree--Poisson system~\eqref{eq:SPintro-hartree-re}.

\end{Cor}
We recall that $f_0 \in L^\infty_{v,\langle v\rangle^p} (H^k_x \cap \dot{H}^{-1}_x)$ means that $\langle v\rangle^p f_0 \in L^\infty_{v} (H^k_x \cap \dot{H}^{-1}_x)$. We refer to Theorem~\ref{thm:semiclassical} and Corollary~\ref{cor:semi} for a detailed version of these statements.

\bigskip

Our last class of results proven in Chapter \ref{Part2-VP} concerns \textbf{nonlinear Lyapunov instability} for some systems studied in the present monograph, mainly assuming that they admit spatially homogeneous stationary profiles being linearly unstable. This study includes the \textbf{Vlasov--Poisson} system for electrons or ions. Our analysis is performed on the flat torus $\T^d$ for any dimension $d$ and covers profiles that are only measures in velocity. Therefore, we extend the results by~\cite{GS,HKH,HKNsima} where the stationary profiles are smooth, as well as the ones of~\cite{CGG} where only two fluids are considered, in spatial dimension $d=1$. As a matter of fact, the multiphasic framework allows to 
unify all aforementioned results. 

The abstract set-up we will develop, and implement in Chapter \ref{Part2-VP} Section \ref{Section-FrameworkInstab},  strongly relies on the method devised by Grenier in \cite{Gr00}, as was already the case for \cite{CGG,HKH,HKNsima}: we will exploit this "linear to nonlinear instability" method to obtain instability results on general multiphasic systems \eqref{eq:general_system_intro}. This will essentially be performed under the assumptions \textbf{(H1)} and \textbf{(H2)}, with a suitable spectral instability assumption (along with a mild assumption on the growth of the associated semigroup). In particular, it will cover the case of general forces $\FF$ in the model \eqref{eq:general_system_intro}, with possibly nonlinear dependence in terms of $(\rhorho, \vv)$. This will prove useful and robust enough to treat the case of electrostatic interaction between ions. 

\medskip

To ease readability and to give a simple flavour of the typical results of Chapter \ref{Part2-VP}, Section \ref{sec:insta-appli}, we state a theorem dealing with the case of Vlasov--Poisson for electrons around rough homogeneous equilibria that are spectrally unstable.
Here, spectral instability refers to the existence of a root of the associated dispersion relation (often referred to as the dielectric function in the physics literature). Typical instabilities are associated  with two-stream equilibria \cite{BohmGross,PierceH,Bun}.

In what follows and until the end of this subsection, everything is restricted to the torus $\T^d$. The instability is first obtained for the multiphasic version of \eqref{eq:VPintro}, which, we recall, reads as
\begin{equation}
\label{eq:VP-multiphaseINTRO}
\left\{ 
\begin{aligned}
\partial_t \rho^\alpha + \Div ( \rho^\alpha v^\alpha) &= 0,\\
\partial_t v^\alpha + (v^\alpha \cdot \nabla) v^\alpha &= -\nabla U,\\
-\Delta U&= \int \rho^\alpha \, \mathrm{d}\mu(\alpha) - \iint  \rho^\alpha  \, \mathrm{d}\mu(\alpha) \, \mathrm{d}x  
,\\
\rho^\alpha |_{t = 0} = \rho^\alpha_0, \ v^\alpha |_{t = 0} &= v^\alpha_0,
\end{aligned}
\right.
\end{equation}
for any measure space $(I, \mu)$. Any density--velocity pair $(1,\mathcal{V}(\alpha))_{\alpha \in I}$ is  a stationary solution to that system, corresponding to homogeneous equilibria $\varphi(v)$  at the kinetic level, by the formula \eqref{eq:f_intro}. The precise definition of the Penrose instability condition in that context (and which depends on $\mu$) has already been spelt out by the first author in \cite{Bar}. For the sake of conciseness, we do not explicitly state it here and postpone the details to  Chapter \ref{Part2-VP} Section \ref{sec:insta-appli}.

\begin{Thm}\label{thmInstab-loose-multiphase}
 \noindent \textbf{Nonlinear instability for multiphasic pressureless Euler--Poisson, informal statement}. 
Let $(1,\mathcal{V}(\alpha))_{\alpha \in I}$ be a profile that satisfies the Penrose instability condition associated with the measure space $(I, \mu)$. Let $s_0 \in \N$.  There exists $\kappa>0$ and $\delta_0>0$ such that for all $s\in \N$ and for all $\delta \in (0, \delta_0)$, there exists an initial condition $(\rho^\alpha_0, v^\alpha_0)_{\alpha \in I}$ satisfying
\begin{equation}
\esssup_{\alpha \in I} \left\| (\rho^\alpha_0-1, v^\alpha_0- \mathcal{V}(\alpha)) \right\|_{W^{s,\infty}_x} \leq \delta,
\end{equation} 
with an associated solution $(\rho^\alpha, v^\alpha)_{\alpha \in I}$ to \eqref{eq:VP-multiphaseINTRO} and a positive time $T_\delta = O(|\log \delta|)$ such that the electrostatic potential $U$ satisfies
\begin{equation}
  \| U(T_\delta) \|_{W^{-k_0,1}_x}  \geq \kappa.
\end{equation} 
\end{Thm}
In this result, for all $p \in [1,\infty]$ and $s \in \Z$, the Sobolev space  $W^{k,p}_x$ stands for the space of $L^p$ functions such that $\| \varphi\|_{W^{k,p}_x} \vcentcolon = \sum_{\ell\leq k} \| \DDD^\ell \varphi\|_{L^p_x} <+\infty$. The space $W^{-s,p}_x$ is the dual of $W^{s,p'}_x$, where $\frac{1}{p} + \frac{1}{p'}=1$.
 \index{W@$W^{k,p}$: Sobolev space based on $L^p$}

From this theorem, we infer a nonlinear instability result that holds at the level of the Vlasov equation, around rough homogeneous profiles in velocity. As before, we defer more general statements regarding a larger class of nonlinear Vlasov equations to Chapter \ref{Part2-VP}.

\begin{Thm}\label{thmInstab-loose-kin}
 \noindent \textbf{Nonlinear instability for Vlasov--Poisson, informal statement}.
Let $\mu \in \mathcal{P}_2(\R^d)$ be a probability measure (with finite variance) satisfying the Penrose instability condition. Let $k_0\in \N$. There exists $\kappa>0$ and $\delta_0>0$ such that for all $s\in \N$, for all $\delta \in (0, \delta_0)$ and for all  $\varphi \in \mathscr{C}^\infty_c(\R^d)$, there exists a measure-valued initial condition $f_{0}$ satisfying
\begin{equation}
\left\| \langle f_{0} - \mu, \varphi\rangle \right\|_{W^{s,\infty}_{x}}, \leq \delta,
\end{equation} 
with an associated measure-valued solution $f(t)$ to the Vlasov--Poisson system \eqref{eq:VPintro}  and a positive time $T_\delta = O(|\log \delta|)$ such that
\begin{equation}
\sup_{[0,T_\delta]} \| U[f] (t) \|_{W^{-k_0,1}_x}  \geq \kappa,
\end{equation} 
where $U[f] $ denotes the electrostatic potential associated with $f$.
\end{Thm}

A series of remarks are in order.

\begin{Rem}
 By leveraging this long-time instability statement, one can deduce some results on the invalidity of the so-called quasineutral limit for Vlasov--Poisson systems. We refer to Remark~\ref{rem:quasi} in Chapter~\ref{Part2-VP} for some details and references. 
\end{Rem}

\begin{Rem}We focus exclusively on homogeneous equilibria, which are the simplest equilibria of the Vlasov--Poisson system. We have not tried to apply the multiphasic framework to study the inhomogeneous cases, namely the so-called BGK equilibria (named after Bernstein, Greene, Kruskal), whose study is challenging and which have been the focus of several recent works \cite{GS-BGK,Lin1,Lin2,GL,BGHP}.
\end{Rem}

\begin{Rem} While the multiphasic framework is particularly efficient to treat nonlinear instability issues, it is fair to acknowledge that we have not been able to use it to prove stability results for (VP)--type systems. One clear issue is the possible blow-up of the multiphasic system which corresponds to the failure of the multiphasic representation; it is not clear whether such a representation is suitable to describe long-time stability results for (VP)--type systems, such as stability of the null solution in $\R^3$ (see \cite{BardosDegond,HRV,Smu,CK-VP,IPWW-scat,BVR}), or Landau damping around stable homogeneous equilibria, which has been the object of several recent works (see \cite{MouhotVillani,BMM,GNR,IPWW} in the torus case).
\end{Rem}

\subsection{Overview of Chapter~\ref{Part3-VNS}: application to Vlasov--Navier--Stokes type systems}

In a third direction, we focus on fluid-particle models, namely \textbf{Vlasov--Navier--Stokes} type systems, set on the flat torus.

In Chapter~\ref{Part3-VNS}, we apply the abstract results from Chapter~\ref{Part1-LWP} to the local well-posedness theory for Vlasov--Navier--Stokes type systems. For the sake of readability, we focus in this overview on the standard Vlasov--Navier--Stokes equation (for an incompressible fluid) introduced in \eqref{eq:VNS_real}:
	\begin{equation}
		\label{eq:VNSintro}
		\left\{
		\begin{gathered}
		\partial_t f + v \cdot \nabla_x f +\Div_v\big( (U-v)f\big) = 0,\\
		\partial_t U + (U \cdot \nabla_x) U + \nabla_x P - \Delta_x U = \int (v-U) f \D v ,\\
		\Div_x(U) = 0,\\
		f|_{t=0} = f_0, \quad U|_{t=0} = U_0,
		\end{gathered}
		\right.
		\end{equation}
but some important variants are also treated in Chapter~\ref{Part3-VNS}.

The analogue of Theorem~\ref{thm:LWP-VP-intro} for the Vlasov--Poisson equation is the following.

\begin{Thm}
\label{thm:LWP-VNS-intro}
Let $k \in \N$ such that $k>d/2$. Let $U_0 \in H^{k+1}_x$ with $\Div_x(U_0)=0$  
and $f_0 \in L^1_{v,\langle v \rangle} H^k_x  $. There exist $T>0$ and a unique weak solution to the Vlasov--Navier--Stokes system~\eqref{eq:VNSintro} on $[0,T]$, with initial condition $(f_0,U_0)$.

\end{Thm}

We recall that $f_0 \in L^1_{v,\langle v \rangle} H^k_x $ means that $\langle v \rangle f_0 \in L^1_v H^k_x$.
The same remarks as in the Vlasov--Poisson case are still relevant. Note that in the Vlasov--Navier--Stokes, we need to ensure higher integrability so as to define the first two moments $\rho_f$ and $j_f$, but no  higher integrability is needed.  In contrast, standard well-posedness results for the incompressible Vlasov--Navier--Stokes system are summarized in Theorem~\ref{thm:vns-classical} from Chapter~\ref{Part3-VNS}.

\medskip

In the second part of Chapter \ref{Part3-VNS}, \textbf{asymptotic concentration in velocity} is proven for {Vlasov--Navier--Stokes} type systems: namely, we build global in time solutions with initial conditions sufficiently close to equilibrium.
We aim at proposing a direct approach that unravels the mechanism of convergence towards a monokinetic profile, which results from the interplay between friction and dissipation. The multiphasic formulation leveraged in this article turns out to be particularly adapted to tackle this problem.

\medskip

For expository purposes, we again focus here on the case of the Vlasov--Navier--Stokes system \eqref{eq:VNSintro}: given a probability space $(I, \mu)$, its multiphasic formulation is 
\begin{equation}
\label{eq:VNS-multiphaseINTRO}
\left\{ 
\begin{aligned}
\partial_t \rho^\alpha + \Div ( \rho^\alpha v^\alpha) &= 0,\\
\partial_t v^\alpha + (v^\alpha \cdot \nabla) v^\alpha &= U - v^\alpha,\\
\partial_t U + (U\cdot \nabla ) U + \nabla P -  \Delta U &= \int \rho^\alpha \big( v^\alpha - U )  \D \mu(\alpha)
,\\
\Div(U) &= 0, \\
\rho^\alpha |_{t = 0} = \rho^\alpha_0, \ v^\alpha |_{t = 0} &= v^\alpha_0, \ U |_{t = 0} = U_0.
\end{aligned}
\right.
\end{equation}
Compared with the kinetic version \eqref{eq:VNSintro}, the expected alignment  in velocity, previously described in \eqref{def:alignmentVNS}, is particularly apparent: it corresponds to the convergence of $U(t)$ and $\vv(t)$ towards the same common constant velocity as $t \rightarrow + \infty$. The main observation is that for any constant speed $\mathcal{V} \in \R^d$, the family $(\rho^\alpha, v^\alpha, U)=(1, \mathcal{V}, \mathcal{V})$ (thus constant for every phase $\alpha$) is a steady state of \eqref{eq:VNS-multiphaseINTRO}, as it directly cancels out the drag terms $v^\alpha-U$.
This paves the way for a perturbative approach around this specific class of equilibria\footnote{We emphasize the fact that we consider only constant equilibria here, leaving aside the interesting case of more general velocity profiles, such as shear flows.}. We will prove that small regular perturbations give rise to global in time solutions, with velocities $U$ and $\vv$ converging towards another  common asymptotic velocity close to $\mathcal{V}$.

\medskip

In what follows, the spatial domain is always the periodic torus $\T^3$. One of the main results of Chapter \ref{Part3-VNS} is the following. 
\begin{Thm}\label{thmVNS-loose-multiphase}
 \noindent \textbf{Asymptotic stability for multiphasic pressureless Euler--Navier--Stokes, informal statement}.
Let $\mathcal{V} \in \R^3$ and $k \in \N$ such that $k>2+3/2$. If
    \begin{align*}
 \Vert \rhorho_0-1 \Vert_{L^{\upinfty}_\alpha H^{k-1}_x}+ \Vert \vv_0-\mathcal{V} \Vert_{L^{\upinfty}_\alpha H^{k}_x}+ \Vert U_0-\mathcal{V} \Vert_{H^k_x} +\Vert \langle U_0 - \vv_0\rangle \Vert_{L^{\upinfty}_\alpha} &\ll 1
    \end{align*}
then the system \eqref{eq:VNS-multiphaseINTRO} admits a unique global in time solution $(\rhorho, \vv, U)$ such that
\begin{align*}
\rhorho \in L^\infty(\R^+; L^\infty_\alpha H^{k-1}_x)  , \ \ \vv \in L^\infty(\R^+;  L^\infty_\alpha H^{k}_x), \ \ U \in L^\infty(\R^+; H^k_x ),
\end{align*}
and with initial data $(\rhorho_0, \vv_0, U_0)$. Furthermore, there exists $\lambda>0$, an asymptotic velocity $\mathcal{W}_{\infty} \in \R^3$ and an asymptotic profile $\bm{\rho}_{\infty} \in L^{\infty}_\alpha L^{\infty}_x$ such that for all $t \geq 0$,
\begin{align*}
\esssup_{\alpha \in (I, \mu)} \mathrm{W}_1 \Big(\rho^{\alpha}(t),\rho^{\alpha}_{\infty}(x-t \mathcal{W}_{\infty}  ) \Big)+\Vert \vv(t) - \mathcal{W}_{\infty} \Vert_{L^{\upinfty}_\alpha L^2_x} &+ \Vert U(t) - \mathcal{W}_{\infty} \Vert_{L^2_x} \\
&\lesssim e^{-\lambda t},
\end{align*}
where $\mathrm{W}_1$ is the $1-$Wasserstein distance on $\T^3$ (recall Definition \ref{def:Wasserstein1}).
\end{Thm}

As a consequence, we recover under different assumptions the main result of~\cite{HKMM} for the Vlasov--Navier--Stokes system.
\begin{Thm}\label{thmVNS-loose}\noindent \textbf{Asymptotic stability for Vlasov--Navier--Stokes, informal statement.}
Let $\mathcal{V} \in \R^3$ and $k \in \N$ such that $k>2+3/2$. If $U_0 \in H^k_x$ and
\begin{equation*}
	f_0(x,v) = \int_I \rho^\alpha_0(x) \otimes \delta_{v= v^\alpha_0(x)} \, \D \mu(\alpha), \ \ \text{with} \ \ \rhorho_0 \in L^\infty_\alpha H^{k-1}_x, \vv_0 \in L^\infty_\alpha H^{k}_x,
	\end{equation*}
satisfy
\begin{align*}
  \Vert \rhorho_0-1 \Vert_{L^{\upinfty}_\alpha H^{k-1}_x}+ \Vert \vv_0^\alpha-\mathcal{V} \Vert_{H^{k}_x}+ \Vert U_0-\mathcal{V} \Vert_{H^k_x}+ \Vert \langle U_0 - \vv_0\rangle \Vert_{L^{\upinfty}_\alpha} &\ll 1,
    \end{align*}
then the following holds. There exists a unique global in time solution $(f, U)$ to the incompressible Vlasov--Navier--Stokes system \eqref{eq:VNSintro} such that for all $T>0$
\begin{align*}
    &U \in L^\infty([0,T]; H^k_x) \cap L^2([0,T]; H^{k+1}_x), \\
&\forall \psi \in \mathscr{C}_b(\R^3), \ \ (t,x)  \mapsto \int_{\R^d} \psi(v) f(t,x,\D v) \in L^\infty([0,T]; H^{k-1}_x).
\end{align*}
Furthermore, we have 
\begin{align*}
\mathrm{W}_1 \left(f(t),\left( \int_I \, \rho^\alpha_\infty(x - t \mathcal{W}_\infty)\mathrm{d}\mu(\alpha) \right) \otimes\delta_{v=\mathcal{W}_\infty}  \Big) \right) + \Vert U(t)- \mathcal{W}_\infty  \Vert_{L^2_x} \underset{t \rightarrow + \infty}{\longrightarrow}0,
\end{align*}
for some asymptotic profile $\rhorho_\infty$ and some asymptotic velocity $\mathcal{W}_\infty \in \R^3$. Here $\mathrm{W}_1$ is the $1-$Wasserstein distance on $\T^3 \times \R^3$ (recall Definition \ref{def:Wasserstein1}).
\end{Thm}

Several comments are now in order - in particular in light of some comparisons with the result for (VNS) on the torus recalled in Theorem \ref{thm:VNSlargetime}.
\begin{itemize}
    \item We refer to the more complete statements of Theorems \ref{thm-globalconvergence} and \ref{thm-cinetique-VNSincomp} in Chapter \ref{Part3-VNS}, Section \ref{Section:monokinetic-stab-theorem}, in which we weaken the assumptions on the initial perturbation (allowing for non-compactly supported distribution functions) and actually obtain refined rates of convergence in higher Sobolev norms.
    \item Our result should be compared with the one obtained in \cite{HKMM}  (see Theorem \ref{thm:VNSlargetime}), where the authors work directly on \eqref{eq:VNSintro}: in Theorem \ref{thmVNS-loose}, no regularity assumptions are required in $v$ for the distribution function $f$, that may be merely a measure. However, we ask for stronger regularity in $x$ compared to \cite{HKMM}, where the solutions roughly belong to the natural energy space \eqref{def:energyspaceVNS}.
    \item It is also important to emphasize that Theorem \ref{thmVNS-loose} only holds when the spatial domain is the periodic torus: at the heart of our method is an explicit spectral gap at the linearized multiphasic level (despite the apparent  partially dissipative nature of the system), which cannot be directly captured by the methods from \cite{HKMM}. 
    Hence it does not apply on the directly to the whole space $\R^3$. Note also that in the periodic case, and because of conservation laws, the asymptotic velocity $\mathcal{W}_\infty$ is \textit{different} from the initial constant velocity $\mathcal{V}$ that we perturb (see the statement of Theorem \ref{thm-globalconvergence}).  

\end{itemize}

\begin{Rem}
Our strategy 
is robust enough to handle the case  of the Vlasov equation  coupled with a compressible isentropic Navier--Stokes system. We refer to Chapter \ref{Part3-VNS} for more details and precise statements -- see in particular Theorem \ref{thm-cinetique-VNScomp}. 
\end{Rem}

\section{Equivalence between the multiphasic and kinetic formulations}
	\label{sec:equivalence_models}

As said before, we systematically work on the torus or on the whole space for the spatial variable: in what follows, $\Omega$ stands for either $\T^d$ or $\R^d$. This section is dedicated to general considerations about the multiphasic formulation introduced in Section \ref{Section-introMultiphaseFormulation}.

\begin{Def}
    We call a \emph{set of labels} a measure space $(I,\mu)$ where $\mu$ is a nonnegative finite measure.
    \index{I@$(I,\mu)$: set of labels}
\end{Def}

	\emph{Throughout this work, we systematically restrict to the case where the force field is Lipschitz continuous.}
 The aim of this preliminary section is to study, given a set of labels $(I,\mu)$,  the link between  multiphasic equations~\eqref{eq:general_system_intro}-\eqref{eq:link_force_Vlasov_multiphase}, namely
 	\begin{equation}
 	\label{eq:multiphasic-equivalence}
	\left\{
	\begin{gathered}
	\partial_t \rho^\alpha + \Div( \rho^\alpha v^\alpha) =0, \\
	\partial_t v^\alpha + ( v^\alpha \cdot \nabla ) v^\alpha = F^\alpha,\\
	\rho^\alpha |_{t = 0} = \rho^\alpha_0, \ v^\alpha |_{t = 0} = v^\alpha_0,
	\end{gathered}
	\right.
	\end{equation}
and Vlasov~\eqref{eq:Vlasov_general_intro} equations, namely
 	\begin{equation}
 	 	\label{eq:vlasov-equivalence}
 	\partial_t f + v \cdot \nabla_x f + \mathrm{div}_v \big( F f \big) = 0,
	\end{equation}
    where the precise relation between the forces $(\mathrm{F}^\alpha)$ and $F$ needs to be specified depending on the context.  The key point is the operation described in~\eqref{eq:f_intro} which maps families of type $(\rhorho,\vv) = (\rho^\alpha,v^\alpha)_{\alpha \in I}$ into measures $f = \mathcal I (\rhorho,\vv)$\footnote{To be precise, this map $\mathcal I$ depends on the set of labels $(I,\mu)$, but we forget about this detail in this paragraph.} on the phase space $\Omega \times \R^d$.

When the family $(v^\alpha)_{\alpha \in I}$ enjoys Lipschitz regularity, we will first prove that this operation maps solutions of~\eqref{eq:multiphasic-equivalence} to solutions of~\eqref{eq:vlasov-equivalence}. Second, we will show that this operation has a right inverse $\mathcal J$, and that calling $f$ a solution of~\eqref{eq:vlasov-equivalence} starting from $f_0$, and $(\rhorho,\vv)$ the solution of~\eqref{eq:multiphasic-equivalence} starting from $\mathcal J(f_0)$ (which will be shown to exist for small time under some assumptions that are coherent with our upcoming general framework), then for all $t$ for which $(\rhorho,\vv)$ is well defined, $f = \mathcal I(\rhorho,\vv)$: that is, $f$ has a multiphasic representation at least for small times.

	The section is organized as follows. In Subsection~\ref{subsec:def_I}, we give a proper definition of the operation described in~\eqref{eq:f_intro} and show that it admits a right inverse. Then, we focus on the uncoupled case, that is, the case when the force field does not depend on the solution itself. In Subsection~\ref{subsec:def_solutions}, we define what we mean by solutions of~\eqref{eq:vlasov-equivalence} and~\eqref{eq:multiphasic-equivalence}, respectively.  In Subsection~\ref{subsec:characterisics}, we provide the main argument for the proof of equivalence: both systems admit as characteristics the solutions of Newton's law of classical mechanics. We conclude the equivalence of the two systems in the uncoupled case in Subsection~\ref{subsec:equivalence}. Finally, in Subsection~\ref{subsec:coupled}, we apply our results to the coupled case, that is where the force field depends on the solution itself. 
	
	\subsection{From the multiphasic formalism to densities in the phase space}
	\label{subsec:def_I}
	The first aim is to give a precise sense to the formula~\eqref{eq:f_intro}. We start by defining a notion of {\bf regular} multiphasic distributions.
		\begin{Def}
		Let $(I,\mu)$ be a set of labels.     
        \index{r@$(\rhorho,\vv)$: multiphasic distribution}
A multiphasic distribution labeled by $(I,\mu)$ is a measurable family $(\rhorho,\vv) = (\rho^\alpha, v^\alpha)_{\alpha \in I}$, well defined up to $\mu$-negligible sets, such that for $\mu$-almost all $\alpha$, $\rho^\alpha$ is a nonnegative Radon measure on $\Omega$, and $v^\alpha : \Omega \to \R^d$ is a measurable vector field. We further say that $(\rhorho,\vv)$ is {\bf regular} provided
		\begin{equation*}
		\esssup_{\alpha} \|  \DDD v^\alpha \|_{L^\infty_x} < + \infty.
		\end{equation*}
	\end{Def}
	The following straightforward proposition corresponds to the rigorous definition behind formula~\eqref{eq:f_intro}. In what follows, $T \pf \nu$ denotes the pushforward of a measure $\nu$ by a mapping $T$.
	\begin{Prop}\label{Prop-meaningI-mu-rho-v}
		Let $(I,\mu)$ be a set of labels and $(\rhorho,\vv)$ be a multiphasic distribution. The map
		\begin{equation*}
		A \in \mathcal B(\Omega \times \R^d) \longmapsto \int_I(\Id, v^\alpha)\pf \rho^\alpha (A) \D \mu(\alpha) \in [0, + \infty],
		\end{equation*}
		 is a measure on $\Omega \times \R^d$, that we denote by $\mathcal I_\mu(\rhorho,\vv)$.     \index{I@$\mathcal I_\mu(\rhorho,\vv)$: distribution function associated with the multiphasic distribution $(\rhorho,\vv)$}
 Setting $f \vcentcolon =  \mathcal I_\mu(\rhorho,\vv)$, we find that for all nonnegative measurable function $\psi$,
		\begin{equation*}
		\int \psi(x,v) \D f(x,v) = \iint \psi(x,v^\alpha(x)) \D \rho^\alpha(x) \D \mu(\alpha) .
		\end{equation*}
	\end{Prop}
	
	\begin{Rem}
		Actually, under mild assumptions on $(\rhorho,\vv)$, the distribution $\mathcal I_\mu(\rhorho,\vv)$ in phase space has some regularity with respect to $x$. For instance, provided for $\mu$-almost all $\alpha$, $\rho^\alpha$ and $v^\alpha$ are continuous functions of $x$ and
		\begin{equation*}
		 \int \| \rho^\alpha \|_{L^\infty(\Omega)} \D \mu(\alpha) < + \infty,
		\end{equation*}
		then $f = f(x, \D v)$ can be viewed as a continuous function of $x$ with values in the set of nonnegative finite measures on $\R^d$, meaning that for any test function $\psi \in \mathscr{C}_b(\R^d)$, the function
		\begin{equation*}
		m_\psi(f): \qquad x \mapsto \int \psi(v) f(x,\D v) \vcentcolon = \int \psi(v^\alpha(x)) \rho^\alpha(x) \D \mu(\alpha)
		\end{equation*}
		is finite and continuous -- this is a consequence of the dominated convergence theorem. Additional uniform bounds on the derivatives of $(\rhorho,\vv)$ would imply higher regularity of $m_\psi(f)$: we refer to Lemma \ref{prop-reg-momentsf} for a more precise statement in the framework of Chapter \ref{Part1-LWP}.  
	\end{Rem}

	A crucial point for showing that any solution of the Vlasov equation has a multiphasic representation is that this operation $\mathcal I$ has a right inverse.  This is a consequence of a disintegration lemma (cf \cite[III-70]{dellacherie1978probabilities}). Moreover, this right inverse takes values in the set of~\emph{regular} multiphasic distributions. 
	
	\begin{Prop}
		\label{prop:disintegration}
		Let $f$ be a nonnegative Radon measure on $\Omega \times \R^d$. Then, there exists  a set of labels $(I,\mu)$ and a multiphasic distribution $(\rhorho,\ww)$  such that $\mathcal I_\mu(\rhorho,\ww) = f$. In addition, we can enforce that $\mu$ is finite, and that for $\mu$-almost all $\alpha \in I$, $w^\alpha$ is constant. In particular, this choice of $(\rhorho,\ww)$ makes it a regular multiphasic distribution.
	\end{Prop}
	There is a lot of freedom in the construction of this right inverse. We propose three constructions. The first two ones are very natural but restricted to specific cases, and the last one is fully general but rather exotic.
	\begin{proof}
		
		\noindent \underline{Case 1: $f$ is a finite measure.}
In this case, we can take $I = \R^d$, $\mu \vcentcolon = \pi_2 {}\pf f$, where $\pi_2$ is the canonical projection from $\Omega \times \R^d$ to $\R^d$, for $\mu$-almost all $v \in \R^d$, $w^v \equiv v$, and for $\mu$-almost all $v \in \R^d$, 
$$\rho^v(x)= f(x | \pi_2 = v),$$ 
where $f(\cdot | \pi_2 = v)$ is the corresponding conditional law  given by the disintegration lemma. With this choice, starting from the disintegration lemma, we have for all test function $\varphi \in \mathscr{C}_b(\R^d \times \R^d)$:
\begin{equation*}
\int \varphi(x,v) \D f(x,v) = \int \left( \int \varphi(x,v) \D \rho^v(x) \right) \D \mu(v).
\end{equation*}
This last quantity rewrites
\begin{align*}
\int \left( \int \varphi(x,v) \D \rho^v(x) \right) \D \mu(v) &= \int \left( \int \varphi(x, w^v(x)) \D \rho^v(x) \right) \D \mu(v) \\
&= \int \left( \int \varphi \D \, (\Id, w^v) \pf \rho^v  \right)\D \mu(v) = \int \varphi \D \mathcal I_\mu(\rhorho,\vv).
\end{align*}
So indeed, we find $\mathcal I(\rhorho,\vv) = f$.

		\bigskip
		
		\noindent \underline{Case 2: $f$ is in $L^1_{\mathrm{loc}}(\Omega \times \R^d)$.}
		In this case, we can choose $I = \R^d$ as well, for $\mu$ any $L^1$ function on $\R^d$ that is positive almost everywhere, for $\Leb$-almost all $v \in \R^d$, $w^v \equiv v$ and for $\Leb$-almost all $(x,v)$ in $\Omega \times \R^d$, 
		\begin{equation*}
		\rho^v(x) = \frac{f(x,v) }{\mu(v)}.
		\end{equation*}
		This formula provides a $L^1_{\mathrm{loc}}$ distribution on $\Omega$ for $\Leb$ and hence $\mu$-almost all $v$, and therefore a Radon measure as well.
	This time as well, the identity $\mathcal I_\mu(\rhorho,\ww) = f$ is easily checked as for all test function $\varphi \in \mathscr{C}_c(\R^d \times \R^d)$, identifying measures with their densities with respect to the Lebesgue measure,
		\begin{align*}
		\int \varphi(x,v) \D f(x,v) &= \iint \varphi(x,v) f(x,v) \D x \D v \\
        &= \int \left( \int \varphi(x,v) \rho^v(x)  \D x \right) \mu(v) \D v = \int \varphi \D \mathcal I_\mu(\rhorho,\vv),
		\end{align*}			
		where the last equality follows from the same lines as in the previous case.

		\bigskip
		
		\noindent \underline{Case 3: General data.}
		In this case, we take $I \vcentcolon= \N^* \times \R^d$. For all $(n,v) \in I$, we define $w^{n,v} \equiv v$ and $\ww \vcentcolon = (w^{n,v})_{(n,v) \in I}$. We also set $\bar \ww \vcentcolon = (\overline{w}^{v})_{v \in \R^d}$ with $\bar w^v \equiv v$. Let $(C_n)_{n\in\N^*}$ be a countable partition of $\Omega \times \R^d$ by bounded sets. As $f$ is Radon, for all $n \in \N^*$, $f_n \vcentcolon = \1_{C_n} f$ is finite. Therefore, the construction of Case~1 applies and leads to $\mu_n$ and $(\rhorho_n)$ such that $f_n = \mathcal I_{\mu_n}(\rhorho_n, \bar \ww)$. Let us define $\mu$ as the finite nonnegative measure $\mu$ on $\N^* \times \R^d$ that satisfies for all Borel set $A \subset \R^d$ and $n \in \N^*$:
		\begin{equation*}
		\mu( \{n\} \times A) \vcentcolon = \frac{\mu_n(A)}{2^n \mu_n(\R^d)}.
		\end{equation*}
		Now, for $\mu$-almost all $(n,v)$ (that is, for all $n \in \N^*$ and $\mu_n$-almost all $v$), let us set
		\begin{equation*}
		\rho^{n,v} \vcentcolon = 2^n \mu_n(\R^d) \rho^v_n.
		\end{equation*}
		For $\mu$-almost all $(n,v)$, this is a finite measure on $\Omega$. Finally, for all nonnegative test function $\varphi$ on $\Omega \times \R^d$, we have
		\begin{align*}
		\int_I \int_\Omega \varphi(x,w^{n,v}(x)) &\D \rho^{n,v} (x) \D \mu(n,v) \\
        &= \sum_n \frac{1}{2^n \mu_n(\R^d)} \int_{\R^d} \int_\Omega \varphi(x,v) 2^n \mu_n(\R^d)\D \rho_n^v(x)\D \mu_n(v)\\
		&= \sum_n \int_{\R^d} \int_\Omega \varphi(x,v) \D \rho_n^v(x)\D \mu_n(v)\\
		&= \sum_n \int_{\Omega \times \R^d} \int \varphi \D f_n(x,v) = \int \varphi \D f.
		\end{align*}
		Therefore, $\I_\mu(\rhorho,\vv) = f$, which completes the proof of the proposition.
	\end{proof}

In practice, one should keep in mind the decompositions introduced in Cases 1 and 2 from the latter proof, which are handy in several practical situations. Note also that in these cases, $w^\alpha= \alpha$ and is thus constant in space, in particular not integrable when $\Omega=\R^d$. This fact will influence the abstract framework developed in the monograph.

In Chapter~\ref{Part1-LWP}, we will develop a Cauchy theory for multiphasic systems, which requires Sobolev smoothness in $x$ (in a sense to be made precise later) for an initial multiphasic distribution $(\rhorho_0,\ww_0)$. In view of applications to Vlasov equations, which are associated to an initial condition $f_0=\mathcal I_\mu(\rhorho_0,\ww_0)$, according to the constructions of Proposition~\ref{prop:disintegration}, this  means, loosely speaking, that we can consider initial distribution functions $f_0(x,v)$ that are measures in $v$ but are Sobolev with respect to $x$.

\begin{Rem} The multiphasic representation provided in the proof of Proposition~\ref{prop:disintegration} is defintely not unique.
For instance, for the case of a sum of $n$ Dirac masses in velocities, that is 
     $$f=\sum_{j=1}^n \rho_j \otimes  \delta_{v=v_j},$$ 
recall that it is natural to take
$$I=[\![1, n]\!], \quad \mathrm{d}\mu(\alpha)=\frac{1}{n}\sum_{j=1}^n \delta_{\alpha=j},$$ 
whereas the proof of Proposition~\ref{prop:disintegration} yields $I=\R^d$.

\end{Rem}

Let us now explain how to ``concatenate'' a sum of two multiphasic representations, which is useful to treat ``mixed'' data such as those involving an integrable part, and a Dirac part, recall~\eqref{def:Nkinsol} in Motivation I.
\begin{Lem}
\label{lem:conc}
    For $i=1,2$, let $(I_i,\mu_i)$ be a set of labels and let $(\rhorho_i,\ww_i)$ be a multiphasic distribution. Consider the distribution function 
    $$
    f = \mathcal I_{\mu_1}(\rhorho_1,\ww_1) + \mathcal I_{\mu_2}(\rhorho_2,\ww_2).
    $$
    Then there exist a set of labels $(I,\mu)$ and a family $(\rhorho,\ww)$ such that $f=\mathcal I_\mu(\rhorho,\ww)$.
\end{Lem}

\begin{proof} We could apply Proposition~\ref{prop:disintegration} but our purpose is to give a practical way to concatenate two sets of labels. Consider the extended set of labels 
$$I\vcentcolon = \{1\}\times I_1 \cup \{2\}\times I_2,$$ 
still denoting $\mu_i$ the extension of the measure on $\{i\}\times I_i$, we endow $I$ with the canonical measure $\mu$ for the union. We can then define  $(\rhorho,\ww)$ as
$$
(\rho^\alpha, w^\alpha) = (\rho_i^{\alpha'}, w_i^{\alpha'})   \quad \text{for} \quad \alpha=(i,\alpha'),
$$
and it follows that 
\begin{align*}
\mathcal I_\mu(\rhorho,\ww) &= \int_I \rho^\alpha \otimes \delta_{v=w^\alpha} \D \mu(\alpha) \\
&= \int_{I_1} \rho_1^\alpha \otimes \delta_{v=w_1^\alpha} \D \mu_1(\alpha) + \int_{I_2} \rho_2^\alpha \otimes \delta_{v=w_2^\alpha} \D \mu_2(\alpha) \\
&= \mathcal I_{\mu_1}(\rhorho_1,\ww_1) + \mathcal I_{\mu_2}(\rhorho_2,\ww_2) = f,
\end{align*}
hence the lemma.
\end{proof}

\begin{Rem}

To conclude this subsection, let us discuss another type of multiphasic representation from those of the proof of Proposition~\ref{prop:disintegration}, that is less general but perhaps somehow more concrete. Fix some measure space $(\mathscr{N},\nu)$ and a family of  functions $(\psi^y)_{y\in \mathscr{N}}$ from $\R^d$ to $\R$ which can be understood as a sort of orthonormal basis. 
Assume we have the general decomposition for all ``nice'' functions $f$:
	\begin{equation}
	\label{eq:fdecomp}
	    f(x,v) =\int_\mathscr{N} m^y_f (x) \, \psi^y(v) \, \D \nu(y),	\end{equation}
        for a certain measure space $(\mathscr{N},\nu)$.
	    In practice, one may assume that there exists a family of functions $\left(\widetilde{\psi}^y \right)_{y\in \mathscr{N}}$ also from $\R^d$ to $\R$ (the ``dual'' family of $(\psi_y)_{y \in  \mathscr{N}}$) such that 
	    $$
	    m^y_f (x) =\left\langle f(x,\cdot), \widetilde{\psi}_y \right\rangle,
	    $$
	    so that $(m^y_f)_{y \in  \mathscr{N}}$ are actually moments of $f$.
For instance, assuming that the Fourier inversion formula holds for all $f(x, \cdot)$, we have
$$
f(x,v) =\int_\xi  \mathcal{F}_v f (x, \xi) e^{i \xi \cdot v} \, \mathrm{d} \xi, \qquad \mathcal{F}_v f (x, \xi)= \frac{1}{(2\pi)^d}\int_v f(x,v) e^{-ix\cdot \xi} \, \D v.
$$
In this example, $\mathscr{N}=\R^d$ and $m_f^y(x)= \mathcal{F}_v f (x, y)$.
Another example that comes to mind is based on the use of an orthonormal basis $(\psi_k)_{k\in\N}$ of a given weighted-$L^2$ space to write
$$
f(x,v) =\sum_{k \in\N}\langle f(x,\cdot), \psi_k\rangle \,   \psi_k (v),
$$
in which case $\mathscr{N}=\N$.

It turns out that the decomposition \eqref{eq:fdecomp} formally gives rise to a multiphasic representation $f=\mathcal I_\mu(\rhorho,\ww) $, with $I = \R^d \times \mathscr{N}$. Indeed, given $\varphi$ a real-valued test function, we have
\begin{align*}
    \int_{\R^{2d}} f(x,v) \varphi (x,v) \, \mathrm{d} v \, \mathrm{d} x &=  \int_{\R^{2d}} \int_\mathscr{N} m_f^y(x)  \psi^y(v)  \, \D \nu(y)   \varphi (x,v)  \,\mathrm{d} v \,  \mathrm{d} x \\
    &=  \int_{\R^d} \int_{\R^d \times \mathscr{N}}  m_f^y(x) \psi^y(v)  \varphi(x,v)   \,  \D v \,  \D \nu(y) \,  \mathrm{d} x .
\end{align*}
Fix $a \in L^1(\R^d)$ and $b \in L^1(\mathscr{N})$. We can define
\begin{align*}
    I\vcentcolon=\R^{d} \times \mathscr{N}, \ \ \alpha\vcentcolon=(v,y) \in I, \ \ \mathrm{d}\mu(\alpha)\vcentcolon= a(v) b(y) \, 
    \mathrm{d} v  \, \mathrm{d} \nu(y), 
\end{align*}
so that $\mu$ is finite, and
\begin{align*}
\rho^\alpha(x) \vcentcolon= \frac{\psi^y(v)}{a(v) b(y)} m_f^y(x) , \ \ v^\alpha(x)\vcentcolon=v.
\end{align*}
Therefore, we end up with the multiphasic decomposition 
\begin{align*}
    f=\int_{I} \rho^{\alpha} \otimes \delta_{v=v^\alpha} \, \mathrm{d}\mu(\alpha).
\end{align*}

We conclude this long  remark with two additional observations. 
\begin{itemize}
    \item 
In the following, as already alluded to, we will systematically require Sobolev regularity for $\rho^\alpha$ and $v^\alpha$. One may thus believe that thanks to this approach, as opposed to that of Proposition~\ref{prop:disintegration}, only regularity of the moments $(m^y_f)$ is required, and not of $f$ itself. However because of the weight $\frac{1}{b(y)}$ in $\rho^\alpha(x)$, one can apply the dominated convergence theorem to~\eqref{eq:fdecomp} to obtain Sobolev regularity in $x$ for the whole distribution function $f$. 

\item This discussion incidentally shows that the approach for studying Vlasov equations that is  based on using a decomposition in a Hermite-like basis, that is popular in numerical analysis, see e.g. \cite{BF1,BF2,BDFV} and references therein, can actually enter the framework of the multiphasic formulation put forward in this monograph.

\end{itemize}

\end{Rem}

	\subsection{Notions of solutions of the uncoupled problem}
	\label{subsec:def_solutions}
	In this subsection, we consider  a given force field $F:[0,T] \times \Omega \times \R^d \to \R^d$ with the following regularity:
	\begin{equation}
	\label{eq:def_lipschitz_norm}
	\| F \|_{\mathcal L,T} \vcentcolon = \int_0^T \Big\{ |F(t,0)| + \Lip(F(t,\cdot ))  \Big\} \D t < + \infty,
	\end{equation}
	where $\Lip(\cdot)$ denotes the Lipschitz semi-norm. This level of regularity, i.e. Lipschitz force fields,  was already alluded to in the beginning of the section and is chosen to ensure global existence and uniqueness for Newton's law starting from any initial condition in the phase space.

 We first give a precise definition of solutions to the Vlasov equation~\eqref{eq:vlasov-equivalence} and to the multiphasic system~\eqref{eq:multiphasic-equivalence} 
 in the \textit{uncoupled case}: namely, we study the multiphasic system
\begin{equation}
 	\label{eq:multiphasic-equivalence2}
	\left\{
	\begin{gathered}
	\partial_t \rho^\alpha + \Div( \rho^\alpha v^\alpha) =0, \\
	\partial_t v^\alpha + ( v^\alpha \cdot \nabla ) v^\alpha = F^\alpha,\\
		F^\alpha(t,x)\vcentcolon = F \left(t,x,v^\alpha(t,x) \right), \\
	\rho^\alpha |_{t = 0} = \rho^\alpha_0, \ v^\alpha |_{t = 0} = v^\alpha_0,
	\end{gathered}
	\right.
	\end{equation}
and the Vlasov equation
\begin{equation}
\left\{
	\begin{gathered}
 	 	\label{eq:vlasov-equivalence2}
 	\partial_t f + v \cdot \nabla_x f + \mathrm{div}_v \big( F f \big) = 0, \\
	f |_{t = 0} = f_0,
    \end{gathered}
	\right.
	\end{equation}
by focusing on the link between their solutions.

	In the case of the Vlasov equation~\eqref{eq:vlasov-equivalence2}, we simply consider distributional solutions, in the spirit of~\cite{ambrosio2014continuity}.
	\begin{Def}[Weak solutions of the Vlasov equation]
		\label{def:solution_Vlasov}
		We call a solution of~\eqref{eq:vlasov-equivalence2} up to time $T>0$ any measurable curve $f = f(t)$, $t \in [0,T]$, valued in the space of nonnegative Radon measures of $\Omega \times \R^d$ such that for all compact $K \subset \Omega \times \R^d$, the map $t \mapsto f(t,K)$ belongs to $L^\infty([0,T])$, and which satisfies~\eqref{eq:vlasov-equivalence2} in the sense of distributions. 
	\end{Def}
	
		In such a case, standard arguments show that up to modifications on negligible sets of time, $f$ is weak-$\star$ continuous with respect to time in the dual of $\mathscr{C}_c(\Omega \times \R^d)$, and hence, it is legitimate to talk of an initial condition $f_0$ (see~\cite[Remark~1.3]{ambrosio2014continuity}). In the following, we systematically consider the continuous version of such solutions.

        \medskip

	For the multiphasic system~\eqref{eq:multiphasic-equivalence2}, we ask more regularity with respect to space for the family of vector fields $\vv$. 
	
	\begin{Def}[Weak solutions of the multiphasic system]
		\label{def:solution_multiphasic}
		Let $(I,\mu)$ be a set of labels. Let $(\rhorho,\vv) = (\rhorho(t),\vv(t))$, $t \in [0,T]$, be a measurable family of multiphasic distributions. We say that $(\rhorho,\vv)$ is a solution of \eqref{eq:multiphasic-equivalence2} provided for $\mu$-almost all $\alpha$:
		\begin{itemize}
			\item For all compact $K \subset \Omega$, the map $t \in [0,T] \mapsto \rho^\alpha(t,K)$ belongs to $L^\infty$.
			\item The vector field $v^\alpha$ satisfies
			\begin{equation}
			\label{eq:lip_bound_v}
			\int_0^T \Lip (v^\alpha(t,\cdot)) \,  \D t < + \infty,
			\end{equation}
			and $v^\alpha \in L^\infty_{\mathrm{loc}}([0,T] \times \Omega)$.
			\item The pair $(\rho^\alpha, v^\alpha)$ solves both equations of~\eqref{eq:multiphasic-equivalence2} in the sense of distributions.
		\end{itemize}
		
		In this case, we can check that for $\mu$-almost all $\alpha$, $t \mapsto \rho^\alpha(t)$ is weak-$\star$ continuous in the dual of $\mathscr{C}_c(\Omega \times \R^d)$ and that $t \mapsto v^\alpha(t)$ is continuous in the topology of local uniform convergence (up to modifying them on a negligible set of times). Therefore, it is legitimate to talk of initial conditions $(\rhorho_0, \vv_0)$ for~\eqref{eq:multiphasic-equivalence2} as well.
	\end{Def}
	
	\begin{Rem}
		This definition is relevant for the following reasons. First, the first point (concerning $\rho^\alpha$) and the Lipschitz bound for $v^\alpha$ are somewhat standard conditions to give a distributional meaning to the continuity equation, that is the first equation in~\eqref{eq:multiphasic-equivalence2}. But it is not enough to give meaning to every term in the velocity equation. Indeed, with these properties, $v^\alpha$ is only $L^1_t L^\infty_{\mathrm{loc},x}$ and $\DDD v^\alpha \in L^1_t L^\infty_x$. Hence, $(v^\alpha \cdot \nabla) v^\alpha$ is not well defined. For the same reason,
		\begin{equation*}
		|F^\alpha(t,x)| = |F(t,x,v^\alpha(t,x))| \leq |F(t,0)| + \Lip (F(t, \cdot))  \Big\{ |x| + |v^\alpha(t,x)|\Big\}
		\end{equation*}
		is not well defined due to time integrability. 
		
		On the other hand, the condition $v^\alpha \in L^\infty_{\mathrm{loc}}([0,T] \times \Omega)$ solves both problems and is natural due to the following a priori estimate deduced from the second equation in~\eqref{eq:multiphasic-equivalence2}
		\begin{equation*}
		\frac{\D}{\D t} |v^\alpha(t,x)| \leq |F(t,0)| + \Lip (F(t, \cdot )) |x| + \Big\{ \Lip v^\alpha(t, \cdot) + \Lip (F(t, \cdot))  \Big\} |v^{\alpha}(t,x)|.
		\end{equation*}
	\end{Rem}
	
	\subsection{Existence and uniqueness for the uncoupled system \emph{via} the method of characteristics}
	\label{subsec:characterisics}
	Let $F: [0,T] \times \Omega \times \R^d \to \R^d$ with $\| F \|_{\mathcal L, T}< + \infty$ (see  \eqref{eq:def_lipschitz_norm}). Associated with $F$, we can define thanks to the Cauchy--Lipschitz theorem the solution of the Newton laws associated with $F$, that is following ODE starting from $(x,v) \in \Omega \times \R^d$, for all $t \in [0,T]$:
	\begin{equation*}
	\left\{  
	\begin{aligned}
	\dot X_t(x,v) &= V_t(x,v),& \dot V_t(x,v) &= F(t,X_t(x,v), V_t(x,v)),\\
	X_0(x,v) &= x
	,& V_0(x,v) &= v.
	\end{aligned}
	\right.
	\end{equation*}
	For  $(t,x,v) \in [0,T] \times \Omega \times \R^d$, we also introduce the notation $$\Phi_t(x,v) \vcentcolon = (X_t(x,v), V_t(x,v)).$$
    In the following two propositions, we will show that $(\Phi_t)_{t \in [0,T]}$ plays the role of characteristics for both equations~\eqref{eq:vlasov-equivalence2} and~\eqref{eq:multiphasic-equivalence2}.

	First, in the case of the Vlasov equation~\eqref{eq:vlasov-equivalence2}, this is easily deduced from the general theory of continuity equations.
	\begin{Prop}
		\label{prop:characteristics_Vlasov}
		Let $f_0$ be a Radon measure on $\Omega \times \R^d$. Equation~\eqref{eq:vlasov-equivalence2} admits a unique solution starting from $f_0$ in the sense of Definition~\ref{def:solution_Vlasov}, and for all $t \in [0,T]$ we have
		\begin{equation*}
		f(t, \cdot ) = \Phi_t {}\pf f_0.
		\end{equation*}
	\end{Prop}
	\begin{proof}
		The main observation is that~\eqref{eq:vlasov-equivalence2} rewrites exactly as the following continuity equation in the phase space~$\Omega \times \R^d$:
		\begin{equation}
        \label{eq:reduced_vlasov}
		\partial_t f + \Div_{x,v} (f \mathsf V) = 0,
		\end{equation}
		where $\mathsf V$ is the smooth and Lipschitz vector field given by
		\begin{equation*}
		\mathsf V(t,x,v) = \begin{pmatrix}
		v\\ F(t,x,v)
		\end{pmatrix}.
		\end{equation*}
        
        The result follows from standard arguments that can be found in~\cite{ambrosio2014continuity}. Let us briefly sketch the arguments for the sake of completeness.

        First, $t \mapsto f(t,\cdot) := \Phi_t {}\pf f_0$ is a solution of~\eqref{eq:vlasov-equivalence2} in the sense of Definition~\ref{def:solution_Vlasov}. Indeed, given a compact $K \subset \Omega \times \R^d$, let us call
        \begin{equation*}
            \overline K := \{ \Phi_t^{-1}(x,v); \, t \in [0,T], \, (x,v) \in K \}.
        \end{equation*}
        This set $\overline K$ is a compact because it is the image of the compact $[0,T] \times K$ by the continuous map $(t,x,v) \mapsto \Phi_t^{-1}(x,v)$. Therefore, 
        \begin{equation*}
            f(t,K) = f_0(\Phi_t^{-1}(K)) \leq f_0(\overline K)
        \end{equation*}
        is bounded uniformly in $t \in [0,T]$. Moreover, if $\varphi = \varphi(t,x,v) \in C_c^\infty((0,T) \times \Omega \times \R^d)$, we have
        \begin{align*}
            \frac{d}{dt} \int \varphi(t,x,v) &f(t,\D x, \D v) = \frac{d}{dt} \int \varphi(t,\Phi_t(x,v) )f_0(\D x, \D v)\\ &= \int \{ \partial_t \varphi(t,\Phi_t(x,v)) + \nabla \varphi(t,\Phi_t(x,v)) \cdot \dot \Phi_t(x,v) \} f_0(\D x, \D v)\\
            &= \int \{ \partial_t \varphi(t,\Phi_t(x,v)) + \nabla \varphi(t,\Phi_t(x,v)) \cdot \mathsf V(t,\Phi_t(x,v)) \} f_0(\D x, \D v)\\
            &= \int \{ \partial_t \varphi(t,x,v) + \nabla \varphi(t,x,v) \cdot \mathsf V (t,x,v) \} f(t,\D x, \D v).
        \end{align*}
        Hence, integrating this computation with respect to time, we find that $f$ is a distributional solution of~\eqref{eq:vlasov-equivalence2} rewritten as~\eqref{eq:reduced_vlasov}. 
        
		Let us finally prove uniqueness. Let $g$ be another solution of~\eqref{eq:vlasov-equivalence2} in the sense of Definition~\ref{def:solution_Vlasov} and $\varphi \in C_c^\infty(\Omega \times \R^d)$. The proof consists in showing that, the map
        \begin{equation*}
           a:  t \mapsto \int \varphi(\Phi_{t}^{-1}(x,v)) g(t,\D x, \D v)
        \end{equation*}
        is constant. If this is true, the equality $a(t) = a(0)$ exactly means that $f_0 = {(\Phi_t^{-1})} \pf g(t, \cdot)$, or differently stated, $g(t,\cdot) = \Phi_{t \, \#}  f_0$, that is $g = f$. To this end, we can derive $a$:
        \begin{align*}
            &\frac{d}{dt} a(t) =  \int \{\partial_t (\varphi(\Phi_t^{-1}(x,v)) +\nabla(\varphi(\Phi_t^{-1}(x,v)) \cdot \mathsf V(t,x,v) \} g(t, \D x, \D v),\\
            &= \int \left\{ \frac{d}{dt} \Phi_t^{-1}(x,v) + \DDD \Phi_t^{-1}(x,v) V(t,x,v) \right\} \cdot \nabla \varphi(\Phi_t^{-1}(x,v))g(t, \D x, \D v)\\
            &= \int \left\{ \frac{d}{dt} \Phi_t^{-1}(\Phi_t(x,v)) + \DDD \Phi_t^{-1}(\Phi_t(x,v)) V(t,\Phi_t(x,v)) \right\} \\
            &  \qquad \qquad \qquad \qquad \qquad \qquad \qquad \qquad \qquad \qquad  \cdot \nabla \varphi(x,v)(\Phi_t^{-1}) {}\pf g(t, \D x, \D v)\\
            &= \int \left\{ (\partial_t\Phi_t^{-1})(\Phi_t(x,v)) + \DDD \Phi_t^{-1}(\Phi_t(x,v)) \dot \Phi_t(x,v) \right\} \\
            & \qquad \qquad \qquad \qquad \qquad \qquad \qquad \qquad \qquad \qquad  \cdot \nabla \varphi(x,v)(\Phi_t^{-1}) {}\pf g(t, \D x, \D v),
        \end{align*}
        where for the first equality, we used the fact that $g$ is a distributional solution of~\eqref{eq:reduced_vlasov}. But, for all $(x,v)$,
        \begin{equation*}
            (\partial_t\Phi_t^{-1})(\Phi_t(x,v)) + \DDD \Phi_t^{-1}(\Phi_t(x,v)) \dot \Phi_t(x,v) = \frac{d}{dt} \Phi_t^{-1} \circ \Phi_t(x,v) = \frac{d}{dt} (x,v) = 0,
        \end{equation*}
        and the result follows.
	\end{proof}

	In the case of the multiphasic system, the link is slightly more intricate because of the possible blow-up of the multiphasic representation. The result goes as follows.
	\begin{Prop}
		\label{prop:characteristics_multiphasic}
		Let $(\rhorho_0, \vv_0)$ be a regular multiphasic distribution and set
		\begin{equation*}
		L \vcentcolon = \esssup_{\alpha \in I}  \| \DDD v_0^\alpha \|_{L^{\upinfty}_x}.
		\end{equation*}
		There exist $T^* \in (0,T]$ only depending on $L$ and $\|F\|_{\mathcal L, T}$ and a solution $(\rhorho,\vv)$ of~\eqref{eq:multiphasic-equivalence2} starting from $(\rhorho_0, \vv_0)$ in the sense of Definition~\ref{def:solution_multiphasic} up to time $T^*$. Moreover, this solution is unique and satisfies the following additional regularity property:
		\begin{equation}
		\label{eq:uniform_lip_bound}
		\esssup_{\alpha} \sup_{t \in T^*} \Big\{|v^\alpha(t,0)|  +  \| \DDD v^\alpha(t) \|_{L^\infty_x} \Big\} < + \infty.
		\end{equation}
		Finally, for a given $\alpha$, $t \in [0,T^*]$ and $x \in \Omega$, set
		\begin{equation*}
		\psi^\alpha(t,x) \vcentcolon = X_t(x,v_0^\alpha(x)).
		\end{equation*}
		Then, we have for $\mu$-almost all $\alpha$, for all $t \in [0,T^*]$ and all $x \in \Omega$:
		\begin{equation*}
		\rho^\alpha_t = \psi^\alpha(t, \cdot) {}\pf \rho_0^\alpha, \qquad \mbox{and} \qquad v^\alpha(t,\psi^\alpha(t,x)) = V_t(x,v_0^\alpha(x)).
		\end{equation*}
	\end{Prop}
	\begin{proof}
		The proof follows very standard arguments. We still provide them for the sake of completeness and to keep track of uniformity of the different quantities with respect to $\alpha$.	
		
		\noindent\underline{Step 1}: $\psi^\alpha$ remains a bi-Lipschitz homeomorphism for a positive time which is uniform in $\alpha$.
		
		To ease readability, we write $\psi_t^\alpha=\psi^\alpha(t, \cdot)$. It is easy to check that for a given $\alpha$ for which $v_0^\alpha$ is Lipschitz continuous (which is the case for $\mu$-almost all $\alpha$), $\psi^\alpha$ has enough regularity to justify the following computations. For such an $\alpha$ and given $x\neq y \in \Omega$ and $t \in [0,T]$ sufficiently small, let us call
		\begin{equation*}
		a(t) \vcentcolon = \log | \psi^\alpha(t,y) - \psi^\alpha(t,x) | \quad \mbox{and} \quad b(t) \vcentcolon = \frac{ | \partial_t \psi^\alpha(t,y) - \partial_t \psi^\alpha(t,x) |}{ | \psi^\alpha(t,y) - \psi^\alpha(t,x) |}.
		\end{equation*}
By standard computations we have
		\begin{equation*}
		\dot a(t) \leq b(t) \quad \mbox{and} \quad \dot b(t) \leq \Lip (F(t, \cdot))(1 + b(t)) + b(t)^2.
		\end{equation*}
		As $b(0) \leq \| D v^\alpha_0 \|_{L^\infty_x} \leq L$, by a  comparison principle and the Cauchy--Lipschitz theorem, there is $T^* \in (0,T]$ only depending on $L$ and $\|F\|_{\mathcal L, T}$ such that for all $t \leq T^*$, $b(t) \leq 2 L$. Therefore, we deduce that for $\mu$-almost all $\alpha$, we have for all $x,y \in \Omega$ and all $t \in [0,T^*]$:
		\begin{equation*}
		| y-x | \exp(-2Lt) \leq |\psi^\alpha(t,y) - \psi^\alpha(t,x)| \leq | y-x | \exp(2L t).
		\end{equation*}
		Hence $\psi^\alpha(t,\cdot)$ is injective, continuous and proper. By the invariance of domain theorem, we can conclude that for $\mu$-almost all $\alpha$, for all $t \in [0,T^*]$, $\psi^\alpha(t, \cdot)$ is a bi-Lipschitz homeomorphism of~$\R^d$. Let us call $\phi^\alpha(t,\cdot)$ the corresponding inverse. Note that this $\phi^\alpha(t, \cdot)$ is Lipschitz continuous.
		
		\bigskip
		
		\noindent \underline{Step 2}: Existence.
		
		Let us define for such an $\alpha$, $t \in [0,T^*]$ and $x \in \Omega$:
		\begin{equation*}
		v^\alpha(t,x) \vcentcolon =V_t\Big(\phi^\alpha(t,x), v^\alpha_0(\phi^\alpha(t,x))\Big) = \partial_t \psi^\alpha(t,\phi^\alpha(t,x)). 
		\end{equation*}
		This vector field is clearly locally Lipschitz continuous in $t$ and $x$. From this expression, it is classical to check that $v^\alpha$ satisfies almost everywhere and in the sense of distributions
		\begin{equation*}
		\partial_t v^\alpha(t,x) + (v^\alpha(t,x) \cdot \nabla)v^\alpha(t,x) = F(t,x,v^\alpha(t,x)).
		\end{equation*}
		The uniform bound~\eqref{eq:uniform_lip_bound} is a direct consequence on the bound we found for $b$ in the previous step.
		
		Eventually, as $\psi^\alpha$ is the flow associated with the Lipschitz vector field $v^\alpha$ (that is, $\partial_t \psi^\alpha(t,x) = v^\alpha(t,\psi^\alpha(t,x))$, the method of characteristics (for instance presented in the proof of Proposition~\ref{prop:characteristics_Vlasov}) shows that $\rho^\alpha(t) \vcentcolon = \psi^\alpha(t) \pf \rho^\alpha_0$ is a distributional solution of the continuity equation
		\begin{equation*}
		\partial_t \rho^\alpha + \Div(\rho^\alpha v^\alpha)=0,
		\end{equation*}
		hence the existence statement.
        
		\bigskip
		
		\noindent \underline{Step 3}: Uniqueness.
		
		Let $(\varrho^\alpha, w^\alpha )$ be another solution starting from $(\rhorho_0, \vv_0)$ in the sense of Definition~\ref{def:solution_multiphasic}. Let $\alpha$ be such that $\psi^\alpha$ is a bi-Lipschitz homeomorphism up to time $T^*$ and $w^\alpha$ satisfies~\eqref{eq:lip_bound_v}. Then, we have for all $x \in \Omega$ and $t \in [0,T^*)$,
		\begin{align*}
		\partial_t | v^\alpha(t&,\psi^\alpha(t,x)) - w^\alpha (t,\psi^\alpha(t,x)) | \\
		&\leq | F(t,\psi^\alpha(t,x), v^\alpha(t,\psi^\alpha(t,x))) - F(t,\psi^\alpha(t,x), w^\alpha(t,\psi^\alpha(t,x)))\\
		&\qquad + ((w^\alpha(t,\psi^\alpha(t,x)) - v^\alpha(t,\psi^\alpha(t,x)))\cdot \nabla)w^\alpha(t,\psi^\alpha(t,x)) |\\
		&\leq \Big( \Lip (F(t, \cdot))  + \Lip (w^\alpha(t, \cdot)) \Big) | v^\alpha(t,\psi^\alpha(t,x)) - w^\alpha (t,\psi^\alpha(t,x)) |.
		\end{align*} 
		Therefore, we must have $w^\alpha = v^\alpha$. Once this is done, $\varrho^\alpha = \rho^\alpha$  as they satisfy the same continuity equation (with Lipschitz velocity field, see the arguments in the proof of Proposition~\ref{prop:characteristics_Vlasov}). 
	\end{proof}
	
	\subsection{The link between the Vlasov and multiphasic equations: uncoupled case}
	\label{subsec:equivalence}
	
	We can finally state precisely the main result of this section, namely, the equivalence between~\eqref{eq:vlasov-equivalence2} and~\eqref{eq:multiphasic-equivalence2}, in the case where the force field $F$ is given and satisfies $\|F\|_{\mathcal L, T}< + \infty$.
	
	\begin{Thm}
		\label{thm:equivalence_models}
		\begin{enumerate}
			\item Let $(I,\mu)$ be a set of labels and $(\rhorho,\vv)$ a solution  of the multiphasic system~\eqref{eq:multiphasic-equivalence2} in the sense of Definition~\ref{def:solution_multiphasic} up to time $T$ with $(\rhorho(0), \vv(0))$ regular. If $\mathcal I_\mu(\rhorho_0, \vv_0)$ is a Radon measure, then $\mathcal I_\mu(\rhorho,\vv)$ is a solution of the Vlasov equation~\eqref{eq:vlasov-equivalence2} up to time $T^*$ in the sense of Definition~\ref{def:solution_Vlasov}.
			\item Let $f$ be a solution of the Vlasov equation~\eqref{eq:vlasov-equivalence2} in the sense of Definition~\ref{def:solution_Vlasov} up to time $T$. Then there exists $T^*\leq T$, $(I,\mu)$ a set of labels and $(\rhorho,\vv)$ a solution of the multiphasic system~\eqref{eq:multiphasic-equivalence2} in the sense of Definition~\ref{def:solution_multiphasic}, satisfying~\eqref{eq:uniform_lip_bound}, such that for all $t \leq T^*$, $f(t) = \mathcal I_\mu(\rhorho(t), \vv(t))$.
		\end{enumerate}
	\end{Thm}
	\begin{proof}
		\underline{First point}. Let $(\Phi_t)$ be defined as in Subsection~\ref{subsec:characterisics}. By the uniqueness part of Proposition~\ref{prop:characteristics_multiphasic}, $(\rhorho,\vv)$ is the solution starting from $(\rhorho_0, \vv_0)$ given in this proposition. Setting $f(t)\vcentcolon = \mathcal I_\mu(\rhorho(t), \vv(t))$, we have for all nonnegative test function $\varphi$ and all $t \leq T$, using the notation of Proposition~\ref{prop:characteristics_multiphasic}:
		\begin{align*}
		\int \varphi(x,v) f(t, \D x, \D v) &= \iint \varphi(x, v^\alpha(t,x))  \rho^\alpha(t,\D x)\D \mu(\alpha)\\
		&=\iint \varphi\big(\psi^\alpha(t,x), v^\alpha(t,\psi^\alpha(t,x)) \big)  \rho^\alpha_0(\D x)\D \mu(\alpha)\\
		&=\iint \varphi \big(X_t(x, v_0^\alpha(x)), V_t(x, v_0^\alpha(x)) \big)  \rho^\alpha_0(\D x)\D \mu(\alpha)\\
		&= \iint \varphi\circ\Phi_t(x, v_0^\alpha(x))  \rho^\alpha_0(\D x)\D \mu(\alpha)\\
		&= \int \varphi \circ \Phi_t \, f_0(\D x, \D v).
		\end{align*}
		Hence, $f(t) = \Phi_t {}\pf f_0$. Provided $f_0$ is a Radon measure, this is the unique solution of~\eqref{eq:vlasov-equivalence2} in virtue of Proposition~\ref{prop:characteristics_Vlasov}.
		
		\bigskip
		
		\noindent \underline{Second point}. Consider $(I,\mu)$ a set of labels and $(\rhorho_0, \vv_0)$ a regular multiphasic distribution such that $\mathcal I_\mu(\rhorho_0,\vv_0) = f_0$ as given by Proposition~\ref{prop:disintegration}. Let $T^*\leq T$ and $(\rhorho,\vv)$ be the solution of~\eqref{eq:multiphasic-equivalence2} in the sense of Definition~\ref{def:solution_multiphasic} starting from $(\rhorho_0, \vv_0)$ given by Proposition~\ref{prop:characteristics_multiphasic}. By the first point, as $\mathcal I_\mu(\rhorho_0, \vv_0)$ is Radon, $\I_\mu(\rhorho,\vv)$ is a solution of~\eqref{eq:vlasov-equivalence2} in the sense of Definition~\ref{def:solution_Vlasov} up to time $T^*$, starting from $f_0$. By the uniqueness part of Proposition~\ref{prop:characteristics_Vlasov}, this solution must coincide with $f$.
	\end{proof}
	
	\subsection{The link between the Vlasov and multiphasic equations: coupled case}
	\label{subsec:coupled}
	In this subsection, we gather consequences of the results obtained previously, in the case when $F$ depends on the solution of the Vlasov equation itself: that is, for all $t\geq 0$, there is a map $\Psi_t$ sending time dependent distribution $f = f(s)$, $s \in [0,t]$, to a vector field $F(s,\cdot)$.
	Note that this corresponds to both~\eqref{eq:VP_real} and~\eqref{eq:VNS_real}. This is not the point of this subsection to be precise on the properties of this map $\Psi_t$
	(this will be done in full details in the multiphasic case in the next section), but rather to compare properties of existence and uniqueness for the Vlasov and multiphasic systems. Namely, we want to compare the Vlasov equation
	\begin{equation}
	\label{eq:coupled_general_Vlasov}
	\left\{ 
	\begin{gathered}
	\partial_t f + v \cdot \nabla_x f + \mathrm{div}_v \big( F f \big) = 0,\\
	F(t,\cdot) = \Psi_t(f|_{[0,t]}),\\
	f|_{t=0} = f_0,
	\end{gathered}
	\right.
	\end{equation}
	and, given a set of labels $(I,\mu)$, the multiphasic system
	\begin{equation}
	\label{eq:coupled_general_multiphasic}
	\left\{
	\begin{gathered}
	\partial_t \rho^\alpha + \Div( \rho^\alpha v^\alpha) =0, \\
	\partial_t v^\alpha + ( v^\alpha \cdot \nabla ) v^\alpha = F(\cdot, \cdot , v^\alpha),\\
	F(t,\cdot) = \Psi_t(\mathcal I_\mu(\rhorho,\vv)|_{[0,t]}),\\
	\rho^\alpha |_{t = 0} = \rho^\alpha_0, \ v^\alpha |_{t = 0} = v^\alpha_0,
	\end{gathered}
	\right.
	\end{equation}
	under the consistency condition
	\begin{equation}
	\label{eq:consistency_initial_condition}
	f_0 = \mathcal I_\mu(\rhorho_0, \vv_0).
	\end{equation}
	
	With the help of Definition~\ref{def:solution_Vlasov} and Definition~\ref{def:solution_multiphasic},  we can precisely define the notion of solutions of systems~\eqref{eq:coupled_general_Vlasov} and~\eqref{eq:coupled_general_multiphasic} that we want to consider.
	\begin{Def}
		\label{def:solution_coupled}
		\begin{itemize}
			\item We call a solution of the Vlasov equation~\eqref{eq:coupled_general_Vlasov} up to time $T>0$ any measurable curve $f = f(t)$, $t \in [0,T]$ valued in the space of nonnegative Radon measures of $\Omega \times \R^d$ such that the force field $F$ defined by the second equation in~\eqref{eq:coupled_general_Vlasov} satisfies $\| F \|_{\mathcal L, T} < + \infty$ and $f$ solves the first equation and satisfies the initial condition in the sense of Definition~\ref{def:solution_Vlasov}.
			\item Let $(\rhorho,\vv) = (\rhorho(t),\vv(t))$, $t \in [0,T]$ be a measurable family of multiphasic distributions. We say that $(\rhorho,\vv)$ is a solution of~\eqref{eq:coupled_general_multiphasic} up to time $T>0$ provided the force field $F$ defined by the third equation satisfies $\| F \|_{\mathcal L, T} < + \infty$ and $(\rhorho,\vv)$ solves the first two equations and satisfies the initial condition in the sense of Definition~\ref{def:solution_multiphasic}.
		\end{itemize}
	\end{Def}
	Let us emphasize the fact that the Lipschitz regularity assumption on the force field is enforced in the definition of solutions for the systems~\eqref{eq:coupled_general_Vlasov} and~\eqref{eq:coupled_general_multiphasic}.
	
	\medskip
	
	First, existence for one system implies existence for the other, and vice versa. 
    
	\begin{Prop}[Existence results for the two systems are equivalent]
		\label{prop:existence_coupled_case}
		Let $f_0$ be a nonnegative Radon measure on $\Omega \times \R^d$ and $(\rhorho_0, \vv_0)$ be a regular multiphasic distribution, such that~\eqref{eq:consistency_initial_condition} holds.
		\begin{enumerate}
			\item If there exists a solution $f$ of~\eqref{eq:coupled_general_Vlasov}, then there exists $0<T^*\leq T$ and $(\rhorho,\vv)$ a solution of~\eqref{eq:coupled_general_multiphasic} up to time $T^*$. In addition, for all $t \leq T^*$, $f(t) = \I_\mu(\rhorho(t),\vv(t))$.
			\item If there exists a solution $(\rhorho,\vv)$ a solution of~\eqref{eq:coupled_general_multiphasic}, then $f \vcentcolon = \I_\mu(\rhorho,\vv)$ is a solution of~\eqref{eq:coupled_general_Vlasov} up to time $T$.
		\end{enumerate}
	\end{Prop}
	\begin{proof}
		\underline{First point}. Let us define $F$ as in the second line of~\eqref{eq:coupled_general_Vlasov}. Let $0<T^* \leq T$ and $(\rhorho,\vv)$ be the solution of~\eqref{eq:multiphasic-equivalence} with this specific $F$, given by Proposition~\ref{prop:characteristics_multiphasic}. By Theorem~\ref{thm:equivalence_models}, $\bar f \vcentcolon = \mathcal I_\mu(\rhorho,\vv)$ is a solution of~\eqref{eq:vlasov-equivalence} up to $T^*$ with this choice of $F$. But there is only one of these solutions in virtue of Proposition~\ref{prop:characteristics_Vlasov}, and therefore, $\bar f = f$. In other words, $f = \mathcal I_\mu(\rhorho,\vv)$, and therefore for all $t \leq T^*$, $F(t, \cdot) = \psi_t(\mathcal I_\mu(\rhorho,\vv)|_{[0,t]})$, and $(\rhorho,\vv)$ is a solution of~\eqref{eq:coupled_general_multiphasic}.
		
		\bigskip
		
		\noindent \underline{Second point}. Defining $F$ as in the third line of~\eqref{eq:coupled_general_multiphasic}, $f \vcentcolon = \mathcal I_\mu(\rhorho,\vv)$ is a solution of the first equation of~\eqref{eq:coupled_general_Vlasov} as a consequence of Theorem~\ref{thm:equivalence_models}. Therefore, $f$ is a solution of \eqref{eq:coupled_general_Vlasov}. 
	\end{proof}

	Similarly, uniqueness for one system implies uniqueness for the other, and vice versa. 
	
	\begin{Prop}[Uniqueness results for the two systems are equivalent]
		\label{prop:uniqueness_coupled_case}
		Let $(I,\mu)$ be a set of labels. Let $f_0$ be a nonnnegative Radon measure on $\Omega \times \R^d$ and $(\rhorho_0, \vv_0)$ be a regular multiphasic distribution, such that~\eqref{eq:consistency_initial_condition} holds and let $T>0$. 
		\begin{enumerate}
			\item If for all $0\leq t \leq T$, there is at most one solution of the Vlasov equation~\eqref{eq:coupled_general_Vlasov} starting from $f_0$ up to time $t$, then for all $0\leq t \leq T$, there is at most one solution of the multiphasic system~\eqref{eq:coupled_general_multiphasic} starting from $(\rhorho_0, \vv_0)$ up to time $t$.
			\item If the multiphasic system~\eqref{eq:coupled_general_multiphasic} admits a solution $(\rhorho,\vv)$ starting from $(\rhorho_0, \vv_0)$ up to time $T$ such that for all $t \leq T$, $(\rhorho,\vv)|_{[0,t]}$ is the unique solution of~\eqref{eq:coupled_general_multiphasic} up to time $t$, then for all $t \leq T$, $\I_\mu(\rhorho,\vv)|_{[0,t]}$ is the unique solution of the Vlasov equation~\eqref{eq:coupled_general_Vlasov} starting from $f_0$ up to time $t$.
		\end{enumerate}
	\end{Prop}
\begin{proof}
	\underline{First point}. Let us consider $0\leq t \leq T$ and two solutions $(\rhorho,\vv)$ and $(\boldsymbol \varrho, \ww)$ of~\eqref{eq:coupled_general_multiphasic} starting from $(\rhorho_0, \vv_0)$ up to time $t$. Let us call respectively $F$ and $G$ the corresponding force fields defined by the third equation in~\eqref{eq:coupled_general_multiphasic}. Then, by Proposition~\ref{prop:existence_coupled_case}, $\I_\mu(\rhorho,\vv)$ and $\I_\mu(\boldsymbol \varrho,\ww)$ are two solutions of~\eqref{eq:coupled_general_Vlasov} up to time $t$, of respective force field $F$ and $G$. Therefore, by assumption, they are the same and hence $F = G$. But by Proposition~\ref{prop:characteristics_multiphasic}, once  the force field is fixed, solutions of~\eqref{eq:multiphasic-equivalence} are unique and therefore, $(\rhorho,\vv) = (\boldsymbol \varrho, \ww)$. The result follows.
	
	\bigskip
	
	\noindent \underline{Second point}. Let $0\leq t \leq T$ and $f$ be a solution of~\eqref{eq:coupled_general_Vlasov} starting from $f_0$, up to time $t$. Relying on Proposition \ref{prop:existence_coupled_case}, let $T^* \leq t$ and $(\boldsymbol \varrho, \ww)$ be the solution of~\eqref{eq:coupled_general_multiphasic} starting from $(\rhorho_0, \vv_0)$ up to time $t$, such that for all $s \leq t$, $f(s) = \I_\mu(\boldsymbol \varrho(s), \ww(s))$. By assumption, for all $s \leq t$, $(\rhorho(s), \vv(s)) = (\boldsymbol \varrho(s), \ww(s))$, and hence $f(s) = \I_\mu(\rhorho(s), \vv(s))$, and the conclusion follows.
\end{proof}

	\section{Summary of the structure of the monograph}
To conclude this introductory chapter, we provide a quick summary of the content of the three main chapters of this monograph.

\medskip

\textbf{Chapter \ref{Part1-LWP}} is organized as follows. In Sections \ref{sec:abstract_existence}--\ref{sec:conseq-rho}, we start by developing an abstract local well-posedness theory in Sobolev spaces for a general class of multiphasic systems.
Consequences at the level of kinetic equations are detailed in Section~\ref{sec:conseq-kin}; we also detail an application to the study of monokinetic limits and provide a comparison with classical (and less classical) results of local in time propagation of regularity for Vlasov equations. In the last section of the Chapter, namely Section~\ref{section:Extension}, we generalize the functional framework in two directions: weighted Sobolev spaces and Besov spaces.

\medskip

 \textbf{Chapter \ref{Part2-VP}} focuses on Vlasov--Poisson type systems, which are first introduced in Section~\ref{Section-PresentationVPoissonandco}. The abstract local well-posedness theory of Chapter~\ref{Part1-LWP} is  applied to the different models, in Sections~\ref{Subsec:VPelectrons}, \ref{Subsec:VPions} and \ref{Subsec:VMongeAMpere} (in the latter, regularity and stability results for the Monge--Ampère equation are also proved see Subsection~\ref{Appendix-MongeAmpère}). In Section~\ref{sec:semiclassical}, we consider the semiclassical Schr\"odinger--Poisson system for mixed states and provide a multiphasic derivation of the Vlasov--Poisson equation. The last sections of the chapter are dedicated to instability issues. First, an abstract framework is provided in Section~\ref{Section-FrameworkInstab}; it is then thoroughly applied to the actual systems in Section~\ref{sec:insta-appli}.

\medskip

\textbf{Chapter \ref{Part3-VNS}} concerns Vlasov--Navier--Stokes type systems, which are first introduced in Section~\ref{Section-PresentationVNSandco}. The abstract local well-posedness theory of Chapter~\ref{Part1-LWP} is  applied to the different models, in Sections~\ref{SubsecVNSmultiphaseLWP}, \ref{SubsecVNScomp-multiphaseLWP}, and \ref{subsec:Stokes}. The last part of the chapter is then dedicated to the proof of global existence and stability around constant monokinetic profiles, in Sections~\ref{Section:monokinetic-stab-theorem} (for the statements and general strategy) and Sections \ref{sectionVNS-Conslaws-Decay}, \ref{Subsection-preuveStab-VNS}, and \ref{Section:VNS-comp} for the actual proofs.

\chapter{Abstract local well-posedness theory}\label{Part1-LWP}

In this chapter, we study the local well-posedness for the multiphasic formulation of a large class of Vlasov equations.  The outline of the chapter is the following.

\begin{itemize}
    \item

In Section~\ref{sec:abstract_existence}, we introduce an abstract framework for local well-posedness, which is put forward in light of the Vlasov--Poisson and the Vlasov--Navier--Stokes systems. The first idea is to focus solely on the momentum equation, seeing the forcing as a function of the momentum (and of possible other fixed parameters). With this point of view, the density equation is thus seen as an auxiliary equation. As an underlying principle,  one of the main structural properties we ask the associated Vlasov equation to satisfy is  that the force field gains one derivative in space with respect to the regularity of the family of densities from the multiphasic decomposition. As a consequence, this essentially translates into asking that the force has the same regularity in space as the momentum. 

In this section, local well-posedness is achieved in Sobolev spaces (based on $L^2_x$) of high enough regularity index, we refer to Theorem~\ref{thm:existence}; the proof is an adaption to the multiphasic setting of the standard procedure for local conservation laws.

\item The density equation in the multiphasic framework is then studied in Section~\ref{sec-continuity-eq}. This allows to leverage the general results of Section~\ref{sec:abstract_existence} to prove local well-posedness at the level of the density--velocity multiphasic formulation; this is the object of Section~\ref{sec:conseq-rho}.

\item Consequences on the associated Vlasov equations are drawn in Section~\ref{sec:conseq-kin}, where we state local existence and uniqueness results, with propagation of regularity at the level of moments in velocity. Relevant comparison with the existing literature is provided. Along the way, we obtain an application to the study of the so-called monokinetic limit for general Vlasov equations, that is developed in Section~\ref{sec:mono}.

\item In Section~\ref{section:Extension}, we extend these abstract results from  Sobolev spaces to weighted Sobolev spaces and to Besov spaces. The latter in particular allow to reach sharper results in terms of regularity.

\item Finally, in Section~\ref{subsec:summary-ass}, we provide a summary of the different sets of abstract assumptions (which pertain to the different frameworks) which lead to well-posedness results, as provided in this chapter.

\end{itemize}

\section{General framework and main result}\label{sec:abstract_existence}

	\subsection{The abstract system under study}\label{subsec:abstract-setting}
	The purpose of this section is to introduce a general abstract setting encompassing at least the multiphasic formulation \eqref{eq:multiphasic-equivalence} of the Vlasov--Poisson system~\eqref{eq:VP_real} and of the Vlasov--Navier--Stokes system~\eqref{eq:VNS_real}. Given a set of labels $(I,\mu)$, we shall focus on the following equation of unknown $\vv = (v^\alpha)_{\alpha \in I}$:
	\begin{equation}
	\label{eq:general_system}
	\left\{
	\begin{aligned}
	\partial_t v^\alpha + ( v^\alpha \cdot \nabla ) v^\alpha &= F^\alpha[\vv,X].\\
	v^\alpha |_{t = 0} &= v^\alpha_0,
	\end{aligned}
	\right.
	\end{equation}
	In this system, the force field $\boldsymbol F=(F^\alpha)_{\alpha\in I}$ is assumed to depend on some extra variable $X$ which lies in some space $\mathcal X$ (the role of this $X$ will be explained shortly after) and is computed using \emph{the whole family} $\vv$, in a possibly nonlinear and nonlocal way, in time and space, with some assumptions to be detailed later. Note that contrary to the multiphasic system \eqref{eq:multiphasic-equivalence} (which is most of the time the true system we want to study), there is no reference to the densities $\rhorho = (\rho^\alpha)_{\alpha \in I}$ and to the continuity equation. As a matter of fact, the effect of the latter is encoded in the extra quantity $X$. We made this choice because the velocity field is somehow pivotal in the local existence result, as it is the most convenient quantity to use in order to write a fixed-point argument.

 \medskip

\begin{Rem}[The Vlasov--Poisson and the Vlasov--Navier--Stokes cases] Let us explain how to compute the value of $\boldsymbol F$ at time $t\geq 0$ knowing $\vv$ in the set of times $[0,t]$ in the Vlasov--Poisson case -- see \eqref{eq:VP-multiphaseINTRO}. Given $\vv$, the value of $\boldsymbol F$ at time $t$ is in this case computed by:
	\begin{enumerate}
		\item Solving the continuity equation for each phase in order to get the family 
	$\rhorho(t, \cdot)$ of densities at time~$t$;
	\item Solving the Poisson equation in order to get the potential $U(t, \cdot)$ at time $t$;
	\item Setting $F^\alpha[\vv,X](t, \cdot) = - \nabla U(t, \cdot)$, for all $\alpha \in I$.
	\end{enumerate}
However, there is still some data missing, as to solve the continuity equation, we need the initial condition for the densities. This is where the parameter $X$ comes into play. In the Vlasov--Poisson case, we consider precisely $X \vcentcolon = \rhorho_0$ in some relevant functional space.

Similarly, in the Vlasov--Navier--Stokes case -- see \eqref{eq:VNS-multiphaseINTRO} -- computing $\boldsymbol F$ at time $t$ from the knowledge of $\vv$ in the set of times $[0,t]$ is done by:
	\begin{enumerate}
	\item Solving the continuity equation for each phase in order to get the family 
	$\rhorho$ in the set of times $[0,t]$;
	\item Solving the incompressible Navier--Stokes equation using the values of $(\rhorho,\vv)$ in the set of times $[0,t]$ in order to get the velocity of the fluid $U(t, \cdot)$ at time $t$;
	\item Setting $F^\alpha[\vv,X](t, \cdot) = v^\alpha(t, \cdot) - U(t, \cdot)$, for all $\alpha \in I$.
\end{enumerate}
This time, some other information is missing: the initial condition of the fluid. In this case, we take $X \vcentcolon = (\rhorho_0,U_0)$ in some relevant space. The other models that we study are treated similarly. 
\end{Rem}

	The main goal of this section is to derive a local existence and uniqueness theory of classical solutions for \eqref{eq:general_system}.
     In view of the equivalence results of Section~\ref{sec:equivalence_models}, it would seem natural to work with Sobolev spaces based on $L^\infty$; however the models we have in mind (namely Vlasov--Poisson and Vlasov--Navier--Stokes) involve operators that do not behave well with respect to $L^\infty$. For this reason, instead, we consider
    data in Sobolev spaces of high regularity, based on $L^2$. This choice is mainly made for the sake of simplicity; we shall later develop a generalization to other functional spaces, such as weighted Sobolev spaces or Besov spaces. 
    
    To do so, we need to give some assumptions on the family of force fields $\FF$. Loosely speaking, we assume that for appropriate $X$, \emph{$\vv$ and $\FF[\vv,X]$ have the same spatial regularity}. However, because we work with force fields computed in a non-local way with respect to time, we cannot take as an assumption that a family of velocity fields in $H^k_x$ provides a family of force fields in $H^k_x$ \emph{for a fixed time}. Instead, we will assume that a family of  velocity fields which is integrable in time with values in $H^k_x$ provides a family of force fields having the same property. Once this assumption is properly stated, the proof will be rather standard, following more or less the classical strategy for local conservation laws, see~\cite[Chapter~2]{majda1984compressible}. The main novelty is that we treat infinitely many phases at once, and hence integrability with respect to $\alpha$ needs to be taken into account.
	
	In general, in the applications, when dealing with multiphasic solutions, we do not want to impose $\vv = (v^\alpha)_{\alpha \in I}$ to be uniformly bounded with respect to the label $\alpha$ (see Remark~\ref{rem:bounded_velocity}). We want to allow some phases to have very large velocities, and only impose some integrability with respect to\ $\alpha$ for the quantity~$\| v^\alpha \|_{L^2_x}$. On the other hand, because of the convective term $(v^\alpha \cdot \nabla) v^\alpha$, blow-ups may occur when some shock appears, that is, when the quantity $\|\DDD v^\alpha \|_{L^{\upinfty}_x}$ becomes infinite for \emph{one phase}. As a result, contrarily to $\| v^\alpha \|_{L^2_x}$, the latter needs to be controlled uniformly with respect to~$\alpha$.
	
	For this reason, we introduce the following norms for families of functions $\vv = (v^\alpha)_{\alpha \in I}$.
	
	\begin{Def}
    \index{H@$\mathcal H^{k,p}$: multiphasic Sobolev space with regularity index $k$ and integrability parameter $p$}
	Given an integer $k \in \N^*$ and an integrability parameter $p \in [1,+\infty]$, we set for $\vv = (v^\alpha)_{\alpha \in I}$,
	\begin{equation}
	\label{def:norm}
	\Vert \vv \Vert_{\mathcal{H}^{k,p}} \vcentcolon = \| \vv \|_{L^p_\alpha L^2_x} + \| \DDD \vv \|_{L^{\upinfty}_\alpha H^{k-1}_x},
	\end{equation}
and
\begin{equation}
	\label{def:space}
	\mathcal H^{k,p} \vcentcolon = \left\{ \vv, \, \Vert \vv \Vert_{\mathcal{H}^{k,p}} <+\infty \right\}.
	\end{equation}
	\end{Def}
	Note that the $\mathcal{H}^{k,p}$ norm  is meant to control derivatives up to order $k$, hence the shift $k-1$ for the Sobolev norm bearing on $\DDD \vv$. The key ingredient for local existence will be an a priori estimate in the norm $\mathcal{H}^{k,p}$ for $k$ large enough.
\begin{Rem}
 Another possible route we could have followed, but which we do not develop in this work, would be to work with the space  $L^p_\alpha L^\infty_x$ instead of $L^p_\alpha L^2_x$  in \eqref{def:norm}. This would correspond to a multiphasic version of the Zhidkov spaces (for which we refer to \cite{Zhi,Gallo} and \cite{Carles-book} for a complete overview of their use for NLS equations and their hydrodynamic formulation).
\end{Rem}

\bigskip

Concerning the set of parameters $\mathcal X$, we want to cover a large variety of possible applications. Therefore, we assume that $\mathcal X$ is a set on which are defined:
\begin{itemize}
	\item a family of nonnegative map $\mathsf N_k : \mathcal X \to \R_+$ for all $k \in \N$, such that $\mathsf N_{k_1} \lesssim_{k_1,k_2} \mathsf N_{k_2}$ for $k_1\leq k_2$;
	\item a metric $\mathsf d$.
\end{itemize}
In the applications, $\mathcal X$ will be a functional space or a subspace of it, the $\mathsf N_k$ will be norms or semi-norms measuring high regularity, and $\mathsf d$ will be a distance measuring low regularity.

	\subsection{Assumptions}\label{Subsec:AssumptionF}
	
As already said in the introduction, the spatial domain for the $x$ variable is always  $\T^d$ or $\R^d$. As usual in the context of systems of conservation laws, the proof of local existence consists in introducing a scheme of fixed-point type, and then to show two key properties of this scheme: a \textbf{uniform bound in high regularity} in order to get compactness, and  \textbf{contraction at low regularity} in order to show that the scheme has a unique limit, which is therefore necessarily a fixed-point, that is, a solution of our system.
	
	To be able to perform such a proof, we will hence have to give two kinds of assumptions on the force fields $\FF$, namely, a bound in high regularity, and uniform regularity in low regularity, with an associated modulus of continuity. This is what we provide in this subsection. But first, let us define properly  the nature of this family of force fields. In this definition, the sign $\infty$ will denote a point "at infinity" in the functional space we will consider, and will hence be assumed to have infinite norm.
	
	\begin{Def}[Families of force fields]
		\label{def:force_fields}
        \index{F@$\boldsymbol{F}$: family of force fields}
		The object $\boldsymbol F = \boldsymbol F[\vv, X](t)$ is called a family of force fields if there exists an integer $k_0$ and $p \in [1, +\infty]$ such that for all $X \in \mathcal X$ and $t \in \R_+$, $\boldsymbol F [\cdot, X](t)$ is a map
		\begin{equation*}
		L^1([0,t]; \mathcal H^{k_0,p}) \longrightarrow \mathcal H^{k_0,p} \cup \{ \infty \}.
		\end{equation*}
		In the following, with a slight abuse of notation, if $0 \leq t \leq T$, $X \in \mathcal X$ and $\vv \in L^1([0,t]; \mathcal H^{k_0,p})$, we will denote
		\begin{equation*}
		\boldsymbol F\left[\vv|_{[0,t]},X\right](t) \qquad \mbox{by} \qquad \boldsymbol F[\vv,X](t). 
		\end{equation*}
	\end{Def} 
\begin{Rem}
	In the terminology of stochastic processes theory, $\boldsymbol F[\cdot,X]$ can be viewed as an adapted process defined on $L^1_{\mathrm{loc}}(\R_+, \mathcal H^{k_0,p})$ endowed with its natural filtration. It means that $\boldsymbol F [\vv, X](t)$ only depends on the value of $\vv$ up to time $t$, or in other words that if $\vv,\ww \in L^1([0,T];\mathcal H^{k_0,p})$ are such that $\vv|_{[0,t]} = \ww|_{[0,t]}$ for some $t \leq T$, then $\boldsymbol F [\vv, X](t) = \boldsymbol F [\ww, X](t)$. This property is essential in order to restrict solutions, that is, to ensure that if $\vv$ is a solution up to time $T>0$ and if $t \in [0, T]$, then $\vv|_{[0,t]}$ is a solution as well.
\end{Rem}
	
	As all the estimates can be nonlinear in a way that we do not want to specify, we introduce a notion of \emph{estimate functions} as follows:
	
	\begin{Def}
	\label{def:estimate_func}
		Let $\F:  \R_+ \times \R_+ \times \R_+\to \R_+ \cup \{ + \infty \}$ be a given function. We say that $\F$ is an \emph{estimate function} provided:
		\begin{enumerate}[label=(\roman*),leftmargin=50pt]
			\item $\F$ is non-decreasing with respect to each variable.
			\item The domain of $\F$ defined by $$\Omega_\F \vcentcolon = \{ (t,M,P) \mbox{ s.t.\ }\F(t,M,P) < + \infty \}$$ is open and $(0,0,0) \in \Omega_\F$.
			\item $\F$ is smooth on its domain $\Omega_\F$.
			\item If $P\in \R_+$ is such that $(0,0,P) \in \Omega_{\mathcal F}$, then $\F(0,0,P) = 0$.
		\end{enumerate}
		We say that it is a \emph{stability estimate function} if in addition to the previous points, $\F$ satisfies:
		\begin{enumerate}[label=(\roman*),leftmargin=50pt,resume]
			\item For all $t\geq 0$ such that $(t,0,0) \in \Omega_{\F}$, we have $\F(t,0,0)=0$.
		\end{enumerate}
	\end{Def} 

	With this notion at hand, we are able to state the general assumptions for our abstract framework. The first one asserts that if $\vv$ and $\mathsf N(X)$ are not too large, then $\FF[\vv,X]$ is well defined and has roughly the same regularity (in space) as $\vv$.
	
	\begin{Ass}{A1}[Bound in high regularity]
		\label{ass:bound_high_regularity}
        \index{A1@Assumptions \ref{ass:bound_high_regularity}--\ref{ass:stability}: for the general $\vv$ framework ($H^k$ spaces)}
		Let $\boldsymbol F$ be a family of force fields and $k_0,p$ be as in Definition~\ref{def:force_fields}. We assume that for all $k \geq k_0$, there exists an estimate function $\F$ such that for all $T>0$, for all $\vv \in L^1([0,T]; \mathcal H^{k,p})$, and all parameters $X \in \mathcal X$,
		we have
		\begin{equation}
		\label{eq:ass_estimate_high_reg}
		\int_0^T \Vert \FF[\vv,X](t)\Vert_{\mathcal{H}^{k,p}} \D t \leq \F\left(T, \int_0^T \Vert \vv(t)\Vert_{\mathcal{H}^{k,p}} \D t, \mathsf N_k(X) \right).
		\end{equation}
	\end{Ass}

	\begin{Rem}
		Recall that in Definition~\ref{def:force_fields}, we allowed $\FF[\vv,X](t)$ to take the value $+\infty$. However, in estimate~\eqref{eq:ass_estimate_high_reg}, if the right-hand side is finite, then necessarily, $\FF[\vv,X](t) \neq \infty$ for almost all $t \in [0,T]$. In practice, we might be able to obtain a control for the left-hand side of \eqref{eq:ass_estimate_high_reg} only for small times $T>0$, however an estimate of the form of  \eqref{eq:ass_estimate_high_reg} will be available as the estimate function $\mathcal{F}$ is allowed to take the value $+\infty$.
	\end{Rem}
\begin{Rem}\label{rem-bound-high-ref-functionQ}
    As will be seen from the coming proof, we can also allow for a natural generalization of Assumption \eqref{ass:bound_high_regularity}, replacing \eqref{eq:ass_estimate_high_reg} by the estimate
    \begin{equation}
		\int_0^T \Vert \FF[\vv,X](t)\Vert_{\mathcal{H}^{k,p}} \D t \leq \F\left(T, \int_0^T \mathcal{Q}(\Vert \vv(t)\Vert_{\mathcal{H}^{k,p}}) \D t, \mathsf N_k(X) \right),
		\end{equation}
        where $\mathcal{Q}:\R^+ \mapsto \R^+$ is a smooth non-decreasing function such that $\mathcal{Q}(0)=0$.
\end{Rem}
	
	The second assumption asserts that if $\vv_1, \vv_2$ and $X_1,X_2$ are bounded in a strong sense and close in weak sense, then $\FF[\vv_1, X_1]$ and $\FF[\vv_2, X_2]$ are close at low regularity as well.
	\begin{Ass}{A2}[Stability at low regularity]
		\label{ass:stability}
		Let $\boldsymbol F, k_0,p$ be as in Definition~\ref{def:force_fields}, and let $R>0$. There exists a stability estimate function $\G_R$ such that for all $T>0$, families $\vv_1,\vv_2 \in L^1([0,T]; \mathcal H^{k_0,p})$, and parameters $X_1$ and $X_2$, if
		\begin{equation}
		\label{eq:uniform_bound_velocities} \max\left( \sup_{t \in [0,T]} \Vert \vv_1(t) \Vert_{\mathcal{H}^{k_0,p}},\sup_{t \in [0,T]} \Vert \vv_2(t) \Vert_{\mathcal{H}^{k_0,p}}, \mathsf N_{k_0}(X_1), \mathsf N_{k_0}(X_2)  \right) \leq R,
		\end{equation}
	then the following estimate holds:
		\begin{equation}
		\label{eq:stability_assumption}
        \begin{aligned}
		\int_0^T \big\|\FF[\vv_1,X_1](t) &- \FF[\vv_2,X_2](t)\big\|_{L^p_\alpha L^2_x} \D t \\
        &\leq \G_R\left(T, \int_0^T \big\|\vv_1(t) - \vv_2(t) \big\|_{L^p_\alpha L^2_x} \D t, \mathsf d(X_1,X_2)\right).
         \end{aligned}
		\end{equation}
	\end{Ass}
\begin{Rem}
	In Assumption~\ref{ass:stability}, we set by convention $\| \infty - \boldsymbol f \|_{L^p_\alpha L^2_x} = + \infty$ for all $\boldsymbol{f} \in \mathcal H^{k_0,p} \cup \{\infty \}$. In particular, in estimate~\eqref{eq:stability_assumption}, if the right-hand side is finite, then $\FF[\vv_1,X_1](t) \neq \infty$ and $\FF[\vv_2,X_2](t) \neq \infty$ for almost all $t \in [0,T]$.
\end{Rem}

	\subsection{Definition of a solution and statement of the main theorem}
	Our main abstract result about local well-posedness for \eqref{eq:general_system} in the sense of Hadamard is the following.
	\begin{Thm}
		\label{thm:existence}
		Under Assumptions \ref{ass:bound_high_regularity} and \ref{ass:stability} for given $k_0 \in \N$ and $p \in [1,+\infty]$, there exists $k_1 > \max(1+d/2, k_0)$ such that for all $k \geq k_1$, given $\vv_0\in \mathcal H^{k,p}$ and $X \in \mathcal X$ such that $(0,0, \mathsf N(X))\in \Omega_\F$, there exist $T>0$ and a unique solution $\vv \in L^\infty([0,T]; \mathcal H^{k,p})$ of~\eqref{eq:general_system} up to time $T$. In addition, we have:
		\begin{itemize}
			\item \underline{Time regularity}. This solution satisfies $\vv \in \mathscr{C}^0([0,T]; L^p_\alpha L^2_x)$ and for $\mu$-almost every $\alpha \in I$, $v^\alpha \in \mathscr{C}^0([0,T]; H^k_x)$.
			\item \underline{Stability}. If $\vv_1$ and $\vv_2$ are two such solutions, defined up to time $T$, and corresponding to initial conditions $\vv_{1,0}$ and $\vv_{2,0}$ and parameters $X_1$ and $X_2$ respectively, then setting 
			\begin{equation*}
			R \vcentcolon = \max \left(\sup_{t \in [0,T]}\Vert \vv_1(t) \Vert_{\mathcal{H}^{k,p}},\sup_{t \in [0,T]} \Vert \vv_2(t) \Vert_{\mathcal{H}^{k,p}}, \mathsf N_{k_0}(X_1), \mathsf N_{k_0}(X_2)  \right),
			\end{equation*}
there holds for all $t \in [0,T]$
\begin{multline}
\label{eq:stability-vv}
	\| \vv_2(t) - \vv_1(t) \|_{L^p_\alpha L^2_x} \\
    \leq  \exp(Rt)  \left\{ \| \vv_{2,0} - \vv_{1,0} \|_{L^p_\alpha L^2_x} + \mathcal{G}_{R} \left(t, 2Rt, \mathsf d(X_1,X_2) \right)\right\}.
	\end{multline}
			\item \underline{Blow-up}. There exists a non-increasing function $T_{\mathrm{inf}}: \R^+ \to \R^+$ such that the following holds. Let $\vv \in L_{\mathrm{loc}}^\infty([0,T^{\star}); \mathcal H^{k,p})$ be a maximally extended solution in $\mathcal H^{k,p}$ (that is, such that $\vv$ is not the restriction to $[0,T^*)$ of a solution $\ww \in L_{\mathrm{loc}}^\infty([0,T]; \mathcal H^{k,p})$ for some $T> T^*$). Then  $T^{\star} \geq T_{\mathrm{inf}}(\Vert \vv_0 \Vert_{\mathcal{H}^{k,p}})$ and we have
			\begin{equation*}
			T^{\star} = + \infty \quad \text{or} \quad  \esssup_{\alpha}\int_0^{T^{\star}} \Big\{ \Vert \nabla  v^{\alpha}(s) \Vert_{L^{\upinfty}_x} + \Vert F_\alpha[\vv(s), X] \Vert_{H^k_x} \Big\}\D s=+\infty.
			\end{equation*}
		\end{itemize}
	\end{Thm}

	\begin{Rem}
    \label{Rem:timeuniform}
	As a by-product of our proof, we actually show the following refined statement related to Theorem \ref{thm:existence}: for all $M>\Vert \vv_0 \Vert_{\mathcal{H}^{k,p}}$, there exists $T(M)>0$ such that the solution of the previous theorem lives at least until $T(M)$ and the following bound holds:
	$$
	 \| \vv\|_{L^\infty([0,T(M)]; \mathcal H^{k,p})} \leq 2M.
	$$
    We emphasize that the time $T(M)$ only depends on $M$ and  not on any kind of precise property of the set of labels $(I,\mu)$. In addition, the mapping $t \mapsto \Vert \vv(t) \Vert_{\mathcal{H}^{k,p}}$ is continuous on $[0,T(M)]$. 
	\end{Rem}

\begin{Rem}
In full generality, it is not clear that one can ensure that $\lim_{R\to0} T_{\mathrm{inf}}(R) =+\infty$ under the sole Assumption \ref{ass:bound_high_regularity} on the force field. However, when applied on specific models, our general framework and result will often yield that this property holds (see for instance Section \ref{Subsec:VPelectrons} where, for the multiphasic Vlasov--Poisson system for electrons, one can ensure $T_{\mathrm{inf}}(R) \sim R^{-1}$).
\end{Rem}

\begin{Rem}
 We refer to Section \ref{section:Extension} for some variants of this result, replacing the Sobolev based regularity for our functional spaces by weighted Sobolev or Besov regularity.
\end{Rem}

	\subsection{Remarks about time regularity}
	The time regularity stated in Theorem~\ref{thm:existence} does not seem to be very strong (it is only pointwise in $\alpha$ when it comes to derivatives). However, it seems to be a good notion as it has the following implications.
	\begin{Cor}
    \label{cor:timecontinuity}
		In the context of Theorem~\ref{thm:existence}, the solution $\vv$ satisfies the following additional time regularity properties:
		\begin{itemize}
			\item For all $\eps>0$, there exists a measurable set $A \subset I$ such that $\mu(I\setminus A) \leq \eps$ and the family of curves $(v^\alpha)_{\alpha \in A} \in \mathscr{C}^0([0,T]; H^k_x)^A$ is uniformly equicontinuous.
			\item  For all $\psi \in L^1(\mu)$, we have
			\begin{equation*}
			\lim_{\delta \to 0} \sup_{\substack{s,t \in [0,T], \\ |t-s| \leq \delta}} \int \psi(\alpha) \| \DDD v^\alpha(t) - \DDD v^\alpha(s) \|_{H^{k-1}_x} \D \mu(\alpha) = 0.
			\end{equation*}
			\item For all $q \in [1, + \infty)$, we have $\DDD \vv \in \mathscr{C}^0([0,T]; L^q_\alpha H^{k-1}_x)$.
		\end{itemize}
	\end{Cor}
	\begin{proof}
		For the first item, let us call $\Phi: I \to \mathscr{C}^0([0,T]; H^k_x)$ the function that maps $\alpha$ to the curve $(t \mapsto v^\alpha(t))$. It is measurable, and therefore the push-forward of $\mu$ by $\Phi$ is a Borel finite nonnegative measure on the Polish space $\Omega \vcentcolon = \mathscr{C}^0([0,T]; H^k_x)$. Hence, it is Radon, and for all $\eps$, there exists $K$ a compact of $\Omega$ with $\Phi \pf \mu(\Omega \backslash K) \leq \eps$. Calling $A \vcentcolon = \Phi^{-1}(K)$, we get the result since by the Ascoli theorem, all curves of $K$ are uniformly equicontinuous.
		
		The second item is a consequence of the first one, of the uniform boundedness in time of $\DDD \vv \in L^\infty_\alpha H^{k-1}_x$ and of the  integrability of $\psi$. 
				
		The third one is a direct consequence of the dominated convergence theorem, using once again the uniform boundedness in time of $\DDD \vv \in L^\infty_\alpha H^{k-1}_x$ and the continuity of $\DDD v^\alpha$ with values in $H^{k-1}_x$ for $\mu$-almost every $\alpha$. 
	\end{proof}
	
	It is not possible to prove a lot more. Indeed, in the context of Theorem~\ref{thm:existence}, $\DDD v$ is in general  not continuous with values in $L^\infty_\alpha L^2_x$ as the following counter-example with $\FF \equiv 0$ shows.
	\begin{Ex} 
		Let $d = 1$, $p \in [1, + \infty)$ (it is not clear what happens when $p = + \infty$), $I=\R$, and $\mu$ any measure of full support (this assumption is made for simplicity, the crucial point being that this support is unbounded) such that $\int \alpha^p \D \mu(\alpha)$. Let $k \geq k_1$ as given by Theorem~\ref{thm:existence}. Let $v_0 \in H^k_x$ be non-constant and $v$ be the corresponding classical local solution of the Burgers equation
		\begin{equation*}
		\left\{\begin{gathered}
		\partial_t v + v \partial_x v = 0,\\
		v|_{t=0} = v_0.
		\end{gathered} \right.
		\end{equation*}
		Let us call $T>0$ a time such that $v$ is defined up to time $T$. Now, it is clear that the only classical solution of the multiphasic system
		\begin{equation*}
		\left\{\begin{gathered}
		\partial_t v^\alpha + v^\alpha \partial_x v^\alpha = 0,\\
		v^\alpha|_{t=0} = \alpha + v_0,
		\end{gathered} \right.
		\end{equation*}
		is defined for $\mu$-almost all $\alpha$ and all $(t,x) \in [0,T] \times \T$ by 
		\begin{equation*}
		v^\alpha(t,x) = \alpha + v(t,x-\alpha t),
		\end{equation*}
		and hence 
		\begin{equation*}
		\partial_x v^\alpha(t,x) = \partial_x v(t, x-\alpha t).
		\end{equation*}
		On the other hand, it is clear by construction that the family $\vv_0$ belongs to $\mathcal H^{k},p$, and indeed, it is easy to check that $\vv$ satisfies the time regularity stated in Theorem~\ref{thm:existence}. On the other hand, let us show that $\vv$ is not continuous with values in $\mathcal H^{k,p}$ by showing that
		\begin{equation*}
		\liminf_{t \to 0} \esssup_{\alpha \in (I, \mu)} \| \partial_x v^\alpha (t) - \partial_x v^\alpha(0)\|_{L^2_x} >0.
		\end{equation*}
		To do so, let us consider $c \in \R$ be such that $ \partial_x v_0 \neq \tau_c \partial_x v_0$, where $\tau_c$ is the translation operator defined for all functions $f$ by  $\tau_c f(\cdot) = f(\cdot-c)$. Note that $\tau_c$ is an isometry on $L^2_x$. Now, for a given $t \in (0,T]$, we have
		\begin{align*}
		\| \partial_x v^\alpha (t) - \partial_x v^\alpha(0)\|_{L^2_x} &= \| \tau_{\alpha t} \partial_x v(t) - \partial_x v(0) \|_{L^2_x} \\
		&\geq \| \tau_{\alpha t} \partial_x v(0) - \partial_x v(0) \|_{L^2_x} - \| \tau_{\alpha t} ( \partial_x v(t) - \partial_x v(0) ) \|_{L^2_x}\\
		&= \| \tau_{\alpha t} \partial_x v(0) - \partial_x v(0) \|_{L^2_x} - \|  \partial_x v(t) - \partial_x v(0)  \|_{L^2_x}.
		\end{align*}
		The second term of the last line does not depend on $\alpha$ and tends to $0$ as $t \to 0$. The first one is continuous in $\alpha$, and hence because $\mu$ is of full support, choosing $\alpha \vcentcolon = c/t$, 
		\begin{align*}
		\liminf_{t \to 0} \esssup_{\alpha \in (I, \mu)} \| \partial_x v^\alpha (t) - \partial_x v^\alpha(0)\|_{L^2_x} &= \liminf_{t \to 0} \sup_{\alpha \in \R} \|  L_{\alpha t} \partial_x v(0) - \partial_x v(0)\|_{L^2_x} \\
		&\geq \| L_c \partial_x v(0) - \partial_x v_0 \|_{L^2_x} \neq 0,
		\end{align*}
		therefore our claim follows.
	\end{Ex}

	\subsection{Proof of Theorem~\ref{thm:existence}}\label{subsec-proofTHM}
	
	Let $k\geq k_0$ and $p$ as in Assumption~\ref{ass:bound_high_regularity} and consider $\vv_0$ and $X$ as in the statement of the theorem. For given $T>0$ and $\ww \in L^1([0,T]; \mathcal H^{k,p})$, we call $\widehat{\ww}  \vcentcolon= \PHI_T[\ww,\vv_0,X]$ the unique classical solution up to time $T$ to
	\begin{equation}\label{eq:defPHI}
	\left\{
	\begin{gathered}
	\partial_t \widehat{w}^\alpha +( w^\alpha \cdot \nabla) \widehat{w}^\alpha = F^\alpha[\ww,X],\\
	\widehat{w}^\alpha|_{t=0} = v^\alpha_0.
	\end{gathered}
	\right.
	\end{equation}
	Under the assumption that $\ww \in L^1([0,T]; \mathcal H^{k,p})$ for $k$ large enough for $\Vert \cdot \Vert_{\mathcal{H}^{k,p}}$ to control the Lipschitz bound of $w^\alpha$ (uniformly in $\alpha$ by definition), $\PHI_T[\ww,\vv_0,X]$ is well defined using the method of characteristics. The main steps of the proof consist first in deriving a bound for 
	\begin{equation*}
\sup_{t \in [0,T]} \Vert \PHI_T[\ww, \vv_0,X](t)\Vert_{\mathcal{H}^{k,p}}
	\end{equation*}
	only implying the $ L^1([0,T]; \mathcal H^{k,p})$ bound on $\ww$, and then to find a modulus of continuity for $\PHI_T$ in $L^p_\alpha L^2_x$ on  bounded sets of $L^\infty([0,T]; \mathcal H^{k,p})$. As we will see, these two steps will be enough to perform a classical fixed-point proof in a suitable space, as a solution to~\eqref{eq:general_system} is nothing but a fixed-point of $\PHI_T$.

\paragraph{Preliminaries}	
We collect here with some standard properties of Sobolev spaces (set on $\T^d$ or $\R^d$) that will prove crucial in the rest of this article. We refer for instance to \cite[Chapter 2]{BCD} and \cite{MB} for the proofs.
\begin{Prop}\label{prop-Sobolev}
The Sobolev spaces satisfy the following properties.
\begin{enumerate}
\item  {\bf Sobolev embedding.} For $m\in \N$, $k >d/2$, we have the continuous embedding $H^{k+m} \hookrightarrow \mathscr{C}^{m}$.
\item  {\bf Product laws.} If $k \leq \min(k_1,k_2)$, $k+d/2<k_1 + k_2$, we have
\begin{align*}
\Vert uv \Vert_{H^k_x} \leq \Vert u \Vert_{H^{k_1}_x} \Vert v \Vert_{H^{k_2}}.
\end{align*}
\item {\bf Tame estimates.} If $k >0$, then for all $w_1,w_2 \in H^k \cap L^{\infty}$, we have
\begin{align*}
\LRVert{ uv}_{H^k_x} \lesssim \LRVert{u}_{L^{\upinfty}_x}\LRVert{v}_{H^k_x}+\LRVert{u}_{H^k_x}\LRVert{v}_{L^{\upinfty}_x} .
\end{align*}
\item {\bf Commutator laws.} For any integer $k>0$ and for any functions $b$ and $u$, we have 
\begin{align}\label{commut-law}
    \sum_{\vert \alpha \vert =k} \Vert [\partial_x^\alpha, b]u \Vert_{L^2_x} \lesssim \Vert \DDD b \Vert_{L^{\upinfty}_x} \sum_{\vert \beta \vert =k-1} \Vert \partial^\beta u \Vert_{L^2_x}+\Vert u \Vert_{L^{\upinfty}_x}\sum_{\vert \beta \vert =k} \Vert  \partial^\beta b \Vert_{L^2_x} .
\end{align}
In particular, the following estimates on the convection operator hold: 
for any integer $k>0$ and any vector fields $a,b$, we have
		\begin{align}\label{eq:convective}
		\begin{split}
		\left\langle   a, (b\cdot \nabla )a  \right\rangle_{H^k_x} &\lesssim  \|  a \|_{H^k_x} \Big\{  \| \DDD b \|_{L^\infty_x} \|  \DDD a \|_{H^{k-1}_x} + \|  \DDD b \|_{H^{k-1}_x}\| \DDD a \|_{L^{\upinfty}_x}  \Big\}, \\ 
		\langle  \DDD a, \DDD \big( (b\cdot \nabla )a \big) \rangle_{H^{k-1}_x} &\lesssim  \| \DDD a \|_{H^{k-1}_x} \\
        &  \qquad\times \Big\{  \| \DDD b \|_{L^{\upinfty}_x} \| \DDD a \|_{H^{k-1}_x} + \| \DDD b \|_{H^{k-1}_x}\| \DDD a \|_{L^{\upinfty}_x}  \Big\}.
		\end{split}
		\end{align}
\item {\bf Composition rules.} Let $s \in \N\setminus\{0\}$ with and $\Gamma \in \mathscr{C}^{\infty}(\R)$ satisfying $\Gamma(0)=0$. Then for all $\mathrm{g} \in H^s \cap L^\infty$,  we have 
\begin{align*}
    \left\Vert \Gamma(\mathrm{g}) \right\Vert_{ H^s} &\leq \mathrm{C}\left(\Gamma'\right) (1+ \Vert \mathrm{g} \Vert_{L^{\upinfty}_x})^{\lceil s \rceil}  \Vert \mathrm{g} \Vert_{H^s}.
\end{align*}
Furthermore, if $\Gamma'(0)=0$ then for all $\mathrm{g} \in H^s \cap L^\infty$, we have 
\begin{align*}
    \left\Vert \Gamma(\mathrm{g}) \right\Vert_{H^s} &\leq \mathrm{C}\left(\Gamma''\right) (1+ \Vert \mathrm{g} \Vert_{L^{\upinfty}_x})^{\lceil s \rceil}  \Vert \mathrm{g} \Vert_{L^{\upinfty}_x} \Vert \mathrm{g} \Vert_{H^s}.
\end{align*}
Above, $\mathrm{C}\left(\Gamma'\right)$ and $\mathrm{C}\left(\Gamma''\right)$ are two continuous nondecreasing functions.

\end{enumerate}
\end{Prop}
We now record a standard lemma about linear transport equation with a source.
	\begin{Lem}\label{LM:existence-vitesseSOBOLEV}
Let $k > 0$. Let $T>0$, $w_0 \in H^k_x$ and $b,S \in L^\infty([0,T]; H^k_x)$ be some vector fields. Then the unique classical solution $w$ to
\begin{align*}
    \partial_t w + (b \cdot \nabla_x) w = S, \qquad w|_{t=0}= w_0,
\end{align*}
in $\mathscr{C}([0,T]; H^k_x)$ satisfies the following: there exists $C=C_k>0$ such that for all $t \in [0,T]$, we have
   
   \begin{align}
	\label{boundL2:existence-vitesseSOBOLEV} \| w(t) \|_{L^2_x} &\leq \| w(0) \|_{L^2_x} + C \int_0^t \big( \| \DDD b \|_{L^{\upinfty}_x} \| w \|_{L^2_x}  + \| S \|_{L^2_x} \big) \, \mathrm{d} s, \\
	\label{boundHk:existence-vitesseSOBOLEV}  \|  w(t) \|_{H^{k}_x} & \leq  \|  w(0) \|_{H^{k}_x} \\
    &\qquad +  C \int_0^t \big(\|  \DDD b \|_{L^{\upinfty}_x} \|  w \|_{H^{k}_x}+ \|   \DDD b \|_{H^{k-1}_x} \|  \DDD w \|_{L^{\upinfty}_x}+  \|  S \|_{H^{k}_x} \big) \, \mathrm{d} s,
\end{align}
    and if $k>1$,
    \begin{equation}
            \label{boundD-Hk-1:existence-vitesseSOBOLEV} 
    \begin{aligned}
	&\| \DDD w(t) \|_{H^{k-1}_x} \\
    &\leq  \| \DDD w(0) \|_{H^{k-1}_x} \\
    &\qquad+   C \int_0^t   \big(   \| \DDD b \|_{L^{\upinfty}_x} \| \DDD w \|_{H^{k-1}_x} + \| \DDD b \|_{H^{k-1}_x}\| \DDD w \|_{L^{\upinfty}_x} + \| \DDD S \|_{H^{k-1}_x}  \big) \, \mathrm{d} s.
	\end{aligned}
        \end{equation}
\end{Lem}
\begin{proof}
 The proof relies on direct \textit{a priori} estimates. First, an energy estimate in $L^2_x$ in the equation on $w$ yield
	\begin{equation*}
	\frac{1}{2}\frac{\D}{\D t} \| w \|_{L^2_x}^2 \leq \| \DDD b \|_{L^{\upinfty}_x} \| w \|_{L^2_x}^2 + \| S \|_{L^2_x} \| w \|_{L^2_x},
	\end{equation*}
	and hence simplifying by $\| w \|_{L^2_x}$, we infer \eqref{boundL2:existence-vitesseSOBOLEV}. On the other hand, computing the $H^{k}_x$ scalar product of this equation with $w$ and using the estimate \ref{eq:convective} on the convective term, we find
	\begin{equation*}
	\frac{1}{2}\frac{\D}{\D t} \|  w \|^2_{H^{k}_x} \leq C \|  w \|_{H^{k}_x} \Big\{  \| \DDD b \|_{L^{\upinfty}_x} \|  w \|_{H^{k}_x} +\|   \DDD b \|_{H^{k-1}_x} \| \DDD w \|_{L^{\upinfty}_x}  \Big\} + \|  S \|_{H^{k}_x} \|   w\|_{H^{k}_x},
	\end{equation*}
	for some $C>0$. Simplifying by $\| w^\alpha \|_{H^k_x}$, we get \eqref{boundHk:existence-vitesseSOBOLEV}
    . The last inequality \eqref{boundD-Hk-1:existence-vitesseSOBOLEV} is deduced in the same way by relying on \ref{eq:convective}, first applying the operator $\mathrm{D}$ to the equation and computing the $H^{k-1}_x$ scalar product of the equation with $\mathrm{D} w$.
   Existence and uniqueness is then standard, using the method of characteristics.
\end{proof}

	\medskip
	
Let us now proceed with the actual proof of Theorem \ref{thm:existence}.
\begin{proof}[Proof of Theorem \ref{thm:existence}]
	
	In the what follows, the constant $C$ refers to a large number possibly growing from line to line, but only depending on the dimension and $k$.
	
	\paragraph{Bound in high regularity.} Let us assume that $k>1+d/2$ so that $H^{k-1}_x \hookrightarrow L^\infty_x$ and  $H^{k-1}_x$ is an algebra. Let $T>0
$ and $\ww \in L^1([0,T]; \mathcal H^{k,p})$. Let us consider $\widehat{\ww}  = \PHI_T[\ww, \vv_0,X]$ as above, that is the solution to \eqref{eq:defPHI}. Appealing to Lemma \ref{LM:existence-vitesseSOBOLEV}, 
we have
\begin{align*}
  \| \widehat{w}^\alpha(t) \|_{L^2_x} &\leq  \| \widehat{w}^\alpha \|_{L^2_x} +  C \int_0^t\left( \| \DDD v_0^\alpha \|_{L^{\upinfty}_x} \| \widehat{w}^\alpha \|_{L^2_x}  + \| F^\alpha[\ww,X] \|_{L^2_x}\right) \, \mathrm{d} s ,
\end{align*}
and 
\begin{align*}
   \| \DDD &\widehat{w}^\alpha \|_{H^{k-1}_x} \\
   &\leq
    \| \DDD v_0^\alpha \|_{H^{k-1}_x} +C \int_0^t \Big(   \| \DDD w^\alpha \|_{L^{\upinfty}_x} \| \DDD \widehat{w}^\alpha \|_{H^{k-1}_x} \\
    & \qquad \qquad \qquad \qquad  \qquad \qquad  + \| \DDD w^\alpha \|_{H^{k-1}_x}\| \DDD \widehat{w}^\alpha \|_{L^{\upinfty}_x} + \| \DDD F^\alpha[\ww,X] \|_{H^{k-1}_x}\Big)\, \mathrm{d} s \\
    &\leq  \| \DDD v_0^\alpha \|_{H^{k-1}_x}  + C
   \int_0^t \left(   \| \DDD w^\alpha \|_{H^{k-1}_x}\| \| \DDD \widehat{w}^\alpha \|_{H^{k-1}_x}  + \| \DDD F^\alpha[\ww,X] \|_{H^{k-1}_x}\right)\, \mathrm{d} s ,
\end{align*}
the former estimates being pointwise in $\alpha$. By Gronwall lemma, and having respectively taken the $L^p_\alpha$ and the $L^\infty_\alpha$ norm for each family of estimates, we get
\begin{align}\label{eq:estimate_bound_high_regularity_0_derivative}
\begin{split}
    &\|\widehat{\ww} (t)\|_{L^p_\alpha L^2_x} \\ 
    &\leq \exp\left[ C \int_0^t \|\mathrm{D} \ww(\tau) \|_{L^{\upinfty}_\alpha L^{\upinfty}_x} \D \tau \right] \| \hat \vv_0\|_{L^p_\alpha L^2_x}  \\
    & \quad + C \int_0^t  \exp\left[ C \int_{\tau}^t \|\mathrm{D} \ww(s) \|_{L^{\upinfty}_\alpha L^{\upinfty}_x} \D s \right] \| \FF[\ww,X](\tau)\|_{L^p_\alpha L^2_x} \D \tau ,
\end{split}
   \end{align}
   and 
\begin{align}\label{eq:estimate_bound_high_regularity_k_derivatives}
\begin{split}
   &\|\mathrm{D} \widehat{\ww}(t)\|_{L^{\upinfty}_\alpha H^{k-1}_x} \\
   &\leq\exp\left[ C \int_0^t \|\mathrm{D} \ww(\tau) \|_{L^{\upinfty}_\alpha H^{k-1}_x} \D \tau \right] \| \mathrm{D}  \vv_0\|_{L^{\upinfty}_\alpha H^{k-1}_x}  \\
   &  \quad +C \int_0^t  \exp\left[ C \int_{\tau}^t \|\mathrm{D} \ww(s) \|_{L^{\upinfty}_\alpha H^{k-1}_x} \D s \right] \| \mathrm{D} \FF[\ww,X](\tau)\|_{L^{\upinfty}_\alpha H^{k-1}_x} \D \tau.
   \end{split}
\end{align}

Recalling the definition of the norm $\mathcal{H}^{k,p}$ in \eqref{def:norm}, we can sum the two previous estimates and use Sobolev embedding to get for all $t\geq 0$ (the right-hand side  being $+\infty$ whenever $\FF$ is not defined up to time $t$)
	\begin{align*}
	\Vert \widehat{\ww}(t) \Vert_{\mathcal{H}^{k,p}} &\leq \Vert \vv_0 \Vert_{\mathcal{H}^{k,p}} \exp\left( C \int_0^t \Vert \ww(s) \Vert_{\mathcal{H}^{k,p}} \D s \right) \\
    &\quad+ \int_0^t \exp\left( C \int_s^t \Vert \ww(\tau) \Vert_{\mathcal{H}^{k,p}} \D \tau \right)\Vert\FF[\ww,X](s)\Vert_{\mathcal{H}^{k,p}} \D s\\
	&\leq \exp\left( C \int_0^t \Vert \ww(s) \Vert_{\mathcal{H}^{k,p}} \D s \right) \left( \Vert \vv_0 \Vert_{\mathcal{H}^{k,p}} + \int_0^t \Vert \FF[\ww,X](s)\Vert_{\mathcal{H}^{k,p}} \D s \right).
	\end{align*}
	Using Assumption~\ref{ass:bound_high_regularity}, we find that for all $t \leq T$,
	\begin{equation}
	\label{eq:bound_highregLWP}
    \begin{aligned}
	\big\Vert\PHI_T[\ww, \vv_0,X](t)\big\Vert_{\mathcal{H}^{k,p}} &\leq \exp\left( C \int_0^t \Vert \ww(s) \Vert_{\mathcal{H}^{k,p}}  \D s \right)\\
    &\quad \times\left( \Vert \vv_0 \Vert_{\mathcal{H}^{k,p}} + \F\left(t, \int_0^t \Vert \ww(s) \Vert_{\mathcal{H}^{k,p}} \D s, \mathsf N_k(X) \right) \right).
    \end{aligned}
	\end{equation}

	\paragraph{Stability at low regularity.} We consider the same exponents $k$ as before. Let us choose $T,R>0$, initial conditions $\vv_{1,0},\vv_{2,0} \in \mathcal H^{k,p}$, associated families of vector fields $\ww_1,\ww_2 \in L^\infty([0,T]; \mathcal{H}^{k,p})$ and parameters $X_1,X_2 \in \mathcal X$ such that
	\begin{equation}
	\label{eq:bounded_w}
\max \Big(\sup_{t \in [0,T]}\Vert \ww_1(t) \Vert_{\mathcal{H}^{k_0,p}},\sup_{t \in [0,T]} \Vert \ww_2(t) \Vert_{\mathcal{H}^{k_0,p}}, \mathsf N_{k_0}(X_1), \mathsf N_{k_0}(X_2)  \Big) \leq R,
	\end{equation}
	and let us call $\widehat{\ww}_i\vcentcolon =\PHI_T[\ww_i, \vv_{i,0},X_i]$ for $i=1,2$. The differences $\widehat{\ww}\vcentcolon =\widehat{\ww}_2-\widehat{\ww}_1$ and $\ww\vcentcolon = \ww_2- \ww_1$  satisfy
	\begin{equation*}
	\left\{ 
	\begin{gathered}
	\partial_t \widehat{w}^{\alpha} + (w^{\alpha}_2 \cdot \nabla ) \widehat{w}^{\alpha} =-( w^{\alpha} \cdot \nabla) \widehat{w}^{\alpha}_1 +  F^{\alpha}[\ww_2,X_2] - F^{\alpha}[\ww_1,X_1]  ,\\
	\widehat{w}^{\alpha}|_{t=0} = v^{\alpha}_{2,0} - v^{\alpha}_{1,0}.
	\end{gathered}
	\right.
	\end{equation*}
By Lemma \ref{LM:existence-vitesseSOBOLEV}, we have
\begin{equation*}
\begin{aligned}
	\frac{\D}{\D t} \| \widehat{w}^{\alpha} \|_{L^2_x} &\leq \| \DDD w^{\alpha}_2 \|_{L^{\upinfty}_x}  \| \DDD \widehat{w}^{\alpha} \|_{L^2_x}  
	+ \Vert \DDD \widehat{w}^{\alpha}_1 \Vert_{L^{\upinfty}_x}
	\| w^{\alpha} \|_{L^2_x}\\
    & \quad+ \| F^\alpha[\ww_2,X_2] - F^\alpha[\ww_1, X_1] \|_{L^2_x}.
    \end{aligned}
	\end{equation*}
	Therefore, taking the $L^p_\alpha$ norm gives for all $t \leq T$:
	\begin{multline*}
	\frac{\D}{\D t} \| \widehat{\ww} \|_{L^p_\alpha L^2_x} \leq R \| \widehat{\ww} \|_{L^p_\alpha L^2_x} + \Vert \widehat{\ww}_1 \Vert_{\mathcal{H}^{k,p}}
	\| \ww \|_{L^p_\alpha L^2_x}\\ 
    + \| \FF^\alpha[\ww_2,X_2] - \FF^\alpha[\ww_1,X_1] \|_{L^p_\alpha L^2_x}.
	\end{multline*}
	by definition of $\Vert \cdot \Vert_{\mathcal{H}^{k,p}}$, the bound~\eqref{eq:bounded_w} and Sobolev embedding. In view of \eqref{eq:bound_highregLWP} and the bound~\eqref{eq:bounded_w}, we also have for all $t \leq T$
	\begin{equation*}
	\Vert\widehat{\ww}_1 \Vert_{\mathcal{H}^{k,p}}=\big\Vert\PHI_T[\ww_1, \vv_{1,0},X_1](t)\big\Vert_{\mathcal{H}^{k,p}} \leq e^{KTR}\left( R + \F\left(T, TR, R\right) \right).
	\end{equation*} 
	We get for all $t \leq T$
	\begin{equation*}
    \begin{aligned}
	\frac{\D}{\D t} \| \widehat{\ww} \|_{L^p_\alpha L^2_x} &\leq R \| \widehat{\ww} \|_{L^p_\alpha L^2_x} + e^{CRt}\left( R + \F\left(t, tR, \mathsf R \right) \right)
	\| \ww \|_{L^p_\alpha L^2_x}\\
    &\quad + \| \FF^\alpha[\ww_2,X_2] - \FF^\alpha[\ww_1,X_1] \|_{L^p_\alpha L^2_x}.
    \end{aligned}
	\end{equation*}
	Using the Gronwall lemma together with Assumption~\ref{ass:stability}, we infer for all $t \leq T$,
	\begin{equation}
	\label{eq:bound_stablowLWP}
	\begin{aligned}
	&\| \PHI_T[\ww_2, \vv_{2,0},X_2](t) - \PHI_T[\ww_1, \vv_{1,0},X_1](t) \|_{L^p_\alpha L^2_x} \leq  \exp(Rt)  \\
    &\quad \times\left\{ \| \vv_{2,0} - \vv_{1,0} \|_{L^p_\alpha L^2_x} + \mathcal{G}_{R} \left(t, \int_0^t \| \ww_2(s) - \ww_1(s) \|_{L^p_\alpha L^2_x} \D s, \mathsf d(X_1,X_2) \right)\right\},
	\end{aligned}
	\end{equation}
	where $\mathcal{G}_{R}$ is a stability estimate function.

	\paragraph{fixed-point procedure.} For given $T,R>0$, we define
	\begin{equation*}
	B_{T,R}\vcentcolon = \Big\lbrace \ww \in \mathscr{C}^0([0,T]; L^p_\alpha L^2_x) \cap L^{\infty}([0,T];\mathcal{H}^{k,p}) \mid \Vert \ww \Vert_{L^{\infty}([0,T];\mathcal{H}^{k,p})} \leq R \Big\rbrace.
	\end{equation*}
	Endowed with the metric $\Vert \cdot \Vert_{L^\infty([0,T]; L^p_\alpha L^2_x)}$, $B_{T,R}$ is a nonnempty and complete metric space. We will show that for all initial condition $\vv_0 \in \mathcal H^{k,p}$ and parameter $X \in \mathcal X$ such that $(0,0, \mathsf N(X)) \in \Omega_\F$ (see Assumption~\ref{ass:bound_high_regularity}), setting
	\begin{equation}
	\label{eq:def_R}
	R \vcentcolon = \max \Big( 2 \Vert \vv_0 \Vert_{\mathcal{H}^{k,p}}, \mathsf N(X) \Big),
	\end{equation}
	there exists $T>0$ such that $\PHI_T[B_{T,R}, \vv_0,X] \subset B_{T,R}$ and $\PHI_T[\cdot, \vv_0,X]$ is a contraction on $B_{T,R}$. Solutions of~\eqref{eq:general_system} of initial condition $\vv_0$ and parameter $X$ being precisely the fixed points of $\PHI_T[\cdot, \vv_0,X]$, the existence of a unique solution will follow.
	
	\bigskip 
	
	\noindent \underline{For $T$ small enough, $B_{T,R}$ is stable by $\PHI_T[\cdot, \vv_0,X]$}. Let $T>0$ and $R$ defined as in~\eqref{eq:def_R}. By~\eqref{eq:bound_highregLWP}, for all $\ww \in B_{T,R}$, as $\F$ is an estimate function in the sense of Definition~\ref{def:estimate_func}, we have:
	\begin{align*}
	&\qquad\|\PHI_T [\ww, \vv_0,X]\|_{L^\infty([0,T]; \mathcal H^{k,p})}  \\
    &\leq \exp\left( K \int_0^T \Vert\ww(s) \Vert_{\mathcal{H}^{k,p}} \D s \right)\left\{ \Vert \vv_0 \Vert_{\mathcal{H}^{k,p}} + \F\left(T, \int_0^T \Vert\ww(s) \Vert_{\mathcal{H}^{k,p}} \D s, \mathsf N(X) \right) \right\}\\
	&\leq \exp(CTR) \left\{ \frac{R}{2} + \F\left(T, TR, R \right) \right\}.
	\end{align*}
	This last quantity converges to $R/2$ as $T \to 0$, therefore there exists $T_1>0$ such that for all $T \leq T_1$ and $\ww \in B_{T,R}$,
	\begin{equation*}
		\|\PHI_T[\ww, \vv_0,X] \|_{L^\infty([0,T]; \mathcal H^{k,p})} \leq R.
	\end{equation*}
	Once this is ensured, it means that for such $\ww$, $\boldsymbol F[\ww, X] \in L^1([0,T]; \mathcal H^{k,p})$. Therefore, continuity of~$\PHI_T[\ww,\vv_0,X]$ in $L^p_\alpha L^2_x$ is straightforward to check, using for all $t \leq T$ the following estimate directly derived from~\eqref{eq:defPHI}:
	\begin{multline*}
	\left\| \partial_t \PHI_T[\ww,\vv_0, X](t) \right\|_{L^p_\alpha L^2_x} \\
    \leq \|\DDD \PHI_T[\ww, \vv_0,X](t)\|_{L^{\upinfty}_\alpha L^{\upinfty}_x} \| \ww(t) \|_{L^p_\alpha L^2_x} + \| \boldsymbol F[\ww,X](t) \|_{L^p_\alpha L^2_x}.
	\end{multline*}
	The claim follows.
	
	\bigskip
	
	\noindent \underline{Contraction estimate}. Let $R$ be defined by~\eqref{eq:def_R}, $T \leq T_1$ and $\ww_1,\ww_2 \in B_{T,R}$. By~\eqref{eq:bound_stablowLWP}, as $\mathcal G_R$ is an estimate function, we have
	\begin{align}
	\notag \| \PHI_T[\ww_2,\vv_0,X]&-\PHI_T[\ww_1,\vv_0,X] \|_{L^{\infty}([0,T];L^p_\alpha L^2_x)}\\ 
	\notag &\leq \exp(RT) \mathcal G_R\left( T, \int_0^T \| \ww_2(s) - \ww_1(s) \|_{L^p_\alpha L^2_x} \D s,0 \right)\\
	\label{eq:estim_contraction} &\leq \exp(RT)  \mathcal G_R\Big( T, T \| \ww_2 - \ww_1 \|_{L^\infty([0,T];L^p_\alpha L^2_x)} \D s,0 \Big).
	\end{align} 
	Notice that as $\ww_1,\ww_2 \in B_{T,R}$, we have
	\begin{equation}
	\| \ww_2 - \ww_1 \|_{L^\infty([0,T];L^p_\alpha L^2_x)} \leq 2R.
	\end{equation}
	Now, let $T_2 \leq T_1$ be such that:
	\begin{itemize}
		\item The rectangle $[0,T_2] \times [0,2RT_2] \times \{ 0 \}$ is a subset of $\Omega_{\G_R}$. 
		\item For all $(T,M) \in [0,T_2] \times [0,2RT_2]$, 
		\begin{equation*}
		\partial_2 \mathcal G_R(T,M,0) \leq \frac{1}{2T_2\exp(RT_2)}.
		\end{equation*}
	\end{itemize} 
For $T \leq T_2$ and $\ww_1,\ww_2 \in B_{T,R}$, we have by the mean value theorem applied to~\eqref{eq:estim_contraction}
	\begin{align*}
	\| \PHI_T[\ww_2,\vv_0,X] &-\PHI_T[\ww_1,\vv_0,X] \|_{L^{\infty}([0,T];L^p_\alpha L^2_x)} \\
    &\leq \frac{\exp(RT)}{2T_2\exp(RT_2)} \times T \| \ww_2 - \ww_1 \|_{L^\infty([0,T];L^p_\alpha L^2_x)}\\
	&\leq \frac{1}{2} \| \ww_2 - \ww_1 \|_{L^\infty([0,T];L^p_\alpha L^2_x)}.
	\end{align*}
	We conclude that $\PHI_T[\cdot, \vv_0, X]$ is a contraction of $B_{T,R}$, and the existence of a fixed point follows.
	
	\paragraph{Time regularity.} The continuity with values in $L^p_\alpha L^2_x$ is directly given by the fixed-point argument (it follows from the completeness of $B_{T,R}$). It remains to prove that for $\mu$-almost all $\alpha$, $w^\alpha \in \mathscr{C}^0([0,T]; H^k_x)$. First, observe that for $\mu$-almost all $\alpha$ the force field $F^\alpha[\ww,X] \in L^1([0,T]; H^k_x)$. Indeed, as $\boldsymbol F[\ww,X] \in L^1([0,T]; \mathcal H^{k,p})$, we have in particular (the measure $\mu$ being of finite mass):
\begin{align*}
	\int \left( \int_0^T \| F^\alpha[\ww,X](t) \|_{L^2_x}\D t \right) \D \mu(\alpha) &\lesssim \int_0^T \left(  \int \| F^\alpha[\ww,X](t) \|^p_{L^2_x} \D \mu(\alpha) \right)^{1/p} \D t \\
    &< + \infty, \\
	\esssup_{\alpha} \int_0^T  \| \DDD F^\alpha[\ww,X](t) \|_{H^{k-1}_x} \D t &\leq  \int_0^T \esssup_{\alpha} \| \DDD F^\alpha[\ww,X](t) \|_{H^{k-1}_x} \D t \\
    &< + \infty.
	\end{align*}
	The claim follows. Therefore, we have that for $\mu$-almost all $\alpha$, 	\begin{itemize}
		\item the initial condition $v_0^\alpha$ is in $H^k_x$;
		\item the velocity field $w^\alpha$ satisfies
		\begin{equation}
		\label{eq:uniform_bound_valpha}
		\sup_{t \in [0,T]} \| \DDD w^\alpha(t) \|_{H^{k-1}_x} < + \infty \qquad \mbox{and hence} \qquad \sup_{t \in [0,T]} \| \DDD w^\alpha(t)\|_{L^{\upinfty}_x} < + \infty;
		\end{equation}
		\item the force field $F^\alpha[\ww,X]$ is in $L^1([0,T]; H^k_x)$.
	\end{itemize}
This is everything we need to conclude with standard arguments, that we only sketch here: 
\begin{enumerate}
\item $\partial_t w^\alpha \in L^1([0,T]; L^2_x)$ by direct estimates, so $w^\alpha \in \mathscr{C}^0([0,T]; L^2_x)$; 
\item $w^\alpha \in L^\infty([0,T]; H^k_x)$ by the first point and~\eqref{eq:uniform_bound_valpha}, so $w^\alpha$ is also weakly continuous with values in $H^k_x$; 
\item the map $t \in [0,T] \mapsto \|w^\alpha(t)\|_{H^k_x}$, which is a priori lower semicontinuous because $\|\cdot \|_{H^k_x}$ is itself lower semicontinuous with respect to the $L^2_x$ topology on bounded sets of $H^k_x$, is actually continuous thanks to energy estimates in $H^k_x$ for~\eqref{eq:defPHI}. Note that by considering the (differential form of) the estimate \eqref{eq:bound_highregLWP} and integrating between two times $t_1$ and $t_2$, we also obtain the fact that $t \mapsto \Vert \vv(t) \Vert_{\mathcal{H}^{k,p}}$ is continuous.
\end{enumerate}

	\paragraph{Stability.} The stability result is a direct consequence of~\eqref{eq:bound_stablowLWP} applied to two solutions and of Gronwall lemma. 
	
	 \paragraph{Blow-up criterion.} 
	Since we now have a solution $\ww$ of 
	\begin{equation}
	\left\{
	\begin{gathered}
	\partial_t w^\alpha + ( w^\alpha \cdot \nabla ) w^\alpha = F^\alpha[\ww,X].\\
	w^\alpha |_{t = 0} = v^\alpha_0,
	\end{gathered}
	\right.
	\end{equation}
	we can argue as in the beginning of the proof regarding bounds in high regularity: on $[0,T]$ by appealing to the tame estimates of Lemma \ref{LM:existence-vitesseSOBOLEV}:
	\begin{align*}
    \| w^\alpha(t) \|_{L^2_x} &\lesssim \| w^\alpha |_{t = 0} \|_{L^2_x}+ \int_0^t \big(\| \DDD w^\alpha \|_{L^{\upinfty}_x} \| w^\alpha \|_{L^2_x}  + \| F^\alpha[\ww,X] \|_{L^2_x} \big) \, \mathrm{d}s, \\
	    \| \DDD w^\alpha(t) \|_{H^{k-1}_x} &\lesssim \| \DDD w^\alpha |_{t = 0} \|_{H^{k-1}_x}\\
        &\quad+ \int_0^t \big(\| \DDD w^\alpha \|_{L^{\upinfty}_x} \| \DDD w^\alpha \|_{H^{k-1}_x} + \| \DDD F^\alpha[\ww,X] \|_{H^{k-1}_x} \big) \, \mathrm{d}s.
	\end{align*}
The desired blow-up criterion in $\mathcal{H}^{k,p}$ follows.
\end{proof}

    To conclude, let us discuss a slight generalization relevant for the whole space case

\begin{Rem}
\label{rem:Rd-alpha}
Let $\mathbf{\lambda}: (I,\mu) \rightarrow \R^d$ be a measurable map and replace the general multiphasic equation~\eqref{eq:general_system} by:
	\begin{equation}
	\label{eq:general_system_Rd}
	\left\{
	\begin{gathered}
	\partial_t w^\alpha + ( (w^\alpha+ \lambda^\alpha) \cdot \nabla ) w^\alpha = G^\alpha[\ww,X].\\
	w^\alpha |_{t = 0} = w^\alpha_0.
	\end{gathered}
	\right.
	\end{equation}
We readily check that only the  \emph{space gradient} of the vector fields $w^\alpha+ \lambda^\alpha$ matters in the estimates of the previous subsection and $\lambda^\alpha$ does not depend on space: hence, all results pertaining to~\eqref{eq:general_system} (in particular the main Theorem~\eqref{thm:existence}), apply {\it mutatis mutandis} to \eqref{eq:general_system_Rd}.

This slight generalization is relevant for the case when the spatial domain is $\R^d$, in view of applying the theory to treat Vlasov equations with an initial condition that is an $L^1$ function in velocity, see the proof of Proposition~\ref{prop:disintegration}, in which case we recall that it is relevant to take 
$I=\R^d$ and $v_0^\alpha(x)= \alpha$. We recall that this initial condition is not $L^2(\R^d)$ integrable in space and thus Theorem~\ref{thm:existence} does not apply directly. However, this is savable by considering the unknown $w^\alpha\vcentcolon = v^\alpha-\alpha$; it is then relevant to consider~\eqref{eq:general_system_Rd} with $$\lambda^\alpha=\alpha, \quad G^\alpha[\ww,X]= F^\alpha[\ww+ \boldsymbol{\alpha},X], \quad w^\alpha_0=0.$$

\end{Rem}

\section{General estimates on the multiphasic continuity equation}\label{sec-continuity-eq}

In this subsection, we derive standard estimates that we will use for the densities solving the continuity equation in the multiphasic systems. This will prove useful in the rest of this article. In what follows, $(I, \mu)$ is a fixed probability space. We also recall the notation $L^q_{\alpha, \Theta}(\mathcal{B})$ for weighted Lebesgues spaces with values in $\mathcal{B}$ -- see \eqref{eq:def_weighted_Lp}.

\paragraph{Bound in high regularity.}
We start by the following elementary result about high Sobolev estimates for general continuity equations, that we then apply to the multiphasic case.

\begin{Lem}\label{LM:SobEstimRHO}
Let $k \in \N$ such that $k>d/2+1$ and $p,q \in [1,+\infty]$. There exists $K>0$ such that the following holds for any $T>0$ an any set of labels $({I}, \mu)$.  Let $\vv$ be a family velocity fields in $L^\infty([0,T]; \mathcal{H}^{k,p}$) 
such that $\DDD v^\alpha \in \mathscr{C}^0([0,T]; H^{k-1}_x)$ for $\mu$-almost all $\alpha \in {I}$.
Then for all measurable function $\Theta: I \rightarrow \R^+$ and for all families $\rhorho_0$ in $L^q_{\alpha, \Theta} H^{k-1}_x$, there exists a unique family of densities 
$$\rhorho \in L^\infty([0,T]; L^q_{\alpha, \Theta} H^{k-1}_x) \cap  \mathscr{C}^0([0,T]; L^{\tilde{q}}_{\alpha, \Theta} H^{k-1}_x)$$ 
where $\tilde{q}=q$ if $q <+\infty$, and $\tilde{q} \in [1,+\infty)$ if $q=+\infty$, 
satisfying the continuity equation
    \begin{align*}
        \partial_t \rho^\alpha+\mathrm{div}(\rho^\alpha v^\alpha)=0,
    \end{align*} 
    with $\rho^\alpha_{\mid t=0} =\rho^\alpha_0$. Furthermore, for all $t \geq 0$,
    \begin{align}
        \label{estim:rhoalphaSob-sup} \Vert \rhorho(t)\Vert_{L^q_{\alpha, \Theta} H^{k-1}_x} &\leq  \Vert \rhorho_0\Vert_{L^q_{\alpha, \Theta} H^{k-1}_x} \exp\left( K \int_0^t \Vert \nabla \vv(s) \Vert_{L^{\upinfty}_\alpha H^{k-1}_x } \D s \right), \\
        \label{estim:rhoalphaSob-integ} \int \| \rho^\alpha(t) - 1\|_{H^{k-1}_x} \D \mu(\alpha) &\leq \left(1 +  \int \| \rho^\alpha_0 - 1 \|_{H^{k-1}_x} \D \mu(\alpha) \right)\\
       \notag &\qquad\times \exp\left( K \int_0^t \Vert \nabla \vv(s) \Vert_{L^{\upinfty}_\alpha H^{k-1}_x } \D s \right) - 1.
    \end{align}
\end{Lem}

\begin{proof}

    Let us first consider a smooth density--velocity pair $(\rho,v)$ satisfying the continuity equation
    \begin{align*}
        \partial_t \rho+\mathrm{div}(\rho v)=0.
    \end{align*} 
For any  $k \in \N {\setminus \lbrace 0 \rbrace}$, there holds 
    \begin{align}\label{estim:rhoalphaSob}
    \dfrac{\mathrm{d}}{\mathrm{d}t}\|\rho\|_{H^{k-1}_x} \lesssim \Vert \nabla v \Vert_{L^{\upinfty}_x} \|\rho\|_{H^{k-1}_x}+ \| \nabla v\|_{H^{k-1}_x}\|\rho\|_{L^{\upinfty}_x}.
    \end{align}
Indeed, applying $\partial_x^\gamma=\partial_{x_1}^{\gamma_1} \cdots \partial_{x_d}^{\gamma_d}$ to the equation, with $0 \leq \vert \gamma \vert \leq k-1$,  and writing
\begin{align*}
\partial_t \partial_x^\gamma\rho+\sum_{j=1}^d v_j \partial_{x_j}\partial^\gamma_x \rho +   \sum_{j=1}^d [\partial_{x_j}\partial^\gamma_x, v_j] \rho=0.
\end{align*}
Multiplying by $\partial_x^\gamma\rho$, integrating by parts, and applying the commutator estimate~\eqref{commut-law}, we thus obtain
\begin{multline*}
\dfrac{\mathrm{d}}{\mathrm{d}t}\| \partial_x^\gamma\rho \|_{L^2_x}^2 \\ 
\lesssim \Vert \nabla v \Vert_{L^{\upinfty}_x} \| \rho \|_{H^{k-1}_x}^2+ \left(\| \nabla v \|_{L^{\upinfty}_x} \| \rho \|_{H^{k-1}_x}+\| \nabla v \|_{H^{k-1}_x}\| \rho \|_{L^{\upinfty}_x} \right) \| \rho \|_{H^{k-1}_x}.
\end{multline*}
Summing on $0 \leq \vert \gamma \vert  \leq k-1$, we obtain the claimed estimate \eqref{estim:rhoalphaSob}.

 We can apply \eqref{estim:rhoalphaSob} for any single label $\alpha \in I$, and rely on Sobolev embedding to obtain the following classical a priori estimate on the continuity equation satisfied by $\rho^\alpha$:
\begin{equation*}
\dfrac{\mathrm{d}}{\mathrm{d}t}\|\rho^\alpha\|_{H^{k-1}_x} \lesssim \| \nabla v^\alpha\|_{H^{k-1}_x} \|\rho^\alpha\|_{H^{k-1}_x},
\end{equation*}
therefore by Gronwall lemma, we get
\begin{equation}\label{eq:dracolosse}
\|\rho^\alpha(t)\|_{H^{k-1}_x} \leq \|\rho_0^\alpha \|_{H^{k-1}_x} \exp\left( K\int_0^t \| \nabla v^\alpha(s) \|_{H^{k-1}_x} \D s \right),
\end{equation}
for some $K>0$. Taking the $L^q_{\alpha, \Theta}$ norm, we obtain the estimate \eqref{estim:rhoalphaSob-sup}.  Thanks to this a priori estimate, the existence and uniqueness of a family of solutions $\rhorho \in L^\infty([0,T]; L^q_{\alpha, \Theta} H^{k-1}_x)$ follows by the method of characteristics: 
defining the characteristics curves associated with the vector field $v^{\alpha}$, i.e \begin{align*} \dfrac{\mathrm{d}}{\mathrm{d}t} \mathrm{X}_{s;t}^\alpha(x)=v^{\alpha}(s,\mathrm{X}_{s;t}^\alpha(x)), \ \ \mathrm{X}_{t;t}^\alpha(x)=x, \end{align*} 
we have 
 \begin{align} \label{eq:charac-densité-gen}
 \rho^{\alpha}(t,x)= \rho^{\alpha}_0(\mathrm{X}_{0;t}^\alpha(x))\exp\left(-\int_0^t \mathrm{div}(v^{\alpha})(\tau, \mathrm{X}_{\tau;t}^\alpha(x)) \mathrm{d} \tau\right).\end{align}
 Classical arguments using the continuity in time of the velocity fields also ensure that for $\mu$-almost all $\alpha \in {I}$, we have $\rho^\alpha \in \mathscr{C}^0([0,T]; H^{k-1}_x)$. To show that $(\rho^{\alpha})_{\alpha \in {I}} \in \mathscr{C}^0([0,T]; L^{\tilde{q}}_{\alpha, \Theta} H^{k-1}_x)$ (where $\tilde{q}=q$ if $q <+\infty$, and $\tilde{q} \in [1,+\infty)$ otherwise), we can use Lebesgue's domination theorem exactly as in Item 3. of Corollary~\ref{cor:timecontinuity}, relying on the bound \eqref{eq:dracolosse}.

The last estimate \eqref{estim:rhoalphaSob-integ} is inferred in the same way: indeed, the equation satisfied by $r^\alpha\vcentcolon=\rho^\alpha-1$ is
$$
\partial_t r^\alpha+\mathrm{div}(r^\alpha v^\alpha)=-\mathrm{div}(v^\alpha),
\qquad
r^\alpha_{\mid t=0}=\rho_0^\alpha-1,
$$
and we can argue as before (applying Gronwall lemma to $\|r^\alpha(t)\|_{H^{k-1}_x}$).
\end{proof}

\paragraph{$L^p$ and pointwise bounds.}
We shall also use standard $L^p$ estimates on the solution to the continuity equation, that follows from a direct application of the method of characteristics (see \eqref{eq:charac-densité-gen} in the proof above).
\begin{Lem}\label{LM:pointwiseEstimRHO}
Assume that $(\rho,v)$ are smooth density and velocity satisfying
    \begin{align*}
        \partial_t \rho+\mathrm{div}(\rho v)=0,
    \end{align*} 
    with $\rho_{\mid t=0} =\rho_0$. Then for all $t \geq 0$ there holds
\begin{align*}
\Vert \rho(t) \Vert_{L^{1}_x}&= \Vert \rho_0 \Vert_{L^{1}_x},\\
\forall p \in (1, \infty], \ \ \Vert \rho(t) \Vert_{L^{p}_x}& \lesssim \Vert \rho_0 \Vert_{L^{p}_x} \exp\left(\int_0^t \Vert \mathrm{D} v(\tau) \Vert_{L^{\upinfty}_x} \mathrm{d} \tau\right), \\
\underset{x}{\inf} \, \rho(t,x) &\geq \left( \underset{x}{\inf} \,  \rho_0(x) \right) \exp\left(-\int_0^t\Vert \mathrm{D} v(\tau) \Vert_{L^{\upinfty}_x} \mathrm{d} \tau\right).
\end{align*}
In particular, if $\rho_0$ is nonnegative then $\rho(t)$ is nonnegative for all $t \geq 0$.
\end{Lem}

\paragraph{Stability at low regularity.}

The following lemma will prove useful when Assumption \ref{ass:stability} (about stability of the force field in $L^2_x$) is checked in practice.
\begin{Lem}\label{LM:estimStabContinuity}
Let $k \in \N$ and $q \in [1, \infty]$. There exists $K>0$ such that the following holds for any set of labels $(I, \mu)$. Let $(\rhorho_1,\vv_1)$ and $(\rhorho_2,\vv_2)$ be two families of smooth densities and velocity fields solutions to 
    \begin{align*}
        \partial_t \rho_i^\alpha+\mathrm{div}(\rho_i^\alpha v_i^\alpha)=0, \ \ i=1,2,
    \end{align*} 
    with initial data ${\rho^\alpha_{i}}_{\mid t=0} =\rho^\alpha_{i,0}$ and $\textbf{v}_i\in \mathcal H^{k,q}$ ($i=1,2$). Then for all measurable function $\Theta: I \rightarrow \R^+$ and all $t \geq 0$ 
	\begin{multline}
	\label{eq:rho_H-1}
	\| \rhorho_2(t) - \rhorho_1(t) \|_{L^q_{\alpha,\Theta} \dot H^{-1}_x} \leq \exp\left( K \int_0^t (\| \nabla \vv_1(s) \|_{L^{\upinfty}_\alpha L^{\upinfty}_x} + \| \nabla \vv_2(s) \|_{L^{\upinfty}_\alpha L^{\upinfty}_x}) \D s \right) \\ \left\{ \| \rhorho_{2,0} - \rhorho_{1,0} \|_{L^q_{\alpha,\Theta}  \dot H^{-1}_x} + \| \rhorho_{1,0}\|_{L^{\upinfty}_{\alpha,\Theta} L^{\upinfty}_x}\int_0^t \| \vv_2(s) - \vv_1(s) \|_{L^q_\alpha  L^2_x} \D s \right\}.
	\end{multline}
\end{Lem}
\begin{Def}
The homogeneous $\| \cdot \|_{\dot H^{-1}_x}$ (semi-)norm is defined as
$$
\| \varphi \|_{\dot H^{-1}_x} \vcentcolon = \sup_{\|\nabla \psi \|_{L^2_x}\leq 1}  \int \varphi \psi \, \D x.
$$
    \index{H@$\dot H^{-1}_x$: homogeneous Sobolev space of order $-1$ in $x$}
        \end{Def}

    \begin{Rem}By considering $\rhorho_{2}\equiv 0$ and $\vv_2\equiv 0$, one also obtains 
    that the solution    to
    \begin{align*}
        \partial_t \rho^\alpha+\mathrm{div}(\rho^\alpha v^\alpha)=0,
    \end{align*} 
    with $\rho^\alpha_{\mid t=0} =\rho^\alpha_0$, satisfies
        \begin{equation}
	\label{eq:rho_H-1-easy}
	\| \rhorho(t) \|_{L^q_\alpha \dot H^{-1}_x} \leq \exp\left( K \int_0^t \| \nabla \vv(s) \|_{L^{\upinfty}_\alpha L^{\upinfty}_x}  \D s \right)  \| \rhorho_{0} \|_{L^q_\alpha  \dot H^{-1}_x},
	\end{equation}
    an estimate we shall often refer to when invoking Lemma~\ref{LM:estimStabContinuity}.
    \end{Rem}
\begin{proof}[Proof of Lemma~\ref{LM:estimStabContinuity}]
We argue at fixed $\alpha$, then taking the $L^q_{\alpha,\Theta}$ norm in the end. As a consequence, we suppress the dependency in $\alpha$. For a given $t>0$, let $\phi(t)$ be the solution of the equation
 	\begin{equation*}
	- \Delta \phi(t) \vcentcolon = \rho_2(t) - \rho_1(t).
	\end{equation*}
	We then have
	\begin{equation*}
	\| \rho_2(t) - \rho_1(t) \|_{\dot H^{-1}_x} = \| \nabla \phi(t) \|_{L^2_x}.
  	\end{equation*}
	Therefore, dropping the time dependency for the sake of readability and denoting the $L^2_x$ scalar product with $\langle \cdot, \cdot \rangle$, we have
	\begin{align*}
	\frac{1}{2} \frac{\D}{\D t} \| \rho_2 - \rho_1 \|^2_{\dot H^{-1}_x}  = \cg  \nabla \phi, \partial_t \nabla \phi \cd = - \cg \phi, \partial_t \Delta \phi \cd = \cg \phi, \partial_t (\rho_2 - \rho_1) \cd, 
	\end{align*}
that is, using the continuity equation 
	\begin{align*}
	\frac{1}{2} \frac{\D}{\D t} \| \rho_2 - \rho_1 \|^2_{\dot H^{-1}_x}&= -\cg \phi , \Div(\rho_2 v_2 - \rho_1 v_1) \cd 
	=  \cg \nabla \phi_t , \rho_2 v_2 - \rho_1 v_1 \cd\\
	&=  \cg \nabla \phi, (\rho_2 - \rho_1)v_2 + (v_2 - v_1) \rho_1 \cd \\
	&= -\cg \Delta \phi \nabla \phi , v_2 \cd + \cg \nabla \phi, \rho_1 (v_2 - v_1) \cd,
	\end{align*}
hence
\begin{align}\label{goldorak}
	\frac{1}{2} \frac{\D}{\D t} \| \rho_2 - \rho_1 \|^2_{\dot H^{-1}_x}=-\cg \Div (\nabla \phi \otimes \nabla \phi) - \frac{1}{2} \nabla |\nabla \phi|^2, v_2 \cd + \cg \nabla \phi, \rho_1 (v_2 - v_1) \cd,
	\end{align}
	and then
	\begin{align*}
	 &\qquad   \frac{\D}{\D t} \| \rho_2 - \rho_1 \|^2_{\dot H^{-1}_x} \\
     &\leq  K \| \DDD v_2 \|_{L^{\upinfty}_x} \| \nabla \phi \|^2_{L^2_x} + \| \rho_1 \|_{L^{\upinfty}_x} \| \nabla \phi \|_{L^2_x} \| v_2 - v_1 \|_{L^2_x}\\
	 &= K \| \DDD v_2 \|_{L^{\upinfty}_x} \| \rho_2 - \rho_1 \|_{\dot H^{-1}_x}^2 + \| \rho_1 \|_{L^{\upinfty}_x} \| \rho_2 - \rho_1 \|_{\dot H^{-1}_x} \| v_2 - v_1 \|_{L^2_x}. 
	\end{align*}
The result follows by simplifying by $\| \rho_2 - \rho_1 \|_{\dot H^{-1}_x}$ and relying on the pointwise estimate of Lemma \ref{LM:pointwiseEstimRHO}.
\end{proof}

Given an initial data and a fixed family of velocity field, and equipped with the \textit{a priori} estimates from Lemma~\ref{LM:SobEstimRHO} and from Lemma~\ref{LM:estimStabContinuity}, one can now classically solve the continuity equations for the densities $(\rho^\alpha)$. We record it in the following definition.
\begin{Def}\label{def:rho-given-v_and_data}
Let $k>d/2+1$, $p\in [1,+\infty]$, $T>0$ and $\vv \in L^\infty([0,T]; \mathcal{H}^{k,p})$. Let  $\rhorho_0 \in L^q_{\alpha, \Theta} H^{k-1}_x$ with $\Theta :I \rightarrow \R^+$ a measurable map. We denote by $\rhorho[\rhorho_0, \vv]$ the unique family $\rhorho \in L^\infty([0,T];L^\infty_{\alpha, \Theta} H^{k-1}_x)$ of solutions to the family of continuity equations on $[0,T]$:
	\begin{equation}
	\left\{
	\begin{gathered}
	\partial_t \rho^\alpha + \Div( \rho^\alpha v^\alpha) =0, \\
	\rho^\alpha |_{t = 0} = \rho^\alpha_0.
	\end{gathered}
	\right.
	\end{equation}
\end{Def}
For such a family of solution $\rhorho[\rhorho_0, \vv]$, all related estimates they satisfy are directly inferred from~\eqref{estim:rhoalphaSob-sup} in Lemma~\ref{LM:SobEstimRHO} for the bound in high regularity, and~\eqref{eq:rho_H-1} in Lemma~\ref{LM:estimStabContinuity} for the stability in $\dot{H}^{-1}_x$.

\begin{Rem}
\label{rem:Rd-alpha2}
As a final comment, following Remark~\ref{rem:Rd-alpha}, we note that all Lemmas~\ref{LM:SobEstimRHO}, \ref{LM:pointwiseEstimRHO} and \ref{LM:estimStabContinuity} also hold when considering the slightly more general continuity equation
$$
\partial_t \rho^\alpha + \Div( \rho^\alpha (w^\alpha+\lambda^\alpha)) =0,
$$
in which we recall $\lambda: (I,\mu) \rightarrow \R^d$ is a measurable map which does not depend on $t$ or $x$.
\end{Rem}

\section{Consequences for multiphasic density--velocity formulations}
\label{sec:conseq-rho}

In this section, we deduce consequences of the general local well-posedness statement from Theorem \ref{thm:existence}: first, we lift this result to the complete multiphasic system \eqref{eq:multiphasic-equivalence}  on density--velocity fields, namely
	\begin{equation}
	\label{eq:general_system_Chapter2}
	\left\{
	\begin{gathered}
	\partial_t \rho^\alpha + \Div( \rho^\alpha v^\alpha) =0, \\
	\partial_t v^\alpha + ( v^\alpha \cdot \nabla ) v^\alpha = F^\alpha,\\
	\rho^\alpha |_{t = 0} = \rho^\alpha_0, \ v^\alpha |_{t = 0} = v^\alpha_0,
	\end{gathered}
	\right.
	\end{equation}
and then we deduce a related statement on the Cauchy theory at the level of Vlasov equations. 

\medskip

 Even though the family of densities $(\rho^\alpha)_{\alpha \in I}$ solution to the first equation of \eqref{eq:general_system_Chapter2} is entirely determined by the family of velocities $(v^\alpha)_{\alpha \in I}$ and initial conditions $(\rho^\alpha_{0})_{\alpha \in I}$ (see Definition \ref{def:rho-given-v_and_data}), it turns out convenient in practice to see the force as a function of $\rhorho$ and $\vv$. 
 In this setting, the parameter is $X=(\rhorho_0,Y)$ where $Y \in \mathcal{Y}$ is another extra-parameter, and we write the family of forces $\FF[\vv,\rhorho_0,Y]$ as $\widetilde\FF[\vv, \rhorho,Y]$. 
 
 In this setting, with $k_0 > 1+ d/2$ 
 and $1\leq p \leq +\infty$ being given, the Assumptions~\ref{ass:bound_high_regularity} and~\ref{ass:stability} reformulate as follows:

	\begin{Ass}{B1}[Bound in high regularity]
		\label{ass:bound_high_regularity-rho}
          \index{A2@Assumptions \ref{ass:bound_high_regularity-rho}--\ref{ass:stability-rho}: for the  $(\rhorho,\vv)$ framework ($H^k$ spaces)}
	We assume that for all $k \geq k_0$, there exists an estimate function $\F$ such that for all $T>0$, for all $\vv \in L^1([0,T]; \mathcal H^{k,p})$, and all parameters $Y \in \mathcal Y$,
		we have
		\begin{multline}
		\label{eq:ass_estimate_high_reg-rho}
		\int_0^T \left\Vert\widetilde{\FF}[\vv, \rhorho,  Y](t)\right\Vert_{\mathcal{H}^{k,p}} \D t \\
        \leq \F\left(T, \int_0^T \left\Vert\vv(t)\right\Vert_{\mathcal{H}^{k,p}} \D t,  \| \rhorho\|_{L^\infty([0,T]; L^{\upinfty}_\alpha H^{k-1}_x)}, \| \rhorho\|_{L^\infty([0,T]; L^{\upinfty}_\alpha \dot{H}^{-1}_x)},  \mathsf N_k(Y) \right).
		\end{multline}
	\end{Ass}

	\begin{Ass}{B2}[Stability at low regularity]
		\label{ass:stability-rho}
	Let $R>0$. There exists a stability estimate function $\G_R$ such that for all $T>0$, families $\vv_1,\vv_2 \in L^1([0,T]; \mathcal H^{k_0,p})$, 
	$\rhorho_1,\rhorho_2 \in L^\infty([0,T]; L^{\upinfty}_\alpha H^{k_0-1}_x)$,
	and parameters $Y_1$ and $Y_2$, if
		\begin{multline}
		\label{eq:uniform_bound_velocities-rho} \max\Bigg( \sup_{t \in [0,T]} \Vert \vv_1(t) \Vert_{\mathcal{H}^{k_0,p}},\sup_{t \in [0,T]} \Vert \vv_2(t) \Vert_{\mathcal{H}^{k_0,p}}, \\
        \| \rhorho_1 \|_{L^\infty([0,T]; L^{\upinfty}_\alpha H^{k_0-1}_x)},\mathsf N_{k_0}(Y_1), \mathsf N_{k_0}(Y_2)  \Bigg) 
		\leq R,
		\end{multline}
	then the following estimate holds:
		\begin{multline}
		\label{eq:stability_assumption-rho}
		\int_0^T \left\Vert\widetilde{\FF}[\vv_1, \rhorho_1, Y_1](t) - \widetilde{\FF}[\vv_2, \rhorho_2,Y_2](t)\right\Vert_{L^p_\alpha L^2_x} \D t \\
		\leq \G_R\left(T, \int_0^T \big\|\vv_1(t) - \vv_2(t) \big\|_{L^p_\alpha L^2_x} \D t, \|\rhorho_1 - \rhorho_2\|_{L^\infty([0,T];L^p
_\alpha \dot{H}^{-1}_x)}, \mathsf d(Y_1,Y_2)\right).
		\end{multline}
	\end{Ass}

\begin{Rem}\label{rem-bound-high-ref-functionQ-withrho}
    In light of Remark \ref{rem-bound-high-ref-functionQ}, we actually allow for a generalization of \eqref{eq:ass_estimate_high_reg-rho} where we require an estimate of of the form 
    \begin{multline}
		\int_0^T \left\Vert\widetilde{\FF}[\vv, \rhorho,  Y](t)\right\Vert_{\mathcal{H}^{k,p}} \D t \\
        \leq \F\left(T, \int_0^T \mathcal{Q}(\left\Vert\vv(t)\right\Vert_{\mathcal{H}^{k,p}}) \D t,  \| \rhorho\|_{L^\infty([0,T]; L^{\upinfty}_\alpha H^{k-1}_x)}, \| \rhorho\|_{L^\infty([0,T]; L^{\upinfty}_\alpha \dot{H}^{-1}_x)},  \mathsf N_k(Y) \right),
		\end{multline}
        where $\mathcal{Q}:\R^+ \mapsto \R^+$ is a smooth non-decreasing function such that $\mathcal{Q}(0)=0$.
\end{Rem}

\begin{Rem}The $\dot{H}^{-1}_x$ control for $\rhorho$ in Assumption~\ref{ass:bound_high_regularity-rho} is  relevant only  for the whole space case (and not for the torus case) for certain equations, as it provides an additional {\it low frequency} control. Typically this is important for the Vlasov--Poisson system, however such a control is not required for  the Vlasov--Navier--Stokes system. 

In practice, according to the Sobolev embedding, the initial $\dot{H}^{-1}_x$ bound for
densities may be ensured by assuming boundedness in a Lebesgue space $L^r_x$ for $r \in [1,2d/(d+2)]$ when $d\geq 3$.
\end{Rem}

We can now infer the following theorem.
\begin{Thm}\label{thm-LWPmultiphase-densityvelocity}
Let $(I, \mu)$ be a set of labels. Let $k>1+d/2$ and $p \in [1, \infty]$.  Suppose that there exists a well defined mapping $(\textbf{w}, \mathbf{\eta}) \mapsto \widetilde{\FF}(\textbf{w}, \mathbf{\eta},Y)$ satisfying the Assumptions~\ref{ass:bound_high_regularity-rho} and~\ref{ass:stability-rho}. If $(\rhorho_0, \vv_0) \in L^\infty_\alpha (H^{k-1}_x \cap \dot{H}^{-1}_x) \times \mathcal{H}^{k,p}$ is a regular multiphasic distribution, there exists $T>0$ and a unique solution $(\rhorho, \vv)$ to the multiphasic system with force field $\widetilde{\FF}(\vv,\rhorho^,Y)$ on $[0,T]$ in the sense of Definition \ref{def:solution_coupled} and such that 
\begin{align*}
    \rhorho \in L^\infty([0,T]; L^\infty_\alpha (H^{k-1}_x\cap\dot{H}^{-1}_x) ), \ \ \vv \in L^\infty([0,T]; \mathcal{H}^{k,p}).
\end{align*}
Furthermore, the following continuity in time  properties hold: we have
\begin{align*}
    \vv \in {\mathscr{C}^0}([0,T]; L^p_\alpha L^2_x) \ \ \text{and} \ \ \DDD \vv \in {\mathscr{C}^0}([0,T]; L^q_\alpha H^{k-1}_x)  \ \ \text{for all} \ \  q \in [1,\infty), \\
    \rhorho \in   \mathscr{C}^0([0,T]; L^{q}_\alpha H^{k-1}_x) \ \ \text{for all} \ \  q \in [1,\infty).
\end{align*}
\end{Thm}

\begin{Rem}\label{rem-continuity-thmgendensityvelocity}
    In view of Theorem \ref{thm:existence}, we also obtain that for $\mu$-almost every $\alpha \in I$, $v^\alpha \in \mathscr{C}^0([0,T]; H^k_x)$.  Furthermore, in view of Remark \ref{Rem:timeuniform}, we also have the continuity on $[0,T]$ of the map $t \mapsto \Vert \vv(t) \Vert_{\mathcal{H}^{k,p}}$.
    
\end{Rem}

\begin{Rem} In the case when the control by $\| \rhorho\|_{L^\infty([0,T]; L^\infty_\alpha \dot{H}^{-1}_x)}$ is not required in the argument of the estimate function $\mathcal F$ in Assumption~\ref{ass:bound_high_regularity-rho}, Theorem~\ref{thm-LWPmultiphase-densityvelocity} still holds by removing all occurrences of $\dot H^{-1}_x$.    
\end{Rem}

\begin{Rem}
\label{rem:smoothkernel} There are situations where in  Assumptions~\ref{ass:bound_high_regularity-rho}--\ref{ass:stability-rho}, only a control of the norm $\| \rhorho\|_{L^\infty([0,T]; L^\infty_\alpha L^q_x)}$ (for some $q \in [1,2]$) is needed in the estimates~\eqref{eq:ass_estimate_high_reg-rho} and \eqref{eq:uniform_bound_velocities-rho} (we leave \eqref{eq:stability_assumption-rho} unchanged). This happens in particular for a class of nonlinear Vlasov equations arising from particle systems with \emph{smooth} interaction kernels, namely
$$
\partial_t f + v\cdot \nabla_x f -( \nabla_x  K \star_x \rho_f) \cdot \nabla_v f =0,
$$
where the interaction kernel $K(x)$ is supposed to be in $W^{k+1,q^\star}(\R^d)$ for $k >1+ d/2$, where $q^\star$ satisfies $1+ \frac{1}{2}= \frac{1}{q}+\frac{1}{q^\star}$. The associated multiphasic system of this nonlinear Vlasov equation is
	\begin{equation}
	\label{eq:multiphasic-smookernel}
	\left\{
	\begin{gathered}
	\partial_t \rho^\alpha + \Div( \rho^\alpha v^\alpha) =0, \\
	\partial_t v^\alpha + ( v^\alpha \cdot \nabla ) v^\alpha = -\nabla_x  K \star_x \left[\int \rho^\alpha \, \mathrm{d}\mu(\alpha) \right] ,\\
	\rho^\alpha |_{t = 0} = \rho^\alpha_0, \ v^\alpha |_{t = 0} = v^\alpha_0.
 	\end{gathered}
	\right.
	\end{equation}
It is clear that by Young inequality for convolutions, the force in this equation, which we denote by $F$, satisfies
$$
\| F \|_{H^{k}_x} \lesssim \| K\|_{W^{k+1,q^\star}_x} \| \rhorho\|_{L^1_\alpha L^q_x},
$$
hence the simplified Assumptions~\ref{ass:bound_high_regularity-rho}--\ref{ass:stability-rho}.
    
In this framework, a variant of Theorem~\ref{thm-LWPmultiphase-densityvelocity} holds without asking $\rhorho_0$ to be this regular; namely if $\rhorho_0  \in L^\infty_\alpha L^q_x$ then the same statement holds with the propagation of this norm:
$$
  \rhorho \in L^\infty([0,T]; L^\infty_\alpha L^q_x ).
$$

\end{Rem}
We refer to Section \ref{section:Extension} for extensions of this result to other functional spaces.

\begin{proof}[Proof of Theorem \ref{thm-LWPmultiphase-densityvelocity}]
    First, let us recall the notation $\rhorho[\rhorho_0,\vv]$ from Definition \ref{def:rho-given-v_and_data}, as well as the associated estimates~\eqref{estim:rhoalphaSob-sup} and~\eqref{eq:rho_H-1}. We aim at applying Theorem~\ref{thm:existence}. To this end, it is enough to check that $\FF[\vv,\rhorho_0,Y] \vcentcolon = \widetilde\FF[\vv,\rhorho[\vv, \rho_0],Y]$ satisfies  Assumptions~\ref{ass:bound_high_regularity} and~\ref{ass:stability}, with $$X= (\rhorho_0,Y),\qquad \mathcal{X}= L^\infty_\alpha H^{k_0-1}_x \cap \dot{H}^{-1}_x \times \mathcal{Y}, $$
    $$N_k(\rhorho_0, Y)= \| \rhorho_0\|_{L^{\upinfty}_\alpha (H^{k-1}_x \cap \dot{H}^{-1}_x)} + N_k(Y),$$
    and
        $$
        d((\rhorho_{1,0}, Y_1),(\rhorho_{2,0}, Y_2)) = \|  \rhorho_{1,0}- \rhorho_{2,0}\|_{L^p_\alpha \dot{H}^{-1}_x} + d (Y_1,Y_2).
        $$
        This will follow directly from  estimates for continuity equations as developed in the previous section, namely Lemmas~\ref{LM:SobEstimRHO} and \ref{LM:estimStabContinuity}.

\medskip
    
\noindent{\bf Assumption \ref{ass:bound_high_regularity-rho} implies Assumption \ref{ass:bound_high_regularity}.}  Assume that Assumption \ref{ass:bound_high_regularity-rho} holds.  By~\eqref{eq:ass_estimate_high_reg-rho}, one has
\begin{multline*}
    		\int_0^T \left\Vert\FF[\vv,\rhorho_0,Y](t)\right\Vert_{\mathcal{H}^{k,p}} \D t \\
            \leq \F\left(T, \int_0^T \left\Vert\vv(t)\right\Vert_{\mathcal{H}^{k,p}} \D t,  \| \rhorho[\vv, \rho_0]\|_{L^\infty([0,T]; L^{\upinfty}_\alpha (H^{k-1}_x)\cap \dot{H}^{-1}_x)},  \mathsf N_k(Y) \right).
\end{multline*}
Thanks to~\eqref{estim:rhoalphaSob-sup} and \eqref{eq:rho_H-1-easy}, one can thus change $\mathcal{F}$ into another estimate function $\widetilde{\mathcal{F}}$ such that
\begin{equation*}
    		\int_0^T \left\Vert\FF[\vv,\rhorho_0,Y](t)\right\Vert_{\mathcal{H}^{k,p}} \D t \leq \widetilde{\mathcal{F}}\left(T, \int_0^T \left\Vert\vv(t)\right\Vert_{\mathcal{H}^{k,p}} \D t,\mathsf N_k(\rhorho_0, Y) \right).
\end{equation*}
This shows that Assumption \ref{ass:bound_high_regularity} is satisfied.

\medskip

\noindent{\bf Assumption \ref{ass:stability-rho} implies Assumption \ref{ass:stability}.}   Assume that Assumption \ref{ass:stability-rho}  holds. Let us assume that we have
$$ \max\left( \sup_{t \in [0,T]} \Vert\vv_1(t)\Vert_{\mathcal{H}^{k_0,p}},\sup_{t \in [0,T]} \Vert\vv_2(t)\Vert_{\mathcal{H}^{k_0,p}}, \mathsf N_{k_0}(X_1), \mathsf N_{k_0}(X_2)  \right) \leq R,
$$
where $X_i = (\rhorho_{0,i},Y_i)$, for some $R>0$. By~\eqref{estim:rhoalphaSob-sup} we thus have
\begin{multline*}
    \max\Bigg( \sup_{t \in [0,T]} \Vert\vv_1(t) \Vert_{\mathcal{H}^{k_0,p}},\sup_{t \in [0,T]} \Vert\vv_2(t) \Vert_{\mathcal{H}^{k_0,p}}, \\
    \| (\rhorho[\vv_1,\rhorho_{1,0}],\rhorho[\vv_2,\rhorho_{2,0}]) \|_{L^\infty([0,T]; L^{\upinfty}_\alpha H^{k_0-1}_x)},
    \mathsf N_{k_0}(Y_1), \mathsf N_{k_0}(Y_2)  \Bigg) 
\leq \Lambda(T, R),
\end{multline*}
for a continuous function $\Lambda$. Since this is an estimate of the form~\eqref{eq:uniform_bound_velocities-rho},  we are in position to apply Assumption \ref{ass:stability-rho}, which yields that
	\begin{multline*}
		\int_0^T \left\Vert{\FF}[\vv_1, \rhorho_{1,0}, Y_1](t) - {\FF}[\vv_2, \rhorho_{2,0},Y_2](t)\right\Vert_{L^p_\alpha L^2_x} \D t \\
		\leq \G_{R}\left(T, \int_0^T \big\|\vv_1(t) - \vv_2(t) \big\|_{L^p_\alpha L^2_x} \D t, \|\rhorho_1 - \rhorho_2\|_{L^\infty([0,T];L^p
_\alpha \dot{H}^{-1}_x)}, \mathsf d(Y_1,Y_2)\right).
		\end{multline*}
Owing  to the stability estimate in $\dot{H}^{-1}_x$ for a family of multiphasic continuity equation, namely  \eqref{eq:rho_H-1}, we deduce that the latter turns into
	\begin{multline*}
		\int_0^T \left\Vert{\FF}[\vv_1, \rhorho_{1,0}, Y_1](t) - {\FF}[\vv_2, \rhorho_{2,0},Y_2](t)\right\Vert_{L^p_\alpha L^2_x} \D t \\
        \leq \widetilde{\mathcal{G}_R}\left(T, \int_0^T \big\|\vv_1(t) - \vv_2(t) \big\|_{L^p_\alpha L^2_x} \D t, \mathsf d((\rhorho_{1,0}, Y_1),(\rhorho_{2,0}, Y_2))\right),
		\end{multline*}
 for yet another stability estimate function $\widetilde{\mathcal{G}_R}$. We can therefore conclude that Assumption \ref{ass:stability} holds.

\end{proof}

In the wake of Remarks~\ref{rem:Rd-alpha}--\ref{rem:Rd-alpha2}, Theorem~\ref{thm-LWPmultiphase-densityvelocity} has a direct analogue for the system
	\begin{equation}
	\label{eq:general_system_Chapter2-re}
	\left\{
	\begin{gathered}
	\partial_t \rho^\alpha + \Div( \rho^\alpha (v^\alpha+\lambda^\alpha) =0, \\
	\partial_t w^\alpha + ( (w^\alpha+\lambda^\alpha) \cdot \nabla ) w^\alpha = G^\alpha[\ww,\rhorho,Y],\\
	\rho^\alpha |_{t = 0} = \rho^\alpha_0, \ w^\alpha |_{t = 0} = w^\alpha_0,
	\end{gathered}
	\right.
	\end{equation}
    for $(\rhorho_0, \ww_0) \in L^1_\alpha (H^{k-1}_x \cap \dot{H}^{-1}_x) \times \mathcal{H}^{k,\infty}$ and in which we recall $\lambda: (I,\mu) \rightarrow \R^d$ is a measurable map which does not depend on $t$ or $x$.
    As explained in  Remark~\ref{rem:Rd-alpha}, this is particularly relevant for the whole space case in view of an application to Vlasov equations, that is when we write a distribution function $f_0$ as
$$f_0(x, \cdot)=\int_{I} \rho^\alpha_0(x)  \otimes \delta_{v=v^\alpha_0(x)} \, \mathrm{d}\mu(\alpha).$$
Compared to the assumptions of Theorem~\ref{thm-LWPmultiphase-densityvelocity}, we note that we have upgraded the integrability assumption for $\rhorho_0$ from $L^\infty_\alpha$ to $L^1_\alpha$, which is better when translated to kinetic distribution functions. On the other hand,  we have restricted to $p=\infty$ in the $\mathcal{H}^{k,p}$ space for $\ww$: but as the velocity vector field of interest is $\lambda^\alpha+\ww^\alpha$, the latter assumption {\it in fine} does not imply that we restrict to bounded velocities on the kinetic side.

Let us state precisely the assumptions and result for~\eqref{eq:general_system_Chapter2-re}. 

	\begin{Ass}{C1}[Bound in high regularity]
		\label{ass:bound_high_regularity-rho-re}
          \index{A3@Assumptions \ref{ass:bound_high_regularity-rho-re}--\ref{ass:stability-rho-re}: for the $(\rhorho,\vv)$ framework, around constant phases ($H^k$ spaces)}
	 We assume that for all $k \geq k_0$, there exists an estimate function $\F$  and a measurable map $\Theta: (I,\mu) \to \R_+ $ such that for all $T>0$, for all $\ww \in L^1([0,T]; \mathcal H^{k,\infty})$, and all parameters $Y \in \mathcal Y$,
		we have
		\begin{multline}
		\label{eq:ass_estimate_high_reg-rho-re}
		\int_0^T \left\Vert \GG [\ww, \rhorho,  Y](t)\right\Vert_{\mathcal{H}^{k,\infty}} \D t \\
        \leq \F\left(T, \int_0^T \Vert \ww(t) \Vert_{\mathcal{H}^{k,\infty}} \D t,  \| \rhorho\|_{L^\infty([0,T]; L^1_{\alpha,\Theta} (H^{k-1}_x\cap \dot{H}^{-1}_x)},  \mathsf N_k(Y) \right).
		\end{multline}
	\end{Ass}

    \begin{Rem}
    The point of the weight $\Theta$ is to handle models involving moments in $v$ of high order; typically when $I=\R^d$ one may take $\Theta(\alpha) = (1+|\alpha|^2)^{p/2}$ if moments of order $p$ come into play. This will be handy for the Vlasov--Navier--Stokes equations which involves moments of order $1$. In addition, this will prove useful to derive some regularity of the moments in $v$ of the kinetic distribution $f=\mathcal{I}_\mu(\rhorho, \vv)$ -- see Section \ref{sec:conseq-kin}.
    
    \end{Rem}

	\begin{Ass}{C2}[Stability at low regularity]
		\label{ass:stability-rho-re} 
	Let $R>0$. There exists a stability estimate function $\G_R$ such that for all $T>0$, families $\ww_1,\ww_2 \in L^1([0,T]; \mathcal H^{k_0,\infty})$, 
	$\rhorho_1,\rhorho_2 \in L^\infty([0,T]; L^1_{\alpha,\Theta} H^{k_0-1}_x)$ (where $\Theta: (I,\mu) \to \R_+$ is  associated with Assumption~\ref{ass:bound_high_regularity-rho-re}),
	and parameters $Y_1$ and $Y_2$, if
		\begin{multline}
		\label{eq:uniform_bound_velocities-rho-re} \max\Bigg( \sup_{t \in [0,T]} \Vert\ww_1(t)\Vert_{\mathcal{H}^{k_0,\infty}},\sup_{t \in [0,T]} \Vert \ww_2(t)\Vert_{\mathcal{H}^{k_0,\infty}}, \\
        \| \rhorho_1\|_{L^\infty([0,T]; L^1_{\alpha,\Theta} H^{k_0-1}_x)},
        \mathsf N_{k_0}(Y_1), \mathsf N_{k_0}(Y_2)  \Bigg) 
		\leq R,
		\end{multline}
	then the following estimate holds:
		\begin{multline}
		\label{eq:stability_assumption-rho-re}
		\int_0^T \big\|{\GG}[\ww^1, \rhorho_1,Y_1](t) - {\GG}[\ww^2, \rhorho_2,Y_2](t)\big\|_{L^{\upinfty}_\alpha L^2_x} \D t \\
		\leq \G_R\left(T, \int_0^T \big\|\ww_1(t) - \ww_2(t) \big\|_{L^{\upinfty}_\alpha L^2_x} \D t, \right. \\
        \left. \|\rhorho_1 - \rhorho_2\|_{L^\infty([0,T];L^1
_{\alpha,\Theta} \dot{H}^{-1}_x)}, \mathsf d(Y_1,Y_2)\right).
		\end{multline}
	\end{Ass}

\begin{Rem}
    Note that, in the case $\Theta=1$, since  $L^\infty_\alpha$ embeds in $L^1_\alpha$ (the measure $\mu$ being of finite mass), the previous Assumptions \ref{ass:bound_high_regularity-rho-re}--\ref{ass:stability-rho-re}, when satisfied, imply the same estimates with $L^1_{\alpha,1}$ replaced everywhere by $L^\infty_\alpha$. 
\end{Rem}
The result is now the following:
\begin{Thm}\label{thm-LWPmultiphase-densityvelocity-re}
Let $(I, \mu)$ be a set of labels and let $\Theta: (I,\mu) \to \R_+ $ be the measurable map. Let $k>1+d/2$.  Suppose that there exists a well defined mapping $(\textbf{w}, \mathbf{\eta}) \mapsto {\GG}(\textbf{w}, \mathbf{\eta},Y)$ satisfying the Assumptions~\ref{ass:bound_high_regularity-rho-re} and~\ref{ass:stability-rho-re}. If $(\rhorho_0, \ww_0) \in L^1_{\alpha,\Theta} (H^{k-1}_x \cap \dot{H}^{-1}_x) \times \mathcal{H}^{k,\infty}$ is a regular multiphasic distribution, there exists $T>0$ and a unique solution $(\rhorho, \ww)$ to the multiphasic system \eqref{eq:general_system_Chapter2-re} with force field ${\GG}(\vv,\rhorho,Y)$ on $[0,T]$ in the sense of Definition \ref{def:solution_coupled} and such that 
\begin{align*}
    \rhorho \in L^\infty([0,T]; L^1_{\alpha,\Theta} (H^{k-1}_x\cap\dot{H}^{-1}_x) ), \ \ \ww \in L^\infty([0,T]; \mathcal{H}^{k,\infty}).
\end{align*}
\end{Thm}

\begin{Rem}
    Time continuity properties similar to the ones stated in Theorem \ref{thm-LWPmultiphase-densityvelocity} also hold in that case, with the suitable addition of a weight $\Theta(\alpha)$. In addition, we can replace all the spaces $L^1_{\alpha,\Theta}$ by $L^\infty_\alpha$ in the statement, with the same replacements in  Assumptions~\ref{ass:bound_high_regularity-rho-re} and~\ref{ass:stability-rho-re}.
\end{Rem}

\begin{proof}
    We argue exactly as in the proof of Theorem~\ref{thm-LWPmultiphase-densityvelocity}, using Remarks~\ref{rem:Rd-alpha}--\ref{rem:Rd-alpha2}. To obtain $L^1_{\alpha,\Theta}$ type estimates for densities, we note that $\tilde{\rho}^\alpha(t,x) \vcentcolon = \Theta(\alpha) \rho^\alpha(t,x)$ satisfies the same continuity equation as $\rho^\alpha$ and thus we can apply Lemmas~\ref{LM:SobEstimRHO}--\ref{LM:estimStabContinuity}.
\end{proof}

\section{Consequences for kinetic formulations}
\label{sec:conseq-kin}

\subsection{Local well-posedness for the Vlasov equation}
In this section, we draw direct consequences of the multiphasic abstract results obtained in Section \ref{sec:conseq-rho} at the level of the nonlinear Vlasov equation.

First, we derive the precise regularity of the moments in velocity of the multiphasic ansatz $\mathcal{I}_\mu(\rhorho, \vv)$ associated to the solutions obtained in Theorem \ref{thm-LWPmultiphase-densityvelocity}. 
\begin{Lem}\label{prop-reg-momentsf}
 Let $k>1+d/2$ and $p \in [1, \infty]$. Let $(I, \mu)$ be a set of labels such that $(\rhorho, \vv)$ is an associated multiphasic density/velocity family satisfying
\begin{align*}
    \rhorho \in L^\infty_\alpha H^{k-1}_x, \ \ \vv \in \mathcal{H}^{k,p}.
\end{align*}
 Then the following holds for the distribution function $f\vcentcolon=\mathcal{I}_\mu(\rhorho, \vv)$:
\begin{itemize}
    \item for any test function $\psi \in 
		\mathscr{C}^\infty_b(\R^d)$, the function
		\begin{equation*}
		m_\psi(f): \qquad x \mapsto \int \psi(v) f(x,\D v)  = \int \psi(v^\alpha(x)) \rho^\alpha(x) \D \mu(\alpha)
		\end{equation*}
        satisfies
        \begin{align}\label{reg-moment-mult-psi}
            m_\psi(f) \in H^{k-1}_x.
        \end{align}
    \item the moments in velocity
    \begin{align*}
    m_j(f) : \qquad x \mapsto \int \langle v \rangle^j f(x,\D v)=\int \langle v^\alpha(x) \rangle^j \rho^\alpha(x) \D \mu(\alpha)
    \end{align*}
    for $j=0,\ldots,\lfloor p\rfloor$, satisfy
    \begin{align}\label{reg-moment-mult-vj}
        m_j(f) \in H^{k-1}_x.
    \end{align}
\end{itemize}
\end{Lem}
\begin{proof}
     Let $\psi \in \mathscr{C}^\infty_b(\R^d)$. By Leibniz rule, it is enough to estimate terms of the form $\partial_x^\gamma \rho^\alpha \partial_x^\beta [\psi(v^\alpha)]$ in $L^2_x$ where $ \beta,\gamma \in \N^d$ with $|\beta|+|\gamma|=k-1$. 
     Assume first that $|\beta|= k-1- |\gamma|>d/2$. Then by Sobolev embedding, we bound $\Vert \partial^\gamma_x    \rhorho\Vert_{ L^{\upinfty}_\alpha L^{\upinfty}_x} \lesssim \Vert \rhorho \Vert_{L^{\upinfty}_\alpha H^{k-1}_x}$ as $\vert \gamma \vert +d/2 < k-1$.
     Since $\psi \in \mathscr{C}^\infty_b(\R^d)$, owing to the Fa\`a di Bruno formula, the term $\| \partial^\beta_x \psi(v^\alpha) \|_{L^2_x} $ can be estimated by a sum of terms of the form
     \begin{equation}
     \label{eq:termesbeta}
    \| \partial^{\beta_1}_x v^\alpha \cdots   \partial^{\beta_\ell}_x v^\alpha \|_{L^2_x}, \qquad \beta_1+ \ldots + \beta_\ell= \beta, \quad |\beta_1|\geq \ldots \geq |\beta_\ell|>0,
     \end{equation}
     for some integer $\ell \geq 1$. Consider
     \begin{align}\label{def:n-bourbakiste}
         n\vcentcolon=\max \left\lbrace j \in [\![1, \ell]\!] \mid \vert \beta_j\vert \geq k-d/2 \right\rbrace
     \end{align}
     if the former set is non-empty, and $n=0$ else. If $n \neq 0$ then $\vert \beta_j \vert \geq k-d/2$ for $1 \leq j \leq n$ and $\vert \beta_j \vert > k-d/2$ for $n< j \leq \ell$.
     By Sobolev embedding, we bound for $j>n$
    \begin{align*}
    \| \partial^{\beta_j}_x v^\alpha \|_{L^{\upinfty}_x} \lesssim \| \DDD v^\alpha \|_{H^{k-1}_x}.    
    \end{align*}
     We apply H\"older's inequality to treat the remaining terms of the product, that is for $1 \leq j \leq n$. Given $q_i\geq 2$ such that 
     $\frac{1}{2} = \sum_{j=1}^n \frac{1}{q_j}$, we have
     $$
         \| \partial^{\beta_1}_x v^\alpha \cdots   \partial^{\beta_n}_x v^\alpha \|_{L^2_x} \leq     \| \partial^{\beta_1}_x v^\alpha \|_{L^{q_1}_x}\cdots   \| \partial^{\beta_n}_x v^\alpha \|_{L^{q_n}_x}.
     $$
 In order to control the Lebesgue norms, we can rely on the Sobolev embedding $H^{k-\vert \beta_j\vert }_x \hookrightarrow L^{q_j}_x $ that holds since $0<k-\vert \beta_j \vert<d/2$ for $1 \leq j \leq n$ and provided that the $q_j$ satisfy 
 \begin{align*}
      \frac{1}{2} \geq \frac{1}{q_j} \geq \frac{1}{2}-\frac{k-\vert \beta_j \vert}{d}, \ \ 1 \leq j \leq n.
 \end{align*} 
To pick such exponents $q_j$, it is sufficient to check that we have
\begin{align*}
    \frac{1}{2}> \sum_{j=1}^n \left(\frac{1}{2} - \frac{k-|\beta_j|}{d}\right).
\end{align*}
     As $\displaystyle \sum_{j=1}^n |\beta_j|\leq |\beta|\leq k-1$, this holds since the former right-hand side is equal to
\begin{align*}
    \frac{n}{2}-\frac{nk}{d}+\frac{1}{d}\sum_{j=1}^n \vert \beta_j \vert \leq \frac{n}{2}-\frac{nk}{d}+ \frac{k-1}{d}=\frac{1}{2}+(n-1)  \left[\frac{1}{2}-\frac{k}{d}\right]-\frac{1}{d}<\frac{1}{2},
\end{align*}
since $n \geq 1$ and $k>d/2$. We end up with
      \begin{equation}
      \label{eq:ouf0}
         \| \partial^{\beta_1}_x v^\alpha \cdots   \partial^{\beta_\ell}_x v^\alpha \|_{L^2_x} \lesssim    \| \DDD v^\alpha \|_{H^{k-1}_x}^\ell.
     \end{equation}
     If $n=0$, then we have $\vert \beta_j \vert <k-d/2$ for all $j \in [\![1, \ell]\!]$, and we proceed as before, putting all the factor in $L^\infty_x$ except one that we keep in $L^2_x$.

   Assume otherwise that $|\beta|= k-1- |\gamma|\leq d/2$. We argue exactly as before, except that this time we study
       \begin{equation}
     \label{eq:termesbeta-re}
    \|\partial^\gamma_x \rho^\alpha \partial^{\beta_1}_x v^\alpha \cdots   \partial^{\beta_\ell}_x v^\alpha \|_{L^2_x}, \qquad \beta_1+ \ldots + \beta_\ell= \beta, \quad |\beta_1|\geq \ldots \geq |\beta_\ell|>0.
     \end{equation}
     Again, by introducing $n$ as previously in \eqref{def:n-bourbakiste}, that is such that for all $j> n$, $|\beta_j| < k-d/2$, then by Hölder inequality  it is sufficient to have
     $$
     \frac{1}{2}\geq \frac{1}{2} - \frac{k-1-|\gamma|}{d} +  \sum_{j=1}^n \left(\frac{1}{2} - \frac{k-|\beta_i|}{d}\right).
     $$
     As $|\gamma|+\sum_{i=1}^n |\beta_i|\leq |\gamma|+|\beta|=k-1$, this holds since $k>d/2$.
We thus obtain
\begin{equation}
\label{eq:ouf}
    \|\partial^\gamma_x \rho^\alpha \partial^{\beta_1}_x v^\alpha \cdots   \partial^{\beta_\ell}_x v^\alpha \|_{L^2_x} 
    {\lesssim}\|  \rho^\alpha \|_{H^{k-1}_x}  \| \DDD v^\alpha \|_{ H^{k-1}_x}^\ell.
\end{equation}
    
    The proof of regularity for the mass $m_0$ is straightforward. That for the higher moments $m_j$ follows from the algebra property and Sobolev embedding $H^{k}_x \hookrightarrow L^\infty_x $ for $k>d/2$ that is used to bound $\Vert \partial_x^\ell (\rho^\alpha \langle v^\alpha \rangle^j) \Vert_{L^2_x}  \lesssim \Vert \rho^\alpha \Vert_{H^{k-1}_x}(1+\Vert v^\alpha \Vert_{H^k_x}^j)$,  and the proof can be performed in the same fashion as before.

\end{proof}
In the context of Theorem \ref{thm-LWPmultiphase-densityvelocity-re}, we also have the following version for the regularity of the moments in velocity, which requires the use of the weighted space $L^1_{\alpha,\Theta}$ for the densities $\rhorho$.
\begin{Lem}\label{prop-reg-momentsf-weights}
 Let $(I, \mu)$ be a set of labels such that $(\rhorho, \ww)$ is an associated multiphasic density/velocity family satisfying
\begin{align*}
    \rhorho \in L^1_{\alpha,\Theta} H^{k-1}_x, \ \ \ww \in \mathcal{H}^{k,\infty},
\end{align*}
for some $k>1+d/2$ and some measurable map $\Theta: I  \mapsto \R^+$ Then for any measurable map $\lambda: I  \mapsto \R^d$ satisfying $\langle \lambda^\alpha \rangle^p \lesssim \Theta(\alpha)$ for some $p \in [1, \infty)$, the conclusions of Lemma \ref{prop-reg-momentsf} hold for the distribution function $f\vcentcolon=\mathcal{I}_\mu(\rhorho, \vv)$ with $\vv=\bm{\lambda}+\ww$.
\end{Lem}
\begin{proof}
    The assertion \eqref{reg-moment-mult-psi} follows directly, as $\ww \in L^\infty_\alpha H^{k}_x$. For the assertion \eqref{reg-moment-mult-vj}, we use the fact that for $j=0,\ldots,\lfloor p\rfloor$ we have the pointwise estimate
    \begin{align*}
        \langle v^\alpha \rangle^j=\langle \lambda^\alpha+w^\alpha \rangle^p \lesssim \langle \lambda^\alpha \rangle^p \langle w^\alpha \rangle^p \lesssim \Theta(\alpha) \langle w^\alpha \rangle^p.
    \end{align*}
    Using Leibniz rule and Faa di Bruno formula, we can proceed as in the proof of  Lemma \ref{prop-reg-momentsf} to obtain the fact that there exists a continuous increasing function $\mathrm{C}_{k,j}$ such that 
    \begin{align*}
        \Vert \langle v^\alpha \rangle^j \rho^\alpha \Vert_{H^{k-1}_x} \lesssim \Theta(\alpha) \mathrm{C}_{k,j}(\Vert w^\alpha \Vert_{H^k_x}) \Vert \rho^\alpha \Vert_{H^{k-1}_x} \lesssim \Theta(\alpha) \mathrm{C}_{k,j}( \Vert\ww \Vert_{\mathcal{H}^{k,\infty}}) \Vert \rho^\alpha \Vert_{H^{k-1}_x},
    \end{align*}
    and the last expression is integrable in $\alpha$, which concludes the proof.
\end{proof}

\medskip

Combining Propositions \ref{prop:existence_coupled_case}--\ref{prop:uniqueness_coupled_case} and Theorem~\ref{thm-LWPmultiphase-densityvelocity}, we can finally state the following local well-posedness result on the original kinetic formulation of the Vlasov equation
\begin{equation}
\label{eq:vlasov-chapter2}
	\left\{
	\begin{gathered}
 	\partial_t f + v \cdot \nabla_x f + \mathrm{div}_v \big( F f \big) = 0, \\
 	f|_{t=0}= f_0.
	\end{gathered}
	\right.
\end{equation}
Assume the force in~\eqref{eq:vlasov-chapter2} can be written as
\begin{equation}
    \label{eq:force-vlasov-gen}
F(\cdot,\cdot,v)=
\mathrm{F}\left(v, \left(\displaystyle \int_{\R^d}\psi_j(w)f(\cdot,\cdot, \mathrm{d}w) \right)_{j \in J}, Y\right)n
\end{equation}
where $(\psi_j(v))_{j \in J}$ is a family smooth functions and $Y\in \mathcal{Y}$ stand for parameters. In other words, we assume that the force field mainly depends on moments in velocity of the distribution function. We shall impose conditions on $\mathrm{F}$ such that Theorem~\ref{thm-LWPmultiphase-densityvelocity} applies; we will check later in the upcoming chapters that various Vlasov--Poisson and  Vlasov--Navier--Stokes type equations satisfy the required assumptions.

\begin{Thm}\label{thm-LWPkinetic}
Let $k>1+d/2$ and $p \in [1, \infty]$. Suppose that the initial distribution function $f_0$ can be written as
\begin{align}\label{eq:compatib-assump-kinetic}
    f_0=\mathcal{I}_\mu(\rhorho_0, \vv_0)
\end{align}
for some set of labels $(I, \mu)$ and some regular multiphasic distribution $(\rhorho_0, \vv_0) \in L^\infty_\alpha (H^{k-1}_x \cap \dot{H}^{-1}_x) \times \mathcal{H}^{k,p}$.
            Assume that the family of forces
$$
\widetilde{\FF}(\vv,\rhorho, Y) \vcentcolon = \mathrm{F}\left(\vv, \left(\displaystyle \int_I \psi_j(v^\alpha) \rho^\alpha  \, \mathrm{d}\mu(\alpha) \right)_{j \in J},Y \right)
$$
satisfies the Assumptions~\ref{ass:bound_high_regularity-rho} and~\ref{ass:stability-rho}. 
Then there exist a time $T>0$ and a unique solution $f$ on $[0,T]$ to the Vlasov equation~\eqref{eq:vlasov-chapter2} with initial data $f_0$ and force field \eqref{eq:force-vlasov-gen} in the sense of Definition \ref{def:solution_coupled}. Furthermore, we have
$$\forall t \in [0,T], \ \ f(t)= \mathcal{I}_\mu( \rhorho(t), \vv(t)),$$
where $(\rhorho, \vv)$ is the unique multiphasic solution from Theorem \ref{thm-LWPmultiphase-densityvelocity} with force field $\widetilde{\FF}(\vv,\rhorho, Y) $. The distribution function $f(t)$ satisfies the following regularity:
\begin{align}\label{eq:weak-continuity-kinetic}
\begin{split}
    &f \in \mathscr{C}^0([0,T]; L^\infty_x (w^*-\mathcal{M}_v)), \\
    &\forall t \in [0,T], \forall \psi \in \mathscr{C}^{k-1}_b(\R^d), \ \ \lim_{s\to t} \left\|  \langle f(t)-f(s), \psi \rangle   \right\|_{H^{k-1}_x} = 0,
    \end{split}
\end{align}
and the following regularity of the moments holds:
\begin{align*}
    (t,x) \mapsto \int_{\R^d} \vert v \vert^j f(t,x,\D v) \in L^\infty([0,T]; H^{k-1}_x), \ \ j=0,\ldots,\lfloor p\rfloor.
\end{align*}
\end{Thm}

\begin{Rem} As before, in the case when the control by $\| \rhorho\|_{L^\infty([0,T]; L^\infty_\alpha \dot{H}^{-1}_x)}$ is not required in the argument of the estimate function $\mathcal F$ in Assumption~\ref{ass:bound_high_regularity-rho}, Theorem~\ref{thm-LWPkinetic} holds for $\rhorho_0 \in L^\infty_\alpha H^{k-1}_x$. 

\end{Rem}

\begin{proof}
The proof consists in considering the multiphasic formulation
	\begin{equation*}
	\left\{
	\begin{gathered}
	\partial_t \rho^\alpha + \Div( \rho^\alpha v^\alpha) =0, \\
	\partial_t v^\alpha + ( v^\alpha \cdot \nabla ) v^\alpha = \widetilde{F}^\alpha,\\
	\rho^\alpha |_{t = 0} = \rho^\alpha_0, \ v^\alpha |_{t = 0} = v^\alpha_0,
	\end{gathered}
	\right.
	\end{equation*}
    and, as announced, to apply Theorem~\ref{thm-LWPmultiphase-densityvelocity} (as Assumptions~\ref{ass:bound_high_regularity-rho} and~\ref{ass:stability-rho} are satisfied) and Propositions \ref{prop:existence_coupled_case}--\ref{prop:uniqueness_coupled_case}. The regularity of moments is a consequence of Proposition~\ref{prop-reg-momentsf}. It therefore remains to prove the continuity in time \eqref{eq:weak-continuity-kinetic}. For the first property, let $\phi \in \mathscr{C}_b(\R^d)$ and $t,s \in [0,T]$. We have
    \begin{align*}
        &\langle f(t,x)-f(s,x), \phi \rangle \\
        &\quad=\int_I \rho^\alpha(t,x) \phi(v^\alpha(t,x))\, \mathrm{d}\mu(\alpha)-\int_I \rho^\alpha(s,x) \phi(v^\alpha(s,x))\, \mathrm{d}\mu(\alpha) \\
        &\quad=\int_I \big(\rho^\alpha(t,x)-\rho^\alpha(s,x) \big) \phi(v^\alpha(t,x))\, \mathrm{d}\mu(\alpha) \\
        &\qquad \qquad \qquad \qquad -\int \rho^\alpha(s,x) \big(\phi(v^\alpha(s,x))-\phi(v^\alpha(t,x)) \big)\, \mathrm{d}\mu(\alpha) \\
        &\quad\vcentcolon=\mathrm{I}_{t,s}(x)+ \mathrm{II}_{t,s}(x).
    \end{align*}
    For the first term, we have by Sobolev embedding that
    \begin{align*}
        \vert \mathrm{I}_{t,s}(x) \vert  \lesssim \Vert \phi \Vert_{L^\infty} \int _I  \Vert \rho^\alpha(t)-\rho^\alpha(s) \Vert_{ H^{k-1}_x} \, \mathrm{d}\mu(\alpha),
    \end{align*}
which tends to zero as $t \rightarrow s$ since $\rhorho \in \mathscr{C}([0,T]; L^1_\alpha H^{k-1}_x)$ by Theorem \ref{thm-LWPmultiphase-densityvelocity}. For the second term, we write
\begin{align*}
    \vert \mathrm{II}_{t,s}(x) \vert \leq  \int_I \vert \rho^\alpha(t,x) \vert \big\vert \phi(v^\alpha(t,x))-\phi(v^\alpha(s,x)) \big\vert \, \mathrm{d}\mu(\alpha).
\end{align*}
As the integrand is bounded by $2 \Vert \phi \Vert_{L^\infty}\vert  \Vert \rho^\alpha \Vert_{L^\infty([0,T];L^{\upinfty}_\alpha L^\infty_x)} $ and since from Remark \ref{rem-continuity-thmgendensityvelocity},  we know that   for $\mu$-almost every $\alpha \in I$, $v^\alpha \in \mathscr{C}^0([0,T]; H^k_x)$, as well as $\rho^\alpha \in \mathscr{C}^0([0,T]; H^{k-1}_x)$, we deduce that $\vert \mathrm{II}_{t,s}(x) \vert$ converges to zero as $t \rightarrow s$ thanks to the dominated convergence theorem (as $\mu(I)<\infty)$.

    Let us now prove the second continuity property in \eqref{eq:weak-continuity-kinetic}, and let $\psi \in \mathscr{C}^{k-1}_b(\R^d)$,  $t,s \in [0,T]$. 
       We apply $\partial^\zeta_x$ to the second line of the decomposition of $\langle f(t,x)-f(s,x), \phi \rangle$ above, for $\zeta \in \N^d$, $|\zeta| \leq  k-1$ and $\ell \leq k-1$ and argue as in the proof of Proposition~\ref{prop-reg-momentsf}.
       According to the Leibniz rule and to the Fa\`a di Bruno formula, it is enough to estimate terms of the form 
       \begin{multline*}
           J_{\gamma,(\beta_i)}\vcentcolon= \\
           \int_I \left( \partial_x^\gamma \rho^\alpha \partial^{\beta_1}_x v^\alpha \cdots   \partial^{\beta_\ell}_x v^\alpha  \phi^{(\ell)}(v^\alpha) (t) - \partial_x^\gamma \rho^\alpha \partial^{\beta_1}_x v^\alpha \cdots   \partial^{\beta_\ell}_x v^\alpha  \phi^{(\ell)}(v^\alpha) (s)   \right) \mathrm{d}\mu(\alpha) 
            \end{multline*}
       in $L^2_x$ where  $ \beta_i,\gamma \in \N^d$ with $\sum_{i=1}^\ell |\beta_i|+|\gamma| \leq k-1$. 
       According to the estimates~\eqref{eq:ouf0}--\eqref{eq:ouf}, we have
       \begin{align*}
       &\|J_{\gamma,(\beta_i)}\|_{L^2_x} \\
       &\lesssim \| \DDD \vv\|_{L^\infty([0,T]; L^{\upinfty}_\alpha H^{k-1}_x)}^{\ell}\int_I  \| \rho^\alpha(t)- \rho^\alpha(s)\|_{H^{k-1}_x}   \mathrm{d}\mu(\alpha) 
       \\ 
       & \ +\| \rho^\alpha\|_{L^\infty([0,T]; L^{\upinfty}_\alpha H^{k-1}_x)}\| \DDD \vv\|_{L^\infty([0,T]; L^{\upinfty}_\alpha H^{k-1}_x)}^{\ell-1}\int_I  \| \DDD v^\alpha(t) - \DDD v^\alpha(s)\|_{H^{k-1}_x} \mathrm{d}\mu(\alpha) \\
       & \ + \| \rho^\alpha\|_{L^\infty([0,T]; L^{\upinfty}_\alpha H^{k-1}_x)} \| \DDD \vv\|_{L^\infty([0,T]; L^{\upinfty}_\alpha H^{k-1}_x)}^{\ell} \\
       &\qquad \qquad \qquad \times \int_I   \|\psi^{(\ell)}(v^\alpha(t)) - \psi^{(\ell)}(v^\alpha(s))\|_{L^{\upinfty}_x} \mathrm{d}\mu(\alpha).
         \end{align*}
To deduce the convergence to zero as $s \rightarrow t$, we apply Theorem~\ref{thm-LWPmultiphase-densityvelocity} and the last item of Corollary~\ref{cor:timecontinuity} for the required time continuity properties of respectively $\vv$ and $\rhorho$, combined with the dominated convergence theorem. Note that to handle the last term, we again use $\mu(I)<+\infty$. This concludes the proof.
\end{proof}

This theorem invites us to put forward a specific notion for time continuity.

\begin{Def}[Continuous in time, strongly in space and weakly in velocity]
		\label{def:solution_Vlasov+timecontinuity}
        Let $m \in \N$.
		A measurable curve $f = f(t)$, $t \in [0,T]$ valued in the space of nonnegative Radon measures of $\Omega \times \R^d$  is said to be continuous in time, {\it $H^m_x$--strongly in space  and weakly in velocity}, if it satisfies the following time continuity property:
        \begin{align}\label{eq:weak-continuity-kinetic-def}
\begin{split}
    &f \in \mathscr{C}^0([0,T]; L^\infty_x (w^*-\mathcal{M}_v)), \\
    &\forall t \in [0,T], \forall \psi \in \mathscr{C}^{m}_b(\R^d), \ \ \lim_{s\to t} \left\|  \langle f(t)-f(s), \psi \rangle   \right\|_{H^{m}_x} = 0,
    \end{split}
\end{align}
	\end{Def}

\medskip
\begin{Ex}\label{eq:example-decompoLinfty}

Let us write down some concrete examples of kinetic initial data for which the initial decomposition  \eqref{eq:compatib-assump-kinetic}  of the initial data $f_0$ of Theorems \ref{thm-LWPkinetic} is satisfied. 
\begin{itemize}
    \item 
The first example is
$$f_0=\sum_{j=1}^n m_j(x) \otimes \delta_{v=v_j(x)},$$ 
with $n \in \N$, for some given sets of  positive masses $(m_j)_{1 \leq j \leq n}$ and velocities $(v_j)_{1 \leq j \leq n}$. Here, one simply takes $$I=[\![1, n]\!], \ \ \mathrm{d}\mu(\alpha)=\frac{1}{n}\sum_{j=1}^n \delta_{\alpha=j},$$
and  
$$\rho^\alpha_0(x)=n m_\alpha(x), \ \ v^\alpha_0(x)=v_\alpha(x),$$ 
so that the decomposition \eqref{eq:compatib-assump-kinetic}  holds. This was of course the prototypical example at the root of our multiphase framework. If $\mathbf{m} \in L^\infty_\alpha H^{s-1}_x$ then we have $\rhorho_0 \in L^\infty_\alpha H^{s-1}_x$.
\item Consider the spatial periodic case $\Omega=\T^d$. Let us assume that  $ f_0 \in L^\infty_{v,\Psi^{-1}}(H^s_x \cap \dot{H}^{-1}_x)$ for some positive $\Psi \in L^1(\R^d)$ having a total mass $1$, and such that $\vert v\vert^p \Psi \in L^1(\R^d)$ with $p \in [1, \infty)$. In light of the proof of Proposition \ref{prop:disintegration} (Case 2.), \eqref{eq:compatib-assump-kinetic} holds for some set of labels $(I, \mu)$. Indeed, set
        $$I=\R^d, \quad \D \mu(\alpha)={\Psi(\alpha) }\mathrm{d}\alpha,$$
and 
$$ \rho_0^\alpha(x)=\frac{f_0(x,\alpha)}{\Psi(\alpha)}, \ \ v^\alpha_0(x)=\alpha.$$
As $x \in \T^d$, we observe that $\vv_0 \in \mathcal{H}^{s+1,p}$. Then  $f_0=\mathcal{I}_\mu(\rhorho_0, \vv_0)$ and since $f_0 \in L^\infty_{v,\Psi^{-1}}(H^s_x \cap \dot{H}^{-1}_x)$ we get that $\rhorho_0 \in L^\infty_\alpha (H^{s}_x \cap \dot H^{-1}_x)$.
\end{itemize}

\end{Ex}

Let us also state the kinetic version of Theorem~\ref{thm-LWPmultiphase-densityvelocity-re}.

\begin{Thm}
\label{thm-LWPkinetic-L1}
Let $k>1+d/2$. Suppose that the initial distribution function $f_0$ can be written as
\begin{align}\label{eq:compatib-assump-kinetic-L1}
    f_0=\mathcal{I}_\mu(\rhorho_0, \ww_0+\boldsymbol{\lambda}),
\end{align}
for some set of labels $(I, \mu)$, some measurable function $\lambda: (I,\mu)\to \R^d$, and some regular multiphasic distribution $(\rhorho_0, \ww_0) \in L^1_{\alpha,\Theta} H^{k-1}_x \cap \dot{H}^{-1}_x \times \mathcal{H}^{k,\infty}$ with a measurable weight $\Theta : I \mapsto \R^+$.
            Assume that the family of forces
$$
\widetilde{\FF}(\ww,\rhorho, Y) \vcentcolon = \mathrm{F}\left(\ww+\boldsymbol{\lambda}, \left(\displaystyle \int_I \psi_j(w^\alpha+\lambda^\alpha) \rho^\alpha  \, \mathrm{d}\mu(\alpha) \right)_{j \in J},Y \right)
$$
satisfies the Assumptions~\ref{ass:bound_high_regularity-rho-re} and~\ref{ass:stability-rho-re}. 
Then there exist a time $T>0$ and a unique solution $f$ on $[0,T]$ to the Vlasov equation~\eqref{eq:vlasov-chapter2} with initial data $f_0$ and force field \eqref{eq:force-vlasov-gen} in the sense of Definition \ref{def:solution_coupled}. Furthermore, we have
$$\forall t \in [0,T], \ \ f(t)= \mathcal{I}_\mu( \rhorho(t), \vv(t)+\boldsymbol{\lambda}),$$
where $(\rhorho, \ww)$ is the unique multiphasic solution from Theorem \ref{thm-LWPmultiphase-densityvelocity-re} with force field $\widetilde{\FF}(\ww,\rhorho, Y)$.

Furthermore, if $1 \lesssim \Theta(\alpha)$ then the distribution function $f(t)$ satisfies the continuity property \eqref{eq:weak-continuity-kinetic} and if $\langle \lambda^\alpha \rangle^p \lesssim \Theta(\alpha)$ for some $p \in [1, \infty)$ then the regularity of moments stated in Theorem~\ref{thm-LWPkinetic} holds.

\end{Thm}
\begin{proof}
   As in the proof of Theorem \ref{thm-LWPkinetic}, the result is a consequence of Theorem \ref{thm-LWPmultiphase-densityvelocity-re}. The only thing to check is the continuity in time and the propagation of regularity for the moments. For the latter, we can  directly appeal to Lemma \ref{prop-reg-momentsf-weights}, since $1 \leq \langle \lambda^\alpha \rangle^p \lesssim \Theta(\alpha)$. Let us detail the continuity in time. We use the same notation as in the proof of  Theorem \ref{thm-LWPkinetic}: as $1 \lesssim \Theta(\alpha)$, we have  
   \begin{align*}
        \vert \mathrm{I}_{t,s}(x) \vert  &\lesssim \Vert \phi \Vert_{L^\infty}\int _I  \Theta(\alpha)\Vert  \rho^\alpha(t)-\rho^\alpha(s) \Vert_{ H^{k-1}_x} \, \mathrm{d}\mu(\alpha), \\
    \vert \mathrm{II}_{t,s}(x) \vert &\lesssim  \int_I \Theta(\alpha)\vert \rho^\alpha(t,x) \vert \big\vert \phi(v^\alpha(t,x))-\phi(v^\alpha(s,x)) \big\vert \, \mathrm{d}\mu(\alpha).
\end{align*}
The term $\mathrm{I}_{t,s}$ tends to zero as $t \mapsto s$ since $\rhorho \in \mathscr{C}^0([0,T]; L^{1}_{\alpha, \Theta} H^{k-1}_x)$. For the term $\mathrm{II}_{t,s}$, we can rely on the dominated convergence theorem, bounding the integrand by $2 \Vert \phi \Vert_{L^\infty}\Theta(\alpha)\Vert \rho^\alpha \Vert_{L^\infty([0,T]; L^\infty_x)}$: from the proof of Lemma \ref{LM:SobEstimRHO} (see \eqref{eq:dracolosse}), we readily see that the former upper bound is integrable in $\alpha$, which gives the first  time continuity property. For the second one, we can treat the term $J_{\gamma,(\beta_i)}$ from the proof of  Theorem \ref{thm-LWPkinetic}, by writing
\begin{align*}
       \|J_{\gamma,(\beta_i)}\|_{L^2_x} 
       &\lesssim \| \DDD \vv\|_{L^\infty([0,T]; L^{\upinfty}_\alpha H^{k-1}_x)}^{\ell}\int_I \Theta(\alpha) \| \rho^\alpha(t)- \rho^\alpha(s)\|_{H^{k-1}_x}   \mathrm{d}\mu(\alpha) 
       \\ 
       &+\| \DDD \vv\|_{L^\infty([0,T]; L^{\upinfty}_\alpha H^{k-1}_x)}^{\ell-1} \\ 
       &\quad \times \int_I \Theta(\alpha)\| \rho^\alpha\|_{L^\infty([0,T]; H^{k-1}_x)} \| \DDD v^\alpha(t) - \DDD v^\alpha(s)\|_{H^{k-1}_x} \mathrm{d}\mu(\alpha) \\
       &+\| \DDD \vv\|_{L^\infty([0,T]; L^{\upinfty}_\alpha H^{k-1}_x)}^{\ell} \\
       &\quad  \times \int_I \Theta(\alpha)\| \rho^\alpha\|_{L^\infty([0,T]; H^{k-1}_x)}   \|\psi^{(\ell)}(v^\alpha(t)) - \psi^{(\ell)}(v^\alpha(s))\|_{L^{\upinfty}_x} \mathrm{d}\mu(\alpha),
         \end{align*}
            and we can conclude by applying the dominated convergence theorem as above.    
  
\end{proof}

\begin{Ex}\label{eq:example-decompoL1}
We can also provide some examples of kinetic initial data for which the decomposition \eqref{eq:compatib-assump-kinetic-L1} of the initial data $f_0$ in Theorem \ref{thm-LWPkinetic-L1} holds.
\begin{itemize}
    \item  Suppose that $f_0 \in L^1_{v,\langle v \rangle^p} (H^s_x \cap \dot{H}^{-1}_x)$ with $s >d/2$ and $p \in [1, \infty)$. Then \eqref{eq:compatib-assump-kinetic-L1} holds for some set of labels $(I, \mu)$ : this time, we set
        $$I=\R^d, \quad \D \mu(\alpha)=\Psi(\alpha) \mathrm{d}\alpha,$$
and 
$$ \rho_0^\alpha(x)=\frac{f_0(x,\alpha)}{\Psi(\alpha)}, \ \ w^\alpha_0(x)=0, \ \ \lambda^\alpha=\alpha,$$
for some positive function $\Psi \in L^1(\R^d)$ having total mass equal to $1$. Then  $f_0=\mathcal{I}_\mu(\rhorho_0, \ww_0+\boldsymbol{\lambda})$  and we have $\ww_0 \in \mathcal{H}^{s+1,\infty}$ as well as  $\rhorho_0 \in L^1_{\alpha,\langle \alpha \rangle^p}(H^{s}_x \cap \dot H^{-1}_x)$, since  $f_0 \in L^1_{v, \langle v \rangle^p} (H^s_x \cap \dot{H}^{-1}_x)$.
\item Suppose now that $\Omega=\R^d$ and $ f_0 \in L^\infty_{v,\Psi^{-1}}(H^s_x \cap \dot{H}^{-1}_x)$ with $s>d/2$ for some positive $\Psi \in L^1(\R^d)$ with total mass $1$. By repeating the same construction for the set of labels and density/velocity families from the previous bullet point, we obtain $(\rhorho_0, \ww_0) \in L^\infty_{\alpha}(H^{s}_x \cap \dot H^{-1}_x)\times \mathcal{H}^{s+1,\infty}$. Then we readily observe that Theorem \ref{thm-LWPkinetic-L1} still holds true with this initial data, replacing $L^{1}_{\alpha, \Theta}$ in the statement by $L^\infty_\alpha$.
    \item The latter constructions can actually be generalized. From the proof of Proposition \ref{prop:disintegration} (Case 1) we can handle more probability measures $f_0 \in \mathscr{P}(\Omega \times \R^d)$, $\Omega =\T^d$ or $\R^d$, satisfying the following assumption. We set $I=\R^d$ and $\mu= \pi_2 {}\pf f$, where $\pi_2$ is the canonical projection from $\Omega \times \R^d$ to $\R^d$.  We ask that 
\begin{equation}
\label{eq:hyp-cond}
    \rho^\alpha_0(w) = f(x | \pi_2 = \alpha) \in L^1_{\alpha, \langle \alpha \rangle^p}(H^s\cap \dot{H}^{-1}), \ \ p \in, [1, \infty],
\end{equation} 
where $f(\cdot | \pi_2 = v)$ is the conditional law  given by the disintegration lemma. 
We also take $w^\alpha_0(x)=0$ and $\lambda^\alpha=\alpha$.   Then $f_0=\mathcal{I}_\mu(\rhorho_0, \ww_0+\boldsymbol{\lambda})$ and we have $\ww_0 \in \mathcal{H}^{s+1,\infty}$ as well as  $\rhorho_0 \in L^1_{\alpha,\langle \alpha \rangle^p}(H^{s}_x \cap \dot H^{-1}_x)$.

However, the condition~\eqref{eq:hyp-cond} may be difficult to check in practice, for general measures.
\end{itemize}
Note that similar considerations hold when replacing $H^s_x$ by $\mathscr{C}^{0,\theta}_x$, the set of $\theta-$H\"olderian continuous functions, for any $\theta \in (0,1]$.
\end{Ex}

\subsection{Application to the monokinetic limit}
\label{sec:mono}
As an application, let us study the problem of the monokinetic limit, which was part of Motivation I from Section \ref{motivation1}. Consider again the Vlasov equation
\begin{equation}
\label{eq:vlasov-mono}
	\left\{
	\begin{gathered}
 	\partial_t f + v \cdot \nabla_x f + \mathrm{div}_v \big( F f \big) = 0, \\
 	f|_{t=0}= f_0,
	\end{gathered}
	\right.
\end{equation}
where the force in~\eqref{eq:vlasov-chapter2} has the form
$$
F(\cdot,\cdot,v)=
\mathrm{F}\left(v, \left(\displaystyle \int_{\R^d}\psi_j(w)f(\cdot,\cdot, \mathrm{d}w) \right)_{j \in J}, Y\right).
$$
This leads to the study of the multiphasic system
\begin{equation}\label{eq:multiphaseGEN-monokin}
\left\{
\begin{aligned}
	\partial_t \rho^\alpha+ \Div( \rho^\alpha v^\alpha) =0, \\
	\partial_t v^\alpha + ( v^\alpha \cdot \nabla ) v^\alpha =\widetilde{\FF}(\vv,\rhorho,Y),\\
	\rho^\alpha |_{t = 0} = \rho^\alpha_{\mathrm{in}}, \ v^\alpha |_{t = 0} = v^\alpha_{\mathrm{in}}.
 \end{aligned}
 \right.
\end{equation}
where the family of forces
$$
\widetilde{\FF}(\vv,\rhorho, Y) \vcentcolon = \mathrm{F}\left(\vv, \left(\displaystyle \int_I \psi_j(v^\alpha) \rho^\alpha  \, \mathrm{d} \mu(\alpha) \right)_{j \in J},Y \right)
$$
is asked to satisfy the Assumptions~\ref{ass:bound_high_regularity-rho} and~\ref{ass:stability-rho}.

\begin{Thm}\label{thm-monokinetic-conv}
    Let $(I,\mu)$ be a set of labels, assumed to be a probability space. Let $k>\max(k_0,1+d/2)$ and $p \in [2, \infty]$. Consider $(\rhorho_{n,\mathrm{in}}, \vv_{n,\mathrm{in}})_{n \in \N}$ a sequence of initial densities and velocities assumed to be uniformly bounded in $L^\infty_\alpha (H^{k-1}_x \cap \dot{H}^{-1}_x) \times \H^{k,p}$, normalized in the sense that $\int_I \int_{\Omega} \rho^\alpha_{n,\mathrm{in}}(x) \, \mathrm{d}x \D \mu(\alpha) =1$ and with $\inf_{\alpha, x ,n} \rho^\alpha_{n,\mathrm{in}}(x) \geq m_{\mathrm{in}}$ for some $m_{\mathrm{in}}>0$.

    Assume that 
    \begin{equation}
        \label{eq:assumption-mono}
\begin{aligned}
\rho_{n, \mathrm{in}}^{\alpha}-\rho_\mathrm{in}^\alpha &\underset{n \rightarrow + \infty}{\longrightarrow} 0, \ \ \text{in}\ \  L^p_\alpha \dot H^{-1}_x, \\
v_{n,\mathrm{in}}^{\alpha} &\underset{n \rightarrow + \infty}{\longrightarrow} v_\mathrm{in}, \ \ \text{in} \ \ L^p_\alpha L^2_x,
\end{aligned}
    \end{equation}
for a family of densities $\rhorho_{\mathrm{in}} \in L^\infty_\alpha (H^{k-1}_x \cap \dot{H}^{-1}_x) $ satisfying $\int_I \int_{\Omega} \rho^\alpha_{\mathrm{in}}(x) \, \mathrm{d}x \D \mu(\alpha) =1$ and a single velocity field $ v_{\mathrm{in}} \in H^{k}_x$.  Let $(\rhorho_n, \vv_n)$ be the unique maximal solution on $[0,T^\star_n)$ to \eqref{eq:multiphaseGEN-monokin} with initial data $(\rhorho_{n,\mathrm{in}}, \vv_{n,\mathrm{in}})$ such that
\begin{align*}
    (\rhorho_n, \vv_n) \in L^\infty([0,T_n^\star]; L^\infty_\alpha H^k_x) \times L^\infty([0,T_n^\star];  \H^{k,p}).
\end{align*}
Then there exists $T^\star>0$ independent of $n$ such that
\begin{align}\label{hyp:unifTIME-monokin}
    T^\star_n > T^\star, \ \ \text{for all } \ n \in \N, 
\end{align}
and the distribution function $f_n(t)\vcentcolon= \mathcal{I}_\mu(\rhorho_n(t), \vv_n(t))$, which is the unique solution  to the Vlasov equation~\eqref{eq:vlasov-mono} associated to the initial  condition $f_\mathrm{in} \vcentcolon= \mathcal{I}_\mu(\rhorho_{n,\mathrm{in}}, \vv_{n,\mathrm{in}})$
converges weakly towards the monokinetic profile 
$$f(t)=\rho(t) \otimes \delta_{v=v(t)}$$ 
on $[0,T^\star]$ as $n \rightarrow + \infty$,
where $(\rho, v)$ is the unique solution of the system 
\begin{equation}
\label{eq:mono-multi}
\left\{
\begin{aligned}
	\partial_t \rho+ \Div( \rho v) =0, \\
	\partial_t v + ( v \cdot \nabla ) v =\widetilde F(v,\rho,Y),\\
	\rho |_{t = 0} = \int_I \rho^\alpha_{\mathrm{in}}\, \D \mu(\alpha), \ v |_{t = 0} = v_{\mathrm{in}},
 \end{aligned}
 \right.
\end{equation}
where the meaning of $\widetilde F(v,\rho,Y)$ is given below. Moreover, $f(t)$ is the unique solution  to the Vlasov equation~\eqref{eq:vlasov-mono} associated with the initial  condition $ \left(\int_I \rho^\alpha_{\mathrm{in}}\, \D \mu(\alpha) \right) \otimes \delta_{v=v_{\mathrm{in}}}$.
\end{Thm}

\begin{Rem}
    The convergence above can actually be made more precise and quantified in terms of the initial assumption \eqref{eq:assumption-mono} -- see the proof below. In particular, we rely on the classical $2$-Wasserstein distance for probability measures. 
\end{Rem}

We note that since $(\rho_n^\alpha, v_n^\alpha)$ is a  solution to the multiphasic system on $[0,T^\star]$, we know that 
\begin{align*}
    f_{n}(t)=\mathcal{I}_\mu(\rhorho_n(t), \vv_n(t))=\int_I \rho^\alpha_n(t) \otimes  \delta_{v=v^\alpha_n(t)} \, \mathrm{d}\mu(\alpha), \ \ 0 \leq t \leq T^\star,
\end{align*}
is a distributional solution to the Vlasov equation on $[0,T^\star]$ (see Section \ref{sec:abstract_existence}). 

We actually have to be a bit more precise about the meaning of \eqref{eq:mono-multi}, especially about the involved force field $F(v,\rho,Y)$. Let us first consider $(\rho^\alpha,v)_{\alpha \in I}$ the unique solution to
\begin{equation}
\label{eq:mono-multi-intermediate}
\left\{
\begin{aligned}
	\partial_t \rho^\alpha+ \Div( \rho^\alpha v) =0, \\
	\partial_t v + ( v \cdot \nabla ) v =\widetilde \FF(v,\rhorho,Y),\\
	\rho^\alpha |_{t = 0} = \rho^\alpha_{\mathrm{in}}, \ v |_{t = 0} = v_{\mathrm{in}},
 \end{aligned}
 \right.
\end{equation}
where
$$\widetilde \FF(v,\rhorho,Y)=\mathrm{F}\left(v, \left(\psi_j(v)\displaystyle \int_I  \rho^\alpha  \, \mathrm{d} \mu(\alpha) \right)_{j \in J},Y \right).$$
We then define 
\begin{align*}
    \rho\vcentcolon=\int_I  \rho^\alpha  \, \mathrm{d} \mu(\alpha).
\end{align*}
Since $v$ does not depend on the label $\alpha$, we can average the first equation from \eqref{eq:mono-multi-intermediate} in $\alpha$ to obtain the final system \eqref{eq:mono-multi} satisfied by $(\rho,v)$, where the meaning of $F(v,\rho,Y)$ is now clarified as
\begin{align*}
    F(v,\rho,Y)=\mathrm{F}\left(v, \left(\psi_j(v)\displaystyle \int_I  \rho^\alpha  \, \mathrm{d} \mu(\alpha) \right)_{j \in J},Y \right)=\mathrm{F}\left(v, \left(\rho\psi_j(v) \right)_{j \in J},Y \right).
\end{align*}

\begin{Rem}
    We mention that following our framework, \cite{LemarieMultiphaseVNS} obtains a similar monokinetic convergence for the incompressible Vlasov--Navier--Stokes system posed on the whole space;  similar assumptions are put forward.
\end{Rem}

\begin{proof}[Proof of Theorem \ref{thm-monokinetic-conv}]  

For all $n$, the maximal solution $(\rhorho_n, \vv_n)$ to the system \eqref{eq:multiphaseGEN-monokin} exists according to Theorem~\ref{thm-LWPmultiphase-densityvelocity}. Moreover, since by assumption there is $ M_{\mathrm{in}}>0$ such that
\begin{align*}
\sup_n \Vert \rhorho_{n, \mathrm{in}} \Vert_{L^{\upinfty}_\alpha H^{k-1}_x} + \Vert \vv_{n, \mathrm{in}} \Vert_{k,p} \leq M_{\mathrm{in}},
\end{align*}
the estimates of Theorem~\ref{thm-LWPmultiphase-densityvelocity} also become uniform with respect to $n$ and there is $T^\star>0$ independent of $n$ such that $T^\star_n > T^\star$ for all $n$ (see Remark~\ref{Rem:timeuniform}).

Next we want to show that $f_{n}(t) \rightharpoonup f(t)=\rho(t) \otimes \delta_{v=v(t)}$ on $[0,T^\star]$ when $n \rightarrow + \infty$, where $(\rho,v)$ solves \eqref{eq:mono-multi}.
As explained earlier, there holds
\begin{align*}
    f(t)=\rho(t) \otimes \delta_{v=v(t)}=\left(\int \rho^\alpha(t)  \, \mathrm{d}\mu(\alpha) \right)\otimes \delta_{v=v(t)},
\end{align*}
where $(\rho^\alpha,v)_{\alpha \in I}$ is solution to \eqref{eq:mono-multi-intermediate}. As such, since
\begin{align*}
    f_n(t)=\int \rho^\alpha_n(t) \otimes \delta_{v=v^\alpha_n(t)} \, \mathrm{d}\mu(\alpha),
\end{align*} 
it becomes readily transparent that $(\rho^\alpha_n,v^\alpha_n)_{\alpha \in I}$ and $(\rho^\alpha, v)_{\alpha \in I}$ solve the \textit{same} system \eqref{eq:multiphaseGEN-monokin}. This paves the way for an approach based on the stability estimates for the multiphasic system we have developed earlier in this manuscript. Note that up to reducing $T^\star$ if necessary, we can assume that these two solutions are defined on the same interval of time $[0,T^\star]$.

We shall actually prove the more precise estimate
\begin{align*}
\underset{t \in [0,T^\star]}{\sup}  \mathrm{W}_2\left(f_{n}(t) , \int \rho^\alpha(t) \otimes\delta_{v=v(t)} \, \mathrm{d}\mu(\alpha)\right)  \overset{n \rightarrow + \infty}{\longrightarrow} 0,
\end{align*}
where $\mathrm{W}_2$ refers to the standard $2$-Wasserstein distance on $\Omega \times \R^d$: for $\mu,\nu$ two probability measures,
\index{W2@$\mathrm{W}_2$: $2$-Wasserstein distance}
$$
\mathrm{W}_2(\mu,\nu) = \inf_{\gamma\in \Pi(\mu,\nu)} \left(\int_{(\Omega \times \R^d)^2} d(z,z')^2 \, \D \gamma(z,z')\right)^{1/2},
$$
where $\Pi(\mu,\nu)$ is the set of probability measures on $(\Omega \times \R^d)^2$ with first marginal $\mu$ and second marginal $\nu$, 
see e.g. \cite{Santambrogio} for a comprehensive introduction.

Equipped with uniform Sobolev estimates, and in view of the assumption at the initial time, we apply the stability estimates~\eqref{eq:stability-vv} and \eqref{eq:rho_H-1} which yield
\begin{equation}
    \label{eq:mono-conv}
\begin{aligned}
&\underset{[0,T^\star]}{\sup} \Vert v_n^\alpha - v\Vert_{L^p_\alpha L^2_x} \underset{n \rightarrow + \infty}{\longrightarrow} 0, \\
&\underset{[0,T^\star]}{\sup} \Vert \rho_n^\alpha - \rho^\alpha \Vert_{L^{\upinfty}_\alpha \dot{H}^{-1}_x} \underset{n \rightarrow + \infty}{\longrightarrow} 0.
\end{aligned}
\end{equation}
We have
\begin{multline*}
\mathrm{W}_2(f_{n}(t) , \rho(t)\otimes \delta_{v=v(t)})\\
=\mathrm{W}_2 \left( \int_{I} \rho^\alpha_n(t) \otimes\delta_{v=v^\alpha_n(t)} \ \mathrm{d}\mu(\alpha),\int \rho^\alpha(t)\otimes \delta_{v=v(t)} \, \mathrm{d}\mu(\alpha)\right).
\end{multline*}
Note that this makes sense since all the arguments are probability measures on $\Omega \times \R^d$. Indeed, they either have total conserved  mass $\int_{I} \int_{\Omega} \rho^{\alpha}_n(t,x)  \D x  \D \mu(\alpha)=\int_{I} \int_{\Omega} \rho^{\alpha}_{n, \mathrm{in}}(x)  \D x  \D \mu(\alpha)=1 $ and $\int_{I} \int_{\Omega} \rho^{\alpha}(t,x) \D x  \D \mu(\alpha)= \int_{I} \int_{\Omega} \rho^{\alpha}_{\mathrm{in}}(x)  \D x  \D \mu(\alpha)=1$ by assumption.
Therefore, by the triangle inequality,
\begin{align*}
&\mathrm{W}_2\left(f_{n}(t) , \int \rho^\alpha(t) \otimes \delta_{v=v(t)} \, \mathrm{d}\mu(\alpha)\right) \\
&\leq \mathrm{W}_2 \left( \int_{I} \rho^\alpha(t) \otimes\delta_{v=v^\alpha_n(t)} \ \mathrm{d}\mu(\alpha),\int_{I} \rho^\alpha(t) \otimes\delta_{v=v(t)} \ \mathrm{d}\mu(\alpha)\right) \\
 &\quad +  \mathrm{W}_2 \left( \int_{I} \rho^\alpha_n(t)\otimes \delta_{v=v^\alpha_n(t)} \ \mathrm{d}\mu(\alpha),\int_{I} \rho^\alpha(t) \otimes\delta_{v=v^\alpha_n(t)} \, \mathrm{d}\mu(\alpha) \right)
\\
&\vcentcolon = \mathrm{I_n}(t)+\mathrm{II_n}(t).
\end{align*}

To estimate $\mathrm{I_n}$, we use the definition of the $2$-Wasserstein distance on $\Omega \times \R^d$ above, by choosing the admissible coupling
\begin{align*}
    \gamma_t=\int_I \int_{\Omega} \rho^\alpha(t,x) \delta_{((x, v_n^\alpha(t,x)),(x,v(t,x)))} \, \mathrm{d}x \, \mathrm{d}\mu(\alpha).
\end{align*}
It first marginal and second marginal are respectively the first and second argument in $\mathrm{I_n}(t)$. As a consequence, we get
\begin{align*}
\mathrm{I_n}(t)& \leq \left(\int_{(\Omega \times \R^d)^2} \vert (y,w)-(y',w') \vert^2\, \mathrm{d}\gamma_t((y,w), (y',w')) \right)^{1/2}\\
&= \left(\int_{I} \int_{\Omega} \rho^\alpha(t) | v^\alpha_n(t) - v(t)|^2   \D x  \D \mu(\alpha)\right)^{1/2}
 \\
 &\lesssim\Vert \rhorho(t) \Vert_{L^\infty_\alpha L^\infty_x} \| \vv_n(t) - v(t)\|_{L^p_\alpha L^2_x},
\end{align*}
since $p \geq 2$ and $(I, \mu)$ is of finite mass.

For $\mathrm{II_n}$, we construct a suitable coupling as follows: consider the optimal coupling $\widetilde{\gamma}^\alpha_t \in \Pi(\rho^\alpha_n(t) \, \mathrm{d}x,\rho^\alpha(t) \, \mathrm{d}x)$ in the $2$-Wasserstein distance sense, and define
$$\Gamma_t=\int_{I} (G^\alpha_n) \pf \widetilde{\gamma}_t^\alpha \, \mathrm{d}\mu(\alpha), \ \ G^\alpha_n(x,y)=((x,v^\alpha_n(x)),(y, v^\alpha_n(y))).$$
One can check that the marginals of $\Gamma_t$ are exactly the arguments in $\mathrm{II}_n(t)$. Hence, we obtain
\begin{align*}
    \mathrm{II_n}(t)^2 &\leq \int_{(\Omega \times \R^d)^2} \vert (y,w)-(y',w')\vert^2\, \mathrm{d}\Gamma_t ((y,w),(y',w')) \\
    &=\int_{I} \int_{\Omega \times \Omega} (\vert x-y \vert^2 + \vert v^\alpha_n(t,x)-v^\alpha_n(t,y) \vert^2)\, \mathrm{d}\widetilde{\gamma}_t^\alpha(x,y) \, \mathrm{d}\mu(\alpha) \\
    & \lesssim (1+\Vert \nabla \vv_n(t) \Vert_{L^\infty_\alpha L^\infty_x}^2) \int_{I} \int_{\Omega \times \Omega} \vert x-y \vert^2 \, \mathrm{d}\widetilde{\gamma}_t^\alpha(x,y) \, \mathrm{d}\mu(\alpha),
\end{align*}
and so
\begin{align*}
   \mathrm{II_n}(t) \lesssim  (1+\Vert \nabla \vv_n(t) \Vert_{L^\infty_\alpha L^\infty_x})\left(\int_{I} \mathrm{W}_2(\rho^\alpha_n(t) \, \mathrm{d}x, \rho^\alpha(t) \, \mathrm{d}x)^2 \, \mathrm{d}\mu(\alpha) \right)^{1/2}. 
\end{align*}
 Thanks to \cite[Lemma 5.34]{Santambrogio}, which relates $\mathrm{W}_2$ metric to the $\dot{H}^{-1}_x$ norm for absolutely continuous measures with uniformly bounded density, we have
 \begin{align*}
     \mathrm{W}_2(\rho^\alpha_n(t) \, \mathrm{d}x, \rho^\alpha(t) \, \mathrm{d}x)^2 \lesssim m_t^{-1} \left\|   \rho^\alpha_n(t) - \rho^\alpha(t) \right\|_{\dot{H}^{-1}_x}^2 
 \end{align*}
 provided that we have a uniform lower bound $m_t>0$ on the densities $\rho^\alpha_n(t)$ and $\rho^\alpha(t)$ on $[0,T^\star]$. This is obtained by Lemma \ref{LM:pointwiseEstimRHO} that entails
 \begin{align*}
     \underset{x}{\inf} \, \rho^\alpha_n(t,x) &\geq \left( \underset{x}{\inf} \,  \rho^\alpha_{\mathrm{in}}(x) \right) \exp\left(-\int_0^t \Vert \nabla \vv_n(\tau) \Vert_{L^{\upinfty}_\alpha L^{\upinfty}_x} \mathrm{d} \tau\right), \\
     \underset{x}{\inf} \, \rho(t,x) &\geq \left( \underset{x}{\inf} \,  \rho_{\mathrm{in}}(x) \right) \exp\left(-\int_0^t \Vert \nabla v(\tau) \Vert_{L^{\upinfty}_x} \mathrm{d} \tau\right),
 \end{align*}
 so that, recalling the uniform lower bound  $m_{\mathrm{in}}$ on the initial densities, we can take $$m_t\vcentcolon=m_{\mathrm{in}} \min\left(e^{-\int_0^t \Vert \nabla \vv _n(\tau) \Vert_{L^{\upinfty}_\alpha L^{\upinfty}_x} \mathrm{d} \tau},e^{-\int_0^t \Vert \nabla v(\tau) \Vert_{L^{\upinfty}_x} \mathrm{d} \tau}\right).$$
 We then obtain
 \begin{align*}
     \mathrm{II_n}(t) &\lesssim  (1+\Vert \nabla \vv_n(t) \Vert_{L^\infty_\alpha L^\infty_x})m_t^{-1/2}\left(\int_{I}  \left\|   \rho^\alpha_n(t) - \rho^\alpha(t) \right\|_{\dot{H}^{-1}_x}^2 \, \mathrm{d}\mu(\alpha) \right)^{1/2} \\
     & \lesssim (1+\Vert \nabla \vv_n(t) \Vert_{L^\infty_\alpha L^\infty_x})m_t^{-1/2} \Vert \rhorho_n(t)-\rho(t) \Vert_{L^{p}_\alpha \dot{H}^{-1}_x}.
 \end{align*}
since $p \geq 2$. Thanks to the uniform estimates on the solutions on $[0,T^\star]$, we thus deduce thanks to~\eqref{eq:mono-conv} that $$\underset{[0,T^\star]}{\sup}\mathrm{I_n} + \underset{[0,T^\star]}{\sup} \mathrm{II_n} \rightarrow 0, \ \ n \rightarrow +\infty.$$
 This concludes the proof.
\end{proof}

\begin{Rem}

\label{rem:hyp-mono}
The key assumption~\eqref{eq:assumption-mono}, which ensures the monokinetic behaviour, means in practice that the sequence of initial conditions must be ``well-prepared'' . Typically, on the torus $\Omega=\T^d$, it holds for initial data of the form
$$
f_{n,\mathrm{in}} (x,v) = \rho_{n,\mathrm{in}} (x)  \varphi_n\left(v-v_{n,\mathrm{in}}(x)\right),
$$
where $\varphi_n(v) = n^d \varphi(n v)$, with $\varphi\geq 0, \int \varphi \D v =1$, is an approximation of unity. We furthermore assume the higher integrability condition
$$
\int \langle v\rangle^p \varphi \D v <+\infty,
$$
for some $p \in [2,\infty)$.
Then we can take in this case 
$$(I,\mu)= (\R^d, \varphi (\alpha)  \mathrm{d} \alpha),$$ 
and  
$$
\rho^\alpha_{n,\mathrm{in}} (x)= \rho_{n,\mathrm{in}} (x), \quad v^\alpha_{n,\mathrm{in}} (x) =v_{n,\mathrm{in}}(x) + n^{-1} \alpha.
$$
The condition~\eqref{eq:assumption-mono} is then enforced if $\rho_{n, \mathrm{in}} \rightarrow \rho_{\mathrm{in}}$ in $\dot H^{-1}_x$ and $v_{n, \mathrm{in}} \rightarrow v_{\mathrm{in}}$ in $L^2_x$.

On the other hand, it is unclear if the multiphasic approach is relevant to treat general monokinetic convergence, for instance when the sequence of initial data converges to a Dirac mass in the sense of a Wasserstein distance. In that case, better results may be achieved, but they depend on the precise properties of the model at hand (whereas Theorem~\ref{thm-monokinetic-conv} applies to a large class of equations): for instance for the Vlasov--Poisson equation, it seems more relevant to rely on Wasserstein stability estimates \cite{loeper2006uniqueness,holding2018uniqueness,Iacobelli} to study the monokinetic limit.

\end{Rem}

\subsection{Comparison with standard well-posedness results in  Sobolev spaces}

In this section, $\Omega = \T^d$. As already alluded to in the introduction, under the same assumptions on the mapping $F$ as in Theorem~\ref{thm-LWPkinetic} (in particular $F$ has the form~\eqref{eq:force-vlasov-gen}), the Vlasov equation 
\begin{equation}
\label{eq:vlasov-WP}
	\left\{
	\begin{gathered}
 	\partial_t f + v \cdot \nabla_x f + \mathrm{div}_v \big( F f \big) = 0, \\
 	f|_{t=0}= f_0,
	\end{gathered}
	\right.
\end{equation}
is locally well-posed in Sobolev spaces with weights in velocity, with large enough indices.
We begin by recalling the following definition of the weighted (in velocity) Sobolev spaces.
    \begin{Def}
\label{def:weighted-v-sobolev}
    Let $n \in \N$ and $r\in \R$. The weighted Sobolev space $ \mathrm{H}^{n}_{r}$ is
    \begin{multline*}
    \mathrm{H}^{n}_{r} \vcentcolon= \Bigg\{f \in L^2(\T^d \times \R^d), \, \\
\| f \|_{\mathrm{H}^{n}_{r}}^2 \vcentcolon= \sum_{|\alpha| + |\beta| \leq n} \int_{\T^d} \int_{\R^d} (1+ |v|^2)^{r} |\partial^\alpha_x \partial^\beta_v f|^2 \, \D v \D x < +\infty\Bigg\}.
    \end{multline*}
    \index{H@$\mathrm{H}^k_r$: Sobolev space in $(x,v)$ with fixed polynomial weight in $v$}
\end{Def}
More precisely, we have the following standard result (for instance proved in \cite[Prop 3.2]{HKapde}):

\begin{Prop}
\label{prop:Hnr}
Assume the force $F(\cdot,f)$ satisfies the same assumptions as in Theorem~\ref{eq:compatib-assump-kinetic}.
Then for $n>1+d/2$ and $r$ large enough, the Vlasov equation~\eqref{eq:vlasov-WP} is locally well-posed in $ \mathrm{H}^{n}_{r}$.
\end{Prop}

\begin{Rem}
    For the Vlasov--Poisson equation, the constraint $n>1+d/2$ has been recently lowered, using averaging lemmas, in \cite{JT}.
\end{Rem}

Theorem~\ref{thm-LWPkinetic} can thus be interpreted as a variant of Proposition~\ref{prop:Hnr}, for initial data with no regularity with respect to velocity and without weights in velocity -- of course only regularity in $x$ is propagated in that case. Both improvements are achieved thanks to the multiphasic approach.

Conversely, it turns out that Proposition~\ref{prop:Hnr} can be somehow recovered from the multiphasic formulation, in the sense that regularity in $v$ may also be propagated.
Assume for simplicity that there is not dependence in $v$ for the force in \eqref{eq:vlasov-WP}. 
This means that in the corresponding multiphasic equations
	\begin{equation}
 	\label{eq:multiphasic-WP}
	\left\{
	\begin{gathered}
	\partial_t \rho^\alpha + \Div( \rho^\alpha v^\alpha) =0, \\
	\partial_t v^\alpha + ( v^\alpha \cdot \nabla ) v^\alpha = F,\\
	\rho^\alpha |_{t = 0} = \rho^\alpha_0, \ v^\alpha |_{t = 0} = v^\alpha_0,
	\end{gathered}
	\right.
	\end{equation}
the forcing $F$ does not depend on $\alpha$.
Though the multiphasic formulation is \textit{a priori} not made for dealing with regularity in $\alpha$, the regularity in $x$ and $v$ of the initial condition $f_0$ can be propagated for the solution $f(t)$ of the Vlasov equation.
Following (the proof of) Proposition~\ref{prop:disintegration}, let us consider the case where $$I=\R^d, \quad \D \mu(\alpha)=\Psi(\alpha) \mathrm{d}\alpha,$$
for a smooth probability density $\Psi \in \mathcal{P}_p(\R^d)$ (for some $p>1$), and 
$$v^\alpha_0(x)=\alpha, \qquad \rho^\alpha_0(x)= \frac{f_0(x, \alpha)}{\Psi(\alpha)}.$$
Finally, suppose that there exists a well defined mapping $(\rhorho, \vv) \mapsto \widetilde{\FF}(\rhorho,\vv,Y)=\vcentcolon F$ satisfying the Assumptions~\ref{ass:bound_high_regularity-rho} and~\ref{ass:stability-rho}.
\begin{Prop}
\label{prop:reg-alpha}
    Assume $\rhorho_0 \in L^\infty_\alpha (H^{k-1}_x\cap \dot{H}^{-1}_x)$. Then there exist $T>0$ and a unique solution $ \rhorho \in L^\infty([0,T]; L^\infty_\alpha (H^{k-1}_x \cap \dot{H}^{-1}_x)), \vv \in L^\infty([0,T]; \mathcal{H}^{k,p})$ to~\eqref{eq:multiphasic-WP}, such that the map $\alpha \mapsto v^\alpha(t,x)$ is a $\mathscr{C}^1$ diffeomorphism, uniformly in $(t,x)$, with inverse $v \mapsto w^v (t,x)$. Moreover, the unique solution to the Vlasov equation~\eqref{eq:vlasov-WP}
    satisfies
    \begin{equation}\label{eq:formule-f}
f(t,x,v) = \rho^{w^v(t,x)}(t,x) \Psi(w^v(t,x))  |\det \DDD_v w^v(t,x)|.
\end{equation}
\end{Prop}

\begin{proof}
We apply Theorem~\ref{thm-LWPmultiphase-densityvelocity} to get a solution
\begin{align}
\label{eq:estimtrucdefou}
    \rhorho \in L^\infty([0,T]; L^\infty_\alpha (H^{k-1}_x \cap \dot{H}^{-1}_x)), \ \ \vv \in L^\infty([0,T]; \mathcal{H}^{k,p}),
\end{align}
to the multiphasic equations~\eqref{eq:multiphasic-WP}.
We can then prove that for short times, for all $x \in \Omega$,  $\alpha \mapsto v^\alpha(t,x)$ is a small Lipschitz perturbation of the identity, and thus a $\mathscr{C}^1$ diffeomorphism, uniformly in $(t,x)$. Indeed $\DDD_\alpha v^\alpha$ satisfies the following forced transport equation
	\begin{equation}\label{eq:defDw}
	\left\{
	\begin{gathered}
	\partial_t  \DDD_\alpha v^\alpha +( v^\alpha \cdot \nabla)  \DDD_\alpha v^\alpha = - \DDD_\alpha v^\alpha : \nabla   v^\alpha ,\\
	  \DDD_\alpha v^\alpha|_{t=0} = \mathrm{I},
	\end{gathered}
	\right.
	\end{equation}
    since $v^\alpha_{t=0}=\alpha$. Owing to product estimates, using~\eqref{eq:estimtrucdefou}, there is $\widetilde{T}>0$ small enough such that
$$
\|\DDD_\alpha \vv -\mathrm{I} \|_{L^{\upinfty}([0,\widetilde{T}]; L^{\upinfty}_\alpha H^{k-1}_x)} \leq 1/2.
$$
Therefore the map $\alpha \mapsto v^\alpha(t,x)$ is a small Lipschitz perturbation of the identity and is thus a $\mathscr{C}^1$ diffeomorphism. Let us denote by $v \mapsto w^v(t,x)$ its inverse.  Now, for $t \in [0,\widetilde{T}]$, by Theorem~\ref{thm-LWPkinetic}, we know that
$$
f(t)= \mathcal{I}_\mu( \rhorho(t), \vv(t))
$$
is the unique solution to the Vlasov~\eqref{eq:vlasov-WP} and
for all test functions $\varphi(x,v)$, by the change of variables (in $\alpha$) $\beta\vcentcolon = w^\alpha(t,x)$,
we get 
\begin{align*}
\int f(t,x,v)  &\phi(x,v) \, \mathrm{d}v \, \mathrm{d}x \\
&=  \int  \int \rho^\alpha(t,x)  \varphi(x,v^\alpha(t,x)) \Psi(\alpha) \, \mathrm{d} \alpha  \, \mathrm{d} x \\
&=  \int  \int \rho^{w^v(t,x)}(t,x)   |\det \DDD_v w^v(t,x)|   \Psi(w^v(t,x))  \varphi(x,v)\, \mathrm{d} v \, \mathrm{d} x,
\end{align*}
and we thus deduce that for $t \in [0,\widetilde{T}]$,
\begin{equation*}
f(t,x,v) = \rho^{w^v(t,x)}(t,x) \Psi(w^v(t,x))  |\det \DDD_v w^v(t,x)|,
\end{equation*}
hence the result.
\end{proof}
Though we shall not dwell on this aspect, analogues of this result in other settings, in particular that of Assumptions~\ref{ass:bound_high_regularity-rho-re} and~\ref{ass:stability-rho-re}.
\begin{Rem}
Higher regularity estimates in $\alpha$ for $\DDD_\alpha v^\alpha$ can be obtained, as for  $j=2, \ldots, k$, $\DDD_\alpha^j v^\alpha$ satisfies
\begin{equation}\label{eq:defDw-higher}
	\left\{
	\begin{aligned}
	\partial_t  \DDD_\alpha^j v^\alpha +( v^\alpha \cdot \nabla)  \DDD_\alpha^j v^\alpha &=- ( \DDD_\alpha^j v^\alpha \cdot \nabla)v^\alpha \\
    & \quad + \eta_j ( (\DDD_\alpha^i v^\alpha,\DDD_\alpha^{i'} \nabla v^\alpha)_{i,i'\leq j-1, \, i+i'=j}),\\
	  \DDD_\alpha^j v^\alpha|_{t=0} &= 0,
	\end{aligned}
	\right.
	\end{equation}
where $\eta_j$ is a certain bilinear function. Therefore, using again product estimates and Theorem~\ref{thm-LWPmultiphase-densityvelocity}, we obtain by induction that for all $j$, there is $\widetilde{T}_j>0$ such that 
$$ 
\left\|\DDD_\alpha^j \vv\right\|_{L^\infty([0,\widetilde{T}_j]; L^{\upinfty}_\alpha H^{k-j}_x)} \leq 1/2.
$$
We note that this implies that 
\begin{equation}
    \label{eq:wbeta}
    \|\DDD_\beta^j w^v \|_{L^\infty(0,\widetilde{T}_j; L^\infty_v H^{k-j}_x)} \lesssim 1.
\end{equation}
This implies higher regularity in $\alpha$ for $\rho^\alpha$ as well, as similarly, for $j\leq k$, $\DDD_\alpha^j \rho^\alpha$ satisfies
\begin{equation}\label{eq:defDw-higherhigher}
	\left\{
	\begin{aligned}
	\partial_t  \DDD_\alpha^j \rho^\alpha  +( v^\alpha \cdot \nabla)  \DDD_\alpha^j \rho^\alpha &=- \Div v^\alpha  \DDD_\alpha^j \rho^\alpha \\
    & \quad + \eta'_j \left( \left(\DDD_\alpha^i   v^\alpha, \DDD_\alpha^{i'}   \nabla \rho^\alpha\right)_{i,i'\leq j-1, i+i'=j}\right) \\
 &\quad + \eta''_j \left( \left(\DDD_\alpha^i \nabla  v^\alpha, \DDD_\alpha^{i'}   \rho^\alpha\right)_{i,i'\leq j-1, i+i'=j}\right), \\
	  \DDD_\alpha^j \rho^\alpha|_{t=0} &= \DDD_\alpha^j \rho_0^\alpha(x),
	\end{aligned}
	\right.
	\end{equation}
where $\eta'_j,\eta''_j$ are certain bilinear functions. Here, the regularity of $\DDD_\alpha^j \rho^\alpha$ in $\alpha$ is limited by its initial regularity, that is to say by the regularity in $x$ and $v$ of the initial distribution function 
$f_0(x,v)$. Namely, if we assume for instance that $f_0(x,v)/ \Psi(v) \in W^{k-1,\infty}_v H^{k-1}_x$, we obtain that for $j=1,\ldots, k-1$,
\begin{equation}
    \label{eq:rhoalpha}
\|\DDD_\alpha^j \rho^\alpha \|_{L^\infty(0,\widetilde{T}_j; L^{\upinfty}_\alpha H^{k-j}_x)} \lesssim 1.
\end{equation}
According to~\eqref{eq:formule-f}, regularity in $x$ and $v$ of $f(t)$ can be inferred from that of $\rho^\alpha(t,x)$ in~\eqref{eq:rhoalpha} and of $w^v(t,x)$ in~\eqref{eq:wbeta}. We shall not dwell on the precise regularity, as details can get tedious.
\end{Rem}

\subsection{Comparison with results of propagation of space regularity}

In \cite{HKapde}, the question of propagation of initial higher regularity in $x$ is studied; it is proved that such regularity can propagated into higher regularity for the moments in velocity of the solution (at least locally in time), for a class of nonlinear Vlasov equations set on $\T^d \times \R^d$.

Though the main matter of \cite{HKapde} focuses on a wider class of Vlasov equations than the one considered in the present work (including for instance Vlasov--Maxwell), it deals in particular with a special class, referred to as the  \emph{transport/elliptic case},  which is very close to the class of force fields we are able to consider in this work. The main example of the transport/elliptic case is the Vlasov--Poisson system, on which we focus to fix ideas.

The work \cite{HKapde} requires the use of weighted (in $v$) Sobolev spaces $\mathrm{H}^{n}_{r}$, that were defined in Definition~\ref{def:weighted-v-sobolev}.

\begin{Thm}[Theorem 9.1 from  \cite{HKapde}]
\label{thm:hkapde}
Let $n\in \N$ and $r>0$ be large enough.
Let $n'>n$ be an integer such that
\begin{equation}
\label{eq:constraint-n-n'}
n> \left\lfloor \frac{n'}{2} \right\rfloor+{d} + 1.
\end{equation}
Assume that $f_0 \in \mathrm{H}^{n}_{r}$. 
Assume furthermore that the initial condition $f_0$ satisfies the following higher space regularity:
\begin{equation}
\label{eq-thm2}
\begin{aligned}
\partial^\alpha_x f_0 \in \mathrm{H}^{0}_{r}, \qquad\forall |\alpha|= n'.
\end{aligned}
\end{equation}
Then there is $T>0$ such that the following holds. There exists a unique solution $f(t)$ with initial data $f_0$ to the Vlasov--Poisson system, such that $f(t) \in \mathscr{C}([0,T]; \mathrm{H}^{n}_r)$. 

Moreover, for all test functions $\psi \in L^\infty([0,T]; {W}^{n',\infty})$, we have
\begin{equation}
\int f \psi \,\D v  \in L^2([0,T]; H^{n'}_x).
\end{equation}

\end{Thm}

We therefore note that the multiphasic approach developed in this work, see Theorem~\ref{thm-LWPkinetic}, improves the result of Theorem~\ref{thm:hkapde} on several aspects.
First of all, neither weight, nor regularity in velocity is required to propagate regularity for moments. Furthermore, any amount of regularity in space can be propagated on moments, which reveals that
the constraint~\eqref{eq:constraint-n-n'} is not necessary (and thus shown to be only a technical restriction in Theorem~\ref{thm:hkapde}).

\subsection{From moment regularity to anisotropic Sobolev regularity}

 To conclude this section, assume to fix ideas that $\Omega= \T^d$. As a matter of fact, the regularity for moments \eqref{eq:weak-continuity-kinetic} obtained in Theorem~\ref{thm-LWPkinetic} yields regularity for the solution itself (see \cite{Gerard} for a microlocal version of this fact in the context of averaging lemmas and for a very similar statement in \cite[Corollary 2.5]{HKapde}) in anisotropic Sobolev spaces (as defined by H\"ormander in \cite[Chapter II, Section 2.5]{Horm}), that we first introduce.

\begin{Def}Let $m, n \in \R$. The anisotropic Sobolev space $H^{m,n}_{x,v}$ is defined as 
$$
H^{m,n}_{x,v}: = \Bigg\{ g \in \mathscr{S}'(\T^d \times \R^d), \, (1+|\xi|^2)^{m/2} (1+|\eta|^2)^{n/2} \widehat{g}(\xi,\eta) \in L^2( \Z^d \times \R^d) \Bigg\},
$$
where $\widehat{g}$ stands for the Fourier transform of $g$ on $\T^d \times \R^d$, that is for Schwartz function $g(x,v)$
\begin{align*}
    \widehat{g}(\xi, \eta)=\frac{1}{(2\pi)^{2d}}\int_{\T^d \times \R^d} e^{-i(\xi\cdot x+\eta \cdot v)}g(x,v) \, \mathrm{d}x \, \mathrm{d}v,  \ \ (\xi, \eta) \in \Z^d \times \R^d,
\end{align*}
then extended by duality on $\mathscr{S}'(\T^d \times \R^d)$. We also denote 
$$H^{m,-\infty}_{x,v} \vcentcolon= \bigcup_{r \in \R} H^{m,r}_{x,v}.$$

\end{Def}

\begin{Cor}
\label{cor}
Consider the same assumptions and notation as in Theorem~\ref{thm-LWPkinetic}.
We have 
$$f \in L^\infty([0,T]; H^{k-1,-\infty}_{x,v}).$$
\end{Cor}

\begin{proof}
We borrow the presentation of \cite{HKapde}; following \cite[Proposition 5.2]{Gerard}, we apply Theorem~\ref{thm-LWPkinetic} with the test function 
\begin{align*}
    \psi_\eta (v) = e^{-i v \cdot \eta} \in \mathscr{C}^\infty_b(\R^d),
\end{align*}
where $\eta \in \R^d$ should be viewed as the Fourier variable in velocity. By a view of the proof of Lemma~\ref{prop-reg-momentsf}, the conclusion of Theorem~\ref{thm-LWPkinetic} can be slightly refined into the statement
\begin{equation}
\label{eq:proof-cor}
\forall \eta \in \R^d, \qquad \left \|\int  \psi_\eta(v) \, f(t,x,\D v)  \right\|_{L^\infty([0,T]; H^{k-1}_x)} \leq \Lambda ( \| \psi_\eta\|_{{W}^{k-1,\infty}_{v}}), 
\end{equation}
where $\Lambda$ is a polynomial function.
Moreover, $ \| \psi_\eta\|_{{W}^{k-1,\infty}_{v}} \lesssim \widetilde{\Lambda}(|\eta|)$, where $\widetilde{\Lambda}$ is also a polynomial function.
If $\mathcal{F}_v$ stands for the partial Fourier transform on the variable $v$, we observe that
\begin{align*}
\int \psi_\eta(v) \, f(t,x,\D v)  = (2\pi)^d \mathcal{F}_v {f}(t,x, \eta),
\end{align*}
and we can deduce from~\eqref{eq:proof-cor} that for some $r>0$ taken large enough,
\begin{align*}
\left \|  \widehat{f}(t,\xi, \eta)  (1+ |\xi|^2)^{(k-1)/2} (1+ |\eta|^2)^{-r/2} \right\|_{L^\infty([0,T]; L^2(\Z^d \times \R^d))} < +\infty,
\end{align*}
which precisely means  that $f \in L^\infty([0,T]; H^{k-1,-r}_{x,v})$.
\end{proof}

	\section{Extensions to other functional spaces}\label{section:Extension}
	In this section, we explain how one can extend the abstract theory developed in Section \ref{subsec-proofTHM} for Sobolev spaces to weighted Sobolev spaces (in space) and Besov spaces.

	\subsection{Weighted Sobolev spaces }\label{subsec-weightedSob}
	 
We first explain
the extension to weighted Sobolev spaces. This functional set-up is particularly relevant when dealing with problems set in the whole space $\R^d$. Notably it allows the inversion of elliptic operators and provides good energy estimates (see for instance \cite[Appendix II]{choquetbruhat}). 
From the multiphasic point of view , it seems particularly adapted to handle equations where there is an auxiliary elliptic or parabolic equation for the force field, with possibly a complicated coupling. In particular, it is a way to bypass the use of negative (homogeneous) Sobolev spaces (for the densities) in order to handle low-frequencies on the whole space, and should therefore be compared with what we have done in Section \ref{sec:conseq-rho} -- see Assumption \ref{ass:bound_high_regularity-rho}.
We therefore believe this extended functional framework could be of interest to treat certain interesting models. We shall give an example in Section \ref{subsec:Stokes}, where we will consider the Vlasov--Stokes system posed in $\R^3$.

	 \medskip

	 Following \cite[Appendix I]{choquetbruhat}, we first introduce the weighted Sobolev spaces of interest (note that the weights increase with the number of derivatives in the norm). Here, the spatial domain is $\R^d$. 
	 
	 \begin{Def}
	The weighted Sobolev spaces $H^{k}_\delta$ with $k \in \N$ and $\delta \in \R$ as the completion of $\mathscr{C}^\infty_c(\R^d)$ under the norm
\begin{align}
    \Vert u \Vert_{H^{k}_\delta}&\vcentcolon = \left(\sum_{0 \leq \vert \alpha \vert \leq k } \int_{\R^d} \langle x \rangle^{2(\delta+ \vert \alpha \vert)} \vert \partial_x^\alpha u \vert^2 \, \mathrm{d}x  \right)^{1/2}, \ \ \langle x \rangle\vcentcolon = (1+\vert x \vert^2)^{1/2}.
\end{align}
\index{Hx@$H^k_\delta$: Sobolev space with varying polynomial weight in $x$}
	 \end{Def}
We will use the notation $L^{2}_\delta=H^{0}_\delta$. We also consider the weighted spaces $\mathscr{C}^{k}_\delta$ with $k \in \N$ and $\delta \in \R$ as the completion of $\mathscr{C}^\infty_c(\R^d)$ under the norm
\begin{align}
    \Vert u \Vert_{\mathscr{C}^{k}_\delta}&\vcentcolon = \sum_{0 \leq \vert \alpha \vert \leq k } \left\Vert \langle \cdot \rangle^{\delta+ \vert \alpha \vert}  \partial_x^\alpha u \right\Vert_{L^\infty(\R^d)}, \ \ \langle x \rangle\vcentcolon = (1+\vert x \vert^2)^{1/2}.
\end{align}
Let us collect some useful facts on these weighted spaces (see  \cite[Appendix I.3]{choquetbruhat}).
For the sake of simplicity, we decided not to include tame estimates.
\begin{Prop}\label{prop-Sobweight}
The weighted Sobolev spaces satisfy the following properties.
\begin{enumerate}
\item {\bf Continuous embeddings.} For all $\delta \geq \delta'$ and $k \geq k'$, we have  $H^{k}_\delta \hookrightarrow H^{k'}_{\delta'}$ and  $\mathscr{C}^{k}_\delta \hookrightarrow \mathscr{C}^{k'}_{\delta'}$.
\item  {\bf Sobolev embedding.} For $m\in \N$, $k >d/2$ and $\delta \geq \delta'-d/2$, we have the continuous embedding $H^{k+m}_{\delta} \hookrightarrow \mathscr{C}^{m}_{\delta'}$.
\item  {\bf Product laws.} If $k \leq \min(k_1,k_2)$, $k+d/2<k_1 + k_2$ and $\delta<\delta_1 + \delta_2+d/2$,
\begin{align*}
\Vert uv \Vert_{H^{k}_\delta} \leq \Vert u \Vert_{H^{k_1}_{\delta_1}} \Vert v \Vert_{H^{k_2}_{\delta_2}}.
\end{align*}
In particular, the space  $H^{k, \delta}$ is an algebra for $k>d/2$ and $\delta >-d/2$.
\item {\bf Poincaré inequality.} If $d \geq 3$, for all $\delta>-d/2$, we have
\begin{align*}
    \Vert u \Vert_{L^{2}_\delta} \lesssim \Vert \DDD u \Vert_{L^{2}_{\delta+1}}.
\end{align*}
\item {\bf Commutator laws.} For $\vert \alpha \vert \leq m$, $m>1+d/2$, $\delta>-d/2$ and $\delta'>-d/2$, we have 
\begin{align}\label{commut-law-weighted}
    \big\Vert [\partial_x^\alpha, b]u \big\Vert_{L^{2}_{\delta}} \lesssim \Vert \DDD b \Vert_{H^{m-1}_{\delta'}} \Vert u \Vert_{H^{m-1}_{\delta}}.
\end{align}
In particular, for $k$ large enough and $\delta >-d/2$, we have 
\begin{align*}
\left\langle  \DDD a, \DDD \big( (b\cdot \nabla )a \big) \right\rangle_{H^{k-1}_{\delta}} \lesssim  \| \DDD b \|_{H^{k}_{\delta}} \| \DDD a \|^2_{H^{k-1}_{\delta}}.
\end{align*}
\end{enumerate}
\end{Prop}
In practice, we will only rely on this lemma when the dimension $d$ is greater than $3$, hence the above Poincaré inequality holds true. With the former properties at hand, we first gather the version of the general Lemma \ref{LM:existence-vitesseSOBOLEV} for the momentum equation in the context of weighted Sobolev spaces. 
\begin{Lem}\label{LM:existence-vitesseSOBOLEV-weighted}
Let $d\geq 3$. Let $k\in \N$ be large enough and $\delta>-d/2$. Let $T>0$, $w_0 \in H^{k}_\delta$ and $v,S \in L^\infty([0,T]; H^{k}_\delta)$ be some vector fields. Then there exists a unique classical solution $w \in \mathscr{C}([0,T]; H^{k}_\delta)$ to
\begin{align*}
    \partial_t w + (v \cdot \nabla_x) w = S, \qquad w|_{t=0}= w_0,
\end{align*}
 and there exists $C>0$ such that, for all $t \in [0,T]$, $w$ satisfies the estimates
\begin{multline*}
   \|w(t)\|_{L^{2}_\delta}  \leq \exp\left[ C \int_0^t \| \DDD v(\tau) \|_{ H^{k-1}_{\delta+1}} \D \tau \right] \| w_0\|_{L^{2}_\delta} \\
   + C \int_0^t  \exp\left[ C \int_{\tau}^t \| \DDD v(s) \|_{ H^{k-1}_{\delta+1}} \D s \right] \| S(\tau)\|_{L^{2}_\delta} \D \tau ,
   \end{multline*}
   and
   \begin{multline*}
   \| \DDD w(t) \|_{H^{k-1}_{\delta+1}} \leq\exp\left[ C \int_0^t \| \DDD v(\tau) \|_{H^{k-1}_{\delta+1}} \| \D \tau \right] \| \DDD w_0 \|_{H^{k-1}_{\delta+1}}  \\
     + C \int_0^t  \exp\left[ C \int_{\tau}^t \| \DDD v(s) \|_{H^{k-1}_{\delta+1}} \| \D s \right] \| \DDD S(\tau) \|_{H^{k-1}_{\delta+1}} \D \tau.
\end{multline*}
\end{Lem}
\begin{proof}
We focus on the derivation of \textit{a priori} estimates, the existence and uniqueness part following from the method of characteristics. 
Multiplying the equation for $w$ by $\langle x \rangle^{2\delta} w$ and integrating in space, we get by integration by parts and Cauchy--Schwarz inequality
\begin{align*}
    \dfrac{\mathrm{d}}{\mathrm{d}t}\| w \|_{L^{2}_\delta}^2 & \lesssim \int \mathrm{div}(\langle x \rangle^{2\delta} v) \vert w \vert^2 +   \Vert S\Vert_{L^{2}_\delta}\Vert w \Vert_{L^{2}_\delta} \\
    &\lesssim \int \vert \mathrm{div} (v)\vert  \vert \langle x \rangle^{\delta}  w \vert^2
    +\int \langle x \rangle^{-1} \vert  v \vert  \vert \langle x \rangle^{\delta}  w \vert^2 +   \Vert S\Vert_{L^{2}_{\delta}}\Vert w \Vert_{L^{2}_{\delta}} \\
    &\lesssim \Vert v \Vert_{\mathscr{
    C}^{1}_{-1}}\Vert w \Vert_{L^{2}_{\delta}}^2+   \Vert S\Vert_{L^{2}_\delta}\Vert w \Vert_{L^{2}_\delta}.
\end{align*}
By Sobolev embedding from Proposition \ref{prop-Sobweight}, we have $\Vert v \Vert_{\mathscr{C}^{1}_{-1}} \lesssim \Vert v \Vert_{H^{k'}_{\delta'}}$ for $k'>1+d/2$ and $\delta'>-1-d/2$.
We then deduce by Poincaré inequality (hence imposing $d \geq 3$ and $\delta'>-d/2$ -- see Proposition \ref{prop-Sobweight}) that
\begin{align*}
    \Vert v \Vert_{\mathscr{
    C}^{1}_{-1}} &\lesssim  \Vert v \Vert_{L^2_{\delta'}}+ \Vert  \nabla v \Vert_{H^{k-1}_{\delta'+1}} \lesssim  \Vert \nabla v \Vert_{L^2_{\delta'+1}}+ \Vert  \nabla v \Vert_{H^{k-1}_{\delta'+1}}
\end{align*}
which means that by taking $\delta'=\delta$, we get
\begin{equation*}
\frac{\D}{\D t} \| w \|_{L^{2}_{\delta}} \leq \| \DDD v \|_{ H^{k-1}_{\delta+1}} \| w \|_{L^{2}_\delta}  + \| S \|_{L^{2}_\delta},
\end{equation*} 
which is enough to obtain the first estimate.
For the second one, we take one spatial derivative in the equation of $w$ and by taking the $H^{k-1}_{\delta+1}$ scalar product of this equation with $\DDD w$, we get
\begin{align*}
\frac{1}{2}\frac{\D}{\D t} \| \DDD w \|^2_{H^{k-1}_{\delta+1}} &\leq \vert \langle \DDD  ((v \cdot \nabla) w) ,\DDD w \rangle_{H^{k-1}_{\delta+1}} \vert  + \| \DDD S \|_{H^{k-1}_{\delta+1}} \| \DDD  w\|_{H^{k-1}_{\delta+1}} \\
& \leq K\| \DDD v \|_{H^{k-1}_{\delta+1}} \| \DDD w \|^2_{H^{k-1}_{\delta+1}}+\| \DDD S \|_{H^{k-1}_{\delta+1}} \| \DDD  w\|_{H^{k-1}_{\delta+1}},
\end{align*}
for some $K>0$, thanks to the commutator estimates of Proposition~\ref{prop-Sobweight}. Like before, we obtain 
\begin{equation*}
\frac{\D}{\D t} \| \DDD w \|_{H^{k-1}_{\delta+1}} \leq K  \| \DDD v \|_{H^{k-1}_{\delta+1}} \| \DDD w \|_{H^{k-1}_{\delta+1}} + \| \DDD S \|_{H^{k-1}_{\delta+1}}.
\end{equation*}
We conclude the proof thanks to Gronwall lemma.
\end{proof}

The previous lemma is enough to obtain a corresponding weighted version of the well-posedness, time regularity and stability from Theorem \ref{thm:existence} in $\R^d$ for $d\geq 3$. To state it, we introduce the following norm on the families of velocities: given a measure space $(I,\mu)$, we define
\begin{equation}
	\label{def:space-weighted}
	\mathcal H^{k,p}_\delta \vcentcolon = \left\{ \uu, \, \Vert \uu \Vert_{H^{k,p}_\delta} <+\infty \right\},
	\end{equation}
        \index{Hx@$\mathcal H^{k,p}_\delta$: multiphasic Sobolev space with varying polynomial weight in $x$}
	where
\begin{equation}
	\label{def:norm-weighted}
	\Vert \uu \Vert_{H^{k,p}_\delta} \vcentcolon = \| \uu \|_{L^p_\alpha L^{2}_ \delta(\R^d)} + \| \DDD \uu \|_{L^{\upinfty}_\alpha H^{k-1}_{\delta+1}(\R^d)},
	\end{equation}
 for $k \in \N$ large enough, $p \in [1, \infty]$ and $\delta>-d/2$. The key Assumptions \ref{ass:bound_high_regularity} and \ref{ass:stability} are then modified accordingly. More precisely:
 \begin{Ass}{D1}[Bound in high regularity]
		\label{ass:bound_high_regularity-weighted}
          \index{A4@Assumptions \ref{ass:bound_high_regularity-weighted}--\ref{ass:stability-weighted}: for the general $\vv$ framework (weighted $H^k$ spaces)}
		Let $d>3$. Let $\boldsymbol F$ be a family of force fields and $k_0,p$ be as in Definition~\ref{def:force_fields}. We assume that for all $k \geq k_0$ and $\delta>-d/2$, there exists an estimate function $\F$ such that for all $T>0$, for all $\vv \in L^1([0,T]; \mathcal H^{k,p}_{\delta})$, and all parameters $X \in \mathcal X$,
		we have
		\begin{equation}
		\label{eq:ass_estimate_high_reg-weighted}
		\int_0^T \Vert \FF[\vv,X](t) \Vert_{H^{k,p}_\delta} \D t \leq \F\left(T, \int_0^T \Vert\vv(t) \Vert_{H^{k,p}_\delta} \D t, \mathsf N_k(X) \right).
		\end{equation}
	\end{Ass}
	\begin{Ass}{D2}[Stability at low regularity]
		\label{ass:stability-weighted}
		Let $d>3$ and $\delta>-d/2$. Let $\boldsymbol F, k_0,p$ be as in Definition~\ref{def:force_fields}, and let $R>0$. There exists a stability estimate function $\G_R$ such that for all $T>0$, families $\vv_1,\vv_2 \in L^1([0,T]; \mathcal H^{k_0,p}_{\delta})$, and parameters $X_1$ and $X_2$, if
		\begin{equation}
		\label{eq:uniform_bound_velocities-weighted} \max\left( \sup_{t \in [0,T]} \Vert\vv_1(t)\Vert_{H^{k_0,p}_\delta},\sup_{t \in [0,T]} \Vert\vv_2(t)\Vert_{H^{k_0,p}_\delta}, \mathsf N_{k_0}(X_1), \mathsf N_{k_0}(X_2)  \right) \leq R,
		\end{equation}
	then the following estimate holds:
		\begin{multline}
		\label{eq:stability_assumption-weighted}
		\int_0^T \big\|\FF[\vv_1,X_1](t) - \FF[\vv_2,X_2](t)\big\|_{L^p_\alpha L^2_{\delta}} \D t \\
        \leq \G_R\left(T, \int_0^T \big\|\vv_1(t) - \vv_2(t) \big\|_{L^p_\alpha L^2_\delta} \D t, \mathsf d(X_1,X_2)\right).
		\end{multline}
	\end{Ass}
 This leads to the following result.
 \begin{Thm}\label{thm:existence-WEIGHTED}
Under the Assumptions \ref{ass:bound_high_regularity-weighted} and \ref{ass:stability-weighted}, the statement of Theorem \ref{thm:existence} still holds, replacing everywhere $L^2_x$ by $L^2_\delta$ and $\mathcal{H}^{k,p}$ by $\mathcal{H}^{k,p}_\delta$ for $k \in \N$ large enough and $\delta>-d/2$.
\end{Thm}
\begin{proof}
Lemma \ref{LM:existence-vitesseSOBOLEV-weighted} is the key result to prove a version of  \eqref{eq:estimate_bound_high_regularity_0_derivative}, \eqref{eq:estimate_bound_high_regularity_k_derivatives} in the context of weighted Sobolev spaces. The proof then follows from the exact same fixed-point as in Section \ref{subsec-proofTHM}. The required adaptation to the proof of Theorem \ref{thm:existence} are directly obtained by taking the suitable $L^p_\alpha$ or $L^\infty_\alpha$ norms respectively in the two estimates of Lemma \ref{LM:existence-vitesseSOBOLEV-weighted}. Lastly, the proof of \eqref{eq:bound_stablowLWP} can be performed in a similar way.
\end{proof}

We can also obtain some versions of Lemmas \ref{LM:SobEstimRHO}--\ref{LM:estimStabContinuity} for the continuity equation in the context of weighted Sobolev spaces. The following results is in the same spirit as Lemma \ref{LM:existence-vitesseSOBOLEV-weighted} before, and be obtained as in the proof of  Lemma \ref{LM:SobEstimRHO}. Details are therefore left to the reader.
\begin{Lem}\label{LM:SobEstimRHO-weighted}
Let $d \geq 3$ and $({I}, \mu)$ be a set of labels. Let $p \in [1,+\infty]$, $k \in \N$  sufficiently large and $\delta, \delta'>-d/2$. There exists $K>0$ such that the following holds. Let $(\rho^{\alpha},v^{\alpha})_{\alpha \in I}$ be smooth densities and velocity fields satisfying
    \begin{align*}
        \partial_t \rho^\alpha+\mathrm{div}(\rho^\alpha v^\alpha)=0,
    \end{align*} 
    with $\rho^\alpha_{\mid t=0} =\rho^\alpha_0$. Then for all $t \geq 0$
\begin{align}
     \Vert \rhorho(t)\Vert_{L^{\upinfty}_\alpha H^{k-1}_\delta} \leq  \Vert \rhorho_0 \Vert_{L^{\upinfty}_\alpha H^{k-1}_\delta} \exp\left( K \int_0^t\Vert\vv(s)\Vert_{H^{k,p}_{\delta'}} \D s \right).
    \end{align}
\end{Lem}
We now state an analogue of Lemma \ref{LM:estimStabContinuity}. 
\begin{Lem}\label{LM-estimWeightNegativeregRho}
Let $d \geq 3$ and $(I, \mu)$ be a set of labels. Let $p \in [1, \infty]$, $k \in \N$ sufficiently large , $\delta \in (-d/2,0)$ and $\delta'>-d/2$. There exists $K>0$ such that the following holds. Let $(\rhorho_1,\vv_1)= (\rho^{\alpha}_1,v^{\alpha}_1)_{\alpha \in I}$ and $(\rhorho_2,\vv_2)=(\rho^{\alpha}_2,v^{\alpha}_2)_{\alpha \in I}$ be two families of smooth densities and velocity fields solutions to 
    \begin{align*}
        \partial_t \rho^\alpha+\mathrm{div}(\rho^\alpha v^\alpha)=0,
    \end{align*} 
    with initial data ${\rho^\alpha_{1}}_{\mid t=0} =\rho^\alpha_{1,0}$ and ${\rho^\alpha_{2}}_{\mid t=0} =\rho^\alpha_{2,0}$, where $\textbf{v}_1,\textbf{v}_2 \in \mathcal H^{k,p}_\delta$ are fixed vector fields. Then for all $T>0$, $t \in [0,T]$, we have
\begin{align*}
	&\| \rhorho_2(t) - \rhorho_1(t) \|_{L^p_\alpha \dot H^{-1}_x} \\
    &\leq \exp\left( K \int_0^t (\| \DDD \vv_1(s) \|_{L^{\upinfty}_\alpha L^{\upinfty}_x} + \| \DDD \vv_2(s) \|_{L^{\upinfty}_\alpha L^{\upinfty}_x}) \D s \right) \\ & \ \times \Bigg\{ \| \rhorho_{2,0} - \rhorho_{1,0} \|_{L^p_\alpha  \dot H^{-1}_x} \\ & \qquad + \Vert  \rhorho_{1,0}\Vert_{L^{\upinfty}_\alpha H^{k-1}_{-\delta-d/2}}  \exp\left( K\Vert \vv_1 \Vert_{L^1([0,T];\mathcal{H}^{k,p}_{\delta'})}\right) \int_0^t \| \vv_2(s) - \vv_1(s) \|_{L^p_\alpha L^2_\delta} \D s \Bigg\}.
	\end{align*}
\end{Lem}
\begin{proof}
 A careful inspection of the end of the proof of Lemma \ref{LM:estimStabContinuity} reveals that we can perform the same steps, by only modifying the estimate for the last term of \eqref{goldorak}: introducing negative weights in the estimates, we can write for all $\delta \in \R $,
    \begin{align*}
        \cg \nabla \phi, \rho_1^\alpha (v_2^\alpha - v_1^\alpha) \cd 
        &\lesssim \| \langle x \rangle^{-\delta} \rho_1^\alpha \|_{L^{\upinfty}_x} \| \nabla \phi \|_{L^2_x} \| \langle x \rangle^{\delta} (v_2^\alpha - v_1^\alpha) \|_{L^2_x} \\
        &\lesssim \|  \rho_1 \|_{\mathscr{C}^0_{-\delta}} \| \nabla \phi \|_{L^2_x} \|  v_2 - v_1 \|_{L^2_\delta}.
    \end{align*}
Using the Sobolev embedding $\|  \rho_1^\alpha \|_{\mathscr{C}^0_{-\delta}} \lesssim \Vert \rho_1^\alpha \Vert_{H^{k-1}_{\delta''}}$ for all $\delta'' \geq -\delta-d/2$ from Proposition \ref{prop-Sobweight} combined with the high-order estimate from Lemma  \ref{LM:SobEstimRHO-weighted}, we then obtain for $k$ large enough, for all $\delta \in (-d/2, 0)$ and $\delta' >-d/2$
 \begin{multline*}
     \cg \nabla \phi, \rho_1^\alpha (v_2^\alpha - v_1^\alpha) \cd \\
     \lesssim \Vert  \rho_{1,0}^\alpha\Vert_{H^{k-1}_{-\delta-d/2}} \exp\left( K \Vert \vv_1 \Vert_{L^1([0,T];\mathcal{H}^{k,p}_{\delta'})}\right) \| \nabla \phi \|_{L^2_x} \|  v_2^\alpha - v_1^\alpha \|_{L^2_\delta}
 \end{multline*}
 on $[0,T]$, hence there exists $K>0$ such that for all $t \in [0,T]$
\begin{multline*}
	\| \rho_2^\alpha(t) - \rho_1^\alpha(t) \|_{\dot H^{-1}_x} \leq \exp\left( K \int_0^t (\| \DDD \vv_1(s) \|_{L^{\upinfty}_\alpha L^{\upinfty}_x} + \| \DDD \vv_2(s) \|_{L^{\upinfty}_\alpha L^{\upinfty}_x}) \D s \right) \\ \Big\{ \| \rho_{2,0}^\alpha - \rho_{1,0}^\alpha \|_{\dot H^{-1}_x} + \Vert  \rho_{1,0}^\alpha\Vert_{H^{k-1}_{-\delta-d/2}}   
    \exp\left( K \Vert \vv_1 \times \Vert_{L^1([0,T];\mathcal{H}^{k,p}_{\delta'})}\right)\\ \times \int_0^t \| v_2^\alpha(s) - v_1^\alpha(s) \|_{L^2_\delta} \D s \Big\},
		\end{multline*}
 for the same range of exponents. The $L^p_\alpha$ version of the inequality then follows.
\end{proof}
Combining Theorem \ref{thm:existence-WEIGHTED} with Lemmas~\ref{LM:SobEstimRHO-weighted}-- \ref{LM-estimWeightNegativeregRho}, we can now infer a version of Theorem \ref{thm-LWPmultiphase-densityvelocity} for the multiphasic system \eqref{eq:general_system_Chapter2} on  density--velocity fields, in the context of weighted Sobolev spaces. The Assumptions~\ref{ass:bound_high_regularity-rho} and \ref{ass:stability-rho} on the force field $\widetilde\FF[\vv, \rhorho,Y]$  must now read as follows (let us recall that $k_0$ and $p$ are fixed):

	\begin{Ass}{E1}[Bound in high regularity]
		\label{ass:bound_high_regularity-rho-weighted}
          \index{A5@Assumptions \ref{ass:bound_high_regularity-rho-weighted}--\ref{ass:stability-rho-weighted}: for the $(\rhorho,\vv)$ framework (weighted $H^k$ spaces)}
	We assume that for all $k \geq k_0$, there exists an estimate function $\F$ such that for all $T>0$, for all $\vv \in L^1([0,T]; \mathcal H^{k,p}_\delta)$, and all parameters $X \in \mathcal X$,
		we have
		\begin{multline}
				\label{eq:ass_estimate_high_reg-rho-weighted}
		\int_0^T \left\Vert\widetilde{\FF}[\vv,\rhorho,X](t) \right\Vert_{H^{k,p}_\delta} \D t \\
        \leq \F\left(T, \int_0^T \Vert\vv(t)\Vert_{H^{k,p}_\delta} \D t,  \| \rhorho\|_{L^\infty([0,T];L^{\upinfty}_\alpha H^{k-1}_{\delta+1})},  \mathsf N_k(X) \right).
		\end{multline}
	\end{Ass}

	\begin{Ass}{E2}[Stability at low regularity]
		\label{ass:stability-rho-weighted}
	There exists a stability estimate function $\G_R$ such that for all $T>0$, families $\vv_1,\vv_2 \in L^1([0,T]; \mathcal{H}^{k_0,p}_{\delta})$, 
	$\rhorho_1,\rhorho_2 \in L^\infty([0,T]; L^{\upinfty}_\alpha H^{k_0-1}_{-\delta-d/2})$,
	and parameters $Y_1$ and $Y_2$, if
		\begin{multline}
		\label{eq:uniform_bound_velocities-rho-weighted} \max\Bigg( \sup_{t \in [0,T]} \Vert\vv_1(t)\Vert_{H^{k_0,p}_\delta},\sup_{t \in [0,T]} \Vert\vv_2(t)\Vert_{H^{k_0,p}_\delta},  \\
        \| \rhorho_1\|_{L^\infty([0,T];L^{\upinfty}_\alpha H^{k_0-1}_{-\delta-d/2})}, \mathsf N_{k_0}(Y_1), \mathsf N_{k_0}(Y_2)  \Bigg) 
		\leq R,
		\end{multline}
	then the following estimate holds:
		\begin{multline}
		\label{eq:stability_assumption-rho-weighted}
		\int_0^T \left\|\widetilde{\FF}[\vv_1,\rhorho_1,Y_1](t) - \widetilde{\FF}[\vv_2,\rhorho_2,Y_2](t)\right\|_{L^p_\alpha L^2_\delta} \D t \\
		\leq \G_R\left(T, \int_0^T \big\|\vv_1(t) - \vv_2(t) \big\|_{L^p_\alpha L^2_\delta} \D t, \|\rhorho_1 - \rhorho_2\|_{L^\infty([0,T];L^p_\alpha \dot{H}^{-1}_x)}, \mathsf d(Y_1,Y_2)\right).
		\end{multline}
	\end{Ass}
 With these Assumptions at hand, we can now directly infer a version of the local well-posedness result from  Theorem \ref{thm-LWPmultiphase-densityvelocity} in the framework of weighted Sobolev spaces.
\begin{Thm}\label{thm-LWPmultiphase-densityvelocity-WEIGHTED}
Let $d \geq 3$. Let $(I, \mu)$ be a set of labels. Let $k \in \N$ large enough, $\delta>-d/2$ and $p \in [1, \infty]$. Suppose that there exists a well defined mapping $(\textbf{w}, \mathbf{\eta}) \mapsto \widetilde{\FF}(\textbf{w}, \mathbf{\eta},Y)$ satisfying the Assumptions~\ref{ass:bound_high_regularity-rho-weighted} and~\ref{ass:stability-rho-weighted}. If $(\rhorho_0, \vv_0) \in L^\infty_\alpha H^{k-1}_{\delta} \times \mathcal{H}^{k,p}_\delta$ is a regular multiphasic distribution, there exists $T>0$ and a unique solution $(\rhorho, \vv)$ to the multiphasic system with force field $\widetilde{\FF}(\vv,\rhorho^,Y)$ on $[0,T]$ in the sense of Definition \ref{def:solution_coupled} and such that 
\begin{align*}
    \rhorho \in L^\infty([0,T]; L^\infty_\alpha H^{k-1}_\delta), \ \ \vv \in L^\infty([0,T]; \mathcal{H}^{k,p}_\delta).
\end{align*}

\end{Thm}

\begin{Rem}    
Let us notice that, compared to the Assumptions \ref{ass:bound_high_regularity-rho}--\ref{ass:stability-rho} that we put forward in the case of Theorem \ref{thm-LWPmultiphase-densityvelocity} for standard Sobolev spaces (without weights), the former estimates \eqref{eq:ass_estimate_high_reg-rho-weighted} and  \eqref{eq:stability_assumption-rho-weighted} in the weighted context do not enforce any $\dot{H}^{-1}_x$ regularity on the initial densities (thanks to Lemma \ref{LM-estimWeightNegativeregRho} applied to two different families of densities having same initial value).

\end{Rem}

\begin{Rem}
    The analogue of Theorem~\ref{thm-LWPmultiphase-densityvelocity-re} which pertains to System~\ref{eq:general_system_Chapter2-re} holds as well for the weighted $H^k_\delta$ spaces, adapting Assumptions~\ref{ass:bound_high_regularity-rho-re} and~\ref{ass:stability-rho-re} accordingly, but we shall not dwell on this aspect. Likewise, the continuity in time properties obtained in the unweighted Sobolev case are still valid.
\end{Rem}

	\subsection{Besov spaces}	\label{Section-Besov}
We now turn to the extension of the abstract result to Besov spaces, another natural generalization of Sobolev spaces.
We begin with some preliminary material related to the classical Littlewood--Paley dyadic decomposition. For all $R,R'>0$, we consider $\mathcal{C}_{R,R'} = \{ \xi \in \R^d, \, R\leq |\xi| \leq R'\}$.
 Pick a radial function $\varphi \in \mathscr{D}(\R^d\setminus\{0\})$, with values in $[0,1]$, such that $\mathrm{supp} \varphi \subset \mathcal{C}_{2/3,7/3}$, that satifies
\begin{equation}
\label{eq:phi}
\forall \xi \neq 0, \quad \sum_{j<0} \varphi(2^{-j} \xi) + \sum_{j \geq 0}  \varphi(2^{-j} \xi)  =1,
\end{equation}
 and is such that $\operatorname{supp} \varphi(2^{-j}\cdot) \cap \operatorname{supp}  \varphi(2^{-j'}\cdot) = \emptyset$, for $|j-j'|\geq 2$ (see \cite[Chapter 2, Section 2.2]{BCD}).
We set $\chi(\xi) = \sum_{j<0} \varphi(2^{-j} \xi)  \in \mathscr{C}^\infty(\R^d\setminus\{0\})$. We can then define the localization operators $\dot{\Delta}_j$ (which correspond, loosely speaking, to localization in frequency space in $2^j \mathcal{C}_{2/3,7/3}$ ) and ${S}_j$ (which correspond, loosely speaking, localization in frequency space in the ball $B(0,{2^j})$). 

For all $j \in \Z$, define for all $a \in \mathcal{S}'(\R^d)$,
\begin{align}
\dot{\Delta}_j a &\vcentcolon= \mathcal{F}^{-1}( \varphi(2^{-j} \xi) \mathcal{F} a), \\
{S}_j a &\vcentcolon=  \mathcal{F}^{-1}( \chi(2^{-j} \xi) \mathcal{F} a),
\end{align}
where $\mathcal{F}$ stands for the Fourier transform in $\R^d$.
We can finally set 
$$
\Delta_j  = \left\{\begin{array}{ll} 0, &\qquad j \leq -2, \\ S_0, &\qquad j=-1, \\  \dot{\Delta}_j, &\qquad j \geq 0.\end{array}\right.
$$
\begin{Def}
Let $s \in \R$ and $q,r \in [1, \infty]$. For $\Omega= \R^d$ or $\T^d$, the Besov space $B^s_{q,r}(\Omega)$ is defined as the set of tempered distributions $u \in \mathscr{S}'(\R^d)$ (identifying functions on $\T^d$ with $(2\pi\Z)^d$ periodic functions) such that
\index{B@$B^s_{q,r}$: Besov space}
$$
\| u \|_{B^s_{q,r}(\Omega)}
 \vcentcolon= \left\| (2^{js} \| \Delta_j u \|_{L^q(\Omega)})_{j \in \Z}\right\|_{\ell^r(\Z)} <+\infty, $$
 endowed with the norm $\| \cdot \|_{B^s_{q,r}(\Omega)}$.
\end{Def}
Remarkably, Besov spaces allow to encode some well-known interesting function spaces.
\begin{itemize}
    \item For $q=r=2$ and $s \in \Z$,  $B^s_{2,2}$ coincides with the Sobolev space $H^s_x$; for general $s\in \R$ this corresponds to the usual fractional Sobolev space $H^s_x$:
    $$
H^s_x \vcentcolon=\left\{ \varphi \in \mathscr{S}', \, \| \varphi \|^2_{H^s_x} \vcentcolon= \int \langle\xi\rangle ^{2s} |\widehat \varphi(\xi)|^2 \D \xi <+\infty \right\}.
$$
    \item For $q=r=\infty$ and $s \in \R_+ \setminus \N$, $B^s_{\infty,\infty}$ coincides with the Hölder space $\mathscr{C}^{\lfloor s \rfloor,s- \lfloor s \rfloor}$.
\end{itemize}
Therefore the extension of Theorem~\ref{thm:existence} to Besov spaces would in particular imply the generalization to fractional Sobolev spaces and to Hölder spaces. We refer to \cite[Chapter 2]{BCD} for a comprehensive introduction to properties and applications of Besov spaces. In particular, the analogue of Proposition~\ref{prop-Besov} reads as follows for Besov spaces.
\begin{Prop}\label{prop-Besov}
Besov spaces satisfy the following properties.
\begin{enumerate}
\item  {\bf Sobolev--Besov embedding.} 
For $m \in \N$, $s > \frac{d}{p}$, we have the continuous embedding
$B^{s+m}_{p,r} \hookrightarrow \mathscr{C}^{m}$.

\item {\bf Tame estimates.} 
If $s>0$, then for all $u,v \in B^s_{q,r} \cap L^{\infty}$,
$$
\|uv\|_{B^s_{q,r}} \lesssim \|u\|_{L^{\upinfty}_x} \|v\|_{B^s_{q,r}} + \|u\|_{B^s_{q,r}} \|v\|_{L^{\upinfty}_x}.
$$
\item {\bf Commutator laws.} 
For $s > 1+ d/p$,  for  all vector fields $a,b$, there exist sequences $c_j$ with $\|c_j\|_{\ell^r}=1$, such that 
\begin{equation}
\label{eq:convective-besov}
\|[\Delta_j, b\cdot \nabla] a\|_{L^q_x} \lesssim 2^{-sj} c_j  \|\nabla b\|_{B^{s-1}_{q,r}}\|a\|_{B^{s}_{q,r}}.
\end{equation}
\end{enumerate}
\end{Prop}

Let us now explain the modification of the definition of the norm $\mathcal{H}^{k,p}$. We can almost keep the same notation, up to to replacing $k$ by a multi-index $(s,q,r)$ with $s \geq 0$, $q,r \in [1,\infty]$, the $L^2_x$ norm by the $B^{0}_{q,r}$ norm when there is no derivatives, and $k-1$ by  $s-1$ at the level of the derivatives. Namely, given a measure space $(I, \mu)$, we  define the space
\begin{equation}
	\label{def:space-Besov}
	\mathcal{B}^{(s,q,r),p} \vcentcolon = \left\{ \uu, \, \Vert \uu \Vert_{\mathcal{B}^{(s,q,r),p}} <+\infty \right\},
	\end{equation}
        \index{B@$\mathcal{B}^{(s,q,r),p} $: multiphasic Besov space}
	where
\begin{equation}\label{def:normBesov}
	\Vert\uu \Vert_{\mathcal{B}^{(s,q,r),p}} \vcentcolon = \| \uu \|_{L^p_\alpha  L^2_x} + \| \DDD \uu \|_{L^{\upinfty}_\alpha B^{s-1}_{q,r}}.
	\end{equation}

\begin{Rem}
    In the definition of this norm, the choice of the space $L^p_\alpha  L^2_x$ is made for simplicity. It indeed essentially allows to use similar (low regularity) stability estimates of in the Sobolev framework of Section~\ref{sec:abstract_existence}. It is possible to replace $L^2_x$ by an abstract functional space $E$ (with for instance $E=L^q_x, B^0_{q,r},\ldots$), modifying the required stability estimates accordingly.
\end{Rem}
    
The Assumptions \ref{ass:bound_high_regularity} and \ref{ass:stability} are then modified accordingly, as follows. 	Fix $q,r \in [1,\infty]$.

	\begin{Ass}{F1}[Bound in high regularity]
		\label{ass:bound_high_regularity-besov}
          \index{A6@Assumptions \ref{ass:bound_high_regularity-besov}--\ref{ass:stability-besov}: for the general $\vv$ framework (Besov spaces)}
	 Let $\boldsymbol F$ be a family of force fields and $s_0,p$ be as in Definition~\ref{def:force_fields}. We assume that for all $s \geq s_0$, there exists an estimate function $\F$ such that for all $T>0$, for all $\vv \in L^1([0,T]; \mathcal \mathcal{B}^{(s,q,r),p} )$, and all parameters $X \in \mathcal X$,
		we have
		\begin{equation}
		\label{eq:ass_estimate_high_reg-besov}
		\int_0^T \Vert \FF[\vv,X](t)\Vert_{\mathcal{B}^{(s,q,r),p}} \D t \leq \F\left(T, \int_0^T \Vert \vv(t)\Vert_{\mathcal{B}^{(s,q,r),p}} \D t, \mathsf N_s(X) \right).
		\end{equation}
	\end{Ass}

	\begin{Ass}{F2}[Stability at low regularity]
		\label{ass:stability-besov}
		Let $\boldsymbol F, s_0,p$ be as in Definition~\ref{def:force_fields}, and let $R>0$. There exists a stability estimate function $\G_R$ such that for all $T>0$, families $\vv_1,\vv_2 \in L^1([0,T]; \mathcal \mathcal{B}^{(s_0,q,r),p})$, and parameters $X_1$ and $X_2$, if
		\begin{equation}
		\label{eq:uniform_bound_velocities-besov} \max\Big( \sup_{t \in [0,T]} \Vert \vv_1(t) \Vert_{\mathcal{B}^{(s_0,q,r),p}},\sup_{t \in [0,T]} \Vert \vv_2(t)\Vert_{\mathcal{B}^{(s_0,q,r),p}}, \mathsf N_{s_0}(X_1), \mathsf N_{s_0}(X_2)  \Big) \leq R,
		\end{equation}
	then the following estimate holds:
		\begin{multline}
		\label{eq:stability_assumption-besov}
		\int_0^T \big\|\FF[\vv_1,X_1](t) - \FF[\vv_2,X_2](t)\big\|_{L^p_\alpha  L^2_x} \D t \\
        \leq \G_R\left(T, \int_0^T \big\|\vv_1(t) - \vv_2(t) \big\|_{L^p_\alpha L^2_x} \D t, \mathsf d(X_1,X_2)\right).
		\end{multline}
	\end{Ass}
We therefore obtain
\begin{Thm}\label{thm:existence-besov}
Under Assumptions \ref{ass:bound_high_regularity-besov} and \ref{ass:stability-besov}, the statement of Theorem \ref{thm:existence} still holds, replacing everywhere  $\mathcal{H}^{k,p}$ by $\mathcal{B}^{(s,q,r),p}$ for  $s>1+d/p$.
\end{Thm}

 In order to prove Theorem~\ref{thm:existence-besov}, we mainly need an version of Lemma \ref{LM:existence-vitesseSOBOLEV} adapted to the Besov setting, that allows for an application of the same fixed-point strategy. With Proposition~\ref{prop-Besov} at hand, we obtain  the following result (see also \cite[Theorem 3.14 p.133]{BCD}).
 Again, we do not state a tame version of this result, for the sake of simplicity.
 
\begin{Lem}\label{LM:existence-vitesseSOBOLEV-Besov}
Let $q,r \in [1,\infty]$, $s > d/q+1$, and $s' \in [0,s]$. Let $T>0$, $w_0 \in B^{s'}_{q,r}$ and $v,S \in L^\infty([0,T]; B^s_{q,r})$. Then there exists a unique classical solution $w$ to
\begin{align*}
    \partial_t w + v \cdot \nabla_x w = S, \qquad w|_{t=0}= w_0,
\end{align*}
in $L^\infty([0,T]; B^{s'}_{q,r}) \cap \mathscr{C}^0([0,T); B^{s'-\eps}_{q,r})$ for all $\eps>0$. Furthermore, if $r<+\infty$ or $s'<s$, then $w  \in \mathscr{C}^0([0,T); B^{s'}_{q,r})$.  Finally there is $C>0$ such that for all $t \in [0,T]$
\begin{multline*}
   \|w(t)\|_{B^0_{q,r}}  \leq \exp\left[ C \int_0^t \| \DDD v(\tau) \|_{ B^{s-1}_{q,r}} \D \tau \right] \| w_0\|_{B^0_{q,r}}  \\
   + C \int_0^t  \exp\left[ C \int_{\tau}^t \| \DDD v(s) \|_{ B^{s-1}_{q,r}} \D s \right] \| S(\tau)\|_{B^0_{q,r}} \D \tau ,
  \end{multline*}
   and
   \begin{multline*}
   \| \DDD w(t) \|_{B^{s'-1}_{q,r}} \leq\exp\left[ C \int_0^t \| \DDD v(\tau) \|_{B^{s-1}_{q,r}} \| \D \tau \right] \| \DDD w_0 \|_{B^{s'-1}_{q,r}}  \\
   + C \int_0^t  \exp\left[ C \int_{\tau}^t \| \DDD v(s) \|_{B^{s-1}_{q,r}} \| \D s \right] \| \DDD S(\tau) \|_{B^{s'-1}_{q,r}} \D \tau.
\end{multline*}
\end{Lem}

Now we provide estimates in Besov spaces for the solution to continuity equations. Applying Lemma~\ref{LM:existence-vitesseSOBOLEV-Besov}, we first obtain
\begin{Lem}\label{LM:SobEstimRHO-besov}
Let $(\mathcal{I}, \mu)$ be a set of labels. Let $p,q,r \in [1,\infty]$ and $s > d/q+1$. There exists $K>0$ such that the following holds. Let $(\rho^{\alpha},v^{\alpha})_{\alpha \in I}$ be smooth densities and velocities field satisfying
    \begin{align*}
        \partial_t \rho^\alpha+\mathrm{div}(\rho^\alpha v^\alpha)=0,
    \end{align*} 
    with $\rho^\alpha_{\mid t=0} =\rho^\alpha_0$. Then for all $t \geq 0$
\begin{align}
     \Vert \rhorho(t) \Vert_{L^{\upinfty}_\alpha B^{s-1}_{q,r}} \leq  \Vert \rhorho_0 \Vert_{L^{\upinfty}_\alpha B^{s-1}_{q,r}} \exp\left( K \int_0^t[\vv(\tau)]_{(s,q,r),p} \D \tau \right).
    \end{align}
\end{Lem}

Combining Theorem \ref{thm:existence-besov} and Lemmas~\ref{LM:SobEstimRHO-besov}, we can now infer a version of Theorem \ref{thm-LWPmultiphase-densityvelocity} for the multiphasic system \eqref{eq:general_system_Chapter2} on  density--velocity fields, in the context of Besov spaces. 
The Assumptions~\ref{ass:bound_high_regularity-rho} and \ref{ass:stability-rho} must read in this framework as follows:

	\begin{Ass}{G1}[Bound in high regularity]
		\label{ass:bound_high_regularity-rho-besov}
           \index{A7@Assumptions \ref{ass:bound_high_regularity-rho-besov}--\ref{ass:stability-rho-besov}: for the  $(\rhorho,\vv)$ framework (Besov spaces)}
	We assume that for all $s \geq s_0$, there exists an estimate function $\F$ such that for all $T>0$, for all $\vv \in L^1([0,T]; \mathcal{B}^{(s,q,r),p})$, and all parameters $X \in \mathcal X$,
		we have
		\begin{multline}
		\label{eq:ass_estimate_high_reg-rho-besov}
		\int_0^T \left\Vert\widetilde{\FF}[\vv,\rhorho,X](t)\right\Vert_{\mathcal{B}^{(s,q,r),p}} \D t \\
        \leq \F\left(T, \int_0^T \Vert\vv(t)\Vert_{\mathcal{B}^{(s,q,r),p}} \D t,  \| \rhorho\|_{L^\infty([0,T];L^{\upinfty}_\alpha (B^{s-1}_{q,r}\cap\dot{H}^{-1}_x))},  \mathsf N_s(Y) \right).
		\end{multline}
	\end{Ass}

	\begin{Ass}{G2}[Stability at low regularity]
		\label{ass:stability-rho-besov}
	There exists a stability estimate function $\G_R$ such that for all $T>0$, families $\vv_1,\vv_2 \in L^1([0,T]; \mathcal{B}^{(s_0,q,r),p})$, 
	$\rhorho_1,\rhorho_2 \in L^\infty([0,T]; L^{\upinfty}_\alpha B^{s_0-1}_{q,r})$,
	and parameters $Y_1$ and $Y_2$, if
		\begin{multline}
		\label{eq:uniform_bound_velocities-rho-besov} \max\Big( \sup_{t \in [0,T]} \Vert\vv_1(t)\Vert_{\mathcal{B}^{(s_0,q,r),p}},\sup_{t \in [0,T]} \Vert \vv_2(t)\Vert_{\mathcal{B}^{(s_0,q,r),p}}, \\
        \| (\rhorho_1,\rhorho_2)\|_{L^\infty([0,T];L^{\upinfty}_\alpha B^{s_0-1}_{q,r})}, \mathsf N_{s_0}(Y_1), \mathsf N_{s_0}(Y_2)  \Big) 
		\leq R,
		\end{multline}
	then the following estimate holds:
		\begin{multline}
		\label{eq:stability_assumption-rho-besov}
		\int_0^T \big\|\widetilde{\FF}[\vv_1,\rhorho_1,Y_1](t) - \widetilde{\FF}[\vv_2,\rhorho_2,Y_2](t)\big\|_{L^p_\alpha L^2_x} \D t \\
		\leq \G_R\left(T, \int_0^T \big\|\vv_1(t) - \vv_2(t) \big\|_{L^p_\alpha L^2_x} \D t, , \|\rhorho_1 - \rhorho_2\|_{L^\infty([0,T];L^p
_\alpha \dot{H}^{-1}_x)},  \mathsf d(Y_1,Y_2)\right).
		\end{multline}
	\end{Ass}

\begin{Thm}\label{thm-LWPmultiphase-densityvelocity-BESOV}
Consider $(I, \mu)$ a set of labels. Let $s>0, p,q,r\in [1,+\infty]$ with $s>1+d/q$. Suppose that there exists a well defined mapping $\textbf{w} \mapsto \widetilde\FF(\textbf{w}, \rhorho[\textbf{w}] )$ satisfying  Assumptions~\ref{ass:bound_high_regularity-rho-besov}--\ref{ass:stability-rho-besov}. 

If $(\rhorho_0, \vv_0) \in L^\infty_\alpha (B^{s-1}_{q,r}\cap \dot{H}^{-1}) \times \mathcal{B}^{(s,q,r),p}$ is a regular multiphasic distribution, there exists $T>0$ and a unique solution $(\rhorho, \vv)$ to the multiphasic system with force field $\FF(\rhorho[\vv], \vv)$ on $[0,T]$ in the sense of Definition \ref{def:solution_coupled} and such that 
\begin{align*}
    &\rhorho \in L^\infty([0,T]; L^\infty_\alpha (B^{s-1}_{q,r}\cap \dot{H}^{-1}_x))\cap L^\infty([0,T]; L^\infty_\alpha (B^{s-1-\eps}_{q,r}\cap \dot{H}^{-1}_x)), \\
    &\vv \in L^\infty([0,T]; \mathcal{B}^{(s,q,r),p})\cap L^\infty([0,T]; \mathcal{B}^{(s-\eps,q,r),p}),
\end{align*}
for any $\eps>0$. Moreover if $r<+\infty$,
$$
  \rhorho \in L^\infty([0,T]; L^\infty_\alpha (B^{s-1}_{q,r}\cap \dot{H}^{-1}_x)), \ \ \vv \in L^\infty([0,T]; \mathcal{B}^{(s,q,r),p}),
$$
\end{Thm}

\begin{Rem}
\label{rem:besov-lambda}
    The analogue of Theorem~\ref{thm-LWPmultiphase-densityvelocity-re} which pertains to System~\ref{eq:general_system_Chapter2-re} also holds for Besov spaces, adapting Assumptions~\ref{ass:bound_high_regularity-rho-re} and~\ref{ass:stability-rho-re} accordingly, but we shall not dwell on this aspect. Similar time continuity properties compared to the Sobolev case can also be inferred. 
\end{Rem}

\begin{Rem}
    We have chosen to keep things as simple as possible and have only developed a basic Besov spaces theory; it is however clear to us that a sharper theory could be obtained following more advanced tools such as those introduced in \cite{BCD}, see \cite{LemarieMultiphaseVNS} for the case of the incompressible Vlasov--Navier--Stokes system.
\end{Rem}

\section{Summary of the assumptions for abstract well-posedness}
\label{subsec:summary-ass}
In this chapter, we have introduced different sets of assumptions pertaining to ``abstract'' multiphasic systems, that ensure their well-posedness. These assumptions, which bear on the forcing term in the equations for velocities, always come in pair; loosely speaking, one asks for
\begin{itemize}
\item a bound in {\bf high regularity},
\item a stability estimate in {\bf low regularity}, assuming some control of the different quantities that appear in the equations.
\end{itemize}
Several such pairs of assumptions appear in this chapter, in order to cover several possible situations, for instance the functional spaces at stake.

In the following of this monograph, we intend to apply the abstract results of this chapter to concrete equations such as Vlasov--Poisson and Vlasov--Navier--Stokes. To this end, at each application, we will need  to check the relevant pair of assumptions. To ease readability, we summarize below the different sets of assumptions of this chapter and indicate where they will be relevant.

\begin{itemize}
\item {\bf Assumptions} \ref{ass:bound_high_regularity}--\ref{ass:stability} pertain to the most general multiphasic system
	\begin{equation}
	\label{eq:general_system-end}
	\partial_t v^\alpha + ( v^\alpha \cdot \nabla ) v^\alpha = F^\alpha[\vv,X].
	\end{equation}
	The functional setting is that of $H^k_x$ Sobolev spaces. 
	
	This set of assumptions will be convenient to use in Chapter \ref{Part2-VP}, for the study of the well-posedness of the Vlasov--Monge--Amp\`ere equation (Section \ref{Subsec:VMongeAMpere}), and for that of the semiclassical limit from Hartree with mixed states to Vlasov--Poisson (Section \ref{sec:semiclassical}).

\item {\bf Assumptions} \ref{ass:bound_high_regularity-rho}--\ref{ass:stability-rho} pertain to the multiphasic system
	\begin{equation}
	\label{eq:general_system_Chapter2-end}
	\left\{
	\begin{gathered}
	\partial_t \rho^\alpha + \Div( \rho^\alpha v^\alpha) =0, \\
	\partial_t v^\alpha + ( v^\alpha \cdot \nabla ) v^\alpha = F^\alpha[\vv,\rhorho,Y],\\
	\rho^\alpha |_{t = 0} = \rho^\alpha_0, \ v^\alpha |_{t = 0} = v^\alpha_0.
	\end{gathered}
	\right.
	\end{equation}
		The functional setting is also that of $H^k_x$ Sobolev spaces.

\item {\bf Assumptions} \ref{ass:bound_high_regularity-rho-re}--\ref{ass:stability-rho-re}
are close to Assumptions \ref{ass:bound_high_regularity-rho}--\ref{ass:stability-rho}, the main difference being that they pertain to a modification of \eqref{eq:general_system_Chapter2-end}, namely
	\begin{equation}
	\label{eq:general_system_Chapter2-re-end}
	\left\{
	\begin{gathered}
	\partial_t \rho^\alpha + \Div( \rho^\alpha (v^\alpha+\lambda^\alpha) =0, \\
	\partial_t w^\alpha + ( (w^\alpha+\lambda^\alpha) \cdot \nabla ) w^\alpha = G^\alpha[\ww,\rhorho,Y],\\
	\rho^\alpha |_{t = 0} = \rho^\alpha_0, \ w^\alpha |_{t = 0} = w^\alpha_0.
	\end{gathered}
	\right.
	\end{equation}

		The two sets of assumptions \ref{ass:bound_high_regularity-rho}--\ref{ass:stability-rho} and  \ref{ass:bound_high_regularity-rho-re}--\ref{ass:stability-rho-re} will be the most relevant in the following, namely in Chapter \ref{Part2-VP}, for the study of well-posedness of the Vlasov--Poisson system, for electrons (Section \ref{Subsec:VPelectrons} ) or for ions (Section \ref{Subsec:VPions}) and in Chapter \ref{Part3-VNS}, for the study of well-posedness of the Vlasov--Navier--Stokes system, in the incompressible case (Section \ref{SubsecVNSmultiphaseLWP}) and in the compressible cas (Section \ref{SubsecVNScomp-multiphaseLWP}).

\item {\bf Assumptions} \ref{ass:bound_high_regularity-weighted}--\ref{ass:stability-weighted} (resp. {\bf Assumptions} \ref{ass:bound_high_regularity-rho-weighted}--\ref{ass:stability-rho-weighted}) pertain to \eqref{eq:general_system_Chapter2-end} (resp. \eqref{eq:general_system_Chapter2-re-end}). 

		The functional setting is that of  weighted $H^k_x$ Sobolev spaces.  
{\bf Assumptions} \ref{ass:bound_high_regularity-rho-weighted}--\ref{ass:stability-rho-weighted} will only be used once in the following, namely in Chapter \ref{Part3-VNS} for the study of the Vlasov--Stokes equation (Section \ref{subsec:Stokes}).

\item {\bf Assumptions} \ref{ass:bound_high_regularity-besov}--\ref{ass:stability-besov} (resp. {\bf Assumptions} \ref{ass:bound_high_regularity-rho-besov}--\ref{ass:stability-rho-besov})  pertain to \eqref{eq:general_system_Chapter2-end} (resp. \eqref{eq:general_system_Chapter2-re-end}). 

	The functional setting is that of Besov spaces.  
{\bf Assumptions} \ref{ass:bound_high_regularity-rho-besov}--\ref{ass:stability-rho-besov} will only be used once in the following, namely in Chapter \ref{Part2-VP} for the study of the Vlasov--Poisson equation for electrons (Section \ref{Subsec:VPelectrons}).

\end{itemize}

	\chapter{Application to Vlasov--Poisson type systems}\label{Part2-VP}

In this chapter, we build on the general multiphase framework developed in Chapter~\ref{Part1-LWP} to prove several results for various Vlasov--Poisson type systems. 
We consider a selection of kinetic models from plasma and gravitational physics, namely the Vlasov--Poisson equation for electrons, the Vlasov--Poisson equation for ions, and the Vlasov--Monge--Ampère equations, which are presented in Section \ref{Section-PresentationVPoissonandco}.

The chapter is divided in three main parts:

\begin{itemize} 

\item The first part focuses on {\bf local well-posedness} results.
Concretely, we thoroughly apply the abstract local well-posedness theory from Chapter \ref{Part1-LWP}; to this purpose we show that the various Vlasov--Poisson type systems satisfy the abstract assumptions,  hence inferring several original statements for their multiphase and kinetic formulations. This set of results and the corresponding proof for each system  are respectively displayed in Sections \ref{Subsec:VPelectrons}--\ref{Subsec:VPions}--\ref{Subsec:VMongeAMpere}. For the sake of presentation, and as its algebraic structure is simpler to handle, more complete results are given for the model for electrons.

\item The second part considers a seemingly different problem, that is the derivation of Vlasov--Poisson as the {\bf semiclassical limit} of a system of Schr\"odinger--Poisson equations for mixed states. We explain how to leverage the multiphase framework to obtain some results that seem new in the literature (this is performed in Section \ref{sec:semiclassical}).

\item The final part is dedicated to {\bf instability} properties of some possibly rough homogeneous equilibria (e.g. sums of Dirac masses in velocity) for Vlasov--Poisson type systems.
We  develop in Section \ref{Section-FrameworkInstab} an abstract robust theory to prove nonlinear instability results, that builds on the multiphasic framework. We eventually explain in Section \ref{sec:insta-appli} why the Vlasov--Poisson equation for electrons and for ions satisfy the required abstract assumptions, thus obtaining nonlinear instability results around non-smooth unstable Penrose homogeneous equilibria.

\end{itemize}

\section{Presentation of the models}\label{Section-PresentationVPoissonandco}
To ease readability, we first introduce the three main kinetic models at stake, together with their multiphasic formulation.

\paragraph{Vlasov--Poisson system for electrons.} 

The standard Vlasov--Poisson equation for electrons, set on $\Omega \times \R^d$ with $\Omega= \T^d$ or $\R^d$, reads as 
\begin{equation}
		\label{eq:VPkin_sectionLWP}
		\left\{
		\begin{gathered}
		\partial_t f + v \cdot \nabla_x f - \nabla_x U\cdot \nabla_v f = 0,\\
		-\Delta_x U = \left\{
    \begin{array}{ll}
         \displaystyle \int f\D v - \displaystyle \iint f\D v \D x, \qquad & \mbox{if }  \ \ \Omega=\T^d, \\
        \displaystyle \int f\D v, \qquad & \mbox{if} \ \  \Omega=\R^d.
    \end{array}
\right.
		\end{gathered}
		\right.
		\end{equation}
        This equation arises when describing the dynamics of electrons,  assuming  ions are so heavy that they remain still at a constant density.
The Cauchy problem for \eqref{eq:VPkin_sectionLWP} has a long history and is by now rather well understood. Let us state right away two emblematic theorems concerning the Cauchy problem for the Vlasov--Poisson system \eqref{eq:VPkin_sectionLWP} set on $\R^3 \times \R^3$.
\begin{Thm}
\label{thm:vp-classical}
    The following global existence results hold for the Vlasov--Poisson equation.
    \begin{itemize}
        \item (Weak solutions, from \cite{arsenev1975existence,BattRein-weak}) Let $f_0 \in L^1_{x,v} \cap L^\infty_{x,v}$ with finite energy, that is
        $$\frac{1}{2}\int \vert v \vert^2 f_0\, \mathrm{d}x \, \mathrm{d}v+\int \vert \nabla_x U_{f_0} \vert^2 \, \mathrm{d}x<+\infty,$$
        where $U_{f_0}$ is the electric potential generated by $f_0$. Then there exists a global weak solution $f \in L^\infty(\R^+; L^1 \cap L^\infty_{x,v}) \cap \mathscr{C}(\R^+; w-L^1_{x,v})$ to \eqref{eq:VPkin_sectionLWP} with initial data $f_0$ and such that the energy (with $f_0$ replaced by $f(t)$) stays bounded.
        \item (Strong solutions, from \cite{Pfaffelmoser,lions1991propagation}) Let $f_0 \in \mathscr{C}^1_{x,v}$ be nonnegative with compact support or with sufficient (polyomial) decay in velocity. Then there exists a unique global solution $f \in \mathscr{C}^1_{t,x,v}$ to \eqref{eq:VPkin_sectionLWP} with initial data $f_0$.
    \end{itemize}
\end{Thm}
For weak solutions, uniqueness is in general not known, unless the density $\rho_f= \int f \D v$ is uniformly bounded \cite{loeper2006uniqueness} (see also \cite{Miot2016uniqueness,CISS} for improvements of this criterion). A class of very weak global solutions, named Lagrangian solutions, which allows to consider initial data $f_0$ only with finite mass and energy, has been obtained in \cite{ACF}.
The theory of strong solutions was settled early in dimension $d=2$ in \cite{UO}; in dimension $d=3$ local strong solutions were also constructed in this reference.
In dimension $d=3$, small data global solutions were first obtained in \cite{BardosDegond}; then global solutions without size restriction were obtained  simultaneously in \cite{Pfaffelmoser}, which was able to treat compactly supported data with a Lagrangian approach, and in \cite{lions1991propagation}, which was able to treat data with high enough moments in velocity with an approach based on the propagation of these moments; various further improvements are due to \cite{Schaeffer,BattRein-strong,Horst,pallard2014-VPtorus}. 
In recent papers, local existence in Sobolev or Besov spaces with improvements on the indices of regularity has been explored in \cite{HC1,HC2,HC3,JT,Tae}.  In the recent \cite{Ngu}, yet another class of unique local solutions has been put forward, asking for initial densities with sufficient integrability in $x$ and H\"older regularity in $v$. 
The mathematical literature related to~\eqref{eq:VPkin_sectionLWP} is huge; the references given here are very far from exhaustive.

\medskip

As already explained in the introduction, the multiphasic formulation of \eqref{eq:VPkin_sectionLWP} reads as follows, with the ansatz $f=\mathcal{I}_\mu(\rhorho,\vv)$:
\begin{equation}
\label{eq:VP}
\left\{ 
\begin{gathered}
\partial_t \rho^\alpha + \Div ( \rho^\alpha v^\alpha) = 0,\\
\partial_t v^\alpha + (v^\alpha \cdot \nabla) v^\alpha = - \nabla U,\\-\Delta U = \left\{
    \begin{array}{ll}
         \displaystyle \int \rho^\alpha  \D \mu(\alpha)-\iint \rho^\alpha   \D \mu(\alpha) \D x, \qquad & \mbox{if }  \ \ \Omega=\T^d, \\
        \displaystyle \int \rho^\alpha  \D \mu(\alpha), \qquad & \mbox{if} \ \  \Omega=\R^d.
    \end{array}
\right.
\end{gathered}
\right.
\end{equation}
As mentioned before, for $\Omega= \T^d$, System \eqref{eq:VP} is exactly the set of equations studied by Grenier in~\cite{Gr96}. Note that in the particular case where $\mu$ is a sum of two Dirac masses, this system corresponds to two-fluid Euler--Poisson models for two fluids, similar to the ones considered in~\cite{CGG}.  Finally, it is possible to generalize~\eqref{eq:VPkin_sectionLWP} and \eqref{eq:VP} to the case of electrons with other ion species which also follow a Vlasov dynamics, we shall come back to this point later.

\paragraph{Vlasov--Poisson system for ions.} 
In dimension $d\leq 3$\footnote{This restriction is only for the sake of simplicity. Higher dimension could be treated as well, at the expense of more technicalities.}, we also consider the following Vlasov--Poisson system set on the torus $\T^d$:
\begin{equation}\label{eq:VPions-kinetic}
\left\{ 
\begin{gathered}
\partial_t f + v \cdot \nabla_x f - \nabla_x U\cdot \nabla_v f = 0,\\
		e^{U}-\Delta_x U = \int f\D v.\\
\end{gathered}
\right.
\end{equation}
This equation arises when describing the dynamics of ions, assuming electrons are so light and fast that they have reached their thermodynamic equilibrium, so that their density $\rho_e$ follow a Maxwell-Boltzmann law $\rho_e=e^U$(see \cite{BardosGolseNguyenSentis} for a formal derivation).  This model has also attracted interest in the mathematical literature, as an alternative\footnote{We mention that the model with a linearized Poisson equation ${U}-\Delta_x U = \int f\D v$ is also often considered, see e.g. \cite{HKR,IRW}.} to the more standard Vlasov--Poisson equation for electrons~\eqref{eq:VPkin_sectionLWP}. Despite the nonlinear nature of the Poisson equation, several results pertaining to~\eqref{eq:VPkin_sectionLWP} were adapted for~\eqref{eq:VPions-kinetic}; in particular, we have the analogue of Theorem~\ref{thm:vp-classical}.
\begin{Thm}
\label{thm:vp-classical-ions}
    The following global existence results hold for the Vlasov--Poisson equation for ions.
    \begin{itemize}
        \item (Weak solutions, adapted from \cite{bouchut1991global}) Let $f_0 \in L^1_{x,v} \cap L^\infty_{x,v}$ with finite energy, that is
        $$\frac{1}{2}\int \vert v \vert^2 f_0\, \mathrm{d}x \, \mathrm{d}v+\int \vert \nabla_x U_{f_0} \vert^2 \, \mathrm{d}x+ \int (U_{f_0}-1) e^{U_{f_0}} \, \mathrm{d} x <+\infty ,$$
        where $U_{f_0}$ is the electric potential generated by $f_0$. Then there exists a global weak solution $f \in L^\infty(\R^+; L^1 \cap L^\infty_{x,v}) \cap \mathscr{C}(\R^+; w-L^1_{x,v})$ to \eqref{eq:VPions-kinetic} with initial data $f_0$ and such that the energy (with $f_0$ replaced by $f(t)$) stays bounded.
        \item (Strong solutions, from \cite{GPI}) Let $f_0 \in \mathscr{C}^1_{x,v}$ be nonnegative with compact support or with sufficient (polynomial) decay in velocity. Then there exists a unique global solution $f \in \mathscr{C}^1_{t,x,v}$ to \eqref{eq:VPions-kinetic} with initial data $f_0$.
    \end{itemize}
\end{Thm}
 The case of the whole space $\R^d$ could be considered as well, at the cost of replacing the electron density $e^U$ by $n(x)e^U$ where $n \in L^1(\R^d)$, but for the sake of simplicity we will not do so in this work; this was done in dimension $d=3$  in \cite{bouchut1991global} for weak solutions and in \cite{griffin2020globalwholespace} for strong solutions. Let us be elusive and only cite without description the recent references \cite{HKI-ions1d,GPI2,HNX,GagnebinIacobelli,CKS} which are devoted to the study of~\eqref{eq:VPions-kinetic}.

  \medskip

The multiphasic formulation of the Vlasov--Poisson equations for ions reads as follows, with the ansatz $f=\mathcal{I}_\mu(\rhorho, \vv)$:
\begin{equation}
\label{eq:VPions}
\left\{ 
\begin{gathered}
\partial_t \rho^\alpha + \Div ( \rho^\alpha v^\alpha) = 0,\\
\partial_t v^\alpha + (v^\alpha \cdot \nabla) v^\alpha = - \nabla U,\\
e^U-\Delta U = \int \rho^\alpha  \D \mu(\alpha).
\end{gathered}
\right.
\end{equation}
In the monokinetic case, this system has recently received a lot of attention, see e.g. \cite{choi2025critical, BaeKimKwon-singularity, BCK, bae2025emergence, song2026h-highMach, EJK}.

\paragraph{Vlasov--Monge--Ampère system
for gravitation.} 
We also address the analysis of the following Vlasov--Monge--Ampère equation, set on the torus $\T^d$:

\begin{equation}\label{eq:VM-kinetic}
\left\{ 
\begin{gathered}
\partial_t f + v \cdot \nabla_x f - \nabla_x U\cdot \nabla_v f = 0,\\
		\mathrm{det}\left( \mathrm{I}+ \mathrm{D}^2_x U\right) = \int f \, \mathrm{d}v,\\
\end{gathered}
\right.
\end{equation}
where we assume implicitly that the total mass of particle is one, that is
\begin{equation*}
    \iint f \D v \D x = 1,
\end{equation*}
a property formally preserved along the flow.

This model has been introduced in~\cite{BreLoe} as a geometric, fully nonlinear counterpart of the Vlasov--Poisson system for gravitation, that is, the system~\eqref{eq:VPkin_sectionLWP} with the opposite sign for $U$. Indeed linearizing the elliptic equation of \eqref{eq:VM-kinetic} around the zero potential yields the following:
\begin{align*}
    \Delta_x U = \int f\D v - 1,
\end{align*}
which is up to the sign the Poisson equation solved by the potential in~\eqref{eq:VPkin_sectionLWP}.

This system is said to be \emph{geometric} because it is deeply linked with the theory of Optimal Transport (see e.g. \cite{Santambrogio}). Indeed, the Monge--Ampère equation
$$
\mathrm{det}\left( \mathrm{I}+ \mathrm{D}^2_x U\right) = \int f \, \mathrm{d}v
$$
asserts that the map
\begin{equation}
\label{eq:transport_map}
    x \in \T^d \mapsto x + \nabla_x U(x)
\end{equation}
transports the density $\int f \D v$ to the Lebesgue measure. Actually, such a map is not unique and it is usual to assume implicitly that the map considered is the optimal Monge map given by quadratic optimal transport, which means assuming that for all $x \in \T^d$,
\begin{equation*}
    \mathrm D^2_x U(x) \geq - \Id
\end{equation*}
in the sense of symmetric matrices. With this choice, it follows that at a given time, a particle located at position $x \in \T^d$ undergoes the force
\begin{equation*}
    x - T(x)
\end{equation*}
where $T$ is the map obtained by solving the above quadratic optimal transport problem.

Local existence of classical solution as well as global existence of weak solutions have been proved in~\cite{BreLoe}. More precisely, we have the following result.
\begin{Thm}[From \cite{BreLoe}]\label{thmBrenierLoeper-VMA} The following holds for the Vlasov--Monge--Ampère equation.
    \begin{itemize}
        \item (Weak solutions) Let $f_0 \in  L^\infty_{x,v}$ be nonnegative with compact support.
        Then there exists a global weak solution $f \in L^\infty(\R^+; L^\infty_{x,v}) \cap \mathscr{C}(\R^+; w-L^1_{x,v})$ to \eqref{eq:VM-kinetic} with initial data $f_0$.
        \item  (Strong solutions) Let $f_0 \in  W^{1,\infty}_{x,v}$ with compact support and such that $\rho_{f_0}(x)=\int f_0(x,v) \D v$ is bounded away from zero. Then there exists $T>0$ and a solution $f  \in  W^{1,\infty}([0,T]; W^{1,\infty}_{x,v})$ to \eqref{eq:VM-kinetic} with initial data $f_0$.
    \end{itemize}
\end{Thm}
Since then, the interest of this model has been justified several times. On the one hand, its monokinetic version has the form of a Hamiltonian ODE in the Wasserstein space, where the Hamiltonian is the Wasserstein distance to the Lebesgue measure~\cite{ambrosio2008hamiltonian}. On the other hand, it seems to emerge from pure noise as the optimality condition of a large deviation principle for a large population of indistinguishable stochastic particles~\cite{brenier2016double,ABB,levy2024monge,leonard2026monge}.

\begin{Rem}[Why restricting to periodic boundary conditions?]
We chose to work only with periodic boundary conditions. Indeed, a finer analysis of the Monge--Ampère equation reveals that the map in~\eqref{eq:transport_map} actually transports the density $\int f \D v$ to the Lebesgue measure \emph{restricted to its image}, which has to be a domain of Lebesgue measure $1$. Somehow, this targeted domain is another unknown of the system, always chosen more or less artificially. On the other hand on $\T^d$, there is only one domain of Lebesgue measure $1$ -- the whole $\T^d$ -- and hence the choice of $U$ is more natural.
\end{Rem}

The multiphasic formulation of the Vlasov--Monge--Ampère system reads as follows, with the ansatz $f=\mathcal{I}_\mu(\rhorho, \vv)$:
\begin{equation}
\label{eq:VM}
\left\{ 
\begin{gathered}
\partial_t \rho^\alpha + \Div ( \rho^\alpha v^\alpha) = 0,\\
\partial_t v^\alpha + (v^\alpha \cdot \nabla) v^\alpha = - \nabla U,\\
\mathrm{det}\left( \mathrm{I}+ \mathrm{D}^2_x U\right) = \int \rho^\alpha \D \mu(\alpha).
\end{gathered}
\right.
\end{equation}
To our knowledge, it has never been studied, except in the monokinetic case, see for instance~\cite{loeper2005quasi,tadmor2022critical}.

\section{Local well-posedness for the Vlasov--Poisson equations for electrons}\label{Subsec:VPelectrons}

In this section, we study the Cauchy problem for the Vlasov--Poisson system for electrons. 
Its algebraic structure is arguably the ``simplest'' considered in this work, as it is the only one where the force field depends linearly on the density. As a result, it  turns out fairly straightforward to apply the abstract results of Chapter~\ref{Part1-LWP}. We have chosen to describe local well-posedness results in Sobolev spaces in Section~\ref{sec:Sobolev-VPE}, and then in Besov spaces in Section~\ref{sec:Besov-VPE}.

\subsection{Sobolev local well-posedness}
\label{sec:Sobolev-VPE}

We start with the study of local wellposedness in Sobolev spaces. To obtain this result, we only need to check that it enters the general framework of Chapter~\ref{Part1-LWP}. In the statement that follows, we put forward two sets of assumptions under which local well-posedness can be obtained. They correspond to the two abstract settings from 
\begin{itemize}
    \item on the first hand, Theorem \ref{thm-LWPkinetic}  with the multiphasic decomposition from Example \ref{eq:example-decompoLinfty}; roughly speaking, we handle in this case  densities $\rhorho$ that are bounded with respect to $\alpha$;
    \item on the other hand, Theorem \ref{thm-LWPkinetic-L1} with the multiphasic decomposition from Example \ref{eq:example-decompoL1}; roughly speaking, we it accounts  in this case  for densities $\rhorho$ that are integrable with respect to $\alpha$.
\end{itemize}
The first assumption turns out to be lighter in terms of presentation. The second assumption is in some sense the most natural one  when considering measurable initial distributions functions, since, once translated as the level of the initial distribution function $f_0$ as in Example \ref{eq:example-decompoLinfty} (see Corollary \ref{coro-cin-VPe} below), it merely requires integrability in velocity. Both cases  allow for propagation of moments in velocity.

\begin{Thm}\label{thm-VPappli-VPE}
Let $k\in \N$ such that $k>1+d/2$. Let us suppose that the initial distribution function $f_0$ can be written as
$$f_0(x, \cdot)=\int_{I} \rho^\alpha_0(x) \otimes  \delta_{v=v^\alpha_0(x)} \, \mathrm{d}\mu(\alpha),$$
for some set of labels $(I, \mu)$ and some regular multiphasic distribution $(\rhorho_0, \vv_0)$. 
\begin{enumerate}

    \item If $(\rhorho_0, \vv_0) \in L^\infty_\alpha (H^{k-1}_x\cap \dot{H}^{-1}_x) \times \mathcal{H}^{k,p}$ with $p \in [1, \infty]$, there exist $T>0$ and a unique solution $(\rhorho, \vv)$ to the multiphasic system \eqref{eq:VP} with initial condition $(\rhorho_0, \vv_0)$, such that
    \begin{align*}
    \rhorho \in L^\infty([0,T]; L^\infty_\alpha (H^{k-1}_x \cap \dot{H}^{-1}_x )), \ \ \vv \in L^\infty([0,T]; \mathcal{H}^{k,p}).
\end{align*}

\item  Let $\boldsymbol\Theta: I \to \R_+$ and $ \boldsymbol\lambda : I \to \R^d$  be measurable maps satisfying $\langle \lambda^\alpha \rangle^p \lesssim \Theta(\alpha)$ with $p \in [0, \infty)$.   If $(\rhorho_0, \vv_0= \ww_0 + \boldsymbol{\lambda})$ is such that\footnote{Recall the notation of weighted $L^p$ spaces introduced at formula~\eqref{eq:def_weighted_Lp}.}  $(\rhorho_0, \ww_0 ) \in L^1_{\alpha,\Theta}(H^{k-1}_x \cap \dot{H}^{-1}_x) \times \mathcal{H}^{k,\infty}$, there exist $T>0$ and a unique solution $(\rhorho, \vv)$ to the multiphasic system \eqref{eq:VP} with initial condition $(\rhorho_0, \vv_0)$, such that
    \begin{align*}
    \rhorho \in L^\infty([0,T]; L^1_{\alpha,\Theta}(H^{k-1}_x \cap \dot{H}^{-1}_x )), \ \ \vv=\boldsymbol{\lambda} + \ww, \ \ \ww \in L^\infty([0,T]; \mathcal{H}^{k,\infty}).
\end{align*}
\end{enumerate}
 In both cases, if the maximal time of existence $T^\star$ is finite, then
\begin{equation*}
\int_0^{T^\star}   \lVe \nabla_x \vv \rVe_{L^{\upinfty}_\alpha L^{\upinfty}_x}  \, \D s = +\infty.
\end{equation*}
Moreover,
$$
f(t,x,\cdot) \vcentcolon= \int_{I} \rho^\alpha(t,x) \otimes  \delta_{v=v^\alpha(t,x)} \, \mathrm{d}\mu(\alpha)
$$
is the unique weak solution on $[0,T]$ of
 \eqref{eq:VPkin_sectionLWP}
  with initial data $f_0$, in the sense of Definition \ref{def:solution_coupled}. The distribution function $f$ is continuous in time, $H^{k-1}_x$--strongly in space and weakly in velocity in the sense of Definition~\ref{def:solution_Vlasov+timecontinuity}, 
and higher order moments in velocity enjoy the following regularity:
\begin{align*}
  &(t,x) \mapsto \int_{\R^d} \vert v \vert^j f(t,x,\D v) \in L^\infty([0,T]; H^{k-1}_x), \ \ j=0, \ldots,\lfloor p\rfloor.  \end{align*}

\end{Thm}

 In the Case 2. of Theorem~\ref{thm-VPappli-VPE}, it is possible to take $\boldsymbol\Theta \equiv 1$. The assumption with a weight is only useful when the initial distribution has high moments in velocity.

\begin{Rem}
    Note that in Theorem~\ref{thm-VPappli-VPE}, we have obtained an improved blow-up criterion compared to that of Theorem~\ref{thm:existence}.
\end{Rem}
Before proving Theorem~\ref{thm-VPappli-VPE}, let us first draw some concrete consequences
for the resolution of the Vlasov--Poisson system. They will be consequences of suitable choices of initial data that allow for a multiphasic decomposition at time $t=0$ that is compatible with one of the assumptions of Theorem~\ref{thm-VPappli-VPE}.

\begin{Cor}\label{coro-multiEuler-VPe}
Let $k\in \N$ such that $k>d/2$. Let $(m^n_0,v^n_0)_{n\in \N}$ be a family of density--velocity fields satisfying the bound
\begin{equation}
    \label{eq:coro1-vp}
 \| \langle n \rangle^\beta m^n_0 \|_{\ell^{\upinfty}_n (H^k_x \cap \dot{H}^{-1}_x)} +\| \langle n \rangle^{-\beta/p} v^n_0 \|_{\ell^p_n L^2_x} + \| \DDD v^n_0 \|_{\ell^{\upinfty}_n H^{k}_x} < +\infty,
\end{equation}
for some $\beta>1$ and $p \in [1,\infty)$. There exist $T>0$ and a unique weak solution $f$ to the Vlasov--Poisson system~\eqref{eq:VPkin_sectionLWP} on $[0,T]$, with initial condition 
$$f_0 = \sum_{n=0}^{+\infty} m^n_0 \otimes \delta_{v= v^n_0},$$ 
in the sense of Definition~\ref{def:solution_coupled}.  It is continuous in time, $H^{k}_x$--strongly in space and weakly in velocity, in the sense of Definition~\ref{def:solution_Vlasov+timecontinuity}. Moreover, higher order moments in velocity enjoy the following regularity:
\begin{align*}
  &(t,x) \mapsto \int_{\R^d} \vert v \vert^j f(t,x,\D v) \in L^\infty([0,T]; H^{k}_x), \ \ j=0, \ldots,\lfloor p\rfloor.
\end{align*}
\end{Cor}

\begin{proof}
We argue following Example~\ref{eq:example-decompoLinfty}, by taking 
$$I=\N, \ \ \mathrm{d}\mu(\alpha)=\sum_{n=0}^{+\infty} \frac{1}{\langle n \rangle^\beta} \delta_{\alpha=n},$$
and  
$$\rho^\alpha_0(x)={\langle \alpha \rangle^\beta} m_0^\alpha(x),
$$ 
so that
$$
f_0 = \sum_{n=0}^{+\infty} m^n_0(x) \otimes \delta_{v= v^n_0(x)}= \int_I \rho_0^\alpha \otimes \delta_{v= v^\alpha_0} \, \mathrm{d} \mu(\alpha). 
$$
According to~\eqref{eq:coro1-vp}, $\rhorho_0 \in L^\infty_\alpha (H^k_x \cap \dot{H}^{-1}_x)$ and $\vv_0 \in \mathcal{H}^{k+1,p}$. We can thus apply Theorem~\ref{thm-VPappli-VPE} with the first set of assumptions, which yields the result.
\end{proof}

\begin{Cor}\label{coro-cin-VPe}
     Let $k\in \N$ such that $k>d/2$. Let $f_0 \in  L^1_v (H^k_x\cap \dot H^{-1}_x)  $. There exist $T>0$ and a unique weak solution $f$ to the Vlasov--Poisson system~\eqref{eq:VPkin_sectionLWP} on $[0,T]$, with initial condition $f_0$ in the sense Definition~\ref{def:solution_coupled}.  It is continuous in time, $H^k_x$--strongly in space and weakly in velocity, in the sense of Definition~\ref{def:solution_Vlasov+timecontinuity}. 
If furthermore we assume that $f_0 \in  L^1_{v, \langle v \rangle^p}(H^k_x\cap \dot H^{-1}_x) $ for some $p \in [1, \infty)$, then higher order moments in velocity enjoy the following regularity:
\begin{align*}
  &(t,x) \mapsto \int_{\R^d} \vert v \vert^j f(t,x,\D v) \in L^\infty([0,T]; H^{k}_x), \ \ j=0, \ldots,\lfloor p\rfloor.
\end{align*}
\end{Cor}

\begin{proof}
We argue following Example~\ref{eq:example-decompoL1}, by taking 
$$I=\R^d, \ \ \mathrm{d}\mu(\alpha)=\frac{1}{\langle \alpha\rangle^\beta} \mathrm{d} \alpha,$$
with $\beta>d$, and  
$$\rho^\alpha_0(x)={\langle \alpha \rangle^\beta} f_0(x,\alpha), \ \ v^\alpha_0(x)=\alpha,
$$ 
so that
$$
f_0 =\int_I \rho_0^\alpha \otimes \delta_{v= v^\alpha_0} \, \mathrm{d} \mu(\alpha). 
$$
According to the assumptions on $f_0$, we have $\rhorho_0 \in L^1_\alpha (H^k_x \cap \dot{H}^{-1}_x)$. We can thus apply Theorem~\ref{thm-VPappli-VPE} with the second set of assumptions (with $\boldsymbol{\lambda}= \boldsymbol{\alpha}$, $\ww_0=0$, and $\boldsymbol{\Theta}=1$ or $\langle \boldsymbol{\alpha} \rangle^p$), which yields the result.

\end{proof}

As a final instance of application, we solve the Vlasov--Poisson for ``mixed'' data, as suggested in Motivation I in the introduction of the monograph.

\begin{Cor}\label{coro-cinVPe-mixed}  Let $k\in \N$ such that $k>d/2$.  Let $(m^n_0,v^n_0)_{n\in \N}$ be a family of density--velocity fields satisfying the bound
$$
 \|  m^n_0 \|_{\ell^{1}_n (H^k_x \cap \dot{H}^{-1}_x)} +\| v^n_0 - \lambda^n \|_{\ell^{\upinfty}_n H^{k+1}_x} +  < +\infty,
$$ 
for some sequence $\lambda: \N \to \R^d$, and $g_0 \in  L^1_v (H^k_x\cap \dot H^{-1}_x)$.
There exist $T>0$ and a unique weak solution $f$ to the Vlasov--Poisson system~\eqref{eq:VPkin_sectionLWP} on $[0,T]$, with initial condition 
$$
        f_0= \sum_{n=1}^\infty m^n_0(x)  \otimes \delta_{v=v^n_0(x)} + g_0,
$$
in the sense of Definition~\ref{def:solution_coupled}.  It is continuous in time, $H^k_x$--strongly in space and weakly in velocity, in the sense of Definition~\ref{def:solution_Vlasov+timecontinuity}.

\end{Cor}

\begin{proof}
 We write the same multiphasic decomposition separately for the Dirac part, and the regular part, as in the previous corollaries. Then we argue as in  Lemma~\ref{lem:conc} to concatenate the two measure spaces and the two multiphasic decompositions. The result follows from an application of Theorem~\ref{thm-VPappli-VPE} with the second set of assumption.
\end{proof}

The end of this subsection is dedicated to the proof of Theorem~\ref{thm-VPappli-VPE}.

\begin{proof}[Proof of Theorem~\ref{thm-VPappli-VPE}]

Thanks to Theorem~\ref{thm-LWPmultiphase-densityvelocity} in Chapter~\ref{Part1-LWP}, to prove Theorem~\ref{thm-VPappli-VPE} with the first set of assumptions on the initial condition, it is sufficient to check that the two main Assumptions \ref{ass:bound_high_regularity-rho} and \ref{ass:stability-rho} hold for the multiphasic system~\eqref{eq:VP}.
Similarly, to prove the theorem with the second set of assumptions on the initial condition, thanks to Theorem~\ref{thm-LWPmultiphase-densityvelocity-re} in Chapter~\ref{Part1-LWP}, it is sufficient to check that the two main Assumptions \ref{ass:bound_high_regularity-rho-re} and \ref{ass:stability-rho-re} are satisfied.

As these tasks are essentially the same, we focus on the verification of Assumptions \ref{ass:bound_high_regularity-rho} and \ref{ass:stability-rho} in the following.

\paragraph{Bound in high regularity.}
In the case of the Vlasov--Poisson equation, the force field does not depend on the phase $\alpha$, so controlling $\Vert \FF[\vv] \Vert_{L^1_T \mathcal{H}^{k,p}}$ is the same as controlling $\| \nabla U \|_{L^1_T  H^{k}_x}$. Since $U$ solves the Poisson equation, we have
\begin{align*}
   - \Delta \nabla U(t)=\nabla \int \rho^\alpha(t)  \D \mu(\alpha), \ \ t>0,
\end{align*}
and therefore, we have
\begin{equation*}
\| \nabla U(t) \|_{H^{k}_x} \lesssim  \int (\|\rho^\alpha(t) \|_{H^{k-1}_x}+\|\rho^\alpha(t) \|_{\dot H^{-1}_x}) \D \mu(\alpha).
\end{equation*}
This bound, which is pointwise in time, clearly implies 
that Assumption \ref{ass:bound_high_regularity-rho} is satisfied, by integrating with respect to time.

\paragraph{Stability at low regularity.}
We consider two families of velocity fields $\vv_1,\vv_2$ satisfying~\eqref{eq:uniform_bound_velocities}, and $\rhorho_1,\rhorho_2$ the corresponding solutions to the continuity equation. As the force field does not depend on $\alpha$, we just need to control the difference of force fields induced by two families of velocities in $L^2_x$. Let us consider two solutions $(U_i, \rhorho_i)$ (for $i=1,2$) of the Poisson equation in \eqref{eq:VP}. We have 
\begin{align*}
-\Delta (U_2 - U_1) = \int(\rho_2^{\alpha}(t) - \rho_1^{\alpha}(t)) \D \mu(\alpha),
\end{align*}
therefore 
\begin{align*}
\| \nabla U_2(t) - \nabla U_1(t) \|_{L^2_x} &= \left\| \int(\rho_2^{\alpha}(t) - \rho_1^{\alpha}(t)) \D \mu(\alpha) \right\|_{\dot H^{-1}_x} \\
&\lesssim  \int \left\| \rho_2^{\alpha}(t) - \rho_1^{\alpha}(t) \right\|_{\dot H^{-1}_x}\D \mu(\alpha)   .
\end{align*}
By integrating in time, we therefore obtain 
that Assumption~\ref{ass:stability-rho} is satisfied. 

\begin{Rem}
We note that an a priori estimate such as~\eqref{eq:uniform_bound_velocities-rho} is not needed in the Poisson case to obtain stabitity. This will not be the case for other models treated in this work.
\end{Rem}

\paragraph{Blow-up criterion.}
Let us now infer a refined continuation criterion for the multiphasic pressureless Euler--Poisson system for electrons. 

From (the proof of) Lemma \ref{LM:SobEstimRHO} on the continuity equation and Lemma \ref{LM:existence-vitesseSOBOLEV} on the velocity equation, we have 
\begin{align*}
\dfrac{\mathrm{d}}{\mathrm{d}t}\|\rho^\alpha\|_{H^{k-1}_x} &\lesssim \| \nabla v^\alpha \|_{L^{\upinfty}_x}\|\rho^\alpha\|_{H^{k-1}_x}  + 
\| \nabla v^\alpha\|_{H^{k-1}_x} \| \rho^\alpha\|_{L^{\upinfty}_x} ,
\end{align*}
together with
\begin{align*}
    \dfrac{\mathrm{d}}{\mathrm{d}t} \Vert v^\alpha \Vert_{L^2_x} &\lesssim \Vert \nabla v^\alpha \Vert_{L^{\upinfty}_x} \Vert v^\alpha \Vert_{L^2_x} + \Vert \nabla U \Vert_{L^2_x} \\
    &\lesssim \Vert \nabla v^\alpha \Vert_{L^{\upinfty}_x} \Vert v^\alpha \Vert_{L^2_x}+ \int \|\rho^\alpha \|_{H^{k-1}_x} \D \mu(\alpha),
\end{align*}
and 
\begin{align*}
	\dfrac{\mathrm{d}}{\mathrm{d}t}  \| \nabla v^\alpha \|_{H^{k-1}_x} \lesssim \| \nabla v^\alpha \|_{L^{\upinfty}_x} \| \nabla w \|_{H^{k-1}_x} + \int \|\rho^\alpha \|_{H^{k-1}_x} \D \mu(\alpha).
	\end{align*}
    In particular, we get
    \begin{align*}
    \dfrac{\mathrm{d}}{\mathrm{d}t}\|\rhorho \|_{L^{\upinfty}_\alpha H^{k-1}_x} &\lesssim \| \nabla \vv \|_{L^{\upinfty}_\alpha L^{\upinfty}_x}\|\rhorho\|_{L^{\upinfty}_\alpha H^{k-1}_x}  + 
\| \nabla \vv \|_{L^{\upinfty}_\alpha  H^{k-1}_x} \| \rhorho\|_{L^{\upinfty}_\alpha 
 L^{\upinfty}_x} ,
        \\
        \dfrac{\mathrm{d}}{\mathrm{d}t} \Vert \vv \Vert_{\mathcal{H}^{k,p}} &\lesssim \| \nabla \vv \|_{L^{\upinfty}_\alpha L^{\upinfty}_x} \Vert \vv \Vert_{\mathcal{H}^{k,p}}+ \|\rhorho \|_{L^{\upinfty}_\alpha H^{k-1}_x} .
	\end{align*}
Summing the two previous estimates yields 
\begin{multline*}
    \dfrac{\mathrm{d}}{\mathrm{d}t} \left(\|\rhorho \|_{L^{\upinfty}_\alpha H^{k-1}_x}+\Vert \vv \Vert_{\mathcal{H}^{k,p}} \right) 
    \lesssim \left(1+\| \rhorho\|_{L^{\upinfty}_\alpha 
 L^{\upinfty}_x}+\| \nabla \vv \|_{L^{\upinfty}_\alpha L^{\upinfty}_x}\right)\\
    \times \left( \|\rhorho \|_{L^{\upinfty}_\alpha H^{k-1}_x}+\Vert \vv \Vert_{\mathcal{H}^{k,p}} +\|\rhorho \|_{L^{\upinfty}_\alpha H^{k-1}_x}\right).
\end{multline*}
From the general pointwise estimate from \eqref{LM:pointwiseEstimRHO}, we infer
\begin{multline*}
    \dfrac{\mathrm{d}}{\mathrm{d}t} \left(\|\rhorho \|_{L^{\upinfty}_\alpha H^{k-1}_x}+\Vert \vv \Vert_{\mathcal{H}^{k,p}} \right) \\
    \lesssim \left(1+ \Vert \rhorho_0 \Vert_{L^{\upinfty}_\alpha L^{\upinfty}_x} \exp\left(\int_0^t \Vert \nabla \vv(\tau) \Vert_{L^{\upinfty}_\alpha L^{\upinfty}_x} \mathrm{d} \tau\right)+\| \nabla \vv \|_{L^{\upinfty}_\alpha L^{\upinfty}_x}\right) \\
    \times \left( \|\rhorho \|_{L^{\upinfty}_\alpha H^{k-1}_x}+\Vert \vv \Vert_{\mathcal{H}^{k,p}} +\|\rhorho \|_{L^{\upinfty}_\alpha H^{k-1}_x}\right).
\end{multline*}
Now using Gronwall lemma, we get the estimate for all $0 \leq t \leq T$
\begin{multline*}
 \Vert \vv(t) \Vert_{\mathcal{H}^{k,p}} 
\leq  \left (\|\rhorho_0 \|_{L^{\upinfty}_\alpha H^{k-1}_x}+\Vert \vv_0 \Vert_{\mathcal{H}^{k,p}} \right) \\
\times \exp\left(C T +  \Psi\left( \int_0^t  \Vert \nabla  \vv(\tau) \Vert_{L^{\upinfty}_\alpha L^{\upinfty}_x} \mathrm{d} \tau \right) \right),
\end{multline*}
for some $C>0$ and where $\Psi(r)=C\left(x +  \Vert \rhorho_0 \Vert_{L^{\upinfty}_\alpha L^{\upinfty}_x} e^x\right)$. This allows to conclude the proof of the blow-up criterion stated in Theorem \ref{thm-VPappli-VPE}.

\end{proof}

\subsection{Besov local well-posedness}
\label{sec:Besov-VPE}

The Besov local well-posedness result for the Vlasov--Poisson system for electrons reads as follows.

\begin{Thm}
\label{thm:vp-besov}
Let $s\in \R$, $q \in [2,+\infty], \, r \in [1,+\infty]$, with $s>1+d/q$, and $p \in [1, \infty]$. Let us suppose that the initial distribution function $f_0$ can be written as
$$f_0(x, \cdot)=\int_{I} \rho^\alpha_0(x) \otimes  \delta_{v=v^\alpha_0(x)} \, \mathrm{d}\mu(\alpha),$$
for some set of labels $(I, \mu)$ and some regular multiphasic distribution $(\rhorho_0, \vv_0)$.

\begin{enumerate}

\item If $(\rhorho_0, \vv_0) \in L^\infty_\alpha (B^{s-1}_{q,r} \cap {\dot{H}^{-1}_x}) \times \mathcal{B}^{(s,q,r),p} $, there exist $T>0$ and a unique solution $(\rhorho, \vv)$ to the multiphasic system \eqref{eq:VP} with initial condition $(\rhorho_0, \vv_0)$, such that
    \begin{align*}
    \rhorho \in L^\infty([0,T]; L^\infty_\alpha (B^{s-1}_{q,r}\cap {\dot{H}^{-1}_x}) ), \ \ \vv \in L^\infty([0,T]; \mathcal{B}^{(s,q,r),p}).
\end{align*}

\item  Let $\Theta: I \to \R_+$ and $ \lambda : I \to \R^d$  be measurable maps. If $(\rhorho_0, \vv_0= \ww_0 + \boldsymbol{\lambda})$ is such that $(\rhorho_0, \ww_0 ) \in L^1_{\alpha,\Theta} (B^{s-1}_{q,r} \cap \dot{H}^{-1}_x) \times \mathcal{B}^{(s,q,r),\infty}$, there exist $T>0$ and a unique solution $(\rhorho, \vv)$ to the multiphasic system \eqref{eq:VP} with initial condition $(\rhorho_0, \vv_0)$, such that
    \begin{align*}
    \rhorho \in L^\infty([0,T]; L^1_{\alpha,\Theta} (B^{s-1}_{q,r}\cap {\dot{H}^{-1}_x}) ), \ \vv=\boldsymbol{\lambda} + \ww, \, \ww\in L^\infty([0,T]; \mathcal{B}^{(s,q,r),\infty}).
\end{align*}

\end{enumerate}

 In both cases, if the maximal time of existence $T^\star$ is finite, then
\begin{equation}
\int_0^{T^\star}   \| \nabla_x \vv \|_{L^{\upinfty}_\alpha L^{\upinfty}_x}  \, \D s = +\infty.
\end{equation}

Moreover,
$$
f(t,x,\cdot) \vcentcolon= \int_{I} \rho^\alpha(t,x) \otimes  \delta_{v=v^\alpha(t,x)} \, \mathrm{d}\mu(\alpha)
$$
is the unique weak solution on $[0,T]$ 
to \eqref{eq:VPkin_sectionLWP}
  with initial data $f_0$, in the sense of Definition \ref{def:solution_coupled}.

\end{Thm}

From Theorem~\ref{thm:vp-besov}, we obtain as a corollary the following kinetic result.

\begin{Cor}\label{coro-cin-VPe-BESOV}
        Let $s\in \R$, $q \in [2,+\infty], \, r \in [1,+\infty]$, with $s>d/q$. 
   Let $f_0 \in  L^1_v (B^s_{q,r}\cap \dot H^{-1}_x)  $. There exist $T>0$ and a unique weak solution to the Vlasov--Poisson system~\eqref{eq:VPkin_sectionLWP} on $[0,T]$, with initial condition $f_0$ in the sense of Definition~\ref{def:solution_Vlasov}. 
\end{Cor}

 A few remarks are in order.

\begin{Rem}
Corollary~\ref{coro-cin-VPe-BESOV} in particular extends the former Sobolev result of Theorem~\ref{thm-VPappli-VPE}.

\begin{itemize}

\item  As $B^{s}_{2,2}$ identifies with $H^s_x$, the fractional $L^2_x$ Sobolev space
for any $s \in \R$, Corollary~\ref{coro-cin-VPe} extends to fractional Sobolev spaces, that is to say for initial data 
$$f_0(x,v) \in L^1_v (H^s_x \cap \dot{H}^{-1}_x),$$ 
for all $s>d/2$. 

\item As $B^\alpha_{\infty,\infty}$ identifies with $\mathscr{C}^{0,\alpha}_x$ for $\alpha \in (0,1)$, Corollary~\ref{coro-cin-VPe-BESOV} also in particular applies to initial data
$$f_0(x,v) \in L^1_v (\mathscr{C}^{0,\alpha}_x \cap \dot{H}^{-1}_x),$$
with  $\alpha \in (0,1)$. 

\item Similarly to Theorem~\ref{thm-VPappli-VPE} and Corollary~\ref{coro-cin-VPe}, one can also infer some suitable time continuity properties and regularity properties. The adaptation to the case of Besov spaces is left to the reader.

\end{itemize}

\end{Rem}

\begin{Rem}
\label{rem:infinitemass}
     We observe that Theorem~\ref{thm:vp-besov}  allows to consider initial data $f_0$ which are not necessarily $L^p_x$ integrable for $p \in [1,\infty)$.
In particular the initial data may be of infinite mass.
Therefore  Theorem~\ref{thm:vp-besov} complements the results of \cite{J00,CCM1,CCM2} which investigated this precise question for the Vlasov--Poisson equation in $\R^3$. More precisely, we note that in \cite{CCM2}, it is assumed that the initial condition $f_0$ satisfies
$$
0\leq f_0(x,v) \leq e^{-\lambda|v|^2}  g(|x|),
$$
for some $\lambda>0$ and $g$ is a bounded continuous function such that there exists $C>0$ and $\eps \in (1/15,1)$  for which for all $n \in \mathbb{Z}^3$,
$$
\int_{|n-x| \leq 1} g (|x|) \, \D x \leq C \frac{1}{|n|^{2+\eps}}.
$$
The result of \cite{CCM2} is the existence of an associated weak solution on some interval $[0,T]$, with uniqueness in an adequate class.
Let us therefore note that  \cite{CCM2} asks that the initial condition
decays like a Gaussian in velocity, an assumption that we do not require; on the other hand, \cite{CCM2} allows merely bounded data, while we require slightly more, namely Hölder continuity; finally \cite{CCM2} asks for a certain good control  of the space integral over balls with centers in $\Z^3$ (allowing for instance for piecewise constant initial densities, with a prescribed decay), while we ask for additional $\dot{H}^{-1}_x$ regularity. 

Finally, our approach  has the merit to show that this kind of infinite mass result is not restricted to the Vlasov--Poisson equation, but holds for a large class of Vlasov equations, satisfying the abstract assumptions of Theorem \ref{thm-LWPmultiphase-densityvelocity-BESOV}.

\end{Rem}

 \begin{proof}[Proof of Theorem \ref{thm:vp-besov}]  Let us focus on the first set of assumptions. We aim at applying  the abstract result of Theorem~\ref{thm-LWPmultiphase-densityvelocity-BESOV}, applied in the Vlasov--Poisson case. To this end, we must check 
   Assumptions \ref{ass:bound_high_regularity-rho-besov} and \ref{ass:stability-rho-besov}
   The verification of stability in low regularity, namely Assumption~\ref{ass:stability-rho-besov}, is exactly the same as in the Sobolev case (see the proof of Theorem~\ref{thm-VPappli-VPE}).
   The high regularity assumption follows from the general inequality
\begin{equation}
\label{eq:truc-besov}    
\| \nabla \Delta^{-1} \varphi \|_{B^{s}_{q,r}} \lesssim \|\varphi\|_{B^{s-1}_{q,r}} + \| \varphi \|_{\dot{H}^{-1}_x}.
\end{equation} 
This directly yields that
\begin{equation*}
\| \nabla U(t) \|_{B^s_{q,r}} \lesssim  \int (\|\rho^\alpha(t) \|_{B^s_{q,r}}+\|\rho^\alpha(t) \|_{\dot H^{-1}_x}) \D \mu(\alpha).
\end{equation*}
Indeed for $j\geq 1$, by \cite[Lemma 2.2]{BCD}, it holds
$$
\|  \nabla \Delta^{-1}  \Delta_j  \varphi \|_{L^q_x} \lesssim 2^{-j} \|\Delta_j  \varphi \|_{L^q_x}.
$$
The term $\| \varphi \|_{\dot{H}^{-1}_x}$ is useful  
to handle  the low frequency contribution. We write that $\varphi= \operatorname{div} \psi$, with $\psi \in L^2_x$ and $\| \psi \|_{L^2_x}= \| \varphi \|_{\dot{H}^{-1}_x}$. Since $\Delta_{-1}  \varphi$ is supported in frequency in a ball, by the Bernstein lemma (see \cite[Lemma 2.1]{BCD}) and the Plancherel theorem,
$$
\|  \nabla \partial_i \partial_{i'} \Delta^{-1}  \Delta_{-1}  \varphi \|_{L^q_x} \lesssim \|  \nabla \Delta^{-1}  \Delta_{-1}   \operatorname{div} \psi \|_{L^2_x}\lesssim \| \psi \|_{L^2_x}= \| \varphi \|_{\dot{H}^{-1}_x} .
$$
This concludes the proof of~\eqref{eq:truc-besov} and thus of the verification of   Assumption \ref{ass:bound_high_regularity-rho-besov}. 

For the second set of assumptions of the theorem, we can argue similarly, relying on Remark~\ref{rem:besov-lambda}.

\end{proof}

\subsection{The case of dynamical electrons and ions}
\label{sec:severalspecies}

As a matter of fact, it is straightforward to generalize all these results to the case of a plasma of electrons and (moving) ions, that is to say, in lieu of~\eqref{eq:VPkin_sectionLWP} and \eqref{eq:VP} , we consider
\begin{equation}
		\label{eq:VPkin-several}
		\left\{
		\begin{gathered}
     \partial_t f_{i} + v \cdot \nabla_x f_{i}- \nabla_x U\cdot \nabla_v f_{i} = 0,\\
		\partial_t f_e + v \cdot \nabla_x f_e + \nabla_x U\cdot \nabla_v f_e = 0,\\
		-\Delta_x U =   \int f_{i}\D v -  \int f_{e}\D v ,
        \end{gathered}
		\right.
		\end{equation}
where $f_i$ (resp. $f_e$) is the distribution function for ions (resp. electrons), with the neutrality constraint 
$$
\iint f_i|_{t=0} \D v \D x =\iint f_e|_{t=0} \D v \D x,
$$
which is propagated by the dynamics. It is also possible to consider several species of ions with barely any modification in the analysis but at the expense of heavier notation.
Consider two sets of labels $(I_i, \mu_i)$ (for ions) and $(I_e, \mu_e)$ (for electrons) and consider the extended set of labels 
$$I\vcentcolon = \{i\}\times I_i \cup \{e\}\times I_e.$$ 
Still denoting $\mu_i$ (resp. $\mu_e$) the extension of the measure on $\{i\}\times I_i$ (resp. $\{e\}\times I_e$), we endow $I$ with the canonical measure $\mu$ for the union. Also define the charge function $\operatorname{ch}: I \to \{-1,1\}$ as
$$
\operatorname{ch}(\alpha)\vcentcolon = \left\{
    \begin{array}{ll}
         \displaystyle 1, \qquad & \mbox{if }  \ \ \alpha \in \{i\}\times I_i, \\
        \displaystyle -1, \qquad & \mbox{if }  \ \ \alpha \in \{e\}\times I_e.
    \end{array}
    \right.
$$
Then the multiphasic version reads
\begin{equation}
\label{eq:VP-several}
\left\{ 
\begin{gathered}
\partial_t \rho^\alpha + \Div ( \rho^\alpha v^\alpha) = 0,\\
\partial_t v^\alpha + (v^\alpha \cdot \nabla) v^\alpha = - \operatorname{ch}(\alpha) \nabla U,
\\-\Delta_x U =  \int \operatorname{ch}(\alpha) \rho^\alpha  \D \mu(\alpha).
\end{gathered}
\right.
\end{equation}
It appears clearly that the multiphasic system~\eqref{eq:VP-several} is structurally very close to the multiphasic system~\eqref{eq:VP} for electrons. Therefore, all local well-posedness statements (Sobolev and Besov) can then be straightforwardly adapted to these systems.

Let us only write one slightly different result where interestingly the low frequency assumption (namely the $\dot{H}^{-1}$ control) only bears on the difference 
$f_i|_{t=0}-f_e|_{t=0}$. 
\begin{Thm}
         Let $k\in \N$ such that $k>d/2$ and let $m\in L^1_{v,\langle v\rangle}(\R^d)$ be a positive weight.  Let $f_{i,0}, f_{e,0} \in  L^\infty_{v,m^{-1}} H^k_x $ and $(f_{i,0} - f_{e,0}) \in  L^\infty_{v,m^{-1}} \dot{H}^{-1}_x$. There exist $T>0$ and a unique weak solution $(f_i,f_e)$ to the Vlasov--Poisson system~\eqref{eq:VPkin-several} on $[0,T]$, with initial condition $(f_{i,0}, f_{e,0})$ in the sense of Definition \ref{def:solution_coupled}.  Moreover, $f_{i,0}, f_{e,0}$ are continuous in time, $H^k_x$--strongly in space and weakly in velocity, in the sense of Definition~\ref{def:solution_Vlasov+timecontinuity}.
\end{Thm}

\begin{proof}
The proof uses the multiphasic formulation and is largely the same as that of Corollary~\ref{coro-multiEuler-VPe} or~\ref{coro-cin-VPe}. 
We rely on Case 1. of Theorem~\ref{thm-VPappli-VPE}. First of all, following Example \ref{eq:example-decompoLinfty}, since $f_{i,0}, f_{e,0} \in  L^\infty_{v,m^{-1}} H^k_x $ we note that we can take the same set of labels $(J,\nu)$ for each species. More precisely, we can take        
$$J=\R^d, \quad \D \nu(\alpha)={m(\alpha) }\mathrm{d}\alpha,$$
and 
$$ \rho_{0}^{(\ell,\alpha)}(x)=\frac{f_{\ell,0}(x,\alpha)}{m(\alpha)}, \ \ v^\alpha_0(x)=\alpha,$$
for $\ell= i,e$. Following the discussion before the statement of the theorem, we can then define $I\vcentcolon = \{i\}\times J \cup \{e\}\times J$.  We are thus led to the study of~\eqref{eq:VP-several} but will not repeat the full arguments.
However, there is only one significant difference in the propagation of the $\dot{H}^{-1}_x$ regularity, which we intend to explain.
The argument with $\dot{H}^{-1}_x$ regularity appears notably to control the $L^2$ norm in the ``bound in high regularity part'' of the proof of Theorem~\ref{thm-VPappli-VPE}. We have
\begin{align*}
-\Delta U  &= \int_I \operatorname{ch}(\alpha)\rho^{\alpha}  \D \mu(\alpha) =\int_J \left[\rho^{(i,\alpha')}-\rho^{(e,\alpha')}\right] \D \nu(\alpha'),
\end{align*}
and therefore, 
\begin{equation*}
\| \nabla U \|_{H^{k}_x} \lesssim  \int_J \left(\left\|\rho^{(i,\alpha')}-\rho^{(e,\alpha')} \right\|_{H^{k-1}_x}+ \left\|\rho^{(i,\alpha')}-\rho^{(e,\alpha')} \right\|_{\dot H^{-1}_x} \right) \D \nu(\alpha').
\end{equation*}
We therefore note that it is sufficient to control $\left\|\rhorho^{(i)}-\rhorho^{(e)} \right\|_{L^\infty_{\alpha'}\dot H^{-1}_x}$; but to this end, 
we can apply Lemma~\ref{LM:estimStabContinuity} which yields
	\begin{multline}
	\label{eq:rho_H-1-multi}
	\left\| \rho^{(i,\alpha')}-\rho^{(e,\alpha')} \right\|_{L^\infty_{\alpha'} \dot H^{-1}_x} \leq \exp\left( K \int_0^t \| \nabla \vv(s) \|_{L^{\upinfty}_\alpha L^{\upinfty}_x}  \right) \\ \left\{ \left\| \rho^{(i,\alpha')}_0-\rho^{(e,\alpha')}_0 \right\|_{L^\infty_{\alpha'}  \dot H^{-1}_x} + \left\| \rhorho_0 \right\|_{L^\infty_{\alpha'} L^{\upinfty}_x} + 2 \int_0^t \| \vv(s)\|_{L^{\upinfty}_\alpha  L^2_x} \D s \right\},
	\end{multline}
    and we conclude noting that 
    $$
    \left\| \rho^{(i,\alpha')}_0-\rho^{(e,\alpha')}_0\right\|_{L^\infty_{\alpha'}  \dot H^{-1}_x} \lesssim \left\| f_{i,0} - f_{e,0} \right\|_{L^\infty_{v,m^{-1}} \dot{H}^{-1}_x},
    $$
  hence the proof of the theorem.
\end{proof}

\section{Local well-posedness for the Vlasov--Poisson equations for ions}\label{Subsec:VPions}
In what follows, we restrict to the dimensions $d=1,2,3$.
We obtain, for the Vlasov--Poisson system for ions, the analogue of local well-posedness results in Sobolev for the Vlasov--Poisson for electrons. Contrary to the latter, we restrict to integer Sobolev indices and do not treat Besov regularity, mainly for the sake of simplicity. Furthermore, we only treat the case of the spatial periodic torus to lighten the exposition.

\begin{Thm}\label{thm-VPappli-VPI}
Let $k\in \N$ such that $k>1+d/2$ and $p \in [1, \infty]$. Let us suppose that the initial distribution function $f_0$ can be written as
$$f_0(x, \cdot)=\int_{I} \rho^\alpha_0(x) \otimes  \delta_{v=v^\alpha_0(x)} \, \mathrm{d}\mu(\alpha),$$
for some set of labels $(I, \mu)$ and some regular multiphasic distribution $(\rhorho_0, \vv_0)$.

\begin{enumerate}

  \item If $(\rhorho_0, \vv_0) \in L^\infty_\alpha (H^{k-1}_x \times \mathcal{H}^{k,p})$, there exist $T>0$ and a unique solution $(\rhorho, \vv)$ to the multiphasic system \eqref{eq:VPions} with initial condition $(\rhorho_0, \vv_0)$, such that
    \begin{align*}
    \rhorho \in L^\infty([0,T]; L^\infty_\alpha (H^{k-1}_x )), \ \ \vv \in L^\infty([0,T]; \mathcal{H}^{k,p}).
\end{align*}

\item  Let $\Theta: I \to \R_+$ and $ \lambda : I \to \R^d$  be a measurable maps satisfying $\langle \lambda^\alpha \rangle^p \lesssim \Theta(\alpha)$ with $p \in [1, \infty)$.  If $(\rhorho_0, \vv_0= \ww_0 + \boldsymbol{\lambda})$ is such that  $(\rhorho_0, \ww_0 ) \in L^1_{\alpha,\Theta}(H^{k-1}_x ) \times \mathcal{H}^{k,\infty}$, there exist $T>0$ and a unique solution $(\rhorho, \vv)$ to the multiphasic system \eqref{eq:VPions} with initial condition $(\rhorho_0, \vv_0)$, such that
    \begin{align*}
    \rhorho \in L^\infty([0,T]; L^1_{\alpha,\Theta} (H^{k-1}_x  )), \ \ \vv=\boldsymbol{\lambda} + \ww, \ \ \ww \in L^\infty([0,T]; \mathcal{H}^{k,\infty}).
\end{align*}

\end{enumerate}
In  both cases,
$$
f(t,x,\cdot) \vcentcolon= \int_{I} \rho^\alpha(t,x) \otimes  \delta_{v=v^\alpha(t,x)} \, \mathrm{d}\mu(\alpha)
$$
is the unique weak solution on $[0,T]$ of \eqref{eq:VPions}
  with initial data $f_0$, in the sense of Definition \ref{def:solution_coupled}. Moreover, $f$ is continuous in time, $H^{k-1}_x$--strongly in space and weakly in velocity, in the sense of Definition~\ref{def:solution_Vlasov+timecontinuity}, 
and  higher order moments in velocity enjoy the following regularity:
\begin{align*}
  &(t,x) \mapsto \int_{\R^d} \vert v \vert^j f(t,x,\D v) \in L^\infty([0,T]; H^{k-1}_x), \ \ j=0, \ldots,\lfloor p\rfloor.  \end{align*}
\end{Thm}

Exactly as in the case of the Vlasov-Poisson system for electrons from Section \ref{Subsec:VPelectrons}, we obtain the following set of results at the kinetic level.
\begin{Cor}\label{coro-cin-VPi}
Corollaries \ref{coro-multiEuler-VPe}--\ref{coro-cin-VPe}--\ref{coro-cin-VPi} hold true for the Vlasov-Poisson system \eqref{eq:VPions-kinetic}.
\end{Cor}

\begin{proof}
    [Proof of Theorem~\ref{thm-VPappli-VPI}]
Thanks to Theorem~\ref{thm-LWPmultiphase-densityvelocity} and Theorem~\ref{thm-LWPmultiphase-densityvelocity-re} in Chapter~\ref{Part1-LWP}, it is sufficient to check the main Assumptions \ref{ass:bound_high_regularity-rho} and \ref{ass:stability-rho} for the Vlasov--Poisson system for electrons (the verification of Assumptions  \ref{ass:bound_high_regularity-rho-re} and \ref{ass:stability-rho-re} being almost identical).

\paragraph{Bound in high regularity.} For a given $(\rho^\alpha)$, let us assume that we have a solution $U \in H^1(\T^3)$ of the elliptic equation 
$$e^U-\Delta U = \int \rho^\alpha  \D \mu(\alpha).$$
The existence and uniqueness of a such a weak solution $U \in H^1_x$ of the elliptic equation (with a given $\rhorho \in L^\infty_\alpha L^2_x$) follows for instance from \cite[Propositions 3.1 and 3.5]{GPI}, that comes from a minimization procedure in $H^1$ (and yields a bound like $\Vert U \Vert_{H^1} \lesssim 1+ \Vert \rhorho \Vert_{L^{\upinfty}_\alpha L^2_x}$). The proof then proceeds in two steps:
\begin{itemize}
    \item first, we upgrade the regularity of $U$ and show that $U \in L^\infty_T H^2_x$ if $\rhorho \in L^\infty_\alpha L^2_x$. This follows from a standard argument with the method of finite differences \textit{à la} Nirenberg, using crucially the fact that $(x-y)(e^x-e^y) \geq 0$ for all $x,y \in \R$ to get rid of the contribution of the semilinear term (it can also be obtained by a direct energy estimate when one tests the equation against $\Delta U$ itself): we obtain the
    \begin{align*}
        \Vert U \Vert_{H^2_x} \lesssim 1+\Vert \rhorho \Vert_{L^{\upinfty}_\alpha L^2_x}.
    \end{align*}
    \item second, we obtain a bound in high regularity. Since $H^2_x(\T^3) \hookrightarrow L^\infty_x(\T^3)$, we already know that $U \in L^\infty_T L^\infty_x$ and therefore, by the composition rule in Sobolev spaces from Proposition \ref{prop-Sobolev}, we infer that for all $s>0$
    \begin{align*}
        \Vert e^U \Vert_{H^s_x} \leq C_s(\Vert U \Vert_{L^{\upinfty}_x})(1+\Vert U \Vert_{H^s}),
    \end{align*}
    where $C_s : \R^+ \rightarrow \R^+$ is a continuous increasing function. Writing the elliptic equation as
    \begin{align*}
        -\Delta  U= -e^U+  \int  \rho^\alpha \, \mathrm{d}\mu(\alpha),
    \end{align*}
and seeing the right-hand side as a source term, the standard elliptic theory and the previous bounds yield for $k>2$
    \begin{align*}
        \Vert \nabla U \Vert_{H^k_x} &\lesssim \Vert e^U \Vert_{H^{k-1}_x}+ \Vert \rhorho \Vert_{L^{\upinfty}_\alpha H^{k-1}_x} \\
        & \lesssim C_{k-1}(\Vert\rhorho \Vert_{L^{\upinfty}_\alpha L^2_x})(1+\Vert U \Vert_{H^{k-1}_x})+\Vert \rhorho\Vert_{L^{\upinfty}_\alpha H^{k-1}_x} \\
        &  \lesssim C_{k-1}(\Vert \rhorho \Vert_{L^{\upinfty}_\alpha L^2_x})(1+\Vert \nabla U \Vert_{H^{k-2}_x})+\Vert \rhorho \Vert_{L^{\upinfty}_\alpha H^{k-1}_x},
    \end{align*}
thanks to the initial bound for $U$ in $H^1$.  An induction argument then allows to conclude that Assumption \ref{ass:bound_high_regularity-rho} is satisfied.
\end{itemize}

\paragraph{Stability at low regularity.}
Let us now consider two solutions $(U^i, (\rho_i^{\alpha}))$ (for $i=1,2$) of the Poisson equation, that is
\begin{align*}
    e^{U_i}-\Delta U_i = \int \rho_i^\alpha  \D \mu(\alpha), \ \ i=1,2.
\end{align*}
We have
\begin{align*}
e^{U_2}-e^{U_1}-\Delta (U_2 - U_1) = \int(\rho_2^{\alpha}(t) - \rho_1^{\alpha}(t)) \D \mu(\alpha),
\end{align*}
therefore, by using again that $(x-y)(e^x-e^y) \geq 0$ for all $x,y \in \R$, we get by testing this identity with $U_2-U_1$: 
\begin{equation*}
\| \nabla U_2(t) - \nabla U_1(t) \|_{L^2_x} \leq   \left\| \int(\rho_2^{\alpha}(t) - \rho_1^{\alpha}(t)) \right\|_{\dot H^{-1}_x}\D \mu(\alpha) .
\end{equation*}
Hence Assumption \ref{ass:stability-rho} holds. The proof of Theorem~\ref{thm-VPappli-VPI} is finally complete.
\end{proof}

\section{Local well-posedness for Vlasov--Monge--Ampère equations}\label{Subsec:VMongeAMpere}

This section is devoted to the study of the Vlasov--Monge--Ampère system. We restrict in this section exclusively to the spatial domain $\T^d$. 

\subsection{Application of Theorem~\ref{thm:existence} to the Vlasov--Monge--Ampère system}

\begin{Thm}\label{thm-VPappli-VM}
Let $k\in \N$ be such that $k>2+d/2$ and $p \in [1, \infty]$. There exists $\eps_0>0$ only depending on $d$ and $k$ with the following property. Let us suppose that the initial distribution function $f_0$ can be written as
$$f_0(x, \cdot)=\int_{I} \rho^\alpha_0(x) \otimes  \delta_{v=v^\alpha_0(x)} \, \mathrm{d}\mu(\alpha),$$
for some set of labels $(I, \mu)$ and some regular multiphasic distribution $(\rhorho_0, \vv_0) \in L^1_\alpha H^{k-1}_x \times \mathcal{H}^{k,p}$ satisfying
\begin{equation*}
    \int \| \rho^\alpha_0 - 1 \|_{H^{k-1}_x} \D \mu(\alpha) \leq \eps_0.
\end{equation*}
Then there exist a time $T>0$ and a unique weak solution $f$ on $[0,T]$  to \eqref{eq:VM-kinetic},
 with initial data $f_0$, in the sense of Definition \ref{def:solution_coupled}. 
\end{Thm}
This theorem is a direct consequence of the following multiphasic statement.
\begin{Thm}
\label{thm-VPappli-multi-VM}
Let $(I, \mu)$ be a set of labels. Let $k\in \N$ be such that $k>2+d/2$ and $p \in [1, \infty]$. There exists $\eps_0>0$ only depending on $k$ and $p$ such that the following holds.  If $(\rhorho_0, \vv_0) \in L^1_\alpha H^{k-1}_x  \times \mathcal{H}^{k,p}$ is a regular multiphasic distribution, and if 
\begin{equation*}
    \int \| \rho^\alpha_0 - 1 \|_{H^{k-1}_x} \D \mu(\alpha) \leq \eps_0,
\end{equation*}
there exist $T>0$ and a unique solution $(\rhorho, \vv)$ to the multiphasic system \eqref{eq:VM} with
$$ \rhorho \in L^\infty([0,T]; L^1_\alpha (H^{k-1}_x \cap \dot{H}^{-1}_x )), \ \ \vv \in L^\infty([0,T]; \mathcal{H}^{k,p}).$$
\end{Thm}
We also have
\begin{Thm}
\label{thm-VPappli2-VM}
Theorems~\ref{thm-VPappli-VM} and \ref{thm-VPappli-multi-VM} hold for a regular multiphasic distribution $(\rhorho_0, \ww_0 + \boldsymbol{\lambda})$ with $(\rhorho_0, \ww_0 ) \in L^1_\alpha (H^{k-1}_x \cap \dot{H}^{-1}_x) \times \mathcal{H}^{k,\infty}$.  
\end{Thm}
 
\begin{Cor}\label{coro-cin-VMA}
     Let $k\in \N$ such that $k>1+d/2$. There exists $\eps_0>0$ only depending on $d$ and $k$ with the following property. Let $f_0 \in  L^1_v H^k_x $ such that  
     \begin{equation}
     \label{eq:smallness-VMA}
         \int  \left\| {f_0(\cdot,v)} -\left(\int f_0(x',v) \D x' \right) \right\|_{H^k_x}  \D v \leq \eps_0.
     \end{equation}
     There exist $T>0$ and a unique weak solution $f$ to the Vlasov--Monge-Ampère system~\eqref{eq:VM-kinetic} on $[0,T]$, with initial condition $f_0$ in the sense of Definition~\ref{def:solution_coupled}. It is continuous in time, $H^k_x$--strongly in space and weakly in velocity, Definition~\ref{def:solution_Vlasov+timecontinuity}.
     If furthermore we assume that $f_0 \in  L^1_{v, \langle v \rangle^p}(H^k_x\cap \dot H^{-1}_x) $ for some $p \in [1, \infty)$, then higher order moments in velocity enjoy the following regularity:
\begin{align*}
  &(t,x) \mapsto \int_{\R^d} \vert v \vert^j f(t,x,\D v) \in L^\infty([0,T]; H^{k}_x), \ \ j=0, \ldots,\lfloor p\rfloor.
\end{align*}
\end{Cor}
\begin{proof}
    We apply Theorem \ref{thm-VPappli-VM} with a suitable choice of initial data. We set
   \begin{align*}
       I=\R^d, \ \ \mathrm{d}\mu(\alpha)=\left(\int f_0(x',\alpha) \D x' \right) \mathrm{d}\alpha,
   \end{align*}
   and
   \begin{align*}
       \rho^\alpha_0(x)=\frac{f(x,\alpha)}{\int f_0(x',\alpha) \D x'}, \ \ v^\alpha(x)=\alpha.
   \end{align*}
Then, in view of \eqref{eq:smallness-VMA}, one has 
\begin{align*}
    \int \| \rho^\alpha_0 - 1 \|_{H^{k}_x} \D \mu(\alpha)=\int  \left\| {f_0(x,v)} -\left(\int f_0(x,v) \D x \right) \right\|_{H^k_x}  \D v  \leq \eps_0,
\end{align*}
together with $\rhorho_0 \in L^1_\alpha H^k_x$. 
\end{proof}
\begin{Rem}
Note that the assumption \eqref{eq:smallness-VMA} on the initial density is stronger than the lower bound requirement from Theorem \ref{thmBrenierLoeper-VMA}: indeed, assuming that $\int f_0(x,v) \D x \D v=1$, one has \begin{align*}
    \Vert \rho_{f_0}-1 \Vert_{L^\infty_x} \lesssim \int  \left\| {f_0(x,v)} -\left(\int f_0(x,v) \D x \right) \right\|_{H^k_x} \, \mathrm{d}v.
\end{align*}
There is presumably a way to replace our smallness assumption by a below bound for the total density, which is the crucial point to get elliptic estimates in the Monge-Ampère equation. Such an approach would be closer to the one of Brenier and Loeper~\cite{BreLoe}, but we did not pursue in this direction.
\end{Rem}
Again, thanks to Theorem~\ref{thm-LWPmultiphase-densityvelocity} and Theorem~\ref{thm-LWPmultiphase-densityvelocity-re} in Chapter~\ref{Part1-LWP}, it is sufficient to check the main Assumptions \ref{ass:bound_high_regularity-rho} and \ref{ass:stability-rho} for the Vlasov--Monge--Ampère equation (the verification of Assumptions  \ref{ass:bound_high_regularity-rho-re} and \ref{ass:stability-rho-re} being identical). In fact, for reasons that will become clear in a few lines, we will instead check Assumptions~\ref{ass:bound_high_regularity} and~\ref{ass:stability}, and apply Theorem~\ref{thm:existence}.

\paragraph{Bounds in high regularity for the Monge-Ampère equation.}
We focus here on the Monge--Ampère equation of unknown $u: \T^d \to \R$, with source term $g$, namely
 \begin{equation}
 \label{eq:MA}
 \left\{
 \begin{gathered}
 \det (\Id + \DDD^2 u ) = 1 + g,\\
 \Id + \DDD^2 u \geq 0.
 \end{gathered}
 \right.
 \end{equation}

The following theorem of existence and regularity for a continuous source is well-known since the founding results of Caffarelli~\cite{caffarelli1990interior} (here written in the periodic domain, in the same spirit as~\cite{loeper2005quasi}).
 \begin{Thm}
 	\label{thm:caf}
 	Suppose that $g \in \mathscr{C}^\theta(\T^d)$ for some $\theta \in (0,1)$ and satisfies $\| g \|_{\mathscr{C}^\alpha} < 1$. Then there exists a unique solution $u \in \mathscr{C}^{2,\theta}(\T^d)$ to~\eqref{eq:MA} up to an additive constant, and $\|\DDD^2 u\|_{\mathscr{C}^{\theta}}$ is bounded by a constant only depending on the dimension and $\|g\|_{\mathscr{C}^\theta}$.
 \end{Thm}


Therefore, at \emph{low} regularity, the Monge--Ampère equation induces a gain of two derivatives. In this theorem, the assumption $\| g \|_{\mathscr{C}^\alpha} < 1$ is not really here to ensure integrability for $g$, but rather to ensure ellipticity for the Monge--Ampère equation. For this reason, when applying the results of this section to the study of the Vlasov--Monge--Ampère system, we will need to restrict ourselves to the case when the density remains close to $1$.

We need a result stating that  the Monge--Ampère equation also induces a gain of two derivatives at \emph{high} regularity, and to be coherent with the rest of this work, in Sobolev spaces. That is, we need to show that if  $g \in H^k_x$ then $u \in H^{k+2}_x$ for $k$ large, under the ellipticity condition $\| g \|_{\mathscr{C}^\alpha} < 1$.

 Observe that by Sobolev embeddings, choosing $k_0>d/2$, there is $\eps_0$ such that $\|g \|_{H^{k_0}_x} \leq \eps_0$ ensures that $\| g \|_{\mathscr{C}^\alpha} < 1$ for some $\alpha$ uniformly. We choose to write the ellipticity condition in that way in order to ensure that the constant in Theorem~\ref{thm:caf} only depends on the dimension. 

\begin{Thm}
	\label{thm:boundMA}
Let $k_0 > d/2$. There exists $\eps_0>0$ such that for all $k \in \N^*$, there exists a polynomial $F_k = F_k[d,k_0, \eps_0]$ cancelling at $0$ and such that for all $g \in H^k(\T^d)$ with $\| g \|_{H^{k_0}} \leq \eps_0$, the solution $u$ to \eqref{eq:MA} satisfies
\begin{equation}
\label{eq:estim_MA}
\| u \|_{H^{k+2}_x} \leq F_k( \| g \|_{H^k_x}).
\end{equation}
\end{Thm}
The proof is postponed to the end of the section.

\bigskip

The additional difficulty with respect to the Vlasov--Poisson equation is that in order to use Theorem~\ref{thm:boundMA}, we need to ensure that the total density is close to one at some fixed regularity $H^{k_0}$ with $k_0>d/2$. We use the flexibility granted by the additional parameter $X$ for this purpose. This is why it is more convenient to try to apply Theorem~\ref{thm:existence} instead of Theorem~\ref{thm-LWPmultiphase-densityvelocity}, and thus to rather check Assumptions \ref{ass:bound_high_regularity} and \ref{ass:stability} in this setting.

Let us fix $k_0 > 1+ d/2$. We set
\begin{equation*}
    \mathcal X \vcentcolon= \Big\{ \rhorho_0 \in L^1_\alpha H^{k_0-1}_x \mbox{ such that } \int \| \rho_0^\alpha - 1  \|_{H^{k_0-1}_x} \D \mu(\alpha) < \eps_0 \Big\},
\end{equation*}
where $\eps_0$ is the parameter given in Theorem~\ref{thm:boundMA} for $k_0 - 1$ in place of $k_0$. For $\rhorho_0 \in \mathcal X$, we shall denote, for $k \geq k_0$,
\begin{equation*}
    N_k(\rhorho_0) \vcentcolon= \| \rhorho_0 - 1 \|_{L^{\upinfty}_\alpha H^{k-1}_x}. 
\end{equation*}
Fix  $k \geq k_0$ and $p \in [1,+\infty]$. 

\paragraph{Bound in high regularity for the Vlasov--Monge--Ampère system.} With these considerations at hand, let us check that the bound at high regularity -- namely, Assumption \ref{ass:bound_high_regularity} -- holds for the multiphasic Vlasov--Monge--Ampère equation \eqref{eq:VM}.

Let $T>0$,  $\vv \in L^1([0,T]; \mathcal H^{k,p})$, and $\rhorho_0 \in \mathcal{X}$. Let $\rhorho \in L^\infty([0,T];L^\infty_\alpha H^{k-1}_x)$ be the associated solution to the continuity equation on $[0,T]$, with vector field $\vv$ and initial condition $\rhorho_0$. 
Let  $U$ be the associated solution to the Monge--Ampère equation, that is
	 \begin{equation*}
 	\det (\Id + \DDD^2 U ) = 1+ \left(\int_I \rho^\alpha \, \mathrm{d}\mu(\alpha) -1\right).
 	\end{equation*}
Using the estimate~\eqref{estim:rhoalphaSob-integ} in Lemma~\ref{LM:SobEstimRHO} for the continuity equation, we find for all $t \in [0,T]$
  \begin{equation}
      \label{eq:rho-1VMA}
    \begin{aligned}
       &\left\| \int \rho^\alpha (t)\D \mu(\alpha) - 1 \right\|_{H^{k-1}_x} \\
       &\qquad \leq \int \| \rho^\alpha(t) - 1\|_{H^{k-1}_x} \D \mu(\alpha) \\
       &\qquad \leq \left(1 +  \int \| \rho^\alpha_0 - 1 \|_{H^{k-1}_x} \D \mu(\alpha) \right)\exp\left( K \int_0^t\Vert\vv(s)\Vert_{\mathcal{H}^{k,p}} \D s \right)  - 1\\
      &\qquad  \leq (1 + \mathsf N_{k}(\rhorho_0)) \exp\left( K \int_0^T \Vert\vv(s)\Vert_{\mathcal{H}^{k,p}} \, \mathrm{d}s \right) - 1.
    \end{aligned}
      \end{equation}
Therefore, if $(1 + N_{k}(\rhorho_0)) \exp\left( K \int_0^T \Vert\vv(s)\Vert_{\mathcal{H}^{k,p}} \right) <1+ \eps_0$, we can apply Theorem~\ref{thm:boundMA} and this yields
\begin{equation}
\label{eq:ma-F}
    \int_0^T \Vert \nabla U(t) \Vert_{H^{k}_x} \D t \leq \F\left(T, \int_0^T \Vert v(t)\Vert_{\mathcal{H}^{k,p}} \D t, \mathsf N_{k}(\rhorho_0)\right),
\end{equation}
       where the estimate function $\F$ is defined as
       \begin{equation*}
    \F(T,A,P) \vcentcolon= \left\{ \begin{aligned}
    &T F_{k}((1+P)\exp(KA)-1), &&\mbox{if } (1+P)\exp(KA) \leq 1 + \eps_0, \\
    & + \infty, && \mbox{else}.
    \end{aligned} \right.
\end{equation*}
in which $F_k$ comes from Theorem~\ref{thm:boundMA}. This proves that Assumption \ref{ass:bound_high_regularity} is satisfied.

\paragraph{Stability at low regularity for the Vlasov--Monge--Ampère system.}
We rely this time on a stability result for the Monge--Ampère equation \eqref{eq:MA}.
\begin{Thm}\label{stabMA}
 	Let $k_1 > d/2 + 1$. There exists $\eps_1>0$ such that whenever 
 	\begin{equation*}
 	\| g \|_{H^{k_1}_x} \leq \eps_1 \qquad \mbox{and} \qquad \| h \|_{H^{k_1}_x} \leq \eps_1,
 	\end{equation*}
 	if $u$ and $v$ satisfy respectively
 	 \begin{equation*}
 	\det (\Id + \DDD^2 u ) = 1 + g \qquad \mbox{and} \qquad \det (\Id + \DDD^2 v ) = 1 + h,
 	\end{equation*}
 	then
 	\begin{equation*}
\| \nabla u - \nabla v \|_{L^2_x} \leq 2 \| g-h \|_{\dot H^{-1}_x}.
 	\end{equation*}
 \end{Thm}
Again, the proof is postponed to the end of the section. 

\bigskip

Equipped with this stability result, we can check that Assumption \ref{ass:stability} is satisfied for the multiphasic Vlasov--Monge--Ampère equation \eqref{eq:VM}. 
Let $R>0$. Consider two families of velocity fields $\vv_1,\vv_2 \in L^\infty([0,T]; \mathcal H^{k,p})$, and $\rhorho_{1,0}, \rhorho_{2,0} \in \mathcal{X}$
such that~\eqref{eq:uniform_bound_velocities} is satisfied, that is to say
\begin{equation*}
	 \max\left( \sup_{t \in [0,T]} \Vert \vv_1(t)\Vert_{\mathcal{H}^{k_1,p}},\sup_{t \in [0,T]} \Vert \vv_2(t) \Vert_{\mathcal{H}^{k_1,p}}, \mathsf N_{k_1}(\rhorho_{1,0}), \mathsf N_{k_1}(\rhorho_{2,0})  \right) \leq R.
		\end{equation*}
Let $\rhorho_1,\rhorho_2$ be the corresponding solutions to the continuity equation, with velocity field $\vv_1,\vv_2$ and initial condition $\rhorho_{1_0},\rhorho_{2_0}$. We need to control the difference of the force fields, induced by the two families of velocities, in $L^2_x$. Let us consider the two associated solutions $U_1$ and $U_2$ of the Monge--Ampère equation \eqref{eq:MA}, that is to say
\begin{align*}
    \mathrm{det}\left( \mathrm{I}+ \mathrm{D}^2_x U_i\right) = \int \rho_i^{\alpha} \D \mu(\alpha), \ \ i=1,2.
\end{align*}
By~\eqref{eq:rho-1VMA} for $k=k_1+1$, there holds
     \begin{equation}
     \label{eq:small-MA}
\left\| \int \rho^{\alpha}_i \D \mu(\alpha)-1 \right\|_{H^{k_1}_x} \leq \eps_1, \ \ i=1,2,
	 \end{equation}
if  $\max_{i=1,2} (1 + N_{k_1}(\rhorho_{i,0})) \exp\left( K \int_0^T \Vert \vv_i(s)\Vert_{\mathcal{H}^{k,p}} \right) <1+ \eps_1$.
     Thus, according to  Theorem \ref{stabMA}, we obtain that 
$$
\| \nabla U_1 - \nabla U_2 \|_{L^2_x} \leq 2 \left\| \int (\rho_1^{\alpha} -\rho_2^\alpha)\D \mu(\alpha) \right\|_{\dot H^{-1}_x}\leq 2  \int \left\|\rho^{\alpha}_1 -\rho^{\alpha}_2 \right\|_{\dot H^{-1}_x} \D \mu(\alpha).
$$
By integrating in time, applying Lemma~\ref{LM:estimStabContinuity}, we finally obtain a stability estimate of the form
\begin{multline*}
  \int_0^T \| \nabla U_1(t) - \nabla U_2(t) \|_{L^2_x} \, \mathrm{d} t \\
  \leq \G_R  \left(T, \int_0^T \| \vv_1(t)-\vv_2(t)\|_{L^p_\alpha L^2_x} \, \mathrm{d}t,  \| \rhorho_{1,0}- \rhorho_{2,0}\|_{L^1_\alpha \dot{H}^{-1}_x}\right),
\end{multline*}
where
    \begin{equation*}
    \G_R(T,A,P) \vcentcolon= \left\{ \begin{aligned}
    &\exp(2 K TR) \left(P + R A\right), &&\mbox{if } (1+ R)\exp(T R)<1+\eps_1, \\
    & + \infty, && \mbox{else}.
    \end{aligned} \right.
\end{equation*}
Hence Assumption \ref{ass:stability} holds,  and this concludes the proof.

\subsection{Proofs of regularity and stability for the Monge--Ampère equation}\label{Appendix-MongeAmpère}

\subsubsection{Proof of Theorem~\ref{thm:boundMA}}

Our proof is an adaptation of Loeper's proof of~\cite[Theorem~3.4]{loeper2005quasi}.
\begin{proof}
	This theorem is based on an \emph{a priori} estimate. What we mean is that it is enough to prove it for any smooth $g$ and $u$. Also, we fix $d$, $k_0>d/2$ and we pick once for all $\eps_0$ such that $\| g \|_{H^{k_0}_x} \leq \eps_0$ implies a uniform bound $\| g \|_{\mathscr{C}^\theta} < 1$ for some $\theta \in (0,1)$. Therefore, we do not refer to the dependence in $d$, $k_0$ and $\eps_0$ in the constants to come.
	
	The main ingredient of the proof consists in observing that choosing $\eps_0$ as explained above, taking a derivative of~\eqref{eq:MA}, and denoting by $M$ the cofactor matrix of $\Id + \DDD^2 u$ (which is divergence-free, see~\cite[Chapter 8]{evans2022partial}), leads to:
	\begin{equation}
	\label{eq:diff_MA}
	\tr(M \DDD^2 \partial_i u ) = \Div (M \nabla \partial_i u ) = \partial_i g.
	\end{equation}
	Observe that because of Theorem~\ref{thm:caf}, the field of nonnegative symmetric matrices $\Id + \DDD^2 u$ is continuous, hence bounded, and of determinant bounded below. Therefore, its eigenvalues are uniformly bounded below and above, and the same automatically holds for 
    \begin{equation*}
        M = \det (\Id + \DDD^2 u ) (\Id + \DDD^2 u)^{-1}.
    \end{equation*}
    This linear elliptic PDE therefore yields a gain of two derivatives, and we can already conclude that the result holds for $k=1$. However, for higher derivatives, this gain only holds provided $M$ is regular enough, and as $M$ depends on $u$, one needs to be precise in the computations.
	
	The proof relies on an induction argument working as follows. First, we apply $k$ derivatives of~\eqref{eq:diff_MA} and we write the PDE obtained as
	\begin{equation}
	\label{eq:MA_k_deriv}
	\tr(M \DDD^2 \partial_\gamma u ) = \Div( M \nabla \partial_\gamma u ) = \partial_\gamma g + \mathrm{LOT},
	\end{equation}
	where $\gamma \in \{1, \dots , d\}^k$ and $\partial_\gamma$ is a notation for $\partial_{\gamma_1} \cdots \partial_{\gamma_k}$. Then, we observe that $\mathrm{LOT}$ is a sum of terms implying derivatives of $u$ of order $k-1$ at most, that we can bound by induction hypothesis. Finally, we estimate these terms in $L^2_x$ thanks to Sobolev inequalities, and hence conclude that $\partial_\gamma u \in H^2_x$, as announced.
	
	With this sketch of proof in mind, we could try to prove the result using the conclusion of the theorem as induction hypothesis. Actually, this would not be enough. Indeed, in order to prove that the terms in $\mathrm{LOT}$ above are in $L^2_x$, we will have to take benefit of the following fact. Imagine that you want to prove~\eqref{eq:estim_MA} for $k = 2$. Then, we know that for any $i= 1, \dots, d$, $\partial_i u$ solves~\eqref{eq:diff_MA}, but with a right-hand side in $H^1_x$, and hence, by Sobolev embeddings, in some $L^p$ with $p>2$. We conclude that $\partial_i u \in W^{2,p}_x$, which is an improvement with respect to the result at stage $k=1$. 
	
	In order to take this effect into account, for a given $n \in \N$ small enough, we define $p_n$ as the exponent for which by Sobolev embeddings, $H^n_x \hookrightarrow L^{p_n}_x$, that is
	\begin{equation*}
	\frac{1}{p_n} = \frac{1}{2} - \frac{n}{d},
	\end{equation*}
	as long as the quantity in the right-hand side is positive. We call $N(d)$ the largest integer for which this is the case. Note that whenever $d$ is even, $1/2 - (N(d) + 1)/d = 0$, which leads to a critical case in Sobolev embeddings where $H^{N(d) + 1}_x$ is not embedded in $L^\infty_x$, but only in $\cap_{p\geq 1}L^p_x$. To avoid this slight complication, from now on, we assume up to formulating the problem in $\T^{d+1}$ instead of $\T^d$ if necessary that $d$ is odd. In this way, whenever $n > N(d)$, $H^n_x \hookrightarrow L^\infty_x$. For this reason, for all $n > N(d)$, we set, $p_n \vcentcolon= \infty$.
		
	We use the following statement as an induction hypothesis, which holds true for $k=1$ using~\eqref{eq:diff_MA}. 
	\begin{quote}
		$(\mathrm{H}_k)$. For any $n =1, \dots, N(d)$, there exists a polynomial function $F_{k,n}: \R_+ \to \R_+$ canceling at $0$ such that for all smooth $u$ and $g$ satisfying~\eqref{eq:MA} with $\| g \|_{H^{k_0}_x} \leq \eps_0$, we have
		\begin{equation*}
		\| u \|_{W^{k+2,p_n}_x} \leq F_{k,n}( \| g \|_{W^{k,p_n}_x}).
		\end{equation*}
	\end{quote}
	
	Let us assume that $(\mathrm H_{k-1})$ holds for some $k \geq 2$ and let us prove $(\mathrm H_{k})$. We pick $n \in \{ 1, \dots, N(d)\}$, a source term $g \in W^{k,p_n}_x$ with $\| g \|_{H^{k_0}_x} \leq \eps_0$ and $\gamma \in \{1, \dots, d\}^{k}$. We just need to prove that there is a polynomial function\footnote{We use the same notation for every such function, hoping it does not lead to any unnecessary confusion.} $F_{k,n}$ canceling at $0$ such that
	\begin{equation*}
	\| \partial_\gamma \DDD^2 u \|_{L^{p_n}_x} \leq F_{k,n} (\| g \|_{W^{k,p_n}_x}).
	\end{equation*}
	Thanks to standard results in elliptic regularity theory, differentiating $k$ times the Monge--Ampère equation~\eqref{eq:MA} as in~\eqref{eq:MA_k_deriv}, we just need to prove that there exists a polynomial function $F_{k,n}$ canceling at $0$ such that
	\begin{equation*}
	\| \mathrm{LOT} \|_{L^{p_n}_x} \leq F_{k,n}(\| g \|_{W^{k,p_n}_x}).
	\end{equation*}
	Comparing~\eqref{eq:diff_MA} and~\eqref{eq:MA_k_deriv}, we see that
	\begin{equation*}
	\mathrm{LOT} = \tr \left\{  \partial_{\hat \gamma}\big( M \DDD^2 \partial_{\gamma_1} u \big) - M \DDD^2 \partial_\gamma u \right\},
	\end{equation*}
	where $\hat \gamma \vcentcolon= (\gamma_2, \dots, \gamma_k)$. It is a sum of terms of the form
	\begin{equation*}
	\partial_\alpha m \, \partial_\beta (\partial_{ij} u),
	\end{equation*}
	where $m$ is a coefficient of $M$, $i,j \in \{1, \dots, d\}$, and $\alpha$ and $\beta$ are multi-indices of respective sizes $|\alpha|$ and $|\beta|$ satisfying $|\alpha| \geq 1$, $|\beta|\geq1$ and $|\alpha| + |\beta| = k$. 
	Recall that $M$ is the cofactor matrix of $\Id + \DDD^2 u$. Therefore $m$ is itself a product of at most $d-1$ factors of type $\partial_{i'j'} u$. All in all,  $\mathrm{LOT}$ is a sum of terms of the form
	\begin{equation}
	\label{eq:each_term_LOT-MA}
	\prod_{a=1}^{d'} \partial_{\alpha_a} (\partial_{i_aj_a} u),
	\end{equation}
	where $2\leq d' \leq d$, $i_1, j_1, \dots, i_a, j_a \in \{1, \dots, d\}$, and $\alpha_1, \dots, \alpha_{d'}$ are multi-indices with $|\alpha_1| + \cdots + |\alpha_{d'}| = k$ and at least two of these $\alpha_a$, say for $a=1$ and $a=2$ satisfy $|\alpha_a| \geq 1$. In particular, they all satisfy $|\alpha_a|\leq k-1$. From now on, we assume that $(\alpha_a)$ is ordered in such a way that $|\alpha_1| \geq |\alpha_2| \geq \dots \geq |\alpha_{d'}|$. Our goal is hence to analyze precisely the integrability of terms of type $\partial_{\gamma_a}(\partial_{i_aj_a}u)$.
	
	Let us first give a precise account for the integrability of $\partial_\alpha (\partial_{ij} u)$ for some choice of $i,j \in \{1, \dots, d\}$ and $|\alpha| \leq k-1$, using the induction assumption. 
	By Sobolev embedding, we know that if $l \leq k$, then $W^{k,p_n}_x \hookrightarrow W^{k-l,p_{n+l}}_x$. 
	Therefore, we infer that $g \in W^{k,p_n}_x \hookrightarrow W^{|\alpha|, p_{n + k - |\alpha|}}$. Therefore, by induction hypothesis, $\partial_\alpha(\partial_{ij} u) \in L^{p_{n + k - |\alpha|}}_x$. Also, by the same argument, if $|\alpha| \leq k_0$, since $g \in H^{k_0}_x$, $\partial_\alpha(\partial_{ij} u) \in L^{p_{k_0 - |\alpha|}}$. All in all, we conclude that calling $K \vcentcolon= \max(k_0, k+n)$,
	\begin{equation*}
	\partial_\alpha(\partial_{ij} u)  \in L^{p_{K-|\alpha|}},
	\end{equation*}
	with a bound that is polynomial in $\| g \|_{W^{k,p_n}}$ and $\eps_0$ and which cancels when $\| g \|_{W^{k,p_n}}=0$.
	
	Going back to~\eqref{eq:each_term_LOT-MA} and using H\"older's inequality, we deduce that in order to prove
	\begin{equation*}
	\prod_{a=1}^{d'} \partial_{\alpha_a} (\partial_{i_aj_a} u) \in L^{p_n}
	\end{equation*}
	with a polynomial bound, it suffices to show that
	\begin{equation*}
	\sum_{a = 1}^{d'} \frac{1}{p_{K-|\alpha_a|}} \leq \frac{1}{p_n},
	\end{equation*}
	(with convention $1/\infty = 0$). This rewrites
		\begin{equation}
		\label{eq:numerologyMA}
	\sum_{a = 1}^{d'} \left( \frac{1}{2} - \frac{K - |\alpha_a|}{d} \right)_+ \leq \frac{1}{2} - \frac{n}{d}.
	\end{equation}
	Let
	\begin{equation*}
d'' \vcentcolon= \max \{ a =1, \dots, d' \mbox{ such that } K- |\alpha_a| \leq N(d) \},
	\end{equation*}
	with convention $\max\emptyset = 0$. We have
	\begin{equation*}
	\sum_{a = 1}^{d'} \left( \frac{1}{2} - \frac{K - |\alpha_a|}{d} \right)_+ = \sum_{a = 1}^{d''}  \frac{1}{2} - \frac{K - |\alpha_a|}{d}.
	\end{equation*}
	therefore, if $d'' = 0$, the conclusion is obvious. If $d'' \geq 1$, inequality~\eqref{eq:numerologyMA} can be reformulated as
	\begin{equation*}
	\frac{d''-1}{2} + \frac{n + |\alpha_1| + \dots + |\alpha_{d''}|}{d} \leq (d''- 1) \frac{K}{d} + \frac{K}{d}.
	\end{equation*}
	But as $K \geq k_0 > d/2$, we have $K/d \geq 1/2$, and as $K \geq n+k$ and $|\alpha_1| + \dots |\alpha_{d''}| \leq k$, this last inequality holds true, which concludes our induction argument.
\end{proof}

\subsubsection{Proof of Theorem~\ref{stabMA}}

 	For a given matrix $M$, we decompose $ \det (\Id + M )$ as
 	\begin{equation*}
 	1 + \tr( M ) + P(M),
 	\end{equation*}
 	where $P$ is a polynomial whose terms have at least degree $2$. Now, notice that by the algebraic identity
 	\begin{equation*}
 	y_1 \dots y_n - x_1 \dots x_n = \sum_{i=1}^n x_1 \dots x_{i-1} (y_i - x_i) y_{i+1} \dots y_n
 	\end{equation*}
 	applied to the difference of each terms of $P$, we see that if $M$ and $N$ are two given matrices, $P(N) - P(M)$ can be written under the form
 	\begin{equation*}
 	P(N) - P(M) = \tr\big(A(M,N) (N-M)\big), 
 	\end{equation*}
 	where $A(N,M)$ is a matrix valued polynomial with respect to the coefficients of $M$ and $N$, all of whose terms have at least degree $1$ (and is therefore small when $M$ and $N$ are small). Notice that whenever $M$ and $N$ are symmetric, we can replace $A$ by $\frac{1}{2}(A + {}^t A)$ and assume that $A$ is symmetric as well. We make this assumption in what follows. With these definitions at hand, taking $u$, $v$, $g$ and $h$ as in the statement of the theorem and calling $w \vcentcolon= v-u$ and $ q\vcentcolon= g-h$, we have
 	\begin{equation*}
 	\tr\Big(\big(\Id + A(\DDD^2 u, \DDD^2 v)\big) \DDD^2w\Big) =q.
 	\end{equation*}
 	Testing this identity against $w$, and after integrating by parts, we get:
 	\begin{multline}
 	\label{eq:estim_stabMA}
 	\int \nabla w \cdot (\Id + A(\DDD^2u, \DDD^2 v))\nabla w \D x \\
    = -\int w \mathrm{div} \Big(  A(\DDD^2 u, \DDD^2 v) \Big) \cdot \nabla w \D x - \cg q, w \cd .
 	\end{multline}
 	Now, because of Theorem~\ref{thm:boundMA}, if $\eps_1$ is sufficiently small, then $\DDD^2 u$ and $\DDD^2 v$ are small enough in $H^{k_1}_x$ so that since $H^{k_1 - 1}_x$ is a Banach algebra included in $L^\infty_x$,
 	\begin{equation*}
 	\|   A(\DDD^2 u, \DDD^2 v) \|_{L^{\upinfty}_x} \leq \frac{1}{4} \qquad \mbox{and} \qquad \Big\| \mathrm{div}\Big( A(\DDD^2 u , \DDD^2 v)\Big) \Big\|_{L^{\upinfty}_x} \leq \frac{1}{4C},
 	\end{equation*} 
 	where $C$ is the Poincaré constant of the torus $\T^d$.
 	From these estimates and the identity~\eqref{eq:estim_stabMA}, we finally find
 	\begin{equation*}
 	\frac{3}{4} \| \nabla w \|_{L^2_x} \leq \frac{1}{4C} \| w \|_{L^2_x} + \| q \|_{\dot H^{-1}_x}.
 	\end{equation*}
	The result follows using the Poincaré inequality. \qed

\section{Application to the justification of semiclassical limits}
\label{sec:semiclassical}

In this section, we study the semiclassical limit from Schr\"odinger--Poisson equation for mixed states to the Vlasov--Poisson system. The starting point is the Hartree--Poisson equation  (recall Motivation IV in the introduction) which is a generalization for mixed states of the Schr\"odinger--Poisson (which is relevant for pure states):
\begin{equation}
\label{eq:SP-hartree}
	\left\{
	\begin{aligned}
i \varepsilon \partial_t \gamma^\eps &=\left[ -\frac{\varepsilon^2}{2} \Delta_x + V_\eps,\gamma^\eps\right], \\
-\Delta_x V^\eps &=  \rho^\eps, \quad \rho^\eps(t,x)= \gamma^\eps(t,x,x), \\ 
\gamma^\eps|_{t=0} &=  \gamma^\eps_0,
\end{aligned}
\right.
\end{equation}
where $\gamma^\eps \in \mathscr{L}(L^2(\R^d))$ is a nonnegative trace-class operator.  
From now on, we restrict to dimensions $d\geq 3$.

The semiclassical limit $\eps\to 0$ for~\eqref{eq:SP-hartree} has been the object of several works, the main goal being the derivation of the Vlasov--Poisson equation. Let us summarize in an informal, somewhat loose, statement, some known results in this direction. Beforehand it is convenient to introduce the so-called Wigner distribution associated to an operator.

\begin{Def}
Given an operator $\gamma^\eps$, the Wigner distribution $W[\gamma^\eps]$ associated to $\gamma^\eps$ is 
\begin{equation}
    \label{def:wigner-true}
   W[\gamma^\eps](x,v)\vcentcolon= 
   \frac{1}{(2\pi)^d} \int_{\R^d} \gamma^\eps\left(x+\eps\frac{y}{2}, x- \eps\frac{y}{2}\right)  e^{-i v\cdot y} \,\, \mathrm{d} y.
   \end{equation}

\end{Def}

\begin{Thm}
\label{thm:class-semiclass}
The following holds, in the semiclassical limit $\eps\to 0$.
\begin{itemize}
    \item (Weak solutions, from \cite{LP}) Given $f_0 \in L^1\cap L^2(\R^d)$ with finite kinetic energy, there exists a sequence of initial data $(\gamma_0^\eps)$ such that $W[\gamma_0^\eps]$ approximates $f_0$, and there exists an associated sequence of solutions $(\gamma^\eps)$ to the Hartree equation~\eqref{eq:SP-hartree} such that $W[\gamma^\eps]$ weakly converges, up to a sequence, to a global weak solution $f$ to the Vlasov--Poisson equation, associated to the initial condition $f_0$.
    
    \item (Strong solutions, from \cite{LS,CLS,LS2}) Given $f_0$ in a weighted Sobolev space $\mathrm{H}^n_r$ (recall Definition~\ref{def:weighted-v-sobolev}), with $n$ and $r$ large enough, there exists a sequence of initial data $(\gamma_0^\eps)$ such that $W[\gamma_0^\eps]$ approximates $f_0$, and there exists an associated sequence of solutions $(\gamma^\eps)$ to the Hartree equation~\eqref{eq:SP-hartree} such that $W[\gamma^\eps]$ strongly converges to the strong solution $f$ to the Vlasov--Poisson equation, associated to the initial condition $f_0$. Moreover the convergence is quantitative.
\end{itemize}
    
\end{Thm}

The results of \cite{Lafleche,IacobelliLafleche} are somehow in between these two classes of results; 
when higher moments in velocity are controlled, they obtain quantitative estimates for the semiclassical limit, in the quantum Wasserstein pseudo-distance of \cite{GP}.

We shall provide, thanks to the multiphasic framework, an alternative class of results for the semiclassical limit, that is somewhat intermediate between the aforementioned weak and strong results.  Rather, we put forward a class of data whose semiclassical limit can be rough in velocity, and with finite regularity in $x$, but which is still uniquely determined.
When the initial condition writes under the form
\begin{equation}
\label{eq:init-hartree-chapter3}
\gamma^\eps_0(x,y) = \int_I  \Psi_0^{\eps,\alpha} (x) \overline{\Psi_0^{\eps,\alpha}}(y) \, \D \mu(\alpha),
\end{equation}
where $(I, \mu)$ is a given measure space, we expect the associated solution to keep the form
\begin{equation}
\label{eq:ansatz-hartree-chapter3}
\gamma^\eps(t,x,y) = \int_I  \Psi^{\eps,\alpha} (t,x) \overline{\Psi^{\eps,\alpha}}(t,y) \, \D \mu(\alpha),
\end{equation}
and 
it turns out relevant to consider the system of Schr\"odinger--Poisson equations set on $\R^d$:
\begin{equation}
\label{eq:SP}
	\left\{
	\begin{aligned}
i \varepsilon \partial_t \Psi^{\eps,\alpha} + \frac{\varepsilon^2}{2} \Delta_x \Psi^{\eps,\alpha} &=   V^\eps \Psi^{\eps,\alpha}, \\
-\Delta_x V^\eps &=  \int_I |\Psi^{\eps,\alpha}|^2 \, \D \mu(\alpha), \\
\Psi^\eps|_{t=0} &=  \Psi^{\eps,\alpha}_0.
\end{aligned}
\right.
\end{equation}
In the following, we use, as everywhere else in this monograph,  bold letters to denote  whole families indexed by $\alpha$.
For such operators, we have
\begin{Def}
\label{def:wigner}
\index{Wp@$W[\boldsymbol\Psi^{\eps}]$: Wigner distribution associated to  $\boldsymbol\Psi^\eps=(\Psi^{\eps,\alpha})_{\alpha \in I}$}
For 
$\gamma^\eps = \int_I  \Psi^{\eps,\alpha} (x) \overline{\Psi^{\eps,\alpha}}(y) \, \D \mu(\alpha)
$, the associated Wigner distribution writes as $   W[\gamma^\eps]=  W[\boldsymbol\Psi^{\eps}]$, where
\begin{equation}
    \label{eq:wigner}
   W[\boldsymbol\Psi^{\eps}](x,v)\vcentcolon= 
   \frac{1}{(2\pi)^d} \int_{\R^d} \int_I \Psi^{\eps,\alpha}\left(x+\eps\frac{y}{2} \right) \overline{\Psi^{\eps,\alpha}}\left(x-\eps\frac{y}{2} \right)  e^{-i v\cdot y} \, \mathrm{d} \mu(\alpha) \, \mathrm{d} y.
   \end{equation}
\end{Def}
 Consider for each $\alpha \in I $ a highly oscillating initial condition of WKB type:
\begin{equation}
\label{eq:multiphasicquant0}
\Psi^{\eps,\alpha}_0(x) = a_0^{\eps,\alpha} (x) \exp \left( \frac{i}{\eps} S^{\eps,\alpha}_0(x) \right),
\end{equation}
for families of amplitudes $\boldsymbol{a}_0^{\eps}$ and phases $\boldsymbol{S}^{\eps}_0$ with values in $\R_+$ and $\R$, respectively.

Following the framework of Theorem~\ref{thm-LWPmultiphase-densityvelocity-re} (that is particularly relevant for the $\R^d$ case as considered here), we focus on phases of the form
$$
S^{\eps,\alpha}_0(x) = x \cdot \lambda^\alpha +  \widetilde{S}^{\eps,\alpha}_0(x),
$$
where $\lambda: I \rightarrow \R^d$ is a given measurable function (typically $\lambda^\alpha \equiv \alpha$ when $I=\R^d$), and $\widetilde{S}^{\eps,\alpha}_0$ has to be seen as an integrable regular perturbation of the linear phase.
The multiphasic formalism developed in Chapter~\ref{Part1-LWP}  entails, in a natural way, the following result.

\begin{Thm}
\label{thm:semiclassical}
Let $k\in \N$ such that $k> 2+d/2$. Let $\boldsymbol{a}^{\eps}_0  \in L^\infty_\alpha H^{k-1}_x\cap \dot{H}^{-1}_x$ and 
$ \boldsymbol{\widetilde{S}}^{\eps}_0  \in \mathcal{H}^{k+1,\infty}$ 
with norms that are uniform with respect to $\eps \in (0,1]$. Then there exists $T>0$ and solutions of the form
$$
\Psi^{\eps,\alpha} (t,x) =  a^{\eps,\alpha} (t,x) \exp \left( \frac{i}{\eps} S^{\eps,\alpha}(t,x) \right)
$$
to the Schr\"odinger--Poisson system~\eqref{eq:SP} on $[0,T]$, with initial conditions~\eqref{eq:multiphasicquant0}. 
Moreover, the family $(\boldsymbol{a}^{\eps}, \nabla \boldsymbol{S}^{\eps}-\boldsymbol\lambda) $  is uniformly bounded in $L^\infty([0,T]; L^\infty_\alpha H^{k-1}_x\cap \dot{H}^{-1}_x) \times L^\infty([0,T];  \mathcal{H}^{k,\infty} )$.

If one assumes that 
$$(\boldsymbol{a}^{\eps}_0, \nabla\boldsymbol{S}^{\eps}_0-\boldsymbol\lambda) \overset{\eps \rightarrow 0}{\longrightarrow} (\sqrt{\rhorho_{0}} , \uu_0), \ \ \text{in} \ \ L^\infty_\alpha (H^{k-3}_x\cap \dot H^{-1}_x ) \times \mathcal{H}^{k-2,\infty}, $$ 
where $(\rhorho_{0}, \uu_0) \in L^\infty_\alpha (H^{k-2}_x\cap \dot H^{-1}_x) \times \mathcal{H}^{k-1,\infty}$ with $\rhorho_0 \geq 0$, 
then we have the convergence
$$(\boldsymbol{a}^{\eps}, \nabla \boldsymbol{{S}}^{\eps}-\boldsymbol\lambda) \overset{\eps \rightarrow 0} {\longrightarrow} (\sqrt{\rhorho}, \vv-\boldsymbol\lambda), \ \ \text{in} \ \ L^\infty([0,T]; L^\infty_\alpha H^{s-1}_x ) \times L^\infty([0,T]; \mathcal{H}^{s,\infty}),$$ 

for all $s<k$, where $(\rhorho, \vv)$  is the unique solution to the multiphasic Vlasov--Poisson system~\eqref{eq:VP} with initial condition $({\rhorho}_{0} , \vv_0= \boldsymbol\lambda+\uu_{0})$ on $[0,T]$.

Moreover, the following quantitative convergences hold:
\begin{multline*}
 \sup_{[0,T]}  \lVe \boldsymbol{a}^{\eps} - \sqrt{\rhorho} \rVe_{L^1_\alpha \dot{H}^{-1}_x} +\sup_{[0,T]} \lVe \nabla \boldsymbol{{S}}^{\eps} - \vv \rVe_{L^p_\alpha L^2_x} \\
 \lesssim  \lVe \boldsymbol{a}^{\eps}_0 - \sqrt{\rhorho_0} \rVe_{L^1_\alpha \dot{H}^{-1}_x} +\lVe \nabla \boldsymbol{{S}}^{\eps}_0 - \vv_0 \rVe_{L^{\upinfty}_\alpha L^2_x} + \eps, 
 \end{multline*}
 and
 \begin{multline*}
 \sup_{[0,T]}  \lVe \boldsymbol{a}^{\eps} - \sqrt{\rhorho} \rVe_{L^1_\alpha H^{j}_x} 
 + \sup_{[0,T]} \lVe \nabla(\nabla \boldsymbol{{S}}^{\eps} - \vv )\rVe_{L^{\upinfty}_\alpha H^j_x}\\
 \lesssim \lVe \boldsymbol{a}^{\eps}_0 - \sqrt{\rhorho}_0 \rVe_{L^1_\alpha H^{j}_x} + \lVe \nabla (\nabla \boldsymbol{{S}}^{\eps}_0 - \vv_0 )\rVe_{L^{\upinfty}_\alpha {H}^j_x}    + \eps, 
\end{multline*}
for all $j=0,\cdots, k-3$.

Finally, if $\mathcal{F}_{v \rightarrow y}$ stands for the partial Fourier transform in the $v$ variable, we have
\begin{multline}
        \label{eq:conv-wkb}
    \sup_{t \in [0,T]} \sup_{ \substack{\varphi \in \mathscr{S}(\R^d_x \times\R^d_y) \\  \left\Vert \langle y \rangle^2 \mathcal{F}_{v\rightarrow y} \varphi\right\Vert_{L^{\upinfty}_x L^1_y} \leq 1}} \left| \int_{\R^d \times\R^d} \left(W[\boldsymbol\Psi^{\eps}](t,x,v) - f(t,x,v)  \right) \varphi(x,v)\, \mathrm{d}x \,  \mathrm{d}v \right| \\
    \lesssim \lVe \boldsymbol{a}^{\eps}_0 - \sqrt{\rhorho}_0 \rVe_{L^1_\alpha (\dot{H}^{-1}_x \cap L^2_x)} + \lVe \nabla \boldsymbol{{S}}^{\eps}_0 - \vv_0 \rVe_{L^{\upinfty}_\alpha H^1_x}   + \eps,
\end{multline}
where 
\begin{itemize}
    \item $W[\boldsymbol\Psi^{\eps}]$ is the Wigner distribution associated with the family $\boldsymbol\Psi^\eps$;

    \item $f(t,x,v)$ is the kinetic distribution associated with the multiphasic data $(\rhorho, \vv)$, that is
    $$
    f(t,x,v) = \int_I \rho^\alpha(t,x) \otimes  \delta_{v=v^\alpha(t,x)}
 \, \mathrm{d} \mu(\alpha).    
 $$
\end{itemize}

\end{Thm}

Loosely speaking,  this therefore justifies the derivation of the Vlasov--Poisson equation, starting from the semiclassical Schr\"odinger--Poisson equation with  mixed WKB states. 
As a corollary we can complement Theorem~\ref{thm:class-semiclass} as follows.

\begin{Cor}
\label{cor:semi}
    Let  $f_0 \in L^\infty_{v,\langle v\rangle^p} (H^k_x \cap \dot{H}^{-1}_x)$ with $k>1+d/2$ and $p>d$, such that $f_0 \geq 0$, 
    and such that $\sqrt{f_0} \in L^\infty_{v,\langle v\rangle^p} (H^k_x \cap \dot{H}^{-1}_x)$.
    Then there exists a sequence of initial data $(\gamma_0^\eps)$ such that $W[\gamma_0^\eps]$ approximates $f_0$, and there exists an associated sequence of solutions $(\gamma^\eps)$ to the Hartree equation~\eqref{eq:SP-hartree} such that $W[\gamma^\eps]$ weakly converges, on a uniform interval of time $[0,T]$, to the unique weak solution $f$ to the Vlasov--Poisson equation, associated to the initial condition $f_0$. Moreover, the weak convergence is quantitative in the sense of~\eqref{eq:conv-wkb}.
\end{Cor}

\begin{Rem} We note that in this result, the condition $\sqrt{f_0} \in L^\infty_{v,\langle v\rangle^p} (H^k_x \cap \dot{H}^{-1}_x)$ is quite demanding on the way $f_0$ cancels. 

Let us also remark that the assumptions of Theorem~\ref{thm:semiclassical}--Corollary~\ref{cor:semi} are not compatible with those of \cite{Lafleche,IacobelliLafleche}; in particular,  Theorem~\ref{thm:semiclassical} does not follow from the latter.
\end{Rem}

\begin{proof}[Proof of Corollary~\ref{cor:semi}]
    For $f_0 \in L^\infty_{v,\langle v\rangle^p} (H^k_x \cap \dot{H}^{-1}_x)$, as already performed several times in this monograph,  we take
    $$
    I= \R^d, \quad \mu = \langle \alpha \rangle^{-p } \mathrm{d} \alpha,
    $$
    and
    $$
    \rho_0^\alpha(x) = f_0(x,\alpha ) \langle \alpha \rangle^p, \quad v_0^\alpha = \alpha,
    $$
    so that 
    $$f_0(x, v)=\int_{I} \rho^\alpha_0(x) \otimes  \delta_{v=v^\alpha_0(x)} \, \mathrm{d}\mu(\alpha).$$
    In particular, $(I, \mu)$ is of finite mass since $p>d$.
    We then take $\boldsymbol{\lambda} \equiv \alpha$ and 
    $$
    \boldsymbol{a}_0^\eps=\sqrt{\rhorho}_0, \quad \boldsymbol{S}_0^\eps = \alpha,
    $$
    and we are in position to apply Theorem~\ref{thm:semiclassical}, which directly yields the result.
\end{proof}

The remaining of this section is dedicated to the proof of the main theorem.

\begin{proof}[Proof of Theorem \ref{thm:semiclassical}]
We closely follow Grenier's approach \cite{Gr98} to the WKB approximation  for the nonlinear Schrödinger equation, which we generalize to the multiphasic setting. Namely, we look for solutions of the form
$$
\Psi^{\eps,\alpha}(t,x) =  a^{\eps,\alpha} (t,x) \exp \left( \frac{i}{\eps} S^{\eps,\alpha}(t,x) \right),
$$
where $a^{\eps,\alpha}$ is complex-valued  and $S^{\eps,\alpha}$ is real-value: we obtain the system
\begin{equation}
\label{eq:prequant}
	\left\{
	\begin{aligned}
&\partial_t S^{\eps,\alpha} + \frac{1}{2} | \nabla S^{\eps,\alpha} |^2= -V^\eps, \\
&\partial_t a^{\eps,\alpha} + \nabla S^{\eps,\alpha}\cdot \nabla a^{\eps,\alpha} + \frac{1}{2} a^{\eps,\alpha} \Div (\nabla S^{\eps,\alpha})  = i \frac{\eps}{2}  \Delta  a^{\eps,\alpha},
\end{aligned}
\right.
\end{equation}
where 
\begin{align*}
    -\Delta V^\eps=\int_I |a^{\eps,\alpha}|^2 \, \D \mu(\alpha).
\end{align*}

\begin{Rem}
    Choosing $a^{\eps, \alpha}$ to be nonnegative (which yields to the so-called Madelung transform) would formally lead, after plugging the above ansatz and identifying of the real and imaginary part, to the different system
    \begin{equation*}
	\left\{
	\begin{aligned}
&a^{\eps,\alpha} \left(\partial_t S^{\eps,\alpha} + \frac{1}{2} | \nabla S^{\eps,\alpha} |^2 +V^\eps \right)=\frac{\eps^2}{2}\Delta a^{\eps,\alpha} , \\
&\partial_t a^{\eps,\alpha} + \nabla S^{\eps,\alpha}\cdot \nabla a^{\eps,\alpha} + \frac{1}{2} a^{\eps,\alpha} \Delta S^{\eps,\alpha}  =0.
\end{aligned}
\right.
\end{equation*}
Relaxing the sign condition for $a^{\eps, \alpha}$ is crucial in the argument from \cite{Gr98}.
\end{Rem}

Differentiating the first equation of \eqref{eq:prequant}, we see that the pair $(a^{\eps,\alpha}, u^{\eps,\alpha}\vcentcolon = \nabla S^{\eps,\alpha}-\lambda^\alpha)$  satisfies
\begin{equation}
\label{eq:quant}
	\left\{
	\begin{aligned}
&\partial_t u^{\eps,\alpha} + (u^{\eps,\alpha} + \lambda^\alpha) \cdot \nabla u^{\eps,\alpha} = -\nabla V^\eps, \\
&\partial_t a^{\eps,\alpha} + (u^{\eps,\alpha}+\lambda^\alpha)\cdot \nabla a^{\eps,\alpha} + \frac{1}{2} a^{\eps,\alpha} \Div (u^{\eps,\alpha})  = i \frac{\eps}{2}  \Delta  a^{\eps,\alpha}, \\
&u^{\eps,\alpha}|_{t=0} = \nabla \widetilde{S}^{\eps,\alpha}_0, \quad a^{\eps,\alpha}|_{t=0} = a^{\eps,\alpha}_0.
\end{aligned}
\right. 
\end{equation}
The idea is to study the Cauchy problem for this system, as in Theorem~\ref{thm-LWPmultiphase-densityvelocity-re}. Then one recovers the phase $S^{\eps,\alpha} $ by integrating the first equation of~\eqref{eq:prequant} (this is the standard eikonal equation). For the study of this Cauchy problem, we rely on the abstract framework of Chapter~\ref{Part1-LWP}.
The main task is to check that a uniform version of Assumptions~\ref{ass:bound_high_regularity-rho} and \ref{ass:stability-rho} are satisfied (in the case $p=\infty$). Namely we intend to show a bound in high regularity and stability in low regularity for the force field; we shall also need an analogue of Lemmas~\ref{LM:SobEstimRHO} and \ref{LM:estimStabContinuity} for the  equation satisfied by  $a^{\eps,\alpha}$. In the process, we need to ensure that all estimates are uniform with respect to $\eps$. We will then be able to apply Theorem~\ref{thm-LWPmultiphase-densityvelocity-re} (complemented by Remark~\ref{rem:Rd-alpha}), and the time of existence together with the control of the norms will also be uniform in $\eps$: namely, we shall obtain that $(\boldsymbol{a}^{\eps}, \nabla \boldsymbol{S}^{\eps}-\boldsymbol\lambda) $  is uniformly bounded in $L^\infty([0,T]; L^\infty_\alpha H^{k-1}_x\cap \dot{H}^{-1}_x) \times L^\infty([0,T];  \mathcal{H}^{k,\infty} )$.

\begin{Lem}
\label{lem:quant} 
Let $k \geq 2+d/2$. There exists $K>0$ such that the following holds for any $T>0$ and $\eps \in (0,1)$.
\begin{itemize}
    \item If $\vv^\eps \in L^\infty([0,T]; \mathcal H^{k,\infty})$ and $\boldsymbol{a}^{\eps}_0 \in L^\infty_\alpha H^{k-1}_x$
then there exists a unique solution $\boldsymbol{a}^{\eps,\alpha} \in L^\infty([0,T]; L^\infty_\alpha H^{k-1}_x)$
to the second equation of \eqref{eq:quant} with initial data $\boldsymbol{a}^{\eps}_0$, such that for all $t \in [0,T]$
\begin{equation}
\label{eq:main-quantum}
\Vert \boldsymbol{a}^{\eps} \Vert_{L^{\upinfty}_\alpha H^{k-1}_x}  \leq  \Vert \boldsymbol{a}^{\eps}_0 \Vert_{L^{\upinfty}_\alpha H^{k-1}_x}   \exp\left( K \int_0^t\Vert\vv^\eps(s)\Vert_{\mathcal{H}^{k,\infty}} \D s \right).
\end{equation}
\item If one considers two vector fields $\vv^\eps, \ww^\eps \in L^\infty([0,T]; \mathcal H^{k,\infty})$ and initial conditions $\boldsymbol{a}^{\eps}_0,\boldsymbol{b}^{\eps}_0$ then their corresponding solutions $\boldsymbol{a}^{\eps,}, \boldsymbol{b}^{\eps}$ to the second equation of \eqref{eq:quant} satisfy for all $t \in [0,T]$
\begin{align}
\label{eq:main-quantum2}
\begin{split}
	&\| \boldsymbol{a}^{\eps}(t)- \boldsymbol{b}^{\eps} (t)\|_{L^1_\alpha \dot H^{-1}_x} 
	 \\
     &\leq \exp\left( K \int_0^t (\| \nabla \vv^\eps(s) \|_{L^{\upinfty}_\alpha H^{k-1}_x} + \| \DDD \ww^\eps(s) \|_{L^{\upinfty}_\alpha H^{k-1}_x}) \D s \right) \\ 
     &  \ \Bigg\{ 	\| \boldsymbol{a}^{\eps}_0- \boldsymbol{b}^{\eps}_0 \|_{L^1_\alpha \dot H^{-1}_x}  
	+\| \boldsymbol{a}^\eps \|_{L^\infty([0,t];L^{\upinfty}_\alpha H^{k-1}_x)}\int_0^t \| \ww^\eps(s) - \vv^\eps(s) \|_{L^1_\alpha  L^2_x} \D s \Bigg\}.
    \end{split}
	\end{align}
\end{itemize}
\end{Lem}

\begin{proof}[Proof of Lemma~\ref{lem:quant}]
Let us start by explaining the a priori estimate~\eqref{eq:main-quantum}: it is obtained by an energy estimate exactly as in the proof of Lemma~\ref{LM:SobEstimRHO}, noticing  that there is no contribution due to the term  $i\eps \Delta a^{\eps,\alpha}$, as 
$$\Re ( i \langle\Delta \varphi, \varphi \rangle_{L^2(\R^d;\C)}) =0.$$ 
Then the existence of a solution to the second equation of  \eqref{eq:quant} relies on this a priori estimate and on an approximation procedure; namely we first solve the linear equation 
\begin{align}
\label{eq:quant2'} &\partial_t a^{\eps,\alpha}_{\delta}  + (u^{\eps,\alpha}+\lambda^\alpha)\cdot \nabla a^{\eps,\alpha}_\delta + \frac{1}{2} a^{\eps,\alpha}_\delta \Div (u^{\eps,\alpha})=   i\frac{\eps}{2} J_\delta^2 \Delta  a^{\eps,\alpha}_{\delta},
\end{align}
 where $J_\delta$ is a smooth approximation of unity, by a Banach fixed-point argument.  We note that~\eqref{eq:main-quantum} still holds for solutions to \eqref{eq:quant2'}, uniformly in $\delta$, and one can pass to the limit $\delta\to0$ and obtain a solution on the whole interval $[0,T]$. Since the equation is linear, uniqueness also follows from \eqref{eq:main-quantum}.
 
 	The stability estimate~\eqref{eq:main-quantum2} can also be obtained following the proof of Lemma~\ref{LM:estimStabContinuity}, noting as well that the contribution of the terms in $i\eps \Delta$ cancels. However, there is a slight difference compared to Lemma~\ref{LM:estimStabContinuity} since the second equation of \eqref{eq:quant} does not have exactly the structure of a continuity equation. 
    
    Let us provide some details, for the sake of completeness.
    Fix the label $\alpha$ and $\eps>0$, omitting the dependency in these parameters for a while. Set
    $$c=a-b, \ \ U=w-v,$$
    that satisfies
    \begin{align}\label{eq:diff-quanteq-plusdidee}
        \partial_t c+ (w+\lambda)\cdot \nabla c +\frac{1}{2}c\mathrm{div}(w)=-U\cdot \nabla a-\frac{1}{2}a \mathrm{div}(U) +i\frac{\eps}{2}\Delta c.
    \end{align}
    Let $\phi$ such that with $-\Delta \phi=c$ and $\Vert c \Vert_{\dot H^{-1}_x}=\Vert \nabla \phi \Vert_{L^2_x}$: we have
    \begin{align*}
        \frac{1}{2}\frac{\mathrm{d}}{\mathrm{d}t}\Vert c \Vert_{\dot H^{-1}_x}^2=\Re\langle \partial_t c, \phi \rangle,
    \end{align*}
    where $\langle \cdot, \cdot \rangle$ denotes the complex $L^2_x$ pairing. Each contribution of the above equation is estimated as follows. First, the new last term in the right-hand side of \eqref{eq:diff-quanteq-plusdidee} is skew-adjoint and makes no contribution since 
    \begin{align*}
    \Re \left \langle i\frac{\eps}{2}\Delta c, \phi \right\rangle =\frac{\eps}{2}\Re i\left \langle c, \Delta \phi \right\rangle=-\frac{\eps}{2}\Re (i \Vert c \Vert_{L^2_x}^2)=0.
    \end{align*}
    Second, similarly to the proof of Lemma~\ref{LM:estimStabContinuity}, the transport term on the left-hand side of \eqref{eq:diff-quanteq-plusdidee} yields
    \begin{align*}
    -\Re \Big\langle (w+\lambda)\cdot \nabla c &+\frac{1}{2}c\mathrm{div}(w), \phi \Big\rangle \\
    &=-\Re \int \Delta \phi \left((w+\lambda) \cdot \nabla \overline{\phi} +\frac{1}{2} \mathrm{div}(w) \overline{\phi}\right)\\
    &=\Re \int \nabla \phi\cdot \nabla w \nabla \overline{\phi}
    +\frac{1}{2} \int (w+\lambda)\cdot \nabla \vert \nabla \phi \vert^2  \\
    & \, +\frac{1}{2} \Re\int \nabla \phi \cdot \nabla(\mathrm{div}(w)) \overline{\phi}+\int \frac{1}{2}\mathrm{div}(w) \vert \nabla \phi \vert^2 \\
    &=\Re \int \nabla \phi\cdot \nabla w \nabla \overline{\phi}+\frac{1}{2} \Re\int \nabla \phi \cdot \nabla(\mathrm{div}(w)) \overline{\phi}.
    \end{align*}
   thanks to successive integration by parts an cancellation between the second and fourth terms. Thanks to Hölder inequality and Sobolev embedding ($1<d/2$), we get 
    \begin{align*}
        &\left\vert \Re \Big\langle (w+\lambda)\cdot \nabla c +\frac{1}{2}c\mathrm{div}(w), \phi \Big\rangle\right\vert \\
        & \lesssim  \Vert \nabla w \Vert_{L^\infty_x} \Vert \nabla \phi \Vert_{L^2_x}^2+ \Vert \nabla(\mathrm{div}(w)) \Vert_{L^d_x} \Vert \nabla \phi \Vert_{L^2_x} \Vert \phi \Vert_{L^{\frac{2d}{d-2}}_x} \\
        & \lesssim \left(\Vert \nabla w \Vert_{L^\infty_x}+ \Vert \nabla(\mathrm{div}(w)) \Vert_{L^d_x}\right)\Vert \nabla \Vert_{L^2_x}^2 \\
        & \lesssim \Vert \nabla w \Vert_{H^{k-1}_x}\Vert \nabla \phi \Vert_{L^2_x}^2,
    \end{align*}
    since $k-1>1+d/2$. Finally, the last source term in \eqref{eq:diff-quanteq-plusdidee} can be rewritten as
    \begin{align*}
        \left\langle U\cdot \nabla a+\frac{1}{2}a \mathrm{div}(U), \phi \right\rangle=\frac{1}{2}\int U \cdot \left( \nabla \overline{\phi}-a \nabla \overline{\phi} \right),
    \end{align*}
    so, by proceeding as before, we get
    \begin{align*}
        \left\vert \Re \left\langle U\cdot \nabla a-\frac{1}{2}a \mathrm{div}(U), \phi \right\rangle \right\vert &\lesssim \left( \Vert \nabla a \Vert_{L^d_x} \Vert \phi \Vert_{L^{\frac{2d}{d-2}}_x}+\Vert a \Vert_{L^\infty_x}\Vert \nabla \phi \Vert_{L^2_x} \right)\Vert U \Vert_{L^2_x} \\
        & \lesssim \Vert a \Vert_{H^{k-1}_x} \Vert U \Vert_{L^2_x} \Vert \nabla \phi \Vert_{L^2_x}.
    \end{align*}
    All in all, we obtain
    \begin{multline*}
        \frac{1}{2}\frac{\mathrm{d}}{\mathrm{d}t}\Vert b^{\eps, \alpha}-a^{\eps, \alpha} \Vert_{\dot H^{-1}_x} \\ \lesssim  \Vert \nabla w^{\eps, \alpha} \Vert_{H^{k-1}_x}\Vert b^{\eps, \alpha}-a^{\eps, \alpha} \Vert_{\dot H^{-1}_x}+\Vert a^{\eps, \alpha}\Vert_{H^{k-1}_x} \Vert v^{\eps, \alpha}-w^{\eps, \alpha} \Vert_{L^2_x}.
    \end{multline*}
    Using Gronwall lemma and then integrating in $\alpha$ yields the result.
    
\end{proof}

\paragraph{Bound in high regularity.} With this lemma at hand, we can argue almost exactly as in the Vlasov--Poisson case to obtain a bound in high regularity. Let $k_0>1+d/2$ and $k\geq k_0$. Consider a vector field $\uu^\eps \in L^1([0,T];\H^{k,\infty})$. We denote by $\boldsymbol{a}^\eps(t)$  the solution to the second equation of \eqref{eq:quant} with vector field $\uu^\eps$ and some initial condition $\boldsymbol{a}^\eps_0$. Consider the force field $-\nabla V^\eps$  solving 
        $$
            -\Delta V^\eps=\int_I |a^{\eps,\alpha}|^2 \, \D \mu(\alpha).
        $$
Owing to \eqref{eq:main-quantum} and since $H^{k-1}_x$ is an algebra (as $k-1>d/2$), we have
\begin{align*}
    \| \nabla V^\eps\|_{H^k_x} &\lesssim \left\Vert \int_I |a^{\eps,\alpha}|^2 \, \D \mu(\alpha)\right\Vert_{H^{k-1}_x}
    + \left\Vert\int_I |a^{\eps,\alpha}|^2 \, \D \mu(\alpha) \right\Vert_{\dot H^{-1}_x}\\
&\lesssim \| |\boldsymbol{a}^{\eps}|^2 \|_{L^1_\alpha H^{k-1}_x}+\Vert \boldsymbol{a}^\eps \Vert_{L^{\upinfty}_\alpha L^{2r}_x}^2\\
&\lesssim \Vert \boldsymbol{a}^{\eps}_0 \Vert_{L^{\upinfty}_\alpha H^{k-1}_x}^2   \exp\left(2 K \int_0^t\Vert \uu^\eps(s)\Vert_{\mathcal{H}^{k,\infty}} \D s \right)+\Vert \boldsymbol{a}^\eps \Vert_{L^{\upinfty}_\alpha L^{2r}_x}^2,
\end{align*}
where $r=2d/(d+2)$, thanks to the embedding $\dot H^{-1}_x \hookrightarrow L^r_x$. Since $2 \leq 2r <\infty$, interpolation and Sobolev embedding shows that $\Vert \boldsymbol{a}^\eps \Vert_{L^{\upinfty}_\alpha L^{2r}_x} \lesssim \Vert \boldsymbol{a}^\eps \Vert_{L^{\upinfty}_\alpha H^{k-1}_x}$,
so that we infer the required estimate for the force field $-\nabla V^\eps$.
Moreover it is clear that they are uniform with respect to $\eps$, provided that the initial data have norms that are uniformly bounded with respect to $\eps \in (0,1]$.
Assumption~\ref{ass:bound_high_regularity-rho} is therefore satisfied.

\paragraph{Stability at low regularity.} For what concerns stability in low regularity, 
we first set $N_{k_0}(\cdot)= \| \cdot\|_{L^{\upinfty}_\alpha H^{k_0-1}_x}$ for $k_0>d/2+2$,
and we assume that we are given $\uu^\eps, \vv^\eps$ and $\boldsymbol{a}^\eps_0, \boldsymbol{b}^\eps_0$ such that
		\begin{equation}
	\label{eq:uniform_bound_velocities-quant} \sup_{\eps \in (0,1]} \max\left( \sup_{t \in [0,T]} \Vert\uu^\eps(t)\Vert_{\mathcal{H}^{k_0,\infty}},\sup_{t \in [0,T]} \Vert\vv^\eps(t)\Vert_{\mathcal{H}^{k_0,\infty}}, \mathsf N_{k_0}(\boldsymbol{a}^\eps_0), \mathsf N_{k_0}(\boldsymbol{b}^\eps_0)  \right) \leq R,
		\end{equation}
        for some $R>0$. We denote by $\boldsymbol{a}^\eps$ (resp. $\boldsymbol{b}^\eps$) the solution to the second equation of \eqref{eq:quant} with vector field $\uu^\eps$ (resp. $\vv^\eps$) and initial condition $\boldsymbol{a}^\eps_0$ (resp. $\boldsymbol{b}^\eps_0$). Then, considering the force fields $-\nabla V^\eps$ and $-\nabla W^\eps$ solving 
        $$
            -\Delta V^\eps=\int_I |a^{\eps,\alpha}|^2 \, \D \mu(\alpha), \qquad     -\Delta W^\eps=\int_I |b^{\eps,\alpha}|^2 \, \D \mu(\alpha),
        $$
        the goal is to show that 
        $$\|\nabla V^\eps-\nabla W^\eps\|_{L^2_x} \lesssim \| |a^{\eps,\alpha}|^2 - |b^{\eps,\alpha}|^2\|_{L^1_\alpha \dot{H}^{-1}_x}.$$
        To this end, we aim at applying the stability estimate~\eqref{eq:main-quantum2}; there is a first small difference compared to the Vlasov--Poisson case, due to the  term  $\| \boldsymbol{a}^\eps \|_{L^\infty([0,t];L^\infty_\alpha H^{k-1}_x)}$  in the stability estimate: this is easily dealt with by combining the estimate \eqref{eq:main-quantum} and the bound \eqref{eq:uniform_bound_velocities-quant}, that yield 
        $$
        \| \boldsymbol{a}^\eps \|_{L^\infty([0,t];L^\infty_\alpha L^\infty_x)} \lesssim_R 1.
        $$
The other difference is that the density is a quadratic function of $a^\eps$, and another argument is thus required for the stability estimate for the force.
 Let $\varphi$ be a  smooth test function. 
 We have
\begin{align*}
 \langle |a^{\eps,\alpha}|^2 - |b^{\eps,\alpha}|^2 ,   \varphi \rangle
 &= \langle a^{\eps,\alpha} - b^{\eps,\alpha},  (\overline{a}^{\eps,\alpha} + \overline{b}^{\eps,\alpha}) \varphi \rangle  \\
 &\leq  \| a^{\eps,\alpha}- b^{\eps,\alpha}\|_{\dot H^{-1}_x} \left( \| C^{\eps,\alpha} \nabla \varphi  \|_{L^2_x} +  \|\nabla C^{\eps,\alpha} \varphi \|_{L^2_x}   \right),
\end{align*}
where we have denoted $C^{\eps,\alpha}\vcentcolon= \overline{a}^{\eps,\alpha} + \overline{b}^{\eps,\alpha}$. By the $H^{k-1}_x$ estimate \eqref{eq:main-quantum}, since $k-1>d/2$, by Sobolev embedding it holds $\| C^{\eps,\alpha} \nabla \varphi  \|_{L^2_x} \lesssim_R \| \nabla \varphi  \|_{L^2_x}$. 

On the other hand, by H\"older's inequality and Sobolev embedding, there holds $\|\nabla C^{\eps,\alpha} \varphi \|_{L^2_x} \leq \|\nabla C^{\eps,\alpha} \|_{L^d_x} \| \nabla \varphi \|_{L^2_x}$. If  $k-2 \geq d/2$, $H^{k-2}_x$ embeds into $L^\infty_x \cap L^2_x$ and thus in particular in $L^d_x$. Otherwise, if $k-2<d/2$ (so that $d-2(k-2)>0$), also by Sobolev embedding, $\| \nabla C^{\eps,\alpha}\|_{L^q_x}$ is uniformly controlled for all $q \in [2, 2d/(d-2(k-2)))$;  we further note that $d < 2d/(d-2(k-2))$ since $k-1>d/2$.
All in all, we deduce that
\begin{equation}
    \label{eq:trucH-1}
\| |a^{\eps,\alpha}|^2 - |b^{\eps,\alpha}|^2 \|_{\dot{H}^{-1}_x} \lesssim_R \| a^{\eps,\alpha}- b^{\eps,\alpha}\|_{\dot H^{-1}_x},
\end{equation}
and the estimate is uniform in $\eps$.
Finally, we can  use
\eqref{eq:main-quantum2} to obtain $L^2_x$ stability for the force field:
\begin{equation}
    \label{eq:trucH-1-re}
\| \nabla V^{\eps}(t)- \nabla W^{\eps}(t) \|_{ \dot{H}^{-1}_x} \lesssim_{t,R}  \| \boldsymbol{a}^{\eps,\alpha}_0- \boldsymbol{b}^{\eps,\alpha}_0\|_{L^1_\alpha\dot H^{-1}_x} + \int_0^t \| \uu^{\eps} - \vv^\eps\|_{L^1_\alpha L^2_x} \, \D s. 
\end{equation}    
Assumption~\ref{ass:stability-rho} follows.

\bigskip

As a consequence, we can apply Theorem~\ref{thm-LWPmultiphase-densityvelocity-re} and as already explained, the statement which we obtain is uniform with respect to $\eps$: this yields the existence part of Theorem~\ref{thm:semiclassical}. It remains to justify the strong convergence to the solution $(\rhorho, \vv)$ of the Vlasov--Poisson equation \eqref{eq:VP} associated with the initial condition $({\rhorho}_{0} , \vv_{0})$. To this end, we use stability estimates\footnote{Without loss of generality, we can always suppose that $(\rhorho, \vv)$ is also defined on $[0,T]$, reducing $T$ if necessary.}.
First note that $\sqrt{\rho^\alpha}$ satisfies the equation
$$
\partial_t \sqrt{\rho^\alpha} + v^\alpha\cdot \nabla \sqrt{\rho^\alpha} + \frac{1}{2} \sqrt{\rho^\alpha} \Div (v^{\alpha})  = 0.
$$
We rely on
\begin{itemize}
    \item a $L^2_x$ stability estimate~\eqref{eq:stability-vv} to control the difference of velocities $$\nabla S^{\eps, \alpha}- v^\alpha,$$ 
    obtained as in the proof of Theorem~\ref{thm:existence};
    \item a $\dot{H}^{-1}_x$ stability estimate to control $a^{\eps,\alpha}- \sqrt{\rho^\alpha}$, similar to~\eqref{eq:main-quantum2} -- there is a slight difference in the treatment of the term $ i\frac{\eps}{2}  \Delta  a^{\eps,\alpha}$ in the ``transport'' equation for $a^{\eps,\alpha}$, which is  seen here as a  uniformly bounded forcing term (whose $L^\infty_\alpha \dot{H}^{-1}_x$ norm is controlled by $\lesssim \eps$, using the uniform estimates for $a^{\eps,\alpha}$);
    \item another $\dot{H}^{-1}_x$ stability estimate to control the difference of densities $$|a^{\eps,\alpha}|^2- {\rho^\alpha},$$ similar to~\eqref{eq:trucH-1}.
\end{itemize}
We deduce
\begin{multline}
\label{eq-arho-quantitatif-0}
\sup_{[0,T]} \lVe \nabla \boldsymbol{{S}}^{\eps} - \vv \rVe_{L^{\upinfty}_\alpha L^2_x} + \sup_{[0,T]}  \lVe \boldsymbol{a}^{\eps} - \sqrt{\rhorho} \rVe_{L^1_\alpha \dot{H}^{-1}_x} \\
\lesssim \lVe \nabla \boldsymbol{{S}}^\eps_0 - \vv_0 \rVe_{L^{\upinfty}_\alpha L^2_x}  + \lVe \boldsymbol{a}^\eps_0 - \sqrt{\rhorho_0} \rVe_{L^1_\alpha \dot{H}^{-1}_x} + \eps.
\end{multline}
Arguing similarly at the level of one derivative in $L^2_x$ for the velocities and $L^2_x$ for the densities, we also obtain
\begin{multline}
\label{eq-arho-quantitatif}
\sup_{[0,T]} \lVe \nabla( \nabla \boldsymbol{{S}}^{\eps} - \vv )\rVe_{L^{\upinfty}_\alpha L^2_x} + \sup_{[0,T]}  \lVe \boldsymbol{a}^{\eps} - \sqrt{\rhorho} \rVe_{L^1_\alpha L^2_x} \\\lesssim \lVe \nabla( \nabla \boldsymbol{{S}}^\eps_0 - \vv_0 )\rVe_{L^{\upinfty}_\alpha L^2_x}  + \lVe \boldsymbol{a}^\eps_0 - \sqrt{\rhorho_0} \rVe_{L^1_\alpha L^2_x} + \eps.
\end{multline}
Arguing by recursion, we get similarly that for all $j=1,\cdots, k-3$,
\begin{multline}
\sup_{[0,T]} \lVe \nabla(\nabla \boldsymbol{{S}}^{\eps} - \vv )\rVe_{L^{\upinfty}_\alpha H^j_x} + \sup_{[0,T]}  \lVe \boldsymbol{a}^{\eps} - \sqrt{\rhorho} \rVe_{L^1_\alpha H^{j}_x} \\\lesssim \lVe \nabla (\nabla \boldsymbol{{S}}^\eps_0 - \vv_0 )\rVe_{L^{\upinfty}_\alpha {H}^j_x}  + \lVe \boldsymbol{a}^\eps_0 - \sqrt{\rhorho_0} \rVe_{L^1_\alpha H^{j}_x} + \eps.
\end{multline}
Note that for what concerns higher regularity, by interpolation between \eqref{eq-arho-quantitatif-0}--\eqref{eq-arho-quantitatif} and the uniform bound for the top-order derivatives,  we also obtain that 
$$
\sup_{[0,T]} \lVe \nabla \boldsymbol{{S}}^{\eps} - \vv \rVe_{L^{\upinfty}_\alpha H^s_x} +\sup_{[0,T]}  \lVe \boldsymbol{a}^{\eps} - \sqrt{\rhorho} \rVe_{L^1_\alpha H^{s-1}_x} \longrightarrow 0,
$$
as $\eps \to 0$, for all $s \in [0,k)$. This proves in particular that the family $(\boldsymbol{a}^{\eps}, \nabla \boldsymbol{{S}}^{\eps}-\boldsymbol\lambda)$ converges in $L^\infty([0,T]; L^\infty_\alpha H^{k-2}_x ) \times L^\infty([0,T]; \mathcal{H}^{k-1,\infty}) $ 
to $(\sqrt{\rhorho}, \vv-\boldsymbol\lambda)$.
The first part of the convergence statement in the theorem follows.

\bigskip

 It remains to prove~\eqref{eq:conv-wkb}. In the following, we dismiss the time variable as it plays the role of a parameter in the estimates (ultimately taking the supremum in time in the forthcoming estimates). Let $\varphi=\varphi(x,v)$ be a smooth test function such that $\|\langle y \rangle^2 \mathcal{F}_{v \rightarrow y} \varphi\|_{L^\infty_x L^1_y} \leq 1$, where we denote by $\mathcal{F}_{v \rightarrow y}$ the partial Fourier transform from $v$ to $y$. Note that by the Fourier inverse formula and the multiphasic representation for $f$, there holds
 \begin{align*}
 \int_{\R^d \times \R^d} f(x,v)   &\varphi(x,v)\, \mathrm{d}x \,  \mathrm{d}v  \\
 &=   \int_{\R^d} \int_I \rho^\alpha(x) \varphi(x,v^\alpha(x))\, \mathrm{d}\mu(\alpha) \, \mathrm{d}x  \\
 &=   \frac{1}{(2\pi)^d} \int_{\R^d} \int_{\R^d} \int_I \rho^\alpha(x) e^{i v^\alpha(x)\cdot y } \mathcal{F}_{v \rightarrow y}\varphi(x,y) \, \mathrm{d}\mu(\alpha) \, \mathrm{d}x  \, \mathrm{d}y.
 \end{align*}
  Therefore we have:
 \begin{align*}
&\int_{\R^d \times\R^d} \left(W[\psi^{\eps,\alpha}](x,v) - f(x,v)  \right) \varphi(x,v)\, \mathrm{d}x \,  \mathrm{d}v \\
&= \frac{1}{(2\pi)^d} \int_{\R^d \times \R^d}  \int_I \mathcal{F}_{v \rightarrow y}\varphi(x,y) \\
&\times \left[a^{\eps,\alpha}(x+\eps y /2)\overline{a^{\eps,\alpha}}(x-\eps y /2)  e^{\frac{i}{\eps}(S^{\eps,\alpha}(x+\eps y/2) - S^{\eps,\alpha}(x-\eps y/2)) } - \rho^\alpha(x) e^{i v^\alpha(x)\cdot y } \right] \\
&\qquad \qquad \qquad \qquad \qquad \mathrm{d}\mu(\alpha) \, \mathrm{d}x \,  \mathrm{d}y \\
&= \frac{1}{(2\pi)^d} \int_{\R^d \times \R^d}  \int_I \mathcal{F}_{v \rightarrow y}\varphi(x,y) \left[a^{\eps,\alpha}(x+\eps y /2)\overline{a^{\eps,\alpha}}(x-\eps y /2) - \rho^\alpha(x) \right] \\
&\qquad \qquad \qquad  \qquad \times e^{\frac{i}{\eps}(S^{\eps,\alpha}(x+\eps y/2) - S^{\eps,\alpha}(x-\eps y/2)) }  \, \mathrm{d}\mu(\alpha) \, \mathrm{d}x \,  \mathrm{d}y \\
&\qquad + \frac{1}{(2\pi)^d} \int_{\R^d \times \R^d}  \int_I \mathcal{F}_{v \rightarrow y}\varphi(x,y) \rho^\alpha(x) \\
&\qquad \qquad \qquad \qquad \times \left[ e^{\frac{i}{\eps}(S^{\eps,\alpha}(x+\eps y/2) - S^{\eps,\alpha}(x-\eps y/2)) } - e^{i v^\alpha(x)\cdot y } \right]\, \mathrm{d}\mu(\alpha)  \, \mathrm{d}x \,  \mathrm{d}y \\
&=\vcentcolon J_1 + J_2.
 \end{align*}
Let us estimate these two terms separately.

\bigskip

 \noindent\textbf{Estimate on $J_1$}. We first write
 \begin{align*}
     J_1&= \frac{1}{(2\pi)^d} \int_{\R^d \times \R^d}  \int_I \mathcal{F}_{v \rightarrow y}\varphi(x,y) \left[\vert a^{\eps,\alpha}(x) \vert^2- \rho^\alpha(x) \right] \\
     &\qquad \qquad \qquad \qquad \qquad \times e^{\frac{i}{\eps}(S^{\eps,\alpha}(x+\eps y/2) - S^{\eps,\alpha}(x-\eps y/2)) }  \, \mathrm{d}\mu(\alpha) \, \mathrm{d}x \,  \mathrm{d}y \\
     &+  \frac{1}{(2\pi)^d} \int_{\R^d \times \R^d}  \int_I \mathcal{F}_{v \rightarrow y}\varphi(x,y) \left[a^{\eps,\alpha}(x+\eps y /2)\overline{a^{\eps,\alpha}}(x-\eps y /2) - \vert a^{\eps,\alpha}(x) \vert^2 \right] \\
      & \qquad \qquad \qquad \qquad \qquad \times e^{\frac{i}{\eps}(S^{\eps,\alpha}(x+\eps y/2) - S^{\eps,\alpha}(x-\eps y/2)) } \, \mathrm{d}\mu(\alpha) \, \mathrm{d}x \,  \mathrm{d}y \\
     & =\vcentcolon J_{1,1}+J_{1,2}.
 \end{align*}
 We obtain directly for $J_{1,1}$ that
 \begin{align*}
     \vert J_{1,1} \vert 
     &\lesssim \int_{\R^d \times \R^d}  \int_I \Vert \mathcal{F}_{v \rightarrow y}\varphi(\cdot,y) \Vert_{L^{\upinfty}_x} \left\vert \vert a^{\eps,\alpha}(x) \vert^2- \rho^\alpha(x) \right\vert   \, \mathrm{d}\mu(\alpha) \, \mathrm{d}x \,  \mathrm{d}y \\
     &\lesssim \Vert \vert \boldsymbol{a}^{\eps} \vert^2-\rhorho \Vert_{L^1_\alpha L^1_x}.
 \end{align*}
For $J_{1,2}$, we first write
\begin{multline*}
    a^{\eps,\alpha}(x+\eps y /2)\overline{a^{\eps,\alpha}}(x-\eps y /2) - \vert a^{\eps,\alpha}(x) \vert^2  \\
    = \big[ a^{\eps,\alpha}(x+\eps y /2)-a^{\eps,\alpha}(x)\big] \overline{a^{\eps,\alpha}}(x-\eps y /2) 
     \quad +a^{\eps,\alpha}(x) \big[ \overline{a^{\eps,\alpha}}(x-\eps y /2)-\overline{a^{\eps,\alpha}}(x)\big] ,
\end{multline*}
and use Cauchy--Schwarz inequality combined  with the mean value inequality 
$$\Vert a^{\eps,\alpha}(\cdot \pm\eps y/2)-a^{\eps,\alpha} \Vert_{L^2_x} \leq \frac{\eps \vert y \vert}{2} \Vert \nabla a^{\eps,\alpha} \Vert_{L^2_x},$$ 
to get for all $y \in \R^d$ that
\begin{align*}
    & \left\Vert a^{\eps,\alpha}(\cdot+\eps y /2)\overline{a^{\eps,\alpha}}(\cdot-\eps y /2) - \vert a^{\eps,\alpha} \vert^2\right\Vert_{L^1_x} 
     \\
     &\qquad\leq \Vert a^{\eps,\alpha}\Vert_{L^2_x}\left( \left\Vert  a^{\eps,\alpha}(\cdot+\eps y /2)-a^{\eps,\alpha}\right\Vert_{L^2_x}  + \left\Vert  a^{\eps,\alpha}(\cdot-\eps y /2)-a^{\eps,\alpha}\right\Vert_{L^2_x}\right) \\
     &\qquad \leq \eps \vert y \vert \Vert \boldsymbol{a}^\eps\Vert_{L^2_x} \Vert \nabla \boldsymbol{a}^{\eps}\Vert_{L^2_x}.
\end{align*}
As a consequence, we get
\begin{align*}
    \vert J_{1,2} \vert &\lesssim \int_{\R^d }  \|\mathcal{F}_{v \rightarrow y}\varphi(\cdot,y) \Vert_{L^{\upinfty}_x} \\
    &\qquad \qquad \times \int_{I\times \R^d}  \left\vert a^{\eps,\alpha}(x+\eps y /2)\overline{a^{\eps,\alpha}}(x-\eps y /2) - \vert a^{\eps,\alpha}(x) \vert^2\right\vert   \, \mathrm{d}\mu(\alpha) \, \mathrm{d}x \,  \mathrm{d}y \\
    & \lesssim \eps \left(\int_{\R^d }  \vert y \vert \Vert \mathcal{F}_{v \rightarrow y}\varphi(\cdot,y) \Vert_{L^{\upinfty}_x} \right) \int_{I} \Vert a^{\eps,\alpha}\Vert_{L^2_x} \Vert \nabla a^{\eps,\alpha}\Vert_{L^2_x} \, \mathrm{d}\mu(\alpha) \\
    & \lesssim \eps \Vert   \boldsymbol{a}^\eps\Vert_{L^{\upinfty}_\alpha H^1_x}^2.
\end{align*}
We have therefore obtained 
\begin{align*}
    \vert J_1 \vert \lesssim \Vert \vert \boldsymbol{a}^\eps \vert^2-\rhorho \Vert_{L^1_\alpha L^1_x}+\eps \Vert   \boldsymbol{a}^\eps\Vert_{L^{\upinfty}_\alpha H^1_x}^2.
\end{align*}

\bigskip

\noindent\textbf{Estimate on $J_2$}. We start by observing that
 \begin{multline*}
     \left\vert e^{\frac{i}{\eps}(S^{\eps,\alpha}(x+\eps y/2) - S^{\eps,\alpha}(x-\eps y/2)) } - e^{i v^\alpha(x)\cdot y } \right\vert \\
     \leq  \left\vert \frac{1}{\eps}(S^{\eps,\alpha}(x+\eps y/2) - S^{\eps,\alpha}(x-\eps y/2))-v^\alpha(x)\cdot y \right\vert,
 \end{multline*}
 and, using the formula $$\frac{1}{\eps}(S^{\eps,\alpha}(x+\eps y/2) - S^{\eps,\alpha}(x-\eps y/2))= \int_{-1/2}^{1/2} y \cdot \nabla S^{\eps,\alpha} (x+\theta \eps y)\, \mathrm{d}\theta ,$$ 
 we can also write 
 \begin{align*}
     &\frac{1}{\eps}(S^{\eps,\alpha}(x+\eps y/2) - S^{\eps,\alpha}(x-\eps y/2))-v^\alpha(x)\cdot y \\
     &= y\cdot \left(\nabla S^{\eps,\alpha}(x)  -v^\alpha(x)\right) + \int_{-1/2}^{1/2} y\cdot \big( \nabla S^{\eps,\alpha} (x+\theta \eps y)-\nabla S^{\eps,\alpha}(x) \big)\, \mathrm{d}\theta. 
 \end{align*}
 As a consequence, we obtain by the mean value inequality
 \begin{multline*}
     \left\vert e^{\frac{i}{\eps}(S^{\eps,\alpha}(x+\eps y/2) - S^{\eps,\alpha}(x-\eps y/2)) } - e^{i v^\alpha(x)\cdot y } \right\vert \\
     \lesssim  \vert y\vert \left\vert \nabla S^{\eps,\alpha}(x)  -v^\alpha(x)\right\vert  +\eps \vert y \vert^2 \Vert \nabla^2 S^{\eps,\alpha} \Vert_{L^{\upinfty}_x} .
 \end{multline*}
Plugging this pointwise estimate into $J_2$, we get with the same principle as for $J_1$:
\begin{align*}
\vert J_2 \vert &\lesssim \int_{\R^d \times \R^d}  \int_I \vert \mathcal{F}_{v \rightarrow y}\varphi(x,y) \vert \vert \rho^\alpha(x) \vert  \\
&\qquad \qquad \qquad \times \left[ \vert y\vert \left\vert \nabla S^{\eps,\alpha}(x)  -v^\alpha(x)\right\vert  +\eps \vert y \vert^2 \Vert \nabla^2 S^{\eps,\alpha} \Vert_{L^{\upinfty}_x}\right]\, \mathrm{d}\mu(\alpha)  \, \mathrm{d}x \,  \mathrm{d}y  \\
& \lesssim \Vert \rhorho \Vert_{L^{\upinfty}_\alpha L^2_x} \Vert \nabla \boldsymbol{S}^\eps  -\vv\Vert_{L^1_\alpha L^2_x}+ \eps \Vert \rhorho \Vert_{L^1_\alpha L^1_x}\Vert \nabla^2 \boldsymbol{S}^\eps \Vert_{L^{\upinfty}_\alpha L^{\upinfty}_x} .
\end{align*}
Gathering the estimates for $J_1$ and $J_2$, we have therefore proven that for all smooth test functions $\varphi(x,v)$ such that $\|\langle y \rangle^2 \mathcal{F}_{v \rightarrow y} \varphi\|_{L^\infty_x L^1_y} \leq 1$, there holds
\begin{align*}
    &\left\vert \int_{\R^d \times\R^d} \left(W[\psi^{\eps,\alpha}](x,v) - f(x,v)  \right) \varphi(x,v)\, \mathrm{d}x \,  \mathrm{d}v  \right\vert \\
    & \lesssim \Vert \vert \boldsymbol{a}^\eps \vert^2-\rhorho \Vert_{L^1_\alpha L^1_x}+\eps \Vert   \boldsymbol{a}^\eps\Vert_{L^{\upinfty}_\alpha H^1_x}^2  \\
    &\quad + \Vert \rhorho \Vert_{L^{\upinfty}_\alpha L^2_x} \Vert \nabla \boldsymbol{S}^\eps  -\vv\Vert_{L^1_\alpha L^2_x}+ \eps \Vert \rhorho \Vert_{L^1_\alpha L^1_x}\Vert \nabla^2 \boldsymbol{S}^\eps \Vert_{L^{\upinfty}_\alpha L^{\upinfty}_x}.
\end{align*}
Moreover, since
$\vert a^{\eps,\alpha} \vert^2-\rho^\alpha = \Re \left[ (a^{\eps,\alpha} - \sqrt{\rho^\alpha} ) (\overline{a}^{\eps,\alpha}+ \sqrt{\rho^\alpha})  \right]$, we have
$$
\Vert \vert \boldsymbol{a}^\eps \vert^2-\rhorho \Vert_{L^1_\alpha L^1_x}
\lesssim (\| \boldsymbol{a}^\eps\|_{L^{\upinfty}_\alpha L^2_x} + \| \rhorho\|_{L^{\upinfty}_\alpha L^1_x}^{1/2}) \Vert  \boldsymbol{a}^\eps -\sqrt{\rhorho} \Vert_{L^1_\alpha L^2_x},
$$
and therefore
\begin{align*}
    &\left\vert \int_{\R^d \times\R^d} \left(W[\psi^{\eps,\alpha}](x,v) - f(x,v)  \right) \varphi(x,v)\, \mathrm{d}x \,  \mathrm{d}v  \right\vert \\
    & \lesssim (\| \boldsymbol{a}^\eps\|_{L^{\upinfty}_\alpha L^2_x} + \| \rhorho\|_{L^{\upinfty}_\alpha L^1_x}^{1/2}) \Vert  \boldsymbol{a}^\eps -\sqrt{\rhorho} \Vert_{L^1_\alpha L^2_x} \\
    & \quad +\eps \Vert   \boldsymbol{a}^\eps\Vert_{L^{\upinfty}_\alpha H^1_x}^2 + \Vert \rhorho \Vert_{L^{\upinfty}_\alpha L^2_x} \Vert \nabla \boldsymbol{S}^\eps  -\vv\Vert_{L^1_\alpha L^2_x}+ \eps \Vert \rhorho \Vert_{L^1_\alpha L^1_x}\Vert \nabla^2 \boldsymbol{S}^\eps \Vert_{L^{\upinfty}_\alpha L^{\upinfty}_x}.
\end{align*}
We now observe that $\boldsymbol{a}^\eps$ is bounded in $L^\infty_\alpha H^1_x$ and that  $\nabla^2 \boldsymbol{S}^\eps=$ is bounded in  $L^\infty_\alpha L^\infty_x$ by Sobolev embedding, uniformly in $\eps$, while the prefactors not depending on $\eps$ are finite by assumption. This is enough to conclude, by using estimates \eqref{eq-arho-quantitatif-0}--\eqref{eq-arho-quantitatif} and taking the supremum over the set of test functions as well as the supremum in time on $[0,T]$.
\end{proof}

\section{Instability of homogeneous profiles for a class of multiphasic systems}\label{Section-FrameworkInstab}

We now turn to the last part of this chapter, that is dedicated to the study of nonlinear instability of homogeneous equilibria to Vlasov--Poisson type systems. We put forward an abstract instability framework, with the aim to apply it to several systems at once.

Following the general philosophy of this work, the actual aim is to establish instability for 
multiphasic systems; namely, we focus in this section on a class of density--velocity fields under the general form~\eqref{eq:general_system_intro}:
\begin{equation}
\label{eq:system_of_study}
\left\{
\begin{aligned}
\partial_t \rho^\alpha + \Div (\rho^\alpha v^\alpha) &= 0,\\
\partial_t v^\alpha + (v^\alpha \cdot \nabla) v^\alpha &= F^\alpha[\rhorho,\vv],
\end{aligned}
\right.
\end{equation}
where the force field $\FF[\rhorho,\vv]$ may depend nonlinearly on the whole families $(\rhorho,\vv)$.  Our main result shows that, under abstract assumptions on $\FF$ encompassing many
models of interest, \emph{linear instability of homogeneous profiles implies nonlinear Lyapunov instability}.

\medskip

The homogeneous profiles under consideration are of the form
$$\rho^\alpha \equiv 1, \qquad  v^\alpha \equiv \mathcal{V}(\alpha),$$ 
and we therefore assume that they are stationary solutions, that is, we impose 
$$\FF[\boldsymbol 1, \bm{\mathcal{V}}] \equiv 0.$$ 
Very loosely speaking, we intend to prove the following type of results: 
\begin{center}
``If there exists a exponentially growing in time solution to the {\bf linearized} equations about the homogeneous profile $(\boldsymbol 1, \bm{\mathcal{V}})$, then there also exists such a solution for the full nonlinear equations~\eqref{eq:system_of_study}''.
\end{center}

\medskip

The abstract assumptions on $\FF$ will be designed to allow nonlocal dependencies in space (e.g.\ through elliptic
couplings) as well as analytic nonlinearities. The motivation directly stems  from the examples of the Vlasov--Poisson system, for electrons and  for ions (recall systems \eqref{eq:VP} and \eqref{eq:VPions}): we will show in Section~\ref{sec:insta-appli} how the abstract result applies in these cases, on their multiphasic formulation.  This first allows to recover classical results on two-fluid plasmas (recall Theorem~\ref{eq:insta-kin-2} from \cite{CGG}) and effortlessly generalize them to an arbitrary number of phases. In a second time, we draw consequences on the original kinetic formulation.

Let us stress that even for the Vlasov--Poisson equation for electrons, the result is original, as our framework allows to deal with stationary profiles that only have the regularity of a measure. 

\medskip

We introduce the abstract setting in Section \ref{Section-assump-instab}, with basically four assumptions to be satisfied by  the force field $\FF$. 
 Loosely speaking, we assume that there holds:
\begin{enumerate}
    \item A \underline{multilinear expansion at any order for the force field $\FF$}, with a suitable control on the remainder at high regularity. In the examples we have in mind, the value of $\FF$ at time $t$ only depends on $\rhorho$ and $\vv$ at time~$t$ but in general the dependence will not be polynomial (this is typical, for instance, for the case of the Vlasov--Poisson system for ions).
    \item A \underline{stability property of the force field $\FF$ at high regularity}. Consistently with Theorem~\ref{thm-LWPmultiphase-densityvelocity},
  controlling $\ell$ derivatives of $\rhorho$ and $\ell+1$ derivatives of $\vv$ yields a control of
  $\ell+1$ derivatives of $\FF[\rhorho,\vv]$ (in the appropriate norms).
    \item  An \underline{exponential instability property} for the linearized system, that is the existence of an exponential growing mode in time solution to this system, with a large enough growing rate with respect to the associated semigroup.
    \item A structural condition ensuring \underline{compatibility with real-valued inputs} (this is mainly a technical point, needed in the nonlinear construction).
\end{enumerate}
Within this framework, the abstract instability result is stated in Theorem~\ref{thm:instability}. The proof, which is given in Section~\ref{sec:instability}, is an adaptation to the multiphasic framework of Grenier's nonlinear instability method \cite{Gr00}.

\subsection{Setting and assumptions}\label{Section-assump-instab}
Let $(I, \mu)$ be a set of labels. From now on, we write $$\rho^\alpha = 1 + r^\alpha, \ \  v^\alpha = \mathcal{V}(\alpha) + w^\alpha,$$ 
where $\mathcal{V}: I \mapsto \R$ is a fixed profile such that $\mathcal{V} \in L^p_\alpha$ with $p \in [1, \infty)$, and recast~\eqref{eq:system_of_study} for the perturbation $(\rr, \ww)$. For simplicity, we denote by $Z^\alpha$ the pair $(r^\alpha, w^\alpha)$, with $$\ZZ \vcentcolon= (\rr,\ww), \ \ \GG[\rr, \ww] \vcentcolon= \FF[\boldsymbol 1 + \rr, \bm{\mathcal{V}}+ \ww],$$
and we thus ask that $\GG[\boldsymbol 0] = 0$ for $(\boldsymbol 1,\bm{\mathcal{V}})$ to be a stationary state. From \eqref{eq:system_of_study}, this leads to the system
\begin{equation}
\label{eq:system_of_study_translated}
\left\{
\begin{aligned}
\partial_t r^\alpha + \Div \Big((1 + r^\alpha)(\mathcal{V}(\alpha) + w^\alpha)\Big) &= 0,\\
\partial_t w^\alpha + \Big((\mathcal{V}(\alpha) + w^\alpha) \cdot \nabla\Big) w^\alpha &= G^\alpha[\rr,\ww].
\end{aligned}
\right.
\end{equation}
We also use the following notation: for all $\ell \in \N$, we set
\begin{equation}\label{def:NORM-Hl-INSTABproof}
\| \ZZ \|_\ell \vcentcolon= \esssup_{\alpha \in I} \| Z^\alpha \|_{H^{\ell}_x \times H^{\ell + 1}_x}.
\end{equation}
We will implicitly use the fact that the previous equation admits a complex extension but that, when restricted on real-valued functions, everything stays real-valued (in particular the force field $\GG$). Our assumptions state as follows.

\begin{enumerate}[label=(insta-H\arabic*)]
	\item \underline{Expansion for the force field}. \label{item:ass_insta_expansion}
	
	For all $N \in \N$, $\GG$ can be decomposed as
	\begin{equation}
	\label{eq:decomposition_G}
	\GG[\ZZ] = \sum_{k=1}^N \GG_k[\ZZ, \dots, \ZZ] + \boldsymbol R_N[\ZZ],
	\end{equation}
	where:
	\begin{enumerate}
		\item For all $k \in \N^*$, $(\ZZ_1, \dots \ZZ_k) \mapsto \GG_k[\ZZ_1, \dots, \ZZ_k]$ is $k$-linear and symmetric, and for all $\vec n = (n_1, \dots n_k) \in (\Z^d)^k$ and all $\boldsymbol{Y_1}, \dots, \boldsymbol{Y_k} \in L^\infty_\alpha$, we can write
		\begin{multline*}
		\GG_k[\boldsymbol{Y_1} \exp(i n_1 \cdot x), \dots , \boldsymbol{Y_k}(\alpha) \exp(i n_k \cdot x)] \\
        = P_{k,\vec n}(\boldsymbol{Y_1}, \dots, \boldsymbol{Y_k}) \exp(i(n_1 + \dots + n_k) \cdot x),
		\end{multline*}
		where $P_{k,\vec n}$ is a continuous $k$-linear and symmetric functional from $(L^\infty_\alpha)^k$ to $L^\infty_\alpha$.
		\item For $\ell \in \N$ large enough, there exists $\eps >0$ such that for all $\ZZ$ satisfying $\|  \ZZ\|_{\ell} \leq \eps$, the remainder $\boldsymbol R_N$ satisfies
		\begin{equation*}
		\esssup_{\alpha \in I}  \| R_N^\alpha[\ZZ] \|_{H^{\ell + 1}_x} \leq C_{N,\ell,\eps} \| \ZZ \|^{N+1}_\ell,
		\end{equation*}
		for some $C_{N, \ell, \eps}>0$.
	\end{enumerate}

\begin{Rem}
    An essential aspect of this assumption, and hence of our result in view of the applications, is that it does not require any uniformity of the continuity of the $k$-linear operator $P_{k,\vec n}$ with respect to $k$ and $\vec n$, nor of $C_{N,\ell,\eps}$ with respect to $N$.
\end{Rem}

\item\label{item:ass_insta_StabHIGHREG} \underline{Stability of the force at high regularity}.

For $\ell \in \N$ large enough, there is a smooth nondecreasing function $\G_\ell : \R_+  \to \R_+ \cup \{+ \infty\}$, finite close to zero, and such that for any smooth $\ZZ_1$ and $\ZZ_2$, we have
\begin{equation}
\esssup_{\alpha \in I}  \| G^\alpha[\ZZ_1] - G^\alpha[\ZZ_2] \|_{H^{\ell + 1}_x} \leq \G_\ell\Big(\max(\| \ZZ_1 \|_\ell, \| \ZZ_2 \|_\ell)\Big) \| \ZZ_1 - \ZZ_2 \|_\ell.
\end{equation}

\item \underline{About the linear system}. \label{item:ass_insta_LINEARIZED}

The linearized system
\begin{equation}
\label{eq:system_of_study_linear}
\left\{
\begin{aligned}
\partial_t r^\alpha + \mathcal{V}(\alpha) \cdot \nabla r^\alpha+ \Div (w^\alpha) &=0, \\
\partial_t w^\alpha + (\mathcal{V}(\alpha) \cdot \nabla) w^\alpha &= G_1^\alpha[\ZZ],
\end{aligned}
\right.
\end{equation}
where $G_1^\alpha$ is defined in~\eqref{eq:decomposition_G} satisfies the following two properties: 
\begin{enumerate}
	\item The semigroup $(S_t)$ associated with the linearized system is well-defined for all $t\geq 0$ on configurations of the form $Z^\alpha(x) = Y_0(\alpha) \exp(i n \cdot x)$ with $n \in \Z^d$ and $Y_0 \in L^\infty_\alpha$. Moreover, calling $Y = Y(t,\alpha)$ the function such that,
	\begin{equation*}
	    S_t( Y_0(\alpha) \exp(in\cdot x)) = Y(t,\alpha) \exp(i n\cdot x),
	\end{equation*}
	there exists $\theta_0 >0$ independent of $n$ and $C_{n} >0$ such that
	\begin{equation*}
	    \| Y(t) \|_{L^{\upinfty}_\alpha} \leq C_n \| Y_0 \|_{L^{\upinfty}_\alpha} \exp(\theta_0 t).
	\end{equation*}
	\item The linearized system admits an exponential growing mode, that is, there exists $n \in \Z^d\backslash\{0\}$ and a non-zero solution of the form
	\begin{equation*}
	X(\alpha) \exp(in\cdot x) \exp( \lambda t), \ \ t>0,
	\end{equation*}
	 of growing rate $\theta \vcentcolon= \Re (\lambda) > \theta_0/2$, and $\boldsymbol{X}\in L^\infty_\alpha \backslash\{ 0 \}$.
\end{enumerate}
\begin{Rem}
    Note that here as well, we do not assume any uniformity of $C_{n}$ with respect to $n$.
\end{Rem}

\item \underline{Real force fields}. \label{item:ass_insta_RealForce}

\begin{enumerate}
    \item When $\ZZ$ has real values, then $\GG[\ZZ]$ has real values.
    \item For all $N \in \N^*$, when $\ZZ$ is real, then $\boldsymbol R_N[\ZZ]$ also is.
    
    \item For all $k \in \N^*$, $\vec n = (n_1, \dots n_k) \in (\Z^d)^k$ and $Y_1, \dots, Y_k \in L^\infty$, 
    \begin{equation*}
        P_{k, \vec n} (\bar Y_1, \dots, \bar Y_k) = \overline{P_{k, \vec n}(Y_1, \dots, Y_k)}.
    \end{equation*}
    In particular, whenever $\ZZ_1, \dots, \ZZ_k$ is real, $\GG_k[\ZZ_1, \dots, \ZZ_k]$ is real.
\end{enumerate}
\begin{Rem}
    The second and third assumptions of \ref{item:ass_insta_RealForce} could be formally deduced from the first one. As they are straightforward to check in the applications, we prefer assuming they are true rather than proving them in full generality.
\end{Rem}
\end{enumerate}

\subsection{Main abstract instability result}
Our main result reads as follows.
\begin{Thm}\label{thm:instability}
	Assume that System~\eqref{eq:system_of_study_translated} satisfies \ref{item:ass_insta_expansion}--\ref{item:ass_insta_StabHIGHREG}--\ref{item:ass_insta_LINEARIZED}--\ref{item:ass_insta_RealForce}. Let  $k_0 \in \N$, let $\boldsymbol{\xi} \in L^1_\alpha$ be a family of vector fields on $\R^d$ such that, considering $\boldsymbol{X}$ as defined in  \ref{item:ass_insta_LINEARIZED}(b) and $P_{1,n}$ defined in \ref{item:ass_insta_StabHIGHREG}(a), there holds
	\begin{equation}\label{eq:def-Lambda(X)-instab-thm}
	\Lambda_\xi(X) \vcentcolon= \int \xi(\alpha) \cdot P_{1,n}(X)(\alpha) \D \mu(\alpha) \neq 0.
	\end{equation}
	 Then there exists $\kappa>0$ such that for all $\delta >0$ small enough and $s \in \N$, there is a real solution $\ZZ$ of~\eqref{eq:system_of_study_translated} well defined up to a time $T_\delta = \mathcal O (|\log \delta|)$ such that
	\begin{equation}
	\label{eq:initial_data_small}
	\esssup_{\alpha \in I} \| Z^\alpha |_{t=0} \|_{W^{s,\infty}_x} \leq \delta
	\end{equation}
	but
	\begin{equation}
	\label{eq:nonlinear_growth}
	\left\| \int \xi(\alpha) \cdot G^\alpha[\ZZ|_{t = T_{\delta}}] \, \mathrm{d}\mu(\alpha) \right\|_{W^{-k_0,1}_x} \geq \kappa.
	\end{equation}
\end{Thm}
\begin{Rem}
    This result is a Lyapunuov instability result because it asserts that if the system~\eqref{eq:system_of_study} is linearly unstable around $0$ in the sense of Assumption~\ref{item:ass_insta_LINEARIZED}(b), then there exists a small neighbourhood of $0$, described by~\eqref{eq:nonlinear_growth}, such that there exist arbitrarily small initial data generating solutions exiting this neighbourhood. Since we are in infinite dimension, we have to choose the topologies of the neighbourhood considered and in which we measure the size of initial data. 
    
    Here, we are able to prove that the initial data can be chosen not only arbitrarily small, but even controlling an arbitrarily large number of derivatives, while the neighbourhood to exit is measured in arbitrarily negative Sobolev regularity. The quantity inside the norm in~\eqref{eq:nonlinear_growth} has to be understood as a moment associated with our solutions. Therefore, a way to formulate informally this result is that linear instability implies that there exist solutions starting as close as we want to zero, even in high Sobolev regularity, but for which one well chosen moment becomes large, even in negative Sobolev regularity. 
\end{Rem}
\subsection{Proof of Theorem \ref{thm:instability}}
\label{sec:instability}
Let us pick $s \in \N$, and consider $\boldsymbol{X}$, $n$ and $\lambda$ as defined in \ref{item:ass_insta_LINEARIZED}(b). We normalize $\boldsymbol{X}$ in such a way that calling 
	\begin{equation}
	\label{eq:def_Y1}
	Y_1^\alpha(t,x) \vcentcolon= \Re \Big( X(\alpha) \exp(i n \cdot x) \exp(\lambda t) \Big),
	\end{equation}
	we have 
	\begin{equation}\label{def-vectorY_1-instab}
		\esssup_{\alpha} \| Y_1^\alpha(0) \|_{W^{s,\infty}_x} = 1.
	\end{equation}
	Note that as a consequence of \ref{item:ass_insta_RealForce}(c), $\YY_1$ is a nontrivial growing real-valued solution of the linear equation~\eqref{eq:system_of_study_linear}.

The whole proof consists in showing that, as soon as $\delta$ is sufficiently small, $\delta \YY_1$, which grows exponentially fast, is a good approximation of the solution $\ZZ_\delta$ of the nonlinear equation~\eqref{eq:system_of_study_translated} starting from the same initial condition
	\begin{equation}
	\label{eq:initial_condition}
	Z_\delta^\alpha(0,x) = \delta Y^\alpha_1(0,x),
	\end{equation}
which satisfies~\eqref{eq:initial_data_small}. Therefore, we first need to measure the growth of $\delta \YY_1$  in the same way that we will estimate the solution itself (see~\eqref{eq:nonlinear_growth}), and then to show that $\delta \YY_1$ is indeed close to $\ZZ_\delta$. Let us start with the first point. By definition of $P_{1,n}$ in \ref{item:ass_insta_expansion}(a), we have
\begin{equation*}
    \GG_1[\YY_1] = \Re\Big(P_{1,n}(X) \exp(in\cdot x) \exp(\lambda t)\Big).
\end{equation*}
Hence, if we give ourselves $k_0$, $\boldsymbol{\xi}$ and $\theta$ as in the statement of the theorem, we have
\begin{align*}
    &\bigg\| \int \xi(\alpha ) \cdot G_1^\alpha[\delta \YY_1|_{t = T}]\D \mu(\alpha) \bigg\|_{W^{-k_0,1}_x} \\
    &\quad= \delta \exp(\theta T) \left\| \Re\left( \int \xi(\alpha ) \cdot P_{1,n}(X(\alpha)) \exp(i (n \cdot x + \Im(\lambda) T) ) \D \mu(\alpha)\right) \right\|_{W^{-k_0,1}_x} \\
	&\quad= \delta \exp(\theta T) \mathfrak{M}_{X,n,k_0},
\end{align*}
where we have set
\begin{align}\label{def:eqM(X,n,k_0)}
    \mathfrak{M}_{X,n,k_0}\vcentcolon = \left\| \Re\Big(\Lambda_\xi(X) \exp(i (n \cdot x + \Im(\lambda) T) ) \Big) \right\|_{W^{-k_0,1}_x},
\end{align}
using the definition \eqref{eq:def-Lambda(X)-instab-thm}. Note that $\mathfrak{M}_{X,n,k_0}>0$.

As a consequence,  at time $ \tau(\delta,A)$ defined as
	\begin{equation}
	\label{eq:def_final_time}
	 \tau(\delta,A) \vcentcolon= \frac{|\log \delta|}{\theta} + A, 
	 \end{equation}
for some $A \in \R$  to be chosen later, and picking $0<\delta<\min(e^{\theta A},1)$ so that $\tau(\delta, A)>0$, we have
\begin{equation}
\label{eq:growing_mode_grew}
        \bigg\| \int \xi(\alpha ) \cdot G_1^\alpha[\delta \YY_1|_{t = \tau(\delta, A)}]\D \mu(\alpha) \bigg\|_{W^{-k_0,1}_x} = 2\kappa_A,
\end{equation}
with
\begin{equation}
\label{eq:def_m}
    \kappa_A \vcentcolon= \frac{1}{2} \exp(\theta A) \frac{|\Lambda_\xi(X)|}{|n|^{k_0}} \int_{\T}|\cos y| \D y.
\end{equation}
Note that by definition of $\tau(\delta,A)$ and $\kappa_A$, we have
\begin{align}\label{albator}
    \bigg\| \int \xi(\alpha ) \cdot G_1^\alpha[\delta \YY_1(t)]&\D \mu(\alpha) \bigg\|_{W^{-k_0,1}_x}=2\kappa_A e^{\theta(t-\tau(\delta,A))}.
\end{align}
In order to slightly simplify the computations, from now on, we assume without loss of generality that $\int |\xi(\alpha)| \D \mu(\alpha) = 1$.

The second point of the proof is to estimate how close is $\delta \YY_1$ from $\ZZ_\delta$. In view of the bound from below
\begin{align}
    	\notag \bigg\| & \int  \xi(\alpha) \cdot G^\alpha[\ZZ_\delta(t)] \, \mathrm{d}\mu(\alpha) \bigg\|_{W^{-k_0,1_x}} \\ 
        \notag	&\geq \bigg\| \int \xi(\alpha ) \cdot G_1^\alpha[\delta \YY_1(t)]\D \mu(\alpha) \bigg\|_{W^{-k_0,1}_x} \\
        \notag &\qquad \qquad \qquad - \bigg\| \int \xi(\alpha ) \cdot (G_1^\alpha[\ZZ_\delta(t)] - G_1^\alpha[\delta \YY_1(t)])\D \mu(\alpha) \bigg\|_{W^{-k_0,1}_x}\\ \notag &\geq \bigg\| \int \xi(\alpha ) \cdot G_1^\alpha[\delta \YY_1(t)]\D \mu(\alpha) \bigg\|_{W^{-k_0,1}_x} \\
        \notag &\qquad \qquad \qquad- \left(\int |\xi(\alpha)| \D \mu(\alpha) \right) \esssup_{\alpha \in I}  \| G^\alpha[\ZZ_\delta(t)] - G_1^\alpha[\delta \YY_1(t)]  \|_{W^{-k_0,1}_x} \\
    	\label{eq:estim_approx_Zdelta_Y} &= \bigg\| \int \xi(\alpha ) \cdot G_1^\alpha[\delta \YY_1(t)]\D \mu(\alpha) \bigg\|_{W^{-k_0,1}} \hspace{-7pt}-  \esssup_{\alpha \in I}  \| G^\alpha[\ZZ_\delta(t)] - G_1^\alpha[\delta \YY_1(t)] \|_{W^{-k_0,1}_x},
\end{align}
it suffices to prove that for $\delta$ sufficiently small, there exists $A \in \R$ such that $\tau(\delta,A)>0$ and for times $t>0$ up to time $\tau(\delta,A)$ as defined in~\eqref{eq:def_final_time}, we have
\begin{align}\label{eq:nonlinear_instability_main_estimate}
    \esssup_{\alpha \in I}  \| G^\alpha[\ZZ_\delta(t)] - G_1^\alpha[\delta \YY_1(t)] \|_{W^{-k_0,1}_x} &\leq \frac{1}{2} \bigg\| \int \xi(\alpha ) \cdot G_1^\alpha[\delta \YY_1(t)]\D \mu(\alpha) \bigg\|_{W^{-k_0,1}_x}.
\end{align}
Indeed, at time $\tau(\delta, A)$, defining $\kappa_A$ as in~\eqref{eq:def_m}, we would conclude using~\eqref{albator}, \eqref{eq:estim_approx_Zdelta_Y} and~\eqref{eq:growing_mode_grew} that
\begin{equation*}
    \bigg\|  \int  \xi(\alpha) \cdot G^\alpha[\ZZ_\delta(\tau(\delta,A)] \, \mathrm{d}\mu(\alpha) \bigg\|_{W^{-k_0,1}_x} \geq 2\kappa_A -\kappa_A = \kappa_A>0,
\end{equation*}
as claimed in \eqref{eq:nonlinear_growth}.

\medskip

We now start the proof of the main estimate \eqref{eq:nonlinear_instability_main_estimate}. Let us first observe that \ref{item:ass_insta_StabHIGHREG} is stronger than the Assumptions \ref{ass:bound_high_regularity-rho} and \ref{ass:stability-rho} needed for the local well-posedness result from Theorem~\ref{thm-LWPmultiphase-densityvelocity} with $p \in [1,\infty)$, since $\mathcal{V} \in L^p_{\alpha}$ by assumption. So we already know that the solution $\ZZ_\delta$ of~\eqref{eq:system_of_study_translated} with initial condition~\eqref{eq:initial_condition} exists for small times. Also, by continuity, we know that given $\ell \in \N$ and $\delta <1$, the time 
	\begin{equation}
	\label{eq:def_T*}
	T_\ell^*(\delta) \vcentcolon= \sup \Big\{ t>0 \mbox{ s.t.~} \ZZ_\delta(t) \mbox{ exists up to time }t \mbox{ and } \esssup_{\alpha \in I} \| \DDD w_\delta^\alpha(t) \|_{H^{\ell}} \leq 1 \Big\}
	\end{equation}
		is positive. Now, following~\cite{Gr00}, the proof goes like this. We will show that for all $N \geq 2$, $\ZZ_\delta$ can be expanded in powers of $\delta$ as:
	\begin{equation}
	\label{eq:decomposition_Z}
	\ZZ_\delta = \delta \YY_1 + \sum_{k=2}^N \delta^k \YY_k + \delta^{N+1}\ZZ_{\delta,N},
	\end{equation}
	where for $k\geq 2$, $\YY_k(0) = \ZZ_{\delta,N}(0)=0$, and:
	\begin{itemize}
		\item for all $k \geq 2$, $\YY_k$ is independent of $\delta$, well defined for all times, and estimated thanks to \ref{item:ass_insta_expansion}(a), \ref{item:ass_insta_LINEARIZED}(b) and the explicit form of $\YY_1$;
		\item for $N$ and $\ell$ large enough, $\ZZ_{\delta,N}$ will be proved to be $\mathcal O(1)$ uniformly in $\delta$ in $H^\ell_x\times H^{\ell + 1}_x$ norm thanks to \ref{item:ass_insta_expansion}(b), \ref{item:ass_insta_StabHIGHREG} and energy estimates, at least up to time $\min(T^*_\ell(\delta), \tau(\delta,A))$ for a small (negative) value of $A$.
	\end{itemize} 
		Thanks to these estimates, we will be able to prove that for some $A<0$, $\delta$ small enough and $\ell$ large enough, there holds $\tau(\delta,A) \leq T^*_\ell(\delta)$   and~\eqref{eq:nonlinear_instability_main_estimate} holds up to time $\tau(\delta,A)$, uniformly in $\delta$.

	\bigskip

	\noindent \underline{Step 1: Definition and estimate for $(\YY_k)$}.
	
\medskip	
	
	Our first task is to write an equation for $\YY_k$, $k \geq 2$, to prove that it is well defined, and to estimate it. The equation solved by $\ZZ_\delta$ can be written formally
	\begin{equation}
		\label{eq:abstract_equation}
	\partial_t \ZZ_\delta + \begin{pmatrix} \mathcal{V} \cdot \nabla \rr_\delta + \Div (\ww_\delta) \\ (\mathcal{V} \cdot \nabla) \ww_\delta\end{pmatrix}  +  \QQ(\ZZ_\delta, \ZZ_\delta) = \begin{pmatrix} 0 \\ \GG[\ZZ_\delta] \end{pmatrix},
	\end{equation}
	where we have used the notation
	\begin{equation*}
	\begin{pmatrix} \mathcal{V} \cdot \nabla \rr + \Div (\ww) \\ (\mathcal{V} \cdot \nabla) \ww\end{pmatrix} \vcentcolon= \begin{pmatrix} \mathcal{V}(\alpha) \cdot \nabla r^\alpha + \Div (w^\alpha) \\ (\mathcal{V}(\alpha) \cdot \nabla) w^\alpha \end{pmatrix}_{\alpha \in I} ,
    \end{equation*}
    and
    \begin{equation*}
\QQ\left(\ZZ, \ZZ' \right) \vcentcolon=  \begin{pmatrix} \operatorname{div} (r'^\alpha w^\alpha)  \\ (w^\alpha\cdot \nabla) w'^\alpha \end{pmatrix}_{\alpha \in I}.
	\end{equation*}
	
	Let us fix $k \geq 2$. In order to derive an equation for $\YY_k$, let us plug~\eqref{eq:decomposition_Z} (for some $N\geq k$) into~\eqref{eq:abstract_equation} and identify all the terms of order $\delta^k$ (these terms are easily seen to be independent of~$N$). 
	The term $\QQ(\ZZ_\delta, \ZZ_\delta)$ is straightforward to expand using its bilinearity, and we find that the term of order $\delta^k$ is
	\begin{equation*}
	\sum_{i=1}^{k-1} \QQ(\YY_i, \YY_{k-i}).
	\end{equation*}
	To expand $\GG[\ZZ_\delta]$, we plug~\eqref{eq:decomposition_Z} into~\eqref{eq:decomposition_G} from \ref{item:ass_insta_expansion}. Doing so, we realize that the term of order $\delta^k$ is
	\begin{equation*}
\sum_{\substack{1\leq j \leq k \\ 1\leq i_1, \dots, i_j \leq k \\ i_1 + \dots + i_j = k}} \GG_j[\YY_{i_1}, \dots, \YY_{i_j}] = \GG_1[\YY_k] + \sum_{\substack{2\leq j \leq k \\ 1\leq i_1, \dots, i_j \leq k-1 \\ i_1 + \dots + i_j = k}} \GG_j[\YY_{i_1}, \dots, \YY_{i_j}].
	\end{equation*}

	Therefore, decomposition~\eqref{eq:decomposition_Z} leads formally to the following equation for $\YY_k$:
	\begin{equation}
	\label{eq:equation_Yk}
\left\{	\begin{aligned}
\partial_t \YY_k &+  \begin{pmatrix} \mathcal{V} \cdot \nabla \rr_k + \Div( \ww_k) \\ (\mathcal{V} \cdot \nabla) \ww_k\end{pmatrix}  + \begin{pmatrix} 0 \\ - \GG_1[\YY_k] \end{pmatrix} \\
&= \sum_{\substack{2\leq j \leq k \\ 1\leq i_1, \dots, i_j \leq k-1 \\ i_1 + \dots + i_j = k}}  \begin{pmatrix}  0 \\ \GG_j[\YY_{i_1}, \dots, \YY_{i_j}]\end{pmatrix} - \sum_{i=1}^{k-1} \QQ(\YY_i, \YY_{k-i}), \\
\YY_k(0) &= 0.
\end{aligned}
\right.
	\end{equation}
	Remark that the right-hand side.~only involves $(\YY_i)_{i<k}$. Therefore, we ask that $\YY_k$ solves \eqref{eq:equation_Yk}, which corresponds to the linearized system~\eqref{eq:system_of_study_linear} with a source term. This paves the way for an induction argument. Let us call $S_t$ the semigroup (formal at this level) associated with~\eqref{eq:system_of_study_linear}. Since $\YY_k(0) = 0$, the Duhamel formula leads for all $t\geq 0$ to
	\begin{multline}
	\label{eq:Duhamel_Yk}
	\YY_k(t)=\\
    \int_0^t S_{t-s} \Bigg( \sum_{\substack{2\leq j \leq k \\ 1\leq i_1, \dots, i_j \leq k-1 \\ i_1 + \dots + i_j = k}}\begin{pmatrix}  0 \\  \GG_j[\YY_{i_1}(s), \dots, \YY_{i_j}(s)]\end{pmatrix} - \sum_{i=1}^{k-1} \QQ(\YY_i(s), \YY_{k-i}(s)) \Bigg) \D s.
	\end{multline}
	
	In order to state the main result of this part of the proof, we need to introduce some functional spaces. For a given $n' \in \Z^d$, we define
	\begin{equation*}
 V_{n'} \vcentcolon= \Big\{ (\alpha, x)\in I \times \T^d \mapsto f(\alpha) \exp(i n' \cdot x), \, f \in L^\infty(I,\mu) \Big\},
	\end{equation*}
 endowed with the norm
 \begin{equation*}
     \| f(\alpha) \exp(i n' \cdot x) \|_{V_n} = \| f \|_{L^{\upinfty}_\alpha}.
 \end{equation*} 
 Also, for a given $k \in \N$, we call
	\begin{equation*}
		 \mathcal U_k \vcentcolon= \bigoplus_{j=-k}^k  V_{j n},
	\end{equation*}
	endowed with the norm defined for all $a_j \in  V_{jn}$ with $j = -k, \dots, k$ by $$\|a_{-k} + \dots + a_k \|_{\mathcal U_k} \vcentcolon= \|a_{-k}\|_{V_{-kn}} + \dots + \|a_k\|_{V_{kn}}.$$
	We recall that $n$ has been given from \ref{item:ass_insta_LINEARIZED}(a). The set $\mathcal U_k$ is the set of multiphasic configurations whose Fourier modes in $x$ have frequencies $-kn, -(k-1)n, \dots, kn$.
	
	\medskip

The main result of this step is the following.
\begin{Lem}
	\label{lem:estimate_Yk}
Let $k \geq 2$. The term $\YY_k(t)$ is well defined in $\mathcal U_k$ for all $t\geq 0$ by the formula~\eqref{eq:Duhamel_Yk}, and satisfies, for some $C_k>0$ only depending on $k$, $\|X\|_\infty$ and the continuity bounds of the different linear or multilinear operators at play:
\begin{equation}
\label{eq:estimYk}
\forall t \geq 0, \ \ \| \YY_k(t) \|_{\mathcal U_{k}}  \leq C_{k} \exp(k \theta t). 
\end{equation}	
Moreover, $\YY_k$ solves~\eqref{eq:equation_Yk}.
\end{Lem}
	\begin{proof}[Proof of Lemma~\ref{lem:estimate_Yk}]
		The proof proceeds by a direct induction argument. For $k=1$, the definition \eqref{eq:def_Y1} of $\YY_1$, which implies that for all $t\geq 0$ and $x\in \T^d$
		\begin{multline*}
		    Y_1^\alpha(t,x) = \\
            \left( \frac{1}{2} X(\alpha)\exp(i \Im(\lambda) t) \exp(in\cdot x) + \frac{1}{2} \bar X(\alpha) \exp(i \Im(\lambda) t) \exp(-in\cdot x) \right)\exp(\theta t),
		\end{multline*}
		reveals that $\YY_1(t) \in \mathcal U_1$ and that
		\begin{equation*}
		    \| \YY_1(t) \|_{\mathcal U_1} = \| \boldsymbol{X} \|_{L^\infty_\alpha}\exp (\theta t).
		\end{equation*}
		
		If the result is true at stage $1, \dots, k-1$ for some $k \geq 2$, let us prove that it is still true at stage~$k$. By \ref{item:ass_insta_expansion}(a) and the explicit definition of $\QQ$, we check using the induction assumption that for all $s \geq 0$, the right-hand side~of~\eqref{eq:equation_Yk}, that is 
        \begin{align*}
            \Xi_k(s)\vcentcolon= \begin{pmatrix}  0 \\ \displaystyle{\sum_{j=2}^k \sum_{\substack{ 1 \leq i_1, \dots, i_j \leq k-1 \\ i_1 + \dots i_j = k}} \GG_j[\YY_{i_1}(s), \dots, \YY_{i_j}(s)]}\end{pmatrix} - \sum_{i=1}^{k-1} \QQ(\YY_i(s), \YY_{k-i}(s))
        \end{align*}
        is well defined and belongs to $\mathcal U_k$, with
		\begin{align*}
		\left\| \Xi_k(s)\right\|_{\mathcal U_k} 
        \leq C_k \exp(k \theta s),
		\end{align*}
		for some $C_k>0$ only depending on $k$, the continuity bound of $Q$ and the $G_j$, $j = 1, \dots, k$.
        
        Therefore, by \ref{item:ass_insta_LINEARIZED}(a), the semigoup $S_t$ can be applied to it. Moreover, also by \ref{item:ass_insta_LINEARIZED}(a), we have for all $t\geq 0$, up to considering a larger $C_k$,
		\begin{align*}
		\left\| S_{t-s}\Xi_k(s)\right\|_{\mathcal U_k}
		 \leq C_k \exp(\theta_0(t-s))\exp(k \theta s).
		\end{align*}
		It is locally integrable with respect to~$s$, so~\eqref{eq:Duhamel_Yk} is a good definition of $\YY_k(t)$ in $\mathcal U_k$, and we have, up to considering a larger $C_k$, 
		\begin{equation*}
		\|\YY_k(t) \|_{\mathcal U_k} \leq C_k \exp(\theta_0 t) \int_0^t \exp((k\theta - \theta_0)s)\D s \leq C_k \exp(k\theta t).
		\end{equation*}
		Note that here, we crucially used $k\theta - \theta_0 >0$, which follows from $\theta_0 < 2 \theta$ (see \ref{item:ass_insta_LINEARIZED}).
	\end{proof}

\bigskip

\noindent \underline{Step 2: Energy estimate for $\ZZ_{\delta,N}$}.

\medskip

Now that we have defined and estimated $(\YY_k)_{k \in \N^*}$ for all times, we introduce the following remainder, following the ansatz~\eqref{eq:decomposition_Z}: for a given $N \in \N^*$ and for all times such that $\ZZ_\delta$ is well defined, we set: 
\begin{equation}\label{def:eq:ZdeltaN}
\ZZ_{\delta, N} =\begin{pmatrix} \rr_{\delta, N} \\ \ww_{\delta,N} \end{pmatrix}\vcentcolon= \frac{1}{\delta^{N+1}} \left( \ZZ_\delta - \sum_{k=1}^N \delta^k \YY_k\right).
\end{equation}
Note that $\ZZ_{\delta, N}(0)=0$ by construction. We expect this remainder to be of size of order $1$ as $\delta \to 0$, so we introduce the following notation for any $\eta \geq 0$:
\begin{equation}
\label{eq:def_TN}
T_{N, \ell}(\delta, \eta) \vcentcolon= \sup \Big\{ t >0 \mbox{ such that } \| \ZZ_{\delta, N}(t) \|_\ell \leq \eta \Big\},
\end{equation}
where we recall the definition \eqref{def:NORM-Hl-INSTABproof} for the norm $\|\cdot\|_\ell $ (for $\ell \in \N$).

The main result of this step will be the following lemma, yielding a suitable energy estimate on the remainder $\ZZ_{\delta, N}$.
\begin{Lem}
	\label{lem:energy_estimate_Grenier}
Let $ N \geq 1$ and $\ell \in \N$ be sufficiently large. There exists $C_\ell, C_{\ell,N}>0$ and $\eta_\ell>0$ such that $\G_\ell(\eta_\ell)< + \infty$ (see \ref{item:ass_insta_StabHIGHREG} for the definition of $\G_\ell$) and $\eta_\ell/2C_{N, \ell}<1/2$ such that for all $\delta \in (0, \eta_\ell/2C_{N, \ell})$, the following holds. Provided that
\begin{equation*}
0 \leq t \leq \min\left( T^*_\ell(\delta) , \tau\left( \delta, \frac{1}{\theta}\log\left(\frac{\eta_\ell}{2C_{N, \ell}} \right)\right), T_{N,\ell}\left( \delta, \frac{\eta_\ell}{2 \delta^{N+1}} \right) \right),
\end{equation*}
(recall the definitions \eqref{eq:def_final_time} and \eqref{eq:def_T*} for $\tau(\delta, A)$ and $T^*_\ell(\delta)$ respectively), we have
\begin{equation*}
\frac{\D}{\D t} \| \ZZ_{\delta, N}(t) \|_\ell \leq C_\ell \| \ZZ_{\delta, N}(t) \|_\ell + C_{N,\ell} \exp((N+1)\theta t).
\end{equation*}
\end{Lem}

\begin{proof}[Proof of Lemma~\ref{lem:energy_estimate_Grenier}]
In order to derive an energy estimate for $\ZZ_{\delta, N}$, let us first write an equation that it satisfies. Using the equation~\eqref{eq:abstract_equation} solved by $\ZZ_\delta$ and the one~\eqref{eq:equation_Yk} solved by $\YY_k$ for all $k$, we find in view of the definition \eqref{def:eq:ZdeltaN} of $\ZZ_{\delta, N}$:
\begin{multline*}
\partial_t \ZZ_{\delta, N} + \begin{pmatrix} \mathcal{V} \cdot \nabla \rr_{\delta, N} + \Div (\ww_{\delta,N}) \\ (\mathcal{V} \cdot \nabla) \ww_{\delta,N}\end{pmatrix} \\
+ \frac{1}{\delta^{N+1}}\left( \QQ(\ZZ_\delta,\ZZ_\delta) - \sum_{k=1}^N \delta^k \sum_{i=1}^{k-1} \QQ(\YY_i, \YY_{k-i}) \right)\\
= \frac{1}{\delta^{N+1}} \begin{pmatrix} 0 \\ \GG[\ZZ_\delta] - \displaystyle{\sum_{k=1}^N \delta^k \sum_{\substack{1\leq j \leq k \\ 1\leq i_1, \dots, i_j \leq k-1 \\ i_1 + \dots + i_j = k}}} \GG_j[\YY_{i_1}, \dots, \YY_{i_j}] \end{pmatrix}.
\end{multline*}
At this stage, it is convenient to introduce the notation:
\begin{equation*}
\YY_{\delta, N} =\begin{pmatrix} \rr'_{\delta, N} \\ \ww'_{\delta,N} \end{pmatrix}\vcentcolon= \sum_{k=1}^N \delta^k \YY_k,
\end{equation*}
so that $\ZZ_\delta=\YY_{\delta,N}+\delta^{N+1} \ZZ_{\delta, N}$ by definition. With this notation at hand, we can simplify by bilinearity on the one hand:
\begin{align*}
&\frac{1}{\delta^{N+1}}\left( \QQ(\ZZ_\delta,\ZZ_\delta) - \sum_{k=1}^N \delta^k \sum_{i=1}^{k-1} \QQ(\YY_i, \YY_{k-i}) \right) \\
&=\frac{1}{\delta^{N+1}}\Bigg(\QQ(\YY_{\delta, N}, \YY_{\delta, N})+ \delta^{N+1}\QQ(\YY_{\delta,N}, \ZZ_{\delta,N})\\
&\qquad+\delta^{N+1}\QQ(\ZZ_{\delta,N}, \YY_{\delta, N})+\delta^{2(N+1)}\QQ(\ZZ_{\delta,N}, \ZZ_{\delta, N})  
-\sum_{k=1}^N \delta^k \sum_{i=1}^{k-1} \QQ(\YY_i, \YY_{k-i}) \Bigg) \\
&= \QQ(\ZZ_\delta, \ZZ_{\delta,N}) + \QQ(\ZZ_{\delta,N}, \YY_{\delta, N}) + \sum_{\substack{1\leq i_1,i_2 \leq N \\ i_1 + i_2 \geq N+1}} \delta^{i_1 + i_2 - (N+1)}  \QQ(\YY_{i_1}, \YY_{i_2}),
\end{align*}
and on the other hand by multilinearity and the expansion of the force field (see \ref{item:ass_insta_expansion}):
\begin{align*}
 &\frac{1}{\delta^{N+1}} \left( \GG[\ZZ_\delta] - \sum_{k=1}^N \delta^k \sum_{\substack{1\leq j \leq k \\ 1\leq i_1, \dots, i_j \leq k-1 \\ i_1 + \dots + i_j = k}} \GG_j[\YY_{i_1}, \dots, \YY_{i_j}]\right) \\
 &= \frac{\GG[\ZZ_\delta] - \GG[\YY_{\delta, N}]}{\delta^{N+1}} \\
 & \quad + \sum_{\substack{1 \leq j \leq N \\ 1 \leq i_1, \dots, i_j \leq N \\ i_1 + \dots + i_j \geq N+1 }} \delta^{i_1 + \dots + i_j - (N+1)} \GG_j[\YY_{i_1}, \dots, \YY_{i_j}] + \frac{\boldsymbol R_N[\YY_{\delta,N}]}{\delta^{N+1}}.
\end{align*}
Gathering everything and using the explicit form of $\QQ$, we find:
\begin{multline*}
\partial_t \ZZ_{\delta,N} + ((\mathcal{V} + \ww_\delta) \cdot \nabla) \ZZ_{\delta,N} + \begin{pmatrix} (1+ \rr'_{\delta,N}) \Div \ww_{\delta, N} +\rr_{\delta,N} \mathrm{div}(\ww_\delta)\\ 0 \end{pmatrix} \\
+ (\ww_{\delta, N} \cdot \nabla) \YY_{\delta,N}\\ = \frac{\GG[\ZZ_\delta] - \GG[\YY_{\delta, N}]}{\delta^{N+1}} + \FF_{\delta, N},
\end{multline*}
where
\begin{multline}\label{def:F-deltaN-INSTAB}
\FF_{\delta,N} \vcentcolon= \begin{pmatrix} 0 \\ \displaystyle{\sum_{\substack{1 \leq j \leq N \\ 1 \leq i_1, \dots, i_j \leq N \\ i_1 + \dots + i_j \geq N+1 }} \delta^{i_1 + \dots + i_j - (N+1)} \GG_j[\YY_{i_1}, \dots, \YY_{i_j}] + \frac{\boldsymbol R_N[\YY_{\delta,N}]}{\delta^{N+1}} }\end{pmatrix} \\ -  \sum_{\substack{1\leq i_1,i_2 \leq N \\ i_1 + i_2 \geq N+1}} \delta^{i_1 + i_2 - (N+1)}  \QQ(\YY_{i_1}, \YY_{i_2}).
\end{multline}
At this stage, by standard Sobolev energy estimates (see Proposition \ref{prop-Sobolev}) and from the definition \eqref{def:NORM-Hl-INSTABproof} of the $\Vert \cdot \Vert_\ell$ norm (controlling $\ell$ derivatives of the first component and $\ell+1$ derivatives of the second component),
and as long as all the objects exist, there exist $C_\ell>0$ such that:
\begin{align*}
&\frac{\D}{\D t} \| \ZZ_{\delta, N} \|_\ell  
\\
&\leq \Bigg( C_\ell\Big\{ 1 + \esssup_{\alpha \in I}  \| \DDD w_\delta^\alpha \|_\ell + \| \YY_{\delta,N} \|_{\ell+1}\Big\}  \G_\ell\big( \max( \|\ZZ_\delta\|_\ell, \|\YY_{\delta,N}\|_\ell ) \big)\Bigg)
 \| \ZZ_{\delta,N}\|_\ell  \\
&\quad +  \| \FF_{\delta,N} \|_\ell.
\end{align*}
Here we have used the stability of the force field at high regularity from \ref{item:ass_insta_StabHIGHREG}.
By definition~\eqref{eq:def_T*}, we also know that if $t \leq T^*_\ell(\delta)$, then 
$$\esssup_{\alpha \in I} \| \DDD w^\alpha_\delta \|_{H^{\ell}} \leq 1.$$ Furthermore, according to the bound $\Vert Y_k \Vert_{\ell+1} \lesssim c_k(\ell) \Vert Y_k \Vert_{\mathcal{U}_k}$ for some $c_k(\ell)>0$, we set $$c_{N,\ell}\vcentcolon=\max_{1 \leq k \leq N} c_k(\ell),$$ 
and choose $A<0$ and $\delta>0$ such that
$$e^{\theta A}< 1/2, \ \ 0<\delta<e^{\theta A},$$
so that $\tau(\delta, A)>0$. We infer from Lemma~\ref{lem:estimate_Yk} that we have for all $0 \leq  t \leq \tau(\delta, A)$,
\begin{align}\label{eq:AVANTestimate_YdeltaN-unif}
    \|\YY_{\delta, N}(t)\|_{\ell+1}
\leq c_{N, \ell} \sum_{k=1}^N  C_k (\delta e^{\theta t})^k \leq c_{N,\ell}\left( \max_{1 \leq k \leq N} C_k \right) \frac{\delta e^{\theta t}}{1-\delta e^{\theta t}},
\end{align}
and therefore
\begin{align*}
    \|\YY_{\delta, N}(t)\|_{\ell+1}
\leq c_{N,\ell}\left( \max_{1 \leq k \leq N} C_k \right) \frac{\delta e^{\theta \tau(\delta,A)}}{1-\delta e^{\theta \tau(\delta,A)}}
=c_{N,\ell}\left( \max_{1 \leq k \leq N} C_k \right)\frac{ e^{\theta A}}{1- e^{\theta A}},
\end{align*}
by the definition \eqref{eq:def_final_time} of $\tau(\delta,A)$. Setting $C_{N, \ell}\vcentcolon= 2c_{N,\ell}  (\max_{1 \leq k \leq N} C_k)$, and using $e^{\theta A}< 1/2$, we infer
\begin{align}\label{eq:estimate_YdeltaN-unif}
 \|\YY_{\delta, N}(t)\|_{\ell+1} 
 \leq C_{N,\ell} e^{\theta A}.
\end{align}
Note that the former former constraint imposes $A<0$, and that without loss of generality, we can always assume that $ C_{N,\ell}>1$.
Also, we have $\| \ZZ_\delta \|_\ell \leq \|\YY_{\delta, N}\|_\ell + \delta^{N+1} \| \ZZ_{\delta, N} \|_\ell$ therefore
\begin{align*} 
\max( \|\ZZ_\delta\|_\ell, \|\YY_{\delta,N}\|_\ell ) \leq  C_{N,\ell}e^{\theta A}+ \delta^{N+1} \| \ZZ_{\delta, N}\|_\ell.
\end{align*}
All in all, if $0 \leq t \leq \min(T^*_\ell(\delta), \tau(\delta, A))$ with the former constraint on $\delta$ and $A$, we find
\begin{align*}
\frac{\D}{\D t} \| \ZZ_{\delta, N} \|_\ell &\leq \left( C_\ell\Big\{ 2 + C_{N,\ell}e^{\theta A}\Big\} + \G_\ell\Big( C_{N,\ell} e^{\theta A} + \delta^{N+1} \|\ZZ_{\delta, N}\|_\ell \Big)\right) \| \ZZ_{\delta,N}\|_\ell \\
& \quad + \| \FF_{\delta,N} \|_\ell.
\end{align*}
Let $0<\eta<1$ be fixed, so that we have $\eta/2C_{N,\ell}<1/2$ since  $C_{N,\ell}>1$.  Choosing $A<0$ such that $e^{\theta A}= \eta / ( 2 C_{N,\ell})<1/2$, we find that whenever
\begin{align*}
t \leq \min \left(T^*_\ell(\delta), \tau\left(\delta, \frac{1}{\theta}\log\left(\frac{\eta_\ell}{2C_{N, \ell}} \right)\right), T_{N, \ell}\left(\delta, \frac{\eta_\ell}{2\delta^{N+1}}\right)\right), \ \ \delta<e^{\delta A}=\frac{\eta_\ell}{2C_{N, \ell}},
\end{align*}
and up to choosing a larger $C_\ell >0$, there holds in view of the definition \eqref{eq:def_TN}
\begin{equation*}
\frac{\D}{\D t} \| \ZZ_{\delta, N} \|_\ell \leq C_\ell \| \ZZ_{\delta,N}\|_\ell + \| \FF_{\delta,N} \|_\ell.
\end{equation*}
It remains to estimate the contribution of the last term $\FF_{\delta,N}$, defined earlier in \eqref{def:F-deltaN-INSTAB}: we simply write
\begin{align*}
 \| \FF_{\delta,N} \|_\ell &\leq \sum_{\substack{1 \leq j \leq N \\ 1 \leq i_1, \dots, i_j \leq N \\ i_1 + \dots + i_j \geq N+1 }} \delta^{i_1 + \dots + i_j - (N+1)} \Vert \GG_j[\YY_{i_1}, \dots, \YY_{i_j}] \Vert_{H^{\ell+1}_x}   \\
 & \quad + \frac{\Vert \boldsymbol R_N[\YY_{\delta,N}] \Vert_{H^{\ell+1}_x}}{\delta^{N+1}} + \sum_{\substack{1\leq i_1,i_2 \leq N \\ i_1 + i_2 \geq N+1}} \delta^{i_1 + i_2 - (N+1)}  \Vert \QQ(\YY_{i_1}, \YY_{i_2})\|_\ell.
\end{align*}
Relying on \eqref{eq:AVANTestimate_YdeltaN-unif}, let us also notice that for $0 \leq t \leq \tau(\delta, A)$, we have
\begin{align}\label{eq:estimate_YdeltaN}
\| \YY_{\delta,N}(t) \|_{\ell} \leq 2C_{N,\ell} \delta \exp(\theta t).
\end{align}
Hence, appealing to Lemma~\ref{lem:estimate_Yk} and \ref{item:ass_insta_expansion}, we infer for $\ell$ sufficiently large and $0 \leq t\leq \tau(\delta, \theta^{-1}\log(\eta_\ell/2C_{N, \ell}))$ with $\delta<\eta_\ell/2C_{N, \ell}<1/2$ that, up to choosing a larger $C_{N,\ell}$, there holds
\begin{equation*}
\| \FF_{\delta, N}(t) \|_\ell \leq C_{N,\ell} \exp((N+1)\theta t).
\end{equation*}
The desired result then follows.
\end{proof}

\noindent \underline{Step 3: Choice of $\ell$, $N$ and $\kappa$}.

\medskip

Frow now on, we choose $\ell_0$ for which Lemma~\ref{lem:energy_estimate_Grenier} holds, $\eta_0 \vcentcolon= \eta_{\ell_0}$ and $N_0$ sufficiently large to have $(N_0+1)\theta - C_{\ell_0} >1$. We also call $C_0 \vcentcolon= C_{N_0,\ell_0}$. These choices are now fixed for the rest of the proof. With these definitions, applying the Gronwall lemma to the conclusion of Lemma~\ref{lem:energy_estimate_Grenier}, we find that whenever
\begin{equation}
\label{eq:first_condition_t}
t \leq \min\left( T^*_{\ell_0}(\delta) , \tau\left( \delta, \frac{1}{\theta}\log\left(\frac{\eta_0}{2 C_0} \right)\right), T_{N_0,\ell_0}\left( \delta, \frac{\eta_0}{2 \delta^{N_0+1}} \right) \right), \ \  0<\delta<\frac{\eta_0}{2C_0},
\end{equation}
we have
\begin{align*}
   \| \ZZ_{\delta,N_0}(t) \|_{\ell_0} &\leq C_0 e^{C_{\ell_0}t} \int_{0}^t e^{((N_0+1)\theta-C_{\ell_0} )s}  \, \mathrm{d}s  \\
   &\leq  \frac{C_0}{(N_0+1)\theta-C_{\ell_0}} \exp((N_0+1)\theta t),
\end{align*}
since $(N_0+1)\theta - C_{\ell_0} >1$, hence
\begin{equation}
\label{eq:estimate_ZdeltaN}
\| \ZZ_{\delta,N_0}(t) \|_{\ell_0} \leq C_0 \exp((N_0+1)\theta t).
\end{equation}

Now, our goal is to simplify~\eqref{eq:first_condition_t} by showing that~\eqref{eq:estimate_ZdeltaN} holds for all $t \leq \tau(\delta,\theta^{-1}\log(\kappa_0))$ for a well-chosen $\kappa_0$ (independent of $\delta$). First, remark that whenever $\kappa \leq \kappa_1\vcentcolon= \eta_0 / (2 C_0)$ and 
\begin{equation*}
t \leq \min\left( T^*_{\ell_0}(\delta) , \tau\left( \delta, \frac{1}{\theta}\log(\kappa) \right), T_{N_0,\ell_0}\left( \delta, \frac{\eta_0}{2 \delta^{N_0+1}} \right) \right),
\end{equation*}
we have by definition~\eqref{eq:def_final_time} of $\tau(\delta,\theta^{-1}\log(\kappa))$, 
\begin{align*}
\| \ZZ_{\delta,N_0}(t) \|_{\ell_0}&\leq C_0 \exp\big((N_0+1)\theta t\big) \\
&\leq  C_0 \exp\big((N_0+1)\theta \tau(\delta,\theta^{-1}\log(\kappa))\big) \leq C_0 \frac{\kappa^{N_0 + 1}}{\delta^{N_0 + 1}}.
\end{align*}
Therefore, if in addition we take $\kappa \leq \kappa_2\vcentcolon= (\eta_0/2C_0)^{1/(N_0+1)}$, then we have
\begin{equation*}
\| \ZZ_{\delta,N_0}(t) \|_{\ell_0}<\frac{\eta_0}{2 \delta^{N_0+1}},
\end{equation*}
and thus $t < T_{N_0,\ell_0}(\delta, \eta_0/(2\delta^{N_0+1}))$ by definition \eqref{eq:def_TN}. We conclude that whenever $\kappa \leq \min(\kappa_1,\kappa_2)$, 
\begin{equation*}
T_{N_0,\ell_0}(\delta, \eta_0/(2\delta)^{N_0 + 1}) \geq \min\left(T^*_{\ell_0}(\delta), \tau\left(\delta,\frac{1}{\theta}\log(\kappa)\right)\right),
\end{equation*}
and consequently that~\eqref{eq:estimate_ZdeltaN} holds for all $t \leq \min(T^*_{\ell_0}(\delta), \tau(\delta,\theta^{-1}\log(\kappa)))$.

With the same kind of argument, let us show that for $\kappa$ sufficiently small, $\tau(\delta,\theta^{-1}\log(\kappa)) \leq T^*_{\ell_0}(\delta)$, so that~\eqref{eq:estimate_ZdeltaN} will then hold for all $t \leq \tau(\delta,\theta^{-1}\log(\kappa))$. So let us assume that $\kappa\leq \min(\kappa_1,\kappa_2)$, and let us consider $$t \leq \min( T^*_{\ell_0}(\delta),\tau(\delta,\theta^{-1}\log(\kappa))).$$ 
By definition, we have
\begin{align*}
\esssup_{\alpha \in I}  \| \DDD w^\alpha_\delta(t) \|_{H^{\ell_0}_x} \leq \| \ZZ_{\delta}(t) \|_{\ell_0} &\leq \| \YY_{\delta,N_0}(t) \|_{\ell_0} + \delta^{N_0+1}\| \ZZ_{\delta,N_0}(t) \|_{\ell_0} \\
&\leq C_0 \Big( \delta \exp(\theta t) + \delta^{N_0+1} \exp\big( (N_0+1) \theta t \big)\Big),
\end{align*}
where we used~\eqref{eq:estimate_YdeltaN} and~\eqref{eq:estimate_ZdeltaN} in the last inequality. Therefore, if in addition $t < \tau(\delta,\theta^{-1}\log(\kappa))$, then
\begin{equation*}
\esssup_{\alpha \in I} \| \DDD w^\alpha_\delta \|_{H^{\ell + 1}_x}< C_0 (\kappa + \kappa^{N_0+1}).
\end{equation*}
Hence, if we consider $\kappa \leq \kappa_3$ with $\kappa_3$ chosen small enough so that $C_0(\kappa_3 + \kappa_3^{N_0+1}) \leq \frac{1}{2}$, there holds $t < T^*_{\ell_0}(\delta)$ in view of the definition \eqref{eq:def_T*} of $T^*_{\ell_0}(\delta)$. We conclude that choosing $\kappa_0 = \min(\kappa_1,\kappa_2,\kappa_3)$, estimate~\eqref{eq:estimate_ZdeltaN} holds for all $t \leq \tau(\delta,\theta^{-1}\log(\kappa_0))$.

\bigskip

\noindent \underline{Step 4: Obtaining the final estimate \eqref{eq:nonlinear_instability_main_estimate}.}

\medskip

By the stability estimate from \ref{item:ass_insta_StabHIGHREG}, we know that for $\ell$ large enough
\begin{align*}
\esssup_{\alpha \in I} \| G^\alpha[\ZZ_{\delta}(t)] &- G^\alpha[\delta \YY_1(t)] \|_{W^{-k_0,1}_x} \\
& \leq \G_\ell\Big(\max(\| \ZZ_{\delta}(t) \|_\ell, \| \delta \YY_1(t) \|_\ell)\Big) \| \ZZ_{\delta}(t) - \delta \YY_1(t) \|_\ell. \\
& \lesssim \G_\ell\Big(\| \ZZ_{\delta}(t)-\delta \YY_1(t) \|_\ell+ \| \delta \YY_1(t) \|_\ell\Big) \| \ZZ_{\delta}(t) - \delta \YY_1(t) \|_\ell.
\end{align*}
 By definition, and thanks to the previous estimates \eqref{eq:estimYk} and \eqref{eq:estimate_ZdeltaN}, we first have for all $t \in [0, \tau(\delta,\theta^{-1}\log(\kappa))]$
\begin{align*}
    \| \ZZ_{\delta}(t) - \delta \YY_1(t) \|_\ell &\leq  \sum_{k=2}^N\delta^k \|\YY_k(t) \|_\ell+ \delta^{N+1}\Vert \ZZ_{\delta,N} (t) \Vert_{\ell}  \\
     & \leq C \frac{(\delta e^{\theta t})^2}{1-\delta e^{\theta t}}+ C\delta^{N+1} e^{(N+1)\theta t},
\end{align*}
so that for all $t \in [0, \tau(\delta,\theta^{-1}\log(\kappa))]$
\begin{align*}
\| \ZZ_{\delta}(t) - \delta \YY_1(t) \|_\ell &\leq C \delta e^{\theta t} \left( \frac{\delta e^{\theta t}}{1-\delta e^{\theta t}}+ \delta^{N} e^{N\theta t} \right) \leq  C \delta e^{\theta t} \left( \frac{2\kappa}{1-2\kappa}+(2\kappa)^N \right).
\end{align*}
Reducing $\kappa$ if necessary, we can now ensure
\begin{align*}
\| \ZZ_{\delta}(t) - \delta \YY_1(t) \|_\ell &\leq \frac{\mathfrak{M}_{X,n,k_0}}{10}\delta e^{\theta t},
\end{align*}
where we recall the definition \eqref{def:eqM(X,n,k_0)} of $\mathfrak{M}_{X,n,k_0}$. Coming back to the estimate given before, we thus infer that 
\begin{align*}
\esssup_{\alpha \in I}  \| G^\alpha[\ZZ_{\delta}(t)] - G^\alpha[\delta \YY_1(t)] \|_{W^{-k_0,1}_x} 
&\lesssim \G_\ell\Big(\frac{\mathfrak{M}_{X,n,k_0}}{10}\delta e^{\theta t}+ \delta e^{\theta t}\Big)  \frac{\mathfrak{M}_{X,n,k_0}}{10}\delta e^{\theta t} \\ &\lesssim \G_\ell\Big(3\kappa\Big) \frac{\mathfrak{M}_{X,n,k_0}}{2}\delta e^{\theta t},
\end{align*}
that is, reducing again $\kappa$ if necessary in the previous step, we end up for all $t \in [0, \tau(\delta,\theta^{-1}\log(\kappa))]$ with 
\begin{align*}
\esssup_{\alpha \in I}  \| G^\alpha[\ZZ_{\delta}(t)] - G^\alpha[\delta \YY_1(t)] \|_{W^{-k_0,1}_x} \leq \frac{\mathfrak{M}_{X,n,k_0}}{10}\delta e^{\theta t}.
\end{align*}
We finally have everything at hand to obtain the final estimate \eqref{eq:nonlinear_instability_main_estimate}. As a matter of fact, in view of \eqref{albator}, it is sufficient to check that for $A=\theta^{-1} \log(\kappa)$ (so that $\kappa_A$ defined in \eqref{eq:def_m} is actually $\kappa_A=\frac{\kappa}{2}\mathfrak{M}_{X,n,k_0}$) there holds
\begin{align*}
\frac{\mathfrak{M}_{X,n,k_0}}{10}\delta e^{\theta t} < \kappa_A e^{\theta(t-\tau(\delta, A))},
\end{align*}
which is precisely equivalent to 
\begin{align*}
\frac{\mathfrak{M}_{X,n,k_0}}{10}<\frac{\kappa_A}{\kappa}=\frac{\mathfrak{M}_{X,n,k_0}}{2}, 
\end{align*}
and the latter is therefore satisfied. This eventually concludes the proof.

\section{Application to nonlinear instability for Vlasov--Poisson type equations}
\label{sec:insta-appli}

We conclude this chapter with the application of the abstract nonlinear instability result from Section \ref{Section-FrameworkInstab} to the study of instabilities in the Vlasov--Poisson models for electrons and ions, both set on the torus $\T^d$. We recall that they respectively read as
\begin{equation*}
\mathrm{(VP_e)}\left\{ 
\begin{gathered}
\partial_t f + v \cdot \nabla_x f - \nabla_x U\cdot \nabla_v f = 0,\\
		-\Delta_x U = \int f\D v- \iint f\D v \D x,
\end{gathered}
\right. 
\end{equation*}
and
\begin{equation*}
\mathrm{(VP_i)}\left\{ 
\begin{gathered}
\partial_t f + v \cdot \nabla_x f - \nabla_x U\cdot \nabla_v f = 0,\\
		e^{U}-\Delta_x U = \int f\D v.
\end{gathered}
\right.
\end{equation*}
We start by introducing the Penrose instability condition \cite{Pen60} which characterizes the spectral instability of a profile $G$. This corresponds to the cancellation of the so-called dielectric function, which determines the dispersion relation of the plasma \cite{vlasov1938vibrational,Landau}. 

\begin{Def}
Let $G$ be a Borel measure. We say that $G$ satisfies the Penrose instability condition if
\begin{itemize}
    \item  For $\mathrm{(VP_e)}$:
there exist $\lambda \in \C$ with $\Re \lambda>0$ and $n \in \Z^d {\setminus \{0\}}$ such that
	\begin{equation}\label{Penrose-kin-elec}
	\int \frac{1}{(\lambda + in \cdot v)^2}\D G(v) = -1.
	\end{equation}
	 \item  For $\mathrm{(VP_i)}$:
there exist $\lambda \in \C$ with $\Re \lambda>0$ and $n \in \Z^d {\setminus \{0\}}$ such that
	\begin{equation}\label{Penrose-kin-ions}
	\frac{|n|^2}{1+|n|^2}\int \frac{1}{(\lambda + in \cdot v)^2}\D G(v) = -1.
	\end{equation}
	\end{itemize}

\end{Def}

The main instability theorem, which pertains to both systems,  is the following.
\begin{Thm}
\label{thm-insta-kin-VP}
Let $G \in \mathcal{P}_p(\R^d)$ (with $p \in [1, \infty)$) be a Borel measure satisfying the Penrose instability condition~\eqref{Penrose-kin-elec} (resp. \eqref{Penrose-kin-ions}).
Let $k_0\in \N$. There exists $M>0$ such that for all $s\in \N$ and for all $\delta>0$, there exist a nonnegative measure-valued initial condition $f_{0,\delta}$, a positive time $T_\delta = O(|\log \delta|)$,  and an associated weak solution to the Vlasov--Poisson system for electrons  $\mathrm{(VP_e)}$ (resp. for ions  $\mathrm{(VP_i)}$)  on $[0,T_\delta]$ denoted by $f_\delta(t)$ such that the following holds. For all  $\varphi \in \mathscr{C}^\infty(\R^d)$ with $\|\varphi\|_{W^{s,\infty}_x} \leq 1$, we have
\begin{equation}
\label{eq:assumptioninsta-VP}
\left\| \langle f_{0,\delta} - G, (1+\vert v \vert^p) \varphi\rangle \right\|_{W^{s,\infty}_x} \leq \delta,
\end{equation} 
 but, if $U[f_\delta] $ stands for the electric potential created by $f_\delta$, then
\begin{equation}
\label{eq:lyap_insta}
\sup_{[0,T_\delta]} \| \nabla U[f_\delta] (t) \|_{W^{-k_0,1}_x}  \geq M.
\end{equation} 

\end{Thm}

\begin{Rem}
    
A few comments are in order.
\begin{itemize}
\item On the first hand this can be viewed as a generalization of the instability result of \cite{HKH} and \cite{HKN} for $\mathrm{(VP_e)}$, about non-smooth equilibria $G$ (in these works, $G$ has to be $\mathscr{C}^k$ with $k$ large);  also note that these works impose a technical condition on $G$ (the so-called $\delta$ (or $\delta'$) condition) to ensure that the initial condition $f_{0,\eps}$ is nonnegative. Theorem~\ref{thm-insta-kin-VP} thus shows that this is not necessary.

\item This result is also a generalization of the ``almost Lyapunov instability'' result of \cite{Bar} (still pertaining to $\mathrm{(VP_e)}$), which also relied on the multiphasic approach but did not reach the complete instability statement, in the sense that the $M$ in the right-hand side of~\eqref{eq:lyap_insta} had to be replaced by $\delta^\alpha$ for any value of $\alpha \in (0,1]$.

\item  In the case of the Vlasov--Poisson system for ions $\mathrm{(VP_i)}$, a result analogous to that of \cite{HKH} was recently obtained in \cite{GPI2}. Theorem~\ref{thm-insta-kin-VP} thus also generalizes this result.

\item Finally, we recover the instability result for two-fluid pressureless Euler--Poisson equation~\eqref{eq:Euler--Poisson-2} of \cite{CGG}, we will come back with more detail to this point.

\end{itemize}

\end{Rem}

\begin{Rem}
\label{rem:perturb-insta} 
   The proof of Theorem~\ref{thm-insta-kin-VP} is based on the multiphasic framework and will follow from an application of the general Theorem~\ref{thm:instability}, whence the following remarks.
\begin{itemize}
    \item In the aforementioned two-fluid case from Item 4 of last Remark, it will readily follow that the unstable solution remains a two-fluid solution as well, that is to say, the nonlinear instability is  not a {\it deus ex machina} kinetic effect.
    \item In the ``genuinely kinetic case'', when the equilibrium $G$ is an integrable function, one can ensure that the unstable solution also remains an integrable function. We send the reader to the upcoming Remark~\ref{rem:perturb-insta2} for some explanations. 
\end{itemize}

\end{Rem}

We can also treat the Vlasov--Poisson equation with both dynamical electrons and ions, which we studied in Subsection~\ref{sec:severalspecies}. We recall the equation reads as
\begin{equation*}
		\mathrm{(VP_{i,e})}
		\left\{
		\begin{gathered}
     \partial_t f_{i} + v \cdot \nabla_x f_{i}- \nabla_x U\cdot \nabla_v f_{i} = 0,\\
		\partial_t f_e + v \cdot \nabla_x f_e + \nabla_x U\cdot \nabla_v f_e = 0,\\
		-\Delta_x U =   \int f_{i}\D v -  \int f_{e}\D v.
        \end{gathered}
		\right.
		\end{equation*}
We note that spatially homogeneous distribution functions $(G_i(v), G_e(v))$ such that
$$
\int G_i(v) \D v = \int G_e(v) \D v
$$
are equilibria for this system.

\begin{Def}
Let $G_i, G_e$ be two Borel measures with equal mass. We say that $(G_i, G_e)$ satisfies the Penrose instability condition for $\mathrm{(VP_{i,e})}$ if 
there exist $\lambda \in \C$ with $\Re \lambda>0$ and $n \in \Z^d {\setminus \{0\}}$ such that
	\begin{equation}\label{Penrose-kin-several}
	\int \frac{1}{(\lambda + in \cdot v)^2}\D G_i(v)  + \int \frac{1}{(\lambda + in \cdot v)^2}\D G_e(v)   = -1.
	\end{equation}

\end{Def}

The nonlinear instability result for $\mathrm{(VP_{i,e})}$  is the following.
\begin{Thm}
\label{thm-insta-kin-VP-several}
Let $G_i, G_e \in \mathcal{P}_p(\R^d)$ (with $p \in [1, \infty)$) be two Borel measures with $(G_i, G_e)$ satisfying the Penrose instability condition~\eqref{Penrose-kin-several}.
Let $k_0\in \N$. There exists $M>0$ such that for all $s\in \N$ and all $\delta>0$, there exist two nonnegative measure-valued initial condition $f_{i,0,\delta}, f_{e,0,\delta}$,  a positive time $T_\delta = O(|\log \delta|)$ and an associated weak solution to the Vlasov--Poisson system  $\mathrm{(VP_{i,e})}$ in $[0,T_\delta]$ denoted by $(f_{i,\delta},f_{e,\delta})$ such that the following holds. For all  $\varphi \in \mathscr{C}^\infty(\R^d)$ with $\|\varphi\|_{W^{s,\infty}_x} \leq 1$, we have 
\begin{equation}
\label{eq:assumptioninsta-VP-several}
\left\| \langle f_{i,0,\delta} - G_i, (1+\vert v \vert^p) \varphi\rangle \right\|_{W^{s,\infty}_x} + \left\| \langle f_{e,0,\delta} - G_e, (1+\vert v \vert^p) \varphi\rangle \right\|_{W^{s,\infty}_x}  \leq \delta,
\end{equation} 
but, if $U[(f_{i,\delta},f_{e,\delta})] $ stands for the electric potential created by $(f_{i,\delta},f_{e,\delta})$, then
\begin{equation}
\label{eq:lyap_insta-several}
\sup_{[0,T_\delta]} \| \nabla U[(f_{i,\delta},f_{e,\delta})] (t) \|_{W^{-k_0,1}_x}  \geq M.
\end{equation} 
\end{Thm}

\begin{Rem}\label{RMQ-instab-GRAVITY}(On gravitational models) 
In principle, we could also include the treatment of systems for gravitation in our instability results. We may effortlessly consider the following gravitational Vlasov--Poisson equation (set on $\T^d$):
\begin{equation}\label{eq:VP-gravitation-INSTAB}
\left\{ 
\begin{gathered}
\partial_t f + v \cdot \nabla_x f - \nabla_x U\cdot \nabla_v f = 0,\\
		\Delta_x U = \int f\D v - \iint f\D v \D x.
\end{gathered}
\right.
\end{equation}
Indeed, that proofs performed for $\mathrm{(VP_e)}$  can be directly adapted to the gravitational case. As a matter of fact, the only difference is the flip of sign in the Poisson equation, hence modifying the associated Penrose condition by a minus sign. We therefore obtain a version of Theorem \ref{thm-insta-kin-VP} for \eqref{eq:VP-gravitation-INSTAB}.

As mentioned in Section \ref{Section-PresentationVPoissonandco}, the Vlasov--Poisson system~\eqref{eq:VP-gravitation-INSTAB} can be viewed as a kind of formal linearization of the Vlasov--Monge--Ampère system
\begin{equation}\label{eq:VMA-INSTAB}
\left\{ 
\begin{gathered}
\partial_t f + v \cdot \nabla_x f - \nabla_x U\cdot \nabla_v f = 0,\\
		\mathrm{det}(\mathrm{Id}+\mathrm{D}^2_x U) = \int f\D v.
\end{gathered}
\right.
\end{equation}
Our abstract result from Theorem \ref{thm:instability} seems to be well-adapted to handle instability issue for such a system, since we developed a general framework allowing for force fields that are genuinely nonlinear in density/velocity.  A result similar to Theorem \ref{thm-insta-kin-VP} should hold for the Vlasov--Monge--Ampère system, and 
the required assumptions of Theorem \ref{thm:instability} are likely to be satisfied for the multiphasic version of \eqref{eq:VMA-INSTAB}, albeit at the cost of much more technical proofs. For the sake of brevity, we prefer not to delve into such refinements. 
\end{Rem}

\subsection{Examples of homogeneous equilibria satisfying the Penrose instability condition}

In this section, we provide a detailed discussion of the patially homogeneous profiles that satisfy the Penrose instability conditions \eqref{Penrose-kin-elec}--\eqref{Penrose-kin-ions}--\eqref{Penrose-kin-several} for the previous Vlasov--Poisson systems.

\medskip

One important idea to keep in mind is that for plasmas described by $\mathrm{(VP)}_{e}$ (or $\mathrm{(VP)}_{i}$), homogeneous equilibria which have the shape of \emph{one single bump} are \textit{spectrally} stable \cite{Pen60,Davidson}, that is do not give rise to some exponentially growing mode in the linearized system. More precisely, in one dimensional plasmas, integrable homogeneous equilibria that have exactly one change of monotonicity are stable; in higher dimension, radially symmetric equilibria are also stable (even better, in dimension greater than or equal to three, all radial equilibria enjoy stability \cite{MouhotVillani}). It turns out that liner and nonlinear asymptotic stability also hold true for such equilibria \cite{MouhotVillani}, in connection with the celebrated Landau damping effect. In particular, Maxwellian equilibria are thus always stable.  Therefore, to see instability, we need to go beyond this class.

As briefly alluded to in the introduction, the most typical instabilities of homogeneous equilibria in plasmas are the so-called two-stream instabilities. Since their discovery \cite{BohmGross,PierceH}, they have been the object of numerous studies in the physics literature, in several settings, be it kinetic (Vlasov) or hydrodynamic (Euler), see e.g. \cite{KrallT,Davidson}. As a matter of fact, the multiphasic framework allows to treat the Vlasov or Euler version of these instabilities in a unified manner.

We have to emphasize that virtually all such studies deal with the equations posed in the whole space $\R^d$ so that all frequencies are allowed, whereas as we work on $\T^d$, we only allow frequencies in the grid $\Z^d$, which imposes some additional constraints compared to the physics literature.

In the following discussion, we most of the time focus on the one-dimensional case $d=1$ and refer to Remark \ref{Rem-Penrose-multi-D} for a straightforward adaptation to the general multi-dimensional case.

\medskip 

\noindent{\bf Hydrodynamic two-stream instability.} The prototypical two-stream equilibrium is 
$$
\mu = a_1 \otimes \delta_{v=v_1}+  a_2 \otimes \delta_{v=v_2},
$$
for $a_1, a_2 \geq 0$, with $a_1 + a_2 = 1$, $v_1, v_2 \in \R^d$, with $v_1 \neq v_2$. It is a classical result in the physics literature that such equilibria are spectrally unstable, which is often referred to as the Buneman instability \cite{Bun}.

This is in this case that we recover the instability result of \cite{CGG} which pertains to the two-phase pressureless Euler--Poisson system. Namely, consider the two-phase system
	\begin{equation}
	\left\{
	\begin{aligned}
	\partial_t \rho^\alpha + \Div( \rho^\alpha u^\alpha) &=0, \qquad \qquad \alpha=1,2, \\
	\partial_t u^\alpha + ( u^\alpha \cdot \nabla ) u^\alpha &= - \nabla U,\\
	-\Delta U &= \rho^1 + \rho^2- \int  (\rho^1 + \rho^2) \D x, \\
	\rho^\alpha |_{t = 0} = \rho^\alpha_0, \ u^\alpha |_{t = 0} &= u^\alpha_0.
	\end{aligned}
	\right.
	\end{equation}
 A nonlinear instability result (namely Theorem~\ref{eq:insta-kin-2} in the introduction) was obtained in \cite{CGG}. This result also follows from our abstract instability framework. The linear instability condition in Theorem~\ref{eq:insta-kin-2} corresponds to the Penrose instability criterion~\eqref{Penrose-kin-elec} written for the equilibrium $\mu$. The fact that $v_1 \neq v_2$ is not sufficient is precisely related to the fact that we work on the torus.

For the one-dimensional case, in  \cite{CGG}, assuming without loss of generality that $v_2>v_1$, for instability to hold, this is proved to be equivalent to asking that $\frac{v_2-v_1}{(a_1^{1/3}+ a_2^{1/3})^{3/2}}$ is strictly less than some universal constant. The same result was later obtained in \cite{Bar} and generalized to higher dimensions.
Considering sums of more than two Dirac masses is also possible, and leads as well to instability.

We finally note that \cite{CGG} puts forward an abstract framework that also relies on the method of Grenier \cite{Gr00} but is actually a bit different from ours. In particular,  \cite{CGG} only treats the 1D case, and the analogue of \ref{item:ass_insta_StabHIGHREG}(b) is delicate to obtain, while this aspect seems bypassed in our approach.
The reason is that \cite{CGG} considers $(\rho_1, \rho_2, \partial_x u_1, \partial_x u_2)$ as the main variable, while we focus here directly on $(\rho_1, \rho_2,  u_1, u_2)$.

For the sake of the exposition, let us briefly show that the prototypical two-stream equilibrium above is indeed unstable both for the electron and the ion case, say for $d=1$ and  assuming for simplicity that $a_1=a_2=\frac{1}{2}$ and $v_1=-v_2=w \neq 0$ (which now replaces $(v_2-v_1)/2$ in what follows). We borrow the presentation from \cite{Bar}. According to \eqref{Penrose-kin-elec}--\eqref{Penrose-kin-ions}, we are looking for $(n , \lambda)$ with and $n \in \Z \setminus \lbrace 0 \rbrace$ and $\Re \lambda>0$ such that
\begin{align*}
    \frac{1}{2}\frac{1}{(\lambda+inw)^2}+\frac{1}{2}\frac{1}{(\lambda-inw)^2}=-e(n), 
    \ \  e(n)\vcentcolon=\begin{cases}
1, & \text{for electrons},\\[1ex]
\dfrac{1+|n|^2}{|n|^2}, & \text{for ions}.
\end{cases}
\end{align*}
It gives after simplification that it is enough to find a root with positive real part to the fourth order polynomial
\begin{align*}
    P(X)=e(n)\left(X^2+(nw)^2\right)^2
+
X^2-(nw)^2.
\end{align*}
Since $P$ has real coefficients and is even, there exists a growing mode if and only if not all four roots of $P$ (counted with multiplicity) are purely imaginary. Hence, it is sufficient to show that the polynomial
\begin{align*}
    Q(X)=P(iX)=e(n)\left(X^2-(nw)^2\right)^2
-
X^2-(nw)^2
\end{align*}
has strictly less than four real roots (counted with multiplicity). One direct study of variation of $Q$ shows that, since $Q(x) \rightarrow +\infty$ as $x \rightarrow \pm \infty$ and $Q(\pm nw)<0$, that $Q$ has only one simple root on $(-\infty, -nw)$ and only one simple root on $(nw, +\infty)$. Since $Q$ is quartic and even, there exists a growing mode if and only if $Q(0)<0$, that is $e(n) (nw)^2<1$, which reads
\begin{align*}   
e(n)n^2 \frac{(v_2-v_1)^2}{4}<1.
\end{align*}
For $n \neq 0$, the minimum of the previous expression is attained for $\vert n \vert=1$, and it is therefore necessary and sufficient that $v_2-v_1$ is small enough.

\medskip

\noindent{\bf Kinetic two-stream instability.} For clarity of exposure, we focus here on the $\mathrm{(VP)}_e$ case. The kinetic case corresponds to the case where the equilibrium $\mu(v)$ is a smooth integrable function.
In the one-dimensional case, for \eqref{Penrose-kin-elec} to hold Penrose put forward a famous instability condition \cite{Pen60} which reads as follows: there exists a local minimum $\overline{v}$ such that 
$$
\int_{\R} \frac{\mu(v)- \mu(\overline{v})}{(v-\overline{v})^2} \D v >0,
$$
which is necessary and sufficient for the whole space case. If the minimum is flat (i.e. it is reached on an interval) then the condition must hold for all points of the interval. On the torus $\T^d$, the Penrose instability condition becomes more demanding: we require that
\begin{equation}
    \label{eq:Pen}
\int_{\R} \frac{\mu(v)- \mu(\overline{v})}{(v-\overline{v})^2} \D v >1,
\end{equation}
and as a matter of fact this is only a necessary, not sufficient, condition.

\medskip

Let us explain in more details this fact in the case where the minimum is not flat. In view of \eqref{Penrose-kin-elec}, setting $z\vcentcolon=\frac{i\lambda}{n}$, one needs to ensure that there exists $n\in \Z$ with $n \neq 0$ such that $n^2$ belongs to the range of the function
$$
z \in \{\Im  z>0\} \mapsto  H(z) \vcentcolon=	\int \frac{\mu(v)}{(v-z)^2}\D v=\int \frac{\mu'(v)}{v-z}\D v.
$$
To this end, we can follow Penrose method \cite{Pen60} (also called the Nyquist diagram method), which is based on the argument principle from complex analysis, and that we explain now .

For $R>0$ and $\eps>0$, let us consider a closed oriented counter-clockwise half-disk $\Omega_{R, \eps}$ in the upper half-plane, made of a segment $L_{R, \eps}$ at height $\eps$ from $-R+i\eps$ to $R+i\eps$ and of a half-circle $C_{R, \eps}$ joining these two endpoints. By the argument principle, since $H$ is holomorphic on $\{\Im z >0\}$, it is enough to prove that, for $R>0$ large enough and $\eps>0$ small enough,  the winding number of $H(\partial \Omega_{R, \eps})$ about some non-zero $n^2$ is positive. Since the zeros of $z \mapsto H(z)-n^2$ are isolated, we can always choose $R$ and $\eps$ such that this function does not vanish on $\partial \Omega_{R, \eps}$

It turns out that the contribution of the large half-circle to the winding number of $C_{R, \eps}$ can actually be discarded. Indeed, for fixed $\eps>0$, we have 
\begin{align*}
    \sup_{z \in C_{R, \eps}} \vert H(z) \vert \leq \frac{4}{R^2} \Vert \mu \Vert_{L^1_v}+\frac{1}{\eps^2}\int_{\vert v \vert>R/2} \vert \mu(v) \vert \, \mathrm{d}v,
\end{align*}
hence for any $\eps>0$, there holds
\begin{align*}
    \sup_{z \in C_{R, \eps}} \vert H(z) \vert  \longrightarrow 0, \ \ R \rightarrow +\infty.
\end{align*}
This proves that the curves $H(C_{R, \eps})$ lies in an arbitrarily small neighbourhood of $0$ for $R$ large enough, and therefore it cannot wind around a positive $n^2>0$. Furthermore, the reasoning above also shows that winding number we are looking for does not depends on $R>0$ for $R$ large enough.

Next, since for any $\eta>0$, we have $H(x+i\eta)\longrightarrow0$ as $x\to\pm\infty $, we are left with the study of the closed curve $x \in \R \mapsto H(x+i\eta)$ for $\eta>0$ small enough. By the classical Plemelj's formula, $x \mapsto H(x+i\eta)$  admits a continuous extension up the boundary when $\eta \rightarrow 0$ and therefore need to to study the curve (still denoted in the same fashion)
$$
x\in \R \mapsto H(x)\vcentcolon =  \mathrm{p.v.} \int  \frac{\mu'(v)}{v-x} \, \mathrm{d}v + i \pi  \mu'(x),
$$
from $-\infty$ to $+\infty$, where $\mathrm{p.v.} \int $ denotes the principal value of the integral. As we have the uniform convergence $\sup_{x \in \R} \vert H(x+i\eta)-H(x) \vert \rightarrow 0$ as $\eta \rightarrow 0$, an homotopy argument shows that it is sufficient to study the winding number of the closed curve $H(\R)$ around any point $n^2>0$. Note that for $x \to \pm \infty$, one has
\begin{align*}
    \frac{\mu'(v)}{v-x}=-\frac{\mu'(v)}{x}-\frac{\mu'(v)v}{x^2}+ \mu'(v)v^2\mathcal{O}(x^{-3}),
\end{align*}
hence after integrating
\begin{align*}
    \mathrm{p.v.} \int  \frac{\mu'(v)}{v-x} \, \mathrm{d}v=-\frac{1}{x^2}\int v \mu'(v) \, \mathrm{d}v+ \mathcal{O}(x^{-3}) 
    =\frac{1}{x^2}\int  \mu(v) \, \mathrm{d}v+ \mathcal{O}(x^{-3}),
\end{align*}
so the curve approaches the origin through the right half-plane.

The curve also crosses the real line at all critical points of $\mu$.
For a point $n^2 \geq 1$ to have positive winding number, it is necessary that there exists $x_0 \in \R$ which is a critical point of $\mu$;  the curve must satisfy $H(x_0)>n^2 \geq 1$ and it must cross the real line from bottom to top, which means this point must be a local minimum. Therefore we find that the Penrose instability condition~\eqref{eq:Pen} is necessary for an instability to exist,
since at a critical point $\overline{v}$ one has
\begin{align*}
    H(\overline{v})= \mathrm{p.v.} \int  \frac{\mathrm{d}}{\mathrm{d}v}(\mu(v)-\mu(\overline{v}))\frac{1}{v-\overline{v}} \, \mathrm{d}v 
    =\int  \frac{\mu(v)-\mu(\overline{v})}{(v-\overline{v})^2} \, \mathrm{d}v,
\end{align*}
so that $H(\overline{v})>1$ is exactly \eqref{eq:Pen}.

\medskip

However it is not sufficient as the curve could have a complex behaviour, crossing the real axis several times and thus preventing the point $n^2$ to have a positive winding number.

Suppose in addition that $\mu$ is even and that $\overline{v}=0$. For any $\eta>0$, one has
\begin{align*}
    H(i\eta)=\int \frac{\mu(v)}{(v-i\eta)^2} \, \mathrm{d}v=-\int \mu(v)\frac{\mathrm{d}}{\mathrm{d}v}\left(\frac{1}{v-i\eta} \right) \, \mathrm{d}v=\int \frac{\mu'(v)}{v-i\eta}  \, \mathrm{d}v,
\end{align*} and therefore, since $\mu'$ is odd, we get
\begin{align*}
    H(i\eta)=\int \frac{v\mu'(v)}{v^2+\eta^2}  \, \mathrm{d}v+ i \eta \int \frac{\mu'(v)}{v^2+\eta^2}  \, \mathrm{d}v=\int \frac{v\mu'(v)}{v^2+\eta^2}  \, \mathrm{d}v \in \R.
\end{align*}
As $\eta \to +\infty$, we have $H(i\eta) \to 0$, while since $\mu'(0)=0$ and $\vert \mu'(v) \vert \lesssim \vert v \vert$, there holds $H(i\eta) \to H(0)$ as $\eta \to 0^+$ with
$$H(0)=\int  \frac{\mu'(v)}{v} \, \mathrm{d}v=\int  \frac{\mu(v)-\mu(0)}{v^2} \, \mathrm{d}v>1,$$ provided that the assumption \eqref{eq:Pen} holds. Hence, by continuity of $\eta \to H(i\eta)$, there exists $\eta_\star>0$ such that $H(i\eta_\star)=1$, which prove that $1^2$ belongs on the range of $H_{\mid \lbrace \mathrm{Im}>0 \rbrace}$. In that case, \eqref{eq:Pen} is therefore a sufficient condition for instability.

Alternatively, one may also follows \cite{GS} and almost explicitly find an unstable eigenmode. Let us explain this computation. First, by continuity, using \eqref{eq:Pen}, one finds $z>0$ such that
$$
\int_{\R} v\frac{\mu'(v)}{v^2+z^2} \D v =1,
$$
and for this $z>0$ one can explicitly define 
$$
r(t,x,v) = e^{z t} \left(v\frac{\mu'(v)}{v^2+z^2} \cos x - z \frac{\mu'(v)}{v^2+z^2} \sin x\right), \qquad \psi(t,x) =  e^{z t} \cos x. 
$$
Since $\mu$ is even, one has
$$
-\Psi''(t,x,v)= \int r(t,x,v) \D v, 
$$
and one can readily check that $r$ satisfies
$$
v\partial_x r - \psi'(x) \mu'(v) = -z r,
$$
which proves that $r$ is indeed an unstable eigenmode for the linearized $\mathrm{(VP)}_e$.

\medskip

\begin{Ex} Denoting by 
$$
\mathcal{M}_{\rho,u,T}(v)\vcentcolon= \frac{\rho}{(2\pi T)^{1/2}} \exp\left( - \frac{|v-u|^2}{2T}\right),
$$ 
the Maxwellian of mass $\rho>0$, velocity $u>0$ and temperature $T>0$, an instance of such unstable equilibrium is thus
$$
\mu(v) =  \mathcal{M}_{1/2,\overline{u},\overline{T}}(v) + \mathcal{M}_{1/2,-\overline{u},\overline{T}}(v),
$$
for $\overline{u}>0$ and $\overline{T}>0$, under the condition tht $\bar u^2 > \bar T$ so that $\mu''(0)>0$. Similar considerations hold for the ion case \eqref{Penrose-kin-ions}, see \cite{HKH}.  
\end{Ex}

\begin{Rem}\label{Rem-Penrose-multi-D}
    The previous reasoning in dimension $d=1$ can be upgraded to cover the general multi-dimensional case: the multidimensional instability criterion can actually be reduced to a one-dimensional criterion along each Fourier direction. Let us detail and explain why.
    
    Define the marginal along an arbitrary unit vector $\sigma \in \mathbb{S}^{d-1}$ of a distribution function $\mu(v)$ with $v \in \R^d$ as
    \begin{align*}
        \mu_\sigma(r)\vcentcolon=\int_{r\sigma + \sigma^\perp} \mu(w) \, \mathrm{d}w=\int_{\sigma^\perp} \mu(r \sigma+w) \, \mathrm{d}w , \ \ v \in \R.
    \end{align*}
    The key observation starting from \eqref{Penrose-kin-elec} is that by the orthogonal decomposition
    \begin{align*}
        v=r\sigma+w, \ \ r \in \R, \ \ w \in \sigma^\perp,
    \end{align*}
    if $\sigma=n /\vert n \vert$ (with $n \in {\Z^d \setminus \lbrace 0 \rbrace}$) then we have $n \cdot v=\vert n \vert r$. Hence, \eqref{Penrose-kin-elec} can be reformulated as
    \begin{align*}
        -1=\int_{\R^d} \frac{\mu(v)}{(\lambda + in \cdot v)^2}  \, \mathrm{d}v =\int_{\R} \int_{\sigma^\perp} \frac{\mu(r\sigma+w)}{(\lambda + i\vert n \vert r)^2}  \, \mathrm{d}r \, \mathrm{d}w=\int_{\R} \frac{\mu_\sigma(r)}{(\lambda + i\vert n \vert r)^2}  \, \mathrm{d}r.
    \end{align*}
    Invoking the analysis performed in the one-dimensional case above with $\mu_\sigma$ instead of $\mu$, we obtain the following: if there exists $n \in {\Z^d \setminus \lbrace 0 \rbrace}$ such that, for $\sigma=n/\vert n \vert$, the distribution $\mu_\sigma$ is even with one local minimum $\overline{r}=0$ and satisfies 
    \begin{align*}
        \int  \frac{\mu_\sigma(r)-\mu_\sigma(0)}{r^2} \, \mathrm{d}r>\vert n \vert^2,
    \end{align*}
    then the Penrose instability condition \eqref{Penrose-kin-elec} is verified. In particular, if for some $1 \leq j \leq d$, the distribution $\mu_{e_j}$ is even with one local minimum $\overline{v}=0$ and satisfies
    \begin{align*}
        \int  \frac{\mu_{e_j}(r)-\mu_{e_j}(0)}{r^2} \, \mathrm{d}r>1,
    \end{align*} then the Penrose instability condition holds.
\end{Rem}

\medskip

 \noindent{\bf Electron-ion two stream instability.}

 For the case of dynamical electrons and ions, that is to say for $\mathrm{(VP_{i,e})}$, on can take
$$
G_i(v) =  \mathcal{M}_{1, 0, T_i}(v), \qquad G_e(v) =  \mathcal{M}_{1, v_e, T_e}(v),
$$
(which, each on their own, would be stable for $\mathrm{(VP)}_{i}$ or $\mathrm{(VP)}_{e}$,  \cite{Pen60,MouhotVillani,GagnebinIacobelli}) and the combination of the dynamics of electrons and ions around these equilibria can give rise to instability in some suitable ranges of parameters \cite{Davidson}. Indeed, the Penrose instability condition~\eqref{Penrose-kin-several} for $\mathrm{(VP_{i,e})}$ is equivalent to the Penrose instability condition for $\mathrm{(VP_{e})}$ for the equilibrium $G= G_i + G_e$. As a result this comes down to a two-stream instability as described above.

In the case of \emph{cold ions}, $G_i$ degenerates to 
$$
G_i = \delta_{v=0},
$$
and \eqref{Penrose-kin-several} turns into finding $\lambda\in \C$ with $\Re \lambda>0$ and a non-zero $n \in \N$ such that
$$
-1=\frac{1}{\lambda^2} + \int \frac{G_e(v)}{(\lambda + in \cdot v)^2} \D v ,
$$  
which we can recast as 
$$
-1=z^2 \left( \int \frac{G_e(v)}{(v - z)^2} \D v -n^2\right) = z^2 \left( \int \frac{G_e'(v)}{v - z} \D v -n^2\right),
$$
with $z\vcentcolon= \frac{i\lambda}{n}$, $\Im z >0$. We are therefore looking for a zero $z \in \C$ with $\Im z >0$ of the function
\begin{align*}
    D(z)\vcentcolon=1+z^2 \left( \int \frac{G_e'(v)}{v - z} \D v -n^2\right).
\end{align*}
Note that $G'_e(v)=-\frac{v-v_e}{T_e} G_e(v)$ and on can assume without loss of generality that $v_e>0$.
We rely on the same Penrose's method \cite{Pen60} and argument principle as above,  with a closed oriented counter-clockwise half-disk $\Omega_{R, \eps}$ in the upper half-plane, that is the union of a segment at height $\eps$ from $-R+i\eps$ to $R+i\eps$ and of a half-circle $C_{R, \eps}$ closing the contour. As before, it is enough to prove that, for $R>0$ large enough and $\eps>0$ small enough,  the winding number of $D(\partial \Omega_{R, \eps})$ about zero is positive.

As $\int \frac{G_e(v)}{(v - z)^2} \D v=z^{-2}+O(z^{-3})$ for $\vert z \vert \to +\infty$ uniformly in $v$, we have $D(z)=2-n^2z^2 + O(z^{-1})$ and hence
$$D(R e^{i\theta})=-n^2 R^2 e^{2i\theta} (1+\eps_R(\theta)), \ \ \sup_{\theta} \vert \eps_R(\theta) \vert \rightarrow 0, \ R \to +\infty.$$
In particular, $\mathrm{arg}(D(R e^{i\theta}))=-\pi+2\theta + o(1)$, which prove that the half-circle defined for $\theta=0$ to $\pi$ has a winding number equal to $1$ as $R\to +\infty$.

Let us now study the boundary value of $D(z)$ from the upper half-plane as $\eps \rightarrow 0^+$. By the Plemelj formula, the extension of $D$ on the real axis is given by the curve $$
x\in \R \mapsto H(x)\vcentcolon =  1+x^2 \left( \mathrm{p.v.} \int \frac{G'_e(v)}{v-x}-n^2\right) + i \pi x^2 G'_e(x).
$$
 As $x\to \pm \infty$, $H \rightarrow - \infty$, and we note that for $x$ near $\pm \infty$ the curve lies in the half-plane $\{ \Re z <0\}$. More precisely, $H(x)=2-n^2 x^2+o(1)$ as $x \to \pm \infty$. Furthermore, we have
 $$\Im H(x)=-\frac{\pi}{T_e}(x-v_e)x^2G_e(x),$$
 so $\Im \, H(x)>0$ for $x<v_e$ and $\Im \, H(x)<0$ for $x>v_e$, with $x \neq 0$. We also have $H(0)=1$ so the curve only touches the real axis at the point 1 without crossing it. The only other point where the imaginary part cancels is $x=v_e$, where
$$
H(v_e) = 1+v_e^2 \left( \int \frac{G_e'(v)}{v-v_e} \D v -n^2\right) = 1-v_e^2 \left( \frac{1}{T_e} + n^2 \right).
$$
As a consequence, $H(v_e)<0$ if and only if the crossing of the real axis (from above) happens on the negative real axis. We conclude that, when $H(v_e)<0$ (resp. $H(v_e) > 0$), the winding number of $H(\R)$ around zero is zero (resp. is $-1$).
Adding the contribution $+1$ from the half-circle detailed previously, we conclude to the following: for instability to occur, it becomes necessary and sufficient to ensure $H(v_e)<0$, that is to say 
$$
v_e^2 \left( \frac{1}{T_e} + n^2 \right) >1.
$$
As we can take $n$ as large as we want, this becomes $v_e >0$, that is $v_e \neq 0$ by symmetry. Thus, in this idealized cold-ion model, any non-zero relative drift
produces an unstable Fourier mode, possibly at arbitrarily high spatial
frequency.

\begin{Rem}[Consequences on the invalidity of the quasineutral limit]
\label{rem:quasi}
The quasineutral limit of the Vlasov--Poisson system (either for electrons or ions) is the limit $\eps \rightarrow 0$ of
\begin{equation*}
\mathrm{(VP_e)^\eps}\left\{ 
\begin{gathered}
\partial_t f_\eps + v \cdot \nabla_x f_\eps - \nabla_x U_\eps\cdot \nabla_v f_\eps = 0,\\
		- \eps^2 \Delta_x U_\eps = \int f_\eps \D v- \iint f_\eps \D v \D x,
\end{gathered}
\right.
\end{equation*}
and
\begin{equation*}
\mathrm{(VP_i)}^\eps\left\{ 
\begin{gathered}
\partial_t f_\eps + v \cdot \nabla_x f - \nabla_x U_\eps \cdot \nabla_v f_\eps = 0,\\
		e^{U_\eps}-\eps^2\Delta_x U_\eps = \int f_\eps \D v,
\end{gathered}
\right.
\end{equation*}
in which the parameter $\eps$ is to be understood as the scaled Debye length of the plasma, which is small in several physical situations so that the approximation $\eps\to 0$ is very common.
This limit has been notably studied from the mathematical point of view in \cite{Gr96,Br00} in the case of $(\mathrm{VP}_{\mathrm{e}})$, and in \cite{HK11,HKR,GPI2} in the case of $(\mathrm{VP}_{\mathrm{i}})$ or its linearized Poisson version. The formal limits, directly obtained when setting $\eps=0$ in the Poisson equation, are \emph{singular} Vlasov equations which display a loss of derivative with respect to the force term, and which are ill-posed in general in the Sobolev framework \cite{HKN,Bar}. Ill-posedness stems for the instability of profiles which satisfy a Penrose instability condition such as \eqref{Penrose-kin-elec} or \eqref{Penrose-kin-ions}. Because of these instabilities, the formal limit does not hold in general, as justified in \cite{HKH}. The latter relies on the combination of a suitable rescaling in $\eps$ with a nonlinear Lyapunov instability result, requiring smoothness of the unstable profile. In view of Theorem \ref{thm-insta-kin-VP}, the multiphasic approach allows to remove this requirement and one can then applies the exact same method as developed \cite{HKH}.

As a consequence, when the Penrose instability condition holds, the formal quasineutral limit  $\eps \rightarrow 0$ fails for (possibly rough) solutions of the Vlasov--Poisson system for electrons and for ions. 
\end{Rem}

\subsection{Multiphasic reformulation}

 In order to prove Theorem \ref{thm-insta-kin-VP}, we shall rely on the multiphasic approach and thus study the multiphasic version of the Vlasov--Poisson systems,
\begin{equation}
\label{eq:VPinstaproof}
\left\{ 
\begin{gathered}
\partial_t \rho^\alpha + \Div ( \rho^\alpha v^\alpha) = 0,\\
\partial_t v^\alpha + (v^\alpha \cdot \nabla) v^\alpha = - \nabla U,\\
-\Delta U = \int \rho^\alpha  \D \mu(\alpha) - \iint \rho^\alpha \D \mu(\alpha) \D x \qquad \text{or} \qquad e^U-\Delta U = \int \rho^\alpha  \D \mu(\alpha),
\end{gathered}
\right.
\end{equation}
where $(I, \mu)$ is a fixed set of labels and  $\mu$ is assumed to be a probability measure. We consider initial data lying in a small neighborhood of stationary solutions of the form
\begin{equation*}
\rho^\alpha \equiv 1, \qquad v^\alpha \equiv \mathcal{V}(\alpha), \qquad U \equiv 0,
\end{equation*}
where $\mathcal{V} :I \to \R^d$ is a function satisfying $\int \vert \mathcal{V}(\alpha) \vert^p \D \mu(\alpha) < + \infty$. In other words we write
\begin{equation}\label{def-perturb-VP}
\rho^\alpha = 1+r^\alpha, \qquad v^\alpha = \mathcal{V}(\alpha)+w^\alpha, \qquad U \equiv \Phi,
\end{equation}
and focus on the equations satisfied by the perturbation $(r^\alpha, w^\alpha, \Phi)$.
In the electron case, we are led to study the following equation on the perturbation $(r^\alpha, w^\alpha, \Phi)$ defined by \eqref{def-perturb-VP}:
\begin{equation}
\label{eq:perturbVP-elec}
\left\{  
\begin{gathered}
\partial_t r^\alpha + \mathcal{V}(\alpha) \cdot \nabla r^\alpha + \Div (w^\alpha) + \Div(r^\alpha w^\alpha) = 0,\\
\partial_t w^\alpha + (\mathcal{V}(\alpha) \cdot \nabla) w^\alpha + (w^\alpha \cdot \nabla) w^\alpha  = - \nabla \Phi,\\
- \Delta \Phi = \int r^\alpha \D \mu(\alpha),
\end{gathered}
\right.
\end{equation}
while on the ion case, the equation for the perturbation reads
\begin{equation}
\label{eq:perturbVP-ions}
\left\{  
\begin{gathered}
\partial_t r^\alpha + \mathcal{V}(\alpha) \cdot \nabla r^\alpha + \Div (w^\alpha) + \Div(r^\alpha w^\alpha) = 0,\\
\partial_t w^\alpha + (\mathcal{V}(\alpha) \cdot \nabla) w^\alpha + (w^\alpha \cdot \nabla) w^\alpha  = - \nabla \Phi,\\
e^{\Phi}-1- \Delta \Phi = \int r^\alpha \D \mu(\alpha).
\end{gathered}
\right.
\end{equation}
 We shall not dwell on the proof of Theorem \ref{thm-insta-kin-VP-several}, for which we can argue similarly. We have indeed already noted that the multiphasic system~\eqref{eq:VP-several} for dynamic electrons and ions is structurally very close to the multiphasic system~\eqref{eq:VPinstaproof} for electrons. As a result, the arguments which will be developed below to treat the electron case can be fairly straightforwardly adapted to the dynamical electron and ion case.

We need to introduce the multiphasic counterpart of the Penrose instability conditions \eqref{Penrose-kin-elec} and \eqref{Penrose-kin-ions}.
\begin{Def}\label{def:Penrose-multi} \begin{itemize}
\item
The stationary solution $(1,\mathcal{V}(\alpha))_{\alpha \in I}$ satisfies the Penrose instability condition for electrons associated to the measure $\mu$ if there exist $\lambda \in \C$ with $\Re \lambda>0$ and $n \in \Z^d {\setminus \{0\}}$ such that
	\begin{equation}
	\label{eq:penrose-multi-VPelec}
	\int \frac{1}{(\lambda + in \cdot \mathcal{V}(\alpha))^2}\D \mu(\alpha) = -1.
	\end{equation}

 \item    The stationary solution $(1,\mathcal{V}(\alpha))_{\alpha \in I}$ satisfies the Penrose instability condition for ions  associated to the measure $\mu$ if there exist $\lambda \in \C$ with $\Re \lambda>0$ and $n \in \Z^d {\setminus \{0\}}$ such that
	\begin{equation}
	\label{eq:penrose-multi-ions}
	\frac{|n|^2}{1+|n|^2}\int \frac{1}{(\lambda + in \cdot \mathcal{V}(\alpha))^2}\D \mu(\alpha) = -1.
	\end{equation}
    \end{itemize}
\end{Def}

Under this condition, the general framework developed in Section \ref{Section-FrameworkInstab} yields the following.
\begin{Prop}\label{prop-check-assump-Instab}
The Assumptions \ref{item:ass_insta_expansion}--\ref{item:ass_insta_StabHIGHREG}--\ref{item:ass_insta_RealForce}  are satisfied for \eqref{eq:perturbVP-elec} (resp. \eqref{eq:perturbVP-ions}). Under the Penrose instability condition~\eqref{eq:penrose-multi-VPelec} (resp. \eqref{eq:penrose-multi-ions}),  Assumption \ref{item:ass_insta_LINEARIZED} is also satisfied for \eqref{eq:perturbVP-elec} (resp. \eqref{eq:perturbVP-ions}). Furthermore, using the same notation from the abstract assumptions, we have
    \begin{align*}
        P_{1,n}(X)(\alpha)=\frac{in}{\vert n \vert^2}\int \sigma \, \mathrm{d}\mu \neq 0,
    \end{align*} 
    where $X(\alpha)=(\sigma(\alpha), \eta(\alpha))$ and $n \in \Z^d {\setminus \lbrace 0 \rbrace}$ correspond to the exponential growing mode from \ref{item:ass_insta_LINEARIZED}(b).
\end{Prop}
From the general Theorem \ref{thm:instability}, we can therefore infer the following nonlinear instability result. 
\begin{Thm}
\label{thm-insta-multi-VP}
Let $(1,\mathcal{V}(\alpha))_{\alpha \in I}$ satisfy the Penrose instability condition from formula~\eqref{eq:penrose-multi-VPelec} (resp.\ \eqref{eq:penrose-multi-VPelec}) associated with the measure $\mu$. Assume also that $\boldsymbol{\mathcal{V}} \in L^p_\alpha$ for some $p \in [1, \infty]$. Let $k_0 \in \N$.  There is $\kappa>0$ such that for all $s\in \N$ and for all $\delta>0$, there exists an initial condition $(\rho^\alpha_0, v^\alpha_0)_{\alpha \in I}$  satisfying
\begin{equation}
\label{eq:assumptioninstamulti-VP}
\esssup_{\alpha \in I} \left\| (\rho^\alpha_0-1, v^\alpha_0- \mathcal{V}(\alpha)) \right\|_{W^{s,\infty}_x} \leq \delta,
\end{equation} 
with associated solution $(\rho^\alpha, v^\alpha)_{\alpha \in I}$ to~\eqref{eq:VPinstaproof}  and a positive time $T_\delta = O(|\log \delta|)$ such that
\begin{equation}
  \| \nabla U[\rhorho](T_\delta) \|_{W^{-k_0,1}_x}  \geq \kappa.
\end{equation} 
\end{Thm}
\begin{proof}
   Here, the force field does not depend on $\alpha$. We can therefore pick $\xi(\alpha)=n$ (independent of $\alpha$), which gives 
    $$\int n \cdot P_{1,n}(X)(\alpha) \D \mu(\alpha) \neq 0.$$
    Thanks to Proposition~\ref{prop-check-assump-Instab} and Theorem \ref{thm:instability}, we obtain
    \begin{align*}
        \left\| \int n \cdot G^\alpha[\ZZ|_{t = T_{\delta}}] \, \mathrm{d}\mu(\alpha) \right\|_{W^{-k_0,1}_x} \geq \kappa,
    \end{align*}
    from which we infer the desired conclusion.
\end{proof}
Let us now justify why Theorem \ref{thm-insta-kin-VP} can be deduced from Theorem \ref{thm-insta-multi-VP}.
\begin{proof}[Proof of Theorem \ref{thm-insta-kin-VP}]
Let $G$ satisfy the assumptions of Theorem \ref{thm-insta-kin-VP}.
Let us take 
$$I=\R^d, \qquad \D \mu=G \D \alpha, \quad \text{and} \quad \mathcal{V}(\alpha)=\alpha.$$
We can therefore apply Theorem~\ref{thm-insta-multi-VP} which gives for $s\in \mathbb{N}$ and $\delta>0$ an initial condition $(\rho^\alpha_0, v^\alpha_0)_{\alpha \in I}$ satisfying \eqref{eq:assumptioninstamulti-VP}. We accordingly set
$$
f_{0,\delta} = \int_I \rho^\alpha_0 \otimes \delta_{v= u^\alpha_0} \, \D \mu(\alpha).
$$
We then observe that \eqref{eq:assumptioninstamulti-VP} implies
$$
\left\| \langle f_{0,\delta} - G, (1+\vert v \vert^2)^{\frac{p}{2}} \varphi\rangle \right\|_{W^{s,\infty}_x} \leq C(\delta),
$$    
where $C$ is a differentiable function such that $C(0)=0$. We thus deduce that up to choosing $\delta$ small enough,
\eqref{eq:assumptioninsta-VP} holds. The conclusion follows.

\end{proof}

\begin{Rem}
\label{rem:perturb-insta2} 
Let us elaborate on the second item of Remark~\ref{rem:perturb-insta}: namely, we want to show that if the equilibrium $G$ is an integrable function, the unstable solution also remains an integrable function. 

     To this end, one may argue as in the proof of Proposition~\ref{prop:reg-alpha} and use in addition the uniform estimates obtained in the proof of Theorem~\ref{thm:instability} in order to reach the relevant times of instability (up to picking $\kappa>0$ not too large in Theorem~\ref{thm:instability}).  We make the same assumption as in Theorem~\ref{thm-insta-multi-VP} (and keep the same notations). 
   For $\ww$ solving~\eqref{eq:perturbVP-elec}, it is sufficient to show that $\| \DDD_\alpha \ww\|_{L^\infty([0,T_\delta]; L^\infty_\alpha L^\infty_x)}<1/2$; we recall that this implies then that the map $\alpha \mapsto \alpha + w^\alpha(t,x)$ is a  $\mathscr{C}^1$ diffeomorphism, uniformly in $(t,x)$. 
It follows from the proof that the unstable solution can then we rewritten as
    \begin{equation}\label{eq:formule-f-again}
f(t,x,v) = \left(1+r^{\zeta(t,x,v)}(t,x)) G(\zeta(t,x,v) \right)  |\det \DDD_v \zeta(t,x,v)|,
\end{equation}
where   $v\mapsto \zeta(t,x,v)$ denotes the inverse of  $\alpha \mapsto \alpha + w^\alpha(t,x)$.

   As for~\eqref{eq:defDw}, $ \DDD_\alpha \ww$ satisfies \begin{equation}\label{eq:defDw-again}
	\partial_t  \DDD_\alpha w^\alpha +( w^\alpha \cdot \nabla)  \DDD_\alpha w^\alpha = - \DDD_\alpha w^\alpha : \nabla   w^\alpha - \nabla w^\alpha,
	\end{equation}
    with $\| \DDD_\alpha w^\alpha|_{t=0}\|_{L^{\upinfty}_x} \lesssim \delta$, according to~\eqref{eq:assumptioninstamulti-VP} and~\eqref{eq:CNS_EM-vpelec}. But arguing as in the beginning of the proof of Theorem~\ref{thm:instability}, one can check that
$$
\int_0^{T_\delta} \| \nabla w^\alpha(t)\|_{L^{\upinfty}_x} \D t 
\leq C \kappa_A,
$$
where the parameter $\kappa_A>0$ can be chosen as small as required (this amounts to choosing $\kappa$ small in Theorem~\ref{thm-insta-multi-VP}).
Therefore, picking $\kappa_A$ small enough, from the Duhamel formula applied for~\eqref{eq:defDw-again}, we end up with the desired estimate $\| \DDD_\alpha \ww\|_{L^\infty([0,T_\delta]; L^\infty_\alpha  L^\infty_x)}<1/2$, and we can conclude.
\end{Rem}

It only remains to prove Proposition~\ref{prop-check-assump-Instab}. To ease readability, we now split out the analysis for the electron and for the ion case.

\subsection{Proof of Proposition~\ref{prop-check-assump-Instab} in the electron case}

\begin{proof}[Proof of Proposition \ref{prop-check-assump-Instab} in the electron case] Let us check the three assumptions \ref{item:ass_insta_expansion}, \ref{item:ass_insta_StabHIGHREG} and \ref{item:ass_insta_LINEARIZED}, the assumption \ref{item:ass_insta_RealForce} being directly obtained by inspection.

\medskip

\noindent

\underline{$\star$ \ref{item:ass_insta_expansion}}: In that case, the force field does not depend on $\alpha$. Since the coupling with the Poisson equation is linear (in $\Phi$) we set for all $N \geq 1$ and $k \geq 2 $:
$$G[\rr, \ww]=G_1[\rr, \ww]\vcentcolon=\nabla \Delta^{-1} \int r^\alpha \, \mathrm{d}\mu(\alpha), \ \ G_k\vcentcolon=0, \ \ R_N\vcentcolon=0,$$
for which it is clear that \ref{item:ass_insta_expansion} is satisfied.

\medskip

\noindent

\underline{$\star$ \ref{item:ass_insta_StabHIGHREG}}: By relying on the linearity of the force field and on elliptic regularity (similarly to what we have done for the bound at high regularity in Section \ref{Subsec:VPelectrons}), we directly obtain that the assumption is satisfied.

\medskip

\noindent

\underline{$\star$ \ref{item:ass_insta_LINEARIZED}}: The linearized system associated with~\eqref{eq:perturbVP-elec}, which corresponds to~\eqref{eq:system_of_study_linear} in this context, is:
	\begin{equation}
	\label{eq:VPlin-elec}
	\left\{  
	\begin{gathered}
	\partial_t r^\alpha + \mathcal{V}(\alpha) \cdot \nabla r^\alpha + \Div(w^\alpha)  = 0,\\
	\partial_t w^\alpha + (\mathcal{V}(\alpha) \cdot \nabla) w^\alpha   = - \nabla \Phi,\\
	- \Delta \Phi = \int r^\alpha \D \mu(\alpha).
	\end{gathered}
	\right.
	\end{equation}
For given $n \in \Z^d$, accordingly with the statement of Assumption~\ref{item:ass_insta_LINEARIZED}, we first need to show that there exists $\theta_0>0$ uniform in $n$ and $C_n \geq 0$ such that for all $(\boldsymbol{\sigma}_0,\boldsymbol{\xi}_0) \in L^\infty(I;\R^{1+d})^2$, there exists a unique solution of~\eqref{eq:VPlin-elec} of the form
\begin{equation*}
r^\alpha(t,x) \vcentcolon = \sigma(t,\alpha) \exp(i n \cdot x), \qquad w^\alpha(t,x) \vcentcolon =  \xi(t,\alpha) \exp(i n\cdot x),
\end{equation*}
with $\boldsymbol{\sigma}(0) = \boldsymbol{\sigma}_0$ and $\boldsymbol{\xi}_0 = \boldsymbol{\xi}(0)$, and that it satisfies
\begin{equation*}
    \| \sigma(t), \xi(t) \|_{L^{\upinfty}_\alpha} \leq C_n \| \sigma_0, \xi_0 \|_{L^{\upinfty}_\alpha}  \exp(\theta_0 t).
\end{equation*}
Formally, $(\boldsymbol{\sigma},\boldsymbol{\xi})$ needs to solve
\begin{multline}
\label{eq:fourier}
\partial_t \begin{pmatrix}
\sigma(t,\alpha) \\ \xi(t,\alpha)
\end{pmatrix}
+ \begin{pmatrix}
in \cdot \mathcal V(\alpha) & in \\
0 & in \cdot \mathcal V(\alpha) \Id
\end{pmatrix} \cdot \begin{pmatrix}
\sigma(t,\alpha) \\ \xi(t,\alpha)
\end{pmatrix} \\
= \begin{pmatrix}
0 \\ \displaystyle{- \frac{in}{|n|^2}\int} \sigma(t,\alpha') \D \mu(\alpha')
\end{pmatrix}.
\end{multline}
To lighten the notation, we denote by
\begin{align*}
L_n(\boldsymbol{\sigma}, \boldsymbol{\xi}) &\vcentcolon = - \begin{pmatrix}
in \cdot \boldsymbol{\mathcal V} & in \\
0 & in \cdot \boldsymbol{\mathcal V} \Id
\end{pmatrix} \cdot \begin{pmatrix}
\boldsymbol{\sigma} \\ \boldsymbol{\xi}
\end{pmatrix}
, \\
K_n(\boldsymbol{\sigma},\boldsymbol{\xi}) &\vcentcolon = \begin{pmatrix}
0 \\ \displaystyle{- \frac{in}{|n|^2}\int} \boldsymbol{\sigma}(\alpha) \D \mu(\alpha)
\end{pmatrix}.
\end{align*}
The operator $L_n$ generates a $\mathscr{C}^0$ semigroup in $L^1(I;\R \times \R^d)$ and $K_n$ has finite range, hence is continuous (even compact) in $L^1(I;\R \times \R^d)$. Therefore, by~\cite[Chapter 3, Theorem 1.1]{pazy2012semigroups}, $L_n + K_n$ generates a $\mathscr \mathscr{C}^0$ semigroup in $L^1(I;\R \times \R^d)$, that we denote by $S_n(t)$.  

Let us call $\bar \gamma_n$ the largest real part of the eigenvalues of $L_n + K_n$. We will use the framework of Shizuta~\cite{shizuta1983classical} to show the following bound for $S_n(t)$ in $L^\infty(I;\R \times \R^d)$: if $\gamma_n>\max(\bar \gamma_n,0)$, then, there exists $C_n$ such that for all $(\boldsymbol{\sigma},\boldsymbol{\xi}) \in L^\infty(I;\R \times \R^d)$ and $t \geq 0$,
\begin{equation}
\label{eq:boundSn}
\| S_n(t)\cdot (\boldsymbol{\sigma},\boldsymbol{\xi}) \|_{L^{\upinfty}_\alpha} \leq C_n \exp(\gamma_n t) \| ( \boldsymbol{\sigma},\boldsymbol{\xi} )\|_{L^{\upinfty}_\alpha}.
\end{equation}
Once this is done, the validity of Assumption~\ref{item:ass_insta_LINEARIZED}(a) follows directly as soon as $\sup_n \bar \gamma_n < +\infty$, with any $\gamma_0>\sup_n \max(\bar \gamma_n,0)$. Let us prove our claim and postpone the proof of $\sup_n \bar \gamma_n < +\infty$.

	A slight difficulty in proving this claim is that $(S_n(t))$ is not a $\mathscr \mathscr{C}^0$ semigroup in $L^\infty(I;\R \times \R^d)$, but only in $L^1(I;\R \times \R^d)$: for a given $(\boldsymbol{\sigma}, \boldsymbol{\xi}) \in L^\infty(I;\R \times \R^d)$, $S_n(t)\cdot (\boldsymbol{\sigma}, \boldsymbol{\xi})$ does not necessarily converge strongly in $L^\infty(I;\R \times \R^d)$ towards $(\boldsymbol{\sigma}, \boldsymbol{\xi})$ as $t \to 0$. Therefore, it is rather unusual to show $L^\infty$ bounds on this semigroup. However, we can bypass this problem by working in the dual: the proof consists in showing the same bound, but for the adjoint semigroup in $L^1(I;\R \times \R^d)$, which is a standard $ \mathscr{C}^0$ semigroup satisfying usual assumptions of perturbation theory.

	Let $L^*_n$ and $K^*_n$ the formal adjoints of $L_n$ and $K_n$ respectively, that is:
	\begin{equation*}
	L^*_n(\boldsymbol{\sigma}, \boldsymbol{\xi}) \vcentcolon =  \begin{pmatrix}
	in \cdot \boldsymbol{\mathcal V} & 0 \\
	in & in \cdot \boldsymbol{\mathcal V} \Id
	\end{pmatrix} \cdot \begin{pmatrix}
	\boldsymbol{\sigma} \\ \boldsymbol{\xi}
	\end{pmatrix}
	, \ \
	K^*_n(\boldsymbol{\sigma},\boldsymbol{\xi}) \vcentcolon =  \begin{pmatrix}
 \displaystyle{\frac{in}{|n|^2} \cdot \int} \sigma(\alpha) \D \mu(\alpha) \\0
	\end{pmatrix}.
	\end{equation*}
	
	For the same reason as for $L_n + K_n$, the operator $L_n^* + K_n^*$ generates a $\mathscr \mathscr{C}^0$ semigroup $(T_n(t))$ on $L^1(I;\R \times \R^d)$, and straightforward arguments show that for all $t$, the restriction of $T_n(t)$ to $L^\infty(I;\R \times \R^d)$ is nothing but $S_n^*(t)$. In other words, $S^*_n(t)$ has a continuous extension to $L^1(I;\R \times \R^d)$, and from now on, we keep the notation $S_n^*(t)$ instead of $T_n(t)$ for this extension.

	By duality, showing~\eqref{eq:boundSn} is equivalent to showing that there is $C_n$ such that for all $t$, $S_n^*(t)$ satisfies 
	\begin{equation}
    \label{eq:estimate_semigroup_adjoint}
	\|  S_n^*(t)\cdot (\boldsymbol{\sigma},\boldsymbol{\xi}) \|_{L^1_\alpha} \leq C_n \exp(\gamma_n t) \| \boldsymbol{\sigma},\boldsymbol{\xi} \|_{L^1_\alpha}.
	\end{equation}
	So let us focus on the proof of this result. It is a direct application of the following theorem of Weyl type, which is a simplification of the results from~\cite{shizuta1983classical}, extending previous results from~\cite{vidav1970}. We do not write the proof of this result and refer directly to~\cite{shizuta1983classical}. Here, we use the assumption that the perturbation is compact and not the notion of "$L$-smoothing" operator as in~\cite{shizuta1983classical}. But it is obvious that a compact operator is $L$-smoothing for any $L$ and we do not need to develop further these considerations.
    \begin{Thm}[Easy version of Theorems~1.1, 1.2 and 1.3 in~\cite{shizuta1983classical}]
    \label{thm:shizuta}
        Let $\mathcal X$ be a Banach space, and let $L$ be an operator generating a $\mathscr \mathscr{C}^0$ semigroup $(\exp(tL))$ on $\mathcal X$. Assume that there exists $\beta_0 \in \R$ such that for all $\beta> \beta_0$, there exists $C = C(\beta)>0$ with
        \begin{equation*}
           \forall X \in \mathcal X, \quad \forall t \geq 0, \qquad \| \exp(tL) \cdot X  \| \leq C \exp(\beta t)\| X\|.
        \end{equation*}

        Let $K$ be a compact operator on $\mathcal X$, and $(S(t))$ be the $\mathscr \mathscr{C}^0$ semigroup generated by $L + K$, as given by~\cite{pazy2012semigroups}. Then:
        \begin{itemize}
            \item For all $\beta > \beta_0$, the spectrum of $L + K$ in the half space $\{ \lambda \in \C \mid \Re\lambda > \beta \}$ consists of a finite number of isolated eigenvalues with finite algebraic multiplicity.
            \item Let $\beta > \beta_0$ be such that $L + K$ has no eigenvalue on the line $\{  \lambda\in \C \mid \Re\lambda = \beta \}$. Let $\lambda_1, \dots, \lambda_N$, $N \in \N^*$ be all the eigenvalues given by the first point, and $\mathcal X_1, \dots, \mathcal X_N$ be the corresponding finite dimensional generalized eigenspaces. There exist $C = C(\beta)$ and continuous projections $P_i: \mathcal X \to \mathcal X_i$, $i =1, \dots, N$ such that for all $X \in \mathcal X$,
            \begin{equation}
            \label{eq:eigenprojection}
              \forall t \geq 0, \qquad  \Big\| S(t) \cdot \Big(X - P_1(X) - \dots -  P_N(X) \Big) \Big\| \leq C \exp(\beta t) \|X\|.
            \end{equation}
        \end{itemize}
    \end{Thm}
    \begin{Rem}
    In particular, with the same notations as in the statement of Theorem~\ref{thm:shizuta}, the estimate~\eqref{eq:eigenprojection} implies
     \begin{equation*}
              \forall t \geq 0, \qquad  \| S(t) \cdot X \| \leq \|S(t) \cdot P_1(X) \| + \dots + \|S(t) \cdot P_N(X) \| + C \exp(\beta t) \|X\|.
            \end{equation*}
        But by finite dimensional considerations relying on Jordan's decomposition, for all $i= 1, \dots, N$, there exists a polynomial $q_i = q_i(t)$ whose order is smaller than $k_i - 1$, being $k_i$ the dimension of $\mathcal X_i$, such that
        \begin{equation*}
            \forall X_i \in \mathcal X_i, \,\forall t\geq 0,\qquad \| S(t) \cdot X_i \| \leq q_i(t) \exp\Big((\Re\lambda_i) t\Big) \|X_i \|. 
        \end{equation*}
        Therefore, we easily deduce from the two above estimates and from the continuity of the projections in play that whenever $\gamma > \beta_0$ is strictly larger than $\max\{\Re\lambda_i, \, i=1, \dots, N\}$, then there exists $C = C(\gamma)$ such that
        \begin{equation}
        \label{eq:consequence_shizuta}
            \forall X \in \mathcal X, \,\forall t\geq 0,\qquad \| S(t) \cdot X \| \leq C \exp(\gamma t) \|X \|. 
        \end{equation}
    \end{Rem}

    \medskip
    
    In our situation, explicit computations show that for all $(\boldsymbol{\sigma}, \xi) \in L^1(I;\R \times \R^d)$,
	\begin{equation*}
	\exp(tL_n^*)\cdot (\boldsymbol{\sigma}, \boldsymbol{\xi}) =  \exp(in \cdot \boldsymbol{\mathcal V}t)\begin{pmatrix}
	1 & 0 \\
	itn &  \Id
	\end{pmatrix} \cdot \begin{pmatrix}
	\boldsymbol{\sigma} \\ \boldsymbol{\xi}
	\end{pmatrix}.
	\end{equation*}
    Therefore, $L^*_n$ satisfies the assumptions of Theorem~\ref{thm:shizuta} with $\beta_0 = 0$. Moreover, $K_n^*$ is of finite range, hence compact. Therefore, Theorem~\ref{thm:shizuta} applies, and~\eqref{eq:estimate_semigroup_adjoint} is an instance of~\eqref{eq:consequence_shizuta}.

    \medskip

Now, it remains to show that $\bar \gamma_0 \vcentcolon = \sup_n \bar \gamma_n < + \infty$, and that under the Penrose instability condition \eqref{eq:penrose-multi-VPelec}, $\bar \gamma_0 >0$ and Assumption~\ref{item:ass_insta_LINEARIZED}(b) holds choosing $\gamma_0 \in (\bar \gamma_0, 2 \bar \gamma_0)$. For a given $n \in \Z^d$, looking for an eigenvalue of $L_n + K_n$ with positive real part is equivalent to looking for solutions to the linearized system~\eqref{eq:VPlin-elec} of the form
\begin{align*}
    (r^\alpha(t,x), w^\alpha(t,x),\Phi(t,x)) = \left(\sigma(\alpha), \eta(\alpha), \frac{1}{|n|^2} \int \sigma \D \mu\right) \exp(\lambda t + in\cdot x),
\end{align*}
for some $\lambda \in \C$ with $\Re \lambda>0$ and $n \in \Z^d {\setminus \{0\}}$, where $(\sigma(\alpha), \eta(\alpha))$ are not zero everywhere. Now plugging this ansatz into~\eqref{eq:VPlin-elec}, we find that $(\sigma,\eta)$ give rise to a solution if and only if for all $\alpha$,
	\begin{gather*}
	(\lambda + in\cdot \mathcal{V}(\alpha)) \sigma(\alpha) + ik \cdot \eta(\alpha) = 0,\\
	(\lambda + in\cdot \mathcal{V}(\alpha)) \eta(\alpha) = -\frac{in}{|n|^2} \int \sigma\D \mu.
	\end{gather*}
	Therefore, using that $\Re\lambda \neq 0$, so that $\lambda + in\cdot \mathcal{V}(\alpha) \neq 0$, we deduce
	\begin{equation}
	\label{eq:CNS_EM-vpelec}
    \begin{aligned}
	\eta(\alpha) &= - \frac{in}{|n|^2} \left(\int \sigma\D \mu \right) \frac{1}{\lambda + in\cdot \mathcal{V}(\alpha)}, \\
    \sigma(\alpha) &= -\left(\int \sigma\D \mu \right) \frac{1}{(\lambda + in\cdot \mathcal{V}(\alpha))^2}.
    \end{aligned}
	\end{equation}
	If $\int \sigma \D \mu = 0$, then $\sigma$ and $\eta$ cancel, which is excluded. So $(\sigma,\eta)$ gives rise to an exponential mode if and only if $\int \sigma \D \mu \neq 0$, and~\eqref{eq:CNS_EM-vpelec} holds. But integrating the second equality with respect to\ $\mu$, we find the Penrose condition~\eqref{eq:penrose-multi-VPelec}. We conclude that $L_n + K_n$ admits $\lambda$ as an eigenvalue if and only if~\eqref{eq:penrose-multi-VPelec} holds.
	
	On the one hand, we deduce that~\eqref{eq:penrose-multi-VPelec} implies that $\bar \gamma_n>0$. On the other hand, we observe that for all $n \in \Z^d$ and $\lambda \in \C$ with $\Re\lambda>0$ such that~\eqref{eq:penrose-multi-VPelec} holds,
    \begin{equation*}
    1 = \left| \int \frac{1}{(\lambda + in\cdot \mathcal V(\alpha))^2} \D \mu(\alpha) \right| \leq \int \frac{1}{|\lambda + i n \cdot \mathcal V(\alpha)|^2}\D \mu(\alpha) \leq \frac{1}{(\Re\lambda)^2}.
    \end{equation*}
    Therefore, $\bar \gamma_n \leq 1$, and in particular, $\bar \gamma_0 < + \infty$. 

    Finally, if the Penrose condition holds for some $(n_0,\lambda_0)$, choosing $\gamma_0 \in (\bar \gamma_0, 2 \bar \gamma_0)$, there exists $n \in \Z^d$ such that $\bar \gamma_n \geq \gamma_0/2$, and therefore, there is a corresponding growing mode of frequency $n$ and eigenvalue $\lambda$ with $\Re\lambda \geq \gamma_0 / 2$. We conclude that Assumption~\ref{item:ass_insta_LINEARIZED} holds in this case. Note that with $X(\alpha)=(\sigma(\alpha), \eta(\alpha))$ defined as in \eqref{eq:CNS_EM-vpelec}, in view of the definition $P_{1,n}(X)(\alpha)=G_1(X(\alpha) e^{i n \cdot x}) e^{-i n \cdot x} $ and thanks to the Penrose condition~\eqref{eq:penrose-multi-VPelec}, we have
    \begin{align*}
        P_{1,n}(X)(\alpha)=\frac{in}{\vert n \vert^2}\int \sigma \, \mathrm{d}\mu \neq 0,
    \end{align*}
 hence $P_{1,n}(X)(\alpha)$ is independent of $\alpha$.
\end{proof}

\subsection{Proof of Proposition~\ref{prop-check-assump-Instab} in the ion case}

\begin{proof}[Proof of Proposition \ref{prop-check-assump-Instab} in the ion case] Let us focus on the three main assumptions.

\medskip

\noindent

\underline{$\star$ \ref{item:ass_insta_expansion} holds}: The force field still does not depend on $\alpha$ but its structure is more involved than in the case for electrons because of the exponential term. For $N \geq 1$ being a fixed integer, we are looking for an expansion of the potential $\Phi$ of the form 
\begin{align}\label{eq:expansion-Potential-ions}
    \Phi=\Phi[\rr]= \sum_{i=1}^N \Phi_i[r]+\mathcal{R}_N,
\end{align}
where $\Phi_i$ is $i$-linear in $\rr$ and $\|R_N \|_{H^\ell_x} \lesssim \| \rr \|_{L^{\upinfty}_\alpha H^\ell_x}^{N+1}$  for $\ell \in \N$ large enough, and so that 
\begin{align*}
    G[\rr, \ww]=-\nabla \Phi=  \sum_{i=1}^N G_i[\rr]+R_N 
\end{align*}
satisfies the desired assumption with $G_i[\rr]=-\nabla \Phi_i$ and $R_N=-\nabla \mathcal{R}_N$. Assuming such an expansion, we can expand $e^\Phi$ as a sum of $i$-linear quantities in $\rr$ plus a remainder:
    \begin{align*}
	e^\Phi  &=1 +  \Phi_1 +  \Big( \Phi_2 + \frac{1}{2} \Phi_1^2 \Big) +  \Big( \Phi_3 + \Phi_1 \Phi_2 + \Phi_1^3 \Big) \\
    & \quad + \dots + \Big( \Phi_N + \Pi_N(\Phi_1, \dots, \Phi_{N-1}) \Big) + S_N \\
    &=\vcentcolon\mathcal E_N(\Phi_1, \dots, \Phi_N)+S_N,
	\end{align*} 
	where $\Pi_N$ is a polynomial where each term is a given coefficient times a product of factors of type $\Phi_i$ for which the sum of labels is $N$, and $\|S_N \|_{H^\ell_x} \lesssim \| \rr \|_{L^{\upinfty}_\alpha H^\ell_x}^{N+1}$.
More precisely, using the multinomial theorem, we have for all $n \in \N$
\begin{align*}
    \left( \sum_{i=1}^N \Phi_i+\mathcal{R}_N\right)^n&=\sum_{k_1+\cdots k_N+k_{N+1}=n}{C}^n_{\textbf{k}_{[N+1]}} 
 \prod_{i=1}^N \Phi_i^{k_i} \mathcal{R}_N^{k_{N+1}} \\&=\sum_{p=1}^N\underbrace{\sum_{\substack{k_1+\cdots +k_N=n \\
    \sum_{i=1}^N i k_i =p}} {C}^n_{\textbf{k}_{[N]}}
    \prod_{i=1}^N \Phi_i^{k_i}}_{\vcentcolon=\text{Main}_{n,p,N}}+ S_{n,N},
\end{align*}
where, if $\textbf{p}_{[m]}=(p_1, \cdots, p_m)$, the coefficient ${C}^n_{\textbf{p}_{[m]}}=\frac{n!}{p_1!\cdots p_m!}$ is the standard multinomial coefficient. One can observe that when $p \geq 2$ and $2 \leq n \leq p$, one has $\text{Main}_{n,p,N}=\text{Main}_{n,p,N}(\Phi_1, \cdots, \Phi_{p-1})$ and we obtain the decomposition 
\begin{align*}
    e^{\Phi}&=1+\sum_{n=1}^{\infty}\frac{1}{n!}\left( \sum_{i=1}^N \Phi_i+\mathcal{R}_N\right)^n\\
    &=1+\sum_{i=1}^N \Phi_i+\mathcal{R}_N +\sum_{n=2}^{\infty}\frac{1}{n!}\left(\sum_{p=1}^N\text{Main}_{n,p,N}(\Phi_1, \cdots, \Phi_N)+S_{n,N}\right)\\
    &=1+\sum_{p=1}^N\left(\Phi_p+\sum_{n=2}^p \frac{1}{n!}\text{Main}_{n,p,N}(\Phi_1, \cdots, \Phi_{p-1})\right) + S_{N}.
\end{align*}
In other words, with the previous notation, we have $$\Pi_p= \sum_{n=2}^p \frac{1}{n!}\text{Main}_{n,p,N}(\Phi_1, \cdots, \Phi_{p-1}).$$ 
Identifying the $i$-linear term in the equation $e^\Phi-1-\Delta \Phi=\int r^\alpha$, we find the following hierarchy of equations that we can solve iteratively in $H^\ell_x$ for $\ell>3/2$:
	\begin{align*}
	\Phi_1 - \Delta \Phi_1 &= \int r^\alpha \, \mathrm{d}\mu(\alpha), \\
	\Phi_p-\Delta \Phi_p&=- \sum_{n=2}^p \frac{1}{n!}\text{Main}_{n,p,N}(\Phi_1, \cdots, \Phi_{p-1}), \ \ 2 \leq p \leq N,
	\end{align*}
so that $\Vert \Phi_i \Vert_{H^\ell_x} \lesssim \Vert \rr \Vert_{L^{\upinfty}_\alpha H^\ell_x}^i$ for all $1 \leq i \leq N$. We can now conclude if we show that the remainder $\mathcal{R}_N$ (recall the expansion 
 \eqref{eq:expansion-Potential-ions}) satisfies $\Vert \mathcal{R}_N \Vert_{H^\ell_x} \lesssim \Vert \rr \Vert_{L^{\upinfty}_\alpha H^\ell_x}^{N+1}$: this can be done by solving the following equation on $\mathcal{R}_N $, obtained by expressing $S_N$ as $e^U - \mathcal E_N$:
	\begin{equation*}
	    - \Delta \mathcal{R}_N =  - \Big( \exp(\Phi_1 + \dots + \Phi_N + \mathcal{R}_N) - \mathcal E_N(\Phi_1, \dots, \Phi_N)\Big),
	\end{equation*}
that is, calling $\Phi_{\leq N}\vcentcolon=\sum_{i=1}^N \Phi_i$,
\begin{equation}\label{eq:closing}
	   e^{\Phi_{\leq N}}(e^{\mathcal{R}_N}-1) - \Delta \mathcal{R}_N = F_N\vcentcolon=\mathcal E_N(\Phi_1, \dots, \Phi_N)-e^{\Phi_{\leq N}}. 
	\end{equation}
We observe that, by construction, the source term $F_N$ satisfies $\Vert F_N \Vert_{H^\ell_x} \lesssim \Vert \rr \Vert_{L^{\upinfty}_\alpha H^\ell_x}^{N+1}$. To verify Assumption \ref{item:ass_insta_expansion}, we may assume that $\Vert \rr \Vert_{L^{\upinfty}_\alpha H^\ell_x} \leq \delta $ with $\delta$ small enough, so that   $\Vert F_N \Vert_{H^\ell_x} \lesssim \Vert \rr \Vert_{L^{\upinfty}_\alpha H^\ell_x}^{N+1} \leq \delta^{N+1}$. This paves the way for a perturbative approach through a semilinear fixed-point procedure at high regularity, that we detail here for the sake of completeness: we rewrite the equation as (dropping the $N$)
\begin{equation*}
	\mathcal{R}- \Delta \mathcal{R} = f+\mathcal{R}-g(e^{\mathcal{R}}-1), \ \ f=O(\delta^{N+1}), \ \ g(x)=e^{\Phi_{\leq N}(x)}.
	\end{equation*}
By introducing $Q(R)=-(e^R-R-1)-(1-g)(e^R-1)$, we can even rewrite it as	
\begin{equation*}
	\mathcal{R}- \Delta \mathcal{R} = f+Q(\mathcal{R}), \ \ f=O(\delta^{N+1}).
	\end{equation*}	
We consider the mapping $\mathcal{T}$ from $B_f=\lbrace R \in H^\ell_x \mid \Vert R \Vert_{H^\ell_x} \leq 2 	\Vert f \Vert_{H^\ell_x} \rbrace$ to $H^\ell_x$ where for all $R \in B_f$, $\mathcal{T}(R)$ is unique solution to
\begin{align*}
(I- \Delta) \mathcal{T}(R) = f+Q(R).
\end{align*}
The following properties hold:
\begin{itemize}
\item \underline{the map $\mathcal{T}$ is well-defined and sends $B_f$ to itself  for $\delta$ small enough}. By inverting $\mathrm{I}-\Delta$, it is enough to estimate $\Vert Q(R) \Vert_{H^\ell_x}$ for $R \in B_f$. By the composition estimates of Proposition~\ref{prop-Sobolev}, we have for some integer $p$ that may change from line to line
\begin{align*}
&\Vert Q(R) \Vert_{H^\ell_x} \\
&\lesssim (1+\Vert R \Vert_{L^{\upinfty}_x})^p \Vert R \Vert_{H^\ell_x}^2+\Vert 1-e^{\Phi_{\leq N}} \Vert_{H^\ell_x}\Vert e^R-1 \Vert_{H^\ell_x} \\
& \lesssim (1+\Vert R \Vert_{L^{\upinfty}_x})^{p} \left(  \Vert R \Vert_{H^\ell_x}^2+\Vert 1-e^{\Phi_{\leq N}} \Vert_{H^\ell_x}\Vert R \Vert_{H^\ell_x} \right) \\
& \lesssim (1+\Vert R \Vert_{L^{\upinfty}_x})^{p} \left(  \Vert R \Vert_{H^\ell_x}^2+(1+ \Vert \Phi_{\leq N} \Vert_{L^{\upinfty}_x})^p \Vert \Phi_{\leq N} \Vert_{H^\ell_x}\Vert R \Vert_{H^\ell_x} \right).
\end{align*}
Using that by construction $\Vert \Phi_{\leq N} \Vert_{H^\ell_x} =O(\delta)$, Sobolev embedding and the fact that $\Vert f \Vert_{H^\ell_x} \leq 1$ for $\delta$ small enough, there exists $C>0$ such that for all $R \in B_f$
\begin{align*}
\Vert f+Q(R) \Vert_{H^\ell_x} \leq  \Vert f \Vert_{H^\ell_x}+C \delta\Vert f \Vert_{H^\ell_x},
\end{align*}
and therefore, reducing $\delta$ if necessary, this proves that $\mathcal{T}$ is well-defined and sends $B_f$ to itself. By elliptic regularity,  the same type of bound holds for $\mathcal{T}(R)$ in $H^{k+2}_x$.
\item \underline{the map $\mathcal{T}$ is a contraction on $B_f$}. We proceed similarly: if $v=\mathcal{T}(R)$ and $v=\mathcal{T}(R')$ with $R,R' \in B_f$, then we obtain
\begin{align*}
(I-\Delta)(v-v')=j(R')-j(R)+(e^{\Phi_{\leq N}}-1)(e^{R'}-e^R),
\end{align*}
where $j(y)=e^y-y-1$. By \cite[Corollary 2.91]{BCD} (that slightly refines the composition estimates of Proposition \ref{prop-Sobolev}), we have
\begin{align*}
\Vert j(R')-j(R) \Vert_{H^\ell_x} \lesssim (\Vert R \Vert_{L^{\upinfty}_x}+\Vert R' \Vert_{L^{\upinfty}_x}) \Vert R'-R \Vert_{H^\ell_x},
\end{align*}
hence there exists $C>0$ such that
\begin{align*}
\Vert j(R')-j(R) \Vert_{H^\ell_x} \leq C \Vert f \Vert_{H^\ell_x}\Vert R'-R \Vert_{H^\ell_x}.
\end{align*}
Likewise, we also have
\begin{align*}
\Vert (e^{\Phi_{\leq N}}-1)(e^{R'}-e^R)\Vert_{H^\ell_x} \leq C \Vert \Phi_{\leq N}  \Vert_{H^\ell_x} \Vert R'-R \Vert_{H^\ell_x}.
\end{align*}
Using $\Vert \Phi_{\leq N} \Vert_{H^\ell_x} =O(\delta)$ and $\Vert f \Vert_{H^\ell_x} =O(\delta^{N+1})$, we conclude that $\mathcal{T}$ is a contraction on $B_f$ provided that $\delta$ is taken small enough
\end{itemize}
By the Banach fixed-point theorem, this proves that $\mathcal{T}$ admits a unique fixed point $\mathcal{R} \in B_f$ for $\delta$ small enough, and we obtain the desired conclusion by elliptic regularity.

\medskip

\noindent

\underline{$\star$ \ref{item:ass_insta_StabHIGHREG} holds}:  We consider $\Phi_i$ for $i=1,2$, the solution to 
\begin{align*}
    e^{\Phi_i}-1- \Delta \Phi_i = \int r^\alpha_i \D \mu(\alpha),
\end{align*}
where the $\rr_i$ are given densities. Our goal is to obtain an estimate of the form  
\begin{align*}
    \Vert \nabla \Phi_2-\nabla \Phi_1 \Vert_{H^{\ell+1}_x} \leq \mathcal{G}_{\ell}\left(\max(\Vert \rr_1 \Vert_{L^{\upinfty}_\alpha H^{\ell}_x},\Vert \rr_2 \Vert_{L^{\upinfty}_\alpha H^{\ell}_x} )\right) \Vert  \rr_2-\rr_1 \Vert_{L^{\upinfty}_\alpha H^{\ell}_x},
\end{align*}
for some smooth function $\mathcal{G}_{\ell}$ and for $\ell$ large enough. In what follows, we will actually denote by $\mathcal{G}_\ell$ any smooth increasing function depending on $\ell$, that is allowed to change from one line to another.

Set $V=\Phi_2-\Phi_1$ and $m=\int (r^\alpha_2-r^\alpha_1) \D \mu(\alpha)$: we have
\begin{align*}
    e^{\Phi_2}-e^{\Phi_1}-\Delta V=m,
\end{align*}
and by Taylor formula we observe that for any $x$
\begin{align*}
     e^{\Phi_2(x)}-e^{\Phi_1(x)}=a(x) V(x), \ \ a(x)=\int_0^1 e^{\Phi_1(x)+t(\Phi_2(x)-\Phi_1(x))} \, \mathrm{d}t,
\end{align*}
hence we obtain the equation  
\begin{align*}
    -\Delta V+a(x) V=m.
\end{align*}
By standard elliptic regularity, we have
\begin{align*}
\Vert \nabla V \Vert_{H^{\ell+1}_x} \leq \mathcal{G}_\ell(\Vert a \Vert_{H^\ell_x}, C^{-1}) \Vert m \Vert_{H^\ell_x} \lesssim \mathcal{G}_\ell(\Vert a \Vert_{H^\ell_x}, C^{-1}) \Vert  \rr_2-\rr_1 \Vert_{L^{\upinfty}_\alpha H^{\ell}_x},
\end{align*}
where $c=\inf \vert a \vert>0$. More precisely, we have the bound
\begin{align*}
\Vert \Phi_1+t(\Phi_2-\Phi_1) \Vert_{L^{\upinfty}_x} \lesssim \Vert \Phi_1+t(\Phi_2-\Phi_1) \Vert_{H^\ell_x} \lesssim M_\ell,
\end{align*}
where $M_\ell\vcentcolon = \max( \Vert \Phi_1 \Vert_{H^\ell_x}, \Vert \Phi_2 \Vert_{H^\ell_x})$.
Hence, there exists $C>0$ such that $e^{\Phi_1(x)+t(\Phi_2(x)-\Phi_1(x))} \geq e^{-C M_\ell}$ and thus
\begin{align*}
c=\inf \vert a \vert \geq e^{-C M_\ell}.
\end{align*}
By analyticity of the exponential function (or composition rule in Sobolev spaces -- see Proposition \ref{prop-Sobolev}), we also have
\begin{align*}
\Vert a \Vert_{H^\ell_x}  \lesssim \mathcal{G}_\ell(M_\ell).
\end{align*}
All in all, we have obtained
\begin{align*}
\Vert \nabla V \Vert_{H^{\ell+1}_x} \leq \mathcal{G}_\ell(M_\ell) \Vert  \rr_2-\rr_1 \Vert_{L^{\upinfty}_\alpha H^{\ell}_x}.
\end{align*}
To conclude, it therefore remains to prove that 
\begin{align*}
M_\ell \lesssim  \mathcal{G}_\ell\left(\max( \Vert \rr_1 \Vert_{ L^\infty_\alpha H^{\ell}_x},\Vert \rr_2 \Vert_{ L^\infty_\alpha H^{\ell}_x} )\right).
\end{align*}
This follows from the high regularity estimate holding for each potential $\nabla \Phi$ (see Section \ref{Subsec:VPions}; the addition of a constant $1$ in the elliptic equation does not change anything to the argument) we have indeed proved that $\Vert \Phi_i \Vert_{H^2_x}+ \Vert \nabla \Phi_i \Vert_{H^{\ell+1}_x} \lesssim \mathcal{G}_\ell(\Vert \rr_i \Vert_{ L^\infty_\alpha H^{\ell}_x})$. This is enough to conclude the proof.
\medskip

\noindent

\underline{$\star$ \ref{item:ass_insta_LINEARIZED} holds}: The associated linearized system is now
\begin{equation}
\label{eq:VPlin-ions}
\left\{  
\begin{gathered}
\partial_t r^\alpha + \mathcal{V}(\alpha) \cdot \nabla r^\alpha + \Div(w^\alpha)  = 0,\\
\partial_t w^\alpha + (\mathcal{V}(\alpha) \cdot \nabla) w^\alpha   = - \nabla \Phi,\\
\Phi- \Delta \Phi = \int r^\alpha \D \mu(\alpha).
\end{gathered}
\right.
\end{equation}
\textit{Modulo} the minor change in the last equation, the analysis of \eqref{eq:VPlin-ions} is similar to that of the linearized case for the electrons, that is system \eqref{eq:VPlin-elec}. By performing the exact same computations as in the proof of \ref{item:ass_insta_LINEARIZED}(a) and \ref{item:ass_insta_LINEARIZED}(b) in Proposition \ref{prop-check-assump-Instab}, and working under the associated  Penrose condition~\eqref{eq:penrose-multi-ions}, the assumption \ref{item:ass_insta_LINEARIZED} holds as well.
\end{proof}

	\chapter{Application to Vlasov--Navier--Stokes type  systems}\label{Part3-VNS}

This chapter is dedicated to the study of important prototypes of fluid-particle models, namely the incompressible or compressible Vlasov--Navier--Stokes system and the Vlasov--Stokes system. These equations are introduced in Section \ref{Section-PresentationVNSandco}).
Again, we build 
through the point of view of the multiphasic framework developed in  Chapter~\ref{Part1-LWP}. 

This chapter is divided in two main parts:

\begin{itemize}
     \item  The first part focuses on {\bf local well-posedness} results. As in Chapter~\ref{Part2-VP} for Vlasov--Poisson type systems, we apply the abstract local well-posedness theory from Chapter~\ref{Part1-LWP}. We justify that the various Vlasov--Navier--Stokes type systems satisfy the abstract assumptions, thus obtaining several original local well-posedness results for their multiphase and kinetic formulations, which are respectively stated in Section \ref{SubsecVNSmultiphaseLWP}--\ref{SubsecVNScomp-multiphaseLWP}--\ref{subsec:Stokes}.

\item In the second part, we study the dynamics of the incompressible or compressible Vlasov--Navier--Stokes system on the periodic torus $\T^3$, near constant monokinetic profiles. Our goal is to establish a global existence result for these system in the multiphasic context, together with {\bf asymptotic stability} properties of such equilibria. Our approach offers a new point of view on the question, unveiling a simple mechanism of alignment for all the phases,  that ultimately leads to concentration in velocity for the kinetic distribution function. The main results are contained in Section \ref{Section:monokinetic-stab-theorem}, while the proofs are displayed through Sections \ref{sectionVNS-Conslaws-Decay}--\ref{Subsection-preuveStab-VNS}--\ref{Section:VNS-comp}--\ref{Section-nonlindrag}.

\end{itemize}

\section{Presentation of the models}\label{Section-PresentationVNSandco}

We first introduce the three main fluid-particle couplings that we will consider in this chapter. To simplify the presentation, and because this is the most relevant physical case, we only deal with the space dimension $d=3$. 
All the models studied in this chapter involves an explicit coupling between the fluid equation and the kinetic equation through  a \textit{drag force}, that is proportional to the relative velocity between the fluid and particles. Strictly speaking, this is only physically meaningful under that form in the tridimensional case,  and is  reminiscent  of the so-called Stokes law used for the computation of the drag force exerted by a viscous fluid on a spherical object -- see e.g.\cite{oro, desvillettes2010-model,DGR-deriv,HoferPHD}). 

For what concerns long-time properties of such systems, we refer to the upcoming Section \ref{Section:monokinetic-stab-theorem}.

\medskip

\paragraph{Incompressible Vlasov--Navier--Stokes system.}
We start with the Vlasov--Navier--Stokes equations \eqref{eq:VNS_real}, which was one of the main motivations of this work. Recall that it reads
\begin{equation}
		\label{eq:VNS_kinLWP}
		\left\{
		\begin{aligned}
		\partial_t f + v \cdot \nabla_x f +\Div_v\big( (U-v)f\big) &= 0,\\
		\partial_t U + (U \cdot \nabla_x) U + \nabla_x P - \Delta_x U &= \int_{\R^3} (v-U) f \, \mathrm{d}v ,\\
		\Div_x(U) &= 0,
		\end{aligned}
		\right.
		\end{equation}
where the fluid velocity is $U(t,x) \in \R^3$, the particle distribution function is $f(t,x,x) \in \R^+$ and $P(t,x) \in \R$ stands for the fluid pressure associated to the incompressibility condition.

The Cauchy theory for \eqref{eq:VNS_kinLWP} is by now well understood. Note that this system intrinsically contains the sole incompressible Navier--Stokes equations (take $(0,U)$ as a solution), and it is therefore not expected that one can obtain a better theory for the fluid part compared to the standard literature on the Navier--Stokes system. In the context of global weak solutions, one can build solutions $(f,U)$ where $f$ is a weak solution to the Vlasov equation, seen as a transport equation in the phase space, while the fluid velocity $U$ is a solution to the Navier--Stokes equations in the sense of Leray. This is the point of view taken in \cite{ABdM,BDGM,BGM}, for finite kinetic energy data.

With a slightly more regular initial velocity field, and under some suitable smallness assumptions, we  refer to \cite{HKMM,danchin2024fujita,HKM3-3D} where more regularity on the velocity field can be propagated globally in time, together with decay in time estimates (see also Section \ref{Section:monokinetic-stab-theorem} below). For what concerns uniqueness, we refer to \cite{HKM3, danchin2024fujita, HKM3-3D}.

In a different direction, one can also construct local in time solutions to \eqref{eq:VNS_kinLWP} in the realm of high order Sobolev spaces. This follows by standard energy estimates, which are available since there is no loss of derivatives in the kinetic equation and since the Navier--Stokes equations display some parabolic effects. We refer to the strategy devised in \cite{HKapde} and on the estimates from \cite{EHK-thick} (which deals with a more singular coupling). We also refer to the regularity estimates from \cite{Dechicha}.

Let us finally summarize the previous discussion  by stating a typical existence result in the context of weak finite energy solutions \textit{versus} strong solutions on $\T^3 \times \R^3$. 
\begin{Thm}\label{thm:vns-classical}  The following existence results hold for the incompressible Vlasov--Navier--Stokes system.
\begin{itemize}
    \item (Weak solutions, from \cite{BDGM}) Let $U_0 \in L^2_x$ with $\mathrm{div}_x(U_0)=0$. Let $f_0 \in L^1_{x,v} \cap L^{\infty}_{x,v}$ such that $f_0 \geq 0$ and $\int_{\T^d \times \R^d} f^{\mathrm{in}} \vert v \vert^2 \, \mathrm{d}x \, \mathrm{d}v<\infty$. There exists a global in time weak solution $(f,U)$ to \eqref{eq:VNS_kinLWP} such that
    \begin{align*}
        f \in L^\infty(\R^+;L^1_{x,v} \cap L^{\infty}_{x,v}), \ \ U \in L^\infty(\R^+, L^2_x) \cap L^2_x(\R^+; H^1_x),
    \end{align*}
    with initial data $(f_0,u_0)$, which satisfies for all $t>0$
\begin{multline*}
 \int\vert U(t) \vert^2 \, \mathrm{d}x+ \int f(t) \vert v \vert^2 \, \mathrm{d}x \, \mathrm{d}v \\
+ 2\int_0^t\left(\int \vert \nabla_x U(\tau) \vert^2 \, \mathrm{d}x+ \int f(\tau) \vert v-U(\tau) \vert^2 \, \mathrm{d}x \, \mathrm{d}v \right) \mathrm{d}\tau \\ \leq  \int\vert U_0 \vert^2 \, \mathrm{d}x+ \int f_0 \vert v \vert^2 \, \mathrm{d}x \, \mathrm{d}v.
\end{multline*}
\item (Strong solutions, from \cite{HKapde}) Let $k \in \N$ with $k>1+3/2$ and $r>0$ large enough. Let $f_0 \in \mathrm{H}^{k-1}_r$ (recall Definition \ref{def:weighted-v-sobolev} of weighted Sobolev spaces) and $U_0 \in H^{k}_x$ with $\mathrm{div}_x(U_0)=0$.
    Then there exist $T>0$ and unique solution $(f,U)$ to \eqref{eq:VNS_kinLWP} with initial data $(f_0,U_0)$ such that $$ f \in \mathscr{C}([0,T]; \mathrm{H}^{k-1}_r), \ \ U \in \mathscr{C}([0,T]; H^{k}_x) \cap L^2([0,T]; H^{k+1}_x).$$
\end{itemize}
\end{Thm}

\medskip

In the multiphasic framework, the multiphasic formulation of \eqref{eq:VNS_kinLWP} is the following:
\begin{equation}
\label{eq:VNS-multiphase}
\left\{ 
\begin{aligned}
\partial_t \rho^\alpha + \Div ( \rho^\alpha v^\alpha) &= 0,\\
\partial_t v^\alpha + (v^\alpha \cdot \nabla) v^\alpha &= U - v^\alpha,\\
\partial_t U + (U\cdot \nabla ) U + \nabla P -  \Delta U &= \int \big( v^\alpha - U ) \rho^\alpha \D \mu(\alpha)
,\\
\Div(U) &= 0,
\end{aligned}
\right.
\end{equation}
where $(I, \mu)$ is a fixed probability space and where  the multiphasic densities and velocities are respectively  $\rho^\alpha(t,x) \in  \R$ and $v^\alpha(t,x) \in \R^3$. 

As mentioned earlier in the introduction, the case of only one phase corresponds to the so-called pressureless Euler/Navier--Stokes system extensively studied in \cite{choi2024revisit, ChoiJungR3,HuangTangZouR3, lemarié2025, zhai, Danchin-ENS2026}. Following the same framework as in this monograph, \cite{LemarieMultiphaseVNS} also considered the multiphasic version \eqref{eq:VNS-multiphase}: there, the focus is on the whole space $\R^d$ in the context of Besov spaces; the main result is the construction of a global in time solution around the null solution within a critical
regularity framework.

\paragraph{Compressible Vlasov--Navier--Stokes system.}
We also consider another variant of the Vlasov--Navier--Stokes system, where the surrounding fluid is assumed to be compressible and barotropic. The fluid-kinetic system on the fluid density $\varrho(t,x) \in \R^+$, fluid velocity $U(t,x) \in \R^3$ and the distribution function $f(t,x,v) \in \R^+$ can be written as
\begin{equation}\label{eq:VNS-comp_kinLWP}
\left\{ 
\begin{aligned}
\partial_t f + v \cdot \nabla_x f +\Div_v\big( (U-v)f\big) &= 0, \\
\partial_t \varrho + \Div_x ( \varrho U) &= 0,\\
\varrho\big(\partial_t U + (U\cdot \nabla_x ) U \big) + \nabla_x P(\varrho) -  \mathcal{A}U &= \int_{\R^3} (v-U) f \, \mathrm{d}v, \\
\end{aligned}
\right.
\end{equation}
where 
\begin{align}\label{def-Laméoperator}
    \mathcal{A}u\vcentcolon = \mu \Delta_x u+ (\lambda+\mu) \nabla_x \mathrm{div}_x u
\end{align}
is the so-called Lamé operator. Here we consider constant coefficients $\mu>0$ and $\lambda \in \R$ so that $\lambda+2\mu>0$. Above, the pressure $P:\R^+ \rightarrow \R $  is a given smooth function. An idealised case study is for instance the case of isentropic law, where $P(\varrho) \sim \rho^{\gamma}$ for some $\gamma>1$.

 Compared to the incompressible case, the Cauchy theory for \eqref{eq:VNS-comp_kinLWP} remains less developed. At the level of weak solutions, the only available result is the one dimensional global existence and uniqueness result from \cite{LiShou}. In the physically relevant case $d=3$, the theory is currently limited to local in time existence and uniqueness of regular solutions; see for instance \cite{ChoiJung-comp,EHK-thick}. The following theorem summarizes the main available results, stated here in the periodic setting for the spatial variable. 
\begin{Thm}
    The following existence results hold for the compressible Vlasov--Navier--Stokes system.
    \begin{itemize}
        \item (Weak solutions, 1d case, from \cite{LiShou}) Let $\varrho_0 \in W^{1, \infty}_x$ with $\inf \varrho_0>0$, $U_0 \in H^1_x$ and $f_0 \in L^1_{x,v}\cap L^\infty_{x,v}$ nonnegative with compact support in velocity. Then there exists a unique global in time weak solution $(f, \varrho, U)$ to \eqref{eq:VNS-comp_kinLWP} with initial data $(f_0, \varrho_0, U_0)$ that satisfies
        \begin{align*}
        &f \in L^\infty(\R^+; L^1_{x,v}\cap L^\infty_{x,v}), \\
        & \varrho \in L^\infty(\R^+; W^{1, \infty}_x), \ \ U \in L^\infty(\R^+; H^1_x) \cap L^2(\R^+; H^2_x).
        \end{align*}
        \item (Strong solutions, 3d case, adapted from \cite{HKapde, EHK-thick}) Let $k \in \N$ with $k>1+3/2$ and $r>0$ large enough. Let $f_0 \in \mathrm{H}^{k-1}_r$ (recall Definition \ref{def:weighted-v-sobolev} of weighted Sobolev space). Let $(\varrho_0,U_0) \in H^{k}_x \times H^{k+1}_x$ with $\inf \varrho_0>0$.
    Then there exists $T>0$ and unique solution $(f,\varrho, U)$ to \eqref{eq:VNS-comp_kinLWP} on $[0,T]$ with initial data $(f_0, \varrho_0, U_0)$ such that \begin{align*}
    &f \in \mathscr{C}([0,T]; \mathrm{H}^{k-1}_r), \\
    &\varrho \in \mathscr{C}([0,T]; \mathrm{H}^{k}_x),  \ \ U \in \mathscr{C}([0,T]; \mathrm{H}^{k}_x) \cap L^2([0,T]; H^{k+1}_x).
    \end{align*}
    \end{itemize}
\end{Thm}

Let us now turn to the multiphasic version of the system \eqref{eq:VNS-comp_kinLWP}. It reads as follows:
\begin{equation}
\label{eq:VNScomp-multiphase}
\left\{ 
\begin{aligned}
\partial_t \rho^\alpha + \Div ( \rho^\alpha v^\alpha) &= 0,\\
\partial_t v^\alpha + (v^\alpha \cdot \nabla) v^\alpha &= U - v^\alpha,\\
\partial_t \varrho + \Div ( \varrho U) &= 0,\\
\varrho\big(\partial_t U + (U\cdot \nabla_x ) U \big) + \nabla P(\varrho) -  \mathcal{A}U &= \int \big( v^\alpha - U ) \rho^\alpha \D \mu(\alpha).
\end{aligned}
\right.
\end{equation}
The case of a single phase coupled with the compressible Navier--Stokes equations has for instance been studied in \cite{GWZ,li2025global,LiShouZhang}.
\paragraph{Vlasov--Stokes system.}
We finally consider a common variant of the incompressible Vlasov--Navier--Stokes case (i.e. system \eqref{eq:VNS_kinLWP}), which is the following coupling between a steady Stokes system and a Vlasov equation, that is:
\begin{equation}
		\label{eq:VSkin}
		\left\{
		\begin{aligned}
		\partial_t f + v \cdot \nabla_x f +\Div_v\big( (U-v)f\big) &= 0,\\
		 - \Delta_x U + \nabla_x P&= \int_{\R^3} (v-U) f \, \mathrm{d}v,\\
		\Div_x(U) &= 0,
		\end{aligned}
		\right.
		\end{equation}
set on the phase space $\R^3 \times \R^3$. It typically arises in the case of low Reynolds regime for the surrounding propellant in which small particles evolve. We refer to \cite{hofer2020sedimentation} for more details. For what concerns the Cauchy problem for \eqref{eq:VSkin}, the construction of global weak solutions can be found in \cite{HoferPHD,HS-VSmeanfield2}.  Let us state the following theorem for well-posedness of \eqref{eq:VSkin}. 
\begin{Thm}[from \cite{HS-VSmeanfield2}]
Let $k>9$ and 
$
f^0 \in \mathcal M (\mathbb R^3 \times \mathbb R^3)
\cap L^\infty(\mathbb R^3 \times \mathbb R^3)
$
with $\int \vert v \vert^k f_0 \, \mathrm{d}v \D x<\infty$. There exists a unique global in time weak solution
$(f,U)$ to \eqref{eq:VSkin} with for all $T>0$
\begin{align*}
    f \in L^\infty((0,T)\times \mathbb R^3\times \mathbb R^3)
\cap \mathscr{C}((0,T);w\text{-}\mathcal M(\mathbb R^3\times \mathbb R^3)), \\
U\in L^\infty((0,T);W^{1,\infty}(\mathbb R^3))
\cap \mathscr{C}((0,T)\times \mathbb R^3),
\end{align*}
and
\begin{equation*}
\sup_{t>0} \int \vert v \vert^k f(t) \, \mathrm{d}v \D x<\infty.
\end{equation*}
\end{Thm}
In the multiphasic framework, the system \eqref{eq:VSkin} is  recast as 
\begin{equation}
\label{eq:VS-multiphase}
\left\{ 
\begin{aligned}
\partial_t \rho^\alpha + \Div ( \rho^\alpha v^\alpha) &= 0,\\
\partial_t v^\alpha + (v^\alpha \cdot \nabla) v^\alpha &= U - v^\alpha,\\
- \Delta U+\nabla P  &= \int \big( v^\alpha - U ) \rho^\alpha \D \mu(\alpha)
,\\
\Div(U) &= 0,
\end{aligned}
\right.
\end{equation}
and is set on $\R^3$. 

\begin{Rem}
As we will explain later on, we will also be able to include a smooth nonlinear drag term in the former models (see in particular Section \ref{sec:nonlindrag-LWP}).
\end{Rem}

\sectionhead
  [Local well-posedness for the incompressible VNS equations]
  {Local well-posedness for the incompressible
   Vlasov--Navier--Stokes equations}\label{SubsecVNSmultiphaseLWP}

We first focus on the multiphasic Vlasov--Navier--Stokes system \eqref{eq:VNS-multiphase}. In this subsection, we consider $x \in \T^3$ or $x \in \R^3$. From now on, we also assume that the initial data $\rho^\alpha_0$, and hence the associated solution of the continuity equation, is nonnegative.

 As for the Vlasov--Poison type systems, see e.g. Section \ref{Subsec:VPelectrons}, we put forward two kinds of assumptions on the initial data that we have identified in our general framework from Chapter \ref{Part1-LWP}: it either corresponds to Theorem \ref{thm-LWPkinetic}  with the initial multiphasic decomposition from Example \ref{eq:example-decompoLinfty} (for densities $\rhorho$ bounded with respect to the label $\alpha$) or to Theorem \ref{thm-LWPkinetic-L1} with the initial multiphasic decomposition from Example \ref{eq:example-decompoL1} (for densities $\rhorho$ integrable with respect to $\alpha$).

Our first result reads as follows.

\begin{Thm}\label{thm:appliVNS-incomp}
Let $k \in \N$ such that $k>1+3/2$. Let $f_0$ be an initial distribution function written
$$f_0(x, \cdot)=\int_{I} \rho^\alpha_0(x) \otimes  \delta_{v=v^\alpha_0(x)} \, \mathrm{d}\mu(\alpha),$$
for some set of labels $(I, \mu)$ and some regular multiphasic distribution $(\rhorho_0, \vv_0)$. Let $ U_0 \in H^{k-1}_x$ with $\Div (U_0)=0$. 

\begin{enumerate}
    \item If $(\rhorho_0, \vv_0) \in L^\infty_\alpha H^{k-1}_x \times \mathcal{H}^{k,p} $ with $p \in [2, \infty]$ and $\rho^\alpha_0 \geq 0$, there exist $T>0$ and a unique solution $(\rhorho, \vv, U)$ on $[0,T]$ to the multiphasic system \eqref{eq:VNS-multiphase} with $\rho^\alpha \geq 0$ and
    \begin{align*}
    &\rhorho \in L^\infty([0,T]; L^\infty_\alpha H^{k-1}_x ), \ \ \vv \in L^\infty([0,T]; \mathcal{H}^{k,p}), \\
    &U \in L^\infty([0,T]; H^{k-1}_x) \cap L^2([0,T]; H^{k}_x).
\end{align*}
    \item Let $\Theta: I \to \R_+$ and $ \lambda : I \to \R^d$ be measurable maps satisfying $\langle \lambda^\alpha \rangle^p \lesssim \Theta(\alpha)$ with $p \in [1, \infty)$. If $(\rhorho_0, \vv_0= \ww_0 + \boldsymbol{\lambda})$ is such that $(\rhorho_0,\ww_0 ) \in L^1_{\alpha,\Theta} H^{k-1}_x  \times \mathcal{H}^{k,\infty}$  with $\rho^\alpha_0 \geq 0$, there exist $T>0$ and a unique solution $(\rhorho, \vv, U)$ on $[0,T]$ to the multiphasic system \eqref{eq:VNS-multiphase} with $\rho^\alpha \geq 0$ and
    \begin{align*}
    &\rhorho \in L^\infty([0,T]; L^1_{\alpha,\Theta} H^{k-1}_x ), \ \ \vv= e^{-t}(\boldsymbol{\lambda}+ \widetilde\ww),\, \widetilde\ww \in L^\infty([0,T]; \mathcal{H}^{k,\infty}), \\ 
    &U \in L^\infty([0,T]; H^{k-1}_x) \cap L^2([0,T]; H^{k}_x).
    \end{align*}
\end{enumerate}
In both cases, if the maximal time of existence $T^\star$ is finite, then
\begin{equation}
\int_0^{T^\star}  \left(  \| \nabla \vv(s) \|_{L^{\upinfty}_\alpha L^{\upinfty}_x} +  \| \nabla U(s) \|_{L^{\upinfty}_x}    \right) \D s = +\infty.
\end{equation}
Moreover, defining
$$
f(t,x,\cdot) \vcentcolon= \int_{I} \rho^\alpha(t,x) \otimes  \delta_{v=v^\alpha(t,x)} \, \mathrm{d}\mu(\alpha),
$$
then $(f,U)$ is the unique weak solution on $[0,T]$ of  \eqref{eq:VNS_kinLWP} with initial data $(f_0,U_0)$, 
in the sense of Definition \ref{def:solution_coupled}. The distribution function $f$ is continuous in time, $H^{k-1}_x$--strongly in space and weakly in velocity in the sense of Definition~\ref{def:solution_Vlasov+timecontinuity},  the fluid velocity $U$ satisfies
\begin{align*}
    &U \in \mathscr{C}^0([0,T]; H^{k-1}_x),
\end{align*}
and  higher order moments in velocity enjoy the following regularity:
\begin{align*}
  &(t,x) \mapsto \int_{\R^d} \vert v \vert^j f(t,x,\D v) \in L^\infty([0,T]; H^{k-1}_x), \ \ j=0, \ldots,\lfloor p\rfloor. 
\end{align*}
Finally, if $U_0 \in H^k_x$ then the whole statement remains true, with the regularity
\begin{align*}
    U \in \mathscr{C}^0([0,T]; H^{ k}_x) \cap L^2([0,T];H^{k+1}_x).
\end{align*}
\end{Thm}

\begin{Rem}The following remarks are in order.

\begin{itemize}
\item Due to the parabolic nature of the Navier--Stokes equations on $U$, we have some flexibility on the regularity of $U_0$. Roughly speaking, we have the following alternative:
\begin{itemize}
    \item either we put $U_0 \in H^{k-1}_x$, at the same level of regularity as the densities $\rhorho_0$. Since $k-1>3/2$, it corresponds to the standard regularity theory for local in time strong solutions of the Navier--Stokes system, see e.g.\cite{RobinsonRodrigoSadowskiNS3};
    \item or we take $U_0 \in H^{k}_x$, at the same level of regularity as the velocities $\vv_0$ (hence requiring that $k+1>3/2$ in view of the framework from Chapter \ref{Part1-LWP}). Since the densities only appear as a source term in the Navier--Stokes system, the parabolic effect enables a shift of one derivative between $U_0$ and $\vv_0$.
\end{itemize}
    \item Contrary to the Vlasov--Poisson type systems studied in Chapter~\ref{Part2-VP}, the $\dot{H}^{-1}_x$ regularity is not required and is thus dispensed with. In addition, in the first set of assumptions,  the constrained $p \geq 2 $ is specific to the Vlasov--Navier--Stokes case: it comes from a technical obstacle in the forthcoming proof of stability at low regularity. We also note that the source term in the Navier--Stokes equations involves in that case a first-order moment in velocity, that is 
    $$\int \rho^\alpha v^\alpha \, \mathrm{d}\mu(\alpha).$$
    In particular, it is well defined under the first set of assumptions since $(\rhorho(t), \vv(t)) \in L^\infty_\alpha H^{k-1}_x \times \mathcal{H}^{k,p}$. However, under the second set of assumptions, it reads
    \begin{align*}
        \int \rho^\alpha v^\alpha \, \mathrm{d}\mu(\alpha)=\int \rho^\alpha \lambda^\alpha \, \mathrm{d}\mu(\alpha)+\int \rho^\alpha w^\alpha \, \mathrm{d}\mu(\alpha).
    \end{align*}
As explained in Remarks~\ref{rem:Rd-alpha}--\ref{rem:Rd-alpha2} and Example \ref{eq:example-decompoL1}, with $\vv=\boldsymbol{\lambda}+\ww$, we typically have in mind $\lambda^\alpha=\alpha$ when it comes to finding a suitable decomposition of the kinetic initial data in the full space case. As a consequence, we enforce the integrability condition $L^1_{\alpha,\Theta}$ for the density $\rhorho$, with $\vert \lambda^\alpha\vert \leq \langle \lambda^\alpha \rangle^p \lesssim \Theta(\alpha)$, which allows to give a meaning to the first term above.

In particular, under both sets of assumptions, the source term in the Navier--Stokes equations lies in $L^2_{\mathrm{loc}}(\R^+; H^{k-1}_x)$.
\item We also refer to the short paragraph \ref{subsec-jaiplusdidee} below about the treatment of the drag force in that context, which explains the presence of the exponential in time factor for the solution in the second set of assumptions.
    \item  As a consequence of the general theory from Chapter \ref{Part1-LWP} -- Remark \ref{Rem:timeuniform},  we also obtain the fact that $t \mapsto \Vert \vv(t) \Vert_{\mathcal{H}^{k,p}}$ (in the first case) and $t \mapsto \Vert \ww(t) \Vert_{\mathcal{H}^{k,\infty}}$ (in the second case) are continuous in $[0,T]$. We shall also obtain a more precise control that implies the above blow-up criterion, see Proposition \ref{prop-VNSblowup} below.
    
    \item We refer to Section \ref{sec:nonlindrag-LWP} where we include the treatment of a smooth nonlinear drag term in the Vlasov--Navier--Stokes system, at the cost of having to impose $p=\infty$ in the former statement, in the first set of assumptions.  In view of Remark \ref{rem:bounded_velocity}, the outcome is that we can only consider kinetic distributions that are compactly supported in velocity. However, let us highlight that taking nonlinear drag forces into account was not, to our knowledge, yet achieved in the mathematical theory of Vlasov--Navier--Stokes systems.
\end{itemize}

\end{Rem}

 As in the Vlasov--Poisson case (see Section \ref{Subsec:VPelectrons}), we first draw direct consequences of Theorem \ref{thm:appliVNS-incomp} on the Vlasov--Navier--Stokes system. As already seen, it mainly stems from a well-chosen set of initial data that fits into the assumptions of Theorem~\ref{thm:appliVNS-incomp}. The initial multiphasic decompositions are exactly the same as in Section \ref{Subsec:VPelectrons} and we do not repeat the proof. 

\begin{Cor}\label{coro-multiEuler-VNSincomp}
Let $k\in \N$ such that $k>1+3/2$. Let $ U_0 \in H^{k-1}_x$ with $\Div (U_0)=0$  and let $(m^n_0,v^n_0)_{n\in \N}$ be a family of density--velocity fields satisfying the bound
\begin{equation}
    \label{eq:coro1-vns}
 \| \langle n \rangle^\beta m^n_0 \|_{\ell^{\upinfty}_n (H^{k-1}_x)} +\| \langle n \rangle^{-\beta/p} v^n_0 \|_{\ell^p_n L^2_x} + \| \DDD v^n_0 \|_{\ell^{\upinfty}_n H^{k-1}_x} < +\infty,
\end{equation}
for some $\beta>1$ and $p \in [2,\infty)$. There exist $T>0$ and a unique weak solution $(f,U)$ to the Vlasov--Navier--Stokes  system~\eqref{eq:VNS_kinLWP} on $[0,T]$, with initial condition $(f_0, U_0)$, with
$$f_0 = \sum_{n=0}^{+\infty} m^n_0(x) \otimes \delta_{v= v^n_0(x)},$$ 
in the sense of Definition~\ref{def:solution_coupled}.  The distribution function $f$ is continuous in time, $H^{k-1}_x$--strongly in space and weakly in velocity in the sense of Definition~\ref{def:solution_Vlasov+timecontinuity} and  the fluid velocity $U$ satisfies 
\begin{align*}
    &U \in \mathscr{C}^0([0,T]; H^{k-1}_x) \cap L^2([0,T]; H^{k}_x).
\end{align*} Moreover, higher order moments in velocity enjoy the following regularity:
\begin{align*}
  &(t,x) \mapsto \int_{\R^d} \vert v \vert^j f(t,x,\D v) \in L^\infty([0,T]; H^{k-1}_x), \ \ j=0, \ldots,\lfloor p\rfloor.
\end{align*}
Finally, if $U_0 \in H^{k}_x$ then the whole statement remains true, with the regularity
\begin{align*}
    U \in \mathscr{C}^0([0,T]; H^{ k}_x) \cap L^2([0,T];H^{k+1}_x).
\end{align*}
\end{Cor}

\begin{Cor}\label{coro:appliVNS-incomp}
    Let $k\in \N$ such that $k>1+3/2$. Let $f_0 \in  L^1_{v, \langle v \rangle} (H^{k-1}_x\cap \dot H^{-1}_x)$ and 
    $U_0 \in H^{k-1}_x$ with $\mathrm{div}_x(U_0)=0$. There exist $T>0$ and a unique weak solution $(f,U)$ to the Vlasov--Navier--Stokes  system~\eqref{eq:VNS_kinLWP} on $[0,T]$, with initial condition $(f_0,U_0)$ in the sense of Definition \ref{def:solution_coupled}. The distribution function $f$ is continuous in time, $H^{k-1}_x$--strongly in space and weakly in velocity in the sense of Definition~\ref{def:solution_Vlasov+timecontinuity}, and the fluid velocity $U$ satisfies 
\begin{align*}
    &U \in \mathscr{C}^0([0,T]; H^{k-1}_x) \cap L^2([0,T]; H^{k}_x).
\end{align*}
If furthermore we assume that $f_0 \in  L^1_{v, \langle v \rangle^p}H^{k-1}_x $ for some $2 \in [1, \infty)$, then higher order moments in velocity enjoy the following regularity:
\begin{align*}
  (t,x) \mapsto \int_{\R^d} \vert v \vert^j f(t,x,\D v) \in L^\infty([0,T]; H^{k-1}_x), \ \ j=0, \ldots,\lfloor p\rfloor.
\end{align*}
Finally, if $U_0 \in H^{k}_x$ then the whole statement remains true, with the regularity
\begin{align*}
    U \in \mathscr{C}^0([0,T]; H^{ k}_x) \cap L^2([0,T];H^{k+1}_x).
\end{align*}
\end{Cor}

\begin{Cor}\label{coro:appliVNS-incomp-Mixed} 
Let $k\in \N$ such that $k>1+3/2$ and $p \in [2, \infty]$.  Let $(m^n_0,v^n_0)_{n\in \N}$ be a family of density--velocity fields satisfying the bound
$$
 \|  m^n_0 \|_{\ell^{1}_n (H^{k-1}_x \cap \dot{H}^{-1}_x)} +\| v^n_0 - \lambda^n \|_{\ell^{\upinfty}_n H^{k}_x} +  < +\infty,
$$
for some sequence $\lambda: \N \to \R^d$ satisfying $\langle \lambda^n \rangle \lesssim \langle n \rangle^p$, a distribution function  $g_0 \in  L^1_{v, \langle v \rangle^p} (H^{k-1}_x\cap \dot H^{-1}_x)$, and $U_0 \in H^{k-1}_x$ with $\mathrm{div}_x(U_0)=0$. There exist $T>0$ and a unique weak solution $(f,U)$ to the Vlasov--Navier--Stokes  system~\eqref{eq:VNS_kinLWP} on $[0,T]$, with initial condition $(f_0, U_0)$
$$
        f_0= \sum_{n=1}^\infty m^n_0  \otimes \delta_{v=v^n_0} + g_0,
$$
 in the sense of Definition~\ref{def:solution_coupled}.  The distribution function $f$ is continuous in time, $H^{k-1}_x$--strongly in space and weakly in velocity in the sense of Definition~\ref{def:solution_Vlasov+timecontinuity} and  the fluid velocity $U$ satisfies
\begin{align*}
    &U \in \mathscr{C}^0([0,T]; H^{k-1}_x) \cap L^2([0,T]; H^{k}_x).
\end{align*}
Moreover, higher order moments in velocity enjoy the following regularity:
\begin{align*}
  &(t,x) \mapsto \int_{\R^d} \vert v \vert^j f(t,x,\D v) \in L^\infty([0,T]; H^{k-1}_x), \ \ j=0, \ldots,\lfloor p\rfloor.
\end{align*}
Finally, if $U_0 \in H^{k}_x$ then the whole statement remains true, with the regularity
\begin{align*}
    U \in \mathscr{C}^0([0,T]; H^{ k}_x) \cap L^2([0,T];H^{k+1}_x).
\end{align*}
\end{Cor}

We now turn to the proof of Theorem~\ref{thm:appliVNS-incomp}. 
Thanks to Theorem~\ref{thm-LWPmultiphase-densityvelocity} in Chapter~\ref{Part1-LWP},for the first set of assumptions on the initial condition, it is enough to check that the Assumptions~\ref{ass:bound_high_regularity-rho} and~\ref{ass:stability-rho}  are satisfied for the system \eqref{eq:VNS-multiphase}. However, for the second set of assumptions on the initial condition, let us emphasize that, contrary to the Vlasov--Poisson case of Chapter~\ref{Part2-VP}, the proof with  does not follow directly from the abstract framework around Theorem~\ref{thm-LWPmultiphase-densityvelocity-re}, but requires a slight adaptation that we explain now.

\paragraph{An adaptation of Theorem~\ref{thm-LWPmultiphase-densityvelocity-re}.}
\label{subsec-jaiplusdidee}

Theorem~\ref{thm-LWPmultiphase-densityvelocity-re} cannot be applied as is to prove Theorem~\ref{thm:appliVNS-incomp}. This is due to the drag force in the Vlasov equation, which entails that when looking for $\vv(t)$ of the form $\ww(t) + \boldsymbol{\lambda}$, the system satisfied by $(\rhorho,\ww)$ is
\begin{equation*}
\left\{ 
\begin{aligned}
\partial_t \rho^\alpha + \Div ( \rho^\alpha (w^\alpha+\lambda^\alpha)) &= 0,\\
\partial_t w^\alpha + (w^\alpha+\lambda^\alpha) \cdot \nabla w^\alpha &= U - w^\alpha- \lambda^\alpha.\\
\end{aligned}
\right.
\end{equation*}
One therefore cannot expect $L^\infty_\alpha$ estimates for $w^\alpha$ due to the term $-\lambda^\alpha$ in the equation for $w^\alpha$; in other words Assumptions~\ref{ass:bound_high_regularity-rho-re} and~\ref{ass:stability-rho-re} cannot hold in this setting. In order to bypass this issue we introduce
$\widetilde{w}^\alpha \vcentcolon= e^{t}(w^\alpha+\lambda^\alpha)$ and the unknown 
$(\rhorho,\widetilde{\ww})$ satisfies the system
\begin{equation*}
\left\{ 
\begin{aligned}
\partial_t \rho^\alpha + \Div ( \rho^\alpha (e^t \widetilde{w}^\alpha+ (1+e^t)\lambda^\alpha)) &= 0,\\
\partial_t \widetilde{w}^\alpha + (e^t \widetilde{w}^\alpha+ (1+e^t)\lambda^\alpha)) \cdot \nabla \widetilde{w}^\alpha &= U.\\
\end{aligned}
\right.
\end{equation*}
Once under this form, the force is reduced to $U$ and Assumptions~\ref{ass:bound_high_regularity-rho-re} and~\ref{ass:stability-rho-re} can be checked as for Assumptions~\ref{ass:bound_high_regularity-rho} and~\ref{ass:stability-rho}. It follows that an analogue of Theorem~\ref{thm-LWPmultiphase-densityvelocity-re} holds in this setting (take $\Theta(\alpha) \gtrsim \langle \lambda^\alpha\rangle$) for any $p \in [1, \infty)$.

\bigskip

\begin{proof}[Proof of Theorem \ref{thm:appliVNS-incomp}]
We fix $k \in \N$ such that $k>1+3/2$. As explained in Subsection \ref{subsec:abstract-setting}, the force field $F^\alpha[\vv]$ acting on phase $\alpha$ in~\eqref{eq:VNS-multiphase} consists of the difference between $U$, a field that does not depend on the phase $\alpha$ and which is computed with the help of the Navier--Stokes equation, and $v^\alpha$ itself. 
We focus on the first set of assumptions, so that we want to check Assumptions~\ref{ass:bound_high_regularity-rho} and~\ref{ass:stability-rho}; to this end, we just need to give estimates for $U$ (the velocities $(v^\alpha)$ being given). Therefore, we are in the same situation as in the Vlasov--Poisson case, namely, that we want to estimate a field that does not depend on $\alpha$.

Even though it will be implicit in our approach,  we will use the fact that given a family of densities $\rhorho \in L^\infty_{\mathrm{loc}}(\R^+; L_{\alpha}^\infty H^{k-1}_x)$ and $\vv \in L^\infty_{\mathrm{loc}}(\R^+; \mathcal H^{k,p})$, and given a fixed velocity field $\overline{U} \in L^\infty_{\mathrm{loc}}(\R^+;  H^{k-1}_x)$, for $k-1>3/2$, one can locally solve  in $L^\infty([0,T];H^{k-1}_x) \cap L^2([0,T];H^k_x)$ the  Navier--Stokes system with forcing
\begin{align*}
    \partial_t U + (U\cdot \nabla ) U + \nabla P -  \Delta U &= \int \big( v^\alpha - \overline{U} ) \rho^\alpha \D \mu(\alpha)
,\\
\Div(U) &= 0,\\
U|_{t=0} &= U_0.
\end{align*}
and with initial data $U_0 \in H^{k-1}_x$. Indeed, by Lemma \ref{prop-reg-momentsf}, the former right-hand side lies in $L^2_{\mathrm{loc}}(\R^+; H^{k-1}_x)$.
One can then rely on a standard iterative scheme for the incompressible Navier--Stokes system. Note that, if $ U_0 \in H^{k}_x$ and $\overline{U} \in L^\infty_{\mathrm{loc}}(\R^+;  H^{k}_x)$, it is also possible to solve the former in  $L^\infty([0,T];H^{k}_x) \cap L^2([0,T];H^{k+1}_x)$, thanks to the parabolic smoothing of the Navier--Stokes equations.

\medskip

Let us now prove that the two main Assumptions \ref{ass:bound_high_regularity-rho} and \ref{ass:stability-rho} holds for \eqref{eq:VNS-multiphase}.

\paragraph{Bound in high regularity, when $U_0 \in H^{k-1}_x$.}
Let $\vv \in L^1([0,T]; \mathcal{H}^{k,p})$ be a family of velocities. Given a family of densities $\rhorho \in L^\infty([0,T]; L^\infty_\alpha H^{k-1}_x)$ that solve the continuity equations associated with $(v^\alpha)$ (in view of the general Lemma \ref{LM:SobEstimRHO}), we need to derive a $L^1([0,T];H^{k}_x)$ on the velocity field $U$ solution to
\begin{align*}
    \partial_t U + (U\cdot \nabla ) U + \nabla P -  \Delta U &= \int \big( v^\alpha - U ) \rho^\alpha \D \mu(\alpha)
,\\
\Div(U) &= 0.
\end{align*}
In what follows, $K>0$ refers to a constant that may change from line to line. An energy estimate in $L^\infty([0,T]; H^{k-1}_x) \cap L^2([0,T]; H^{k}_x)$ for the Navier--Stokes equations with a source term gives
\begin{multline}\label{eq:estimVNSdebase-neverstop}
\frac{1}{2}\frac{\D}{\D t}\|U\|^2_{H^{k-1}_x} +   \| \DDD U \|^2_{H^{k -1}_x} 
 \\\leq  \vert \langle (U\cdot \nabla ) U, U \rangle_{H^{k-1}_x} \vert  + \|U\|_{H^{k-1}_x} \int \| ( v^\alpha - U  )\rho^{\alpha} \|_{H^{k-1}_x} \D \mu(\alpha).
 \end{multline}
Using the fact that $H^{k-1}_x$ is an algebra for $k-1>3/2$, we get by Young inequality
\begin{align*}
    \vert \langle (U\cdot \nabla ) U, U \rangle_{H^{k-1}_x} \vert  &\leq \Vert (U\cdot \nabla ) U \Vert_{H^{k-1}_x} \Vert U \Vert_{H^{k-1}_x} \\
    &\leq K\Vert  \nabla  U \Vert_{H^{k-1}_x} \Vert U \Vert_{H^{k-1}_x}^2 \\
    &\leq \frac{1}{2} \Vert  \nabla  U \Vert_{H^{k-1}_x}^2+ K \Vert U \Vert_{H^{k-1}_x}^4,
\end{align*}
and therefore,  after absorption of the gradient term in the left-hand side, we obtain 
\begin{multline}\label{eq:estimVNSdebase-neverstop2}
\frac{\D}{\D t}\|U\|^2_{H^{k-1}_x} +   \| \nabla U \|^2_{H^{k -1}_x} 
 \\\leq   K\left(\| U \|^4_{H^{k-1}_x} + \Vert \rhorho \Vert_{L^{\upinfty}_\alpha H^{k-1}_x}\|U\|_{H^{k-1}_x}^2 +  \Vert \rhorho \Vert_{L^{\upinfty}_\alpha H^{k-1}_x} \Vert \vv \Vert_{L^{p}_\alpha H^{k}_x}\ |U\|_{H^{k-1}_x}\right).
 \end{multline}
Let us now set
 \begin{align*}
     y(t)\vcentcolon=1+\|U(t)\|^2_{H^{k-1}_x} +  \int_0^t \| \nabla U(s) \|^2_{H^{k -1}_x} \, \mathrm{d}s.
 \end{align*}
 Since $1 \leq y(t)$ and $\|U(t)\|_{H^{k-1}_x} \leq y(t)^{1/2}$, we deduce from \eqref{eq:estimVNSdebase-neverstop2} that
 \begin{align*}
     y' &\leq K\big(y^2+\Vert \rhorho \Vert_{L^{\upinfty}_\alpha H^{k-1}_x}y +  \Vert \rhorho \Vert_{L^{\upinfty}_\alpha H^{k-1}_x} \Vert \vv \Vert_{L^{p}_\alpha H^{k}_x}\ y^{1/2} \big) \\
     & \leq  K(1+\Vert \rhorho \Vert_{L^{\upinfty}_\alpha H^{k-1}_x})(1+\Vert \vv \Vert_{L^{p}_\alpha H^{k}_x})y^2.
 \end{align*}
Hence, there holds
\begin{align*}
    \frac{\mathrm{d}}{\mathrm{d}t}y^{-1}=-y^{-2} y' \gtrsim -K(1+\Vert \rho^\alpha \Vert_{L^{\upinfty}_\alpha H^{k-1}_x})(1+\Vert \ww\Vert_{L^p_\alpha H^{k-1}_x}),
\end{align*}
so that, after integration in time,  we get for all times $t \in [0,T]$
\begin{multline*}
    y^{-1}(t)-y^{-1}(0) \\
    \geq -K(1+\Vert \rho^\alpha \Vert_{L^\infty([0,T];L^{\upinfty}_\alpha H^{k-1}_x)}\left(T+\int_0^T (1+\Vert \vv(s)\Vert_{L^p_\alpha H^{k-1}_x}) \mathrm{d}s\right).
\end{multline*}
Setting for $M_1, M_2, P \geq 0$ 
$$\Lambda(t,M_1,M_2,P)\vcentcolon=\frac{(1+P^2)}{1-K(1+P^2)(1+M_2)(t+M_1)},$$ one gets
\begin{align*}
   y(t) \leq \Lambda\left(T,\int_0^T \Vert \vv(s)\Vert_{L^p_\alpha H^{k-1}_x} \mathrm{d}s,\Vert \rho^\alpha \Vert_{L^\infty([0,T];L^{\upinfty}_\alpha H^{k-1}_x)}, \|U_0\|_{H^{k-1}_x}\right),
\end{align*}
for time $T$ and $t \in [0,T]$ such  that the former denominator is strictly positive. In particular, we get for those times that
\begin{align*}
    &\|U\|^2_{L^\infty([0,T],H^{k-1}_x)} +  \int_0^T \| \nabla U(s) \|^2_{H^{k-1}_x} \, \mathrm{d}s \\
    & \leq \Lambda\left(T,\int_0^T \Vert \vv(s)\Vert_{L^p_\alpha H^{k-1}_x} \mathrm{d}s,\Vert \rho^\alpha \Vert_{L^\infty([0,T];L^{\upinfty}_\alpha H^{k-1}_x)}, \|U_0\|_{H^{k-1}_x}\right).
\end{align*}
For all $T>0$, we end up with 
\begin{align*}
\int_0^T \| U(s) \|_{H^k_x} \D s &\lesssim \int_0^T \big( \|U(s)\|_{L^2_x} + \| \nabla U(s) \|_{H^{k-1}_x} \big) \D s \\
&\lesssim T\|U\|_{L^\infty([0,T],H^{k-1}_x)} +\sqrt{T} \left( \int_0^T \| \nabla U(s) \|_{H^{k-1}_x}^2 \D s \right)^{1/2}\\
&  \leq \F\left(T,\int_0^T \Vert \vv(\tau) \Vert_{\mathcal{H}^{k,p}} \D \tau, \Vert \rhorho \Vert_{L^\infty([0,T];L^\infty_\alpha H^{k-1}_x)}, \|U_0\|_{H^{k-1}_x} \right).
 \end{align*}
 where 
 \begin{multline*}
    \F(t,M_1,M_2,P) \\\vcentcolon= \left\{ \begin{aligned}
    &2(t + \sqrt{t}) \Lambda(t, M,N, P)^{1/2}, &&\mbox{if }  K(1+M_2)(t+M_1)(1+P)<1, \\
    & + \infty, && \mbox{else}.
    \end{aligned} \right.
\end{multline*}
 This show that  Assumption~\ref{ass:bound_high_regularity-rho} is satisfied for the system \eqref{eq:VNS-multiphase}. It thus concludes the proof for the bound in high regularity for small enough times.

\paragraph{Bound in high regularity, when $U_0 \in H^{k}_x$.} We now explain how to adapt the argument above to derive the bound at high-regularity in that context. We will also justify the higher regularity for $U$ compared to the previous case: this will entail the last part of the statement of Theorem \ref{thm-cinetique-VNScomp}. We still consider the density-velocity pair $\rhorho \in L^\infty([0,T]; L^\infty_\alpha H^{k-1}_x)$, $\vv \in L^2([0,T]; \mathcal{H}^{k,p})$.
The main difference with the former proof for the high-regularity estimate is that we assume some $L^2$--integrability in time for the fixed velocity field $\vv$, and not merely some $L^1$--integrability. This actually enters the slight generalization of Assumption \ref{ass:bound_high_regularity-rho} highlighted in Remark \ref{rem-bound-high-ref-functionQ-withrho}.

As before, we need to provide an estimate for $U$ in $L^1([0,T];H^{k}_x)$. We will actually bound its $L^\infty([0,T];H^{k}_x)$ norm directly, as we allow on more derivative control on $U$. The decisive point is indeed to rely on the parabolic effect for the Navier--Stokes equations, by seeing the densities $\rhorho$ as lower-order in terms of regularity for the source term : instead of \eqref{eq:estimVNSdebase-neverstop}-\eqref{eq:estimVNSdebase-neverstop2}, we perform an energy estimate for the Navier--Stokes equations in $H^{k}_x$ and get 
\begin{align*}
    \frac{\mathrm{d}}{\mathrm{d}t} \Vert U \Vert^2_{H^{k}_x} 
    + \Vert \nabla U \Vert_{H^k_x}^2  
     &\lesssim   \| U\|_{H^k_x}^4+ \Vert U \Vert_{L^2_x}\int \Vert \rho^\alpha (v^\alpha-U) \Vert_{L^2_x} \, \mathrm{d}\mu(\alpha)\\
    &\quad + \Vert \nabla U \Vert_{H^k_x}\int \Vert \rho^\alpha \Vert_{H^{k-1}_x} \Vert (w^\alpha-U)  \Vert_{H^{k-1}_x} \, \mathrm{d}\mu(\alpha),
\end{align*}
thanks to one integration by parts in the source term. Note that the treatment of the convective term, which is similar to the argument already given above, would only requires $k>3/2$ here. However, we still relied on $k-1>3/2$ to use the algebra property of $H^{k-1}_x$ on the source term.
By Hölder inequality in $\alpha$ combined with Young inequality, we get after absorbing the last gradient term:
\begin{align}\label{NS-oneceagain-noidea}
    \frac{\mathrm{d}}{\mathrm{d}t} \Vert U \Vert^2_{H^{k}_x} 
    + \Vert \nabla U \Vert_{H^k_x}^2   
    &\lesssim   \| U\|_{H^k_x}^4  + (1+\Vert \rho^\alpha \Vert_{L^{\upinfty}_\alpha H^{k-1}_x}^2 )(\Vert \vv\Vert_{L^p_\alpha H^{k-1}_x}^2+\Vert U \Vert_{H^{k-1}_x}^2).
\end{align}
We now introduce 
$$y(t):=1+\Vert U(t) \Vert_{H^{k}_x}^2+\int_0^t \Vert \nabla U \Vert_{H^k_x}^2(s) \, \mathrm{d}s,$$
and since $\| U\|_{H^k_x}^3  \leq Y^{3/2} $ and $\| U\|_{H^k_x}^2  \leq Y $,
so that the former differential inequality can be rewritten as
\begin{align*}
    \frac{\mathrm{d}}{\mathrm{d}t}y &\lesssim y^{2}+ (1+\Vert \rho^\alpha \Vert_{L^{\upinfty}_\alpha H^{k-1}_x}^2)(\Vert \vv\Vert_{L^p_\alpha H^{k-1}_x}^2  +y) \\
    & \lesssim (1+\Vert \rho^\alpha \Vert_{L^{\upinfty}_\alpha H^{k-1}_x}^2)(1+\Vert \ww\Vert_{L^p_\alpha H^{k-1}_x}^2) y^2.
\end{align*}
We are thus back to the former proof for the bound at high regularity. In particular, we obtain 
$$    U \in L^\infty([0,T]; H^{ k}_x) \cap L^2([0,T];H^{k+1}_x),
$$
as well as an estimate of the form
\begin{align*}
\int_0^T \| U(s) \|_{H^k_x} \D s \leq  \F\left(T,\int_0^T \Vert \vv(\tau) \Vert_{\mathcal{H}^{k,p}}^2 \D \tau, \Vert \rhorho \Vert_{L^\infty([0,T];L^\infty_\alpha H^{k-1}_x)}, \|U_0\|_{H^{k-1}_x} \right).
 \end{align*}
 where, as alluded above, an additional exponent is present in the integral in time inside the second  argument of  the function $\F$.

\begin{Rem}
In the former two cases, the continuity in time of the fluid velocity field $U$ (either with value in $H^{k-1}_x$ or in $H^k_x$) is a standard consequence of the general theory for the Navier--Stokes equations, and we do not detail it further. Furthermore, in the next paragraph, it will be important to notice that  for all $t$ sufficiently small, we also have the pointwise bound
\begin{equation}
\label{eq:estimate_u_Hk}
\|U(t)\|_{H^{k-1}_x} \leq \Lambda\left(t,\int_0^t \Vert \vv(s)\Vert_{L^p_\alpha H^{k-1}_x} \mathrm{d}s,\Vert \rho^\alpha \Vert_{L^\infty([0,t];L^{\upinfty}_\alpha H^{k-1}_x)}, \|U_0\|_{H^{k-1}_x}\right)^{1/2},
\end{equation}     
or a version with $\Vert \vv(s)\Vert_{L^p_\alpha H^{k-1}_x}$ in the second argument above.
\end{Rem}

\paragraph{Stability at low regularity.} In what follows, we will need to assume that $p \geq 2$. Here, we need to find a stability estimate of type~\eqref{eq:stability_assumption-rho} for the VNS system, in order to ensure that Assumption \ref{ass:stability-rho} is satisfied. Once again, this problem reduces to find a $L^1_T L^2_x$ estimate for $U_2 - U_1$, where $U_1$ and $U_2$ are solutions to
\begin{equation*}
\mbox{for }i=1,2,\qquad \left\{
\begin{aligned}
\partial_t U_i + (U_i \cdot \nabla)U_i + \nabla P_i - \Delta U_i &= \int \rho^{\alpha}_i \left(v^{\alpha}_i - U_i \right) \D \mu(\alpha),\\
\Div(U_i) &= 0,\\
U_i|_{t=0} &= U_0,
\end{aligned}
\right. 
\end{equation*}
where $\vv_1,\vv_2$ are two given families of velocity fields satisfying \eqref{eq:uniform_bound_velocities-rho} and $\rhorho_1,\rhorho_2$ are the corresponding solutions to the continuity equation. 
The difference $W \vcentcolon = U_2 - U_1$ satisfies
\begin{align*}
\partial_t W + (U_2 \cdot \nabla)W + \nabla P -  \Delta W &=-(W \cdot \nabla)U_1 \\
&\quad + \int \left[\rho_2^{\alpha} \left(v_2^{\alpha} - U_2 \right)-\rho_1^{\alpha} \left(v_1^{\alpha} - U_1 \right) \right] \D \mu(\alpha), \\
\mathrm{div} (W)&=0.
\end{align*}
Observe that the last term can be rewritten as
\begin{align*}
    &\int \left[\rho_2^{\alpha} \left(v_2^{\alpha} - U_2 \right)-\rho_1^{\alpha} \left(v_1^{\alpha} - U_1 \right) \right] \D \mu(\alpha)
    \\
    &\qquad=\int \left[(\rho_2^{\alpha}v_2^{\alpha}-\rho_1^{\alpha}v_1^{\alpha}) +(\rho_1^{\alpha}-\rho_2^{\alpha})U_1 - \rho_2^{\alpha}W \right] \D \mu(\alpha) \\
    &\qquad=\int \left[\rho_2^{\alpha}(v_2^{\alpha}-v_1^{\alpha})+(\rho_2^{\alpha}-\rho_1^{\alpha})(v_1^{\alpha}-U_1) - \rho_2^{\alpha}W \right] \D \mu(\alpha) \\
    &\qquad\vcentcolon =\int \left[\mathrm{A}^\alpha+\mathrm{B}^\alpha + \mathrm{C}^\alpha \right] \D \mu(\alpha),
\end{align*}
with
\begin{align}\label{def:BrinkmanstabABC}
\mathrm{A}^\alpha\vcentcolon = \rho_2^{\alpha}(v_2^{\alpha}-v_1^{\alpha}), \ \ \mathrm{B}^\alpha\vcentcolon =(\rho_2^{\alpha}-\rho_1^{\alpha})(v_1^{\alpha}-U_1), \ \ \mathrm{C}^\alpha\vcentcolon =-\rho_2^{\alpha}W.
\end{align}
Testing against $W$, we obtain the estimate
\begin{align*}
\frac{1}{2}\frac{\D}{\D t} \| W \|^2_{L^2_x} +  \| \nabla W\|^2_{L^2_x} &\leq 
\langle (W \cdot \nabla)U_1, W \rangle_{L^2_x}+ 
\int \left\langle \mathrm{A}^\alpha+\mathrm{B}^\alpha + \mathrm{C}^\alpha, W  \right\rangle_{L^2_x} \D \mu(\alpha) \\
& \leq \Vert \nabla U_1 \Vert_{L^{\upinfty}_x} \Vert W \Vert^2_{L^2_x}+\int \left\langle \mathrm{A}^\alpha+\mathrm{B}^\alpha + \mathrm{C}^\alpha, W  \right\rangle_{L^2_x} \D \mu(\alpha),
\end{align*}
since $U_2$ is divergence-free. 
For the term $\mathrm{A}^\alpha$, we write
\begin{align*}
    \int \left\langle  \mathrm{A}^\alpha, W  \right\rangle_{L^2_x} \D \mu(\alpha) &\leq \int \Vert \rho_2^{\alpha}(v_2^{\alpha}-v_1^{\alpha}) \Vert_{L^2_x} \Vert W \Vert_{L^2_x} \D \mu(\alpha) \\
    & \leq \frac{1}{2}\int \Vert \rho_2^{\alpha} \Vert_{L^{\upinfty}_x}^2 \Vert v_2^{\alpha}-v_1^{\alpha} \Vert_{L^2_x}^2\D \mu(\alpha) + \frac{1}{2}\Vert W \Vert_{L^2_x}^2,
\end{align*}
by Young inequality, since $(I,\mu)$ is of total mass $1$. For the other two terms, we have
\begin{align*}
\int \left\langle  \mathrm{B}^\alpha, W \right\rangle_{L^2_x} \D \mu(\alpha) &\leq \int \Vert \rho_2^{\alpha}-\rho_1^{\alpha} \Vert_{\dot H^{-1}_x} \Vert \nabla ((v_1^{\alpha}-U_1) \cdot W) \Vert_{L^2_x} \D \mu(\alpha) \\
&\leq  \int (\Vert \nabla v_1^{\alpha} \Vert_{L^{ \infty}}+\Vert \nabla U_1 \Vert_{L^{\upinfty}_x})\Vert \rho_1^{\alpha}-\rho_2^{\alpha} \Vert_{\dot H^{-1}_x}\Vert W \Vert_{L^2_x} \D \mu(\alpha) \\
& \quad + \int (\Vert v_1^{\alpha} \Vert_{L^{ \infty}}+\Vert  U_1 \Vert_{L^{\upinfty}_x})\Vert \rho_1^{\alpha}-\rho_2^{\alpha} \Vert_{\dot H^{-1}_x}\Vert \nabla W \Vert_{L^2_x} \D \mu(\alpha).
\end{align*}
For the second term, let us use the Hölder inequality in $(I, \mu)$: we have, since $p'=p/(p-1) \leq p$ for $p \geq 2$,
\begin{align*}
    \int \Vert v_1^{\alpha} \Vert_{L^{ \infty}_x}\Vert \rho_1^{\alpha}-\rho_2^{\alpha} \Vert_{\dot H^{-1}_x} \D \mu(\alpha)  
    &\leq \Vert  \vv_1 \Vert_{L^{p'}_\alpha L^{ \infty}_x}\Vert \rhorho_1-\rhorho_2 \Vert_{L^{p}_\alpha \dot H^{-1}_x}  \\
    &\lesssim  \Vert  \vv_1 \Vert_{L^{p}_\alpha L^{ \infty}_x}\Vert \rhorho_1-\rhorho_2 \Vert_{L^{p}_\alpha \dot H^{-1}_x}.
\end{align*}
By Young inequality, we infer that there exists a constant $C>0$ such that 
\begin{align*}
\int \left\langle  \mathrm{B}^\alpha, W \right\rangle_{L^2_x} \D \mu(\alpha) &\leq C
  \int (\Vert \nabla v_1^{\alpha} \Vert_{L^{ \infty}_x}^2+\Vert U_1 \Vert_{W^{1,\infty}_x}^2)\Vert v_2^\alpha- v_1^\alpha  \Vert_{\dot H^{-1}_x}^2 \D \mu(\alpha) \\
  & \quad +  C \Vert  \vv_1 \Vert_{L^{p}_\alpha L^{ \infty}_x}^2 \Vert \vv_2-\vv_1 
  \Vert_{L^p_\alpha \dot H^{-1}_x}^2 \\
& \quad + \frac{1}{4}\Vert W \Vert_{L^2_x}^2+ \frac{1}{4}\Vert \nabla W \Vert_{L^2_x}^2.
\end{align*}
For the last term $\mathrm{C}^\alpha$, we have directly
\begin{align*}
    \int \left\langle  \mathrm{C}^\alpha, W  \right\rangle_{L^2_x} \D \mu(\alpha)
    \leq \Vert W \Vert^2_{L^2_x} \left(\int \Vert \rho_2^{\alpha} \Vert_{L^{\upinfty}_x}\D \mu(\alpha) \right).
\end{align*}
Absorbing the terms with $\nabla W$ in the left-hand side, we obtain
\begin{align*}
\begin{split}
\frac{1}{2}\frac{\D}{\D t} \| W \|^2_{L^2_x} +  \| \nabla W\|^2_{L^2_x}  &\lesssim \left(1+\Vert \nabla U_1 \Vert_{L^{\upinfty}_x} + \Vert \rhorho_2 \Vert_{L^1_\alpha L^{\upinfty}_x}\right) \Vert W \Vert^2_{L^2_x} \\
& \quad + (\Vert \nabla \vv_1 \Vert^2_{L^{\upinfty}_\alpha W^{1, \infty}_x}+\Vert U_1 \Vert^2_{W^{1, \infty}})\Vert \rhorho_1-\rhorho_2  \Vert_{L^2_\alpha \dot H^{-1}_x}^2 \\
& \quad +\Vert  \vv_1 \Vert_{L^{p}_\alpha L^{ \infty}_x}^2 \Vert \rhorho_1-\rhorho_2 \Vert_{L^p_\alpha \dot H^{-1}_x}^2 \\
& \quad + \Vert \rhorho _2 \Vert^2_{L^{\upinfty}_\alpha L^{\upinfty}_x} \Vert \vv_2-\vv_1  \Vert_{L^2_\alpha  L^2_x}^2.
\end{split}
\end{align*}
We end up with 
\begin{align*}
\begin{split}
\frac{1}{2}\frac{\D}{\D t} \| W \|^2_{L^2_x} +  \| \nabla W\|^2_{L^2_x}  &\lesssim \left(1+\Vert \nabla U_1 \Vert_{L^{\upinfty}_x} +  \Vert \rhorho_2\Vert_{L^1_\alpha L^{\upinfty}_x}\right) \Vert W \Vert^2_{L^2_x} \\
& \quad + \left(\Vert \vv_1 \Vert_{\mathcal{H}^{k,p}}^2+\Vert U_1 \Vert^2_{W^{1, \infty}} \right)\Vert \rhorho_1-\rhorho_2 \Vert_{L^p_\alpha \dot H^{-1}_x}^2 \\
& \quad + \Vert \rhorho_2 \Vert^2_{L^{\upinfty}_\alpha L^{\upinfty}_x} \Vert \vv_2-\vv_1 \Vert_{L^p_\alpha  L^2_x}^2.
\end{split}
\end{align*}
Using the pointwise estimates from \eqref{eq:estimate_u_Hk}, we can turn the previous estimate on $\Vert W \Vert_{L^2_x}$ into a suitable stability estimate of the form~\eqref{eq:stability_assumption-rho} thanks to the Gronwall lemma. This completes the proof of the stability at low regularity for the system \eqref{eq:VNS-multiphase}.

\paragraph{Blow-up criterion.} 
We now study a continuation criterion for the solution of the multiphasic Vlasov--Navier--Stokes system. To ease the presentation, and because it corresponds to the framework of Section \ref{Section:monokinetic-stab-theorem} about asymptotic stability, we have decided to work under the assumption $U_0 \in H^k_x$: it corresponds to the case where $U$ and $\vv$ have the same level of regularity.   
\begin{Prop}
\label{prop-VNSblowup}
Let $k \in \N$ such that  $k>1+3/2$ and $p \geq 2$. Let $T^\star$ be the maximal time of existence given by Theorem \ref{thm:existence}. Let $T \in (0, T^\star)$. If there exists a finite constant $M>0$ such that 
\begin{equation}
\label{eq:hyp}
\int_0^{T}  \left( \| \nabla \vv(s) \|_{L^{\upinfty}_\alpha L^{\upinfty}_x}+   \| \nabla U(s) \|_{L^{\upinfty}_x}   \right) \D s \leq M,
\end{equation}
then there exists some continuous nonnegative function $\Psi_{k,p}$, nondecreasing in each of its variable, such that for all $t \in [0,T)$, 
\begin{multline}\label{eq:borneBlowUpVNS}
\| \rhorho(t)\|_{L^{\upinfty}_\alpha H^{k-1}_x} +     \Vert \vv(t) \Vert_{\mathcal{H}^{k,p}} + \| U(t)\|_{H^k_x}  \\
\leq  
\left(\Vert \rhorho_0\|_{L^{\upinfty}_\alpha H^{k-1}_x} + \Vert \vv_0 \Vert_{\mathcal{H}^{k,p}}+\| U_0\|_{ H^{k}_x} \right)
\Psi_{k,p} \left(M,T \right).
\end{multline}
In particular, if $T^\star$ is finite, then
\begin{equation}
\int_0^{T^\star}  \left(  \| \nabla \vv (s) \|_{L^{\upinfty}_\alpha L^{\upinfty}_x} +  \| \nabla U(s) \|_{L^{\upinfty}_x}    \right) \D s = +\infty.
\end{equation}
\end{Prop}
\begin{Rem}
    Proposition~\ref{prop-VNSblowup} pertains to the first set of assumptions in Theorem \ref{thm:appliVNS-incomp}. Up to minor adaptations in the proof, we could also obtain the same blow-up criterion for the second set of assumptions.
\end{Rem}
We first need a preliminary lemma on the Vlasov--Navier--Stokes system. 
\begin{Lem}\label{LM:energyCons+Brink+Chemin}
    Let $(I , \mu)$ be a measure space (with finite mass). Consider any smooth solution $(\rhorho, \vv, U)$ to the multiphasic Vlasov--Navier--Stokes system \eqref{eq:VNS-multiphase} associated to an initial data $(\rhorho_0, \vv_0, U_0)$ with $\rhorho_0 \geq 0$. Then the following estimates hold.
    \begin{itemize}
        \item We have the following energy equality: for all $t \geq 0$
        \begin{align*}
    &\int  \rho^\alpha(t,x) \vert v^\alpha(t,x) \vert^2 \, \mathrm{d}x \, \mathrm{d}\mu(\alpha) +\int \vert U(t,x) \vert^2 \, \mathrm{d}x \\
    & \qquad  + \int_0^t \left\lbrace\int  \rho^\alpha(s,x) \vert( v^\alpha-U)(s,x) \vert^2 \, \mathrm{d}x \, \mathrm{d}\mu(\alpha) +\int \vert \nabla U (s,x) \vert^2 \, \mathrm{d}x \right\rbrace \, \mathrm{d}s \\
    &= \int  \rho^\alpha_0(x) \vert v^\alpha_0(x) \vert^2 \, \mathrm{d}x \, \mathrm{d}\mu(\alpha) +\int \vert U_0(x) \vert^2 \, \mathrm{d}x.
\end{align*}
        \item The Brinkman force
        $\int \rho^\alpha (v^\alpha-U) \, \mathrm{d}\mu(\alpha)$ satisfies  for all $T \geq 0$
\begin{multline*}
   \left\Vert \int \rho^\alpha (v^\alpha-U) \, \mathrm{d}\mu(\alpha) \right\Vert_{L^2([0,T]; L^2_x)}^2 \\
   \lesssim \left(\int  \rho^\alpha_0 \vert v^\alpha_0 \vert^2 \, \mathrm{d}x \, \mathrm{d}\mu(\alpha) +\int \vert U_0 \vert^2 \, \mathrm{d}x \right)   \Vert \rhorho \Vert_{L^\infty([0,T]; L^\infty_\alpha L^\infty_x)}.
\end{multline*}
\item For all $T \geq 0$, there holds
    \begin{multline}\label{eq:CheminNS}
        \int_0^T \Vert U (t) \Vert_{L^{\upinfty}_x}\, \mathrm{d}t \\
        \leq \left(\int  \rho^\alpha_0 \vert v^\alpha_0 \vert^2 \, \mathrm{d}x \, \mathrm{d}\mu(\alpha) +\int \vert U_0 \vert^2 \, \mathrm{d}x \right)\left(1+C(T)  \Vert \rhorho \Vert_{L^\infty([0,T]; L^\infty_\alpha L^\infty_x)} \right),
    \end{multline}
    for some continuous and increasing function $C(T)$.
    \end{itemize}
\end{Lem}
\begin{proof}
    We first rewrite the equation on $(\rhorho,\vv,U)$ as
    \begin{align*}
    \partial_t \rho^\alpha+ \mathrm{div}(\rho^\alpha v^\alpha)&=0, \\
        \partial_t (\rho^\alpha v^\alpha) + \mathrm{div}(\rho^\alpha v^\alpha \otimes v^\alpha) &= \rho^\alpha(U - v^\alpha),\\
\partial_t U + (U\cdot \nabla ) U + \nabla P -  \Delta U &= \int \big( v^\alpha - U ) \rho^\alpha \D \mu(\alpha)
,\\
\Div(U) &= 0.
    \end{align*}
The energy balance is nothing but the energy--dissipation identity (see e.g. \cite{BDGM}) written in the multiphasic formalism:
\begin{multline*}
    \frac{\mathrm{d}}{\mathrm{d}t} \left\lbrace \int  \rho^\alpha \vert v^\alpha \vert^2 \, \mathrm{d}x \, \mathrm{d}\mu(\alpha) +\int \vert U \vert^2 \, \mathrm{d}x  \right\rbrace \\
    +\iint \rho^\alpha \vert( v^\alpha-U) \vert^2 \, \mathrm{d}x \, \mathrm{d}\mu(\alpha) +\int \vert \nabla U \vert^2 \, \mathrm{d}x=0.
\end{multline*}
For the estimate on the Brinkman force, we first write by Cauchy--Schwarz inequality (as $\rho^\alpha \geq 0$)
\begin{align*}
    \left\vert \int \rho^\alpha (v^\alpha-U) \, \mathrm{d}\mu(\alpha) \right\vert \leq \left( \int \rho^\alpha \, \mathrm{d}\mu(\alpha) \right)^{1/2} \left(\int   \rho^\alpha  \vert v^\alpha-U \vert^2 \, \mathrm{d}\mu(\alpha)\right)^{1/2},
\end{align*}
so that, since $(I, \mu)$ is of finite measure, there holds
\begin{align*}
    \left\vert \int \rho^\alpha (v^\alpha-U) \, \mathrm{d}\mu(\alpha) \right\vert^2 \lesssim \Vert \rhorho \Vert_{L^{\upinfty}_\alpha L^{\upinfty}_x} \int   \rho^\alpha  \vert v^\alpha-U\vert^2 \, \mathrm{d}\mu(\alpha).
\end{align*}
By integration and the use of the former energy-dissipation inequality, we get
\begin{align*}
& \int_0^{T} \int \left\vert \int \rho^\alpha (v^\alpha-U) \, \mathrm{d}\mu(\alpha) \right\vert^2  \, \mathrm{d}x \, \mathrm{d}t \\
&\qquad \qquad \lesssim \Vert \rhorho  \Vert_{L^\infty([0,T];L^\infty_{\alpha}L^{\infty}_x)} \int_0^{T}\int     \rho^\alpha  \vert v^\alpha-U\vert^2 \, \mathrm{d}\mu(\alpha)\, \mathrm{d}x \, \mathrm{d}t \\
&\qquad \qquad \lesssim \Vert \rhorho \Vert_{L^\infty([0,T];L^\infty_{\alpha}L^{\infty}_x)}\left( \int  \rho^\alpha_0 \vert v^\alpha_0 \vert^2 \, \mathrm{d}x \, \mathrm{d}\mu(\alpha) +\int \vert U_0 \vert^2 \, \mathrm{d}x \right),
\end{align*}
which proves the second estimate. For the last estimate on $\Vert U \Vert_{L^1_T L^{\upinfty}_x}$, we recast the Navier--Stokes equations on $U$ as
\begin{align*}
\partial_t U-\Delta U= -\mathbb{P}[(U \cdot \nabla)U]+ \mathbb{P}S, \ \ S\vcentcolon =\int \rho^{\alpha} \left(v^{\alpha} - U \right) \D \mu(\alpha),
\end{align*}
where $\mathbb{P}$ stands for the Leray projection on the set of divergence-free vector fields. It yields, following \cite[Lemme 3.2]{Chemin-Leray} and using the homogeneous Besov space $\dot{B}^{3/2}_{2,1}$ (see Section \ref{Section-Besov}) that
\begin{align*}
   \int_0^T \Vert U (t) \Vert_{L^{\upinfty}_x}\, \mathrm{d}t \lesssim \int_0^T \Vert U (t) \Vert_{\dot{B}^{3/2}_{2,1}}\, \mathrm{d}t  &\lesssim C_0+  \int_0^T \left\Vert \int_0^t e^{-\Delta(t'-t)} S(t') \D t'  \right\Vert_{\dot{B}^{3/2}_{2,1}} \D t \\
    & \lesssim C_0+C(T)\| S \|_{L^2([0,T]; L^2_x)} \\
    & \lesssim C_0 + C(T)C_0 \Vert \rhorho \Vert_{L^\infty([0,T];L^\infty_{\alpha}L^{\infty}_x)},
\end{align*}
where $C_0=\left( \int  \rho^\alpha_0 \vert v^\alpha_0 \vert^2 \, \mathrm{d}x \, \mathrm{d}\mu(\alpha) +\int \vert U_0 \vert^2 \, \mathrm{d}x \right)$, thanks to the previous estimate on the Brinkman force $S$. This concludes the proof.
\end{proof}
Let us now obtain the blow-up condition, that is Proposition \ref{prop-VNSblowup}.
\begin{proof}[Proof of Proposition \ref{prop-VNSblowup}]
Let $T>0$ and $M>0$ such that~\eqref{eq:hyp} holds. First, we perform estimates that are integrated in $\alpha$, introducing the quantity
\begin{multline*}
    y_{k}(t)\vcentcolon= \| U(t)\|_{H^k_x}^2+\int \| \rho^\alpha(t)\|_{H^{k-1}_x}^2\D \mu(\alpha) \\
     + \int \| v^\alpha(t)\|_{H^{k-1}_x}^2 \D \mu(\alpha)+ \Vert \nabla v^\alpha(t) \Vert_{L^{\upinfty}_\alpha H^{k-1}_x}^2.
\end{multline*}
We now perform tame estimates to bound $y_k$. From Lemma \ref{LM:SobEstimRHO} on the continuity equation for each $\rho^\alpha$, we have
\begin{align*}
\dfrac{\mathrm{d}}{\mathrm{d}t}\|\rho^\alpha\|_{H^{k-1}_x}^2 &\lesssim \| \nabla v^\alpha \|_{L^{\upinfty}_x}\|\rho^\alpha\|_{H^{k-1}_x}^2  + \| \rho^\alpha\|_{L^{\upinfty}_x}
\left( \| \nabla v^\alpha\|_{H^{k-1}_x}^2+ \|\rho^\alpha\|_{H^{k-1}_x}^2 \right),
\end{align*}
thanks to Young inequality. By integration in $\alpha$, we get
\begin{align}\label{eq-blowup:rhoalphaHk-1VNS}
\dfrac{\mathrm{d}}{\mathrm{d}t} \|\rhorho\|_{L^{\upinfty}_\alpha H^{k-1}_x}^2  &\lesssim (\| \rhorho \|_{L^{\upinfty}_\alpha L^{\upinfty}_x} +\| \nabla \vv\|_{L^{\upinfty}_\alpha L^{\upinfty}_x}) y_k. 
\end{align}
Next, Lemma \ref{LM:existence-vitesseSOBOLEV} on the momentum equation for each $v^\alpha$ yields the estimates
\begin{align*}
\frac{\D}{\D t}  \| v^\alpha\|_{H^{k-1}_x}^2+\| v^\alpha\|_{H^{k-1}_x}^2 &\lesssim \| \nabla v^\alpha\|_{L^{\upinfty}_x} \| v^\alpha\|_{H^{k-1}_x}^2
+    \| U\|_{H^{k-1}_x}\| v^\alpha\|_{H^{k-1}_x}, \\
\frac{\D}{\D t} \|  \nabla v^\alpha \|^2_{H^{k-1}_x}+\| \nabla v^\alpha\|_{H^{k-1}_x}^2 &\lesssim \| \nabla v^\alpha\|_{L^{\upinfty}_x} \| \nabla v^\alpha\|_{H^{k-1}_x}^2 + \|  \nabla U \|_{H^{k-1}_x} \|   \nabla v^\alpha\|_{H^{k-1}_x}.
\end{align*}
By Young inequality and integration in $\alpha$ in the first inequality, we infer
\begin{align}\label{eq-blowup:valphaHk-VNS}
\frac{\D}{\D t}  \int \| v^\alpha\|_{H^{k}_x}^2 \D \mu(\alpha)+\frac{\D}{\D t}  \| \nabla v^\alpha\|_{H^{k-1}_x}^2\lesssim (1+\| \nabla \vv\|_{L^{\upinfty}_\alpha L^{\upinfty}_x}) y_k.
\end{align}
We now perform an $H^k_x$ energy estimate the Navier--Stokes equations for $U$, we obtain
\begin{align*}
    \frac{\mathrm{d}}{\mathrm{d}t} \Vert U \Vert^2_{H^{k}_x} 
    + \Vert \nabla U \Vert_{H^k_x}^2  
    &\lesssim  \Vert \nabla U \Vert_{L^{\upinfty}_x} \| U\|_{H^k_x}^2 + \Vert U \Vert_{L^2_x}\int \Vert \rho^\alpha (w^\alpha-U) \Vert_{L^2_x} \, \mathrm{d}\mu(\alpha)\\
    &\quad + \Vert \nabla U \Vert_{H^k_x}\int \Vert \rho^\alpha (w^\alpha-U) \Vert_{H^{k-1}_x} \, \mathrm{d}\mu(\alpha).
\end{align*}
By a tame estimate (see Proposition \ref{prop-Sobolev}), we have
\begin{multline*}
\int \| \rho^\alpha( v^\alpha-U) \|_{H^{k-1}_x} \, \mathrm{d}\mu(\alpha) \\
\lesssim \int \| \rho^\alpha \|_{H^{k-1}_x}\|  v^\alpha-U \|_{L^{\upinfty}_x}\, \mathrm{d}\mu(\alpha)+ \int \|  \rho^\alpha \|_{L^{\upinfty}_x}\|( v^\alpha-U) \|_{H^{k-1}_x} \, \mathrm{d}\mu(\alpha),
\end{multline*}
followed by Young and Cauchy--Schwarz (in $\alpha$) inequalities, we can absorb the factor $\Vert \nabla U \Vert_{H^k_x}$ in the left-hand side and get
\begin{multline*}
    \frac{\mathrm{d}}{\mathrm{d}t} \Vert U \Vert^2_{H^{k}_x} 
    + \Vert \nabla U \Vert_{H^k_x}^2  \\
    \lesssim  \Vert \nabla U \Vert_{L^{\upinfty}_x} \| U\|_{H^k_x}^2  + \left( \Vert \rhorho \Vert_{L^{\upinfty}_\alpha L^{\upinfty}_x}+ \int \Vert \vv \Vert_{L^{\upinfty}_x} \, \mathrm{d}\mu(\alpha)+  \Vert U \Vert_{L^{\upinfty}_x}\right)y_k.
\end{multline*}
Summing the previous inequality to \eqref{eq-blowup:rhoalphaHk-1VNS} and \eqref{eq-blowup:valphaHk-VNS}, we end up with the following differential inequality on $y_k$:
\begin{align*}
&\frac{\mathrm{d}}{\mathrm{d}t} y_k \\
 &\, \lesssim 
  \left (1+ \Vert \nabla U \Vert_{L^{\upinfty}_x}+\Vert \nabla \vv \Vert_{L^{\upinfty}_\alpha L^{\upinfty}_x}+ \Vert \rhorho \Vert_{L^{\upinfty}_\alpha L^{\upinfty}_x}+ \int \Vert \vv \Vert_{L^{\upinfty}_x} \, \mathrm{d}\mu(\alpha)+  \Vert U \Vert_{L^{\upinfty}_x}\right)y_k.
\end{align*}
By Gronwall inequality, we infer the bound
\begin{align*}
    y_k(t)& \lesssim  \left(\Vert \rho^\alpha_0\|_{L^{\upinfty}_\alpha H^{k-1}_x}^2 + \Vert \vv_0\Vert_{\mathcal{H}^{k,p}}^2+\| U_0\|_{ H^{k}_x}^2 \right) \\
    & \qquad \times \exp\Bigg( \int_0^t ( 1+\Vert \nabla U \Vert_{L^{\upinfty}_x}+\Vert \nabla \vv \Vert_{L^{\upinfty}_\alpha L^{\upinfty}_x}+ \Vert \rhorho \Vert_{L^{\upinfty}_\alpha L^{\upinfty}_x} \\
    &\qquad \qquad \qquad \qquad \qquad \qquad \qquad + \Vert \vv \Vert_{L^1_\alpha L^{\upinfty}_x} +  \Vert U \Vert_{L^{\upinfty}_x})  \, \mathrm{d}\tau \Bigg).
\end{align*}
Here, we have used the fact that $\Vert \vv_0\|_{L^1_\alpha H^{k}_x} \lesssim \Vert \vv_0 \Vert_{\mathcal{H}^{k,p}}$, since $(I, \mu)$ is of finite measure. 

To conclude, it is now enough show how to control the three last terms in the exponential factor. Since $\rho^\alpha$ satisfies a continuity equation with velocity field $v^\alpha$, we have by Lemma \ref{LM:pointwiseEstimRHO}
\begin{align*}
\Vert \rhorho(t) \Vert_{L^{\upinfty}_\alpha L^{\upinfty}_x} &\lesssim \Vert \rhorho_0 \Vert_{L^{\upinfty}_\alpha L^{\upinfty}_x} \exp\left(\int_0^t \Vert \nabla \vv(\tau) \Vert_{L^{\upinfty}_\alpha L^{\upinfty}_x} \mathrm{d} \tau\right). 
\end{align*}
Similarly, the method of characteristics on the momentum equation of $v^{\alpha}$, with source term $U-v^{\alpha}$, shows that $\|  \vv(\tau) \|_{L^1_\alpha L^{\upinfty}_x}$ is controlled by $\int_0^{T}  \| U \|_{L^{\upinfty}_x} \D s$ itself on $[0,T]$, that is
\begin{align*}
\|  \vv(\tau) \|_{L^1_\alpha L^{\upinfty}_x} \lesssim \|  \vv_0 \|_{L^1_\alpha L^{\upinfty}_x}+e^{T} \int_0^{T}  \| U \|_{L^{\upinfty}_x} \D s \lesssim \Vert \vv_0 \Vert_{\mathcal{H}^{k,p}}+ e^{T} \int_0^{T}  \| U \|_{L^{\upinfty}_x}.
\end{align*}
Appealing to \eqref{eq:CheminNS} in Lemma \ref{LM:energyCons+Brink+Chemin}, we can therefore control all the remaining terms since
\begin{align*}
 &\qquad \int_0^T \Vert U (t) \Vert_{L^{\upinfty}_x}\, \mathrm{d}t  \\
 &\lesssim \left(\int  \rho^\alpha_0 \vert v^\alpha_0 \vert^2 \, \mathrm{d}x \, \mathrm{d}\mu(\alpha) +\int \vert U_0 \vert^2 \, \mathrm{d}x \right)\left(1+C(T)  \Vert \rho^\alpha \Vert_{L^\infty([0,T]; L^\infty_\alpha L^\infty_x)} \right) \\
 &\lesssim \left(1+\Vert \vv_0 \Vert_{\mathcal{H}^{k,p}}^2+\| U_0\|_{ H^{k}_x}^2 \right)\left( 1+C(T) \right) \Vert \rhorho_0 \Vert_{L^{\upinfty}_\alpha L^{\upinfty}_x} e^{\int_0^t \Vert \nabla \vv(\tau) \Vert_{L^{\upinfty}_\alpha L^{\upinfty}_x} \mathrm{d} \tau},
\end{align*} 
for $p \geq 2$, for some continuous nondecreasing and nonnegative function $C(T)$. All in all, we have therefore proven that 
\begin{multline}\label{estim-key-blowupVNS}
    \| \rhorho(t)\|_{L^{\upinfty}_\alpha H^{k-1}_x}^2+\| U(t)\|_{H^k_x}^2  
     \lesssim  \left(\Vert \rhorho_0\|_{L^{\upinfty}_\alpha H^{k-1}_x}^2 + \Vert \vv_0 \Vert_{\mathcal{H}^{k,p}}^2+\| U_0\|_{ H^{k}_x}^2 \right) \\
     \times \Psi \left( \Vert \rhorho_0 \Vert_{L^{\upinfty}_\alpha L^{\upinfty}_x},T, \int_0^T \| \nabla \vv \|_{L^{\upinfty}_\alpha L^{\upinfty}_x}, \int_0^T \| \nabla U\|_{L^{\upinfty}_x}  \right),
\end{multline}
 for some continuous nonnegative function $\Psi$ that is nondecreasing in each of its variable. In particular, we obtain the same estimate for $\| U(t)\|_{H^k_x}$ on $[0,T]$. Coming back to the equation satisfied by $v^\alpha$ and performing a tame estimate in $\mathcal{H}^{k,p}$ (recall Proposition~\ref{prop-Sobolev}) yields
\begin{align*}
    \frac{\D}{\D t} \| \vv \|_{\mathcal{H}^{k,p}} &\leq \| \nabla  \vv \|_{L^{\upinfty}_\alpha L^{\upinfty}_x} \| \vv \|_{\mathcal{H}^{k,p}}  + \| U \|_{H^k_x},
\end{align*}
hence by Gronwall inequality
\begin{align*}
\| \vv(t) \|_{\mathcal{H}^{k,p}} \lesssim \left(\| \vv_0 \|_{\mathcal{H}^{k,p}}+ \int_0^t \| U(\tau) \|_{H^k_x} \, \mathrm{d}\tau  \right)\exp\left( \int_0^t \| \nabla  \vv(\tau) \|_{L^{\upinfty}_\alpha L^{\upinfty}_x} \, \mathrm{d}\tau\right).
\end{align*}
Modifying the function $\Psi$ above, the same type estimate as \eqref{estim-key-blowupVNS} then holds for $\| \vv(t) \|_{\mathcal{H}^{k,p}}$. This is enough to obtain the claimed estimate \ref{eq:borneBlowUpVNS} and then to conclude the proof of Proposition \ref{prop-VNSblowup}.
\end{proof}

This concludes the proof of Theorem~\ref{thm:appliVNS-incomp} with the first set of assumptions. The case of the second set of assumptions is very similar and thus skipped; keep in mind though that we use the adaption of Theorem~\ref{thm-LWPmultiphase-densityvelocity-re} that is explained above before the start of the proof.

\end{proof}

\subsection{Nonlinear drag force case}\label{sec:nonlindrag-LWP}
In this section, we explain how to treat the case where one rather considers a possible smooth nonlinear drag force term in the Vlasov--Navier--Stokes system \eqref{eq:VNS_kinLWP} (and its multiphasic version \eqref{eq:VNS-multiphase}). We replace all the terms $U-v$ (resp. $U-v^\alpha$) by $U-v+\Gamma(U-v)$ (resp. by
   $ U - v^\alpha+ \Gamma(U-v^\alpha)$)
where $\Gamma$ satisfies
\begin{align*}
    \Gamma\in \mathscr{C}^{\infty}(\R^3;\R^3), \ \  \Gamma(0)=0.
\end{align*}
More concretely, we consider 
\begin{equation}
		\label{eq:VNS_kinLWP-nonlindrag}
		\left\{
		\begin{aligned}
		\partial_t f + v \cdot \nabla_x f +\Div_v\big( (U-v)f+\Gamma(U-v)\big) &= 0,\\
		\partial_t U + (U \cdot \nabla_x) U + \nabla_x P - \Delta_x U = -&\int_{\R^3} \Big((U-v)+\Gamma(U-v)f\Big) f \, \mathrm{d}v ,\\
		\Div_x(U) &= 0,
		\end{aligned}
		\right.
		\end{equation}
and its multiphasic counterpart
\begin{equation}
\label{eq:VNS-multiphase-nonlindrag}
\left\{ 
\begin{aligned}
\partial_t \rho^\alpha + \Div ( \rho^\alpha v^\alpha) &= 0,\\
\partial_t v^\alpha + (v^\alpha \cdot \nabla) v^\alpha &= U - v^\alpha+\Gamma(U-v^\alpha),\\
\partial_t U + (U\cdot \nabla ) U + \nabla P -  \Delta U &= -\int \Big( U-v^\alpha  +\Gamma(U-v^\alpha)\Big) \rho^\alpha \D \mu(\alpha)
,\\
\Div(U) &= 0.
\end{aligned}
\right.
\end{equation}
It turns out that the local well-posedness results  
Theorem \ref{thm:appliVNS-incomp} and Theorem \ref{eq:VNS-multiphase-nonlindrag} still pertain in that context, at the cost of requiring some compact support in velocity for the kinetic distribution function (which translates into taking $\vv \in \mathcal{H}^{k,\infty}$ on the multiphasic side, see Remark \ref{rem:bounded_velocity}). This constraint might be of technical nature and 
        we leave open the question of relaxing this condition.
\begin{Thm}
    Under the above assumption on the $\Gamma$, Theorem \ref{thm:appliVNS-incomp} still holds for \eqref{eq:VNS_kinLWP-nonlindrag} with $p=\infty$, and Corollaries \ref{coro-multiEuler-VNSincomp}--\ref{coro:appliVNS-incomp}--\ref{coro:appliVNS-incomp-Mixed} still holds if $f_0$ is compactly supported in velocity.
\end{Thm}
To our knowledge, such result has not been previously obtained in the literature. It will follow from a modification of the proof of Theorem \ref{thm:appliVNS-incomp} and we therefore strive for 
the equivalent of the bound at high regularity and stability at low regularity for \eqref{eq:VNS_kinLWP-nonlindrag}. For the sake of simplicity, we only consider the case $U_0 \in H^{k-1}_x$ and highlight the essential modifications compared to the previous proof. 

\medskip

$\bullet$ \textbf{Bound in high regularity.} Thanks to the composition estimate from Proposition \ref{prop-Sobolev}, we know that 
\begin{align*}
    \Vert \Gamma(U-v^\alpha) \Vert_{H^{k-1_x}} &\lesssim (1+ \Vert U-v^\alpha \Vert_{L^\infty_x})^{k-1} \Vert U-v^\alpha \Vert_{H^{k-1}_x} \\
    & \lesssim (1+ \Vert U\Vert _{H^{k-1}_x}+\Vert v^\alpha \Vert_{H^{k-1}_x})^{k-1} (\Vert U \Vert_{H^{k-1}_x}+\Vert v^\alpha \Vert_{H^{k-1}_x}) \\
    & \lesssim \Vert v^\alpha \Vert_{H^{k-1}_x}^k+\Vert v^\alpha \Vert_{H^{k-1}_x}+\Vert U\Vert _{H^{k-1}_x}+\Vert U\Vert _{H^{k-1}_x}^k,
\end{align*}
from which we deduce the following energy inequality
\begin{align*}
    \begin{split}
        &\frac{\D}{\D t}\|U\|_{H^{k-1}_x}^2 + \Vert \nabla U \Vert_{H^{k-1}_x}^2 \\ &\leq K  \Bigg( \| U \|^4_{H^{k-1}_x} \\
        & \qquad \qquad  + \Vert \rhorho \Vert_{L^{\upinfty}_\alpha H^{k-1}_x} \left(\|U\|_{H^{k-1}_x}\mathcal{Q}(\Vert \vv \Vert_{\mathcal{H}^{k,p}})+ \|U \|_{H^{k-1}_x}^{k+1}+\|U \|_{H^{k-1}_x}^2\right)\Bigg),
        \end{split}
    \end{align*}
for $k>1+3/2$ and some constant $K>0$, and where we have set $\mathcal{Q}(z)\vcentcolon=z+z^k$.
Once again, we introduce 
$$y(t):=1+\Vert U(t) \Vert_{H^{k-1}_x}^2+\int_0^t \Vert \nabla U \Vert_{H^{k-1}_x}^2(s) \, \mathrm{d}s,$$
and since $k \in \N$, which implies that $k \geq 3$ and $\frac{k+1}{2} \geq 2$, it satisfies the following differential inequality:
\begin{align*}
    y'(t) &\leq K\left( y^2+ \Vert \rhorho \Vert_{L^{\upinfty}_\alpha H^{k-1}_x}(y+y^{\frac{k+1}{2}})+ \mathcal{Q}(\Vert \vv \Vert_{\mathcal{H}^{k,p}})y^{1/2}\right) \\
    & \lesssim K (1+\Vert \rhorho \Vert_{L^{\upinfty}_\alpha H^{k-1}_x})(1+\mathcal{Q}(\Vert \vv \Vert_{\mathcal{H}^{k,p}})) y^{\frac{k+1}{2}}.
\end{align*}
Hence, there holds
\begin{align*}
    \frac{\mathrm{d}}{\mathrm{d}t}y^{-\frac{k-1}{2}}=-\frac{k-1}{2}y^{-\frac{k+1}{2}} y' \gtrsim -K(1+\Vert \rho^\alpha \Vert_{L^{\upinfty}_\alpha H^{k-1}_x})(1+\mathcal{Q}(\Vert \ww\Vert_{L^p_\alpha H^{k-1}_x})),
\end{align*}
so that, after integration in time,  we get for all times $t \in [0,T]$
\begin{multline*}
    y(t) \leq  \\
    \frac{1+ \Vert U_0 \Vert_{H^{k-1}_x}^2}{\left[1-K(1+\Vert U_0 \Vert_{H^{k-1}_x}^2)^{\frac{k-1}{2}}(1+\Vert \rho^\alpha \Vert_{L^{\upinfty}([0,t];L^{\upinfty}_\alpha H^{k-1}_x)})(t+\int_0^t \mathcal{Q}(\Vert \ww(s)\Vert_{L^p_\alpha H^{k-1}_x}) \, \mathrm{d}s) \right]^{\frac{2}{k-1}}}.
\end{multline*}
provided the former denominator is strictly positive. Exactly as in the case of the linear drag, we find an estimate of the form
    \begin{multline*}
		\int_0^T \| U(s)\|_{H^k_x} \, \mathrm{d}s \\
        \leq (T+T^{1/2})\F\left(T, \int_0^T \mathcal{Q}\big(\Vert \vv(t) \Vert_{\mathcal{H}^{k,p}}\big) \D t, \| \rhorho\|_{L^\infty([0,T]; L^\infty_\alpha H^{k-1}_x)}, \mathsf N_k(X) \right),
		\end{multline*}
    where $\mathsf N_k(X)=\|U_0\|_{H^{k-1}_x}$. The estimate that we have obtained here is slightly different from the one we have put forward in Assumption~\ref{ass:bound_high_regularity-rho}, but it actually fits the allowed generalization highlighted in Remark \ref{rem-bound-high-ref-functionQ-withrho}.

        \medskip
        
$\bullet$ \textbf{Stability at low regularity.} Compared to the case of a linear drag force, we now write
    \begin{align*}
        &\int \rho_2^\alpha \Gamma(v_2^\alpha-U_2) \, \mathrm{d}\mu(\alpha)-\int \rho_1^\alpha \Gamma(v_1^\alpha-U_1)\, \mathrm{d}\mu(\alpha) \\
        &=\int (\rho_2^\alpha-\rho_1^\alpha) \Gamma(v_2^\alpha-U_2)\, \mathrm{d}\mu(\alpha)+\int \rho_1^\alpha \big[\Gamma(v_2^\alpha-U_2)-\Gamma(v_1^\alpha-U_1) \big]\, \mathrm{d}\mu(\alpha).
    \end{align*}
    Let us use the notation $W=U_2-U_1$. For the second term we have by the mean value inequality that 
    \begin{align*}
        \int \left\langle  \rho_1^\alpha \big[\Gamma(v_2^\alpha-U_2)-\Gamma(v_1^\alpha-U_1) \big] , W \right\rangle_{L^2_x} \, \mathrm{d}\mu(\alpha)&\lesssim \int \vert \rho_1^\alpha \vert \vert v_2^\alpha-v_1^\alpha\vert \vert W \vert \, \mathrm{d}x \, \mathrm{d}\mu(\alpha) \\
        & \quad + \int \vert \rho_1^\alpha \vert  \vert W \vert^2 \, \mathrm{d}x \, \mathrm{d}\mu(\alpha).
    \end{align*}
    These two terms  can therefore be treated as contribution of the terms $\mathrm{A}^\alpha$ and $\mathrm{C}^\alpha$ in \eqref{def:BrinkmanstabABC}. For the first term, we write thanks to the composition estimate of Proposition \ref{prop-Sobolev} and Sobolev embedding that with $k>1+d/2$
    \begin{align*}
  & \int \langle (\rho_2^\alpha-\rho_1^\alpha) \Gamma(v_2^\alpha-U_2) , W \rangle_{L^2_x} \, \mathrm{d}\mu(\alpha) \\
  &\qquad \lesssim  \int \Vert \rho_2^{\alpha}-\rho_1^{\alpha} \Vert_{\dot H^{-1}_x} \Vert \nabla (\Gamma(v_2^{\alpha}-U_2) \cdot W) \Vert_{L^2_x} \D \mu(\alpha) \\
   &\qquad \lesssim \int \Vert \nabla (\Gamma(v_2^{\alpha}-U_2)) \Vert_{L^{\upinfty}_x} \Vert \rho_1^{\alpha}-\rho_2^{\alpha} \Vert_{\dot H^{-1}_x}\Vert W \Vert_{L^2_x} \D \mu(\alpha)  \\
   &\qquad\qquad+ \int \Vert \Gamma(v_2^{\alpha}-U_2) \Vert_{L^{\upinfty}_x} \Vert \rho_1^{\alpha}-\rho_2^{\alpha} \Vert_{\dot H^{-1}_x}\Vert\nabla W \Vert_{L^2_x} \D \mu(\alpha)
   \\
   &\qquad \lesssim \int \Vert \Gamma(v_2^{\alpha}-U_2) \Vert_{H^k_x} \Vert \rho_1^{\alpha}-\rho_2^{\alpha} \Vert_{\dot H^{-1}_x}(\Vert W \Vert_{L^2_x}+\Vert\nabla W \Vert_{L^2_x} ) \D \mu(\alpha)   \\
   &\qquad \lesssim \int (1+ \Vert v_2^{\alpha}\Vert_{L^{\upinfty}_x}^k +\Vert U_2 \Vert_{L^{\upinfty}_x}^k)  (\Vert v_2^{\alpha} \Vert_{H^k_x}+ \Vert U_2\Vert_{H^k_x}) \Vert \rho_1^{\alpha}-\rho_2^{\alpha} \Vert_{\dot H^{-1}_x}\\
   &\qquad\qquad \qquad \times (\Vert W \Vert_{L^2_x}+\Vert\nabla W \Vert_{L^2_x} ) \D \mu(\alpha) \\
   &\qquad \lesssim \left(1+\Vert \vv_2 \Vert_{\mathcal{H}^{k,\infty}}^{2(k+1)}+ \Vert U_2\Vert_{H^k_x}^{2(k+1)} \right) \int \Vert \rho_1^{\alpha}-\rho_2^{\alpha} \Vert_{\dot H^{-1}_x}^2 \, \mathrm{d}\mu(\alpha)  \\
   &\qquad\qquad + \frac{1}{4}\Vert W \Vert_{L^2_x}^2 +  \frac{1}{4}\Vert \nabla W \Vert_{L^2_x}^2.
    \end{align*}
Note that the nonlinear drag term has compelled us to use the space $\mathcal{H}^{k,\infty}$ for the velocities (that is taking the exponent $p=\infty$ compared to the previous situation $p \geq 2$ in the linear drag case). We can now conclude the proof of stability at low regularity as in the case of the linear drag force.

\sectionhead
  [Local well-posedness for the compressible VNS equations]
  {Local well-posedness for the compressible
   Vlasov--Navier--Stokes equations}\label{SubsecVNScomp-multiphaseLWP}

In this section, we now consider the multiphasic compressible Vlasov--Navier--Stokes system \eqref{eq:VNScomp-multiphase}, still set on the torus $\T^3$ or the whole space $\R^3$.

\begin{Thm}\label{thm:appliVNS-comp}
Let $k \in \N$ such that $k>2+3/2$. Let $f_0$ be an initial distribution function such that
$$f_0(x, \cdot)=\int_{I} \rho^\alpha_0(x) \otimes  \delta_{v=v^\alpha_0(x)} \, \mathrm{d}\mu(\alpha),$$
for some set of labels $(I, \mu)$ and some regular multiphasic distribution $(\rhorho_0, \vv_0)$. 
Let $(\varrho_0, U_0) \in H^{k-1}_x \times H^{k-1}_x$ such that $\inf_{x}\vert \varrho_0(x) \vert>0$.

\begin{enumerate}

\item If $(\rhorho_0, \vv_0) \in L^\infty_\alpha H^{k-1}_x \times \mathcal{H}^{k,p} $ is such that $\rho^\alpha \geq 0$ with $p \in [2, \infty]$, there exist $T>0$ and a unique solution $(\rhorho, \vv, \varrho, U)$ on $[0,T]$ to the multiphasic system \eqref{eq:VNScomp-multiphase} with 
    \begin{align*}
    &\rhorho \in L^\infty([0,T]; L^\infty_\alpha H^{k-1}_x ), \ \ \vv \in L^\infty([0,T]; \mathcal{H}^{k,p}), \\
    &U \in L^\infty([0,T]; H^{k-1}_x) \cap L^2([0,T]; H^{k}_x), \ \ \varrho \in L^\infty([0,T]; H^{k-1}_x).
\end{align*}

\item Let $\Theta: I \to \R_+$ and $ \lambda : I \to \R^d$ be measurable maps satisfying $\langle \lambda^\alpha \rangle^p \lesssim \Theta(\alpha)$ with $p \in [1, \infty)$. If $(\rhorho_0, \vv_0=\ww_0 + \boldsymbol{\lambda})$ with $(\rhorho_0, \ww_0 ) \in L^1_{\alpha, \Theta} H^{k-1}_x  \times \mathcal{H}^{k,\infty}$ is such that $\rho^\alpha_0 \geq 0$, there exist $T>0$ and a unique solution $(\rhorho, \vv, U)$ on $[0,T]$ to the multiphasic system \eqref{eq:VNScomp-multiphase} with
    \begin{align*}
    &\rhorho \in L^\infty([0,T]; L^1_{\alpha, \Theta} (H^{k-1}_x) ), \\
    &\vv= e^{-t}(\boldsymbol{\lambda}+ \widetilde\ww),\, \widetilde\ww \in L^\infty([0,T]; \mathcal{H}^{k,\infty}), \\ 
    &U \in L^\infty([0,T]; H^{k-1}_x) \cap L^2([0,T]; H^{k}_x), \ \ \varrho \in L^\infty([0,T]; H^{k-1}_x).
    \end{align*}
\end{enumerate}

In  both cases, defining
$$
f(t,x,\cdot) \vcentcolon= \int_{I} \rho^\alpha(t,x) \otimes  \delta_{v=v^\alpha(t,x)} \, \mathrm{d}\mu(\alpha),
$$
$(f,\varrho, U)$ is the unique weak solution on $[0,T]$ of \eqref{eq:VNS-comp_kinLWP} with the initial data $(f_0,\varrho_0, U_0)$,
in the sense of Definition \ref{def:solution_coupled}.  
It is continuous in time, $H^{k-1}_x$--strongly in space and weakly in velocity, Definition~\ref{def:solution_Vlasov+timecontinuity}, it holds
  \begin{align*}
    &\varrho \in \mathscr{C}^0([0,T]; H^{k-1}_x), \ \ U \in \mathscr{C}^0([0,T]; H^{k-1}_x) \cap L^2([0,T]; H^{k}_x),
\end{align*}
and  higher order moments in velocity enjoy the following regularity:
\begin{align*}
  &(t,x) \mapsto \int_{\R^d} \vert v \vert^j f(t,x,\D v) \in L^\infty([0,T]; H^{k-1}_x), \ \ j=0, \ldots,\lfloor p\rfloor.  \end{align*}
Finally, if $(\varrho_0,U_0)  \in H^k_x \times H^k_x$ then the whole statement remains true, with the regularity
\begin{align*}
    \varrho \in \mathscr{C}^0([0,T]; H^{ k}_x), \ \ U \in \mathscr{C}^0([0,T]; H^{ k}_x) \cap L^2([0,T];H^{k+1}_x).
\end{align*}
\end{Thm}

\begin{Rem}
    As for the incompressible case, our approach allows us to treat the case of a nonlinear smooth drag term in \eqref{eq:VNScomp-multiphase} by choosing $p=\infty$, in the first set of assumptions. The proof being similar to the one explained in Section \ref{sec:nonlindrag-LWP}, at the cost of heavier technical modifications, we prefer not to include it here in the compressible case. Details are left to the reader.
\end{Rem}

As in the incompressible case, we can now deduce the same consequences at the kinetic level.
\begin{Cor}
    Let $k\in \N$ such that $k> 2+3/2$. Let $(\varrho_0, U_0) \in H^{k-1}_x \times H^{k-1}_x$ such that $\inf_{x}\vert \varrho_0(x) \vert>0$. Then the conclusion of Corollaries \ref{coro-multiEuler-VNSincomp}--\ref{coro:appliVNS-incomp}--\ref{coro:appliVNS-incomp-Mixed} holds true for the compressible Vlasov--Navier--Stokes system \eqref{eq:VNS-comp_kinLWP}, with
    $$\varrho \in \mathscr{C}^0([0,T]; H^{k-1}_x), \ \ U \in \mathscr{C}^0([0,T]; H^{k-1}_x) \cap L^2([0,T]; H^{k}_x).$$
    If $(\varrho_0,U_0)  \in H^k_x \times H^k_x$ then the whole statements remain true, with the regularity
\begin{align*}
    \varrho \in \mathscr{C}^0([0,T]; H^{ k}_x), \ \ U \in \mathscr{C}^0([0,T]; H^{ k}_x) \cap L^2([0,T];H^{k+1}_x).
\end{align*}
\end{Cor}

\begin{proof}[Proof of Theorem \ref{thm:appliVNS-comp}]
We focus on the first set of assumptions. Thanks to Theorem~\ref{thm-LWPmultiphase-densityvelocity} in Chapter~\ref{Part1-LWP}, to prove Theorem~\ref{thm:appliVNS-comp},
 it is enough to check that the Assumptions~\ref{ass:bound_high_regularity-rho} and~\ref{ass:stability-rho}  are satisfied for the system \eqref{eq:VNScomp-multiphase}.
For the second set of assumptions on the initial condition, we can rely again on the adaptation of Theorem~\ref{thm-LWPmultiphase-densityvelocity-re} provided in Section~\ref{SubsecVNSmultiphaseLWP}; details are left to the reader.

The force field of interest is again $ U - v^\alpha$. Let us prove that the two main Assumptions \ref{ass:bound_high_regularity-rho} and \ref{ass:stability-rho} holds for \eqref{eq:VNScomp-multiphase}. For the sake of presentation, we choose $\lambda=\mu=1$ in the definition of the Lamé operator $\mathcal{A}$ (see \eqref{def-Laméoperator}), the general case being a minor adaptation of the proof.

\paragraph{Bounds in high regularity.}
Let $k \in \N$ such that $k>2+3/2$. For the sake of conciseness, we focus on the case $(\varrho_0, U_0) \in H^{k-1}_x \times H^{k-1}_x$, that is we consider that the fluid data has the same level of regularity as the densities $\rhorho_0$. The case where $(\varrho_0, U_0) \in H^{k-1}_x \times H^{k-1}_x$ can be handled by arguments similar to the ones given in the incompressible case, and we leave it to the reader\footnote{Note that the refined energy estimates given in Section \ref{Section:VNS-comp} (for the proof of asymptotic stability) will actually be performed in that case.}.

As for the incompressible case, it is enough to bound $U$ in $L^1([0,T]; H^k_x)$.  First, in view of Lemma \ref{LM:pointwiseEstimRHO}, we can work under the assumption that
\begin{align*}
    m\vcentcolon=\inf \varrho(t)>0,
\end{align*}
at least on a small interval of time, if it holds at time $t=0$ and if we control the Lipschitz norm of $u$ integrated in time. To ease readability, we do not detail such standard procedure. It will anyway be carried on in our perturbative long-time well-posedness proof in Section \ref{Section:VNS-comp}.

\medskip

By dividing by $\varrho(t)$, we therefore write the equation on $U$ as 
\begin{align*}
\partial_t U -I(\varrho) \mathcal{A}U =-(U\cdot \nabla ) U -\nabla \pi (\varrho)   + I(\varrho)\int \big( v^\alpha - U ) \rho^\alpha \D \mu(\alpha),
\end{align*}
where $\pi'(z)=P'(z)/z$ and where $I(z)=1/z$ on $[m,+\infty)$. We can always extend smoothly the functions $\pi$ and $I$ in $\R$ (for instance by $0$ at $z=0$), since there will only be evaluated on the range of $\varrho$. In view of this remark, composition rule in Sobolev spaces from Proposition \ref{prop-Sobolev} applies and we have 
\begin{align*}
\Vert \pi(\varrho) \Vert_{H^{k-1}_x}+ \Vert I(\varrho) \Vert_{H^{k-1}_x} \lesssim (1+ \Vert \varrho \Vert_{L^{\upinfty}_x})^{k-1}\Vert \varrho \Vert_{H^{k-1}_x}.
\end{align*}
As in the incompressible case, we perform an energy estimate in $H^{k-1}_x$, that entails
\begin{align*}
&\dfrac{\mathrm{d}}{\mathrm{d}t} \Vert U \Vert_{H^{k-1}_x}^2-\sum_{\vert \beta \vert \leq k-1} \int I(\varrho) \partial^\beta U \cdot \partial^\beta (\mathcal{A}U)  \\
&\lesssim \Vert  U\Vert_{H^{k-1}_x}^3+ \Vert U\Vert_{H^{k-1}_x}\Vert I(\varrho)\Vert_{H^{k-1}_x} \int \Vert \rho^\alpha (v^\alpha-U ) \Vert_{H^{k-1}_x}\, \mathrm{d}\mu(\alpha)  \\
& \qquad \qquad \qquad \qquad \quad + \langle U, \nabla \pi(\varrho) \rangle_{H^{k-1}_x}+\sum_{\vert \beta \vert \leq k-1 }\Vert [\partial^\beta, I(\varrho)](\mathcal{A}U) \Vert_{L^2_x} \Vert U \Vert_{H^{k-1}_x}.
\end{align*}
Here, we have just treated the convective term thanks to 
\begin{align*}
    \langle (U\cdot \nabla ) U, U \rangle_{H^{k-1}_x} \lesssim \Vert \nabla U \Vert_{L^{\upinfty}_x} \Vert U \Vert_{H^{k-1}_x}^2 \lesssim \Vert U \Vert_{H^{k-1}_x}^3,
\end{align*}
by Sobolev embedding, that requires $k-1>1+3/2$, even though proceeding like in the proof for the incompressible case is also possible.

We need to understand the new term associated to diffusion in the left-hand side and the new last terms in the right-hand side.
\begin{itemize}
\item For the diffusion term, we write  carefully for all $\vert \beta \vert \leq k-1$ by subsequent integrations by parts
\begin{align*}
 -\int I(\varrho) \partial^\beta U \cdot \partial^\beta (\mathcal{A}U)& =-\int I(\varrho) \partial^\beta U\cdot   \partial^\beta (\Delta U+ \nabla \mathrm{div}(U)) \\
 &=\int I(\varrho)\big( \vert \nabla \partial^\beta U \vert^2 +  \vert \mathrm{div} (\partial^\beta U) \vert^2 \big) \\
 & \quad +\int   \partial^\beta U  \cdot (\nabla  (\partial^\beta U) \nabla I(\varrho)) \\
 &\quad +  \int \mathrm{div}(\partial^\beta U)\partial^\beta U \cdot \nabla I(\varrho).
\end{align*}
We will put these two last terms in the right-hand side and estimate them: they are actually suitable for an absorption argument by writing
\begin{multline*}
\int   \partial^\beta U  \cdot (\nabla  (\partial^\beta U) \nabla I(\varrho)) + \int \mathrm{div}(\partial^\beta U)\partial^\beta U \cdot \nabla I(\varrho) \\
\leq \Vert (\nabla,\mathrm{div})U \Vert_{H^{k-1}_x} \Vert U \Vert_{H^{k-1}_x} \Vert \nabla I(\varrho) \Vert_{L^{\upinfty}_x}.
\end{multline*}
\item For $\langle U, \nabla \pi(\varrho) \rangle_{H^{k-1}_x}$, we integrate by parts and get $\langle U, \nabla \pi(\varrho) \rangle_{H^{k-1}_x}=-\langle \mathrm{div}(U), \pi(\varrho) \rangle_{H^{k-1}_x}$, that can be absorbed by Cauchy--Schwarz and Young inequalities.
\item For the commutators  $\Vert [\partial^\beta, I(\varrho)](\mathcal{A}U) \Vert_{L^2_x}$ with $\vert \beta \vert \leq k-1$, we use the commutator law from Proposition \ref{prop-Sobolev} and get
\begin{align*}
\sum_{\vert \beta \vert \leq k-1} \Vert &[\partial^\beta, I(\varrho)](\mathcal{A}U) \Vert_{L^2_x} \\
&\lesssim \Vert \nabla I(\varrho) \Vert_{L^{\upinfty}_x} \sum_{\vert \gamma \vert \leq k-2} \Vert \partial^\gamma \mathcal{A}U \Vert_{L^2_x}
+\Vert \mathcal{A}U \Vert_{L^{\upinfty}_x} \Vert I(\varrho) \Vert_{H^{k-1}_x} \\
&\lesssim \Vert \nabla I(\varrho) \Vert_{L^{\upinfty}_x} \Vert (\nabla, \mathrm{div})(U) \Vert_{H^{k-1}_x}
+\Vert \mathcal{A}U \Vert_{L^{\upinfty}_x} \Vert I(\varrho) \Vert_{H^{k-1}_x}.
\end{align*}
Hence the last term is estimated by
$$\Vert \nabla I(\varrho) \Vert_{L^{\upinfty}_x} \Vert (\nabla, \mathrm{div})(U) \Vert_{H^{k-1}_x}\Vert U \Vert_{H^{k-1}_x}
+\Vert \mathcal{A}U \Vert_{L^{\upinfty}_x} \Vert I(\varrho) \Vert_{H^{k-1}_x}\Vert U \Vert_{H^{k-1}_x},$$
and the first piece is fine for an absorption argument.
\end{itemize}
Gathering all the previous estimates together, we have
\begin{align*}
&\dfrac{\mathrm{d}}{\mathrm{d}t} \Vert U \Vert_{H^{k-1}_x}^2+ \frac{1}{\Vert \varrho \Vert_{L^{\upinfty}_x}} \Vert (\nabla, \mathrm{div}) U\Vert_{H^{k-1}_x}^2 \\
&\qquad \qquad \qquad\lesssim \Vert  U\Vert_{H^{k-1}_x}^3+ \Vert U\Vert_{H^{k-1}_x}\Vert I(\varrho)\Vert_{H^{k-1}_x} \int \Vert \rho^\alpha (v^\alpha-U ) \Vert_{H^{k-1}_x}\, \mathrm{d}\mu(\alpha) \\
&\qquad \qquad \qquad \quad + \Vert \mathrm{div}(U) \Vert_{H^{k-1}_x}  \Vert \pi(\varrho) \Vert_{H^{k-1}_x}   \\
&\qquad \qquad \qquad \quad + \Vert \nabla I(\varrho) \Vert_{L^{\upinfty}_x} \Vert (\nabla, \mathrm{div})(U) \Vert_{H^{k-1}_x}\Vert U \Vert_{H^{k-1}_x} \\
&\qquad \qquad \qquad \quad + \Vert \mathcal{A}U \Vert_{L^{\upinfty}_x} \Vert I(\varrho) \Vert_{H^{k-1}_x}\Vert U \Vert_{H^{k-1}_x}.
\end{align*}
By Sobolev embedding, using $k-1>1+3/2$ to write
\begin{align*}
    \Vert \mathcal{A}U \Vert_{L^{\upinfty}_x} \lesssim \Vert \mathcal{A}U \Vert_{H^{k-2}_x}  \lesssim \Vert (\nabla, \mathrm{div})U \Vert_{H^{k-1}_x}, \ \ \Vert \nabla I(\varrho) \Vert_{L^{\upinfty}_x} \lesssim \Vert I(\varrho) \Vert_{H^{k-1}_x},
\end{align*}
we can absorb all the terms $\Vert (\nabla, \mathrm{div}) U\Vert_{H^{k-1}_x}^2$ in the left-hand side and get
\begin{align*}
&\dfrac{\mathrm{d}}{\mathrm{d}t} \Vert U \Vert_{H^{k-1}_x}^2+ \left(\frac{1}{\Vert \varrho \Vert_{L^{\upinfty}_x}}-\eta\right) \Vert (\nabla, \mathrm{div}) U\Vert_{H^{k-1}_x}^2 \\
&\lesssim \Vert  U\Vert_{H^{k-1}_x}^3+ \Vert U\Vert_{H^{k-1}_x}\Vert I(\varrho)\Vert_{H^{k-1}_x} \int \Vert \rho^\alpha (v^\alpha-U ) \Vert_{H^{k-1}_x}\, \mathrm{d}\mu(\alpha) +C_\eta  \Vert \pi(\varrho) \Vert_{H^{k-1}_x}^2 \\
& \quad  + C_\eta \Vert I(\varrho) \Vert_{H^{k-1}_x}^2 \Vert U \Vert_{H^{k-1}_x}^2  \\
& \lesssim \Vert  U\Vert_{H^{k-1}_x}^3+  (1+ \Vert \varrho \Vert_{L^{\upinfty}_x})^{(k-1)}\Vert \varrho \Vert_{H^{k-1}_x} \Vert U\Vert_{H^{k-1}_x}\int \Vert \rho^\alpha (v^\alpha-U ) \Vert_{H^{k-1}_x}\, \mathrm{d}\mu(\alpha)  \\
& \quad  + C_\eta  (1+ \Vert \varrho \Vert_{L^{\upinfty}_x})^{4(k-1)}\Vert \varrho \Vert_{H^{k-1}_x}^4 +\Vert U \Vert_{H^{k-1}_x}^4.
\end{align*}
for all $\eta>0$, where $C_\eta>0$ is a function diverging as $\eta \rightarrow 0$, thanks to the composition rule recalled above for $I(\varrho)$ and $\pi(\varrho)$. 
On the other hand, we have by (the proof of) Lemma \ref{LM:SobEstimRHO} that
\begin{align*}
   \dfrac{\mathrm{d}}{\mathrm{d}t} \Vert \varrho \Vert_{H^{k-1}_x}^2 &\lesssim \Vert \nabla U \Vert_{L^{\upinfty}_x} \Vert \varrho \Vert_{H^{k-1}_x}^2
   + \Vert \nabla U \Vert_{H^{k-1}_x} \Vert \varrho \Vert_{L^{\upinfty}_x} \Vert \varrho \Vert_{H^{k-1}_x} \\
   &\leq C_{\eta}(\Vert U \Vert_{H^{k-1}_x}^4
   + \Vert \varrho \Vert_{H^{k-1}_x}^4)+\eta \Vert \nabla U \Vert_{H^{k-1}_x}^2,
\end{align*}
by Young inequality and Sobolev embedding. Hence, setting
$$x(t)\vcentcolon=\Vert \varrho(t) \Vert_{H^{k-1}_x}^2+\Vert U(t) \Vert_{H^{k-1}_x}^2,$$ 
we get by summing the two estimates that
\begin{align*}
    &x'+ \left(\frac{1}{\Vert \varrho \Vert_{L^{\upinfty}_x}}-2\eta\right) \Vert (\nabla, \mathrm{div}) U\Vert_{H^{k-1}_x}^2 \\
    &\leq K\Big[x^{3/2}+\Vert \rhorho(t) \Vert_{L^{\upinfty}_\alpha H^{k-1}_x} (1+ x)^{4(k-1)}x  (\Vert \vv \Vert_{\mathcal{H}^{k,p}}+ x^{1/2} ) \\
    &\qquad \qquad \qquad \qquad \qquad \qquad \qquad \qquad + C_\eta (1+ x)^{4(k-1)} x^2 \Big],
\end{align*}
for some constant $K>0$, by Sobolev embedding. By Lemma \ref{LM:pointwise-densité} and Sobolev embedding, there exists $\kappa>0$ and $K>0$ such that we have
$$\frac{1}{\Vert \varrho(t) \Vert_{L^{\upinfty}_x}} \geq \kappa \Vert \varrho_0 \Vert_{H^{k-1}_x}^{-1} \exp\left(-K\int_0^t \Vert  U(\tau) \Vert_{H^{k-1}_x} \mathrm{d} \tau\right).$$
As a consequence, picking $\eta=\kappa \Vert \varrho_0 \Vert_{H^{k-1}_x}^{-1}/8$, there holds (modifying $K$ if needed)
\begin{align*}
    &x' \leq K\Big[x^{3/2}+\Vert \rhorho(t) \Vert_{L^{\upinfty}_\alpha H^{k-1}_x}  (1+ x)^{4(k-1)}x  (\Vert \vv \Vert_{\mathcal{H}^{k,p}}+ x^{1/2} ) \\
    &\qquad \qquad \qquad \qquad \qquad \qquad \qquad \qquad \qquad+ C_\eta (1+ x)^{4(k-1)} x^2 \Big],
\end{align*}
for times $t>0$ such that 
\begin{align}\label{cond-BobMarley-comp}
    \int_0^t \Vert  U(\tau) \Vert_{H^{k-1}_x} \mathrm{d} \tau<K^{-1}\ln(2),
\end{align}
because the factor in front of the former dissipative term is then nonnegative. Now, if we set $y(t)=(x+\delta)^{1/2}$ for $\delta>0$, we have $x(t)^{1/2} \leq y(t)$ and it is enough to bound $y$ and letting $\delta \rightarrow 0$ at the end of the proof. As long as \eqref{cond-BobMarley-comp} holds, we have
\begin{align*}
2 y y'&= x' \leq  K \Big(y^3+\Vert \rhorho(t) \Vert_{L^{\upinfty}_\alpha H^{k-1}_x} (1+ y^2)^{4(k-1)} y^2\left(\Vert \vv \Vert_{\mathcal{H}^{k,p}}+ y \right) \\
& \qquad  \qquad \qquad \qquad  \qquad \qquad \qquad  \qquad  + C_\eta (1+ y^2)^{4(k-1)}y^4 \Big),
\end{align*}
and therefore after simplification by $y$ we end up with 
\begin{align*}
    y' &\leq 
    K \Big(y^2+\Vert \rhorho(t) \Vert_{L^{\upinfty}_\alpha H^{k-1}_x}  (1+ y^2)^{4(k-1)} y\left(\Vert \vv \Vert_{\mathcal{H}^{k,p}}+ y \right) \\
    &\qquad \qquad \qquad \qquad \qquad \qquad \qquad \qquad+ C_\eta (1+ y^2)^{4(k-1)} y^3 \Big).
\end{align*}
By a comparison principle, it is enough to look at  $f(t)$ solution to the differential equation:
\begin{equation*}
\left\{ 
\begin{aligned}
f'(t) &= K f(t) \Big( f(t) + \Vert \rhorho(t) \Vert_{L^{\upinfty}_\alpha H^{k-1}_x} (1+f(t)^2)^{4(k-1)} (\Vert \vv(t) \Vert_{\mathcal{H}^{k,p}}+f(t)) \\
&  \qquad \qquad \qquad \qquad    \qquad    +C_\eta (1+f(t)^2)^{4(k-1)} f(t)^2 \Big)
 \\
 f(0) &= y(0),
 \end{aligned} \right.
\end{equation*}
on its maximal interval of existence, where it is in particular non-decreasing. Note that, if \eqref{cond-BobMarley-comp} holds, and by the definition of $x$ and $y$, we also have 
\begin{align}\label{cond-BobMarley-comp2}
    \int_0^t \Vert U(s) \Vert_{H^{k-1}_x} \, \mathrm{d}s \leq t f(t),
\end{align}
since $f$ is nondecreasing. Setting $A_g=(1+g^2)^{4(k-1)}$ and $B_g=C_\eta (1+g^2)^{4(k-1)} g^2$, a direct integration of the previous ODE gives that, for all times $t \in [0, T_{\max} ]$ such that $f(t) \leq 2 f(0)$, we have 
\begin{multline*}
    f(t) \leq f(0) e^{K(2 f(0)+B_{2f(0)})t} \\
    \times \exp\left( K \Vert \rhorho \Vert_{L^\infty([0,t];L^\infty_\alpha H^{k-1}_x)} A_{2f(0)}\left[  \int_0^t \Vert \vv(s) \Vert_{\mathcal{H}^{k,p}} \, \mathrm{d}s+2t f(0) \right]\right).
    \end{multline*} 
In view of \eqref{cond-BobMarley-comp}--\eqref{cond-BobMarley-comp2} and the previous computation, we finally set for $M,N\geq0$, $P=(P_1,P_2)$, $P_i\geq 0$,
\begin{align*}
    \Lambda(t,M,N,P)=e^{K(2(P_1+P_2)+B_{2(P_1+P_2)})t} \exp\left( K N A_{2(P_1+P_2)} \left[  M+2t (P_1+P_2) \right]\right),
\end{align*}
followed by 
 \begin{equation*}
    \F(t,M,N,P) \vcentcolon= \left\{ \begin{aligned}
    &(P_1+P_2)\Lambda(t,M,N,P), &&\mbox{if }  \Lambda(t,M,N,P)<2 \quad \text{and} \ \\ &\ && t(P_1+P_2)\Lambda(t,M,N,P)  \\ &\ &&\qquad \qquad < K^{-1}\ln(2), \\
    & + \infty, && \mbox{else},
    \end{aligned} \right.
\end{equation*}
we end up with
\begin{align*}
    \int_0^T \Vert U(s) \Vert_{H^{k-1}_x} \, \mathrm{d}s \leq  T \F\left(T,\int_0^T \Vert \vv(s) \Vert_{\mathcal{H}^{k,p}} \, \mathrm{d}s, \Vert \rhorho \Vert_{L^\infty([0,t];L^\infty_\alpha H^{k-1}_x)} ,N(X)\right), 
\end{align*}
where $N(X)=( \|\varrho_0\|_{H^{k-1}_x}, \|U_0\|_{H^{k-1}_x})$. This shows that an estimate of the type \ref{ass:bound_high_regularity-rho} is satisfied.

\paragraph{Stability at low regularity.} We now derive  a stability estimate of type~\ref{ass:stability-rho} for the compressible VNS system. Once again, as in the incompressible case, this question reduces to find a $L^1_T L^2_x$ estimate for $U_2 - U_1$, where $U_1$ and $U_2$ are solutions to
\begin{equation*}
\mbox{for }i=1,2,\quad \left\{
\begin{aligned}
\partial_t \varrho_i+\mathrm{div}(U_i \varrho_i)&=0, \\
\varrho_i\left(\partial_t U_i + (U_i \cdot \nabla)U_i \right) + \nabla P(\varrho_i) - \mathcal{A}U_i &= \int \rho_i^{\alpha} \left(v_i^{\alpha} - U_i \right) \D \mu(\alpha),\\
{\varrho_i|}_{t=0} = \varrho_0, \ \ {u_i|}_{t=0} &= U_0,
\end{aligned}
\right. 
\end{equation*}
where $\vv_1,\vv_2$ are two given families of velocity fields satisfying \eqref{eq:uniform_bound_velocities} and $\rhorho_1,\rhorho_2$ are the corresponding solutions to the continuity equation. The differences $R=\varrho_2-\varrho_1$ and $W=U_2-U_1$ satisfy
\begin{align*}
\partial_t R+\mathrm{div}(R U_2)&=-\mathrm{div}(\varrho_1 W), \\
\partial_t W  -\frac{1}{\varrho_2} \mathcal{A} W 
 &=\sum_{i=1}^5 \mathfrak{S}_k,
\end{align*}
where
\begin{align*}
    &\mathfrak{S}_1=-(U_2 \cdot \nabla)W, \ \  \mathfrak{S}_2=-(W \cdot \nabla)U_1, \\
    &\mathfrak{S}_3=\left(\frac{1}{\varrho_2}-\frac{1}{\varrho_1} \right) \mathcal{A}U_1, \ \ \mathfrak{S}_4=-\left( \nabla \pi(\varrho_2)-\nabla \pi (\varrho_1)\right), \\
    &\mathfrak{S}_5=\frac{1}{\varrho_2}\int \rho_2^{\alpha} \left(v_2^{\alpha} - U_2 \right) \, \mathrm{d}\mu(\alpha)-\frac{1}{\varrho_1}\int\rho_1^{\alpha} \left(v_1^{\alpha} - U_1 \right) \D \mu(\alpha).
\end{align*}
In view of the desired estimate~\eqref{eq:stability_assumption} and the bound at high regularity from the previous subsection, we will denote by $C_A>0$ a generic constant changing from line to line and depending on an upper bound $A>0$ of the right-hand side of \eqref{eq:stability_assumption}. As for the proof of the bounds at high regularity, we can always assume that 
\begin{align*}
    \inf \varrho_1(t)>0, \ \ \inf \varrho_2(t)>0.
\end{align*}
We shall now derive an estimate on 
$$\mathrm{Z}\vcentcolon=\Vert R \Vert_{L^2_x}^2+\Vert W \Vert_{L^2_x}^2.$$ 
First, thanks to an integration by parts, the equation on $R$ yields 
\begin{align*}
    \frac{1}{2}\frac{\mathrm{d}}{\mathrm{d}t} \Vert R \Vert^2_{L^2_x}&=-\int \mathrm{div}(R U_2)R-\int \mathrm{div}(\varrho_1 W)R \\
    &=\frac{1}{2}\int \mathrm{div}(U_2) \vert R \vert^2-\int (\varrho_1 \mathrm{div}(W)+\nabla \varrho_1 \cdot W) R\\
    & \lesssim \Vert \nabla U_2 \Vert_{L^{\upinfty}_x}\Vert R \Vert_{L^2_x}^2+\Vert  \varrho_1 \Vert_{W^{1, \infty}}(\Vert R \Vert_{L^2_x}^2+\Vert W \Vert_{L^2_x}^2)+ \frac{\delta}{2} \Vert \mathrm{div}(W) \Vert_{L^2_x}^2,
\end{align*}
for all $\delta>0$. Next, we have the following direct estimate on $W$:
\begin{multline*}
    \frac{1}{2}\frac{\mathrm{d}}{\mathrm{d}t} \Vert W \Vert^2_{L^2_x}+\int \left(\frac{1}{\varrho_2}-\eta \left\Vert \nabla  \frac{1}{\varrho_2} \right\Vert_{L^{\upinfty}_x} \right)\left( \vert \nabla W \vert^2+\vert \mathrm{div}(W) \vert^2\right) \\
    \leq \sum_{k=1}^5 \int \mathfrak{S}_k \cdot W \, \mathrm{d}x  +\eta^{-1} \Vert W \Vert^2_{L^2_x}, 
\end{multline*}
for all $\eta>0$. We now estimate each term from the sum in the right-hand side. For $\mathfrak{S}_1$ and $\mathfrak{S}_2$, there holds
\begin{align*}
    \int \mathfrak{S}_1 \cdot W \, \mathrm{d}x \lesssim \Vert \nabla W_2 \Vert_{L^{\upinfty}_x} \Vert W \Vert_{L^2_x}^2, \ \ \int \mathfrak{S}_2 \cdot W \, \mathrm{d}x \lesssim \Vert \nabla U_1 \Vert_{L^{\upinfty}_x} \Vert W \Vert_{L^2_x}^2.
\end{align*}
For $\mathfrak{S}_3$, we write
\begin{align*}
    \int \mathfrak{S}_3 \cdot W \, \mathrm{d}x &=\int \frac{\varrho_1-\varrho_2}{\varrho_1 \varrho_2}  \mathcal{A}U_1 \cdot W \, \mathrm{d}x  
    \leq \Vert R \Vert_{L^2_x} \left\Vert\frac{1}{\varrho_1 \varrho_2}  \mathcal{A}U_1 \cdot W\right\Vert_{L^2_x}\\
    & \leq \left\Vert \frac{\mathcal{A}U_1}{\varrho_1 \varrho_2} \right\Vert_{L^{\upinfty}_x}^2 \Vert R \Vert_{L^2_x}^2+ \Vert W \Vert_{L^2_x}^2 \\
    & \leq C_A \Vert R \Vert_{L^2_x}^2+ \Vert W \Vert_{L^2_x}^2,
\end{align*}
while for $\mathfrak{S}_4$, we have by the mean-value theorem 
\begin{align*}
\int \mathfrak{S}_4 \cdot W \, \mathrm{d}x \leq C_A(\delta) \Vert R \Vert_{L^2_x}^2+\frac{\delta}{2} \Vert \nabla W \Vert_{L^2_x}^2,
\end{align*}
for all $\delta>0$. For the last term $\mathfrak{S}_5$, we first write  
\begin{align*}
    \mathfrak{S}_5&=\mathfrak{S}_{5,1}+ \mathfrak{S}_{5,2} \\
    &\vcentcolon=\left( \frac{1}{\varrho_2}-\frac{1}{\varrho_1} \right)\int \rho_2^{\alpha} \left(v_2^{\alpha} - U_2 \right) \, \mathrm{d}\mu(\alpha)+\frac{1}{\varrho_1}\int \left[\mathrm{A}^\alpha+\mathrm{B}^\alpha + \mathrm{C}^\alpha \right] \D \mu(\alpha),
\end{align*}
where $\mathrm{A}^\alpha,\mathrm{B}^\alpha$ and $\mathrm{C}^\alpha$ are defined in \eqref{def:BrinkmanstabABC} in the incompressible case. For the first term $\mathfrak{S}_{5,1}$, we have thanks to composition estimates from Proposition \ref{prop-Sobolev}
\begin{align*}
    \int \mathfrak{S}_{5,1} \cdot W \, \mathrm{d}x &\leq \left\Vert \frac{1}{\varrho_1 \varrho_2} \int \rho_2^{\alpha} \left(v_2^{\alpha} - U_2 \right) \, \mathrm{d}\mu(\alpha)\right\Vert_{L^{\upinfty}_x}^2\Vert R \Vert_{L^2_x}^2+ \Vert W \Vert_{L^2_x}^2 \\
    &\leq C_A \Vert R \Vert_{L^2_x}^2+ \Vert W \Vert_{L^2_x}^2.
    \end{align*}
    For the second term $\mathfrak{S}_{5,2}$, this can be handled, assuming $p \geq 2$, by using straightforward modifications of the strategy from the incompressible case, with the additional prefactor $\Vert 1/ \varrho_1 \Vert_{\mathrm{W}^{1, \infty}}$, but this term is harmless by composition estimates and the bound at high regularity.

All in all, relying on the same estimates used at the end of the argument for the incompressible case, we obtain for all $\eta, \delta>0$
\begin{align*}
    &\frac{\mathrm{d}}{\mathrm{d}t} \mathrm{Z}+\int \left(\frac{1}{\varrho_2}-\eta \left\Vert \nabla  \frac{1}{\varrho_2} \right\Vert_{L^{\upinfty}_x}-\delta  \right)\left( \vert \nabla W \vert^2+\vert \mathrm{div}(W) \vert^2\right) \\
    &\lesssim_{\eta, \delta} (1+C_A)\left(\mathrm{Z}+ \int \Vert \rho_1^{\alpha}-\rho_2^{\alpha} \Vert_{\dot H^{-1}_x}^2 \D \mu(\alpha) +\int \Vert v_2^{\alpha}-v_1^{\alpha} \Vert_{ L^2_x}^2 \D \mu(\alpha) \right).
\end{align*}
Eventually, thanks to the bounds at high regularity, a suitable choice of small enough $\delta>0$ and $\eta>0$ with respect to $C_A$ allows one to ensure that the second term in the previous estimate is nonnegative. Thanks to Gronwall lemma, we therefore obtain that Assumption~\ref{ass:stability-rho} is satisfied, which concludes the proof.

\end{proof}

\section{Local well-posedness for the Vlasov--Stokes equations}\label{subsec:Stokes}

We will now apply the general framework of Chapter \ref{Part1-LWP} to the multiphasic Vlasov--Stokes system \eqref{eq:VS-multiphase}, which is set in the Euclidean space $\R^3$. Let us emphasis the fact that the Stokes system
\begin{align*}
    - \Delta U+\nabla P  &= \int \big( v^\alpha - U ) \rho^\alpha \D \mu(\alpha)
,\\
\Div(U) &= 0,
\end{align*}
only makes sense for $x \in \R^3$ (the above source term having no reason to be average-free in general, the former system therefore cannot be considered on the torus $\T^3$). Consequently, since we work on the whole space and have no direct control on the low-frequency part of the solution $U$, we do not expect to have $U \in L^2(\R^3)$ (the study of the general fundamental solution for the Stokes equations rather produces the so-called \textit{Stokeslet}, that has a tail behaving like $\vert x \vert^{-1}$ -- see e.g. \cite[Chapter 4]{galdi2000introduction}): the general theory for the Stokes system naturally calls for homogeneous Sobolev spaces, typically with $U\in\dot H^1(\R^3)$, as in \cite{HoferPHD,HS-VSmeanfield2}. This framework is particularly well suited to the study of weak solutions.

However, it is not clear that the general theory developed in Chapter \ref{Part1-LWP}, which is used to construct a high-regularity solution $(\rhorho,\vv)$ (being of finite energy) coupled through the force field $U$, is compatible with such homogeneous spaces. Indeed, it requires the intensive use of product estimates and embeddings that are not always behaving well with homogeneous spaces. Instead, it is more natural in the present setting to work with weighted Sobolev spaces with varying weights, as introduced in Section \ref{subsec-weightedSob} from Chapter \ref{Part1-LWP}. Below, we directly use the notation set forward there. This framework will allow us to produce good estimates the Stokes operator above, directly at high regularity. We refer to Lemma \ref{LM:StokesWeightElliptic} below for a more precise statement and some simple heuristics about the use of the weighted Sobolev spaces $H^k_\delta(\R^3)$ for the Stokes problem.

\medskip

Our main result reads as follows.
\begin{Thm}\label{thm:appliVS}
Let $k \in \N$ such that $k>1+3/2$, $p \in [2, \infty]$ and $\delta \in (-3/2,-1/2)$. Suppose that the initial distribution function $f_0$ can be written as
$$f_0(x, \cdot)=\int_{I} \rho^\alpha_0(x) \otimes  \delta_{v=v^\alpha_0(x)} \, \mathrm{d}\mu(\alpha),$$
for some set of labels $(I, \mu)$ and some regular multiphasic distribution $(\rhorho_0, \vv_0)\in L^\infty_\alpha H^{k-1}_{\delta} \times \mathcal{H}^{k,p}_\delta $, where $\rho^\alpha \geq 0$, $\rhorho_0 \in L^\infty_\alpha H^{k-1}_{-\delta-3/2}$ and $\langle x \rangle^{-\delta}\rhorho_0 \in L^\infty_\alpha L^3_x$. Then there exist $T>0$ and a unique solution $(\rhorho, \vv, U)$ on $[0,T]$ to the multiphasic system \eqref{eq:VS-multiphase} with
    \begin{align*}
    \rhorho \in L^\infty([0,T]; L^\infty_\alpha (H^{k-1}_\delta)  ), \ \ \vv \in L^\infty([0,T]; \mathcal{H}^{k,p}_\delta).
\end{align*}
Moreover,
$$
f(t,x,\cdot) \vcentcolon= \int_{I} \rho^\alpha(t,x) \otimes  \delta_{v=v^\alpha(t,x)} \, \mathrm{d}\mu(\alpha)
$$
is the unique weak solution on $[0,T]$ to \eqref{eq:VSkin} with initial condition $f_0$, in the sense of Definition \ref{def:solution_coupled}. 
\end{Thm}

\begin{Rem}
    As for the incompressible and compressible cases, one can also include the case of a smooth nonlinear drag term in \eqref{eq:VS-multiphase}, taking $p=\infty$ (see Section \ref{sec:nonlindrag-LWP}). Details are left to the reader.
\end{Rem}

\begin{Rem}
    The condition $\langle x \rangle^{-\delta}\rhorho_0 \in L^\infty_\alpha L^3_x$ is actually the most stringent among the three polynomial-decay requirements on the densities in Theorem \ref{thm:appliVS}. The constraint on the initial data is not empty and is for instance satisfied if $\rho^\alpha \sim \langle x \rangle^{p}$ with $p>1-\delta$, uniformly in the label $\alpha$. The need for such constraint seems mostly technical, see the proof of the stability estimates below, but we did not try to relax it.  
\end{Rem}
\begin{Cor}
    Let $k \in \N$ such that $k>1+3/2$ and $\delta \in (-3/2,-1/2)$. Let $m \in L^1(\langle v \rangle^2)$ and $f_0 \in  L^\infty_{v,m^{-1}} H^k_\delta  $ with $f_0 \geq 0$,  $f_0 \in L^\infty_{v,m^{-1}} H^{k-1}_{-\delta-3/2}$ and $\langle x \rangle^{-\delta}f_0 \in L^\infty_{v,m^{-1}} L^3_x$.  There exist $T>0$ and a unique weak solution to the system~\eqref{eq:VSkin} on $[0,T]$, with initial condition $f_0$, in the sense of Definition \ref{def:solution_coupled}.
\end{Cor}
Thanks to Theorem~\ref{thm-LWPmultiphase-densityvelocity-WEIGHTED} in Chapter~\ref{Part1-LWP}, to prove Theorem~\ref{thm:appliVS}
 it is enough to check that the Assumptions~\ref{ass:bound_high_regularity-rho-weighted} and~\ref{ass:stability-rho-weighted}  are satisfied for the system \eqref{eq:VS-multiphase}. 

We start with some preliminary estimates (without weights to begin with) on the Stokes equations in dimension 3.
\begin{Lem}
Let $\vv \in L^{\infty}_T\mathcal{H}^{k,p}$ be a family of velocity field and $\rhorho$ be the associated nonnegative solution to the continuity equation. There exists a unique weak solution $U \in  L^\infty_T \dot H^1(\R^3)$ to the Stokes problem
$$- \Delta U+\nabla P  = \int \big( v^\alpha - U ) \rho^\alpha \D \mu(\alpha)
,\qquad
\Div(U) = 0.$$
Furthermore, we have the following pointwise in time estimates:
\begin{align}
\label{estim-Stokes1}
    \Vert \nabla U \Vert_{L^2_x} &\lesssim \int \Vert \rho^\alpha v^\alpha   \Vert_{L^{6/5}_x} \D \mu(\alpha), \\
\label{estim-Stokes2}\Vert  \nabla U \Vert_{L^3_x} &\lesssim \int \Vert \rho^\alpha v^\alpha \Vert_{L^{3/2}_x} \D \mu(\alpha)+ \left( \int \Vert \rho^\alpha v^\alpha   \Vert_{L^{6/5}_x} \D \mu(\alpha)  \right) \int \Vert \rho^\alpha  \Vert_{L^2_x} \D \mu(\alpha), \\
\label{estim-Stokes3}    \Vert  U \Vert_{L^{\upinfty}_x} &\lesssim \int \Vert \rho^\alpha v^\alpha \Vert_{L^2_x} \D \mu(\alpha)+ \left( \int \Vert \rho^\alpha v^\alpha   \Vert_{L^{6/5}_x} \D \mu(\alpha)  \right) \int \Vert \rho^\alpha  \Vert_{L^3_x} \D \mu(\alpha).
\end{align}
\end{Lem}
\begin{proof}
    The proof is mainly extracted from \cite{HS-VSmeanfield2}, with an adaptation to the multiphasic framework. The well-posedness part and the first bound comes from rewriting of the Stokes problem as 
    \begin{align*}
        - \Delta U+ \left(\int \rho^\alpha \D \mu(\alpha) \right)U + \nabla P  = \int \rho^\alpha v^\alpha   \D \mu(\alpha)
,\qquad 
\Div(U) = 0,
    \end{align*}
    from which we infer the \textit{a priori} estimate
        \begin{align*}
    \Vert \nabla U \Vert_{L^2_x}^2 \leq \int_{\R^3} \int \rho^\alpha v^\alpha \cdot U \D \mu(\alpha)  \, \mathrm{d}x \lesssim \Vert \nabla U \Vert_{L^2_x}\int \Vert \rho^\alpha v^\alpha   \Vert_{L^{6/5}_x} \D \mu(\alpha),
    \end{align*}
    thanks to the divergence-free condition, positivity of the densities, Hölder inequality and Sobolev embedding in dimension 3. It allows to apply the Lax-Milgram theorem in $\dot H^1_x$ and to obtain the pointwise estimate \eqref{estim-Stokes1}. For the second estimate, we first use the elliptic regularity for the Stokes problem with the source term $\int \big( v^\alpha - U ) \rho^\alpha \D \mu(\alpha)$ combined with Hölder inequality and Sobolev embedding to get for all $p \in (1, 6]$ and $q$ such that $\frac{1}{p}=\frac{1}{6}+\frac{1}{q}$
    \begin{align*}
        \Vert \mathrm{D}^2 U \Vert_{L^p_x} 
        &\lesssim \int \Vert \rho^\alpha v^\alpha \Vert_{L^p_x} \D \mu(\alpha)+ \Vert U \Vert_{L^6_x} \int \Vert \rho^\alpha  \Vert_{L^q_x} \D \mu(\alpha) \\
        &\lesssim \int \Vert \rho^\alpha v^\alpha \Vert_{L^p_x} \D \mu(\alpha)+ \left( \int \Vert \rho^\alpha v^\alpha   \Vert_{L^{6/5}_x} \D \mu(\alpha)  \right) \int \Vert \rho^\alpha  \Vert_{L^q_x} \D \mu(\alpha),
    \end{align*}
    Here we have used the previous estimate \eqref{estim-Stokes1}. With $p=3/2$ and $q=2$, and using the Sobolev inequality $\Vert \nabla U \Vert_{L^3_x} \lesssim \Vert \mathrm{D}^2 U \Vert_{L^{3/2}_x}$, we get the estimate \eqref{estim-Stokes2}. With $p=2$ and $q=3$, we also obtain $U \in \dot H^2_x$ from which we deduce $U \in L^\infty_x$ with the associated estimate \eqref{estim-Stokes3}. This concludes the proof.
\end{proof}

Weighted elliptic regularity for the Stokes system in dimension $3$ reads as follows (see \cite{Stokesweight}). 
\begin{Lem}\label{LM:StokesWeightElliptic}
    Let $\delta \in (-3/2, -1/2), k \geq 2$ and $a \in H^k_{2}(\R^3)$ with $a \geq 0$. For any $J \in H^k_{\delta+2}(\R^3)$, there exists a unique $(U,P) \in H^{k+2}_\delta(\R^3) \times H^{k+1}_\delta(\R^3)$ such that 
    \begin{align*}
        - \Delta U+a U+\nabla P  = J
,\\
\Div(U) = 0.
    \end{align*}
Moreover, we have 
\begin{align*}\Vert U \Vert_{H^{k+2}_\delta} \leq C(\Vert a \Vert_{H^k_{2}}) \Vert J \Vert_{H^k_{\delta+2}},
\end{align*}
where $C: \R^+ \rightarrow \R^+$ is a continuous nondecreasing function.
\end{Lem}

\begin{Rem}
    As mentioned earlier, it might be expected that the solution $U$ to the Stokes system on $\R^3$ behaves likes $U (x)\sim \vert x \vert^{-1}$ at infinity (see \cite{galdi2000introduction}), so that $\partial^k_x U(x) \sim \vert x \vert^{-1-k}$. Therefore 
    $\langle x \rangle^{2(k+\delta)}\vert \partial^k_x U(x)\vert^2 \sim \vert x \vert^{2(k+\delta)}\vert x \vert^{-2-2k} =\vert x \vert^{2\delta-2}$ hence such $U(x)$ satisfies $U \in H^k_\delta $ if and only if $\delta<-1/2$. On the other hand, we exclude the fact that the solution $U$ behaves like a non-zero constant at infinity, since $\langle x \rangle^\delta \in L^2_x$ if and only $\delta<-3/2$. This heuristics roughly explains the numerology in the weight parameter $\delta$ in Lemma \ref{LM:StokesWeightElliptic}.
\end{Rem}

\begin{proof}[Proof of Theorem \ref{thm:appliVS}]
We now aim at checking that the two main Assumptions \ref{ass:bound_high_regularity-rho-weighted} and \ref{ass:stability-rho-weighted} hold for \eqref{eq:VS-multiphase} (see Section \ref{subsec-weightedSob}).
Similarly to the Vlasov--Navier--Stokes case, the force field of interest is the fluid velocity field $U$ itself, solution to the Stokes equations.

\paragraph{Bound in high regularity.}
We first aim at deriving a suitable $L^1_T H^{k}_\delta $ estimate for $U$ solution to the Stokes equations for $\delta>-3/2$. We write the Stokes equation as
$$- \Delta U+ \left(\int  \rho^\alpha \D \mu(\alpha) \right)U+\nabla P   = \int \rho^\alpha v^\alpha   \D \mu(\alpha).$$
Now appealing to Lemma \ref{LM:StokesWeightElliptic} with
\begin{align*}
    a\vcentcolon = \int  \rho^\alpha \D \mu(\alpha), \ \ J\vcentcolon =\int \rho^\alpha v^\alpha   \D \mu(\alpha),
\end{align*}
we infer that the Stokes problem can be solved for $k$ large enough and $\delta \in (-3/2,-1/2)$ with 
\begin{align*}
\Vert U \Vert_{H^{k}_\delta} &\leq  C(\Vert \rho^\alpha \Vert_{L^{\upinfty}_\alpha H^{k-2}_{2}}) \int \Vert \rho^\alpha v^\alpha \Vert_{H^{k-2}_{\delta+2}} \, \mathrm{d}\mu(\alpha) \\
& \leq  C(\Vert \rho^\alpha \Vert_{L^{\upinfty}_\alpha H^{k-2}_{2}})\Vert \rho^\alpha \Vert_{L^{\upinfty}_\alpha H^{k-1}_{2}} \int \Vert v^\alpha \Vert_{H^{k}_{\delta}} \, \mathrm{d}\mu(\alpha),
\end{align*}
in view of the product law from Proposition \ref{prop-Sobweight}.
This shows that Assumption \ref{ass:bound_high_regularity-rho-weighted} is satisfied and now concludes the proof.

\paragraph{Stability at low regularity.}
We consider
\begin{equation*}
\mbox{for }i=1,2,\qquad \left\{
\begin{gathered}
- \Delta U_i+\nabla P_i  = \int \rho_i^{\alpha} \left(v_i^{\alpha} - U_i \right) \D \mu(\alpha),\\
\Div(U_i) = 0,
\end{gathered}
\right. 
\end{equation*}
where $\vv_1,\vv_2$ are two given families of velocity fields satisfying~\eqref{eq:uniform_bound_velocities-weighted} with an associated $R>0$, and $\rhorho_1,\rhorho_2$ are the corresponding solutions to the continuity equation. We observe that the difference $W \vcentcolon = U_2 - U_1$ satisfies
\begin{align*}
-  \Delta W +\nabla P &=\int \left[\rho_2^{\alpha} \left(v_2^{\alpha} - U_2 \right)-\rho_1^{\alpha} \left(v_1^{\alpha} - U_1 \right) \right] \D \mu(\alpha), \\
\mathrm{div} (W)&=0,
\end{align*}
and we can also rewrite the first line as
\begin{multline*}
    -\Delta W + \nabla P+ \left(\int \rho_2^{\alpha} \, \D \mu(\alpha) \right)W \\
    =\int \left[\rho_2^{\alpha}(v_2^{\alpha}-v_1^{\alpha})+(\rho_2^{\alpha}-\rho_1^{\alpha})(v_1^{\alpha}-U_1) \right] \D \mu(\alpha).
\end{multline*}
Now observe that according to the Poincaré inequality of Proposition \ref{prop-Sobweight}, there holds for all $\delta  \leq -1$
\begin{align}\label{ineq:HölderWeights}
 \Vert W \Vert_{L^2_\delta} \lesssim \Vert \nabla W \Vert_{L^2_x}.
\end{align}
 We can therefore focus on estimating $\Vert \nabla W \Vert_{L^2_x}$:
testing with $W$, and thanks to the divergence-free condition, we obtain
\begin{multline*}
    \Vert \nabla W \Vert_{L^2_x}^2+ \int \left(\int \rho_2^{\alpha} \, \D \mu(\alpha) \right)  \vert W \vert^2 \, \mathrm{d}x 
    \\
    =\int \left\langle \rho_2^{\alpha}(v_2^{\alpha}-v_1^{\alpha})+(\rho_2^{\alpha}-\rho_1^{\alpha})(v_1^{\alpha}-U_1),  W  \right\rangle_{L^2_x} \D \mu(\alpha),
\end{multline*}
from which we infer,  as the densities $\rho_2^{\alpha}$ can always be assumed to be nonnegative, that
\begin{multline}\label{eq:pikachuStokes}
    \Vert \nabla W \Vert_{L^2_x}^2
     \leq \int \Vert \langle x \rangle^{-\delta}\rho_2^{\alpha} \Vert_{L^3_x} \Vert \langle x \rangle^{\delta} (v_2^{\alpha}-v_1^{\alpha}) \Vert_{L^2_x} \Vert   W \Vert_{L^6_x}  \D \mu(\alpha) \\
      + \int \Vert (\rho_2^{\alpha}-\rho_1^{\alpha}) \Vert_{\dot H^{-1}_x} \Vert \nabla ((v_1^{\alpha}-U_1)\cdot  W) \Vert_{L^2_x}   \D \mu(\alpha).
\end{multline}
\textbf{First term in \eqref{eq:pikachuStokes}}: For the first term, we use Sobolev embedding  
to obtain the control 
\begin{multline*}
\int \Vert \langle x \rangle^{-\delta}\rho_2^{\alpha} \Vert_{L^3_x} \Vert \langle x \rangle^{\delta} (v_2^{\alpha}-v_1^{\alpha}) \Vert_{L^2_x} \Vert   W \Vert_{L^6_x}  \D \mu(\alpha) \\
\lesssim \Vert  \nabla W \Vert_{L^2_x} \int \Vert \langle x \rangle^{-\delta}\rho_2^{\alpha} \Vert_{L^3_x} \Vert v_2^{\alpha}-v_1^{\alpha} \Vert_{L^2_\delta}   \D \mu(\alpha).
\end{multline*}
The factor $\Vert \langle x \rangle^{-\delta}\rho_2^{\alpha} \Vert_{L^3_x}$ in the integral is handled thanks to the method of characteristics. Indeed, using the representation formula for the density
\begin{align*}
    \rho^{\alpha}_2(t,x)=\rho^{\alpha}_{2,0}(\mathrm{X}^{0;t}(x))\exp\left(-\int_0^t \mathrm{div}(v_2^{\alpha})(\tau, \mathrm{X}^{\tau;t}(x)) \,  \mathrm{d} \tau\right),
\end{align*}
where $\mathrm{X}^{t;s}$ is the characteristics associated to the vector field $v_2^\alpha$, that is  
$$\frac{\mathrm{d}}{\mathrm{d}s}\mathrm{X}^{t;s}(x)=v^{\alpha}_2(s,\mathrm{X}^{t;s}(x)), \qquad \mathrm{X}^{s;s}(x)=x,$$ 
we obtain
\begin{align*}
    \Vert \langle x \rangle^{-\delta}\rho_2^{\alpha}(t) \Vert_{L^3_x}^3=\int \langle x \rangle^{-3\delta} \vert \rho^{\alpha}_{2,0}(\mathrm{X}^{0;t}(x)) \vert^3 e^{-3\int_0^t \mathrm{div}(v_2^{\alpha})(\tau, \mathrm{X}^{\tau;t}(x))  \, \mathrm{d} \tau} \, \mathrm{d} x.
\end{align*}
Changing variable along the flow, and then integrating the ODE in time, we infer (since $\delta<0$) that for all $0 \leq t \leq T$, 
\begin{align*}
     &\Vert \langle x \rangle^{-\delta}\rho_2^{\alpha}(t) \Vert_{L^3_x}^3 \\
     &\lesssim e^{2T\Vert \nabla v_2^{\alpha} \Vert_{L^{\infty}_T L^{\upinfty}_x}} \int \langle \mathrm{X}^{t;0}(x) \rangle^{-3 \delta} \vert \rho^{\alpha}_{2,0}(x)\vert^3 \, \mathrm{d}x    \\
     & \lesssim e^{2T\Vert \nabla v_2^{\alpha} \Vert_{L^{\infty}_T L^{\upinfty}_x}} \\
     &\qquad \times\left( \int \langle x \rangle^{-3 \delta} \vert \rho^{\alpha}_{2,0}(x)\vert^3 \, \mathrm{d}x+\int \left\langle \int_0^T v_2^\alpha(\mathrm{X}^{\tau;0}(x)) \, \mathrm{d}\tau  \right\rangle^{-3 \delta} \vert \rho^{\alpha}_{2,0}(x)\vert^3 \, \mathrm{d}x \right) \\
     & \lesssim e^{2T\Vert \nabla v_2^{\alpha} \Vert_{L^{\infty}_T L^{\upinfty}_x}} \left(1+T^{-3\delta} \Vert v_2^{\alpha} \Vert_{L^{\infty}_T L^{\upinfty}_x}^{-3\delta} \right)  \Vert \langle x \rangle^{-\delta} {\rho_2^{\alpha}}_0 \Vert_{L^3_x}^3,
\end{align*}
and we therefore deduce the following control (for $\delta<0$): for all $t \in [0,T]$
\begin{align*}
    \Vert \langle x \rangle^{-\delta}\rho_2^{\alpha}(t) \Vert_{L^3_x} \lesssim e^{\frac{2}{3}T\Vert \nabla v_2^{\alpha} \Vert_{L^{\infty}_T L^{\upinfty}_x}} \left( 1+T^{-\delta} \Vert v_2^{\alpha} \Vert_{L^{\infty}_T L^{\upinfty}_x}^{-\delta} \right) \Vert \langle x \rangle^{-\delta} {\rho^{\alpha}_{2,0}} \Vert_{L^3_x}.
\end{align*}
Using the Sobolev embedding and Poincaré inequality, with $\delta>-3/2$ and $k>1+3/2$  from Proposition \ref{prop-Sobweight}, we know that $\Vert \nabla v_2^{\alpha} \Vert_{L^{\infty}_T L^{\upinfty}_x} + \Vert v_2^{\alpha} \Vert_{L^{\infty}_T L^{\upinfty}_x} \lesssim \Vert \vv_2 \Vert_{L^{\infty}_T \mathcal{H}^{k,p}_\delta}$.
This is enough to control the first term in \eqref{eq:pikachuStokes} (having $0>\delta>-3/2$) as
\begin{multline}\label{eq:pikachuStokes2}
    \int \Vert \langle x \rangle^{-\delta}\rho_2^{\alpha} \Vert_{L^3_x} \Vert \langle x \rangle^{\delta} (v_2^{\alpha}-v_1^{\alpha}) \Vert_{L^2_x} \Vert   W \Vert_{L^6_x}  \D \mu(\alpha) \\
    \lesssim C_{R,T} \Vert  \nabla W \Vert_{L^2_x} \int  \Vert v_2^{\alpha}-v_1^{\alpha} \Vert_{L^2_\delta}   \D \mu(\alpha),
\end{multline}
as we have assumed that $\vv_1,\vv_2$ satisfy the bound~\eqref{eq:uniform_bound_velocities-rho-weighted}.

\textbf{Second term in \eqref{eq:pikachuStokes}}: For the second term, we write thanks to Hölder's inequality and Sobolev embedding 
\begin{align*}
    \Vert \nabla ((v_1^{\alpha}-U_1)\cdot  W)& \Vert_{L^2_x} \\
    &\leq \left(\Vert v_1^{\alpha}\Vert_{L^{\upinfty}_x}+ \Vert U_1 \Vert_{L^{\upinfty}_x} \right) \Vert  \nabla W \Vert_{L^2_x}+ \Vert (\nabla v_1^{\alpha}- \nabla U_1)  W \Vert_{L^2_x}\\
    & \leq \left(\Vert v_1^{\alpha}\Vert_{L^{\upinfty}_x}+ \Vert U_1 \Vert_{L^{\upinfty}_x} \right) \Vert \nabla W \Vert_{L^2_x}+ \Vert \nabla v_1^{\alpha}- \nabla U_1 \Vert_{L^3_x} \Vert  W \Vert_{L^6_x} \\
    & \leq \left(\Vert v_1^{\alpha}\Vert_{L^{\upinfty}_x}+\Vert \nabla v_1^{\alpha}\Vert_{L^3_x}+ \Vert U_1 \Vert_{L^{\upinfty}_x}+ \Vert \nabla U_1 \Vert_{L^3_x} \right) \Vert \nabla  W \Vert_{L^2_x}.
\end{align*}
Let us study each term in the prefactor. The first one is estimated using Sobolev embedding, yielding $\Vert v_1^{\alpha}\Vert_{L^{\upinfty}_x} \lesssim  \Vert v_1^{\alpha}\Vert_{H^k_\delta}$,
while for the second one, we use an additional interpolation argument, as follows:
\begin{align*}
    \Vert \nabla v_1^{\alpha}\Vert_{L^3_x} 
    \lesssim 
    \Vert \nabla v_1^{\alpha}\Vert_{L^2_x}^{2/3}\Vert \nabla v_1^{\alpha}\Vert_{L^{\upinfty}_x}^{1/3}  
    \lesssim 
    \Vert \nabla v_1^{\alpha}\Vert_{H^1_{-1}}^{2/3}\Vert v_1^{\alpha}\Vert_{\mathscr{C}^1_0}^{1/3} \lesssim 
    \Vert \nabla v_1^{\alpha}\Vert_{H^{k-1}_{\delta+1 }}^{2/3}\Vert v_1^{\alpha}\Vert_{H^{k}_{\delta}}^{1/3}.
\end{align*}
Integrating in $\alpha$, and in view of the definition of the norm $\Vert \cdot \Vert_{\mathcal{H}^{k,p}_\delta}$ (recall \eqref{def:norm-weighted}) we can therefore bound the former contribution by $R$ (as we have assumed that $\vv_1,\vv_2$ satisfy the bound~\eqref{eq:uniform_bound_velocities-weighted}). For the other two ones, we apply the preliminary estimates \eqref{estim-Stokes2} and \eqref{estim-Stokes3} that give
\begin{align*}
    \Vert  \nabla U^1 \Vert_{L^3_x} 
    &\lesssim \int \Vert \rho_1^{\alpha} v_1^{\alpha} \Vert_{L^{3/2}_x} \D \mu(\alpha)+ \left( \int \Vert \rho_1^{\alpha} v_1^{\alpha}   \Vert_{L^{6/5}_x} \D \mu(\alpha)  \right) \int \Vert \rho_1^{\alpha}  \Vert_{L^2_x} \D \mu(\alpha) \\
    &\lesssim \int \Vert \rho_1^{\alpha} \Vert_{L^{3/2}_x} \Vert v_1^{\alpha} \Vert_{L^{\upinfty}_x} \D \mu(\alpha) \\
    &\qquad + \left( \int \Vert \rho_1^{\alpha} \Vert_{L^{6/5}_x} \Vert v_1^{\alpha} \Vert_{L^{\upinfty}_x} \D \mu(\alpha)  \right) \int \Vert \rho_1^{\alpha}  \Vert_{L^2_x} \D \mu(\alpha),
\end{align*}
and similarly
\begin{align*}
    \Vert  U^1 \Vert_{L^{\upinfty}_x} 
    &\lesssim \int \Vert \rho_1^{\alpha} \Vert_{L^{2}_x} \Vert v_1^{\alpha} \Vert_{L^{\upinfty}_x} \D \mu(\alpha) \\
    &\qquad+ \left( \int \Vert \rho_1^{\alpha} \Vert_{L^{6/5}_x} \Vert v_1^{\alpha} \Vert_{L^{\upinfty}_x} \D \mu(\alpha)  \right) \int \Vert \rho_1^{\alpha}  \Vert_{L^3_x} \D \mu(\alpha).
\end{align*}
Since all the $L^p$ norms of the densities $\rho_1^{\alpha}$ for $p \in (1, \infty)$ satisfy $\Vert \rho_1^\alpha(t) \Vert_{L^p_x} \lesssim e^{RT} \Vert \rho_{1, \mid t=0}^\alpha \Vert_{L^p_x}$ (see Lemma \ref{LM:pointwiseEstimRHO}) and $\Vert v_1^{\alpha} \Vert_{L^{\upinfty}_x} \lesssim \Vert v_1^{\alpha}\Vert_{H^k_\delta}$, we can therefore control these two types of terms by an expression depending on $R,T$ and the initial densities, and therefore infer that the second term in \eqref{eq:pikachuStokes} is bounded by
\begin{multline}\label{eq:pikachuStokes3}
    \int \Vert (\rho_2^{\alpha}-\rho_1^{\alpha}) \Vert_{\dot H^{-1}_x} \Vert \nabla ((v_1^{\alpha}-U_1)\cdot  W) \Vert_{L^2_x}   \D \mu(\alpha) \\
    \leq C_{R,T}\int \Vert (\rho_2^{\alpha}-\rho_1^{\alpha}) \Vert_{\dot H^{-1}_x}   \D \mu(\alpha),
\end{multline}
where $C_{R,T}>0$ implicitly depends on the initial densities. Gathering the two previous estimates \eqref{eq:pikachuStokes2}--\eqref{eq:pikachuStokes3}, and in view of \eqref{ineq:HölderWeights} and \eqref{eq:pikachuStokes}, there exists another constant $C_{R,T}$ of the same type such that
\begin{align*}
    \Vert W \Vert_{L^2_\delta}
     \leq C_{R,T} \left( \int  \Vert v_2^{\alpha}-v_1^{\alpha} \Vert_{L^2_\delta }   \D \mu(\alpha) + \int \Vert (\rho_2^{\alpha}-\rho_1^{\alpha}) \Vert_{\dot H^{-1}_x}   \D \mu(\alpha) \right),
\end{align*}
from which we deduce that Assumption~\ref{ass:stability-rho-weighted} is satisfied.
\end{proof}

\section[Stability of monokinetic profiles in Vlasov--Navier--Stokes  systems]{Stability of monokinetic profiles in Vlasov--Navier--Stokes  systems: main results}\label{Section:monokinetic-stab-theorem}

The second part of this chapter is dedicated to the study of the large time behaviour for the incompressible or compressible Vlasov--Navier--Stokes systems, set on the periodic torus $\T^3$. We aim at studying asymptotic stability of constant monokinetic profiles for these systems and, in doing so, to show that the well-posedness results from Section \ref{SubsecVNSmultiphaseLWP} and Section \ref{SubsecVNScomp-multiphaseLWP} extend to a global-in-time framework. We strive for the simple yet striking idea that, around simple states describing alignment  between the fluid and the particle velocities, the multiphasic framework can sustain global in time solutions as well as their long-time dynamics.

More precisely, observe that for a fixed density $\rho\in \R$ and a fixed  velocity $\mathcal V \in \R^3$, the distribution
\begin{equation}
    \label{eq:mono-VNS-2}
    f = \rho \otimes \delta_{v=\mathcal V}, \quad U = \mathcal V,
\end{equation}
which we refer to as the constant monokinetic profile, is a stationary solution to the incompressible Vlasov--Navier--Stokes system~\eqref{eq:VNS_kinLWP}. It represents one of the simplest equilibrium configurations, in which the particle cloud and the fluid move with the same constant velocity.
Very loosely speaking, we intend to prove the following type of results: 
\begin{center}
``\textit{Constant monokinetic profiles~\eqref{eq:mono-VNS-2} enjoy asymptotic stability properties for the Vlasov--Navier--Stokes equations~\eqref{eq:VNS_kinLWP}.}''
\end{center}

Although such stability statement was already obtained for the incompressible Vlasov--Navier--Stokes system (in its original kinetic formulation), recall Theorem~\ref{thm:VNSlargetime} of \cite{HKMM} in the introduction, we aim at developing an alternative approach based on the multiphasic formulation. 

The key observation is that the monokinetic equilibrium~\eqref{eq:mono-VNS-2} admits an immediate counterpart in the multiphasic formulation. Indeed, for any probability space $(I, \mu)$, the constant family
\begin{equation}
\label{eq:mono-VNS-multiphase}
\rho^\alpha =\rho, \qquad v^\alpha = \mathcal{V} , \qquad U \equiv \mathcal{V}, \qquad \alpha \in I,
\end{equation}
defines a stationary solution of the multiphasic Vlasov--Navier--Stokes system, since the drag force $v^\alpha-U$ identically vanishes. This elementary fact is the starting point of our analysis, allowing us to investigate the dynamics as a perturbation of a particularly simple equilibrium in which all phases are perfectly aligned with the fluid.

 Whereas the kinetic approach of \cite{HKMM} relies on the energy--dissipation method (see Theorem \ref{thm:VNSlargetime}), we shall obtain analogous results by a robust perturbative approach around the former constant solution \eqref{eq:mono-VNS-multiphase}, providing along the way  quantitative information on the rate of exponential convergence. 
Moreover, thanks to the multiphasic framework, the effect of friction forces becomes particularly transparent in the analysis, ultimately leading to global existence and to the alignment of all phases. The corresponding result for the tridimensional compressible Vlasov--Navier--Stokes system~\eqref{eq:VNS-comp_kinLWP} is new.

\medskip

For clarity of exposure, we now (re)set notation and the systems under study, and state the results concerning the incompressible/compressible case in separate subsections.

\subsection{The case of the incompressible Vlasov--Navier--Stokes system}
We start with the case of the incompressible Vlasov--Navier--Stokes system, in the multiphasic formulation:
\begin{equation}
\label{eq:VNS-CHAP4}
\left\{ 
\begin{gathered}
\partial_t \rho^\alpha + \Div ( \rho^\alpha v^\alpha) = 0,\\
\partial_t v^\alpha + (v^\alpha \cdot \nabla) v^\alpha = U - v^\alpha+ \Gamma(U-v^\alpha),\\
\partial_t U + (U\cdot \nabla ) U + \nabla P -  \Delta U = -\int \big( U - v^\alpha+ \Gamma(U-v^\alpha) ) \rho^\alpha \D \mu(\alpha)
,\\
\Div(U) = 0,
\end{gathered}
\right.
\end{equation}
where $(I, \mu)$ is a fixed probability space. Note that we include a possible smooth quadratic nonlinear drag $\Gamma(U-v^\alpha)$, where $\Gamma\in \mathscr{C}^{\infty}(\R^3;\R^3)$ is asked to satisfy $\Gamma(0)=0$ and $\Gamma'(0)=0$ -- see Section \ref{sec:nonlindrag-LWP}. The system~\eqref{eq:VNS-CHAP4} is endowed with initial data
\begin{align*}
    (\rhorho, \vv, U)_{\mid t=0}=(\rhorho_0, \vv_0, U_0).
\end{align*}

\medskip

As explained before, recalling \eqref{eq:mono-VNS-2} and \eqref{eq:mono-VNS-multiphase}, we will study this multiphasic system \eqref{eq:VNS-CHAP4} for initial data in the neighbourhood of a (normalized) stationary solution of the form
\begin{equation}
\label{eq:mono-VNS}
\rho^\alpha \equiv 1, \qquad v^\alpha \equiv \mathcal{V} , \qquad U \equiv \mathcal{V},
\end{equation}
where $\mathcal{V} \in \R^3$ is a constant vector.  We refer to the later Remark \ref{rem-normalizationVNS-longtime-simpl} where we explain why it is enough to consider the value $1$ for the particle densities.

 Given an initial condition, let us introduce the total momentum
\begin{align}
\label{def:WM}
 \mathcal{W}_M\vcentcolon= \frac{1}{1+M}\left\langle U_0+\int \rho^\alpha_0 v_0^\alpha \, \mathrm{d}\mu(\alpha) \right\rangle,
\end{align}
where $M$ is the total mass for the particles
\begin{align}
\label{def:M-VNS}
M\vcentcolon=\int \langle \rho_0^\alpha \rangle \, \mathrm{d}\mu(\alpha),
\end{align}
which is a quantity  formally conserved  by the dynamics of \eqref{eq:VNS-CHAP4}. Our main result, asserting global well-posedness and asymptotic stability for~\eqref{eq:VNS-CHAP4}  around the steady state \eqref{eq:mono-VNS}, reads as follows.
\begin{Thm}\label{thm-globalconvergence}
There exist $\Lambda>0$ and $\eps_0>0$  such that for any integer $k \geq 4$,
the following holds. For any $\mathcal{V} \in \R^3$, if
    \begin{align}
        \label{1st-smallness-cond} \Vert \rhorho_0-1 \Vert_{L^{\upinfty}_\alpha H^{k-1}_x}+ \Vert \vv_0-\mathcal{V} \Vert_{\mathcal{H}^{k,p}}+ \Vert U_0-\mathcal{V} \Vert_{H^k_x} &\leq \eps_0,
    \end{align}
with
\begin{align}\label{cond-exposant-p-thm}
    \left\{
    \begin{array}{ll}
        p \in [2, \infty]  & \mbox{if } \ \ \Gamma=0, \\[2mm]
        p =\infty  & \mbox{if } \ \ \Gamma \neq 0,
    \end{array}
\right.
\end{align}
then the system \eqref{eq:VNS-CHAP4} admits a unique global in time solution $(\rhorho, \vv, U)$ such that
\begin{align*}
\rhorho \in L^\infty(\R^+; L^\infty_\alpha H^{k-1}_x), \ \, \vv \in L^\infty(\R^+; \mathcal{H}^{k,p}), \ \, U \in L^\infty(\R^+; H^k_x ), \ \ \nabla_x U \in  {L^2}(\R^+; H^{k-1}_x ),
\end{align*}
associated  with the initial data $(\rhorho_0, \vv_0, U_0)$. Furthermore, for all $0<\lambda<\Lambda$, there exists $C_\lambda>0$ such that 
\begin{align}\label{final-convergence-velocities}
    \Vert \vv(t) - \mathcal{W}_M \Vert_{L^p_\alpha L^2_x}+ \Vert U(t) - \mathcal{W}_M \Vert_{L^2_x} &\leq C_\lambda e^{-\lambda t}, \ \ t>0,
\end{align}
and there exist $\gamma=\gamma(\Lambda,k) >0$ and $C_\gamma>0$ such that 
\begin{align}\label{final-convergence-velocities2}
     \left\Vert \nabla \vv(t) \right\Vert_{L^{\upinfty}_\alpha H^{k-1}_x}+ \left\Vert \nabla U(t) \right\Vert_{H^{k-1}_x} &\leq C_\gamma e^{-\gamma t}, \ \ t>0.
\end{align}
Finally, there exists some density profile $$\bm{\rho}_{\infty} \in L^{\infty}_\alpha L^{\infty}(\T^3),$$ satisfying for all $\alpha \in I$,
$$\langle {\rho}^\alpha_{\infty}\rangle = \langle {\rho}^\alpha_{0}\rangle,$$ 
and such that for all $0<\lambda<\Lambda$, there exists $C'_\lambda>0$ such that or all $t \geq 0$
\begin{align}\label{final-convergence-densities}
\esssup_{\alpha \in (I, \mu)} & \, \, \mathrm{W}_1 \Big(\rho^{\alpha}(t),\rho^{\alpha}_{\infty}(x-t \mathcal{W}_M  ) \Big) \leq C'_\lambda e^{-\lambda t}, \ \ t>0,
\end{align}
where $\mathrm{W}_1$ is the $1-$Wasserstein distance on $\T^3$ (as defined in \eqref{def:Wasserstein1}).
\end{Thm}
A series of remarks is now in order.

\begin{Rem}
    In the former theorem, the solution satisfy actually  the same time continuity properties as that of Theorem \ref{thm:appliVNS-incomp} and the remark below this statement.
\end{Rem}

\begin{Rem}[On the convergence of the  velocities]
One important observation is that the fluid velocity $U(t)$ and particles velocities $\vv(t)$ do \textit{not} converge back to the perturbed steady state $\mathcal{V}$ (while the family of densities roughly behaves like a travelling wave). We can actually rephrase the final asymptotics \eqref{final-convergence-velocities} in the following way: a small perturbation of a constant density equation to $1$ and a common constant velocity $\mathcal{V}$  by $(\bm{r}_0, \bm{w}_0, u_0)$ makes the velocities align towards $\mathcal{W}_M$, that is also expressed in terms of the perturbation as
\begin{multline*}
\mathcal{W}_M=\mathcal{V}^{\infty}[\bm{r}_0, \bm{w}_0, u_0, \mathcal{V} ] \vcentcolon=\mathcal{V}+\frac{1}{2+\int \langle r_0^\alpha \rangle \, \mathrm{d}\mu(\alpha)}\left\langle u_0+\int r^\alpha_0 w^\alpha_0 \D \mu(\alpha) \right\rangle \\
+\frac{1}{2+\int \langle r_0^\alpha \rangle \, \mathrm{d}\mu(\alpha)}\left\langle \int  w^\alpha_0 \D \mu(\alpha) \right\rangle.
\end{multline*}
This convergence is a specific feature of the periodic torus case. Indeed, the forthcoming analysis based on conservation laws of the equations (see Section \ref{sectionVNS-Conslaws-Decay})  will reveal that  the spatial average of the perturbation of the steady state, that is 
\begin{align*}
    \mathcal{V}(t)=\langle U(t)-\mathcal{V}\rangle, \ \ \mathcal{V}^\alpha(t)=\langle v^\alpha(t)-\mathcal{V}\rangle,
\end{align*}
are convergent in time, but not decaying to zero.

\end{Rem}

\begin{Rem}[On the spectral gap]\label{rmk-partiallydissipative}
As it will be made clearer later on, the constant $\Lambda>0$ in our statement comes from an explicit spectral gap of the associated linearized system (around $(\rho^\alpha, v^\alpha, U)=(1, \mathcal{V}, \mathcal{V})$ for any constant vector $\mathcal{V} \in \R^3$), the existence of which 
is a specific feature of the periodic torus case. At the nonlinear level, we obtain an almost optimal exponential decay rate with respect to the linearized problem. 

The situation is quite different in the whole space $\R^3$ where there is no spectral gap and where polynomial in time decay estimates are expected. When the multiphasic ansatz is reduced to one phase, that is with the coupling between one pressureless Euler equation and the incompressible Navier--Stokes system (see \eqref{eq:Euler-NS}), it has been proven in \cite{choi2024revisit, ChoiJungR3,HuangTangZouR3, lemarié2025, zhai} that small and smooth perturbations of the zero stationary state $(\rho, v,U)=(0,0,0)$ yield global in time solutions for which the velocities converge back algebraically to zero. Let us also mention the recent work \cite{Danchin-ENS2026} (in the spirit of \cite{danchin2024fujita} for the Vlasov--Navier--Stokes system), where, in particular, the smallness of the initial density is relaxed. We finally mention the paper \cite{LemarieMultiphaseVNS} on the stability of the trivial steady state for the multiphasic version of the system on the whole space with critical (Besov) regularity.

Let us finally emphasize that the decay in time of both the fluid velocity and the particle velocities is a genuinely nontrivial consequence of the fluid--particle coupling in the Vlasov--Navier--Stokes system.
 Indeed, explicit dissipation is only present in the Navier--Stokes equations on $U$ (through the viscous Laplacian term) and not in the coupled equations for the $v^\alpha$. From this perspective, the multiphasic VNS system \eqref{eq:VNS-CHAP4} provides a notable example of a partially dissipative system (see in particular the survey \cite{danchin2023-dissipativereview} and references therein).

\end{Rem}

\begin{Rem}[On the well-preparedness assumption]
The key assumption ensuring that the data are initially close to being monokinetic is $\Vert \vv_0^\alpha-\mathcal{V} \Vert_{\mathcal{H}^{k,p}}\ll 1$, in view of \eqref{1st-smallness-cond}. Furthermore, since 
\begin{align*}
    \Vert  \langle  \vv_0- U_0 \rangle \Vert_{L^p_\alpha}\leq \Vert    \vv_0- U_0 \Vert_{L^p_\alpha L^2_x} \leq  \Vert \vv_0-  \mathcal{V} \Vert_{\mathcal{H}^{k,p}}+\Vert U_0-  \mathcal{V}\Vert_{H^k_x} \leq 2\eps_0,
\end{align*}
the hypothesis \eqref{1st-smallness-cond} ensures a well-preparedness assumption between the initial velocities of the fluid and particles phases, that have to be sufficiently close on average. 
\end{Rem}

\begin{Rem}[On the regularity assumptions]
    Our result is coherent with the local well-posedness framework developed in Chapter \ref{Part1-LWP} and applied in Section \ref{SubsecVNSmultiphaseLWP}. More precisely:
    \begin{itemize}
        \item  The regularity threshold that is required for the global result in Theorem \ref{thm-globalconvergence}, that is actually $k \geq 2+3/2$, is not optimal in the sense of the local well-posedness theory from Sections \ref{sec:abstract_existence}--\ref{SubsecVNSmultiphaseLWP}, which requires at least $1+3/2$ derivatives in space. This slight difference appears technical: the proof is indeed based on the exponential decay of a semigroup associated to a linearized system satisfied by gradients of the perturbation. This system involves a source term containing nonlinear terms whose control requires two derivatives in $L^\infty$, explaining the $2+d/2$ threshold. 
 
        \item We keep the total velocity $\vv=\mathcal{V}+\bm{w}$ in the space $\mathcal{H}^{k,p}$, i.e. we have a perturbation $\bm{w} \in \mathcal{H}^{k,p}$ but not necessarily in $L^\infty_\alpha H^k_x$. 
        As a consequence, this allows to consider distribution functions which are not necessarily compactly supported in velocity (see Remark \ref{rem:bounded_velocity}).
        However,  taking into account the additional nonlinear drag $\Gamma(U-v^\alpha)$ seems to require the additional constraint $p=\infty$  (see Section \ref{Section-nonlindrag} at the end of the chapter). Therefore our result in that case is restricted to data that are compactly supported in velocity. Note that this constraint was already present in the local well-posedness theory, see Section \ref{sec:nonlindrag-LWP}. As far as we know, this is nevertheless the first stability result handling an extra nonlinear drag term in fluid-kinetic settings.
    \end{itemize}
\end{Rem}

   \begin{Rem}
By introducing a viscosity parameter $\nu>0$ in front of the term $\Delta U$ and a friction parameter $\kappa >0$ in front of the drag terms $(U-v)$ in the VNS system, our methods  allow to quantify exactly the optimal exponential decay rate $\Lambda(\nu, \kappa)$ in terms of $\nu$ and $\kappa$. 
We can indeed show that
\begin{align*}
    \Lambda(\nu, \kappa)=\kappa+\frac{\nu}{2}-\frac{1}{2}\sqrt{\nu^2+4\kappa^2}>0.
\end{align*}
\end{Rem}

\begin{Rem}[On the renormalisation for the densities]\label{rem-normalizationVNS-longtime-simpl} 
    The assumption in Theorem~\ref{thm-globalconvergence} that the densities $\rhorho_0$ are small perturbations of the constant state $1$ is merely a normalisation, and we could rather assume, instead of \eqref{1st-smallness-cond}, that there holds initially
    $$\left\Vert \frac{\rhorho_0}{\mathrm{C}}-1 \right\Vert_{L^{\upinfty}_\alpha H^{k-1}_x}+ \Vert \vv_0-\mathcal{V} \Vert_{\mathcal{H}^{k,p}}+ \Vert U_0-\mathcal{V} \Vert_{H^k_x} \leq \eps_0$$
    for some constant $\mathrm{C}>0$. The only minor modification stems from the forthcoming spectral linearized analysis -- see Remark \ref{rem-linsystem-incom-massnot1}. Note that this more general initial smallness assumption allows for a little bit more room for the initial data in the subsequent kinetic counterpart of the multiphasic result -- see Remark \ref{rem-wellpreparedness-VNSlargetime}. Note also that this construction also shows that it is always legitimate to consider a probability space $(\mathcal{I}, \mu)$ in Theorem~\ref{thm-globalconvergence}.
\end{Rem}

Theorem \ref{thm-globalconvergence} has a kinetic counterpart, that reads as follows.
\begin{Thm}\label{thm-cinetique-VNSincomp}
Under the assumptions of Theorem \ref{thm-globalconvergence}, let us consider the global solution $(\rhorho, \vv, U)$ to the multiphasic system \eqref{eq:VNS-CHAP4} and set
\begin{align*}
    f(t,x,v)=\int_{I} \rho^\alpha(t,x) \otimes \delta_{v=v^\alpha(t,x)} \, \D \mu(\alpha).
\end{align*}
Then $(f, U)$ is the unique global in time solution of the incompressible Vlasov--Navier--Stokes system 
\begin{equation}
\label{eq:VNS-kinetic}
\left\{ 
\begin{gathered}
\partial_t f + v \cdot \nabla_x f+\mathrm{div}_v[(U-v)f+\Gamma(U-v)f] =0,\\
\partial_t U + (U\cdot \nabla_x ) U + \nabla_x P - \Delta_x U = -\int_{\R^3} \big( U - v+ \Gamma(U-v) ) f \, \mathrm{d}v
,\\
\Div_x(U) = 0,
\end{gathered}
\right.
\end{equation}
in the sense of Definition \ref{def:solution_coupled}. Furthermore, we have the following asymptotic behaviour: for all $0 <\lambda <\Lambda$ and $t>0$, there holds
\begin{multline}\label{final-convergence-densities-kin}
\Vert U(t) - \mathcal{W}_M \Vert_{L^2_x}+ \mathrm{W}_1 \left(f(t),\left( \int_I \, \rho^\alpha_\infty(x - t \mathcal{W}_M)\mathrm{d}\mu(\alpha) \right) \otimes\delta_{v=\mathcal{W}_M}  \Big) \right) \\ \lesssim_\lambda e^{-\lambda t},
\end{multline}
where $\mathrm{W}_1$ is the $1-$Wasserstein distance on $\T^3 \times \R^3$ (see Definition \ref{def:Wasserstein1}).
\end{Thm}

\begin{Rem}
In Theorem \ref{thm-cinetique-VNSincomp}, the solution $(f,U)$ satisfies the continuity in time properties and regularity of the moments stated at the end of Theorem \ref{thm:appliVNS-incomp}.
\end{Rem}

\begin{Rem}[On the well-preparedness assumption II]\label{rem-wellpreparedness-VNSlargetime}
    For the assumptions of Theorem \ref{thm-globalconvergence} to be satisfied, the initial data $f_0$ that we can consider must be well-prepared. As in Remark~\ref{rem:hyp-mono}, this means typically that
$$
f_{0} (x,v) = \rho_{0} (x)  \varphi_n\left(v-v_{0}(x)\right),
$$
where $\varphi_n(v) = n^3 \varphi(n v)$ with $\varphi\geq 0, \int \varphi(v) \D v =1$ being an approximation of unity satisfying $\int \vert v \vert^p \varphi(v) \, \mathrm{d}v<\infty$. In this case, we take $(I,\mu)= (\R^d, \varphi (\alpha)  \mathrm{d} \alpha)$, that is a probability space, and  
$$
\rho^\alpha_0(x) = \rho_0(x), \quad v^\alpha_{0} (x) =v_{0}(x) + n^{-1} \alpha.
$$
Then, for a constant velocity $\mathcal{V} \in \R^3$, the assumption \eqref{1st-smallness-cond} of Theorem~\ref{thm-globalconvergence} is satisfied provided that $n$ is taken sufficiently large and that the quantity
\begin{align*}
    \|\rho_0-1\|_{H^{k-1}_x}+\|v_0-\mathcal{V}\|_{H^{k}_x}+\|U_0-\mathcal{V}\|_{H^{k}_x}
\end{align*}
is sufficiently small, since there holds
\begin{multline*}
    \Vert \rhorho_0-1 \Vert_{L^{\upinfty}_\alpha H^{k-1}_x}+ \Vert \vv_0-\mathcal{V} \Vert_{\mathcal{H}^{k,p}}+ \Vert U_0-\mathcal{V} \Vert_{H^k_x} \\
    \lesssim \|\rho_0-1\|_{H^{k-1}_x}+\|v_0-\mathcal{V}\|_{H^{k}_x}+\|U_0-\mathcal{V}\|_{H^{k}_x} +\frac{1}{n}\left(\int \vert v \vert^p \varphi(v) \, \mathrm{d}v  \right)^{1/p}.
\end{multline*} 
Hence, the hypotheses of Theorem~\ref{thm-globalconvergence} correspond to kinetic initial data that are small perturbations of a monokinetic equilibrium. Note finally that, in view of Remark \ref{rem-normalizationVNS-longtime-simpl}, we can also assume that $\rho_0$ is a small perturbation of any constant $\mathrm{C}>0$.

\end{Rem}

\subsection{The case of the compressible Vlasov--Navier--Stokes system}

We now focus on the multiphasic compressible version of the Vlasov--Navier--Stokes system, that is 
\begin{equation}
\label{eq:VNScompressibleCHAP4}
\left\{ 
\begin{aligned}
\partial_t \rho^\alpha + \Div ( \rho^\alpha v^\alpha) &= 0,\\
\partial_t v^\alpha + (v^\alpha \cdot \nabla) v^\alpha &= U - v^\alpha+ \Gamma(U-v^\alpha),\\
\partial_t \varrho + \Div ( \varrho U)&=0, \\
\varrho\left(\partial_t U + (U \cdot \nabla) U \right)+ \nabla P(\varrho) - \mathcal{A}U &= -\int \big( U - v^\alpha+ \Gamma(U-v^\alpha) ) \rho^\alpha \D \mu(\alpha).
\end{aligned}
\right.
\end{equation}
The fluid pressure  is assumed to be $P=P(\varrho)$ for a smooth function $P :\R^+ \rightarrow \R$ applied to the fluid density $\varrho$. This corresponds to the standard barotropic case. We also assume that $P(0)=0$ and $P'(1)>0$ (see Remark \ref{rem-normalizationcompVNS-longtime}).

The so-called Lamé opérator $\mathcal{A}$ is defined as $\mathcal{A}u=\Delta u+ 2 \nabla \Div(u)$ (all coefficients have been normalized for the sake of simplicity). As in the incompressible case, the quadratic nonlinear drag $\Gamma\in \mathscr{C}^{\infty}(\R^3;\R^3)$ satisfies $\Gamma(0)=0$ and $\Gamma'(0)=0$. The system~\eqref{eq:VNScompressibleCHAP4} is endowed with the initial data
\begin{align*}
    (\rhorho, \vv, \varrho, U)_{\mid t=0}=(\rhorho_0, \vv_0, \varrho_0, U_0).
\end{align*}

\medskip

As before, we will study the multiphasic system \eqref{eq:VNScompressibleCHAP4} in the neighborhood of constant stationary solutions of the form
\begin{align*}
\rho^\alpha \equiv 1 , \qquad v^\alpha \equiv \mathcal{V} , \qquad \varrho \equiv 1, \qquad U \equiv \mathcal{V} ,
\end{align*}
where $\mathcal{V} \in \R^3$ is a constant vector. We also introduce 
\begin{align}
\label{def:WmM}
\mathcal{W}_{m,M}\vcentcolon= \frac{1}{m+M}\left\langle \varrho_0 U_0+\int \rho^\alpha_0 v_0^\alpha \, \mathrm{d}\mu(\alpha) \right\rangle,
\end{align}
where $M$ (resp. $m$) is the total mass of particles (resp. of the fluid) defined by
\begin{align}
\label{def:Mm-VNScomp}
M\vcentcolon=\int \langle \rho_0^\alpha \rangle \, \mathrm{d}\mu(\alpha), \ \ m\vcentcolon=\langle \varrho_0 \rangle,
\end{align}
which are conserved quantities.
Our main result for \eqref{eq:VNScompressibleCHAP4}  reads as follows.
\begin{Thm}\label{thm-globalconvergence-comp}
There exist $\Lambda>0$ and $\eps_0>0$ 
such that for any integer $k \geq 5$,
the following holds. For any $\mathcal{V} \in \R^3$, if
    \begin{align}
        \label{1st-smallness-cond-comp} \Vert \rhorho_0-1 \Vert_{L^{\upinfty}_\alpha H^{k-1}_x}+ \Vert \varrho_0-1 \Vert_{H^{k}_x}+ \Vert \vv_0 -\mathcal{V} \Vert_{\mathcal{H}^{k,p}}+ \Vert U_0-\mathcal{V} \Vert_{H^k_x} &\leq \eps_0, \\
    \label{3rd-novacuum-comp}      \inf \varrho_0 &\geq c_0,
    \end{align}
    for some $c_0>0$ and with
\begin{align}\label{cond-exposant-p-thmcomp}
    \left\{
    \begin{array}{ll}
        p \in [1, \infty]  & \mbox{if } \ \ \Gamma=0, \\[2mm]
        p =\infty  & \mbox{if } \ \ \Gamma \neq 0,
    \end{array}
\right.
\end{align}
then the system \eqref{eq:VNScompressibleCHAP4} admits a unique global in time solution $(\rhorho, \vv, \varrho, u)$ such that
\begin{align*}
 &\rhorho \in L^\infty(\R^+; L^\infty_\alpha H^{k-1}_x)  , \ \ \vv \in L^\infty(\R^+; \mathcal{H}^{k,p}), \\ 
 &\varrho \in L^\infty(\R^+; H^{k}_x ), \ \ U \in L^\infty(\R^+; H^k_x ),
\end{align*}
associated with the initial data $(\rhorho_0, \vv_0, \varrho_0, U_0)$. Furthermore, for all $0<\lambda<\Lambda$, there exists $C_\lambda>0$ such that
\begin{align}\label{final-convergence-velocities-comp}
\Vert \nabla \varrho(t)  \Vert_{L^{\upinfty}_x}+ \Vert \vv(t) - \mathcal{W}_{m,M} \Vert_{L^p_\alpha L^2_x}+ \Vert U(t) - \mathcal{W}_{m,M} \Vert_{L^2_x} \leq C_\lambda e^{-\lambda t}, \ \ t>0.
\end{align}
Finally, there exists some density profile $$\bm{\rho}_{\infty} \in L^{\infty}_\alpha L^{\infty}_x,$$ satisfying for all $\alpha \in I$,
$$\langle {\rho}^\alpha_{\infty}\rangle = \langle {\rho}^\alpha_{0}\rangle,$$ 
and such that for all $0<\lambda<\Lambda$, there exists $C'_\lambda>0$ such that
\begin{align}\label{final-convergence-densities-comp}
\esssup_{\alpha \in (I, \mu)} & \, \, \mathrm{W}_1 \Big(\rho^{\alpha}(t),\rho^{\alpha}_{\infty}(x-t \mathcal{W}_{m,M}  ) \Big) \leq C'_\lambda e^{-\lambda t}, \ \ t>0.
\end{align}
\end{Thm}
When $(\rhorho,\vv)\equiv (\boldsymbol{0},\boldsymbol{0})$, the system reduces to the classical compressible Navier--Stokes system. Interestingly, the study of the large time behaviour of this system around constant states in the torus $\T^2$ or $\T^3$ has been the object of several recent works, see in particular \cite{ZZ,DW}.
When only one phase is retained in the pressureless Euler equations of \eqref{eq:VNScompressibleCHAP4}, the system results in the so-called (pressureless) Euler--Navier--Stokes system; we mention the long-time asymptotic results of \cite{GWZ} for solutions around $(\rho, v, \varrho, U)=(1,0,1,0)$ and of \cite{li2025global} for solutions around $(\rho, v, \varrho, U)=(0,0,1,0)$, when the system is set on $\R^3$, and for Sobolev data with an integrability hypothesis.
We may also mention the recent work \cite{LiShouZhang} where, for solutions around $(\rho, v, \varrho, U)=(0,0,1,0)$ in $\R^3$ and in a critical framework, the incompressible and high-friction limit of the system is studied in large times.

\medskip

Remarks similar to the ones made after Theorem \ref{thm-globalconvergence} for the case of an incompressible fluid also pertain for the compressible case, \textit{mutatis mutandis}. 
\begin{Rem}\label{rem-normalizationcompVNS-longtime}
Let us also comment on the differences in terms of assumptions and results that are specific to the compressible Navier--Stokes case:
\begin{itemize}
    \item The assumption \eqref{3rd-novacuum-comp} is a standard hypothesis for what concerns the compressible Navier--Stokes systems, and asserts that there is no vacuum initially.
    \item In the final decay estimate \eqref{final-convergence-velocities-comp}, we only obtain the decay of the gradient of the fluid density $\varrho$. By Poincaré inequality on the torus, and because the average $\langle \varrho(t) \rangle$ is preserved along the time evolution,  we infer a similar decay for $\Vert  \varrho(t)-\langle \varrho_0 \rangle \Vert_{L^{\upinfty}_x}$.
    \item In view of the latter point, it would make sense to rather assume that  $\varrho_0-\langle \varrho_0 \rangle $ is initially small (with $\langle \varrho_0 \rangle>0)$ . In that case, arguing as for the normalisation procedure from Remark \ref{rem-normalizationVNS-longtime-simpl} by rescaling by $\langle \varrho_0 \rangle$, we can  always come back to the case treated in Theorem \ref{thm-globalconvergence-comp}. The only change stem from a modification of the dissipation coefficients in the Lamé operator $\mathcal{A}u=\Delta u+ 2 \nabla \Div(u)$, where new constant coefficients appear (see however Remark \ref{rem:Lamé-coeff} showing that this is harmless for our analysis), and from a new pressure law, where one has to assume that from the beginning that $P'(\langle \varrho_0 \rangle)>0$.
\end{itemize}
\end{Rem}

Theorem \ref{thm-globalconvergence-comp} has the following kinetic counterpart.

\begin{Thm}\label{thm-cinetique-VNScomp}
Under the assumptions of Theorem \ref{thm-globalconvergence-comp}, let us consider the global solution $(\rhorho, \vv, \varrho, U)$ to the multiphasic system \eqref{eq:VNScompressibleCHAP4} and set
\begin{align*}
    f(t,x,v)=\int_{I} \rho^\alpha(t,x) \otimes \delta_{v=v^\alpha(t,x)} \, \D \mu(\alpha).
\end{align*}
Then $(f, \varrho, U)$  is the unique global in time solution of the compressible Vlasov--Navier--Stokes system 
\begin{equation}
\label{eq:VNS-kinetic-comp}
\left\{ 
\begin{gathered}
\partial_t f + v \cdot \nabla_x f+\mathrm{div}_v[(U-v)f+\Gamma(U-v^\alpha)f] =0,\\
\partial_t \varrho + \Div_x ( \varrho U)=0, \\
\varrho\left(\partial_t U + (U \cdot \nabla_x) U \right)+ \nabla_x P(\varrho) - \mathcal{A}U = -\int_{\R^3} \big( U - v+ \Gamma(U-v) ) f \, \mathrm{d}v,
\end{gathered}
\right.
\end{equation}
in the sense of  Definition~\ref{def:solution_coupled}. Furthermore, we have the following convergence: for all $0 <\lambda <\Lambda$ and  for all $t \geq 0$, there holds
\begin{multline}\label{final-convergence-densities-comp-kin}
\Vert \nabla \varrho(t)  \Vert_{L^{\upinfty}_x}+ \Vert U(t) - \mathcal{W}_{m,M} \Vert_{L^2_x}\\+\mathrm{W}_1 \left(f(t),\left( \int_I \, \rho^\alpha_\infty(x-t\mathcal{W}_{m,M})\mathrm{d}\mu(\alpha) \right) \otimes\delta_{v=\mathcal{W}_M}  \Big) \right) \lesssim e^{-\lambda t}. 
\end{multline}
\end{Thm}
To our knowledge, there have been only a few works dealing with the large time behaviour of the compressible Vlasov--Navier--Stokes system.  
The one-dimensional case is settled in \cite{LiShou}: exponential decay to equilibrium is proved for global weak solutions, without requiring initial data close to equilibrium. In higher dimension, a conditional large-time result is obtained in \cite{Choi-comp}, assuming suitable bounds on the solution itself. In principle, the strategy of \cite{HKMM} could be adapted to remove this conditionality in the case of a compressible fluid. However, we are not aware of such a complete proof in the literature.
Theorems \ref{thm-globalconvergence-comp} and \ref{thm-cinetique-VNScomp} somehow fill in this gap, though with a different set of assumptions for the initial condition and a completely different approach.

\subsection{Strategy of proof}\label{section-stratproofVNS}
Let us sketch the main elements for the proof of Theorems \ref{thm-globalconvergence} and \ref{thm-globalconvergence-comp}. We shall study the multiphasic systems \eqref{eq:VNS-CHAP4} and \eqref{eq:VNScompressibleCHAP4} by looking for solutions of the form
\begin{align}
 & \underline{\text{for} \ \ \eqref{eq:VNS-CHAP4}:} \ \ \rho^\alpha = 1+r^\alpha, \qquad v^\alpha =\mathcal{V}+w^\alpha, \qquad U = \mathcal{V}+ u, \label{perturb-incomp}\\
 & \underline{\text{for} \ \ \eqref{eq:VNScompressibleCHAP4}:} \ \ \rho^\alpha = 1+r^\alpha, \qquad v^\alpha = \mathcal{V}+w^\alpha, \qquad \varrho = 1+n, \qquad U = \mathcal{V}+ u. \label{perturb-comp}
\end{align}
In what follows, we focus on the case of the incompressible Navier--Stokes equations, that is the multiphasic system \eqref{eq:VNS-CHAP4}, since the method we will develop contains the main ingredients that can be adapted to treat the case of \eqref{eq:VNScompressibleCHAP4}. We also drop here the nonlinear drag force $\Gamma$ to simplify the presentation.

With the notation from \eqref{perturb-incomp}, the main goal to obtainTheorem \ref{thm-globalconvergence} will be to show that the solution can be decomposed for large times as 
\begin{align}\label{eq:decompo-vitesse-stratVNS}
\begin{split}
    v^\alpha(t,x) &=  \mathcal{V}+\langle w^\alpha(t) \rangle + \mathrm{O}(e^{- \lambda t}), \\
U(t,x) &= \mathcal{V}+  \langle u(t) \rangle  + \mathrm{O}(e^{-\lambda t}),
\end{split}
\end{align}
for some $\lambda>0$ related to an explicit spectral gap obtained at the linearized level (see below). To describe the asymptotic behaviour of the averages in space, we shall also obtain that 
\begin{align*}
v^\alpha(t) - U(t) = \mathrm{O}(e^{-\lambda t}),
\end{align*}
This should eventually yield that for almost every $\alpha \in I$, $\langle w^\alpha(t) \rangle$  has the same limit as $t \to +\infty$ as $\langle u(t) \rangle$. The identification of the common limit, that is expected to be $\mathcal{W}_M-\mathcal{V},$ will be based on identities stemming from conservation of the total momentum.

Overall, our strategy is based on a fairly elementary linearized analysis combined with a bootstrap procedure at the nonlinear level --  a method that has to be put in contrast compared to the one of \cite{HKMM}. One of the main difficulties comes from the fact that on the one hand, no decay is expected to hold for the family of densities $(\rho^\alpha)$ and, on the other hand, a dissipation mechanism on the velocities is \textit{a priori} only coming from the Navier--Stokes part on $U$.
It is the coupling between the two equations through the drag force $U-v^\alpha$ that will yield  decay on the whole family of velocities $(v^\alpha)$ and $U$. It should be noticed that this stabilizing mechanism due to the coupling is somehow reminiscent of the case of the so-called partially dissipative systems (see Remark \ref{rmk-partiallydissipative}).

\paragraph{Spectral analysis.}
First, we perform a linearisation around constant density/common velocity profiles (see \eqref{perturb-incomp}): the main observation is that this procedure entirely decouples the perturbation for the densities and the perturbation for the velocities. It therefore turns out sufficient to understand the following linearized system
\begin{equation*}
\left\{ 
\begin{aligned}
\partial_t w^\alpha + (\mathcal{W} \cdot \nabla) w^\alpha&=  u- w^\alpha,\\
\partial_t u+  (\mathcal{W} \cdot \nabla)u - \Delta u + \nabla P &=  \int (w^\alpha-u) \D \mu(\alpha),\\
\Div(u) &= 0,
\end{aligned}
\right.
\end{equation*}
for a fixed $\mathcal{W} \in \R^3$. At this stage, a direct frequency analysis yields an \textbf{explicit spectral gap} $\Lambda>0$, that leads to exponential decay at rate $\Lambda$ for $\Vert \nabla \ww \Vert_{L^{\upinfty}_\alpha L^2_x}$ and  $\Vert \nabla u \Vert_{L^2_x}$, as well as for the difference of the averages of the velocities, namely $\langle u-w^\alpha \rangle$. Let us highlight that no decay can be expected for each average of the velocities: this will entail, at the nonlinear level, the appearance of the constant in space profiles $\langle w^\alpha(t) \rangle $ and $\langle u(t) \rangle $ in the above decomposition of the solution (see \eqref{eq:decompo-vitesse-stratVNS}).

\paragraph{Bootstrap argument.} Our stability proof for the nonlinear system is based on a bootstrap argument. We will show that exponential type decay bounds in time, inherited from the linearized spectral analysis, can be propagated. At the same time, we shall ensure that 
the $L^1_T Lip_x$ norms of the velocities remain under control, and   a continuation criterion from Section \ref{SubsecVNSmultiphaseLWP}, namely Proposition~\ref{prop-VNSblowup} will allow us to claim that the solution exists globally in time.

More precisely, for any $\lambda \in (0, \Lambda)$ (where $\Lambda>0$ stands for the former spectral gap), and for some $\lambda' \in (0, \lambda)$ and  $\delta \in (0,1)$ to be fixed later, we consider the maximal time of existence $T>0$ (for the perturbed system satisfied by $(r^\alpha, w^\alpha, u)$) such that for all $t \in [0,T]$
\begin{align}\label{eq:bounds-bootstrapSTRAT-VNS}
\begin{split}
    \Vert \nabla \ww(t) \Vert_{L^{\upinfty}_\alpha L^2_x}+\Vert \nabla u(t) \Vert_{L^2_x} + \Vert \langle \ww-u \rangle (t)\Vert_{L^p_\alpha}\leq \delta e^{-\lambda t}, \\
    \Vert  (\nabla \ww, \nabla^2 \ww)(t) \Vert_{L^{\upinfty}_\alpha L^{\upinfty}_x}+\Vert (\nabla u, \nabla^2 u)(t) \Vert_{L^{\upinfty}_x} \leq \delta e^{-\lambda' t}. 
    \end{split}
\end{align}
The former set of bounds \eqref{eq:bounds-bootstrapSTRAT-VNS} will be referred to as the \textit{bootstrap assumption}. Note that the presence of second-order derivatives is only technical, and comes from the fact that the linear analysis yields decay on \textit{derivatives} of the solution, hence requiring a control on two derivatives at the nonlinear level.

The main goal is then to prove that the decay estimates from \eqref{eq:bounds-bootstrapSTRAT-VNS} actually hold on $[0,T]$ by replacing the constant $\delta$ by $C(\eps_0+\delta^2)$, for some constant $C>0$, where $\eps_0$ controls the size of the initial perturbation. By picking $\delta$ and $\eps_0$ small enough so that
$C(\eps_0+\delta^2)<\delta/2$, a continuation argument then classically allows us to close the bootstrap argument, yielding $T=+\infty$ and then concluding the proof. To do so, we first obtain two family of estimates:
\begin{itemize}
    \item on the one hand, \textit{refined exponentially growing bounds for ``high-order'' norms} of the form
    \begin{equation}
    \label{eq:refined-expoSTRAT-VNS}
        \forall t \in [0,T], \ \ \| \ww(t) \|_{\mathcal{H}^{k,p}}+ \| u(t)\|_{H^k_x} \lesssim \eps_0 e^{\sigma^{-1} t},
    \end{equation}
    for all $\sigma>1$. This comes from a refinement of the local well-posedness theory from Section \ref{SubsecVNSmultiphaseLWP}, here directly performed at the level of the perturbation. Making this exponential bound as tame as possible (thanks to a large parameter $\sigma$) will prove essential at the final stage of the argument to obtain almost optimal result in terms of Sobolev regularity.
    \item on the other hand, \textit{improved exponential decay of ``lower-order'' norms} of the form
    \begin{align}\label{eq:decay-expoSTRAT-VNS}
    \forall t \in [0,T], \ \ \Vert \nabla \ww(t) \Vert_{L^{\upinfty}_\alpha L^2_x}+\Vert \nabla u(t) \Vert_{L^2_x} \lesssim (\eps_0 + \delta^2) e^{-\lambda t},
\end{align}
which is obtained by combining the exponential decay estimates from the linearized semigroup and Duhamel formula, treating nonlinearities as source terms. Nonlinearities of quadratic type are well-designed to obtain the $\delta^2$ pre-factor in the previous inequality, thanks to $L^2_x-L^\infty_x$ estimates. In the case of \eqref{eq:VNS-CHAP4}, one also encounters the nonlinear term $\int r^\alpha (w^\alpha-u) \D \mu(\alpha)$, which is not of quadratic nature in terms of decay in time: it is actually handled by relying on the decay of the difference of the velocities in the bootstrap assumption and on the smallness of the initial density perturbation, explaining the pre-factor $\eps_0$ in~\eqref{eq:decay-expoSTRAT-VNS}
    
\end{itemize}

To improve the last bound from the bootstrap assumption \eqref{eq:bounds-bootstrapSTRAT-VNS}, we finally rely on  {interpolation} to write that for any integer $k>1+3/2$, we have
\begin{equation*}
\begin{aligned}
 \| \nabla u \|_{L^{\upinfty}_x}   \lesssim\|  \nabla u  \|_{L^2_x}^{\theta_k}    \| \mathrm{D}^k u\|_{L^2_x}^{1-\theta_k}, \ \ \| \nabla \ww \|_{L^{\upinfty}_\alpha  L^{\upinfty}_x}   \lesssim \|  \nabla \ww\|_{L^{\upinfty}_\alpha L^2_x}^{\theta_k}    \| \mathrm{D}^k \ww \|_{L^{\upinfty}_\alpha  L^2_x}^{1-\theta_k},
\end{aligned}
\end{equation*}
for some $\theta_k \in (0,1)$. This comes from the Gagliardo--Nirenberg--Sobolev inequality, and similar estimates also apply for the second-order derivatives. Combining the previous interpolation inequalities, the exponential bound \eqref{eq:refined-expoSTRAT-VNS} and the decay estimate \eqref{eq:decay-expoSTRAT-VNS} yields
\begin{multline*}
     \forall t \in [0,T], \ \ \| \nabla u(t) \|_{L^{\upinfty}_x}+\| \nabla \ww(t) \|_{L^{\upinfty}_\alpha  L^{\upinfty}_x} \\\lesssim  \big(\eps_0+\eps_0^{\theta_k}\delta^{2(1-\theta_k)} \big) \exp\big(\theta_k \sigma^{-1} t+(1-\theta_k)\lambda t \big).
\end{multline*}
By taking $\eps_0$ sufficiently small compared to $\delta$, and by choosing $\sigma$ large enough, a suitable choice of $\lambda'>0$ in the bootstrap assumption \eqref{eq:bounds-bootstrapSTRAT-VNS} ensures that we can improve the second decay estimate there. Note that this strategy will allow to us to obtain an \textit{almost} optimal global in time well-posedness result in the sense of the local well-posedness framework of Section \ref{SubsecVNSmultiphaseLWP}, since the exponent $k$ is fixed at the beginning of the proof.

A main specificity of the spatial periodic case is the following: since decay in time is only expected for the non-zero modes of the solutions, we have to introduce their average in the estimate and rely on Poincaré inequality to make gradients appear in order to use the bootstrap assumption \eqref{eq:bounds-bootstrapSTRAT-VNS}. Furthermore, since the velocities are not converging back to the initial velocity profile $\mathcal{V}$ that has been perturbed (but rather slightly deviate from it), we naturally introduce the expected asymptotic velocity $\mathcal{W}_M$ from Theorem \ref{thm-globalconvergence} and the difference $ \mathcal{W}_M-\mathcal{V}$. We refer to the subsequent Section \ref{sectionVNS-Stratbootstrap} for more details. 

\paragraph{Conservation laws and asymptotic density profile.}
The final goal is therefore to show that the averages of the perturbation $\langle u(t) \rangle$ and $\langle w^\alpha(t) \rangle$ converge to $\mathcal{W}_M-\mathcal{V}$ when $t \rightarrow + \infty$. This is equivalent to capturing the \textit{asymptotic behaviour of the spatial means} $\langle w^\alpha (t)\rangle$ and $\langle u (t)\rangle$ in the decomposition \eqref{eq:decompo-vitesse-stratVNS}. We crucially rely on the conservation of momentum
\begin{align*}
\frac{\D}{\D t} \int_{\T^3} \left[  \int (1+r^\alpha(t,x))\big(\mathcal{V}+ w^\alpha(t,x)) \D \mu(\alpha) + (\mathcal{V}+u(t,x)\big)  \D \mu(\alpha) \right] \D x = 0,
\end{align*}
written on the perturbation. This last part is reminiscent of the strategy performed in \cite{HKMM} at the level of the kinetic equation. We can then infer the final asymptotic behaviour for the total velocities in view of the decomposition \eqref{eq:decompo-vitesse-stratVNS}. The existence of the asymptotic profile for the densities $(\rho^\alpha)$ is obtained in the same abstract way as in \cite{HKMM}, and we refer to Proposition \ref{Prop-exist-asympto-density} for more details.

\medskip

\textbf{From now on, we fix a constant speed $\mathcal{V} \in \R^3$, and we will write all our solutions according to the decompositions \eqref{perturb-incomp} and \eqref{perturb-comp}. 
}
\section{Conservation laws and linear decay}\label{sectionVNS-Conslaws-Decay}

\subsection{Conservation laws}
We start by deriving some estimates bearing on space averages, that will allow us to capture the limit of the zero-modes of the solution to the perturbed system.
As explained previously in the strategy of proof from Section \ref{section-stratproofVNS}, this is based on conservation laws for the multiphasic Vlasov--Navier--Stokes systems on the torus.
\begin{Lem}\label{smallLemma-moyenne}
Let $(\rhorho, \vv, U)$  (resp. $(\rhorho, \vv, \varrho,U)$) be a smooth solution to the incompressible VNS system \eqref{eq:VNS-CHAP4}  (resp. to the compressible VNS system \eqref{eq:VNScompressibleCHAP4}), and assume that \eqref{perturb-incomp} (resp. \eqref{perturb-comp}) holds. 
The following inequalities hold: we have for all $t \geq 0$ and all $p \in [1, \infty]$: 
\begin{itemize}
    \item for the incompressible VNS system \eqref{eq:VNS-CHAP4}: using the notation \eqref{def:WM} and \eqref{def:M-VNS}, we have
    \begin{align*}
        &\vert \langle u(t) \rangle -(\mathcal{W}_M-\mathcal{V}) \vert \\
        &\qquad \qquad \lesssim_M  \Vert \rhorho(t) \Vert_{L^{\upinfty}_\alpha L^2_x} \left( \Vert \nabla u(t) \Vert_{L^2_x} + \int \Vert w^\alpha(t)-u(t) \Vert_{L^2_x}\, \mathrm{d}\mu(\alpha)\right), \\    
        &\Vert \langle w^\alpha(t) \rangle -(\mathcal{W}_M-\mathcal{V}) \Vert_{L^p_\alpha} \\
        &\qquad \qquad \lesssim_M \Vert \langle w^\alpha(t)-u(t) \rangle \Vert_{L^p_\alpha} \\
        &\qquad \qquad \quad +  \Vert \rhorho(t) \Vert_{L^{\upinfty}_\alpha L^2_x} \left( \Vert \nabla u(t) \Vert_{L^2_x} + \int \Vert w^\alpha(t)-u(t) \Vert_{L^2_x}\, \mathrm{d}\mu(\alpha)\right);
    \end{align*}
   \item for the compressible VNS system 
 \eqref{eq:VNScompressibleCHAP4}: using the notation  \eqref{def:WmM} and \eqref{def:Mm-VNScomp}, we have \begin{align*}
        &\vert \langle u(t) \rangle -(\mathcal{W}_{m,M}-\mathcal{V}) \vert \\
         &\qquad \lesssim_{M,m} \left( \Vert \rhorho(t) \Vert_{L^{\upinfty}_\alpha L^2_x}+ \Vert \varrho(t) \Vert_{L^2_x} \right)  \left( \Vert \nabla u(t) \Vert_{L^2_x} + \int \Vert w^\alpha(t)-u(t) \Vert_{L^1_x}\, \mathrm{d}\mu(\alpha)\right), \\    
        &\Vert \langle \ww(t) \rangle -(\mathcal{W}_{m,M}-\mathcal{V}) \Vert_{L^p_\alpha} \\
        & \qquad \lesssim_{M,m} \Vert \langle \ww(t)-u(t) \rangle \Vert_{L^p_\alpha}\\
        & \qquad \qquad  + \left( \Vert \rhorho(t) \Vert_{L^{\upinfty}_\alpha L^2_x}+ \Vert \varrho(t) \Vert_{L^2_x} \right) \left( \Vert \nabla u(t) \Vert_{L^2_x} + \int \Vert w^\alpha(t)-u(t) \Vert_{L^1_x}\, \mathrm{d}\mu(\alpha)\right).
    \end{align*}
\end{itemize}
\end{Lem}
\begin{proof}
First, we focus on the incompressible case, namely system \eqref{eq:VNS-CHAP4}. We rely on the conservation law
\begin{align*}
\frac{\D}{\D t} \int_{\T^3} \left(  U(t,x)+\int \rho^\alpha(t,x)v^\alpha(t,x) \D \mu(\alpha)  \right) \D x = 0,
\end{align*}
which is obtained by integration in $x$ and $\alpha$ and sum of
\begin{equation*}
\left\{  
\begin{gathered}
\partial_t (\rho^\alpha v^\alpha) + \mathrm{div}(\rho^\alpha v^\alpha \otimes v^\alpha) =  \rho^\alpha(U - v^\alpha+ \Gamma(U-v^\alpha)),\\
\partial_t U + \mathrm{div}( U \otimes U) - \Delta U + \nabla P =  -\int \rho^\alpha(U - v^\alpha+ \Gamma(U-v^\alpha)) \D \mu(\alpha).
\end{gathered}
\right.
\end{equation*}
We infer that for all $t >0$, we have
\begin{align*}
  \left\langle U(t)+\int \rho^\alpha(t) v^\alpha(t) \D \mu(\alpha)   \right\rangle  = \left\langle U_0+\int \rho^\alpha_0 v^\alpha_0 \D \mu(\alpha) \right\rangle  =(1+M) \mathcal{W}_{M},
\end{align*}
recalling the definition of $\mathcal{W}_M$ in~\eqref{def:WM}.
According to the decomposition \eqref{perturb-incomp} of the solution, we now write, using the definition \eqref{def:M-VNS} of $M$:
\begin{align*}
    \left\langle  U+\int \rho^\alpha v^\alpha \D \mu(\alpha)   \right\rangle&=\left\langle \mathcal{V}+u+\int \rho^\alpha (\mathcal{V}+w^\alpha) \D \mu(\alpha) \right\rangle \\
    &=\mathcal{V}+ \langle u \rangle + \mathcal{V}M   + \left\langle\int \rho^\alpha (w^\alpha-u) \D \mu(\alpha) \right\rangle
    \\ &\quad + \left\langle (u-\langle u \rangle) \int \rho^\alpha \D \mu(\alpha) \right\rangle  + M\langle u \rangle,
\end{align*}
where we have dropped the time dependency for simplicity, so that by the identity above, we  get
\begin{align*}
    (1+M)\mathcal{W}_M&=(1+M)\mathcal{V}+(1+M)\langle u\rangle+ \left\langle\int \rho^\alpha (w^\alpha-u) \D \mu(\alpha) \right\rangle
     \\
     & \quad + \left\langle (u-\langle u \rangle) \int \rho^\alpha \D \mu(\alpha) \right\rangle.
\end{align*}
hence by dividing by $1+M$
\begin{align*}
    \langle u(t) \rangle -(\mathcal{W}_M-\mathcal{V}) &=-\frac{1}{1+M}\Bigg(\left\langle\int \rho^\alpha (w^\alpha-u \D \mu(\alpha) \right\rangle \\
    & \quad + \left\langle (u-\langle u \rangle) \int \rho^\alpha \D \mu(\alpha) \right\rangle \Bigg).
\end{align*}
We then obtain the result by Cauchy--Schwarz inequality and Poincaré inequality. The last stated inequality is now straightforward by writing $$\langle w^\alpha(t) \rangle -(\mathcal{W}_M-\mathcal{V})=\langle w^\alpha(t)-u(t) \rangle +\langle u(t) \rangle -(\mathcal{W}_M-\mathcal{V}),$$
and taking the $L^p_\alpha$ norm.
\medskip

The proof follows the same lines for the compressible case, namely system~\eqref{perturb-comp}. One has the conservation law
\begin{align*}
\frac{\D}{\D t} \int_{\T^3} \left( \varrho(t,x) U(t,x)+ \int \rho^\alpha(t,x)v^\alpha(t,x) \D \mu(\alpha)  \right) \D x = 0,
\end{align*}
which is readily obtained from
System~\eqref{perturb-comp} and implies that 
for all $t >0$ we have, using the definition of $\mathcal{W}_{m,M}$ in~\eqref{def:WmM}:
\begin{align*}
  \left\langle \varrho(t) U(t)+\int \rho^\alpha(t) v^\alpha(t) \D \mu(\alpha)   \right\rangle  &= \left\langle \varrho_0 U_0+\int \rho^\alpha_0 v^\alpha_0 \D \mu(\alpha) \right\rangle  \\
  &=(m+M) \mathcal{W}_{m,M}.
\end{align*}
We can now develop with the decomposition \eqref{perturb-comp} and obtain
\begin{align*}
    \left\langle \varrho U+\int \rho^\alpha v^\alpha \D \mu(\alpha)   \right\rangle
    &=\left\langle \varrho(\mathcal{V}+u)+\int \rho^\alpha (\mathcal{V}+w^\alpha) \D \mu(\alpha) \right\rangle \\
    &=m \mathcal{V}+ \langle \varrho (u-\langle u \rangle) \rangle +m \langle u \rangle+ M\mathcal{V} \\
    &\quad + \left\langle\int \rho^\alpha (w^\alpha-u) \D \mu(\alpha) \right\rangle 
    \\ &\quad + \left\langle (u-\langle  u \rangle) \int \rho^\alpha \D \mu(\alpha) \right\rangle  + M\langle u \rangle,
\end{align*}
by definition of $m$, $M$ in \eqref{def:Mm-VNScomp}, so that we get
\begin{align*}
    (m+M)\mathcal{W}_{m,M}&=(m+M)\mathcal{V}+(m+M) \langle u \rangle +\langle \varrho (u-\langle u \rangle) \rangle\\
    &\quad + \left\langle\int \rho^\alpha (w^\alpha-u) \D \mu(\alpha) \right\rangle
     + \left\langle (u-\langle u \rangle) \int \rho^\alpha \D \mu(\alpha) \right\rangle.
\end{align*}
The end of the proof is then similar to the incompressible case.
\end{proof}

\subsection{Normalisation and system for the perturbation}

Let us now introduce the precise perturbative framework that we will use to prove Theorems \ref{thm-globalconvergence} and \ref{thm-globalconvergence-comp}.

\medskip
$\bullet$ \underline{Incompressible case \eqref{eq:VNS-CHAP4}}:
we look for a solution $(\rho^\alpha, v^\alpha, U)$ to \eqref{eq:VNS-CHAP4} of the form
\begin{align*}
\rho^\alpha = 1+r^\alpha, \qquad v^\alpha =\mathcal{V}+w^\alpha, \qquad U = \mathcal{V}+ u.
\end{align*}

\medskip

Recall that the expected asymptotic velocity for all the phases $\mathcal{W}_M$ has been defined in \eqref{def:WM}: from now on, we will use the lightened notation $\mathcal{W}$.
We introduce the difference
$$ \mathcal{Z}\vcentcolon=\mathcal{W}-\mathcal{V},$$ in the equations for the velocities, so that the perturbation 
$(r^\alpha, w^\alpha, u)$ has then to satisfy
\begin{equation}
\label{eq:VNSpertub}
\left\{  
\begin{aligned}
&\partial_t r^\alpha + \mathcal{V} \cdot \nabla r^\alpha + \Div (w^\alpha)  + \Div (r^\alpha w^\alpha)= 0,\\
&\partial_t w^\alpha + (\mathcal{W} \cdot \nabla) w^\alpha   =  u- w^\alpha - ((w^\alpha-\mathcal{Z}) \cdot \nabla) w^\alpha
+\Gamma(u-w^\alpha)
,\\
&\partial_t u + (\mathcal{W} \cdot \nabla) u - \Delta u + \nabla P =  \int  (w^\alpha-u) \D \mu(\alpha) +\int r^\alpha (w^\alpha-u) \D \mu(\alpha) \\
& \qquad \qquad \qquad \qquad \qquad \qquad  \qquad   
- \int \rho^\alpha \Gamma(u-w^\alpha) \, \mathrm{d}\mu(\alpha)
-((u-\mathcal{Z}) \cdot \nabla) u,\\
& \qquad \qquad \qquad \quad   \qquad\Div(u) = 0.
\end{aligned}
\right.
\end{equation}

\medskip

$\bullet$ \underline{Compressible case \eqref{eq:VNScompressibleCHAP4}}: we look for a solution $(\rho^\alpha, v^\alpha, \varrho, U)$ to \eqref{eq:VNScompressibleCHAP4} of the form
\begin{align*}
    \rho^\alpha=1+r^\alpha, \ \ 
    v^\alpha=\mathcal{V}+w^\alpha, \ \     \varrho=1+n, \ \     
    U =\mathcal{V}+u.
\end{align*}
As in the incompressible case, we write $\mathcal{W}=\mathcal{W}_{m, M}$ and, introducing the final velocity \eqref{def:WmM}
and $$ \mathcal{Z}\vcentcolon=\mathcal{W}_{m,M}-\mathcal{V},$$
the perturbation $(r^\alpha, w^\alpha,  n, u)$ satisfies
\begin{equation}
\label{eq:VNScomp-pertub1}
\left\{  
\begin{aligned}
&\partial_t r^\alpha + \mathcal{V} \cdot \nabla r^\alpha + \Div(w^\alpha)  + \Div (r^\alpha w^\alpha)= 0,\\
&\partial_t w^\alpha + (\mathcal{W} \cdot \nabla) w^\alpha   =  u- w^\alpha - ((w^\alpha-\mathcal{Z}) \cdot \nabla) w^\alpha
+\Gamma(u-w^\alpha),
\\
&\partial_t  n + \mathcal{W} \cdot \nabla  n + \Div(u)  + \Div ( n u)= -(u-\mathcal{Z}) \cdot \nabla  n, \\
&(1+ n)\big(\partial_t u + (\mathcal{W} \cdot \nabla) u + (u \cdot \nabla) u \big)- \mathcal{A} u + \nabla [P(1+ n)] =   \int  (w^\alpha-u) \D \mu(\alpha) \\
  &+\int r^\alpha (w^\alpha-u) \D \mu(\alpha) 
  - \int \rho^\alpha \Gamma(u-w^\alpha) \, \mathrm{d}\mu(\alpha)
   -(1+ n)((u-\mathcal{Z}) \cdot \nabla) u.
\end{aligned}
\right.
\end{equation}

\subsection{Spectral gaps and linear decay}

In this section, we prove the existence of a spectral gap and deduce exponential decay in time for the linearized systems, both for the incompressible and compressible versions of the Vlasov--Navier--Stokes system. This decay only hold for the non-zero modes of the linearized solution and will be quantified in $\dot H^1_x$ norm.  Interestingly, all the results are based on a direct and simple spectral study on the Fourier side. 

\paragraph{Incompressible case.} We focus on the linearized version of the system \eqref{eq:VNSpertub}. We note that, dropping the nonlinear terms, the linearized equations yield a subsystem on the velocities $(w^\alpha,u)$ that is decoupled from the densities $\rho^\alpha$.
\begin{Prop}\label{Prop-lineardecayVNS}
There exists a universal constant $\Lambda \in (0,1/2)$ such that the following holds. Let $T>0$ and $\mathcal{\widetilde{W}} \in \R^3$. Let $(\ww , u)$ satisfy on $[0,T]$ the linearized system
\begin{equation}
\label{eq:lin-VNS}
\left\{  
\begin{gathered}
\partial_t w^\alpha + (\mathcal{\widetilde{W}} \cdot \nabla) w^\alpha =  u- w^\alpha+S_1^\alpha,\\
\partial_t u + (\mathcal{\widetilde{W}} \cdot \nabla) u - \Delta u + \nabla P =  \int (w^\alpha-u) \D \mu(\alpha)+S_2,\\
\Div (u) = 0,
\end{gathered}
\right.
\end{equation}
for some given smooth sources $S_1^\alpha, S_2$. Then the following estimate holds for all $t \in [0,T]$: 
\begin{align}\label{Lin-estimate-Nabla-u}
 \| \nabla j(t) \|_{L^2_x} + \| \nabla u(t) \|_{L^2_x}&\lesssim e^{-\Lambda t}  \left(  \Vert \nabla j_0\|_{ L^2_x} 
 + \|  \nabla u_0 \|_{ L^2_x} \right)    \\
 & \quad + \int_0^t e^{-\Lambda(t-s)} \left (\|  \nabla S_1 (s)  \|_{L^2_x}+\| \nabla S_2(s) \|_{L^2_x}  \right)\D s,
\end{align}
where 
$$ j(t,x)\vcentcolon=\int w^\alpha(t,x) \, \mathrm{d}\mu(\alpha), \ \ j_0\vcentcolon=\int w^\alpha_0 \, \mathrm{d}\mu(\alpha), \ \ S_1\vcentcolon=\int S^\alpha_1 \, \mathrm{d}\mu(\alpha),$$
and for any $p \in [1, \infty]$,
\begin{align}
\label{Lin-estimate-Diffvel-integ1}  \Vert  \langle  \ww(t)- u(t) \rangle \Vert_{L^p_\alpha} &\lesssim e^{- t}   \Vert \langle  \ww_0- u_0 \rangle \Vert_{L^p_\alpha} \\
\notag & \quad + \int_0^t e^{-(t-s)} ( \Vert \langle  S^\alpha_1(s) \rangle \Vert_{L^p_\alpha} + \vert \langle S_2(s)  \rangle \vert) \D s, \\[2mm]
\label{Lin-estimate-Diffvel-integ2}\left\vert \left\langle  \int w^\alpha(t) \, \mathrm{d}\mu(\alpha)- u(t) \right\rangle \right\vert  &\lesssim e^{-2 t}  \left\vert \left\langle  \int w^\alpha_0 \, \mathrm{d}\mu(\alpha)- u_0 \right\rangle \right\vert \\
\notag  &  \quad +\int_0^t e^{-2(t-s)}\left( \left(\int | \langle  S^\alpha_1(s) \rangle | \, \mathrm{d}\mu(\alpha) \right)+ | \langle S_2(s)  \rangle |\right)  \D s.
\end{align}
\end{Prop}

\begin{Rem}
   We emphasize that the previous result holds for any $\widetilde{\mathcal{W}}$, which is a consequence of Galilean invariance.
\end{Rem}

\begin{proof}
First, we choose a natural moving frame which accounts for Galilean invariance by introducing the change of variables
    $x\vcentcolon= z+t\mathcal{\widetilde{W}}$,
whose purpose is to remove the transport term. We then set 
\begin{align*}
    \left(\widetilde{w^\alpha}, \widetilde{u}, \widetilde{S_1^\alpha}, \widetilde{S_2} \right)(t,z) \vcentcolon=(w^\alpha, u,S_1^\alpha, S_2)(t,x).
\end{align*}
By invariance by translation of the norms in~\eqref{Lin-estimate-Nabla-u}--\eqref{Lin-estimate-Diffvel-integ1}--\eqref{Lin-estimate-Diffvel-integ2}, this change of variable does not affect any of the estimates we will derive. Removing the tilde superscripts for readability, we now focus on 
\begin{equation}
\label{eq:lin2-VNS}
\left\{  
\begin{aligned}
\partial_t w^\alpha &=  u- w^\alpha+S_1^\alpha,\\
\partial_t u   - \Delta u + \nabla P &=  \int (w^\alpha-u) \D \mu(\alpha)+S_2,\\
\Div (u) &= 0.
\end{aligned}
\right.
\end{equation}

\medskip

\noindent \textbf{\underline{Step 1}: proof of \eqref{Lin-estimate-Nabla-u}}.

\medskip

We start with the case where the sources are set to zero, that is with $S_1^\alpha=S_2=0$. Defining 
$j\vcentcolon=\int w^\alpha  \D \mu(\alpha)
$, since $\mu$ is of mass $1$, we have 
\begin{align*}
\partial_t \begin{pmatrix}
w^\alpha\\
 u
\end{pmatrix}=\mathbb{M}\begin{pmatrix}
 w^\alpha\\
 u
\end{pmatrix}+ \begin{pmatrix}
0\\
 \mathbb{P}j
\end{pmatrix}, \ \ \mathbb{M}=\begin{bmatrix}
- \mathrm{Id} &  \mathrm{Id} &\\
0 &   - \mathrm{Id}+\Delta
\end{bmatrix},
\end{align*}
where $\mathbb{P}$ stands for the standard Leray projection onto divergence-free vector fields. By integrating the first equation in $\alpha$, we obtain that $(\mathbb{P}j, u)$ satisfies the equation
\begin{equation}
\left\{  
\begin{gathered}
\partial_t (\mathbb{P}j)  =   u- \mathbb{P}j,\\
\partial_t u  -  \Delta u =   \mathbb{P}j-u,
\end{gathered}
\right.
\end{equation}
which can be written as
$\partial_t \begin{pmatrix}
 \mathbb{P}j\\
 u
\end{pmatrix}=\mathbb{A}\begin{pmatrix}
 \mathbb{P}j\\
 u
\end{pmatrix}$,
where the differential operator $\mathbb{A}$ is defined as
\begin{align}\label{def-opdiffA-VNS}
\mathbb{A}\vcentcolon=\begin{bmatrix}
- \mathrm{Id} &  \mathrm{Id} &\\
 \mathrm{Id} & - \mathrm{Id} +  \Delta
\end{bmatrix}.
\end{align}
Using Fourier series on the periodic torus, we are thus led to study the following ODE for all $k \in \Z^3$:
\begin{align*}
\frac{\D}{\D t}\begin{pmatrix}
 \widehat{\mathbb{P}j}_k\\
\widehat{u}_k
\end{pmatrix}
=\mathbb{A}_k\begin{pmatrix}
 \widehat{\mathbb{P}j}_k \\
 \widehat{u}_k
\end{pmatrix}, \ \ 
\mathbb{A}_k\vcentcolon=\begin{bmatrix}
-\mathrm{I}_3 &  \mathrm{I}_3 &\\
 \mathrm{I}_3 & -(1+ \vert k \vert^2)\mathrm{I}_3  
\end{bmatrix}\in M_6(\R),
\end{align*}
with a given initial condition $( \widehat{\mathbb{P}j}_k(0), \widehat{u}_k(0))$.
We observe that for all $k \in \Z^3$, the matrix $\mathbb{A}_k$ is real symmetric, hence diagonalizable  in an orthonormal basis. 
We find that
$$ \mathbb{A}_k=Q_k\mathbb{D}_k Q_k^{-1},$$
where $Q_k\in O_6(\R)$, $\mathbb{D}_k=\mathrm{Diag}(\Lambda_k^-\mathrm{I}_3,\Lambda_k^+ \mathrm{I}_3 )$, with
$$\Lambda_k^\pm\vcentcolon=\frac{1}{2}\left(\pm \sqrt{\vert k \vert^4+4} -\vert k\vert^2 -2\right).$$
We observe that for $k \in \Z^3\setminus\{0\}$, we have 
\begin{align*}\mathrm{Sp}(\mathbb{A}_k)= \left\lbrace \Lambda_k^+, \Lambda_k^-   \right\rbrace \subset \left(-\infty, -\Lambda \right], \ \ \Lambda\vcentcolon=\frac{3-\sqrt{5}}{2}   \in \left (0,\frac{1}{2} \right).
\end{align*}
The endpoint value of the spectrum is reached for $\vert k \vert =1$. We infer that $e^{t \mathbb{A}_k}=Q_k e^{t \mathbb{D}_k}Q_k^{-1}$. Since the solution $X_k(t)=(\widehat{\mathbb{P}j}_k(t), \widehat{u}_k(t))$ of $\dot X_k=\mathbb{A}_k X_k$ is given by $X_k(t)=e^{t\mathbb{A}_k }X_k(0)$ and $Q_k$ is an  orthogonal matrix, we deduce that for all $t \geq 0$
\begin{align}\label{decay-semigroup1st}
\vert X_k(t) \vert= \vert e^{ t \mathbb{A}_k }X_k(0) \vert \lesssim e^{-\Lambda t} \vert X_k(0) \vert,
\end{align}
this estimate being uniform in $k \neq 0$. 

Coming back to $(\mathbb{P}j,u)$, we therefore know that any non-zero mode decays as follows:
\begin{align*}
    \forall k \in \Z^3 \setminus \lbrace 0 \rbrace, \ \ \forall t \geq 0 \ \ ,\left\vert \widehat{\mathbb{P}j}_k(t) \right\vert +  \vert \widehat{u}_k(t) \vert \lesssim e^{-\Lambda t} \left(  \left\vert \widehat{\mathbb{P}j}_k(0) \right\vert+ \vert \widehat{u}_k(0) \vert \right).
\end{align*}
We can therefore estimate for all $t \geq 0$
\begin{align*}
\left\Vert \nabla \mathbb{P}j (t) \right\Vert_{L^2_x}^2+ \| \nabla u(t) \|_{L^2_x}^2 &=\sum_{k \in \Z^3 \setminus \lbrace 0 \rbrace}  \vert k \vert^2  \left( \left\vert \widehat{\mathbb{P}j}_k(t) \right\vert^2+\vert \widehat{u}_k(t) \vert^2 \right) \\
&\lesssim e^{-2 \Lambda t} \sum_{k \in \Z^3 \setminus \lbrace 0 \rbrace}  \vert k \vert^2  \left( \left\vert \widehat{\mathbb{P}j}_k(0) \right\vert^2+\vert \widehat{u}_k(0) \vert^2 \right),
\end{align*}
from which we infer 
\begin{align*}
   \left\Vert \nabla \mathbb{P}j (t) \right\Vert_{L^2_x}+  \| \nabla u(t) \|_{L^2_x} \lesssim e^{-\Lambda t} \left( \| \nabla j(0) \|_{L^2_x}+ \| \nabla u(0) \|_{L^2_x} \right).
\end{align*}
When the sources $S_1^\alpha$ and $S_2$ are not zero, the same computations as before yield
\begin{align*}
\partial_t \begin{pmatrix}
 \mathbb{P}j\\
 u
\end{pmatrix}=\mathbb{A}\begin{pmatrix}
 \mathbb{P}j\\
 u
\end{pmatrix}
+\begin{pmatrix}
 \mathbb{P} S_1\\
  \mathbb{P} S_2
\end{pmatrix},
\end{align*}
where $S_1=\int S_1^\alpha \, \mathrm{d}\mu(\alpha)$,
so that by Duhamel formula, we have for any $k \neq 0$ 
\begin{align*}
   \begin{pmatrix}
 \widehat{\mathbb{P}j}_k (t)\\
 \widehat{u}_k(t)
\end{pmatrix} =e^{tA_k} \begin{pmatrix}
 \widehat{\mathbb{P}j}_k(0)\\
 \widehat{u}_k(0)
\end{pmatrix}
+
\int_0^t e^{(t-\tau)A_k}\begin{pmatrix}
 \widehat{\mathbb{P} S_1}_k(\tau)\\
 \widehat{\mathbb{P} S_2}_k(\tau)
\end{pmatrix} \, \mathrm{d}\tau.
\end{align*}
 According the first step of the proof, we know that for all $k \neq 0$, we have 
 \begin{align*}
     \left\lvert e^{t A_k} 
     \begin{pmatrix}
 \widehat{\mathbb{P}j}_k(0)\\
 \widehat{u}_k(0) 
 \end{pmatrix}\right\rvert &\lesssim e^{-\Lambda t} \left(  \left\vert \widehat{\mathbb{P}j}_k(0) \right\vert+ \vert \widehat{\mathbb{P} u}_k(0) \vert \right), \ \ t \geq 0, \\
 \left\lvert e^{(t-\tau) A_k} 
     \begin{pmatrix}
 \widehat{\mathbb{P} S_1}_k(\tau)\\
 \widehat{S_2}_k(\tau) 
 \end{pmatrix}\right\rvert &\lesssim e^{-\Lambda (t-\tau)} \left(  \left\vert  \widehat{\mathbb{P} S_1}_k(\tau) \right\vert+ \vert  \widehat{ \mathbb{P} S_2}_k(\tau) \vert \right), \ \ 0 \leq \tau \leq t.
 \end{align*}
 Similarly, we obtain after summation over the frequencies $k \neq 0$:
 \begin{multline*}
    \Vert \nabla \mathbb{P} j(t) \Vert_{L^2_x}+ \| \nabla u(t) \|_{L^2_x}\lesssim e^{-\Lambda t}  \left(  \Vert \nabla j_0\|_{ L^2_x} + \|  \nabla u_0 \|_{ L^2_x} \right)    
 \\
 + \int_0^t e^{-\Lambda(t-s)} \left (\|  \nabla S_1 (s)  \|_{L^2_x}+\| \nabla S_2(s) \|_{L^2_x}  \right)\D s,
\end{multline*}
which is the desired \eqref{Lin-estimate-Nabla-u}, \textit{modulo} the fact that it holds for $\mathbb{P}j$ instead of $j$ in the left-hand side. However, the orthogonal part $\mathbb{P}^\perp j$ (where $\mathbb{P}^\perp$ is the orthogonal projection on potential vector field) satisfies
\begin{align*}
    \partial_t \mathbb{P}^\perp j=-\mathbb{P}^\perp j+\mathbb{P}^\perp S_1,
\end{align*}
since $u$ is divergence-free, and therefore one directly obtains the estimate
\begin{align*}
    \Vert \nabla \mathbb{P}^\perp j(t) \Vert_{L^2_x} \leq e^{-t} \Vert \nabla j_0 \Vert_{L^2_x}+ \int_0^t e^{-(t-s)} \Vert \nabla S_1(s) \Vert_{L^2_x} \, \mathrm{d}s.
\end{align*}
Since $\Lambda<1$, we eventually infer the full estimate \eqref{Lin-estimate-Nabla-u}.

\medskip

\noindent \textbf{\underline{Step 2}: proof of \eqref{Lin-estimate-Diffvel-integ1} and \eqref{Lin-estimate-Diffvel-integ2}}.

\medskip

For the estimate on the mean $ \langle w^\alpha-u \rangle$, we use the fact that $w^\alpha-u$ is solution to the equation
$$
\partial_t (w^\alpha-u)   + \Delta u - \nabla p =  u- w^\alpha - \int (w^\beta-u) \D \mu(\beta) + S^\alpha_1 - S_2,
$$
and thus the mean $ \langle w^\alpha-u \rangle$ satisfies for any phase $\alpha$
\begin{equation}
\label{eq:mean}
\partial_t \langle w^\alpha-u \rangle  = - \langle  w^\alpha-u \rangle - \int \langle w^\beta-u \rangle \D \mu(\beta) +   \langle S^\alpha_1 \rangle - \langle S_2 \rangle .
\end{equation}
In particular,
$$
\frac{\mathrm{d}}{\mathrm{d}t} \int \langle w^\alpha-u\rangle  \D \mu(\alpha) =- 2 \int \langle w^\alpha-u \rangle \D \mu(\alpha) + \int \left( \langle S^\alpha_1 \rangle - \langle S_2 \rangle \right)\D \mu(\alpha).
$$
Therefore, we obtain
\begin{multline}
\label{eq:estimmeanmu}
\left|\int \langle w^\alpha-u \rangle  \D \mu(\alpha) (t)\right | \lesssim e^{-2 t} \left| \int \langle w^\alpha_0-u_0 \rangle  \D \mu(\alpha)\right|  \\
+ \int_0^t e^{-2(t-s)} \left| \int \left( \langle S^\alpha_1 \rangle - \langle S_2 \rangle \right)\D \mu(\alpha)\right| \D s,
\end{multline}
from which we infer \eqref{Lin-estimate-Diffvel-integ2}. From \eqref{eq:mean}, we find that
\begin{multline*}
    \langle w^\alpha-u \rangle  (t)=e^{-t}\langle w^\alpha_0-u_0 \rangle  \\
    +\int_0^t e^{-(t-s)} \left( - \int \langle w^\beta-u \rangle(s) \D \mu(\beta) +   \langle S^\alpha_1(s) \rangle - \langle S_2(s) \rangle\right)  \, \mathrm{d}s,
    \end{multline*}
and plugging in~\eqref{eq:estimmeanmu}, this gives
\begin{multline}
\label{eq:estimmean}
\left| \langle w^\alpha-u \rangle  (t) \right| \lesssim e^{-t} \left|\langle w^\alpha_0-u_0 \rangle\right| +  e^{- t} \left| \int \langle w^\beta_0-u_0 \rangle  \D \mu(\beta)\right| \\ + \int_0^t e^{-(t-s)}   \left(\left|   \langle S^\beta_1 \rangle - \langle S_2 \rangle \right| + \int \left|   \langle S^\beta_1 \rangle - \langle S_2 \rangle \right| \, \mathrm{d} \mu(\beta) \right)(s) \D s.
\end{multline}
By integrating in $\alpha$ and applying the Hölder inequality, this yields the estimate \eqref{Lin-estimate-Diffvel-integ1}, hence concluding the proof. 

\end{proof}

\begin{Rem}\label{rem-linsystem-incom-massnot1}
  In the case where we rather write $\rho^\alpha=\mathrm{C}(1+r^\alpha)$ with a constant $\mathrm{C}>0$, the equivalent of  the linearized system \eqref{eq:lin-VNS} (after the Galilean transformation) is
    \begin{equation}
\label{eq:lin2-VNS-massC}
\left\{  
\begin{aligned}
\partial_t w^\alpha &=  u- w^\alpha,\\
\partial_t u   - \Delta u + \nabla P &=  \mathrm{C}\int (w^\alpha-u) \D \mu(\alpha),\\
\Div (u) &= 0.
\end{aligned}
\right.
\end{equation}
Introducing the unknown $j^{\mathrm{C}}\vcentcolon=\sqrt{C}\int w^\alpha  \D \mu(\alpha)$ in order so symmetrize the system, we obtain the fact that $(\mathbb{P}j^{\mathrm{C}}, u)$ satisfies for all $k \in \Z^3$: 
\begin{align*}
\frac{\D}{\D t}
\begin{pmatrix}
 \widehat{\mathbb{P}j^{\mathrm{C}}}_k\\
\widehat{u}_k
\end{pmatrix}
=\mathbb{A}_k^{\mathrm{C}}\begin{pmatrix}
 \widehat{{\mathbb{P}j^{\mathrm{C}}}}_k \\
 \widehat{u}_k
\end{pmatrix}, \ \ \mathbb{A}_k^{\mathrm{C}}\vcentcolon=
\begin{bmatrix}
-\mathrm{I}_3 &  -\sqrt{\mathrm{C}}\mathrm{I}_3 \\
 \sqrt{\mathrm{C}}\mathrm{I}_3 & -(\mathrm{C}+ \vert k \vert^2)\mathrm{I}_3  
\end{bmatrix}\in M_6(\R).
\end{align*}
One can check that the argument from above for the case $\mathrm{C}=1$ still holds true. Indeed, the eigenvalue of the symmetric matrix $\mathbb{A}_k^{\mathrm{C}}$ are now
$$\Lambda_k^{\mathrm{C},\pm}\vcentcolon=\frac{1}{2}\left(\pm \sqrt{(\vert k \vert^2+\mathrm{C}-1)^2+4 \mathrm{C})} -\vert k\vert^2 -1-\mathrm{C}\right).$$
and one similarly obtain an explicit uniform spectral gap since
\begin{align*}
    \max_{\vert k \vert \geq 1} \Lambda_k^{\mathrm{C},\pm} =-\frac{2}{\mathrm{C}+2+\sqrt{\mathrm{C}^2+4\mathrm{C}}}<0.
\end{align*} 
The decay of the non-zero modes then follows. For the decay of the difference of the average, we now have (restauring the source term) $$
\frac{\mathrm{d}}{\mathrm{d}t} \int \langle w^\alpha-u\rangle  \D \mu(\alpha) =- (1+\mathrm{C}) \int \langle w^\alpha-u \rangle \D \mu(\alpha) + \int \left( \langle S^\alpha_1 \rangle - \langle S_2 \rangle \right)\D \mu(\alpha).
$$
so the argument from the case $\mathrm{C}=1$ still holds true.
\end{Rem}

\paragraph{Compressible case.}
Next, we turn to the spectral study of the linearized version of \eqref{eq:VNScomp-pertub1}, for which the same type of proof can be performed. Let us recall that, by assumption, we have assumed that $P'(1)>0$.
\begin{Prop}\label{Prop-lineardecayVNS-comp}
There exists a universal constant $\Lambda \in (0,1)$ such that the following holds. Let $T>0$ and $\mathcal{\widetilde{W}} \in \R^3$. Let $(\ww , u)$ satisfy on $[0,T]$ the linearized system
\begin{equation}
\left\{ 
\begin{aligned}
\partial_t w^\alpha + (\mathcal{\widetilde{W}}\cdot \nabla) w^\alpha &= u - w^\alpha + S_1^\alpha,\\
\partial_t  n + \mathcal{\widetilde{W}} \cdot \nabla  n+\Div ( u)&=S_3, \\
\partial_t u + (\mathcal{\widetilde{W}}\cdot \nabla) u+ P'(1)\nabla  n- \mathcal{A}u &= \int \big( w^\alpha - u ) \D \mu(\alpha)+S_2,
\end{aligned}
\right.
\end{equation}
for some given smooth sources $S_1^\alpha, S_2$ and $S_3$, and where we recall that $\mathcal{A}u=\Delta u+ 2 \nabla \Div (u)$. Then the following estimates hold for all $t \in [0,T]$:  
\begin{align}\label{Lin-estimate-Nabla-utheta-comp}
\begin{split}
 &\| \nabla j(t) \|_{L^2_x} +   \| \nabla  n(t) \|_{L^2_x}+\| \nabla u(t) \|_{L^2_x}
  \\
  &\lesssim e^{-\Lambda t}  \left(  \| \nabla  n_0 \|_{L^2_x}+ \Vert \nabla  j_0\|_{ L^2_x} + \|  \nabla u_0 \|_{ L^2_x} \right)    \\
      & \quad + \int_0^t e^{-\Lambda(t-s)} \left (\|  \nabla  S_1 (s)  \|_{L^2_x}+\| \nabla S_2(s) \|_{L^2_x}+\| \nabla S_3(s) \|_{L^2_x}   \right)\D s,
      \end{split}
\end{align}
where
$$  j(t,x)\vcentcolon=\int w^\alpha(t,x) \, \mathrm{d}\mu(\alpha), \ \ j_0\vcentcolon=\int w^\alpha_0 \, \mathrm{d}\mu(\alpha), \ \ S_1\vcentcolon=\int S^\alpha_1 \, \mathrm{d}\mu(\alpha),$$
and for any $p \in [1, \infty]$, 
\begin{align}
\label{Lin-estimate-Diffvel-integ1-comp}  \Vert  \langle  w^\alpha(t)- u(t) \rangle \Vert_{L^p_\alpha} &\lesssim e^{- t}   \Vert \langle  w^\alpha_0- u_0 \rangle \Vert_{L^p_\alpha} \\
\notag &  \quad + \int_0^t e^{-(t-s)} ( \Vert \langle  S^\alpha_1(s) \rangle \Vert_{L^p_\alpha} + \vert \langle S_2(s)  \rangle \vert) \D s, \\[2mm]
\label{Lin-estimate-Diffvel-integ2-comp}\left\vert \left\langle  \int w^\alpha(t) \, \mathrm{d}\mu(\alpha)- u(t) \right\rangle \right\vert  &\lesssim e^{-2 t}  \left\vert \left\langle  \int w^\alpha_0 \, \mathrm{d}\mu(\alpha)- u_0 \right\rangle \right\vert \\
\notag &  \quad +\int_0^t e^{-2(t-s)}\left( \left(\int_I | \langle  S^\alpha_1(s) \rangle | \, \mathrm{d}\mu(\alpha) \right)+ | \langle S_2(s)  \rangle \vert\right)  \D s.
\end{align}
\end{Prop}
\begin{proof}
For the estimates \eqref{Lin-estimate-Diffvel-integ1-comp} and \eqref{Lin-estimate-Diffvel-integ2-comp}  bearing on averages, we observe that $\langle u-w^\alpha \rangle$ satisfies the equation 
\begin{equation}
\label{eq:mean-comp}
\partial_t \langle w^\alpha-u \rangle  = -P'(1)\nabla  n- \langle  w^\alpha-u \rangle - \int \langle w^\beta-u \rangle \D \mu(\beta) +   \langle S^\alpha_1 \rangle - \langle S_2 \rangle,
\end{equation}
and therefore, since $\langle \nabla  n \rangle=0$, the same proof as in Proposition \ref{Prop-lineardecayVNS} for the incompressible case applies.

We now focus on the estimate \eqref{Lin-estimate-Nabla-utheta-comp}. 
Let us mention the fact that the analysis performed below is reminiscent of the spectral study performed in \cite{GWZ} for the (one phase) pressureless-Euler equations coupled with the compressible Navier--Stokes equations (but for the case of the whole space $\R^3$).

\medskip

\noindent \textbf{\underline{Step 1}: Projection on incompressible and compressible part}.

\medskip

By Galilean invariance, we can again assume that $\widetilde{\mathcal{W}}=0$ and we start by  setting the source terms equal to $0$. Denote by $\mathbb{P}$ and $\mathbb{P}^\perp$ the orthogonal projectors on divergence-free and potential vector fields. Since $\mathcal{A}u=\Delta u+ 2 \nabla \Div (u)$, we have 
\begin{equation*}
\left\{ 
\begin{aligned}
\partial_t w^\alpha  &= u - w^\alpha,\\
\partial_t  n  +\Div ( \mathbb{P}^\perp u)&=0, \\
\partial_t \mathbb{P}u - \Delta \mathbb{P}u &= \mathbb{P}\int w^\alpha \D \mu(\alpha)-\mathbb{P}u, \\
\partial_t \mathbb{P}^\perp u + P'(1)\nabla  n- 3\Delta \mathbb{P}^\perp u &=\mathbb{P}^\perp \int w^\alpha  \D \mu(\alpha)-\mathbb{P}^\perp u.
\end{aligned}
\right.
\end{equation*}
Indeed, we have $u=\mathbb{P}u+ \nabla g$  where $\nabla g=\mathbb{P}^\perp u$ so that $$\nabla \, \mathrm{div}(u)=\nabla \, \mathrm{div}(\nabla g)=\nabla \Delta g=\Delta \mathbb{P}^\perp u.$$
By defining 
$$j\vcentcolon=\int w^\alpha  \D \mu(\alpha),
$$
we first obtain the following two decoupled systems on $(\mathbb{P}u, \mathbb{P}j)$ and $( n, \mathbb{P}^\perp u, \mathbb{P}^\perp j)$:
\begin{equation}\label{eq:NSCLin-incomp}
\left\{ 
\begin{aligned}
\partial_t \mathbb{P}j  &= \mathbb{P}u - \mathbb{P}j, \\
\partial_t \mathbb{P}u - \Delta \mathbb{P}u &= \mathbb{P}j-\mathbb{P}u,
\end{aligned}
\right.
\end{equation}
and 
\begin{equation}\label{eq:NSCLin-comp}
\left\{ 
\begin{aligned}
\partial_t  n +\Div ( \mathbb{P}^\perp u)&=0, \\
\partial_t \mathbb{P}^\perp u + P'(1)\nabla  n- 3\Delta \mathbb{P}^\perp u &=\mathbb{P}^\perp j-\mathbb{P}^\perp u, \\
\partial_t \mathbb{P}^\perp j&  = \mathbb{P}^\perp u - \mathbb{P}^\perp j.
\end{aligned}
\right.
\end{equation}

\medskip

\noindent \textbf{\underline{Step 2}: study of the incompressible part}.

\medskip

The incompressible macroscopic part $(\mathbb{P}u, \mathbb{P}j)$  satisfies the same kind of system \eqref{eq:NSCLin-incomp} as previously studied in the incompressible Navier--Stokes case: more precisely, we have
    \begin{align*}
\partial_t \begin{pmatrix}
 \mathbb{P}j\\
 \mathbb{P}u
\end{pmatrix}=\mathbb{A}\begin{pmatrix}
 \mathbb{P}j\\
 \mathbb{P}u
\end{pmatrix},
\end{align*}
where the differential operator $\mathbb{A}$ is defined in \eqref{def-opdiffA-VNS}. As a consequence, thanks to (the proof of) Proposition \ref{Prop-lineardecayVNS}, we infer for all $x \in \R^6$ and all $t \geq 0$
$$
\vert e^{t \mathbb{A}_k}x \vert  \lesssim e^{- t  \Lambda_i}\vert x \vert,
$$
for some $\Lambda_i>0$, uniformly in $k\in\Z^3 \setminus\{0\}$. 

\medskip

\noindent \textbf{\underline{Step 3}: study of the compressible part}.

\medskip

Let us focus on the compressible macroscopic part $( n, \mathbb{P}^\perp u, \mathbb{P}^\perp j)$ which satisfies \eqref{eq:NSCLin-comp}. We now consider the following scalar quantities:
\begin{align*}
\mathfrak{U}\vcentcolon=\vert D \vert^{-1} \mathrm{div} (\mathbb{P}^\perp u), \ \ \mathfrak{J}\vcentcolon=\vert D \vert^{-1} \mathrm{div} (\mathbb{P}^\perp j), \ \ \widehat{(\vert D \vert ^{-1}g)_k}\vcentcolon=\vert k \vert^{-1} \widehat{g}_k.
\end{align*}
Since these new variables are of the same order in terms of number of derivatives, it is enough\footnote{and the same holds in the case where sources are present.} to estimate $( n, \mathfrak{U}, \mathfrak{J})$ instead of $( n, u, j)$. 
Applying $\vert D \vert^{-1} \mathrm{div} $ to the equation on $ \mathbb{P}^\perp u$ and $ \mathbb{P}^\perp j$, we obtain the following system in $\R^3$:
\begin{equation}
\left\{ 
\begin{aligned}
\partial_t  n +\vert D \vert \mathfrak{U}&=0, \\
\partial_t \mathfrak{U}  -P'(1)\vert D \vert  n-3 \Delta \mathfrak{U} &=\mathfrak{J}-\mathfrak{U}, \\
\partial_t \mathfrak{J} &  = \mathfrak{U}-\mathfrak{J}.
\end{aligned}
\right.
\end{equation}
For all $k \in \Z^3$, the $k--$th Fourier-mode satisfy
\begin{align*}
\frac{\mathrm{d}}{\mathrm{d}t} \begin{pmatrix} \widehat{ n}_k \\  \widehat{\mathfrak{U}}_k \\  \widehat{\mathfrak{J}}_k \\ \end{pmatrix}
= \mathbb{C}_k \begin{pmatrix} \widehat{ n}_k \\  \widehat{\mathfrak{U}}_k \\  \widehat{\mathfrak{J}}_k \\ \end{pmatrix}, \ \ 
\mathbb{C}_k=\begin{bmatrix}
0 & -\vert k \vert &  0\\
P'(1)\vert k \vert  & -1-3\vert k \vert^2  & 1 \\
0 & 1 & -1 
\end{bmatrix} \in  \mathrm{M}_{3}(\R).
\end{align*}
Compared to the incompressible part, a difficulty arises due to the fact that the matrix $\mathbb{C}_k$ is not symmetric and it is therefore not so straightforward to obtain some exponential decay estimate that is uniform\footnote{A careful study of the characteristic polynomial of $\mathbb{C}_k$ actually shows that $\mathbb{C}_k$ has three distinct real eigenvalues (for $\vert k \vert >\sqrt{2}$, the case $\vert k \vert=1$ having to be treated separately) and is therefore diagonalizable. However, on has to prove that the associated projectors are uniformly bounded with respect to the frequency $k$, which requires additional argument. One possible strategy (which we do not pursue in this work) would be to rely on perturbation theory at high frequency.} in the frequency $k$. 

\medskip

Hence, we adopt a more direct strategy based on the construction of a suitable Lyapunov function for the previous system of ODEs, in the spirit of hypocoercivity methods (see among others \cite{Herau1,HerauNier,Villani,MouhotNeumann,DolbeaultMouhotSchmeiser}). This procedure is actually quite standard for partially dissipative systems \cite{danchin2023-dissipativereview}. 

In what follows, we fix $k \in \Z^3 \setminus \lbrace 0 \rbrace$ (so that $\vert k \vert \geq 1$). For the sake of readability, and to  symmetrize the problem, we set
\begin{align*}
    (\mathcal{N}, \mathcal{U}, \mathcal{J})=(\sqrt{c}\widehat{ n}_k,\widehat{\mathfrak{U}}_k,\widehat{\mathfrak{J}}_k), \ \ c=P'(1)>0,
\end{align*}
 where, for the sake of readability, we intentionally do not write the dependency with respect to the frequency $k$ in index. We therefore seek exponential decay for $( \mathcal{N}, \mathcal{U}, \mathcal{J})(t)$, which satisfies the following system of ODES:
\begin{equation}
    \label{eq:ODES-vns-comp}
\begin{aligned}
    \frac{\mathrm{d}}{\mathrm{d}t} \begin{pmatrix}  \mathcal{N} \\  \mathcal{U} \\  \mathcal{J} \\ \end{pmatrix}
=\begin{pmatrix} -\sqrt{c} \vert k \vert \mathcal{U} \\  \sqrt{c}\vert k \vert  \mathcal{N}-(1+3 \vert k \vert^2)\mathcal{U} + \mathcal{J}\\  \mathcal{U}-\mathcal{J}.\end{pmatrix}.
\end{aligned}
\end{equation}
Let $\eta \in (0,1)$ be a parameter that will be fixed later on. We introduce the functional
\begin{align*}
    E_k(t)\vcentcolon=\frac{1}{2}\big(\vert  \mathcal{N}(t) \vert^2+\vert \mathcal{U}(t) \vert^2+\vert \mathcal{J}(t) \vert^2 \big)-\eta \gamma_k \Re\left( \mathcal{N}(t) \overline{\mathcal{U}(t)}\right), \ \ t \geq 0,
\end{align*}
where
\begin{align*}
    \gamma_k\vcentcolon=\sqrt{c}\frac{ \vert k \vert}{1+3\vert k \vert^2} \in (0, \sqrt{c}/4]. 
\end{align*}
First, we observe that $E_k$ and $\vert ( \mathcal{N},\mathcal{U},\mathcal{J})\vert^2$ are equivalent, in the sense that for all $t \geq 0$
\begin{align}\label{eq:equiv-norm-VNSlincom}
    \frac{1}{4}\big(\vert  \mathcal{N}(t) \vert^2+\vert \mathcal{U}(t) \vert^2+\vert \mathcal{J}(t) \vert^2 \big) \leq E_k(t) \leq \frac{3}{4} \big(\vert  \mathcal{N}(t) \vert^2+\vert \mathcal{U}(t) \vert^2+\vert \mathcal{J}(t) \vert^2 \big).
\end{align}
Indeed, since $\vert \Re( \mathcal{N} \overline{\mathcal{U}}) \vert  \leq \frac{1}{2}(\vert  \mathcal{N} \vert^2+\vert \mathcal{U} \vert^2)$ and $\eta \gamma_k \leq \eta \sqrt{c}/4 \leq 1/16$ for $\eta \in (0, 1/4\sqrt{c}]$, one has
\begin{align*}
    E_k &\geq \frac{1}{2}\big(\vert  \mathcal{N} \vert^2+\vert \mathcal{U} \vert^2+\vert \mathcal{J} \vert^2 \big)-\eta \gamma_k \frac{1}{2}(\vert  \mathcal{N} \vert^2+\vert \mathcal{U} \vert^2) \\
    & =\frac{1}{2}(1-\eta \gamma_k)(\vert  \mathcal{N} \vert^2+\vert \mathcal{U} \vert^2)+ \frac{1}{2}\vert \mathcal{J} \vert^2 \\
    & \geq \frac{1}{4}\big(\vert  \mathcal{N} \vert^2+\vert \mathcal{U} \vert^2+\vert \mathcal{J} \vert^2 \big),
\end{align*}
and 
\begin{align*}
    E_k \leq \frac{1}{2}\left(1+\frac{1}{2}\eta \gamma_k\right)\big(\vert  \mathcal{N} \vert^2+\vert \mathcal{U} \vert^2+\vert \mathcal{J} \vert^2 \big) \leq \frac{3}{4}\big(\vert  \mathcal{N} \vert^2+\vert \mathcal{U} \vert^2+\vert \mathcal{J} \vert^2 \big).
\end{align*}
The goal is now to obtain, for a suitable choice of $\eta$ (independent of $k$) an energy inequality of form
\begin{align}\label{eq:ineqdiff-hypocercive}
    \frac{\mathrm{d}}{\mathrm{d}t}E_k(t)+\Lambda E_k(t) \leq 0, \ \  t>0,
\end{align}
where $\Lambda>0$ does not depend on $k$. If such an estimate were to hold, then
\begin{align*}
    E_k(t) \leq  e^{-\Lambda t}E_k(0),
\end{align*}
and in view of \eqref{eq:equiv-norm-VNSlincom}, we would obtain
\begin{align}\label{hypoco-ccl1}
    \vert  \mathcal{N}(t) \vert^2+\vert \mathcal{U}(t) \vert^2+\vert \mathcal{J}(t) \vert^2 \lesssim e^{-\Lambda t} \big(\vert  \mathcal{N}(0) \vert^2+\vert \mathcal{U}(0) \vert^2+\vert \mathcal{J}(0) \vert^2 \big).
\end{align}
Differentiating the functional $E_k$, we have
\begin{align*}
    \frac{\mathrm{d}}{\mathrm{d}t}E_k=\Re\left(\frac{\mathrm{d}}{\mathrm{d}t} \mathcal{N} \overline{ \mathcal{N}}\right)+\Re\left(\frac{\mathrm{d}}{\mathrm{d}t}\mathcal{U} \overline{\mathcal{U}}\right)+\Re\left(\frac{\mathrm{d}}{\mathrm{d}t}\mathcal{J} \overline{\mathcal{J}}\right)\\
    -\eta \gamma_k\Re\left(\frac{\mathrm{d}}{\mathrm{d}t} \mathcal{N} \overline{\mathcal{U}}\right) -\eta \gamma_k \Re\left(\frac{\mathrm{d}}{\mathrm{d}t}\overline{\mathcal{U}}  \mathcal{N}\right), 
\end{align*}
and get after simplification, using the system~\eqref{eq:ODES-vns-comp}:
\begin{align*}
    \frac{\mathrm{d}}{\mathrm{d}t}E_k&=\Big[ -(1+3 \vert k \vert^2) + \sqrt{c}\eta \gamma_k \vert k \vert \Big]\vert\mathcal{U} \vert^2 -\vert \mathcal{J} \vert^2 - \sqrt{c} \eta \vert k \vert \gamma_k \vert  \mathcal{N} \vert^2  \\
    & \quad 
    +2 \Re(\mathcal{J} \overline{\mathcal{U}})
    + \sqrt{c}\eta \vert k \vert \Re( \mathcal{N} \overline{\mathcal{U}})
    -\eta \gamma_k \Re( \mathcal{N} \overline{\mathcal{J}}).
\end{align*}
We now handle the three cross-terms thanks to Young inequality: we have
\begin{align*}
    2 \Re(\mathcal{J} \overline{\mathcal{U}}) &\leq 2 \vert \mathcal{U} \vert^2 + \frac{\vert \mathcal{J} \vert^2}{2}, \\
    \sqrt{c}\eta \vert k \vert \Re( \mathcal{N} \overline{\mathcal{U}})& \leq \frac{c\eta^2 \vert k \vert^2}{\sqrt{c}\eta \gamma_k \vert k \vert} \frac{\vert \mathcal{U} \vert^2}{2}+\sqrt{c}\eta \gamma_k \vert k \vert \frac{\vert  \mathcal{N} \vert^2}{2}\\
    &=\eta (1+3 \vert k \vert^2) \frac{\vert \mathcal{U} \vert^2}{2}+\sqrt{c}\eta \gamma_k \vert k \vert \frac{\vert  \mathcal{N} \vert^2}{2}, \\
    -\eta \gamma_k \Re( \mathcal{N} \overline{\mathcal{J}}) & \leq \frac{\sqrt{c}\eta \gamma_k \vert k \vert}{4} \vert  \mathcal{N} \vert^2+ \frac{\eta \gamma_k}{\sqrt{c} \vert k \vert} \vert \mathcal{J} \vert^2 \leq \frac{\sqrt{c}\eta \gamma_k \vert k \vert}{4} \vert  \mathcal{N} \vert^2+ \frac{\eta }{4} \vert \mathcal{J} \vert^2,
\end{align*}
and plugging these estimates back in the former differential equality, we obtain after simplification
\begin{multline*}
    \frac{\mathrm{d}}{\mathrm{d}t}E_k\leq \left[ \left(\frac{\eta}{2}-1 \right)(1+3 \vert k \vert^2) + \sqrt{c}\eta \gamma_k \vert k \vert+2\right]\vert \mathcal{U} \vert^2  +\left[\frac{\eta}{4}-\frac{1}{2} \right]\vert \mathcal{J} \vert^2 \\
    -\frac{\sqrt{c}\eta \gamma_k \vert k \vert}{4} \vert  \mathcal{N} \vert^2.
\end{multline*}
Let us choose $\eta=\frac{1}{8(1+c)} \leq \frac{1}{8}$ (so that $\eta c \leq \frac{1}{8}$): since we have
\begin{align*}
     \sqrt{c}\eta \gamma_k \vert k \vert=\eta c\frac{ \vert k \vert^2}{1+3\vert k \vert^2} \leq \frac{\eta c}{3} \leq \eta c \leq \frac{1}{8}, \ \  \vert k \vert \geq 1,
\end{align*} 
and 
\begin{align*}
    \sqrt{c}\eta \gamma_k \vert k \vert=\eta c\frac{ \vert k \vert^2}{1+3\vert k \vert^2}  \geq \frac{\eta c}{4}=\frac{1}{32}\frac{c}{1+c},
\end{align*}
we obtain for $\vert k \vert \geq 1$
\begin{align*}
    \frac{\mathrm{d}}{\mathrm{d}t}E_k&\leq\left[ -\frac{15}{16}(1+3 \vert k \vert^2) + \frac{1}{8} +2\right]\vert \mathcal{U} \vert^2  -\frac{15}{32} \vert \mathcal{J} \vert^2 
    -\frac{\sqrt{c}\eta \gamma_k \vert k \vert}{4} \vert  \mathcal{N} \vert^2 \\
    & \leq \left[ -\frac{15}{16} \cdot 4 + \frac{1}{8}  + 2 \right] \vert \mathcal{U} \vert^2  -\frac{15}{32} \vert \mathcal{J} \vert^2 
    -\frac{1}{32}\frac{c}{1+c} \frac{1}{4} \vert  \mathcal{N} \vert^2.
\end{align*}
As a consequence, there exists a constant $\lambda>0$ depending on $c$ and independent of $k$ such that
\begin{align*}
    \frac{\mathrm{d}}{\mathrm{d}t}E_k \leq -\lambda \big(\vert  \mathcal{N} \vert^2+\vert \mathcal{U} \vert^2+\vert \mathcal{J} \vert^2 \big),
\end{align*}
and thanks to the equivalence of the norms \eqref{eq:equiv-norm-VNSlincom}, we can eventually infer that an estimate of the type \eqref{eq:ineqdiff-hypocercive} holds. 

In view of \eqref{hypoco-ccl1}, we have therefore proved the following: there exists $\Lambda_c>0$ such that for any $k \in \Z^3 \setminus \lbrace 0 \rbrace$ and $t \geq 0$, we have 
\begin{align}\label{hypoco-ccl1-2}
    \vert  n_k(t) \vert+\vert \mathfrak{U}_k(t) \vert+\vert \mathfrak{J}_k(t) \vert \lesssim e^{-\Lambda_c t} \big(\vert  n_k(0) \vert+\vert \mathfrak{U}_k(0) \vert+\vert \mathfrak{J}_k(0) \vert \big).
\end{align}
In terms of $( n, \mathbb{P}^\perp u, \mathbb{P}^\perp j)$, we therefore have
\begin{align*}
    \left\vert \left(\widehat{ n}, \widehat{\mathbb{P}^\perp u}, \widehat{\mathbb{P}^\perp j}\right)_k (t) \right\vert &\leq C e^{-\Lambda_c t}  \left\vert \left(\widehat{ n}, \widehat{\mathbb{P}^\perp u}, \widehat{\mathbb{P}^\perp j}\right)_k (0) \right\vert.
\end{align*}

\medskip

\noindent \textbf{\underline{Conclusion}: decay estimates}.

\medskip

Coming back to $( n, j,u)$, gathering all pieces together, we have thus proven that by taking $\Lambda=\min( \Lambda_c, \Lambda_i)>0$, the following holds: for all $t \geq 0$, the solution $(\mathbb{P}u, \mathbb{P}j)$ to \eqref{eq:NSCLin-incomp} and the solution $( n, \mathbb{P}^\perp u, \mathbb{P}^\perp j)$ to \eqref{eq:NSCLin-comp} have Fourier modes with the following exponential linear decay for all $k \in \Z^3 \setminus \lbrace 0 \rbrace$
\begin{align*}
\left\vert \left(\widehat{\mathbb{P}u}, \widehat{\mathbb{P}j}\right)_k (t) \right\vert &\leq C e^{-\Lambda t}  \left\vert \left(\widehat{\mathbb{P}u}, \widehat{\mathbb{P}j}\right)_k (0) \right\vert, \\
\left\vert \left(\widehat{ n}, \widehat{\mathbb{P}^\perp u}, \widehat{\mathbb{P}^\perp j}\right)_k (t) \right\vert &\leq C e^{-\Lambda t}  \left\vert \left(\widehat{ n}, \widehat{\mathbb{P}^\perp u}, \widehat{\mathbb{P}^\perp j}\right)_k (0) \right\vert.
\end{align*}
for some universal constant $C>0$. Arguing as in the proof of Proposition \ref{Prop-lineardecayVNS}, we can deduce that for all $t \geq 0$
\begin{align*}
   \left\Vert \nabla \mathbb{P}j (t) \right\Vert_{L^2_x}+  \| \nabla \mathbb{P} u(t) \|_{L^2_x} &\lesssim e^{-\Lambda t} \left( \| \nabla j(0) \|_{L^2_x}+ \| \nabla u(0) \|_{L^2_x} \right),
   \end{align*}
   and 
   \begin{multline*}
       \left\Vert \nabla  n (t) \right\Vert_{L^2_x}+\left\Vert \nabla \mathbb{P}^\perp j (t) \right\Vert_{L^2_x}+  \| \nabla \mathbb{P}^\perp u(t) \|_{L^2_x} \\\lesssim e^{-\Lambda t} \left( \left\Vert \nabla  n(0) \right\Vert_{L^2_x}+\| \nabla j(0) \|_{L^2_x}+ \| \nabla u(0) \|_{L^2_x} \right),
\end{multline*}
therefore 
\begin{multline*}
    \left\Vert \nabla  n (t) \right\Vert_{L^2_x}+\left\Vert \nabla  j (t) \right\Vert_{L^2_x}+  \| \nabla  u(t) \|_{L^2_x} \\
    \lesssim e^{-\Lambda t} \left( \left\Vert \nabla  n(0) \right\Vert_{L^2_x}+\| \nabla j(0) \|_{L^2_x}+ \| \nabla u(0) \|_{L^2_x} \right).
\end{multline*}
When the sources $S_1^\alpha$ and $S_2$ are not zero, we can proceed 
using the Duhamel formula as in the proof of Proposition \ref{Prop-lineardecayVNS}.

\end{proof}

\begin{Rem}\label{rem-linsystem-comp-massnot1}
Remark \eqref{rem-linsystem-incom-massnot1}, i.e. when we rather write $\rho^\alpha=\mathrm{C}(1+r^\alpha)$ for $\mathrm{C}>0$, is still valid in the compressible case. The argument above based on a well-adapted energy function can be adapted, by considering this time the unknown
$$(\mathcal{N}, \mathcal{U}, \mathcal{J})=(\sqrt{c}\widehat{ n}_k,\widehat{\mathfrak{U}}_k,\sqrt{C}\widehat{\mathfrak{J}}_k), \ \ c=P'(1)>0.$$
\end{Rem}

\begin{Rem}\label{rem:Lamé-coeff}
The former proof is also adapted to more general dissipative operators with  constant coefficients in the compressible Navier--Stokes equations. More precisely, the general form is
\begin{align*}
    \mathcal{A}u=\mathrm{div}\big(2\nu D(u)+\lambda \mathrm{div}(u) Id\big)
    =\mu \Delta u+ (\lambda+\mu) \nabla \mathrm{div}u,
\end{align*}
with $\mu>0$ and $\zeta\vcentcolon=\lambda+2\mu>0$. The case we considered before corresponds to $\lambda=\mu=1$, hence $\zeta=3$. In the general case $\zeta>0$, and using the notation from the proof of Proposition  \ref{Prop-lineardecayVNS-comp}, we get
\begin{equation}
\left\{
\begin{aligned}
\partial_t  n +\vert D \vert \mathfrak{U}&=0, \\
\partial_t \mathfrak{U}  -\vert D \vert  n-\zeta \Delta \mathfrak{U} &=\mathfrak{J}-\mathfrak{U}, \\
\partial_t \mathfrak{J} &  = \mathfrak{U}-\mathfrak{J}.
\end{aligned}
\right.
\end{equation}
which leads to the study, for all $k \in \Z^3$ with $k \neq 0$, of the semigroup generated by
$$\mathbb{C}_{\zeta,k}=\begin{bmatrix}
0 & -\vert k \vert &  0\\
P'(1)\vert k \vert  & -1-\zeta\vert k \vert^2  & 1 \\
0 & 1 & -1 
\end{bmatrix}.$$
The same argument as in the proof of Proposition  \ref{Prop-lineardecayVNS-comp} can then be performed, defining  a new  coefficient $\gamma_k=\sqrt{P'(1)}\vert k \vert/(1+\zeta \vert k \vert^2)$ (after symmetrization).
\end{Rem}

\sectionhead
  [Proof of Theorem \ref{thm-globalconvergence} for the incompressible VNS system]
  {Proof of Theorem \ref{thm-globalconvergence} for the incompressible Vlasov--Navier--Stokes system}\label{Subsection-preuveStab-VNS}
We now turn to the proof of Theorem \ref{thm-globalconvergence} and Theorem \ref{thm-globalconvergence-comp}. In this section, we focus on the proof of Theorem \ref{thm-globalconvergence}, that is concerned with  the multiphasic incompressible Vlasov--Navier--Stokes system \eqref{eq:VNS-CHAP4}. The compressible case, that is system \eqref{eq:VNScompressibleCHAP4}, will be treated in Section \ref{Section:VNS-comp} following similar methods and ideas. Moreover, we mainly study the case where $\Gamma=0$. The required adaptation when  the additional term $\Gamma$ is taken into account is deferred in Section \ref{Section-nonlindrag}.

\medskip

In the proof, we will constantly use the fact that 
\begin{align*}
C_{0,k}\vcentcolon=  \| \rr_0\|_{L^{\upinfty}_\alpha H^{k-1}_x}+\Vert \bm{w}_0 \Vert_{\mathcal{H}^{k,p}}+ \| u_0\|_{H^{k}_x} \leq \eps_0,
\end{align*}
can be made small enough. The parameter $\eps_0$ therefore controls the size the  initial perturbation involved in the initial hypothesis \eqref{1st-smallness-cond} and will be appropriately fixed later on. 

\subsection{Refined energy estimates}
First, we refine the blow-up criterion already obtained in Section \ref{SubsecVNSmultiphaseLWP} by focusing on the equation for the perturbation ($r^\alpha, w^\alpha, u)$. Let us recall that it satisfies
\begin{equation}
\label{eq:VNSpertub-lwp}
\left\{  
\begin{aligned}
\partial_t r^\alpha + \mathcal{V} \cdot \nabla r^\alpha + \Div (w^\alpha)  + \Div (r^\alpha w^\alpha)&= 0,\\
\partial_t w^\alpha + ((\mathcal{V}+w^\alpha )\cdot \nabla) w^\alpha  &=  u- w^\alpha,
\\
\partial_t u + ((\mathcal{V}+u) \cdot \nabla) u - \Delta u + \nabla p &=  \int (w^\alpha-u) \D \mu(\alpha) \\
&\quad +\int r^\alpha (w^\alpha-u) \D \mu(\alpha),\\
\Div(u) &= 0.
\end{aligned}
\right.
\end{equation}
We have the following result, providing a refined exponential estimate of the Sobolev norm of the perturbation. 

\begin{Prop}\label{prop-refinedLWP-VNSincomp}
Let $k \in \N$ with $k>1+3/2$ and $ p \in [2, \infty]$. There exists $C_k>0$ such that for all $\sigma>1$ and $T>0$, all smooth solutions to \eqref{eq:VNSpertub-lwp} on $[0,T]$ satisfy the following estimate for all $t \in [0,T]$:
\begin{align*}
    &\| \ww(t) \|_{\mathcal{H}^{k,p}}+ \| u(t)\|_{H^k_x} \\
    &\leq \eps_0 \exp\Bigg(\sigma^{-1}t+ C_k\int_0^t \Big(\| (\nabla \ww, \nabla^2 \ww)(s) \|_{L^{\upinfty}_\alpha L^{\upinfty}_x}+ \|\nabla u(s) \|_{L^{\upinfty}_x}   \\
    &   \qquad \qquad \qquad   +\sigma \| (\ww -u)(s)  \|_{L^2_\alpha L^{\upinfty}_x}^2 +\| \rr(s) \|_{L^{\upinfty}_\alpha L^{\upinfty}_x}+ \| \rr(s) \|_{L^{\upinfty}_\alpha L^{\upinfty}_x}^2 \Big) \, \mathrm{d}s \Bigg).
\end{align*}
\end{Prop}
\begin{proof}
In what follows, the constants $C_k$ are universal positive constants that may change from line to line. For $\sigma>1$, let us introduce the functional
\begin{align*}
    y_{k,\sigma}(t)&=\int \| w^\alpha(t)\|_{H^k_x}^2 \D \mu(\alpha)+ \| u(t)\|_{H^k_x}^2+ \int \| \mathrm{div}(w^\alpha(t))\|_{H^{k-1}_x}^2 \D \mu(\alpha)\\
    & \quad+\sigma^{-1}  \int \| r^\alpha(t)\|_{H^{k-1}_x}^2 \D \mu(\alpha).
\end{align*}
We start by estimating the three first terms. Performing an energy estimate in $H^k_x$ for $w^\alpha$ that we integrate in $\alpha$, which is then summed to an $H^k_x$ energy estimate on $u$, we obtain
\begin{equation}
    \label{eq:trucHk}
\begin{aligned}
    &\frac{1}{2}\frac{\mathrm{d}}{\mathrm{d}t} \left(\int \| w^\alpha\|_{H^k_x}^2 \D \mu(\alpha)+\Vert u \Vert^2_{H^k_x} \right)
    + \Vert \nabla u \Vert_{H^k_x}^2 + \int \Vert u-w^\alpha \Vert^2_{H^k_x} \, \mathrm{d}\mu(\alpha) \\
    &\leq C_k \left(\int \Vert \nabla w^\alpha \Vert_{L^{\upinfty}_x} \| w^\alpha\|_{H^k_x}^2 \D \mu(\alpha) + \Vert \nabla u \Vert_{L^{\upinfty}_x} \| u\|_{H^k_x}^2 \right) \\
    &  \quad + \Vert u \Vert_{L^2_x}\int \Vert r^\alpha (w^\alpha-u) \Vert_{L^2_x} \, \mathrm{d}\mu(\alpha)+ \Vert \nabla u \Vert_{H^k_x}\int \Vert r^\alpha (w^\alpha-u) \Vert_{H^{k-1}_x} \, \mathrm{d}\mu(\alpha),
\end{aligned}
\end{equation}
for some constant $C_k>0$. Here, we have used the fact that the (linear) contribution of the two terms $u-w^\alpha$ and  $\int (w^\alpha-u) \D \mu(\alpha)$ combined in such a way that $w^\alpha \cdot (u-w^\alpha)+u\cdot(w^\alpha-u)=-\vert u-w^\alpha \vert^2$. The same applies if one replaces $u$ by $\partial_x^\ell u$ and $w^\alpha$ by $\partial_x^\ell w^\alpha$ with $\ell \leq k$.
The last term either stems from the $L^2_x$ estimate on $u$ or on its derivatives (in that case, an additional integration by parts is performed). Now we observe that $\mathrm{div}(w^\alpha)$ satisfies the equation
\begin{align*}
    \partial_t \mathrm{div}(w^\alpha)+\mathrm{div}(w^\alpha)=-\mathrm{div}\left( (\mathcal{V}+w^\alpha )\cdot \nabla w^\alpha \right),
\end{align*}
since $u$ is divergence-free. Since 
\begin{align*}
&\left\langle \mathrm{div}\left( (\mathcal{V}+w^\alpha )\cdot \nabla w^\alpha \right), \mathrm{div}(w^\alpha) \right\rangle_{H^{k-1}_x} \\
&=\left\langle (\nabla w^\alpha)^{T} : \nabla w^\alpha , \mathrm{div}(w^\alpha) \right\rangle_{H^{k-1}_x}
+\left\langle (\mathcal{V}+w^\alpha )\cdot \nabla \mathrm{div}(w^\alpha), \mathrm{div}(w^\alpha) \right\rangle_{H^{k-1}_x} \\
&= \sum_{\ell \leq k-1} \left\langle \partial^\ell ((\nabla w^\alpha)^{T} : \nabla w^\alpha) , \partial^\ell \mathrm{div}(w^\alpha) \right\rangle_{L^2_x} \\
& \quad +\sum_{\ell \leq k-1}\left\langle w^\alpha \cdot \partial^\ell \nabla \mathrm{div}(w^\alpha), \partial^\ell\mathrm{div}(w^\alpha) \right\rangle_{L^2_x} \\
& \quad +\sum_{\ell \leq k-1}\left\langle [\partial^\ell, w^\alpha \cdot \nabla]\mathrm{div}(w^\alpha), \partial^\ell\mathrm{div}(w^\alpha) \right\rangle_{L^2_x},
\end{align*}
we can employ tame estimates from Proposition \ref{prop-Sobolev} to get
\begin{align*}
    \left\langle \mathrm{div}\left( ((\mathcal{V}+w^\alpha )\cdot \nabla) w^\alpha \right), \mathrm{div}(w^\alpha) \right\rangle_{H^{k-1}_x} &\leq C_k\Vert \nabla w^\alpha \Vert_{L^{\upinfty}_x} \Vert \nabla w^\alpha \Vert_{H^{k-1}_x} \Vert \mathrm{div}(w^\alpha) \Vert_{H^{k-1}_x}\\
    & \quad   +C_k\Vert \nabla w^\alpha \Vert_{L^{\upinfty}_x} \Vert \mathrm{div}(w^\alpha) \Vert_{H^{k-1}_x}^2  \\
    & \quad  + C_k \Vert \nabla^2 w^\alpha \Vert_{L^{\upinfty}_x} \Vert w^\alpha \Vert_{H^{k-1}_x} \Vert \mathrm{div}(w^\alpha) \Vert_{H^{k-1}_x},
\end{align*}
so that we obtain by Young inequality
\begin{multline*}
    \frac{1}{2}\frac{\mathrm{d}}{\mathrm{d}t} \Vert \mathrm{div}(w^\alpha) \Vert_{H^{k-1}_x}^2+ \Vert \mathrm{div}(w^\alpha)\Vert_{H^{k-1}_x}^2  \\
     \leq C_k (\Vert \nabla w^\alpha \Vert_{L^{\upinfty}_x}+\Vert \nabla^2 w^\alpha \Vert_{L^{\upinfty}_x})(\Vert w^\alpha \Vert_{H^k_x}^2+\Vert \mathrm{div}(w^\alpha) \Vert_{H^{k-1}_x}^2).
\end{multline*}
By integrating in $\alpha$ and summing to~\eqref{eq:trucHk}, we get
\begin{align*}
    &\frac{1}{2}\frac{\mathrm{d}}{\mathrm{d}t} \left(\int \| w^\alpha\|_{H^k_x}^2 \D \mu(\alpha)+ \| u\|_{H^k_x}^2+ \int \| \mathrm{div}(w^\alpha)\|_{H^{k-1}_x}^2 \D \mu(\alpha) \right) \\
    & \qquad \qquad \qquad   + \int \| \mathrm{div}(w^\alpha)\|_{H^{k-1}_x}^2 \D \mu(\alpha) + \Vert \nabla u \Vert_{H^k_x}^2 \\
    &\leq C_k \left(\Vert (\nabla w^\alpha, \nabla^2 w^\alpha) \Vert_{L^{\upinfty}_\alpha L^\infty_x } \int  (\Vert w^\alpha \Vert_{H^k_x}^2+\Vert \mathrm{div}(w^\alpha) \Vert_{H^{k-1}_x}^2) \D \mu(\alpha) \right.\\
    & \qquad \left.+ \Vert \nabla u \Vert_{L^{\upinfty}_x} \| u\|_{H^k_x}^2 \right) + \Vert \rr \Vert_{L^{\upinfty}_\alpha L^{\upinfty}_x} \Vert u \Vert_{L^2_x}\int  \Vert w^\alpha-u  \Vert_{L^2_x} \, \mathrm{d}\mu(\alpha)
      \\
      & \quad +\Vert \nabla u \Vert_{H^k_x}\int \Vert r^\alpha (w^\alpha-u) \Vert_{H^{k-1}_x} \, \mathrm{d}\mu(\alpha).
\end{align*}
We now use Young inequality on the two last terms, with a tame estimate for products (see again Proposition \ref{prop-Sobolev}) in the last one, to get:
\begin{multline*}
    \Vert \rr \Vert_{L^{\upinfty}_\alpha L^{\upinfty}_x} \Vert u \Vert_{L^2_x}\int  \Vert w^\alpha-u  \Vert_{L^2_x} \, \mathrm{d}\mu(\alpha) \\
    \leq \Vert \rr \Vert_{L^{\upinfty}_\alpha L^{\upinfty}_x} \Vert u \Vert_{H^k_x}^2+ \Vert \rr \Vert_{L^{\upinfty}_\alpha L^{\upinfty}_x} \int  \Vert w^\alpha-u  \Vert_{H^k_x}^2 \, \mathrm{d}\mu(\alpha)
    \end{multline*}
    and
    \begin{align*}
    &\Vert \nabla u \Vert_{H^k_x}\int \Vert r^\alpha (w^\alpha-u) \Vert_{H^{k-1}_x} \, \mathrm{d}\mu(\alpha) \\
    &\leq \frac{1}{2} \Vert \nabla u \Vert_{H^k_x}^2 \\
    &\qquad + C_k\left( \int \Vert r^\alpha \Vert_{L^{\upinfty}_x} \Vert w^\alpha-u \Vert_{H^{k-1}_x} \, \mathrm{d}\mu(\alpha)+ \int \Vert r^\alpha \Vert_{H^{k-1}_x} \Vert w^\alpha-u \Vert_{L^{\upinfty}_x} \, \mathrm{d}\mu(\alpha)\right)^2\\
    &\leq \frac{1}{2} \Vert \nabla u \Vert_{H^k_x}^2+ C_k \Vert \rr \Vert_{L^{\upinfty}_\alpha L^{\upinfty}_x}^2 \int  \Vert w^\alpha-u \Vert_{H^k_x}^2 \, \mathrm{d}\mu(\alpha)\\
    & \qquad   + C_k \left(\int \Vert r^\alpha \Vert_{H^{k-1}_x}^2 \, \mathrm{d}\mu(\alpha)\right) \left( \int\Vert w^\alpha-u \Vert_{L^{\upinfty}_x}^2 \, \mathrm{d}\mu(\alpha)\right),
\end{align*}
thanks to the Cauchy--Schwarz inequality in $\alpha$. We can now absorb the gradient term $\Vert \nabla u \Vert_{H^k_x}$ in the left-hand side, which leads to 
\begin{equation}
    \label{eq:trucHk2}
    \begin{aligned}
    &\frac{\mathrm{d}}{\mathrm{d}t} \left(\int \| w^\alpha\|_{H^k_x}^2 \D \mu(\alpha)+ \| u\|_{H^k_x}^2+ \int \| \mathrm{div}(w^\alpha)\|_{H^{k-1}_x}^2 \D \mu(\alpha) \right) \\
    & \qquad  \qquad \qquad  \qquad \qquad \qquad \qquad + \int \| \mathrm{div}(w^\alpha)\|_{H^{k-1}_x}^2 \D \mu(\alpha) + \frac{1}{2}\Vert \nabla u \Vert_{H^k_x}^2 \\
    &\qquad \leq  C_k \left(\Vert (\nabla w^\alpha, \nabla^2 w^\alpha) \Vert_{L^{\upinfty}_\alpha L^\infty_x } \int  (\Vert w^\alpha \Vert_{H^k_x}^2+\Vert \mathrm{div}(w^\alpha) \Vert_{H^{k-1}_x}^2) \D \mu(\alpha) \right. \\
    &\quad \quad  + (\Vert \rr \Vert_{L^{\upinfty}_\alpha L^{\upinfty}_x}+\Vert \nabla u \Vert_{L^{\upinfty}_x}) \| u\|_{H^k_x}^2  \\
    & \quad \quad   + \left. \Vert \rr \Vert_{L^{\upinfty}_\alpha L^{\upinfty}_x}^2 \int  \Vert w^\alpha-u \Vert_{H^{k}_x}^2 \D \mu(\alpha)+  \Vert \ww-u \Vert_{L^2_\alpha L^{\upinfty}_x}^2 \int \Vert r^\alpha \Vert_{H^{k-1}_x}^2 \D \mu(\alpha)   \right),
\end{aligned}
\end{equation}
 Let us now proceed with a tame estimate of $r^\alpha$ in $H^{k-1}_x$ (recall Proposition~\ref{prop-Sobolev}): we get for all  $\sigma>1$
\begin{align}\label{eq:estim-Sobolev-perturb-ralpha}
\begin{split}
    \frac{1}{2} \frac{\mathrm{d}}{\mathrm{d}t} \Vert r^\alpha \Vert^2_{H^{k-1}_x} 
     &\leq C_k \left( \| \nabla w^\alpha \|_{L^{\upinfty}_x} \| r^\alpha\|_{H^{k-1}_x}^2  + 
\| \nabla w^\alpha\|_{H^{k-1}_x} \| r^\alpha\|_{L^{\upinfty}_x}\|r^\alpha\|_{H^{k-1}_x} \right) \\
& \quad +   \Vert r^\alpha \Vert_{H^{k-1}_x}\Vert \mathrm{div}(w^\alpha) \Vert_{H^{k-1}_x} \\
&\leq C_k \left( \| \nabla w^\alpha \|_{L^{\upinfty}_x} \| r^\alpha\|_{H^{k-1}_x}^2  + 
\sigma\| r^\alpha\|_{L^{\upinfty}_x} \| w^\alpha\|_{H^k_x}^2 +\sigma^{-1}\| r^\alpha\|_{L^{\upinfty}_x}\|r^\alpha\|_{H^{k-1}_x}^2 \right) \\
& \quad + \sigma^{-1}  \Vert r^\alpha \Vert_{H^{k-1}_x}^2+\frac{\sigma}{2}\Vert \mathrm{div}(w^\alpha)\Vert_{H^{k-1}_x}^2.
\end{split}
 \end{align}
Dividing by $\sigma>1$ and integrating in $\alpha$, we get
\begin{align*}
     &\frac{\mathrm{d}}{\mathrm{d}t} \sigma^{-1}\int \Vert r^\alpha \Vert^2_{H^{k-1}_x} \D \mu(\alpha)\\
&\leq C_k \Bigg( (\| \nabla \ww \|_{L^{\upinfty}_\alpha L^{\upinfty}_x}+\| \rr\|_{L^{\upinfty}_\alpha L^{\upinfty}_x}) \sigma^{-1}\int \| r^\alpha\|_{H^{k-1}_x}^2 \D \mu(\alpha) \\
&\qquad + 
\| \rr\|_{L^{\upinfty}_\alpha L^{\upinfty}_x} \int \| w^\alpha\|_{H^k_x}^2 \D \mu(\alpha) \Bigg)  \\
& \quad + \sigma^{-2}  \int \Vert r^\alpha \Vert_{H^{k-1}_x}^2 \D \mu(\alpha) +\frac{1}{2}\int \Vert \mathrm{div}(w^\alpha)\Vert_{H^{k-1}_x}^2 \D \mu(\alpha).
\end{align*}
Summing to the previous differential inequality~\eqref{eq:trucHk2}, and finally absorbing the last term above in the left-hand side, we end up with the differential inequality
\begin{multline*}
    \frac{\mathrm{d}}{\mathrm{d}t}y_{k,\sigma} \leq \bigg[\sigma^{-1}+C_k \big(\| (\nabla \ww, \nabla^2 \ww) \|_{L^{\upinfty}_\alpha L^{\upinfty}_x}+ \|\nabla u \|_{L^{\upinfty}_x}\\
    + \sigma \| \ww -u  \|_{L^2_\alpha L^{\upinfty}_x}^2 + \| \rr \|_{L^{\upinfty}_\alpha L^{\upinfty}_x}+ \| \rr \|_{L^{\upinfty}_\alpha L^{\upinfty}_x}^2 \big) \bigg]y_{k,\sigma}.
\end{multline*}
By Gronwall inequality, and since $\| u(t)\|_{H^k_x}^2 \leq y_{k, \sigma}(t)$, we have in particular
\begin{multline}
    \label{eq:boundu}
    \| u(t)\|_{H^k_x}^2 \\ 
    \leq y_{k,\sigma}(0)\exp\Bigg(\sigma^{-1}t+ C_k\int_0^t \Big(\| (\nabla \ww, \nabla^2 \ww)(s) \|_{L^{\upinfty}_\alpha L^{\upinfty}_x}+ \|\nabla u(s) \|_{L^{\upinfty}_x} \\
    + \sigma\| (\ww -u)(s)  \|_{L^2_\alpha L^{\upinfty}_x}^2 + \| \rr(s) \|_{L^{\upinfty}_\alpha L^{\upinfty}_x} +\| \rr(s) \|_{L^{\upinfty}_\alpha L^{\upinfty}_x}^2\Big)\, \mathrm{d}s\Bigg).
\end{multline}
 Note that since $p \geq 2$, we actually have 
\begin{align*}
    y_{k,\sigma}(0) \lesssim \eps_0.
\end{align*}
We now come back to the equation satisfied by $w^\alpha$, treating $u$ as a fixed source term. Namely, we write
\begin{align*}
\partial_t w^\alpha + w^\alpha+ ((\mathcal{V}+w^\alpha )\cdot \nabla) w^\alpha  =  u.
\end{align*}
Relying on (the proof) of Lemma \ref{LM:existence-vitesseSOBOLEV}, we obtain
\begin{align*}
  \frac{\mathrm{d}}{\mathrm{d}t}\| \ww \|_{\mathcal{H}^{k,p}} \leq (C_k \Vert \nabla w^\alpha \Vert_{L^{\upinfty}_x}-1)\| \ww \|_{\mathcal{H}^{k,p}}+ \Vert u \Vert_{H^k_x},
\end{align*}
for some $C_k>0$, from which we infer
\begin{align*}
\| \ww(t) \|_{\mathcal{H}^{k,p}} 
&\leq (\| \ww_0 \|_{\mathcal{H}^{k,p}}+t \Vert u \Vert_{L^\infty_t H^k_x})e^{-t}e^{C_k\int_0^t \Vert \nabla w^\alpha(s) \Vert_{L^{\upinfty}_x} \, \mathrm{d}s} \\
&\leq (\| \ww_0 \|_{\mathcal{H}^{k,p}}+\Vert u \Vert_{L^\infty_t H^k_x})e^{C_k\int_0^t \Vert \nabla w^\alpha(s) \Vert_{L^{\upinfty}_x} \, \mathrm{d}s}.
\end{align*}
By the previous bound~\eqref{eq:boundu} on  $\| u(t)\|_{H^k_x}$, we finally reach the desired conclusion.
\end{proof}

We also need the following pointwise estimates on the perturbation of the density $(r^\alpha)_{\alpha \in I}$, slightly refining what was stated in Lemma \ref{LM:pointwiseEstimRHO}.
\begin{Lem}\label{LM:pointwise-densité}
For $\mu$-almost all $\alpha$, for all $T>0$, all smooth solutions on $[0,T]$ to 
$$\partial_t r^\alpha + \mathcal{V} \cdot \nabla r^\alpha + \Div (w^\alpha)  + \Div (r^\alpha w^\alpha)= 0$$
satisfy, for all $t  \in [0,T]$
    \begin{align*}
\Vert r^\alpha(t) \Vert_{L^{\upinfty}_x} &\lesssim 1+\left( 1+\Vert r^\alpha_0 \Vert_{L^{\upinfty}_x} \right)\exp\left(\int_0^t \Vert \nabla w^{\alpha}(\tau) \Vert_{L^{\upinfty}_x} \mathrm{d} \tau\right), \\
\Vert \nabla r^\alpha(t) \Vert_{L^{\upinfty}_x} &\lesssim \exp\left(C\int_0^T \Vert \nabla w^{\alpha}(\tau) \Vert_{L^{\upinfty}_x} \mathrm{d} \tau\right) \\
&\qquad \qquad\times \left[\Vert \nabla r^\alpha_0 \Vert_{L^{\upinfty}_x} +(1+\Vert  r^\alpha_0 \Vert_{L^{\upinfty}_x})\int_0^t \Vert \mathrm{D}^2 w^\alpha(s) \Vert_{L^{\upinfty}_x} \mathrm{d}s \right].
\end{align*}
\end{Lem}
\begin{proof}
Consider the characteristics associated to the vector field $w^{\alpha}+\mathcal{V}$, i.e
\begin{align*}
\dfrac{\mathrm{d}}{\mathrm{d}s} \mathrm{X}^{s;t}(x)=\mathcal{V}+w^{\alpha}(s,\mathrm{X}^{s;t}(x)), \ \ \mathrm{X}^{t;t}(x)=x.
\end{align*}
As $\rho^\alpha= 1+ r^\alpha$ satisfies the continuity equation
$$
\partial_t \rho^\alpha + \Div((\mathcal{V} + w^\alpha) \rho^\alpha)=0,
$$
we obtain by the method of characteristics that
\begin{align}
    r^{\alpha}(t,x)=-1+\left( 1+r^{\alpha}_0(\mathrm{X}^{0;t}(x))\right)\exp\left(-\int_0^t \mathrm{div}(w^{\alpha})(\tau, \mathrm{X}^{\tau;t}(x)) \mathrm{d} \tau\right).
\end{align}
The first estimate thus readily follows. To estimate $\nabla r^\alpha$, we take the gradient of the latter formula and get
\begin{align*}
    \nabla r^{\alpha}(t,x)=\exp&\left(-\int_0^t \mathrm{div}(w^{\alpha})(\tau, \mathrm{X}^{\tau;t}(x)) \mathrm{d} \tau\right)  \Big[ \mathrm{D}\mathrm{X}^{0;t}(x) \nabla r^{\alpha}_0(\mathrm{X}^{0;t}(x))\\
    &-(1+r^{\alpha}_0(\mathrm{X}^{0;t}(x)))\int_0^t \mathrm{D}\mathrm{X}^{\tau;t}(x) \nabla \mathrm{div}(w^{\alpha})(\tau, \mathrm{X}^{\tau;t}(x)) \mathrm{d} \tau \Big].
\end{align*}
Since
\begin{align*}
\mathrm{D}\mathrm{X}^{s;t}(x)&=\mathrm{I}+\int_t^{s} \mathrm{D}\mathrm{X}^{\sigma;t}(x) \nabla w^\alpha(\sigma, \mathrm{X}^{\sigma;t}(x)) \, \mathrm{d}\sigma,
\end{align*}
we infer from Gronwall lemma that
\begin{align*}
    \vert \mathrm{D}\mathrm{X}^{s;t}(x) \vert \leq \exp\left(\int_0^T \Vert \nabla w^{\alpha}(\tau) \Vert_{L^{\upinfty}_x} \mathrm{d} \tau\right),
\end{align*}
for all $s,t \in [0,T]$. Coming back to $\nabla r^{\alpha}$, we can now estimate
\begin{multline*}
\Vert \nabla r^{\alpha}(t) \Vert_{L^{\upinfty}_x} \lesssim \exp\left(C\int_0^T \Vert \nabla w^{\alpha}(\tau) \Vert_{L^{\upinfty}_x} \mathrm{d} \tau\right) \\
\times\left[\Vert \nabla r^\alpha_0 \Vert_{L^{\upinfty}_x} +(1+\Vert  r^\alpha_0 \Vert_{L^{\upinfty}_x})\int_0^t \Vert \mathrm{D}^2 w^\alpha(s) \Vert_{L^{\upinfty}_x} \mathrm{d}s \right],
\end{multline*}
and this concludes the proof.
\end{proof}

\subsection{Bootstrap}\label{sectionVNS-Stratbootstrap}

\paragraph{Set up of the bootstrap argument.}
Let $k \in \N$ such that $k>2+3/2$ and $p \in [2, \infty]$. For a given initial data $(\rr_0, \ww_0, u_0)$, we introduce the maximal time of existence $T^{\mathrm{max}}_{0}>0$ of the unique strong solution $(\rr,\ww,u)$ to the equation \eqref{eq:VNSpertub}, which is given by Theorem \ref{thm:appliVNS-incomp}, and such that for any $T \in (0, T^{\mathrm{max}}_{0})$, we have
\begin{align*}
    &\rr \in L^\infty([0,T]; L^\infty_\alpha (H^{k-1}_x) ), \ \ \ww \in L^\infty([0,T]; \mathcal{H}^{k,p}), \\
    &u \in L^\infty([0,T]; H^{k}_x) \cap L^2([0,T]; H^{k+1}_x).
\end{align*}
Note that Theorem \ref{thm:appliVNS-incomp} applies straightforwardly to the perturbation $(\rr,\ww,u)$ or to the original solution  $(\rhorho,\ww,U)$). 

\medskip

Let $\lambda \in (0, \Lambda)$ where $\Lambda \in (0,1)$ is the linear decay rate stemming from Proposition \ref{Prop-lineardecayVNS}. Let $\lambda' \in (0, \theta \lambda)$ for some $\theta \in (0,1)$ that will be chosen later on. For $\delta \in (0,1)$ small enough also to be determined later, we introduce
\begin{equation}\label{assump-bootstrapVNS}
\begin{aligned}
T^\star_{\delta} \vcentcolon= \sup\Bigg\{ T \in \left(0, T^{\mathrm{max}}_{0}\right),    \, \forall t \in [0,T], \qquad \qquad \qquad \qquad\qquad \qquad & 
  \\
  \|(\nabla \ww,\nabla^2 \ww)(t)\|_{L^{\upinfty}_\alpha  L^{\upinfty}_x} + \|(\nabla u,\nabla^2 u)(t) \|_{L^{\upinfty}_x}  &\leq \delta e^{-\lambda' t} \\
    \|\nabla \ww(t)\|_{L^{\upinfty}_\alpha  L^{2}_x} + \|\nabla u(t) \|_{L^{2}_x}  &\leq \delta e^{-\lambda t} \\
     \Vert  \langle  \ww(t)- u(t) \rangle \Vert_{L^p_\alpha} &\leq \delta e^{-\lambda t}
\Bigg\}.
\end{aligned}
\end{equation}
This implies in particular that for all $t \in [0, T^\star_{\delta}]$, we have
\begin{align*}
     \|\nabla \ww(t)-\nabla u(t)\|_{L^{\upinfty}_\alpha L^{2}_x} \leq \delta e^{-\lambda t}.
\end{align*}
By the local in time theory from Theorem \ref{thm:appliVNS-incomp} (and the remark below this statement), as soon as the size of the initial conditions $(\bm{w}_0,u_0) \in \mathcal{H}^{k,p} \times H^k_x$ is small with respect to $\delta$,  we have by continuity in time that $T^\star_{\delta} >0$. Note that the continuity property of the true solution $(\vv, U)$ is enough to handle \eqref{assump-bootstrapVNS} on the perturbation $(\ww, u)$, since they only differs by the same constant velocity.

\medskip

The goal is now to prove that 
$T^\star_{\delta}=+\infty$. To do so, we will show that  for $\delta$ small enough and then $\varepsilon_0$ (the size of the initial data -  see the notation of Theorem \ref{thm-globalconvergence}) small enough, the above bootstrap bounds on $[0,T^\star_\delta]$ (see \eqref{assump-bootstrapVNS}) can be improved in the following way: for all $t\in[0,T^\star_\delta]$, there holds
\begin{align*}
\|(\nabla \ww,\nabla^2 \ww)(t)\|_{L^{\upinfty}_\alpha  L^{\upinfty}_x} + \|(\nabla u,\nabla^2 u)(t) \|_{L^{\upinfty}_x} &\leq \frac{\delta}{2}e^{-\lambda' t}, \\
\|\nabla \ww(t)\|_{L^{\upinfty}_\alpha  L^{2}_x} + \|\nabla u(t) \|_{L^{2}_x}&\leq \frac{\delta}{2}e^{-\lambda t}, \\
\Vert  \langle  \ww(t)- u(t) \rangle \Vert_{L^p_\alpha} &\leq \frac{\delta}{2} e^{-\lambda t}.
\end{align*}
By continuity in time, if this improvement holds then it would imply that the inequalities defining $T^\star_\delta$
still hold on $[0,T^\star_\delta+\eta]$ for some $\eta>0$, contradicting the definition of $T^\star_\delta$,
unless $T^\star_\delta = T^{\max}_{0}$. Finally, if $T^{\mathrm{max}}_{0}<\infty$, the blow-up criterion (Proposition~\ref{prop-VNSblowup}) yields
$$
\int_0^{T^{\mathrm{max}}_{0}} \big(\|\nabla \ww(s)\|_{L^{\upinfty}_\alpha L^{\upinfty}_x}
+\|\nabla u(s)\|_{L^{\upinfty}_x}\big) \, \mathrm{d}s=+\infty,
$$
which is impossible in view the exponential in time decay of the integrands. Hence $T^{\mathrm{max}}_{0}=+\infty$, and therefore we will obtain $T^\star_\delta=+\infty$.

Let us assume from now on that $T^\star_{\delta}<T^{\mathrm{max}}_{0}$. In particular, we have $T^\star_{\delta}<\infty$.
We will repeatedly use the notation $C_{0, \delta, \lambda}$ (resp. $C_{\delta, \lambda})$ to denote some constants that may change from line to line and that depend continuously, in an increasing way with respect to each of the parameters, of $\eps_0, \delta, \lambda$ (resp. to the parameters $\delta, \lambda$).

\medskip

From the refined bounds from Proposition \ref{prop-refinedLWP-VNSincomp}, we can first of all obtain a control of high-order derivatives of the velocities. Indeed, using the bootstrap assumption \eqref{assump-bootstrapVNS}, we find that for any $k \in \N$ with $k>1+3/2$, there exists $C_k>0$ such that for all $\sigma>0$ and for all $t \in [0,T_\delta^\star)$
\begin{align*}
    &\| \ww(t) \|_{\mathcal{H}^{k,p}}+ \| u(t)\|_{H^k_x} \\
    &\leq \eps_0 \exp\Bigg(\sigma^{-1}t+ C_k\int_0^t \Big(\| (\nabla \ww, \nabla^2 \ww)(s) \|_{L^{\upinfty}_\alpha L^{\upinfty}_x}+ \|\nabla u(s) \|_{L^{\upinfty}_x}   \\
    &  \qquad \qquad \qquad \qquad    +\sigma \| (\ww -u)(s)  \|_{L^2_\alpha L^{\upinfty}_x}^2 +\| \rr(s) \|_{L^{\upinfty}_\alpha L^{\upinfty}_x}+ \| \rr(s) \|_{L^{\upinfty}_\alpha L^{\upinfty}_x}^2 \Big) \, \mathrm{d}s \Bigg) \\
    & \leq \eps_0 \exp\bigg(\sigma^{-1}t+ C_k\int_0^t \Big(\| (\nabla \ww, \nabla^2 \ww)(s) \|_{L^{\upinfty}_\alpha L^{\upinfty}_x}+ \|\nabla u(s) \|_{L^{\upinfty}_x} \\
    &   \qquad \qquad \qquad \qquad  \qquad \qquad  + C\sigma\| \nabla (\ww -u)(s)  \|_{L^{\upinfty}_\alpha L^{\upinfty}_x}^2  +\sigma \Vert  \langle  \ww- u \rangle(s) \Vert_{L^2_\alpha}^2 \\
    & \qquad \qquad \qquad \qquad \qquad \qquad + \| \rr(s) \|_{L^{\upinfty}_\alpha L^{\upinfty}_x}(1+\| \rr(s) \|_{L^{\upinfty}_\alpha L^{\upinfty}_x}) \Big) \, \mathrm{d}s \bigg) \\
    & \leq \eps_0 \exp \bigg(\sigma^{-1}t+ \delta(1+\sigma) C_{\lambda'}C_k+ C_k \int_0^t \| \rr(s) \|_{L^{\upinfty}_\alpha L^{\upinfty}_x}(1+\| \rr(s) \|_{L^{\upinfty}_\alpha L^{\upinfty}_x}) \, \mathrm{d}s\bigg),
\end{align*}
for some constant $C_{\lambda'}>0$. Furthermore, we have obtained in the proof of Lemma \ref{LM:pointwise-densité} that 
\begin{multline*}
    r^{\alpha}(t,x)=r^{\alpha}_0(\mathrm{X}^{0;t}(x))\exp\left(-\int_0^t \mathrm{div}(w^{\alpha})(\tau, \mathrm{X}^{\tau;t}(x)) \mathrm{d} \tau\right) \\ +\exp\left(-\int_0^t \mathrm{div}(w^{\alpha})(\tau, \mathrm{X}^{\tau;t}(x)) \mathrm{d} \tau\right)-1.
\end{multline*}
Setting $a(t)\vcentcolon=-\int_0^t \mathrm{div}(w^{\alpha})(\tau, \mathrm{X}^{\tau;t}(x)) \mathrm{d} \tau$, we have by the mean-value theorem that
\begin{align*}
    \vert e^{a(t)}-1\vert =\vert e^{a(t)}-e^{a(0)}\vert \leq \sup_{[0,t]}\vert a' e^a\vert t \leq \Vert \nabla w^\alpha(t) \Vert_{L^{\upinfty}_x} e^{\int_0^t \Vert \nabla w^\alpha(s) \Vert_{L^{\upinfty}_x} \, \mathrm{d}s}t ,
\end{align*}
hence we infer from the bootstrap assumption that
\begin{align}\label{bound-perturb-density-bootstrap}
    \Vert \rr(t) \Vert_{L^{\upinfty}_\alpha L^{\upinfty}_x} \leq \Vert \rr_0 \Vert_{L^{\upinfty}_\alpha L^{\upinfty}_x}e^{\delta/\lambda'}+\delta e^{-\lambda' t}e^{\delta/\lambda'}t \leq e^{\delta/\lambda'} C_{\lambda'}(\eps_0+\delta),
\end{align}
for some constant $C_{\lambda'}>0$. Plugging this estimate in the former exponential bounds with the new estimate at hands, we get
\begin{align*}
    &C_k \int_0^t \| \rr(s) \|_{L^{\upinfty}_\alpha L^{\upinfty}_x}(1+\| \rr(s) \|_{L^{\upinfty}_\alpha L^{\upinfty}_x}) \, \mathrm{d}s \\
    &\qquad \leq C_k e^{\delta/\lambda'} C_{\lambda'}(\eps_0+\delta) (1+e^{\delta/\lambda'} C_{\lambda'}(\eps_0+\delta))t \\
    &\qquad  \leq e^{\delta/\lambda'} C_{\lambda'}(\eps_0+\delta)t,
\end{align*}
for a modified constant $C_{\lambda'}$, if one imposes  without loss of generality that $\eps_0+\delta \leq 10$. We end up with the following estimate: for $t \in [0,T_\delta^\star)$ and $\sigma>0$, we have 
\begin{multline}\label{bound-opti-growthexpoVNS}
    \| \ww(t) \|_{\mathcal{H}^{k,p}}+ \| u(t)\|_{H^k_x} \\
     \leq \eps_0 \exp \bigg(\left(\sigma^{-1}+ C_k e^{\delta/\lambda'} C_{\lambda'}(\eps_0+\delta) \right)t+ \delta(1+\sigma) C_{\lambda'}C_k\bigg).
\end{multline}

\subsection{Decay estimates}\label{sectionVNS-decay-estimate}

We now improve the estimates appearing in the bootstrap assumption \eqref{assump-bootstrapVNS}, starting with the decay of the gradients in $L^2_x$. We begin with an improved decay of  $\| \nabla u(t) \|_{L^2_x}$.
\begin{Lem}\label{lem:ImproveLinL2-nabla-u}
There exists $C_{0,\lambda',\delta}>0$ such that all $t \in [0, T^\star_\delta]$, we have
\begin{align*}
\| \nabla u(t) \|_{L^2_x}    \lesssim \left((1+C_{0, \lambda',\delta})\eps_0+(1+C_{0, \lambda',\delta})\delta^2 \right) e^{-\lambda t}.
\end{align*}
\end{Lem}
\begin{proof}
   From the system~\eqref{eq:VNSpertub} satisfied by the perturbation, we aim at applying Proposition~\ref{Prop-lineardecayVNS} with $\widetilde{\mathcal{W}} = \mathcal{W}$ (recall \eqref{def:WM} and $\mathcal{Z}=\mathcal{W}-\mathcal{V}$) and
\begin{align*}
S^\alpha_1&\vcentcolon=- ((w^\alpha-\mathcal{Z}) \cdot \nabla) w^\alpha, \\
S_2 &\vcentcolon= S_{2,a} + S_{2,b},   \\
S_{2,a}&\vcentcolon=- ((u-\mathcal{Z}) \cdot \nabla) u, \ \ S_{2,b}\vcentcolon= \int r^\alpha (w^\alpha-u) \D \mu(\alpha). \ \ 
\end{align*}
We thus obtain, for all $t \in [0, T^\star_\delta]$
\begin{align}
\label{decay-nablau-proof} \|  \nabla u(t) \|_{L^2_x}  &\lesssim \eps_0 e^{-\Lambda t} + \int_0^t e^{-\Lambda (t-s)} \Big( \|  \nabla S^\alpha_1 (s)  \|_{L^p_\alpha L^2_x}+ \|  \nabla S_2(s)\|_{ L^2_x} \Big) \, \mathrm{d}s.
\end{align}
 Here we used the fact that $\mu$ is of finite mass to get a $L^p_\alpha$ norm in the right-hand side.
 We now estimate each term from the sources under the integrals in time, with the following guideline. For terms like  $S_1^\alpha$ and $S_{2,a}$ involving products of velocities with at least one derivative, we perform standard $L^2_x-L^\infty_x$ estimates combined with the bootstrap assumption. For the product term $S_{2,b}$, where no decay in time on $r^\alpha$ or its derivatives is expected, we rely on the initial smallness assumption.

\medskip

$\bullet$ \underline{Term $\|  \nabla S^\alpha_1\|_{L^p_\alpha L^2_x}$}: we have
\begin{align*}
\| \nabla S^\alpha_1 \|_{L^2_x}&=\| \nabla \big(((w^\alpha-\mathcal{Z}) \cdot \nabla) w^\alpha \big) \|_{L^2_x} \\
& \lesssim \| \nabla w^\alpha\|_{L^{2}_x}\| \nabla w^\alpha\|_{L^{\upinfty}_x}+\| w^\alpha-\mathcal{Z} \|_{L^2_x}\| \nabla^2 w^\alpha\|_{L^{\upinfty}_x}\\
& \lesssim \| \nabla w^\alpha\|_{L^{2}_x}\| \nabla w^\alpha\|_{L^{\upinfty}_x} \\
&\quad +\Big(\| w^\alpha-\langle w^\alpha \rangle \|_{L^2_x}+\| \langle w^\alpha \rangle-\langle u \rangle \|_{L^2_x}+\| \langle u \rangle-\mathcal{Z} \|_{L^2_x}\Big)\|\nabla^2 w^\alpha\|_{L^{\upinfty}_x} \\
& \lesssim \| \nabla w^\alpha\|_{L^{2}_x}\| \nabla w^\alpha\|_{L^{\upinfty}_x}+\Big(\| \nabla w^\alpha \|_{L^2_x}
    +| \langle w^\alpha- u \rangle|
    +| \langle u \rangle-\mathcal{Z} |\Big)\|\nabla^2 w^\alpha\|_{L^{\upinfty}_x},
\end{align*}
where we have used Poincaré inequality. Relying on Lemma \ref{smallLemma-moyenne}, we also know that
\begin{align*}
    | \langle u \rangle-\mathcal{Z} | &\lesssim \Vert \rhorho(t) \Vert_{L^{\upinfty}_\alpha L^2_x} \left( \Vert \nabla u(t) \Vert_{L^2_x} + \int \Vert w^\alpha(t)-u(t) \Vert_{L^2_x}\, \mathrm{d}\mu(\alpha)\right).
\end{align*}
Using the exponential decay of $\Vert \nabla \ww \Vert_{L^{\upinfty}_\alpha L^{\upinfty}_x}$ from the bootstrap assumption to bound the density thanks to Lemma \ref{LM:pointwise-densité}, we get
\begin{align*}
    | \langle u \rangle-\mathcal{Z} | \lesssim C_{0, \lambda',\delta} \left( \Vert \nabla u(t) \Vert_{L^2_x} + \int \Vert w^\alpha(t)-u(t) \Vert_{L^2_x}\, \mathrm{d}\mu(\alpha)\right).
\end{align*}
By Poincaré inequality, and since $p \in [1, \infty]$, we obtain
\begin{align}\label{amelio-conservationlaw}
| \langle u \rangle-\mathcal{Z} |\lesssim C_{0, \lambda',\delta} \left( \Vert \nabla u(t) \Vert_{L^2_x} +\Vert \nabla (\ww-u) \Vert_{L^p_\alpha L^2_x} + \Vert \langle  \ww-u \rangle  \Vert_{L^p_\alpha}\right).
\end{align}
All in all, we can thus infer
\begin{align*}
    &\| \nabla S^\alpha_1 \|_{L^2_x}  \\
    &\lesssim
    \| \nabla w^\alpha\|_{L^{2}_x}\| \nabla w^\alpha\|_{L^{\upinfty}_x}+\| \nabla w^\alpha \|_{L^2_x}\|\nabla^2 w^\alpha\|_{L^{\upinfty}_x}
    \\ 
    & \quad +\Big( \vert \langle w^\alpha- u \rangle |+ C_{0, \lambda',\delta} \Vert \nabla u \Vert_{L^2_x}+ C_{0, \lambda',\delta} (\Vert \nabla (\ww-u) \Vert_{L^p_\alpha L^2_x} \\
    & \quad \qquad \qquad \qquad \qquad \qquad \qquad \qquad+ \Vert \langle  \ww-u \rangle  \Vert_{L^p_\alpha})
    \Big)\|\nabla^2 w^\alpha\|_{L^{\upinfty}_x}.
\end{align*}
 By taking the $L^p$ norm in $\alpha$ and using the bootstrap assumption, we get for all $s \in [0,T^\star_\delta]$
\begin{align*}
    \| \nabla S^\alpha_1(s) \|_{L^2_x} \lesssim (1+C_{0, \lambda',\delta})\delta^2 e^{-(\lambda+\lambda') s}.
\end{align*}

$\bullet$ \underline{Term $\|  \nabla S_{2,a}\|_{L^2_x}$}: we have 
\begin{align*}
\| \nabla S_{2,a} \|_{L^2_x}&=\| \nabla \big((u-\mathcal{Z}) \cdot \nabla u \big) \|_{L^2_x} \\
& \lesssim \| \nabla u\|_{L^{2}_x}\| \nabla u\|_{L^{\upinfty}_x}+\| u-\mathcal{Z} \|_{L^2_x}\| \nabla^2 u\|_{L^{\upinfty}_x}\\
& \lesssim \| \nabla u\|_{L^{2}_x}\| \nabla u\|_{L^{\upinfty}_x}+\Big(\| u-\langle u \rangle \|_{L^2_x}+| \langle u \rangle-\mathcal{Z} \vert \Big)\|\nabla^2 u \|_{L^{\upinfty}_x}
,\end{align*}
and we can therefore conclude as before to get for all $s \in [0,T^\star_\delta]$
\begin{align*}
    \| \nabla S_{2,a}(s) \|_{L^2_x} \lesssim (1+C_{0,\delta})\delta^2 e^{-(\lambda+\lambda') s}.
\end{align*}

$\bullet$ \underline{Term $\|   \nabla S_{2,b}\|_{L^2_x}$}: we first write 
\begin{align*}
    \Vert \nabla S_{2,b} \Vert_{L^2_x} &\leq \Vert \rr \Vert_{L^{\upinfty}_\alpha W^{1, \infty}_x}\int \left(\Vert w^\alpha-u \Vert_{L^2_x}+\Vert \nabla(w^\alpha-u) \Vert_{L^2_x} \right)\, \mathrm{d}\mu(\alpha) \\
    & \lesssim \Vert \rr \Vert_{L^{\upinfty}_\alpha W^{1, \infty}_x}\int \left(\vert \langle w^\alpha-u \rangle \vert  +\Vert \nabla(w^\alpha-u) \Vert_{L^2_x} \right)\, \mathrm{d}\mu(\alpha),
\end{align*}
by Poincaré inequality. According to the bootstrap assumption combined with Hölder inequality, we already know that the second factor decays like $\delta e^{-\lambda t}$ on $[0,T_\delta^\star]$. For the first factor, we have now from Lemma \ref{LM:pointwise-densité} and the bootstrap assumption that
\begin{align*}
    \Vert \nabla r^\alpha \Vert_{ L^{\upinfty}_x} \lesssim \exp(C \delta)\big[\Vert \nabla r^\alpha_0 \Vert_{ L^{\upinfty}_x}+(1+\Vert  r^\alpha_0 \Vert_{ L^{\upinfty}_x})\delta  \big].
\end{align*}
Recalling the estimate \eqref{bound-perturb-density-bootstrap} for $\Vert \rr \Vert_{L^{\upinfty}_\alpha L^{\upinfty}_x}$, we obtain for all $s \in [0,T^\star_\delta]$
\begin{align*}
    \Vert \nabla S_{2,b}(s) \Vert_{L^2_x} \lesssim C_{\delta, \lambda'}\left( \eps_0+\delta\right)\delta e^{-\lambda s}.
\end{align*}

\medskip

At this stage, we can improve the decay of $\nabla u(t)$ thanks to \eqref{decay-nablau-proof} and the previous estimates on the sources: indeed, we now have for all $t \in [0,T_\delta^\star)$
\begin{multline*}
    \|  \nabla u(t) \|_{L^2_x}  \lesssim \eps_0  e^{-\Lambda t} + ((1+C_{0,\delta}) \delta^2 \\
    +C_{\delta, \lambda'}\left( \eps_0+\delta\right)\delta)\int_0^t e^{-\Lambda (t-s)}\left( e^{-(\lambda+\lambda')s} + e^{-\lambda s}\right)\, \mathrm{d}s.
\end{multline*}
Note that because 
\begin{align*}
    \int_0^t e^{-\Lambda (t-s)} e^{-(\lambda+\lambda')s} \, \mathrm{d}s \lesssim \left\{ \begin{array}{ll}
          e^{-(\lambda+\lambda')t} & \mbox{if } \Lambda-(\lambda+\lambda')>0, \\
          e^{-\Lambda t} & \mbox{if } \Lambda-(\lambda+\lambda')<0, \\
          t e^{-\Lambda t} & \mbox{if } \Lambda-(\lambda+\lambda')=0,
    \end{array} 
    \right.
\end{align*}
we always get
\begin{align*}
   \|  \nabla u(t) \|_{L^2_x}  & \lesssim \eps_0 e^{-\Lambda t} + ((1+C_{0,\delta})\delta^2+C_{\delta, \lambda'}\left( \eps_0+\delta\right)\delta) e^{- \lambda t} \\
& \leq (\eps_0+ (1+C_{0, \delta})\delta^2+C_{\delta, \lambda'}\left( \eps_0+\delta\right)\delta) e^{- \lambda t},
\end{align*}
since $\lambda < \Lambda$. We therefore obtain the desired decay and this concludes the proof.
\end{proof}

Thanks to the previous lemma, we are now in position to improve the decay of  $\|  \nabla \ww(t) \|_{L^{\upinfty}_\alpha L^2_x}$.
\begin{Lem}\label{lem:ImproveLinL2-nabla-walpha}
There exists $C_{0,\lambda',\delta}>0$ such that all $t \in [0, T^\star_\delta]$, we have
\begin{align*}
&    \| \nabla \ww(t) \|_{L^{\upinfty}_\alpha L^2_x}    \lesssim \left((1+C_{0, \lambda',\delta})\eps_0+(1+C_{0, \lambda',\delta})\delta^2 \right) e^{- \lambda t},
\end{align*}
where $C_{0, \lambda,\delta}$ is a continuous and monotone function in all the parameters.
\end{Lem}
\begin{proof}
    Applying $\nabla$ to the equation
$$\partial_t w^\alpha + (\mathcal{W} \cdot \nabla) w^\alpha   =  u- w^\alpha - ((w^\alpha-\mathcal{Z}) \cdot \nabla) w^\alpha,$$
and performing an energy estimate, we obtain
\begin{align*}
    \frac{1}{2}  \frac{\mathrm{d}}{\mathrm{d}t}\Vert \nabla w^\alpha \Vert_{L^2_x}^2 +\Vert \nabla w^\alpha \Vert_{L^2_x}^2 \leq \Vert \nabla u \Vert_{L^2_x} \Vert \nabla w^\alpha \Vert_{L^2_x}+ \Vert \nabla w^\alpha \Vert_{L^{\upinfty}_x}\Vert \nabla w^\alpha \Vert_{L^2_x}^2 ,
\end{align*}
and therefore 
\begin{align*}
    \frac{1}{2}  \frac{\mathrm{d}}{\mathrm{d}t}\Vert \nabla w^\alpha \Vert_{L^2_x} +\Vert \nabla w^\alpha \Vert_{L^2_x}\leq  \Vert \nabla u \Vert_{L^2_x}  + \Vert \nabla w^\alpha \Vert_{L^{\upinfty}_x}^2 .
\end{align*}
Combining the bootstrap assumption on $\Vert \nabla w^\alpha \Vert_{L^{\upinfty}_x}$ with the decay for $\|  \nabla u(t) \|_{L^2_x}$ obtained in Lemma \ref{lem:ImproveLinL2-nabla-u}, we can infer
\begin{align*}
    \frac{\mathrm{d}}{\mathrm{d}t}\Vert \nabla w^\alpha(t) \Vert_{L^2_x} +\Vert \nabla w^\alpha (t)\Vert_{L^2_x}\lesssim \left((1+C_{0, \lambda,\delta})\eps_0+(1+C_{0, \lambda,\delta})\delta^2 \right) e^{- \lambda t},
\end{align*}
and thus
\begin{align*}
    \Vert \nabla \ww(t) \Vert_{L^{\upinfty}_\alpha L^2_x} &\lesssim \eps_0 e^{-t}+ \left((1+C_{0, \lambda',\delta})\eps_0+(1+C_{0, \lambda',\delta})\delta^2 \right)\int_0^t e^{-(t-s)} e^{-\lambda s} \, \mathrm{d}s.
\end{align*}
Since $\lambda<1$, we end up with 
\begin{align*}
    \Vert \nabla \ww(t) \Vert_{L^{\upinfty}_\alpha L^2_x} \lesssim \left((1+C_{0, \lambda',\delta})\eps_0+(1+C_{0, \lambda',\delta})\delta^2 \right) e^{-\lambda t},
\end{align*}
and this concludes the proof.
\end{proof}

As a consequence of Lemmas \ref{lem:ImproveLinL2-nabla-u}--\ref{lem:ImproveLinL2-nabla-walpha}, we can directly improve the decay of the difference of the gradients in $L^2_x$. 
\begin{Lem}\label{lem:ImproveLinL2-diffnabla}
There exists $C_{0,\lambda',\delta}>0$ such that all $t \in [0, T^\star_\delta]$, we have
\begin{align*}
\| \nabla  \ww(t)-\nabla u(t) \|_{L^{\upinfty}_\alpha L^2_x}    \lesssim \left((1+C_{0, \lambda',\delta})\eps_0+(1+C_{0, \lambda',\delta})\delta^2 \right) e^{-\lambda t}.
\end{align*}
\end{Lem}

Lastly, we improve the decay of the difference in average $\Vert  \langle  \ww(t)- u(t) \rangle \Vert_{L^p_\alpha}$.
\begin{Lem}\label{lem:ImproveLdiffMoyenne}
    There exists $C_{0,\lambda',\delta}>0$ such that all $t \in [0, T^\star_\delta]$, we have
\begin{align*}
\Vert  \langle  \ww(t)- u(t) \rangle \Vert_{L^p_\alpha}    \lesssim \left((1+C_{0, \lambda',\delta})\eps_0+(1+C_{0, \lambda',\delta})\delta^2 \right) e^{- \lambda t}.
\end{align*}
\end{Lem}

\begin{proof}
As in the proof of Lemma \ref{lem:ImproveLinL2-nabla-u}, we rely on Proposition~\ref{Prop-lineardecayVNS} with 
\begin{align*}
S^\alpha_1&\vcentcolon=- ((w^\alpha-\mathcal{Z}) \cdot \nabla) w^\alpha, \\
S_2 &\vcentcolon= S_{2,a} + S_{2,b},   \\
S_{2,a}&\vcentcolon=- ((u-\mathcal{Z}) \cdot \nabla) u, \ \ S_{2,b}\vcentcolon= \int r^\alpha (w^\alpha-u) \D \mu(\alpha),
\end{align*}
which provides for all $t \in [0,T_\delta^\star]$
\begin{align*}
\Vert  \langle  \ww(t)- u(t) \rangle \Vert_{L^p_\alpha} &\lesssim \eps_0 e^{-\Lambda t}   + \int_0^t e^{-\Lambda(t-s)} ( \Vert \langle  S^\alpha_1(s) \rangle \Vert_{L^p_\alpha} + \vert \langle S_2(s)  \rangle \vert) \D s,
\end{align*}    
since $\Lambda<1$. It remains to estimate the source terms under the integral in time. Since on the torus, we have $\vert \langle \cdot  \rangle \vert \leq \Vert \cdot \Vert_{L^2_x}$, we can argue  as in the end of the proof of Lemma \ref{lem:ImproveLinL2-nabla-u}.

\medskip

$\bullet$ \underline{Term $\Vert \langle  S^\alpha_1 \rangle \Vert_{L^p_\alpha}$}: we write 
\begin{align*}
    \vert \langle  S^\alpha_1 \rangle \vert  \leq \| w^\alpha-\mathcal{Z} \|_{L^2_x}\| \nabla w^\alpha\|_{L^{\upinfty}_x}
\end{align*}
and we conclude as in the proof of Lemma \ref{lem:ImproveLinL2-nabla-u} to get 
\begin{align*}
    \Vert \langle  S^\alpha_1 \rangle \Vert_{L^p_\alpha} \lesssim (1+C_{0,\delta})\delta^2 e^{-(\lambda+\lambda') s}.
\end{align*}

\medskip

 $\bullet$ \underline{Term $| \langle S_{2,a} \rangle|$}:  we write 
\begin{align*}
| \langle S_{2,a} \rangle|&=| \langle  (u-\mathcal{Z}) \cdot \nabla u \rangle|  \lesssim \Big(\| u-\langle u \rangle \|_{L^2_x}+\| \langle u \rangle-\mathcal{Z} \|_{L^2_x}\Big)\|\nabla u \|_{L^{\upinfty}_x},
\end{align*}
and we can get as in the proof of Lemma \ref{lem:ImproveLinL2-nabla-u} that 
\begin{align*}
    \vert \langle S_{2,a}(s) \rangle \vert \lesssim (1+C_{0,\delta})\delta^2 e^{-(\lambda+\lambda') s}.
\end{align*}

 $\bullet$ \underline{Term $| \langle S_{2,b} \rangle|$}: similarly, we obtain 
\begin{align*}
    |  \langle S_{2,b}(s) \rangle| \lesssim C_{\delta, \lambda'}(\eps_0+\delta) \delta e^{-\lambda s}.
\end{align*}
Looking at the penultimate inequality obtained in the proof of Lemma \ref{lem:ImproveLinL2-nabla-u}, and using $\lambda<\Lambda$, we can conclude in the same way.

\end{proof}

We are finally in position to conclude the bootstrap argument, thanks to Lemmas \ref{lem:ImproveLinL2-nabla-u}--\ref{lem:ImproveLinL2-nabla-walpha}--\ref{lem:ImproveLinL2-diffnabla}--\ref{lem:ImproveLdiffMoyenne}. At the core of our proof is an interpolation procedure enabling us to upgrade the decay of the velocity gradients in $L^2_x$ to a decay in $L^\infty_x$.
\begin{Lem}\label{lem-bootstrap-amelio} 
    There exist $\delta>0$ and $\eps_0>0$ such that for all $k >2+3/2$ and $\eps \in (0, \eps_0)$, we have for all $t \in [0,T_\delta^\star]$
    \begin{align}
\begin{split}\label{ineq-bootstrap-amelio} 
          \ \  \|(\nabla \ww,\nabla^2 \ww)(t)\|_{L^{\upinfty}_\alpha L^{\upinfty}_x} + \|(\nabla u,\nabla^2 u)(t) \|_{L^{\upinfty}_x}  &\leq \frac{\delta}{2} e^{-\lambda' t}, \\
  \ \  \|\nabla \ww(t)\|_{L^{\upinfty}_\alpha L^{2}_x} + \|\nabla u(t) \|_{L^{2}_x}  &\leq \frac{\delta}{2} e^{-\lambda t}, \\
  \Vert  \langle  \ww(t)- u(t) \rangle \Vert_{L^p_\alpha}
 &\leq \frac{\delta}{2} e^{-\lambda t}.
 \end{split}
    \end{align}
\end{Lem}
\begin{proof}
First, we focus on the estimates for the first-order derivatives of the velocity fields and for their difference. By the Gagliardo--Nirenberg--Sobolev, the following functions being mean-free,  we have 
    \begin{equation}\label{eq:GNS-gradient}
\begin{aligned}
 \| \nabla u \|_{L^{\upinfty}_x}   &\leq C_k  \|   \nabla u  \|_{L^2_x}^{\theta_k}    \| \mathrm{D}^k u\|_{L^2_x}^{1-\theta_k}, \\ 
   \| \nabla w^\alpha \|_{L^{\upinfty}_x} &\leq C_k\|  \nabla w^\alpha \|_{L^2_x}^{\theta_k}  \|   \mathrm{D}^k w^\alpha\|_{L^2_x}^{1-\theta_k},
\end{aligned}
\end{equation}
for some $\theta_k \in (0,1)$ and $C_k>0$, with  $k \in \N$ such that $k>2+3/2$. In what follows, constants appearing like $C_k$ only depend on $k$ and can vary from line to another. We first recall the exponential estimate from \eqref{bound-opti-growthexpoVNS} that reads
\begin{multline*}
    \| \ww(t) \|_{\mathcal{H}^{k,p}}+ \| u(t)\|_{H^k_x} \\
     \leq \eps_0 \exp \bigg(\left(\sigma^{-1}+ C_k e^{\delta/\lambda'} C_{\lambda'}(\eps_0+\delta) \right)t+ \delta(1+\sigma) C_{\lambda, \lambda'}C_k\bigg),
\end{multline*}
for all $\sigma>1$. Combining with Lemmas \ref{lem:ImproveLinL2-nabla-u}--\ref{lem:ImproveLinL2-nabla-walpha}, we thus get for all $t \in [0, T^\star_\delta]$
\begin{align}\label{bound-intermediate-decay-grad}
\begin{split}
    &\| \nabla u (t) \|_{L^{\upinfty}_x} 
+\| \nabla \ww (t) \|_{L^{\upinfty}_\alpha L^{\upinfty}_x}   
\\ 
&\lesssim C_k \Big((1+C_{0, \lambda',\delta})\eps_0+(1+C_{0, \lambda',\delta})\delta^2  \Big)^{\theta_k} \eps_0^{1-\theta_k} e^{(1-\theta_k)\delta(1+\sigma) C_{\lambda, \lambda'}C_k} \\
& \qquad \qquad \qquad \qquad  \times \exp\left((1-\theta_k)\left(\sigma^{-1}+ C_k e^{\delta/\lambda'} C_{\lambda'}(\eps_0+\delta) \right)t -\theta_k\lambda  t\right).
\end{split}
\end{align}
From Lemmas \ref{lem:ImproveLinL2-diffnabla}--\ref{lem:ImproveLdiffMoyenne}, we also get the additional estimate 
\begin{align}\label{bound-intermediate-decay-diff}
\Vert  \langle  \ww(t)- u(t) \rangle \Vert_{L^p_\alpha}
    \leq C \left((1+C_{0, \lambda',\delta})\eps_0+(1+C_{0, \lambda',\delta})\delta^2 \right) e^{-\lambda t},
\end{align}
for all $t \in [0, T^\star_\delta]$.
\medskip

We then proceed as follows. Recall that we have fixed any $\lambda<\Lambda$ and $\lambda' \in (0, \theta \lambda)$, where $\theta \in (0,1)$ was left free. We fix $\theta=\theta_k$. Next, we take $\sigma>0$ large enough and  $\eps_0, \delta$ small enough so that 
\begin{align*}
(1-\theta_k)\left(\sigma^{-1}+ C_k e^{\delta/\lambda'} C_{\lambda'}(\eps_0+\delta) \right) -\theta_k\lambda  &<-\lambda', \\
     C_k \Big((1+C_{0, \lambda',\delta})\eps_0+(1+C_{0, \lambda',\delta})\delta^2  \Big)^{\theta_k} \eps_0^{1-\theta_k} e^{(1-\theta_k)\delta(1+\sigma) C_{\lambda, \lambda'}C_k} &<\frac{\delta}{2}, \\
     C \left((1+C_{0, \lambda',\delta})\eps_0+(1+C_{0, \lambda',\delta})\delta^2 \right) &<\frac{\delta}{2},
\end{align*}
which is possible since $\theta_k \lambda-\lambda'>0$. This can be done while also imposing $100\eps_0 <\delta$, so that the time $T^\star_\delta$ from \eqref{assump-bootstrapVNS} is indeed well-defined,  and eventually ensures the following improvement of \eqref{bound-intermediate-decay-grad}--\eqref{bound-intermediate-decay-diff}: for all $t \in [0, T^\star_\delta]$, there holds
\begin{align*}
         \|\nabla \ww(t)\|_{L^{\upinfty}_\alpha L^{\upinfty}_x} + \|\nabla u(t) \|_{L^{\upinfty}_x}  &\leq \frac{\delta}{2} e^{-\lambda' t}, \\
  \Vert  \langle  \ww(t)- u(t) \rangle \Vert_{L^p_\alpha}
 +\|\nabla \ww(t)-\nabla u(t)\|_{L^{\upinfty}_\alpha L^{2}_x} &\leq \frac{\delta}{2} e^{-\lambda t}.
    \end{align*}
The same procedure applies \textit{mutatis mutandis} for the second-order derivatives of the velocity fields: indeed, similar estimates as \eqref{eq:GNS-gradient} hold for the second-order derivatives on the left-hand side, with other constants $\widetilde{C}_k>0$ and other exponents $\widetilde{\theta}_k \in (0,1)$. It therefore concludes the proof.
\end{proof}

\subsection{Conclusion and final asymptotics}\label{sectionVNS-cclbootstrap}

It is now fairly standard to conclude the bootstrap argument.

\begin{Lem}\label{Lem-cclBootstrap}
There exists $\eps_0>0$ such that for all $k \in \N$ with $k > 2+3/2$ and all $\eps \in (0,\eps_0)$, the perturbed system \eqref{eq:VNSpertub} admits a unique global in time solution $(\rr, \vv, u) $ (in the class of Theorem \ref{thm-globalconvergence}). Furthermore, for all $0<\lambda<\Lambda$, there exists $C_\lambda>0$ such that for all $t \geq 0$
\begin{align}\label{eq:boundfromProp-cclBootstrap}
      \|\nabla \ww(t)\|_{L^{\upinfty}_\alpha L^{2}_x} + \|\nabla u(t) \|_{L^{2}_x}+ \Vert  \langle  \ww(t)- u(t) \rangle \Vert_{L^p_\alpha}  \leq C_\lambda e^{-\lambda t}.
\end{align}
\end{Lem}
\begin{proof}
As explained in the set-up of the bootstrap argument, we assume by contradiction that $T_\delta^\star<T_{[0]}^{\mathrm{max}}$. Invoking Lemma \ref{lem-bootstrap-amelio}, we know that there exist $\delta>0$ and $\eps_0>0$ such that for all  $\eps \in (0, \eps_0)$, the improved estimate \eqref{ineq-bootstrap-amelio} holds for all $t \in [0,T_\delta^\star]$. These estimates contradict the definition of $T_\delta^\star$ and we thus must have $T_\delta^\star=T_{[0]}^{\mathrm{max}}$. In particular, in view of the blow-up criterion and the exponential decay in time of $\|\nabla w^\alpha(t)\|_{L^{\upinfty}_x}$ and $\|\nabla u(t) \|_{L^{\upinfty}_x}$ on $[0,T_{[0]}^{\mathrm{max}}]$, we infer that we must have $T_{[0]}^{\mathrm{max}}=+\infty$.  The claimed estimate follows directly from the definition of $T_\delta^\star$.
\end{proof}
In what follows, let $\lambda \in (0, \Lambda)$ be fixed.
\paragraph{Convergence of the averages.}
In view of \eqref{eq:boundfromProp-cclBootstrap} and Poincaré inequality, we also infer the estimate
\begin{align}\label{decay:u-minus-average}
    \Vert  u(t)-\langle u(t) \rangle  \Vert_{L^2_x} + \Vert  \ww(t)-\langle \ww(t) \rangle   \Vert_{L^{\upinfty}_\alpha L^2_x} \leq &C_\lambda e^{-\lambda t},
\end{align}
It remains to understand the limit of the averages $\langle u(t) \rangle$ and $\langle w^\alpha(t) \rangle$. We rely on the preliminary Lemma \ref{smallLemma-moyenne}: since 
$ \Vert \rhorho \Vert_{L^{\upinfty}_\alpha L^2_x}$ is uniformly bounded for all times (thanks to Lemma \ref{LM:pointwise-densité} and the bootstrap assumption), we have for all $t>0$
\begin{align*}
    &\qquad \vert \langle u(t) \rangle -\mathcal{Z} \vert +
        \Vert \langle \ww(t) \rangle -\mathcal{Z} \Vert_{L^p_\alpha} \\
        &\lesssim \Vert \langle \ww(t)-u(t) \rangle \Vert_{L^p_\alpha}  +\Vert \nabla u(t) \Vert_{L^2_x} + \int \Vert w^\alpha(t)-u(t) \Vert_{L^2_x}\, \mathrm{d}\mu(\alpha) \\
        &\lesssim \Vert \langle w^\alpha(t)-u(t) \rangle \Vert_{L^p_\alpha}  +\Vert \nabla u(t) \Vert_{L^2_x} + \int \vert \langle  w^{\alpha}(t) -u(t)\rangle \vert \, \mathrm{d}\mu(\alpha) \\
        & \quad + \int \Vert \nabla (w^\alpha(t)-u(t)) \Vert_{L^2_x} \, \mathrm{d}\mu(\alpha) \\
        & \lesssim \Vert \langle \ww(t)-u(t) \rangle \Vert_{L^p_\alpha}  +\Vert \nabla u(t) \Vert_{L^2_x} +\|\nabla \ww(t)-\nabla u(t)\|_{L^{\upinfty}_\alpha L^{2}_x},
\end{align*}
by Poincaré  inequality. Again by \eqref{eq:boundfromProp-cclBootstrap}, we deduce that for all $t>0$, we have 
\begin{align}\label{eq:decay-<u>-Z}
     \vert \langle u(t) \rangle -\mathcal{Z} \vert +
        \Vert \langle \ww(t) \rangle -\mathcal{Z} \Vert_{L^p_\alpha} &\leq C_\lambda e^{-\lambda t}.
\end{align}
We thus infer the desired convergence \eqref{final-convergence-velocities} since $\mathcal{W}=\mathcal{V}+\mathcal{Z}$.

\paragraph{Decay of higher order derivatives.}
Let us improve the previous convergences by showing that the following exponential decay holds for our global in time solutions: for any $k >2+3/2$, there exists $C, \gamma >0$ such that for all $t>0$
\begin{align}\label{decay-final-highderivative}
\Vert \ww(t) -\langle \ww(t) \rangle \Vert_{\mathcal{H}^{k,p}}+ \Vert u(t) -\langle u(t) \rangle \Vert_{H^k_x} \leq C e^{-\gamma t}.
\end{align}
Combined with \eqref{eq:decay-<u>-Z} , this will entail \eqref{final-convergence-velocities2}. 

\paragraph{Step 1: decay in $H^{k-1}_x$.} First, coming back to the exponential estimates from \eqref{bound-opti-growthexpoVNS} -- and up to taking $\sigma$ larger, and $\delta$ and $\eps_0$ smaller, if necessary -- we know that we can obtain an exponential control of the top-order derivatives (of order $k$) of the form
\begin{align}\label{bound:expo-toporderderivative}
    \Vert \nabla^k \ww(t) \Vert_{L^{\upinfty}_\alpha L^2_x}+\Vert \nabla^k u(t) \Vert_{L^2_x} \lesssim \exp(\eps t),
\end{align}
for any $\eps>0$. The constants that are involved may depend on $\eps, \eps_0$ or $\delta$ but this is harmless at this stage.

\medskip

Relying on the same interpolation procedure as in \eqref{eq:GNS-gradient}, we know that for all $ 2 \leq i \leq k-1$, we have
\begin{align*}
      \Vert \nabla^i u \Vert_{L^2_x} 
   &\leq \Vert \nabla u \Vert_{L^2_x}^{1-\theta_{i,k}} \Vert \nabla^k u \Vert_{L^2_x}^{\theta_{i,k}},
\end{align*}
with $\theta_{i,k}=\frac{i-1}{k-1}$. Using the exponential decay from Lemma \ref{Lem-cclBootstrap}, we thus obtain
\begin{align*}
  \forall t>0, \ \ \Vert \nabla^i u(t) \Vert_{L^2_x} \leq C \exp((-(1-\theta_{i,k}) \lambda + \theta_{i,k} \eps)t).
\end{align*}
Hence, if we choose $\gamma >0$ such that
\begin{align*}
    \gamma <\frac{1}{k-1}\lambda \leq \frac{k-i}{k-1}\lambda, \ \ 2 \leq i \leq k-1,
\end{align*}
we then have 
\begin{align*}
    0< \frac{k-i}{i-1}\lambda - \frac{k-1}{i-1}\gamma=\frac{1-\theta_{i,k}}{\theta_{i,k}} \lambda -\frac{1}{\theta_{i,k}}\gamma, \ \ 2 \leq i \leq k-1,
\end{align*}
and by picking $\eps>0$ such that
\begin{align*}
    \eps< \underset{2 \leq i \leq k-1}{\inf} \left\lbrace  \frac{1-\theta_{i,k}}{\theta_{i,k}} \lambda -\frac{1}{\theta_{i,k}}\gamma \right\rbrace,
\end{align*}
we end up with $-(1-\theta_{i,k}) \lambda + \theta_{i,k} \eps<\gamma$ for all $2 \leq i \leq k-1$. We obtain for all $t>0$ and  $ 2 \leq i \leq k-1$
\begin{align*}
\Vert \nabla^i u(t) \Vert_{L^2_x} \lesssim e^{-\gamma t}.
\end{align*}
By the same argument, the previous estimate can also be inferred on $\Vert \nabla^i \ww(t) \Vert_{L^{\upinfty}_\alpha L^2_x}$, and then on the difference $\Vert \nabla^i \ww(t)-\nabla^i u(t) \Vert_{L^{\upinfty}_\alpha L^2_x}$ for all $ 2 \leq i \leq k-1$. In addition, note that Lemma \ref{Lem-cclBootstrap} also entails (by Poincaré inequality) the decay
\begin{align*}
    \Vert (u-\ww)(t) \Vert_{L^p_\alpha L^2_x} \lesssim e^{-\lambda t}.
\end{align*}
All in all, we have obtained the existence of $\gamma>0$ such that, for all $t>0$,
    \begin{align}
    \label{eq:decay-diff-H^(k-1)}\Vert (u-\ww)(t) \Vert_{L^p_\alpha H^{k-1}_x} &\lesssim e^{-\gamma t},\\
   \label{eq:decay-velocities-H^(k-1)} \Vert \nabla^i \ww(t) \Vert_{L^{\upinfty}_\alpha L^2_x}+\Vert \nabla^i u(t) \Vert_{L^2_x} &\lesssim e^{-\gamma t}, \ \ 1 \leq i \leq k-1.
\end{align}

\medskip

\paragraph{Step 2: decay for Navier--Stokes equations in $H^k_x$.} Next, we recast the Navier--Stokes equations satisfied by the velocity $u$ as the following linear Stokes equations:
\begin{equation}
\left\{ 
\begin{aligned}
    \partial_t u  +(\mathcal{W} \cdot \nabla) u- \Delta u + \nabla p &= \mathfrak{F}, \\
    \mathrm{div}(u)&=0,
\end{aligned}
\right.
\end{equation}
where the source term is 
\begin{align*}
\mathfrak{F}\vcentcolon=-((u-\mathcal{Z}) \cdot \nabla) u+\int (1+r^\alpha) (w^\alpha-u) \D \mu(\alpha)= \vcentcolon \mathfrak{F}_1+\mathfrak{F}_2. 
\end{align*}
By standard estimates on the linear Stokes system on the torus, we know that $\Vert u-\langle u \rangle \Vert_{H^k_x}$ decays exponentially fast provided that the source term $\mathfrak{F}$ enjoys a pointwise exponential decay in $H^{k-1}_x$.

For $\mathfrak{F}_1$, we write by the algebra property of $H^{k-1}_x$ (since $k-1>d/2$) that
\begin{align*}
    \Vert \mathfrak{F}_1 \Vert_{H^{k-1}_x} \leq \Vert u-\mathcal{Z} \Vert_{H^{k-1}_x}\Vert \nabla u \Vert_{H^{k-1}_x}
     \lesssim \left( \vert \langle u \rangle - \mathcal{Z} \vert +\Vert \nabla u \Vert_{H^{k-2}_x} \right) \Vert  u \Vert_{H^{k}_x}.
\end{align*}
In view of \eqref{eq:decay-velocities-H^(k-1)} and \eqref{eq:decay-<u>-Z}, we know that all the terms in the prefactor enjoy an exponential decay, while the second term is bounded by $e^{\eps t}$ for any $\eps>0$, thanks to \eqref{bound:expo-toporderderivative}: more precisely, we get
\begin{align*}
    \Vert \mathfrak{F}_1(t) \Vert_{H^{k-1}_x} \lesssim e^{(-\gamma+\eps)t}.
\end{align*}
For $\mathfrak{F}_2$, we write in a similar way (since $k-1>3/2$) that
\begin{align*}
    \Vert \mathfrak{F}_2(t) \Vert_{H^{k-1}_x} &\leq \Vert (1+\rr)(t) \Vert_{L^{\upinfty}_\alpha H^{k-1}_x}   \Vert (\ww-u)(t) \Vert_{L^p_\alpha H^{k-1}_x} \\
    &\lesssim \left(1+\Vert \rr(t) \Vert_{L^{\upinfty}_\alpha H^{k-1}_x} \right)e^{-\gamma t}, 
\end{align*}
thanks to the estimate \eqref{eq:decay-diff-H^(k-1)}. Let us now bound the density term involving $r^\alpha$: relying on the first line of \eqref{eq:estim-Sobolev-perturb-ralpha}, we have 
\begin{align*}
    \frac{\mathrm{d}}{\mathrm{d}t}\Vert r^\alpha \Vert_{H^{k-1}_x} \lesssim \Vert \nabla w^\alpha \Vert_{L^{\upinfty}_x}\Vert r^\alpha \Vert_{H^{k-1}_x} +\Vert \nabla w^\alpha \Vert_{H^{k-1}_x}\Vert r^\alpha \Vert_{L^{\upinfty}_x}+ \Vert \mathrm{div}(w^\alpha) \Vert_{H^{k-1}_x}.
\end{align*}
By Gronwall inequality, we get
\begin{align*}
    \Vert \rr(t) \Vert_{L^{\upinfty}_\alpha  H^{k-1}_x} &\leq \left( \Vert \rr_0 \Vert_{L^{\upinfty}_\alpha  H^{k-1}_x}+ \Vert \nabla \ww(t) \Vert_{L^{\upinfty}_\alpha  H^{k-1}_x} \right)\exp\left(\int_0^t \Vert \nabla \ww(s) \Vert_{L^{\upinfty}_\alpha L^{\upinfty}_x} \, \mathrm{d}s \right) \\
    & \lesssim \Vert \rr_0 \Vert_{L^{\upinfty}_\alpha  H^{k-1}_x}+ e^{\eps t},
\end{align*}
for any $\eps>0$, thanks to \eqref{bound:expo-toporderderivative} and the inequalities from \eqref{assump-bootstrapVNS}. Coming back to $\mathfrak{F}_2$, we end up with
\begin{align*}
    \Vert \mathfrak{F}_2(t) \Vert_{H^{k-1}_x}   \lesssim \left(1+\Vert \rr_0 \Vert_{L^{\upinfty}_\alpha  H^{k-1}_x}+ e^{\eps t} \right)e^{-\gamma t}.
\end{align*}
All in all, the previous lines show that for any $t>0$, we have
\begin{align*}
     \Vert \mathfrak{F}(t) \Vert_{H^{k-1}_x} \lesssim e^{(\eps-\gamma)t}.
\end{align*}
Taking $\eps$ small enough to ensure that $\eps<\gamma$, we obtain 
\begin{align*}
    \Vert \mathfrak{F}(t) \Vert_{H^{k-1}_x}   \lesssim e^{-\gamma' t}, 
\end{align*}
for some $\gamma'>0$. Since we already know from \eqref{decay:u-minus-average} that $\Vert u-\langle u \rangle \Vert_{L^2_x}$ displays some exponential decay, we can conclude that there exists $\gamma>0$ such that 
\begin{align*}
    \Vert u(t)-\langle u(t) \rangle \Vert_{H^k_x}  \lesssim e^{-\gamma t}.
 \end{align*}
\paragraph{Step 3: decay for pressureless Euler in $L^\infty_\alpha H^{k}_x$.}  We can now estimate $\Vert \nabla \ww \Vert_{L^{\upinfty}_\alpha H^{k-1}_x}$ by a direct energy estimate. Thanks to Lemma \ref{LM:existence-vitesseSOBOLEV}, we get for $\mu$-almost all $\alpha$ that
\begin{multline*}
    \frac{\mathrm{d}}{\mathrm{d}t} \Vert \nabla w^\alpha \Vert^2_{H^{k-1}_x}+ \Vert \nabla w^\alpha \Vert^2_{H^{k-1}_x} \\
    \lesssim \Vert \nabla u \Vert_{H^{k-1}_x}\Vert \nabla w^\alpha \Vert_{H^{k-1}_x}+ \Vert \nabla w^\alpha \Vert_{L^{\upinfty}_x}\Vert \nabla w^\alpha \Vert_{H^{k-1}_x}^2.
\end{multline*}
According to the decay of $\Vert \nabla u \Vert_{H^{k-1}_x}$ from \textbf{Step 2} and to the decay of $\Vert \nabla  \ww \Vert_{L^{\upinfty}_\alpha L^{\upinfty}_x}$ from the bootstrap assumption \eqref{assump-bootstrapVNS}, we obtain
\begin{align*}
    \frac{\mathrm{d}}{\mathrm{d}t} \Vert \nabla  \ww \Vert_{L^{\upinfty}_\alpha H^{k-1}_x} \lesssim -\big(1-e^{-\gamma t} \big) \Vert \nabla  \ww \Vert_{L^{\upinfty}_\alpha H^{k-1}_x} + e^{-\gamma t}, 
\end{align*}
for some $\gamma \in (0,1)$. From Gronwall lemma, using the fact that 
$$\exp\left( -\int_0^t (1-e^{-\gamma s}) \, \mathrm{d}s \right) \lesssim e^{-t},
$$
we infer 
\begin{align*}
    \Vert \nabla  \ww (t) \Vert_{L^{\upinfty}_\alpha H^{k-1}_x} \lesssim e^{-t} \Vert \nabla  \ww_0 \Vert_{L^{\upinfty}_\alpha H^{k-1}_x}+ \int_{0}^t e^{-\gamma s} e^{-(t-s)} \, \mathrm{d}s \lesssim e^{-\gamma t}.
\end{align*}
Since we already know from \eqref{decay:u-minus-average} that $\Vert \ww-\langle \ww \rangle \Vert_{L^{\upinfty}_\alpha L^2_x}$ decays exponentially fast in time, we can now conclude that the same holds for $\Vert \ww -\langle \ww \rangle \Vert_{\mathcal{H}^{k,p}}$.

\medskip

Combining the conclusions from \textbf{Step 2} and \textbf{Step 3} altogether, we eventually obtain \eqref{decay-final-highderivative}.

\bigskip

\paragraph{Asymptotic density profile.}
To conclude and obtain \eqref{final-convergence-densities}, it remains to investigate the limit of $(\rho^{\alpha})$. 
\begin{Prop}\label{Prop-exist-asympto-density}
There exists an asymptotic profile $\rhorho_{\infty} \in L^{\infty}_\alpha L^{\infty}_x$ such that for $\mu$-almost all $\alpha$,
$
\langle \rho^\alpha_0 \rangle = \langle \rho^{\alpha}_{\infty} \rangle
$
and
\begin{align*}
\esssup_{\alpha \in (I, \mu)} & \, \, \mathrm{W}_1 \Big(\rho^{\alpha}(t),\rho^{\alpha}_{\infty}(x-t \mathcal{W} ) \Big) \underset{t \rightarrow + \infty}{\longrightarrow}  0,
\end{align*}
with exponential rate, where $\mathrm{W}_1$ is the $1-$Wasserstein distance on $\T^3$ (see Definition \eqref{def:Wasserstein1}) and where
$$\mathcal{W}=\frac{1}{2}\left\langle U_0+\int \rho^\alpha_0 v^\alpha_0 \D \mu(\alpha) \right\rangle.$$
\end{Prop}
\begin{proof}
We closely follow the approach of \cite{HKMM}. We rely on Cauchy's criterion in the space $L^{\infty}_\alpha L^{\infty}_x$. For the unperturbed fluid density and velocity $(\rhorho,\vv)$, let us set
\begin{align*}
\overline{\rho}^{\alpha}(t,x)\vcentcolon=\rho^{\alpha}(t,x+t \mathcal{W}), \ \ \overline{v}^{\alpha}(t,x)\vcentcolon=v^{\alpha}(t,x+t \mathcal{W}),
\end{align*}
which satisfy
\begin{align*}
\partial_t \overline{\rho}^{\alpha}=\mathrm{div}_x\Big(\overline{\rho}^{\alpha}(\overline{v}^{\alpha}-\mathcal{W}) \Big).
\end{align*}
Let $\psi \in \mathrm{Lip}(\T^3)$ such that $\Vert \nabla_x \psi \Vert_{L^{\upinfty}_x} \leq 1$. For all $t \geq s \geq 0$, we have
\begin{align*}
\int_{\T^3} \psi \overline{\rho}^{\alpha}(t)-\int_{\T^3} \psi \overline{\rho}^{\alpha}(s)=\int_0^t \int_{\T^3} \nabla_x \psi (\overline{v}^{\alpha}(\tau)-\mathcal{W}) \, \mathrm{d}\tau.
\end{align*}
By using the uniform control on the densities $(\rho^\alpha)$, we get 
\begin{align*}
\left\vert \int_{\T^3} \psi \overline{\rho}^{\alpha}(t)-\int_{\T^3} \psi \overline{\rho}^{\alpha}(s) \right\vert &\lesssim \Vert \rho^\alpha_0 \Vert_{L^{\upinfty}_x} \int_s^t \int_{\T^3}  \left\vert \overline{v}^{\alpha}(\tau)-\mathcal{W} \right\vert \, \mathrm{d}\tau \\
& = \Vert \rho^\alpha_0 \Vert_{L^{\upinfty}_x} \int_s^t \int_{\T^3}  \left\vert v^\alpha(\tau)-\mathcal{W} \right\vert \, \mathrm{d}\tau \\
& \leq \Vert \rhorho_0 \Vert_{L^{\upinfty}_\alpha L^{\upinfty}_x} \int_s^t  \left\Vert \vv(\tau)-\mathcal{W} \right\Vert_{L^{\upinfty}_\alpha L^2_x} \, \mathrm{d}\tau.
\end{align*}
Thanks to Monge--Kantorovitch duality formula (see \eqref{def:Wasserstein1}), this yields
\begin{align*}
\esssup_{\alpha \in (I, \mu)} \mathrm{W}_1 \Big(\overline{\rho}^{\alpha}(t),\overline{\rho}^{\alpha}(s) \Big) \lesssim  \Vert \rhorho_0 \Vert_{L^{\upinfty}_\alpha L^{\upinfty}_x} \int_s^t  \left\Vert \vv(\tau)-\mathcal{W} \right\Vert_{L^{\upinfty}_\alpha L^2_x} \, \mathrm{d}\tau. 
\end{align*}
By~\eqref{eq:decay-<u>-Z}--\eqref{decay-final-highderivative}, namely the exponential decay of $\left\Vert \vv(\tau)-\mathcal{W} \right\Vert_{L^{\upinfty}_\alpha L^2_x}$ when $\tau \rightarrow + \infty$, we can now use Cauchy's criterion (and an additional weak compactness argument based on the uniform bound for the densities) to get the existence of $\rhorho_{\infty} \in L^{\infty}_\alpha L^{\infty}_x$, with 
$
\langle \rho^\alpha_0 \rangle = \langle \rho^{\alpha}_{\infty} \rangle
$
, such that for almost all $\alpha$, $(\overline{\rho}^{\alpha})$ weakly converges towards this profile,
at exponential rate. We obtain the desired estimate by a final change of variable.
\end{proof}

\sectionhead
  [Proof of Theorem \ref{thm-globalconvergence-comp} for the compressible VNS system]
  {Proof of Theorem \ref{thm-globalconvergence-comp} for the compressible Vlasov--Navier--Stokes system}\label{Section:VNS-comp}

In this section, we focus on the compressible case,  that is we study \eqref{eq:VNScompressibleCHAP4}, and detail the proof of Theorem \ref{thm-globalconvergence-comp}. The proof will essentially follow the same lines of proof as for Theorem  \ref{thm-globalconvergence} and our goal is to explain the required modifications in the strategy. As in the former section, we will consider the case where the drag term is linear, that is taking $\Gamma=0$. The case of a nonlinear drag term will be handled in Section \ref{Section-nonlindrag}.

\medskip

Recall that we denote $\ww= (w^\alpha)_{\alpha \in I}$. In what follows, we write the initial smallness assumption as
\begin{align*}
   C_{0,k}\vcentcolon=  \| \rr_0\|_{L^{\upinfty}_\alpha H^{k-1}_x}+\Vert  n_0  \Vert_{H^{k}_x}+\Vert \bm{w}_0 \Vert_{\mathcal{H}^{k,p}}+ \| u_0\|_{H^{k}_x} \leq \eps_0,
\end{align*}
that encodes the size of the initial perturbation (see Theorem \ref{thm-globalconvergence-comp})
\subsection{Refined energy estimates}

We first refine the local well-posedness result for the perturbation satisfying \eqref{eq:VNScomp-pertub1}. We require that we work on a time interval $[0,T]$ such that for all $t \in [0,T]$, the fluid density satisfies
\begin{align}\label{assumption-novoid-theta}
    \inf_{x} \varrho(t) >0.
\end{align}
In view of Lemma \ref{LM:pointwise-densité}, this property will hold for $T$ small enough, since it is true initially -- see Assumption \eqref{3rd-novacuum-comp}. We will propagate it for all times (keeping the perturbation $ n$ small in $L^\infty_x$) in the subsequent bootstrap argument. Dividing by the fluid density $\varrho$ in the momentum equation, we find that the perturbation ($r^\alpha, w^\alpha,  n, u)$ satisfies
\begin{equation}
\label{eq:VNScomp-pertubLWP}
\left\{  
\begin{aligned}
\partial_t r^\alpha + \mathcal{V}\cdot\nabla r^\alpha + &\Div (w^\alpha) + \Div(r^\alpha w^\alpha) = 0,\\
\partial_t w^\alpha + \bigl((\mathcal{V}+w^\alpha)&\cdot\nabla\bigr) w^\alpha = u-w^\alpha,\\
\partial_t  n + (\mathcal{V}+u)\cdot\nabla  n &+ \Div (u) = - n\,\Div(u),\\
\partial_t u + \bigl((\mathcal{V}+u)\cdot\nabla\bigr)u &- \mathcal{A}u + \nabla n
\\
&= \int (w^\alpha-u)\,\mathrm{d}\mu(\alpha)
 + \int r^\alpha(w^\alpha-u)\,\mathrm{d}\mu(\alpha) \\
&\quad - Q( n)\,\nabla n
 - K( n)\,\mathcal{A}u \\
&\quad - K( n)\int (w^\alpha-u)\,\mathrm{d}\mu(\alpha)
 - K( n)\int r^\alpha(w^\alpha-u)\,\mathrm{d}\mu(\alpha),
\end{aligned}
\right.
\end{equation}
where $$K( n)\vcentcolon=\frac{ n}{(1+ n)}, \ \ Q( n)\vcentcolon=\frac{P'(1+ n)-(1+ n)}{1+ n}.$$

Let us remind the reader that we work under the assumption that $P'(1)=1$, hence $K(0)=Q(0)=0$.

\medskip 

In the spirit of Proposition \ref{prop-refinedLWP-VNSincomp} for the incompressible case, the goal is now to derive a refined blow-up criterion for \eqref{eq:VNScomp-pertubLWP}, with tame   exponential growth.

\begin{Prop}\label{prop-refinedLWP-VNScomp}
Let $k \in \N$ with $k>2+3/2$. There exists $C_k>0$ such that for all $\sigma>1$ and for all $T>0$, all smooth solutions to \eqref{eq:VNScomp-pertubLWP} on $[0,T]$ such that
\begin{align}\label{bound-technical-theta-VNScompressible}
    \Vert  n \Vert_{L^\infty([0,T];L^\infty_x)} \leq \min\left( \frac{1}{2C_k}, \frac{1}{2}\right)
\end{align}
satisfy on $[0,T]$  the following exponential bound:
\begin{multline*}
    \| \ww(t) \|_{\mathcal{H}^{k,p}}+ \| u(t)\|_{H^k_x}+ \Vert  n(t) \Vert_{H^k_x} \\
    \leq \eps_0 \exp\left(\sigma^{-1}t  + C_k\int_0^t \mathcal{Q}_{k, \sigma} \Big[ n(s), u(s), \rr(s),\ww(s) \Big] \, \mathrm{d}s  \right),
\end{multline*}
    where
\begin{align}\label{eq:fct-blowup-comp}
\begin{split}
    &\mathcal{Q}_{k, \sigma}\Big[ n(s), u(s), \rr(s),\ww(s) \Big]\\
    &\vcentcolon= 
     \varphi_k(\Vert  n \Vert_{L^{\upinfty}_x}) \bigg(\Vert  n \Vert_{L^{\upinfty}_x}+\Vert ( n, \nabla  n) \Vert_{L^{\upinfty}_x}^2+(1+\Vert \rr \Vert_{L^{\upinfty}_\alpha L^{\upinfty}_x})\Vert \rr \Vert_{L^{\upinfty}_\alpha L^{\upinfty}_x}\\
    & \quad \qquad  \qquad \qquad+\Vert  n \Vert_{L^{\upinfty}_x}^2\Vert \rr \Vert_{L^{\upinfty}_\alpha L^{\upinfty}_x}^2 
     +(\Vert \rr \Vert_{L^{\upinfty}_\alpha L^{\upinfty}_x}^2+\sigma \Vert  n \Vert_{L^{\upinfty}_x}^2)\Vert \ww-u  \Vert_{L^2_\alpha L^{\upinfty}_x}^2 \\ 
     & \qquad \qquad \qquad \quad + \Vert \rr \Vert_{L^{\upinfty}_\alpha L^{\upinfty}_x}\Vert \ww-u  \Vert_{L^2_\alpha L^{\upinfty}_x}
     +\Vert \ww-u  \Vert_{L^1_\alpha L^{\upinfty}_x}^2 \\
     & \quad \qquad \qquad \qquad
      +\Vert (\nabla \ww, \nabla^2 \ww) \Vert_{L^{\upinfty}_\alpha L^{\upinfty}_x}+\Vert (\nabla u, \nabla^2 u) \Vert_{L^{\upinfty}_x} \bigg),
    \end{split}
\end{align}
in which $\varphi_k$ is a continuous positive function.
\end{Prop}
\begin{Rem}
    The assumption \eqref{bound-technical-theta-VNScompressible} will be  automatically satisfied when we will apply Proposition~\ref{prop-refinedLWP-VNScomp} in the subsequent bootstrap procedure.
\end{Rem}
\begin{proof}[Proof of Proposition \ref{prop-refinedLWP-VNScomp}]
First, in view of \eqref{bound-technical-theta-VNScompressible}, since we have $\vert  n \vert \leq \frac{1}{2}$, we easily observe that for all $t \in [0,T]$
\begin{align*}
    \Vert K( n(t)) \Vert_{L^{\upinfty}_x}\lesssim  \Vert  n(t) \Vert_{L^{\upinfty}_x}, \qquad \Vert K'( n(t)) \Vert_{L^{\upinfty}_x} \lesssim 1.
\end{align*}
We will repeatedly use these bounds throughout the proof.
In what follows, the constant $C_k$ (resp. the function $\varphi_k$) will refer to a constant (resp. a continuous positive function) allowed to change from line to line and only depending on $k$. For $\sigma>1$, we introduce
\begin{align*}
    y_{k,\sigma}(t)&=\int \| w^\alpha(t)\|_{H^k_x}^2 \D \mu(\alpha)+\int \| \mathrm{div}(w^\alpha(t))\|_{H^{k-1}_x}^2 \D \mu(\alpha)  \\
    & \quad +  \sigma^{-1}\int \| r^\alpha(t)\|_{H^{k-1}_x}^2 \D \mu(\alpha) + \| u(t)\|_{H^k_x}^2+\|  n(t)\|_{H^k_x}^2.
\end{align*}
We follow the idea of proof from Proposition \ref{prop-refinedLWP-VNSincomp} in the incompressible case, the main differences being as follows.

$\bullet$ It turns out that the energy estimate for the fluid velocity $u$ solution to the momentum equation in the compressible Navier--Stokes system (last equation in \eqref{eq:VNScomp-pertubLWP}) is more intricate compared to the incompressible case. The main differences come from the four last terms in the right-hand side, and the linear contribution from term $\nabla  n$ in the left-hand side. However,  the latter will get compensated by using the coupled equation on $ n$ in \eqref{eq:VNScomp-pertubLWP}, which reflects an underlying hyperbolic structure of the system. 

More precisely, let us perform an $H^k_x$ estimate for $u$ and $w^\alpha$: with an integration by parts in the dissipative operator, we first get as in the proof of Proposition \ref{prop-refinedLWP-VNSincomp} that
\begin{align*}
&\frac{1}{2}\frac{\mathrm{d}}{\mathrm{d}t} \left(\Vert  u \Vert_{H^k_x}^2+\int \| w^\alpha\|_{H^k_x}^2 \D \mu(\alpha) +\right)
     +\Vert (\nabla, \mathrm{div}) u \Vert_{H^k_x}^2 \\
     &\quad + \int \Vert u-w^\alpha \Vert^2_{H^k_x} \, \mathrm{d}\mu(\alpha)
 + \langle \nabla  n, u \rangle_{H^k_x} \\
    & \qquad \qquad \qquad \leq  C_k \left(\int \Vert \nabla w^\alpha \Vert_{L^{\upinfty}_x} \| w^\alpha\|_{H^k_x}^2 \D \mu(\alpha) + \Vert \nabla u \Vert_{L^{\upinfty}_x} \| u\|_{H^k_x}^2 \right)  \\
    & \qquad \qquad  \qquad \quad  + \Vert \rr \Vert_{L^{\upinfty}_\alpha L^{\upinfty}_x} \Vert u \Vert_{L^2_x}\int  \Vert w^\alpha-u  \Vert_{L^2_x} \, \mathrm{d}\mu(\alpha)\\
    & \qquad \qquad  \qquad \quad  +\Vert \nabla u \Vert_{H^k_x}\int \Vert r^\alpha (w^\alpha-u) \Vert_{H^{k-1}_x} \, \mathrm{d}\mu(\alpha) +S,
\end{align*}
where $S\vcentcolon=(\mathrm{I})+(\mathrm{II})+(\mathrm{III})+(\mathrm{IV})$ with \begin{align*}
    (\mathrm{I})&\vcentcolon=-\langle Q( n) \nabla  n , u \rangle_{H^k_x}, \\
    (\mathrm{II}) &\vcentcolon=-\langle K( n) \mathcal{A}u, u \rangle_{H^k_x},   \\
     (\mathrm{III}) &\vcentcolon=-\left\langle K( n)\int (w^\alpha-u) \D \mu(\alpha), u \right\rangle_{H^k_x}, \\ 
     (\mathrm{IV}) &\vcentcolon= -\left\langle K( n)\int r^\alpha (w^\alpha-u) \D \mu(\alpha)   , u \right\rangle_{H^k_x}.
\end{align*}
Arguing as in the proof of Proposition \ref{prop-refinedLWP-VNSincomp}, we can bound the right-hand side of the previous inequality and obtain for all $\sigma>1$ that
\begin{align*}
    &\frac{1}{2}\frac{\mathrm{d}}{\mathrm{d}t} \left(\Vert  u \Vert_{H^k_x}^2+\int \| w^\alpha\|_{H^k_x}^2 \D \mu(\alpha) +\right)
     +\Vert (\nabla, \mathrm{div}) u \Vert_{H^k_x}^2 + \int \Vert u-w^\alpha \Vert^2_{H^k_x} \, \mathrm{d}\mu(\alpha) \\
     & \leq C_k \left(\Vert \nabla w^\alpha \Vert_{L^{\upinfty}_\alpha L^\infty_x }\int  \| w^\alpha\|_{H^k_x}^2 \D \mu(\alpha) + (\Vert \nabla u \Vert_{L^{\upinfty}_x} +\Vert \rr \Vert_{L^{\upinfty}_\alpha L^{\upinfty}_x}) \| u\|_{H^k_x}^2 \right. \\
     &  \qquad \quad   + \Vert \rr \Vert_{L^{\upinfty}_\alpha L^{\upinfty}_x}^2 \int  \Vert w^\alpha-u \Vert_{H^{k}_x}^2 \D \mu(\alpha)\\
     & \qquad \quad   \left. + \sigma \Vert w^\alpha-u \Vert_{L^2_\alpha L^{\upinfty}_x}^2 \sigma^{-1}\int \Vert r^\alpha \Vert_{H^{k-1}_x}^2 \D \mu(\alpha) \right) +S.
\end{align*}

Let us now treat one by one the four terms contained in $S$. For the first one, integrating by parts once, we have the bound
\begin{align*}
    (\mathrm{I}) &\leq \Vert Q( n) \Vert_{L^{\upinfty}_x} \Vert \nabla  n \Vert_{L^2_x}\Vert u \Vert_{L^2_x} + \Vert \nabla u \Vert_{H^k_x} \Vert Q( n) \nabla  n \Vert_{H^{k-1}_x} \\
    & \leq \Vert Q( n) \Vert_{L^{\upinfty}_x} \Vert \nabla  n \Vert_{L^2_x}^2 + \Vert Q( n) \Vert_{L^{\upinfty}_x}\Vert u \Vert_{L^2_x}^2 \\
   & \quad + \eta \Vert \nabla u \Vert_{H^k_x}^2+C_k \eta^{-1}(\Vert Q( n) \Vert_{L^{\upinfty}_x}^2 \Vert  n \Vert_{H^k_x}^2+\Vert Q( n) \Vert_{H^{k-1}_x}^2 \Vert \nabla  n \Vert_{L^{\upinfty}_x}^2),
\end{align*}
for all $\eta>0$, thanks to tame estimates (see Proposition \ref{prop-Sobolev}) and Young inequality.
Since we have $Q(0)=0$, by Taylor expansion we observe that $\vert Q( n) \vert \leq C(\Vert  n \Vert_{L^{\upinfty}_x}) \vert  n \vert$. By the composition rule from Proposition \ref{prop-Sobolev}, we infer
\begin{multline*}
    (\mathrm{I}) \leq  \eta \Vert \nabla u \Vert_{H^k_x}^2 \\
    +  C_k \eta^{-1} \varphi_k(\Vert  n \Vert_{L^{\upinfty}_x})(\Vert  n \Vert_{L^{\upinfty}_x}+\Vert  n \Vert_{L^{\upinfty}_x}^2 + \Vert \nabla  n \Vert_{L^{\upinfty}_x}^2)
    (\Vert u \Vert_{H^k_x}^2 +\Vert  n \Vert_{H^k_x}^2).
\end{multline*}
Next, the term $(\mathrm{II})$ is written as
\begin{align*}
    (\mathrm{II}) &=  -\sum_{|\beta|\le k}\left\langle \partial^\beta\!\bigl(K( n)\mathcal{A}u\bigr),\,\partial^\beta u\right\rangle_{L^2_x}
    = (\mathrm{II})_{\le k-1}+(\mathrm{II})_{k}, \\
    (\mathrm{II})_{\le k-1}
&\vcentcolon= -\sum_{|\beta|\le k-1}\left\langle \partial^\beta\!\bigl(K( n)\mathcal{A}u\bigr),\,\partial^\beta u\right\rangle_{L^2_x}, \\
(\mathrm{II})_{k}
&\vcentcolon= -\sum_{|\beta|=k}\left\langle \partial^\beta\!\bigl(K( n)\mathcal{A}u\bigr),\,\partial^\beta u\right\rangle_{L^2_x}.
\end{align*}
For the first term, we directly get by tame estimates from Proposition \ref{prop-Sobolev} and Young inequality
\begin{align*}
    (\mathrm{II})_{\le k-1} &\leq C_k \Vert u \Vert_{H^{k-1}_x}(\Vert K( n)\Vert_{L^{\upinfty}_x} \Vert \mathcal{A}u\Vert_{H^{k-1}_x}+\Vert K( n) \Vert_{H^{k-1}_x} \Vert \mathcal{A}u\Vert_{L^{\upinfty}_x}), \\
    & \leq \eta \Vert (\nabla, \Div)(u) \Vert_{H^k_x}^2+  C_k \eta^{-1} \varphi_k(\Vert  n \Vert_{L^{\upinfty}_x}) \\
    &\qquad \times \left(\Vert  n \Vert_{L^{\upinfty}_x}^2+ \Vert \nabla^2 u \Vert_{L^{\upinfty}_x}\right) \left(\Vert u \Vert_{H^k_x}^2 +\Vert  n \Vert_{H^k_x}^2\right),
\end{align*}
for all $\eta>0$. For the second term, we decompose  it as
\begin{align*}
    (\mathrm{II})_{k}
&= -\sum_{|\beta|=k}\Big\langle K( n)\,\mathcal{A}\partial^\beta u,\,\partial^\beta u\Big\rangle_{L^2_x}
   -\sum_{|\beta|=k}\Big\langle [\partial^\beta,K( n)]\mathcal{A}u,\,\partial^\beta u\Big\rangle_{L^2_x} \\
   & 
\leq C_k\|K( n)\|_{L^{\upinfty}_x}\|(\nabla,\Div)u\|_{H^k_x}^2 + \eta \|(\nabla,\Div)u\|_{H^k_x}^2+ \eta^{-1}\Vert K'( n) \nabla  n \Vert_{L^{\upinfty}_x}^2\Vert u \Vert_{H^k_x}^2  \\
& \quad +C_k \Big(\Vert K'( n) \nabla  n \Vert_{L^{\upinfty}_x} \Vert \mathcal{A}u \Vert_{H^{k-1}_x}+ \Vert K( n) \Vert_{H^k_x} \Vert \mathcal{A}u \Vert_{L^{\upinfty}_x}\Big)\|u\|_{H^k_x}\\
& \leq (\eta + C_k\| n\|_{L^{\upinfty}_x})\|(\nabla,\Div)u\|_{H^k_x}^2 
+ C_k \eta^{-1}\,\varphi_k(\| n\|_{L^{\upinfty}_x})
\\
&\qquad \times \big(\|\nabla  n\|_{L^{\upinfty}_x}^2+\|\nabla^2 u\|_{L^{\upinfty}_x}\big)\big(\|u\|_{H^k_x}^2+\| n\|_{H^k_x}^2\big),
\end{align*}
thanks to an integration by parts in the first sum, and the use of a commutator estimate for the second sum.

Gathering everything together, we infer the bound
\begin{align*}
    (\mathrm{II}) &\leq (2\eta+C_k \Vert  n\Vert_{L^{\upinfty}_x})\Vert (\nabla, \Div)(u) \Vert_{H^k_x}^2 +\eta^{-1}C_k \varphi_k(\Vert  n \Vert_{L^{\upinfty}_x}) (\Vert ( n, \nabla  n) \Vert_{L^{\upinfty}_x}\\
    &\quad + \Vert \nabla^2 u \Vert_{L^{\upinfty}_x})(\Vert u \Vert_{H^k_x}^2 +\Vert  n \Vert_{H^k_x}^2),
\end{align*}
for all $\eta>0$. For $(\mathrm{III})$ and $(\mathrm{IV})$, we integrate by parts in $(\mathrm{IV})$ and then obtain thanks to tame estimates from Proposition \ref{prop-Sobolev} that
\begin{align*}
    &(\mathrm{III})+(\mathrm{IV}) \\
    &\leq \Vert u \Vert_{H^k_x} \int \Vert K( n) (w^\alpha-u) \Vert_{H^k_x} \, \mathrm{d}\mu(\alpha) \\
    &\quad + \Vert \nabla u \Vert_{H^k_x} \int \left\Vert K( n)r^\alpha (w^\alpha-u) \right\Vert_{H^{k-1}_x}\, \mathrm{d}\mu(\alpha) \\
    & \quad + \Vert \rr \Vert_{L^{\upinfty}_\alpha L^{\upinfty}_x} \Vert \ww-u \Vert_{L^1_\alpha L^{\upinfty}_x} \Vert u \Vert_{L^2_x} \Vert K( n) \Vert_{L^2_x} \\ 
    & \leq C_k\Vert u \Vert_{H^k_x} \int (\Vert K( n)\Vert_{L^{\upinfty}_x} \Vert (w^\alpha-u) \Vert_{H^k_x}+\Vert K( n)\Vert_{H^k_x} \Vert (w^\alpha-u) \Vert_{L^{\upinfty}_x}) \, \mathrm{d}\mu(\alpha)\\
    &
    \quad + C_k\Vert \nabla u \Vert_{H^k_x}\Bigg(\int \left\Vert K( n) \right\Vert_{L^{\upinfty}_x} \Vert r^\alpha \Vert_{L^{\upinfty}_x}  \Vert w^\alpha-u \Vert_{H^{k-1}_x}    
    \\
    & \qquad\qquad\qquad \quad    + \left\Vert K( n) \right\Vert_{L^{\upinfty}_x} \Vert r^\alpha \Vert_{H^{k-1}_x}  \Vert w^\alpha-u \Vert_{L^{\upinfty}_x}  \\
    & \qquad\qquad\qquad \quad + \Vert K( n) \Vert_{H^{k-1}_x} \Vert r^\alpha \Vert_{L^{\upinfty}_x}  \Vert w^\alpha-u \Vert_{L^{\upinfty}_x} \, \mathrm{d}\mu(\alpha) \Bigg) \\
    & \quad + \Vert \rr \Vert_{L^{\upinfty}_\alpha L^{\upinfty}_x} \Vert \ww-u \Vert_{L^2_\alpha L^{\upinfty}_x} \Vert u \Vert_{L^2_x} \Vert K( n) \Vert_{L^2_x}.
\end{align*}
By using the composition rule from Proposition \ref{prop-Sobolev} and rearranging all the terms, we get for all $\eta>0$ and $\sigma>1$
\begin{align*}
    (\mathrm{III})&+(\mathrm{IV}) 
    \leq \eta \Vert \nabla u \Vert_{H^k_x}^2  \\
    &\quad +C_\eta \varphi_k(\Vert  n \Vert_{L^{\upinfty}_x}) \Big(\Vert  n \Vert_{L^{\upinfty}_x}+\Vert \rr \Vert_{L^{\upinfty}_\alpha L^{\upinfty}_x}^2+\Vert  n \Vert_{L^{\upinfty}_x}^2\Vert \rr \Vert_{L^{\upinfty}_\alpha L^{\upinfty}_x}^2 \\
    &\quad \quad+(\Vert \rr \Vert_{L^{\upinfty}_\alpha L^{\upinfty}_x}^2+\sigma \Vert  n \Vert_{L^{\upinfty}_x}^2)\Vert \ww-u  \Vert_{L^2_\alpha L^{\upinfty}_x}^2 + \Vert \rr \Vert_{L^{\upinfty}_\alpha L^{\upinfty}_x}\Vert \ww-u  \Vert_{L^2_\alpha L^{\upinfty}_x} \Big) \\
    & \qquad  \qquad  \qquad \times   \left(\sigma^{-1} \int \Vert r^\alpha \Vert_{H^{k-1}_x}^2+ \int\Vert w^\alpha  \Vert_{H^k_x}^2+\Vert u \Vert_{H^k_x}^2+\Vert  n \Vert_{H^k_x}^2 \right).
\end{align*}
Adding all the terms together, we end up for all $\eta>0$ and $\sigma>1$ with
\begin{equation}
    \label{eq:maincomp}
\begin{aligned}
    &\frac{1}{2}\frac{\mathrm{d}}{\mathrm{d}t} \left(\int \| w^\alpha\|_{H^k_x}^2 \D \mu(\alpha)+\Vert u \Vert^2_{H^k_x} \right) \\
    &\qquad +\left(1-4\eta-C_k\Vert  n \Vert_{L^\infty_{T,x}} \right)\Vert (\nabla, \mathrm{div}) u \Vert_{H^k_x}^2 
 + \langle \nabla  n, u \rangle_{H^k_x}\\
    &\leq  \eta^{-1} \varphi_k(\Vert  n \Vert_{L^{\upinfty}_x}) \\
    &\quad \times \bigg(\Vert  n \Vert_{L^{\upinfty}_x}+\Vert ( n, \nabla  n) \Vert_{L^{\upinfty}_x}^2+\Vert \rr \Vert_{L^{\upinfty}_\alpha L^{\upinfty}_x}+\Vert \rr \Vert_{L^{\upinfty}_\alpha L^{\upinfty}_x}^2+\Vert  n \Vert_{L^{\upinfty}_x}^2\Vert \rr \Vert_{L^{\upinfty}_\alpha L^{\upinfty}_x}^2 \\
    & \qquad  \quad   +(\Vert \rr \Vert_{L^{\upinfty}_\alpha L^{\upinfty}_x}^2+\sigma \Vert  n \Vert_{L^{\upinfty}_x}^2)\Vert \ww-u  \Vert_{L^2_\alpha L^{\upinfty}_x}^2 + \Vert \rr \Vert_{L^{\upinfty}_\alpha L^{\upinfty}_x}\Vert \ww-u  \Vert_{L^2_\alpha L^{\upinfty}_x} \\
    & \qquad \qquad   \quad   +\Vert \ww-u  \Vert_{L^1_\alpha L^{\upinfty}_x}^2+\Vert \nabla \ww \Vert_{L^{\upinfty}_\alpha L^{\upinfty}_x}+\Vert (\nabla u, \nabla^2 u) \Vert_{L^{\upinfty}_x} \bigg)y_{k, \sigma}.
\end{aligned}
\end{equation}

$\bullet$ Since $u$ is not necessarily divergence-free, the equation on $\mathrm{div}(w^\alpha)$ becomes
\begin{align*}
    \partial_t \mathrm{div}(w^\alpha)+\mathrm{div}(w^\alpha)=\mathrm{div}(u)-\mathrm{div}\left( (\mathcal{V}+w^\alpha )\cdot \nabla w^\alpha \right),
\end{align*}
Proceeding as in the proof of Proposition \ref{prop-refinedLWP-VNSincomp}, we get
\begin{multline*}
    \frac{1}{2}\frac{\mathrm{d}}{\mathrm{d}t} \Vert \mathrm{div}(w^\alpha) \Vert_{H^{k-1}_x}^2+ \Vert \mathrm{div}(w^\alpha)\Vert_{H^{k-1}_x}^2  \leq \Vert \mathrm{div}(u) \Vert_{H^{k-1}_x} \Vert \mathrm{div}(w^\alpha) \Vert_{H^{k-1}_x}  \\
    +C_k \Vert (\nabla \ww, \nabla^2 \ww) \Vert_{L^{\upinfty}_\alpha L^{\upinfty}_x}(\Vert w^\alpha \Vert_{H^k_x}^2+\Vert \mathrm{div}(w^\alpha) \Vert_{H^{k-1}_x}^2),
\end{multline*}
which, by Young inequality and an absorption of the term $\Vert \mathrm{div}(w^\alpha) \Vert_{H^{k-1}_x} $ in the left-hand side, yields:
\begin{multline*}
    \frac{1}{2}\frac{\mathrm{d}}{\mathrm{d}t} \Vert \mathrm{div}(w^\alpha) \Vert_{H^{k-1}_x}^2+ \frac{1}{2}\Vert \mathrm{div}(w^\alpha)\Vert_{H^{k-1}_x}^2 \\ \leq \frac{1}{2}\Vert \mathrm{div}(u) \Vert_{H^{k-1}_x}^2+C_k \Vert (\nabla \ww, \nabla^2 \ww) \Vert_{L^{\upinfty}_\alpha L^{\upinfty}_x}(\Vert w^\alpha \Vert_{H^k_x}^2+\Vert \mathrm{div}(w^\alpha) \Vert_{H^{k-1}_x}^2).
\end{multline*}
Integrating in $\alpha$ and  summing up this inequality to~\eqref{eq:maincomp}, and absorbing the term $\Vert \mathrm{div}(u) \Vert_{H^{k-1}_x} $ in the left-hand side, we now obtain 
\begin{equation}
    \label{eq:maincomp2}
\begin{aligned}
    &\frac{1}{2}\frac{\mathrm{d}}{\mathrm{d}t} \left(\int \| w^\alpha\|_{H^k_x}^2 \D \mu(\alpha)+\int \| \mathrm{div}(w^\alpha)\|_{H^{k-1}_x}^2 \D \mu(\alpha)+\Vert u \Vert^2_{H^k_x} \right) \\
    & \quad +\int \| \mathrm{div}(w^\alpha)\|_{H^{k-1}_x}^2 \D \mu(\alpha) \\
    & \quad +\left(1-\eta-C_k\Vert  n \Vert_{L^\infty_{T,x}} \right)\Vert (\nabla, \mathrm{div}) u \Vert_{H^k_x}^2 
 + \langle \nabla  n, u \rangle_{H^k_x}\\
    &\leq C_\eta \varphi_k(\Vert  n \Vert_{L^{\upinfty}_x}) \\
    &\quad \times \bigg(\Vert  n \Vert_{L^{\upinfty}_x}+\Vert ( n, \nabla  n) \Vert_{L^{\upinfty}_x}^2+\Vert \rr \Vert_{L^{\upinfty}_\alpha L^{\upinfty}_x}+\Vert \rr \Vert_{L^{\upinfty}_\alpha L^{\upinfty}_x}^2+\Vert  n \Vert_{L^{\upinfty}_x}^2\Vert \rr \Vert_{L^{\upinfty}_\alpha L^{\upinfty}_x}^2 \\
    & \qquad \quad   +(\Vert \rr \Vert_{L^{\upinfty}_\alpha L^{\upinfty}_x}^2+\sigma \Vert  n \Vert_{L^{\upinfty}_x}^2)\Vert \ww-u  \Vert_{L^2_\alpha L^{\upinfty}_x}^2 + \Vert \rr \Vert_{L^{\upinfty}_\alpha L^{\upinfty}_x}\Vert \ww-u  \Vert_{L^2_\alpha L^{\upinfty}_x} \\
    & \qquad \qquad  \quad   +\Vert \ww-u  \Vert_{L^1_\alpha L^{\upinfty}_x}^2+\Vert (\nabla \ww, \nabla^2 \ww) \Vert_{L^{\upinfty}_\alpha L^{\upinfty}_x}+\Vert (\nabla u, \nabla^2 u) \Vert_{L^{\upinfty}_x} \bigg)y_{k, \sigma}.
\end{aligned}
\end{equation}

$\bullet$ We now perform an estimate of $r^\alpha$ in $H^{k-1}_x$: as in the proof of Proposition \ref{prop-refinedLWP-VNSincomp}, we get for all $\sigma>1$
\begin{align*}
    &\frac{1}{2} \frac{\mathrm{d}}{\mathrm{d}t} (\sigma^{-1}\Vert r^\alpha \Vert^2_{H^{k-1}_x}) 
\\ &\qquad \leq \left(\sigma^{-1}+C_k\| \nabla w^\alpha \|_{L^{\upinfty}_x}   + C_k 
\| r^\alpha\|_{L^{\upinfty}_x}  \right)(\| w^\alpha\|_{H^k_x}^2+\sigma^{-1}\Vert r^\alpha \Vert^2_{H^{k-1}_x}) \\
& \qquad \quad + \frac{1}{2} \Vert \mathrm{div}(w^\alpha)\Vert_{H^{k-1}_x}^2.
 \end{align*}
Integrating in $\alpha$ then entails
\begin{align*}
    &\frac{1}{2} \frac{\mathrm{d}}{\mathrm{d}t} \sigma^{-1}\int \Vert r^\alpha \Vert^2_{H^{k-1}_x} \D\mu(\alpha)
\\ 
& \leq \left( \sigma^{-1}+C_k\| \nabla \ww \|_{L^{\upinfty}_\alpha L^{\upinfty}_x}   + C_k 
\| \rr\|_{L^{\upinfty}_\alpha L^{\upinfty}_x}  \right) \\
&\qquad \times\left(\int \| w^\alpha\|_{H^k_x}^2 \D\mu(\alpha) +\sigma^{-1}\int \Vert r^\alpha \Vert^2_{H^{k-1}_x}) \D\mu(\alpha \right) \\
& \quad + \frac{1}{2} \int \Vert \mathrm{div}(w^\alpha)\Vert_{H^{k-1}_x}^2 \D\mu(\alpha).
 \end{align*}
Similarly, for the perturbation of the fluid density $ n$, we have the following estimate in $H^k_x$, thanks to Proposition \ref{prop-Sobolev}:
\begin{align*}
&\frac{1}{2}\frac{\mathrm{d}}{\mathrm{d}t} \Vert  n \Vert^2_{H^k_x}+\langle   n, \mathrm{div} u\rangle_{H^k_x} \\
&\leq C_k \left( \| \nabla u \|_{L^{\upinfty}_x} \|  n\|_{H^k_x}^2  + 
\| \nabla u\|_{H^k_x} \|  n\|_{L^{\upinfty}_x}\| n\|_{H^k_x} \right)  \\ 
&\leq C_k \left( \| \nabla u \|_{L^{\upinfty}_x} \|  n\|_{H^k_x}^2  + 
(\eta')^{-1}\|  n\|_{L^{\upinfty}_x}^2 \| n\|_{H^k_x}^2 \right) +\eta' \| \nabla u \|_{H^k_x}^2,
\end{align*}
for all $\eta'>0$. Summing up these two estimates to~\eqref{eq:maincomp2}, we infer
\begin{align*}
    &\frac{1}{2}\frac{\mathrm{d}}{\mathrm{d}t} y_{k, \sigma}+\left(1-4\eta -\eta'-C_k\Vert  n \Vert_{L^\infty_{T,x}} \right)\Vert (\nabla, \mathrm{div}) u \Vert_{H^k_x}^2 \\
    &\leq \eta^{-1} \varphi_k(\Vert  n \Vert_{L^{\upinfty}_x}) \\
    &\quad \times \bigg((1+(\eta')^{-1})\Vert  n \Vert_{L^{\upinfty}_x}+\Vert ( n, \nabla  n) \Vert_{L^{\upinfty}_x}^2+(1+\Vert \rr \Vert_{L^{\upinfty}_\alpha L^{\upinfty}_x})\Vert \rr \Vert_{L^{\upinfty}_\alpha L^{\upinfty}_x}  \\
    & \qquad \quad +\Vert  n \Vert_{L^{\upinfty}_x}^2\Vert \rr \Vert_{L^{\upinfty}_\alpha L^{\upinfty}_x}^2  +(\Vert \rr \Vert_{L^{\upinfty}_\alpha L^{\upinfty}_x}^2+\sigma \Vert  n \Vert_{L^{\upinfty}_x}^2)\Vert \ww-u  \Vert_{L^2_\alpha L^{\upinfty}_x}^2  \\
    & \qquad   \quad + \Vert \rr \Vert_{L^{\upinfty}_\alpha L^{\upinfty}_x}\Vert \ww-u  \Vert_{L^2_\alpha L^{\upinfty}_x}   +\Vert \ww-u  \Vert_{L^1_\alpha L^{\upinfty}_x}^2 \\
    & \qquad \quad +\Vert (\nabla \ww, \nabla^2 \ww) \Vert_{L^{\upinfty}_\alpha L^{\upinfty}_x}+\Vert (\nabla u, \nabla^2 u) \Vert_{L^{\upinfty}_x} +\sigma^{-1}\bigg)y_{k, \sigma},
\end{align*}
that holds for all $\sigma>1, \eta,\eta'>0$. Here, we have crucially used the fact that the 
terms $\langle \nabla  n, u \rangle_{H^k_x}$ and $\langle  n, \mathrm{div}(u) \rangle_{H^k_x}$ exactly cancel out. 

It remains to explain how we can drop out the dissipative term in the left-hand side, so that we can then deduce a suitable estimate on $ \| u(t)\|_{H^k_x}^2+\|  n(t)\|_{H^k_x}^2$. Thanks to the assumption \eqref{bound-technical-theta-VNScompressible}, we have $C_k\Vert  n \Vert_{L^\infty_{T,x}} \leq \frac{1}{2}$, hence choosing $\eta=\eta'=1/20$, we can ensure that
\begin{align*}
    1-\eta -\eta'-C_k\Vert  n \Vert_{L^\infty_{T,x}} >0.
\end{align*}
We therefore deduce the desired estimate for $ \| u(t)\|_{H^k_x}^2+\|  n(t)\|_{H^k_x}^2$. The estimate for $\ww$ can finally be obtained exactly as in the proof of Proposition \ref{prop-refinedLWP-VNSincomp}, since it satisfies the same equation in the incompressible and the compressible case. This concludes the proof of Proposition \ref{prop-refinedLWP-VNScomp}.
\end{proof}

\subsection{Bootstrap}
We now set up a bootstrap argument similar to the one introduced in Section \ref{sectionVNS-Stratbootstrap} for the incompressible case.  

\medskip

Let $k \in \N$ such that $k>3+3/2$ and $p \in [2, \infty]$. For a given initial data $(\rr_0, \ww_0,  n_0, u_0)$, let $T^{\mathrm{max}}_{0}>0$ be the maximal time of existence of the unique strong solution $(\rr,\ww,  n, u)$ to the system \eqref{eq:VNScomp-pertubLWP}, obtained from Theorem \ref{thm:appliVNS-comp}, and such that for any $T \in (0, T^{\mathrm{max}}_{0})$, we have
\begin{align*}
    &\rr \in L^\infty([0,T]; L^\infty_\alpha (H^{k-1}_x) ), \ \ \ww \in L^\infty([0,T]; \mathcal{H}^{k,p}), \ \ n \in L^\infty([0,T]; H^{k}_x), \\ 
    &u \in L^\infty([0,T]; H^{k}_x) \cap L^2([0,T]; H^{k+1}_x).
\end{align*}

Recall that in view of Remark \ref{rem-normalizationcompVNS-longtime}, we can always assume that $\langle n_0 \rangle=0$ (enforcing $\varrho_0 =\langle \varrho_0 \rangle+n_0$ since the beginning) and therefore
\begin{align*}
    \langle n(t) \rangle=0, \ \ t \in [0,T],
\end{align*}
since this condition is preserved by the evolution.

As a direct application from the refined energy estimates from Proposition \ref{prop-refinedLWP-VNScomp}, we also know that a blow-up criterion for the evolution of the perturbation is 
\begin{align}\label{eq:blowup-compVNS}
    T_0^{\max}<+\infty \Longleftrightarrow \int_0^{ T_0^{\max}} \mathcal{Q}_{k, \sigma}\Big[ n(s), u(s), \rr(s),\ww(s) \Big] \, \mathrm{d}s =+\infty.
\end{align}
where the quantity $\mathcal{Q}_{k, \sigma}\Big[ n(s), u(s), \rr(s),\ww(s) \Big]$ has been defined in \eqref{eq:fct-blowup-comp}. 

Denote the exponential decay rate from Proposition \ref{Prop-lineardecayVNS-comp} as $\Lambda=\min(\Lambda_i, \Lambda_c) \in (0,1)$ and let $\lambda \in (0, \Lambda)$ and  $\lambda' \in (0, \theta \lambda)$ for some $\theta \in (0,1)$ to be determined later.
We now set up the following bootstrap procedure, which is a variant of that used for the incompressible case in Section \ref{sectionVNS-Stratbootstrap}.

For $\delta \in (0,1)$ small enough to be determined later on, we define 
\begin{align}\label{assumpExtra-bootstrapCOMP}
\begin{split}
&T^\star_{\delta} \vcentcolon= \sup\Bigg\{ T \in (0, T^{\mathrm{max}}_{0}),    \, \forall t \in [0,T], \\
  &  \ \|(\nabla \ww,\nabla^2 \ww)(t)\|_{L^{\upinfty}_\alpha L^{\upinfty}_x} + \|(\nabla u,\nabla^2 u, \nabla ^3 u)(t) \|_{L^{\upinfty}_x}+\|(\nabla  n,\nabla^2  n)(t) \|_{L^{\upinfty}_x}  \leq \delta e^{-\lambda' t}, \\
  &  \quad \qquad \qquad\qquad \qquad \qquad \qquad  \|\nabla \ww(t)\|_{L^{\upinfty}_\alpha L^{2}_x} + \|\nabla u(t) \|_{L^{2}_x} +\|\nabla  n(t) \|_{L^{2}_x} \leq \delta e^{-\lambda t}, \\
 &
 \qquad  \qquad  \qquad   \qquad\qquad \qquad  \qquad \qquad \qquad \qquad \qquad \quad \Vert  \langle  \ww(t)- u(t) \rangle \Vert_{L^p_\alpha} \leq \delta e^{-\lambda t}
\Bigg\},
\end{split}
\end{align}
For small enough initial data with respect to $\delta$, the local well-posedness theory from Theorem~\ref{thm:appliVNS-comp} ensures that $T_\delta^\star>0$. Compared to the bootstrap assumption \eqref{assump-bootstrapVNS}, notice that we have added additional decay estimates for the derivatives of the fluid density and velocity.   Note that these assumptions actually require $k>3+3/2$ (to be in agreement with Sobolev embedding), which explains why the regularity threshold is higher in the compressible than in the incompressible case.

Furthermore,  thanks to Poincaré inequality,  we have for all $t \in (0, T^\star_\delta)$,
\begin{align}\label{estimate-theta-decay-comp}
     \|  n(t)-\langle n(t) \rangle \|_{L^{2}_x}  \leq \delta e^{-\lambda t}, \ \ \|  n(t)-\langle n(t) \rangle \|_{L^{\upinfty}_x}  \leq \delta e^{-\lambda' t}.
\end{align}
In particular, since $\langle  n(t) \rangle =\langle n_0 \rangle$, we get
\begin{align}\label{eq:ronflex}
     \|  n(t)  \|_{L^{\upinfty}_x}  \leq  \vert \langle n_0 \rangle \vert + \delta= \vert \langle \varrho_0-1 \rangle \vert + \delta,
\end{align}
and by choosing $\delta$ small enough combined with the smallness hypothesis \eqref{1st-smallness-cond-comp} (that implies $ \vert \langle \varrho_0-1 \rangle\vert\leq \eps_0$), the Assumptions \eqref{bound-technical-theta-VNScompressible} and  \eqref{assumption-novoid-theta} are satisfied up to time $T^\star_\delta$.

\medskip

\noindent \textbf{Blow-up criterion}. First let us show that, provided that $T_0^{\max}=T_\delta^\star$, the estimates from \eqref{assumpExtra-bootstrapCOMP} are enough to prove the control 
$$\int_0^{ T_0^{\max}} \mathcal{Q}_{k, \sigma}\Big[ n(s), u(s), \rr(s),\ww(s) \Big] \, \mathrm{d}s<\infty.$$
Thanks to \eqref{eq:blowup-compVNS}, it will entail that $T_0^{\max}=T_\delta^\star=+\infty$. By inspection on \eqref{eq:fct-blowup-comp}, we observe that the only contributions that have to be understood are the ones from the term $\Vert \ww-u  \Vert_{L^2_\alpha L^{\upinfty}_x}$, $\Vert \ww-u  \Vert_{L^2_\alpha L^{\upinfty}_x}^2$ and $\Vert \ww-u  \Vert_{L^1_\alpha L^{\upinfty}_x}^2$. By Hölder inequality (in $\alpha)$, it is enough to prove that $\Vert \ww-u  \Vert_{L^2_\alpha L^{\upinfty}_x}$ decays exponentially fast in time. Since $p \geq 2$, we rely on Poincaré inequality and  estimate this term as
\begin{align*}
    \Vert \ww(t)-u(t)  \Vert_{L^2_\alpha L^{\upinfty}_x} \lesssim \Vert  \langle  \ww(t)- u(t) \rangle \Vert_{L^p_\alpha}+\Vert \nabla(\ww(t)-u(t))  \Vert_{L^{\upinfty}_\alpha L^{\upinfty}_x},
\end{align*}
from which we infer the desired exponential decay in time thanks to \eqref{assumpExtra-bootstrapCOMP}.

\medskip

\noindent \textbf{Improved estimates}. It remains to prove that $T_0^{\max}=T_\delta^\star$. As in the incompressible case, a continuity argument ensures that it is sufficient to prove that the estimates  from \eqref{assumpExtra-bootstrapCOMP} hold replacing $\delta$ by $\delta/2$. Working on $[0 ,T^\star_\delta)$, we now rewrite \eqref{eq:VNScomp-pertub1} as 
\begin{equation}
\label{eq:VNScomp-pertub2}
\left\{  
\begin{aligned}
\partial_t r^\alpha + \mathcal{V} &\cdot \nabla r^\alpha + \Div (w^\alpha)  + \Div (r^\alpha w^\alpha)= 0,\\
\partial_t w^\alpha + (\mathcal{W} &\cdot \nabla) w^\alpha   =  u- w^\alpha - ((w^\alpha-\mathcal{Z}) \cdot \nabla) w^\alpha
,
\\
\partial_t  n + \mathcal{W} &\cdot \nabla  n +  \Div(u)  = -   n \mathrm{div}(u) -(u-\mathcal{Z}) \cdot \nabla  n, \\
\partial_t u + (\mathcal{W} \cdot \nabla) u  &+\nabla  n - \mathcal{A} u \\
&=
\int (w^\alpha-u) \D \mu(\alpha)  -K( n) \mathcal{A} u- Q( n)\nabla  n   - ((u-\mathcal{Z}) \cdot \nabla) u \\
& \quad -K( n)\int (w^\alpha-u) \D \mu(\alpha) +\frac{1}{1+ n}\int r^\alpha (w^\alpha-u) \D \mu(\alpha),
\end{aligned}
\right.
\end{equation}
where $$K( n)=\frac{ n}{1+ n}, \ \ Q( n)=\frac{P'(1+ n)-(1+ n)}{1+ n}.$$
Since $K(0)=Q(0)=0$, we have
\begin{align}\label{propfunctionK1}
    \vert  n \vert <\delta<\frac{1}{2} \Longrightarrow   
    \left\{\begin{array}{ll}
    1+  n \geq \frac{1}{2}, \\[1mm]
    \vert K( n) \vert \lesssim \vert  n \vert, \\[1mm]
    \vert Q( n) \vert \leq C_\delta \vert  n \vert, \ \ \vert Q'( n)\vert \leq C_\delta,
    \end{array}\right.
\end{align}
thanks to Taylor formula and direct bounds. Here, the constant $C_\delta$ depends in a continuous and increasing way on $\delta$.

\medskip

Let us label the nonlinear source terms appearing in the right-hand side of \eqref{eq:VNScomp-pertub2}:

$\bullet$ for the equation on $w^\alpha$:
$$S_1^\alpha \vcentcolon=- ((w^\alpha-\mathcal{Z}) \cdot \nabla) w^\alpha.$$

$\bullet$ for the equation on $u$:
$$S_2\vcentcolon=S_{2,a}+S_{2,b}+S_{2,c}+S_{2,d}+S_{2,e},$$
where
\begin{align*}
    S_{2,a}&\vcentcolon=-K( n) \mathcal{A} u, \ \ S_{2,b}\vcentcolon= - Q( n)\nabla  n,  \ \ 
      S_{2,c}\vcentcolon=- ((u-\mathcal{Z}) \cdot \nabla) u ,\\
       S_{2,d}&\vcentcolon= -K( n)\int (w^\alpha-u) \D \mu(\alpha), \ \  S_{2,e}\vcentcolon= \frac{1}{(1+ n)}\int r^\alpha (w^\alpha-u) \D \mu(\alpha).
\end{align*}

$\bullet$ for the equation on $ n$:
$$S_{3}\vcentcolon=S_{3,a}+S_{3,b},$$
where
\begin{align*} 
    S_{3,a}&\vcentcolon=-  n \Div (u), \quad  S_{3,b}\vcentcolon=-(u-\mathcal{Z}) \cdot \nabla  n.
\end{align*}

Note that several terms display a derivative bearing on $ n$ and $u$ but this will be taken into account in the forthcoming bootstrap procedure.

\medskip

First, let us now explain how one can proceed to obtain the analogues of Lemmas \ref{lem:ImproveLinL2-nabla-u}--\ref{lem:ImproveLinL2-nabla-walpha}--\ref{lem:ImproveLinL2-diffnabla}--\ref{lem:ImproveLdiffMoyenne} in the compressible case, combining the bootstrap assumption \eqref{assumpExtra-bootstrapCOMP} with the linear decay from Proposition \ref{Prop-lineardecayVNS-comp}. We recall that the latter entails the following exponential decay:  
\begin{multline*}
    \| \nabla  n(t) \|_{L^2_x}+\| \nabla u(t) \|_{L^2_x} \\
  \lesssim \eps_0 e^{-\Lambda t}  + \int_0^t e^{-\Lambda_c(t-s)} \left (\|  \nabla  S_1^\alpha  (s)  \|_{L^p_\alpha L^2_x}+\| \nabla S_2(s) \|_{L^2_x}+\| \nabla S_3(s) \|_{L^2_x}   \right)\D s.
\end{multline*}
Let us briefly explain how the source terms $S_1^\alpha$, $S_2$ and $S_3$ can be handled. Until the end of the proof,  $C_\delta>0$ will refer to a constant depending in a continuous increasing way on $\delta$ and that may vary from line to line. 
\begin{itemize}
\item \underline{Terms $S_1^\alpha$ and $S_{2,c}$}: we can treat these terms exactly as in the incompressible case.
\item \underline{Terms $S_{2,d}$ and $S_{2,e}$}: compared to the incompressible case, the new factor involving $ n$ is harmless by using the bound \eqref{propfunctionK1}, since $ n$ is  small enough in $L^\infty_x$, see \eqref{eq:ronflex}.
    \item \underline{Term $S_{2,a}$}:  we have
    \begin{align*}
        \Vert \nabla S_{2,a} \Vert_{L^2_x} 
        &=\Vert \nabla(K( n) \mathcal{A} u) \Vert_{L^2_x} \\
        &\lesssim \Vert K(  n) \Vert_{L^2_x}\Vert \nabla \mathcal{A}u \Vert_{L^{\upinfty}_x}+ \Vert \nabla K( n) \Vert_{L^2_x} \Vert \mathcal{A}u \Vert_{L^{\upinfty}_x}.
        \end{align*}
        Using the composition rule from Proposition \ref{prop-Sobolev} followed by \eqref{propfunctionK1}, we obtain
        \begin{align*}
   \Vert \nabla S_{2,a} \Vert_{L^2_x} &\lesssim  \Vert K(  n) \Vert_{L^2_x}\Vert \nabla^3 u \Vert_{L^{\upinfty}_x}+ C(\Vert  n \Vert_{L^{\upinfty}_x})\Vert \nabla  n \Vert_{L^2_x} \Vert \mathcal{A}u \Vert_{L^{\upinfty}_x}      \\
   & \lesssim  C_{\delta}\left(\Vert  n \Vert_{L^2_x}\Vert \nabla^3 u \Vert_{L^{\upinfty}_x}+ \Vert \nabla  n \Vert_{L^2_x} \Vert \mathcal{A}u \Vert_{L^{\upinfty}_x}\right),
   \end{align*}
hence for all $s \in (0,T_\delta^\star)$,
   \begin{align*}
   \Vert \nabla S_{2,a}(s) \Vert_{L^2_x} \leq C_{\delta} \delta^2 e^{-(\lambda+\lambda')s},
        \end{align*} 
in view of the bootstrap assumption \eqref{assumpExtra-bootstrapCOMP}. 
    \item \underline{Term $S_{2,b}$}: 
As for the term $S_{2,a}$, we get
\begin{align*}
    \Vert \nabla S_{2,b} \Vert_{L^2_x}&=\Vert \nabla (Q( n)\nabla  n) \Vert_{L^2_x} \\
    &\lesssim \Vert Q( n) \Vert_{L^2_x} \Vert \nabla^2  n \Vert_{L^{\upinfty}_x}+\Vert \nabla Q( n) \Vert_{L^2_x} \Vert \nabla  n \Vert_{L^{\upinfty}_x} \\
    &\lesssim C_\delta \Vert  n \Vert_{L^2_x} \Vert \nabla^2  n \Vert_{L^{\upinfty}_x}+C_\delta \Vert \nabla  n \Vert_{L^2_x} \Vert \nabla  n \Vert_{L^{\upinfty}_x},
\end{align*}
therefore the bootstrap assumption \eqref{assumpExtra-bootstrapCOMP} yields for all $s \in (0,T_\delta^\star)$
   \begin{align*}
   \Vert \nabla S_{2,b}(s) \Vert_{L^2_x} \leq C_{\delta} \delta^2 e^{-(\lambda+\lambda')s}.
        \end{align*} 
    \item \underline{Term $S_{3,a}$}: we have
    \begin{align*}
        \Vert \nabla S_{3,a} \Vert_{L^2_x} 
        &=\Vert \nabla ( n \Div (u)) \Vert_{L^2_x} 
        \lesssim \Vert \nabla  n \Vert_{L^2_x}\Vert \nabla u \Vert_{L^{\upinfty}_x} + \Vert  n \Vert_{L^2_x} \Vert \nabla^2 u \Vert_{L^{\upinfty}_x},
    \end{align*}
    so as before we get for all $s \in (0,T_\delta^\star)$
    \begin{align*}
        \Vert \nabla S_{3,a}(s) \Vert_{L^2_x} \lesssim \delta^2 e^{-(\lambda+\lambda')s}.
    \end{align*}
\item \underline{Term $S_{3,b}$}: we have
    \begin{align*}
         \Vert \nabla S_{3,b} \Vert_{L^2_x} 
        &=\Vert \nabla \left((u-\mathcal{Z}) \cdot \nabla  n \right) \Vert_{L^2_x} \\
        & \lesssim \Vert \nabla  n \Vert_{L^2_x}\Vert \nabla u \Vert_{L^{\upinfty}_x} + \left(\Vert u-\langle u \rangle \Vert_{L^2_x}+\vert \langle u-\mathcal{Z} \rangle \vert  \right)  \Vert \nabla^2  n \Vert_{L^{\upinfty}_x} \\
        & \lesssim \Vert \nabla  n \Vert_{L^2_x}\Vert \nabla u \Vert_{L^{\upinfty}_x} + \left(\Vert \nabla u \Vert_{L^2_x}+\vert \langle u-\mathcal{Z} \rangle \vert \right)  \Vert \nabla^2  n \Vert_{L^{\upinfty}_x}.
         \end{align*}
         Invoking \eqref{amelio-conservationlaw} that also holds in the compressible case thanks to Lemma \ref{smallLemma-moyenne}, we obtain for all $s \in (0,T_\delta^\star)$
         \begin{align*}
             \Vert \nabla S_{3,b}(s) \Vert_{L^2_x} \lesssim  (1+C_{0,\delta})\delta^2 e^{-(\lambda+\lambda') s}, 
         \end{align*}
thanks to the bootstrap assumption \eqref{assumpExtra-bootstrapCOMP}. 
\end{itemize}
These estimates are enough to obtain the analogues of Lemmas \ref{lem:ImproveLinL2-nabla-u}--\ref{lem:ImproveLinL2-nabla-walpha}--\ref{lem:ImproveLinL2-diffnabla}--\ref{lem:ImproveLdiffMoyenne} in the compressible case. We gather a summarized statement in the following lemma. 
\begin{Lem}\label{lem:ImproveLinL2-COMPRESSIBLE}
There exists $C_{0,\lambda',\delta}>0$ such that all $t \in [0, T^\star_\delta]$, we have
\begin{align*}
\| \nabla  n(t) \|_{L^2_x} +\| \nabla u(t) \|_{L^2_x}  &+ \|  \nabla \ww \|_{L^{\upinfty}_\alpha L^2_x} \\
&\lesssim \left((1+C_{0, \lambda',\delta})\eps_0+(1+C_{0, \lambda',\delta})\delta^2 \right) e^{-\lambda t}, \\
\| \nabla  \ww(t)-\nabla u(t) \|_{L^{\upinfty}_\alpha L^2_x}    &\lesssim \left((1+C_{0, \lambda',\delta})\eps_0+(1+C_{0, \lambda',\delta})\delta^2 \right) e^{-\lambda t}, \\
\Vert  \langle  \ww(t)- u(t) \rangle \Vert_{L^p_\alpha}    &\lesssim \left((1+C_{0, \lambda',\delta})\eps_0+(1+C_{0, \lambda',\delta})\delta^2 \right) e^{- \lambda t},
\end{align*}
where $C_{0, \lambda',\delta}$ is a continuous and monotone function with respect to each of the parameters.

\end{Lem}

We are finally in position to conclude the bootstrap argument, that is improving the bounds appearing in \eqref{assumpExtra-bootstrapCOMP}. Similarly to Lemma \ref{lem-bootstrap-amelio} for the incompressible case, we have the following statement.

\begin{Lem}\label{lem-bootstrap-amelioCOMPRESSIBLE} 
    There exist $\delta>0$ and $\eps_0>0$ such that for all $k \geq 3+3/2$ and $\eps \in (0, \eps_0)$, we have for all $t \in [0,T_\delta^\star]$
    \begin{align}
\begin{split}\label{ineq-bootstrap-amelioCOMPRESSIBLE} 
            \|(\nabla \ww,\nabla^2 \ww)(t)\|_{L^{\upinfty}_\alpha L^{\upinfty}_x} + \|(\nabla u,\nabla^2 u, \nabla ^3 u)(t) \|_{L^{\upinfty}_x}& \\
           +\|(\nabla  n,\nabla^2  n)(t) \|_{L^{\upinfty}_x} &\leq \frac{\delta}{2} e^{-\lambda' t}, \\
  \ \ \|\nabla \ww(t)\|_{L^{\upinfty}_\alpha L^{2}_x} + \|\nabla u(t) \|_{L^{2}_x}+\|\nabla  n(t) \|_{L^{2}_x}  &\leq \frac{\delta}{2} e^{-\lambda t}, \\
  \Vert  \langle  \ww(t)- u(t) \rangle \Vert_{L^p_\alpha}
 &\leq \frac{\delta}{2} e^{-\lambda t}.
 \end{split}
    \end{align}
\end{Lem}
The proof of Lemma~\ref{lem-bootstrap-amelioCOMPRESSIBLE} follows the same arguments as in Lemma~\ref{lem-bootstrap-amelio} in the incompressible case, that is the proof of Lemma~\ref{lem-bootstrap-amelio}. Namely we use the same interpolation procedure, combining the exponential bound from Proposition \ref{prop-refinedLWP-VNScomp} and the improved decay estimates from Lemma \ref{lem:ImproveLinL2-COMPRESSIBLE}, we choose the constant $\theta=\theta_k$ (with $k$ fixed before the bootstrap argument) as in the proof ofLemma~\ref{lem-bootstrap-amelio}, followed by $\sigma>0$ large enough and then $\delta, \eps_0$ small enough.  Details are left to the reader.

\subsection{Conclusion}

We can now end the proof of Theorem \ref{thm-globalconvergence-comp}. The conclusion of the bootstrap and the final convergences are the same as in Section \ref{sectionVNS-cclbootstrap}. Note that the convergence towards zero of $\nabla \varrho=\nabla n$ comes from \eqref{assumpExtra-bootstrapCOMP}.

\section{The case of nonlinear drag force}\label{Section-nonlindrag}
In this last section, we explain how to take into account the nonlinear drag term $\Gamma (U-v^\alpha)$ in the asymptotic stability of Vlasov--Navier--Stokes type systems studied in this chapter, both in the incompressible case  \eqref{eq:VNS-CHAP4} and compressible case \eqref{eq:VNScompressibleCHAP4}. 
More precisely, we are concerned with the adaptation of the  proofs of Theorems \ref{thm-globalconvergence}--\ref{thm-globalconvergence-comp} in the case where $\Gamma \neq 0$. 
 Let us recall that the quadratic nonlinear term $\Gamma\in \mathscr{C}^{\infty}(\R^3;\R^3)$ is asked to satisfy $\Gamma(0)=0$ and $\Gamma'(0)=0$.

\medskip

We can follow the same lines of proof as what we have previously done in this chapter, starting from the perturbed systems \eqref{eq:VNSpertub} in the incompressible case and \eqref{eq:VNScomp-pertub1} in the compressible case and view the new nonlinear term $\Gamma (u-w^\alpha)$ due to the nonlinear drag as a an additional source term. Therefore it only remains to  study the  contributions of the terms
\begin{align*}
\mathscr{S}^\alpha&\vcentcolon=\Gamma(u-w^\alpha), \\
\mathscr{S}_{\mathrm{incomp}}&\vcentcolon=-\int (1+r^\alpha) \Gamma(u-w^\alpha) \, \mathrm{d}\mu(\alpha), \\
\mathscr{S}_{\mathrm{comp}}&\vcentcolon=-\frac{1}{1+ n}\int (1+r^\alpha) \Gamma(u-w^\alpha) \, \mathrm{d}\mu(\alpha),
\end{align*}
where  $\mathscr{S}_{\mathrm{incomp}}$ is for the incompressible case \eqref{eq:VNSpertub} and $\mathscr{S}_{\mathrm{comp}}$ for the compressible case \eqref{eq:VNScomp-pertub1}. The term $\mathscr{S}^\alpha$ appears in the equation for $w^\alpha$, while $\mathscr{S}_{\mathrm{incomp}}$ or $\mathscr{S}_{\mathrm{comp}}$ appears in the equation for $u$.

\medskip

For the sake of conciseness, we focus on the incompressible regime, that is $\mathscr{S}^\alpha$ and $\mathscr{S}_{\mathrm{incomp}}$, the treatment of $\mathscr{S}_{\mathrm{comp}}$ being similar. The main tool we need is the composition rule from Proposition \ref{prop-Sobolev}. They appear in a few instances in the proof:
\begin{itemize}
    \item[$\bullet$] when we rely on the linear analysis to obtain decay, that is in Lemmas \ref{lem:ImproveLinL2-nabla-u} and \ref{lem:ImproveLdiffMoyenne};
    \item[$\bullet$]  when we perform a direct energy estimate for $\nabla w^\alpha$, that is in Lemma \ref{lem:ImproveLinL2-nabla-walpha}.
\end{itemize}
Let us emphasize that in the equations satisfied by $w^\alpha$, the term $\Gamma(u-w^\alpha)$ is not integrated in $\alpha$ and terms such as $\nabla (\Gamma(u-w^\alpha)) \cdot \nabla w^\alpha$ do not have a well-defined sign; we are thus enforced to use the rough estimate
\begin{align*}
    \int \nabla (\Gamma(u-w^\alpha)) \cdot \nabla w^\alpha \, \mathrm{d}x \leq \Vert \nabla (\Gamma(u-w^\alpha))\Vert_{L^2_x} \Vert \nabla w^\alpha \Vert_{L^2_x},
\end{align*}
and then to take the supremum in $\alpha$, since we want to improve the decay of $\Vert \nabla \ww \Vert_{L^{\upinfty}_\alpha L^2_x}$. In view of the composition rules from Proposition \ref{prop-Sobolev}, we will have to deal with some factors depending on $\Vert u-w^\alpha \Vert_{L^{\upinfty}_x}$, hence on
\begin{align}\label{estim-diffvelocity-Linfty}
    \Vert u-w^\alpha \Vert_{L^{\upinfty}_x} \lesssim \vert \langle u-w^\alpha \rangle \vert+\Vert \nabla (u-w^\alpha) \Vert_{L^{\upinfty}_x}.
\end{align}
This justifies why we need to restrict $p=\infty$ when treating nonlinear drags (recall \eqref{cond-exposant-p-thm} in the statement of Theorem \ref{thm-globalconvergence}); in particular we take $p=\infty$ in the Proposition \ref{Prop-lineardecayVNS} about linear decay and in  the definition of $T^\star_\delta$ for the bootstrap argument \eqref{assump-bootstrapVNS}, we demand a $L^\infty_\alpha$ control for the difference of the spatial means of the velocities. 

\medskip

On the other hand, we argue as follows to estimate the new contributions of $\mathscr{S}^\alpha$ and $\mathscr{S}_{\mathrm{incomp}}$, in view of applying Proposition \ref{Prop-lineardecayVNS}. At the end of the day, it amounts to controlling $\Vert\Gamma (u-w^\alpha) \Vert_{H^1_x}$. Thanks to the composition estimates from Proposition \ref{prop-Sobolev} and Poincaré inequality, we know that
\begin{multline*}
    \Vert\Gamma (u-w^\alpha) \Vert_{H^1_x} \\ 
     \lesssim \mathrm{C}(\Gamma'', \Vert u-w^\alpha \Vert_{L^{\upinfty}_x})\Vert u-w^\alpha \Vert_{L^{\upinfty}_x} \big(\vert \langle u-w^\alpha \rangle \vert+\Vert \nabla (u-w^\alpha) \Vert_{L^2_x}\big),
\end{multline*}
for some increasing function $\mathrm{C}$ with respect to all its parameters. In view of \eqref{estim-diffvelocity-Linfty}, we end up with
\begin{multline*}
    \Vert\Gamma (u-w^\alpha) \Vert_{H^1_x} \\
    \lesssim \mathrm{C}(\Gamma'', \vert \langle u-w^\alpha \rangle \vert+\Vert \nabla (u-w^\alpha) \Vert_{L^{\upinfty}_x})\left( \vert \langle u-w^\alpha \rangle \vert^2+\Vert \nabla (u-w^\alpha) \Vert_{L^{\upinfty}_x}^2 \right).
\end{multline*}
By taking the supremum in $\alpha$ and using \eqref{assump-bootstrapVNS} (recall $p=\infty$), we can argue as in the treatment of the other source terms in the proof of Lemmas \ref{lem:ImproveLinL2-nabla-u} and \ref{lem:ImproveLdiffMoyenne}. Details of the whole procedure are left to the reader.

\appendix

 \printindex

 \bibliography{biblio}

@article {BDGM,
    AUTHOR = {Boudin, L. and Desvillettes, L. and Grandmont, C. and Moussa, A.},
     TITLE = {Global existence of solutions for the coupled {V}lasov and
{N}avier-{S}tokes equations},
   JOURNAL = {Differ. Integral Equ.},
  FJOURNAL = {Differential and Integral Equations. An International Journal for
Theory \& Applications},
    VOLUME = {22},
      YEAR = {2009},
    NUMBER = {11-12},
     PAGES = {1247--1271},
      ISSN = {0893-4983},
   MRCLASS = {76D05 (35A01 35D30 35Q35 76T10)},
  MRNUMBER = {2555647 (2011e:76039)},
MRREVIEWER = {Xianpeng Hu}
}

@article {BGM,
    AUTHOR = {Boudin, L. and Grandmont, C.
              and Moussa, A.},
     TITLE = {Global existence of solutions to the incompressible
              {N}avier-{S}tokes-{V}lasov equations in a time-dependent domain},
   JOURNAL = {J. Differ. Equations},
    VOLUME = {262},
      YEAR = {2017},
     PAGES = {1317--1340},
}

@article{ABdM,
  title={The existence of the global generalized solution of the system of equations describing suspension motion},
  author={Anoshchenko, O. and Boutet de Monvel-Berthier, A.},
 FJournal = {Mathematical Methods in the Applied Sciences},
 Journal = {Math. Methods Appl. Sci.},
   volume={20},
  number={6},
  pages={495--519},
  year={1997},
  publisher={Wiley Online Library}
}

@article{HKMM,
  title={Large Time Behavior of the {V}lasov-{N}avier-{S}tokes System on the Torus},
  author={Han-Kwan, D. and Moussa, A. and Moyano, I.},
  journal={Arch. Ration. Mech. Anal.},
  volume={236},
  number={3},
  pages={1273--1323},
  year={2020},
  publisher={Springer}
}

@article {HKH,
    AUTHOR = {Han-Kwan, D. and Hauray, M.},
     TITLE = {Stability issues in the quasineutral limit of the
              one-dimensional {V}lasov-{P}oisson equation},
   JOURNAL = {Comm. Math. Phys.},
  FJOURNAL = {Communications in Mathematical Physics},
    VOLUME = {334},
      YEAR = {2015},
    NUMBER = {2},
     PAGES = {1101--1152},
      ISSN = {0010-3616},
   MRCLASS = {35Q83 (35B35)},
  MRNUMBER = {3306612},
MRREVIEWER = {Calvin Tadmon},
       DOI = {10.1007/s00220-014-2217-4},
       URL = {https://doi.org/10.1007/s00220-014-2217-4},
}

@article{Bar,
  title={Nonlinear instability in {V}lasov type equations around rough velocity profiles},
  author={Baradat, A.},
  journal={Annales de l'Institut Henri Poincar{\'e} C, Analyse non lin{\'e}aire},
  volume={37},
  number={3},
  pages={489--547},
  year={2020},
}

@article {Gr96,
    AUTHOR = {Grenier, E.},
     TITLE = {Oscillations in quasineutral plasmas},
   JOURNAL = {Comm. Partial Differential Equations},
  FJOURNAL = {Communications in Partial Differential Equations},
    VOLUME = {21},
      YEAR = {1996},
    NUMBER = {3-4},
     PAGES = {363--394},
      ISSN = {0360-5302},
   MRCLASS = {82C22 (35Q99 76X05 82D10)},
  MRNUMBER = {1387452},
MRREVIEWER = {Gershon Wolansky},
       DOI = {10.1080/03605309608821189},
       URL = {https://doi.org/10.1080/03605309608821189},
}

@Article{Gr00,
 Author = {Grenier, E.},
 Title = {On the nonlinear instability of {Euler} and {Prandtl} equations.},
 FJournal = {Communications on Pure and Applied Mathematics},
 Journal = {Commun. Pure Appl. Math.},
 ISSN = {0010-3640},
 Volume = {53},
 Number = {9},
 Pages = {1067--1091},
 Year = {2000},
 Language = {English},
 DOI = {10.1002/1097-0312(200009)53:9<1067::AID-CPA1>3.0.CO;2-Q},
 zbMATH = {1591999},
 Zbl = {1048.35081}
}

@article{HKN,
  title={Ill-posedness of the hydrostatic {E}uler and singular {V}lasov equations},
  author={Han-Kwan, D. and Nguyen, T.},
  journal={Arch. Ration. Mech. Anal.},
  volume={221},
  number={3},
  pages={1317--1344},
  year={2016},
  publisher={Springer}
}

@Article{CGG,
 Author = {Cordier, S. and Grenier, E. and Guo, Y.},
 Title = {Two-stream instabilities in plasmas},
 FJournal = {Methods and Applications of Analysis},
 Journal = {Methods Appl. Anal.},
 ISSN = {1073-2772},
 Volume = {7},
 Number = {2},
 Pages = {391--405},
 Year = {2000},
 Language = {English},
 DOI = {10.4310/MAA.2000.v7.n2.a7},
 zbMATH = {1757307},
 Zbl = {1002.82029}
}

@article {GS,
    AUTHOR = {Guo, Y. and Strauss, W. A.},
     TITLE = {Nonlinear instability of double-humped equilibria},
   JOURNAL = {Ann. Inst. H. Poincar\'e Anal. Non Lin\'eaire},
  FJOURNAL = {Annales de l'Institut Henri Poincar\'e. Analyse Non
              Lin\'eaire},
    VOLUME = {12},
      YEAR = {1995},
    NUMBER = {3},
     PAGES = {339--352},
}

@article{Zak,
title={Benney equations and quasiclassical approximation in the method of the inverse problem},
  author={Zakharov, V.E.},
  journal={Functional analysis and its applications},
  volume={14},
  number={2},
  pages={89--98},
  year={1980},
  publisher={Springer}
  }

@article{zakharov1981benney,
  title={On the {B}enney equations},
  author={Zakharov, V.E.},
  journal={Physica D: Nonlinear Phenomena},
  volume={3},
  number={1-2},
  pages={193--202},
  year={1981},
  publisher={North-Holland}
}

@article {HKNsima,
    AUTHOR = {Han-Kwan, D. and Nguyen, T.},
     TITLE = {Nonlinear {I}nstability of {V}lasov--{M}axwell {S}ystems in
              the {C}lassical and {Q}uasineutral {L}imits},
   JOURNAL = {SIAM J. Math. Anal.},
  FJOURNAL = {SIAM Journal on Mathematical Analysis},
    VOLUME = {48},
      YEAR = {2016},
    NUMBER = {5},
     PAGES = {3444--3466},
      ISSN = {0036-1410},
     CODEN = {SJMAAH},
   MRCLASS = {35Q83 (35B35)},
  MRNUMBER = {3553927},
       DOI = {10.1137/15M1028765},
       URL = {http://dx.doi.org/10.1137/15M1028765},
}

@article{Pen60,
        author = {O. {Penrose}},
        title = {{Electrostatic instability of a uniform non-Maxwellian plasma}},
        journal = {Phys. Fluids},
        volume = {3},
        number = {},
        year = {1960},
	pages = {258-265}
}

@article {Br97,
    AUTHOR = {Brenier, Y.},
     TITLE = {A homogenized model for vortex sheets},
   JOURNAL = {Arch. Rational Mech. Anal.},
  FJOURNAL = {Archive for Rational Mechanics and Analysis},
    VOLUME = {138},
      YEAR = {1997},
    NUMBER = {4},
     PAGES = {319--353},
      ISSN = {0003-9527},
     CODEN = {AVRMAW},
   MRCLASS = {76C05 (35Q35 49J40 76M30)},
  MRNUMBER = {1467558 (98m:76030)},
MRREVIEWER = {Vladimir V. Shelukhin},
       DOI = {10.1007/s002050050044},
       URL = {http://dx.doi.org/10.1007/s002050050044},
}

@incollection {Brehyqua,
    AUTHOR = {Brenier, Y.},
     TITLE = {Comparaison du r\'egime hydrostatique des fluides
              incompressibles non-visqueux et du r\'egime quasineutre des
              plasmas},
 BOOKTITLE = {Trends in applications of mathematics to mechanics ({N}ice,
              1998)},
    SERIES = {Chapman \& Hall/CRC Monogr. Surv. Pure Appl. Math.},
    VOLUME = {106},
     PAGES = {285--292},
 PUBLISHER = {Chapman \& Hall/CRC, Boca Raton, FL},
      YEAR = {2000},
   MRCLASS = {76B99 (35Q35 76X05)},
  MRNUMBER = {1734891},
}

@article {Br-var1,
    AUTHOR = {Brenier, Y.},
     TITLE = {The least action principle and the related concept of
              generalized flows for incompressible perfect fluids},
   JOURNAL = {J. Amer. Math. Soc.},
  FJOURNAL = {Journal of the American Mathematical Society},
    VOLUME = {2},
      YEAR = {1989},
    NUMBER = {2},
     PAGES = {225--255},
      ISSN = {0894-0347},
   MRCLASS = {58D05 (35Q10 49H05 58E30 58F11 76A02 76C05)},
  MRNUMBER = {969419},
MRREVIEWER = {Alberto Valli},
       DOI = {10.2307/1990977},
       URL = {http://dx.doi.org/10.2307/1990977},
}

@article {Br-var2,
    AUTHOR = {Brenier, Y.},
     TITLE = {The dual least action problem for an ideal, incompressible
              fluid},
   JOURNAL = {Arch. Rational Mech. Anal.},
  FJOURNAL = {Archive for Rational Mechanics and Analysis},
    VOLUME = {122},
      YEAR = {1993},
    NUMBER = {4},
     PAGES = {323--351},
      ISSN = {0003-9527},
     CODEN = {AVRMAW},
   MRCLASS = {58E30 (58E50 76C05 76M30)},
  MRNUMBER = {1217592},
MRREVIEWER = {J. E. Marsden},
       DOI = {10.1007/BF00375139},
       URL = {http://dx.doi.org/10.1007/BF00375139},
}

@article {Br-var3,
    AUTHOR = {Brenier, Y.},
     TITLE = {Minimal geodesics on groups of volume-preserving maps and
              generalized solutions of the {E}uler equations},
   JOURNAL = {Comm. Pure Appl. Math.},
  FJOURNAL = {Communications on Pure and Applied Mathematics},
    VOLUME = {52},
      YEAR = {1999},
    NUMBER = {4},
     PAGES = {411--452},
      ISSN = {0010-3640},
     CODEN = {CPAMA},
   MRCLASS = {58E10 (37K65 49Q20 58D05 76B03 76M30)},
  MRNUMBER = {1658919},
MRREVIEWER = {Vladimir V. Shelukhin},
       DOI = {10.1002/(SICI)1097-0312(199904)52:4<411::AID-CPA1>3.0.CO;2-3},
       URL =
              {http://dx.doi.org/10.1002/(SICI)1097-0312(199904)52:4<411::AID-CPA1>3.0.CO;2-3},
}

@Book{pazy2012semigroups,
 Author = {Pazy, A.},
 Title = {Semigroups of linear operators and applications to partial differential equations},
 FSeries = {Applied Mathematical Sciences},
 Series = {Appl. Math. Sci.},
 ISSN = {0066-5452},
 Volume = {44},
 Year = {1983},
 Publisher = {Springer, Cham},
 Language = {English},
 DOI = {10.1007/978-1-4612-5561-1},
 zbMATH = {3816434},
 Zbl = {0516.47023}
}

@Article{BHK,
 Author = {Brigouleix, N. and Han-Kwan, D.},
 Title = {The non-relativistic limit of the {Vlasov}-{Maxwell} system with uniform macroscopic bounds},
 FJournal = {Annales de la Facult{\'e} des Sciences de Toulouse. Math{\'e}matiques. S{\'e}rie VI},
 Journal = {Ann. Fac. Sci. Toulouse, Math. (6)},
 ISSN = {0240-2963},
 Volume = {31},
 Number = {2},
 Pages = {545--594},
 Year = {2022},
 Language = {English},
 DOI = {10.5802/afst.1702},
 zbMATH = {7549948},
 Zbl = {1492.35340}
}

@article{BreMoy,
  title={Relaxed solutions for incompressible inviscid flows: A variational and gravitational approximation to the initial value problem},
  author={Brenier, Y. and Moyano, I.},
  journal={Philosophical Transactions of the Royal Society A},
  volume={380},
  number={2219},
  pages={20210078},
  year={2022},
  publisher={The Royal Society}
}

@book {majda1984compressible,
    AUTHOR = {Majda, A.},
     TITLE = {Compressible fluid flow and systems of conservation laws in
              several space variables},
    SERIES = {Applied Mathematical Sciences},
    VOLUME = {53},
 PUBLISHER = {Springer-Verlag, New York},
      YEAR = {1984},
     PAGES = {viii+159},
      ISBN = {0-387-96037-6},
   MRCLASS = {35L65 (76L05 76N10)},
  MRNUMBER = {748308},
MRREVIEWER = {Joel Smoller},
       DOI = {10.1007/978-1-4612-1116-7},
       URL = {https://doi.org/10.1007/978-1-4612-1116-7},
}

@article {Besse2009multi,
    AUTHOR = {Besse, N. and Berthelin, F. and Brenier, Y. and
              Bertrand, P.},
     TITLE = {The multi-water-bag equations for collisionless kinetic
              modeling},
   JOURNAL = {Kinet. Relat. Models},
  FJOURNAL = {Kinetic and Related Models},
    VOLUME = {2},
      YEAR = {2009},
    NUMBER = {1},
     PAGES = {39--80},
      ISSN = {1937-5093,1937-5077},
   MRCLASS = {35F20 (35A35 35Q83 65M60 82C40)},
  MRNUMBER = {2472149},
MRREVIEWER = {Zhaohui\ Huo},
       DOI = {10.3934/krm.2009.2.39},
       URL = {https://doi.org/10.3934/krm.2009.2.39},
}

@article {Bardos2013cauchy,
    AUTHOR = {Bardos, C. and Besse, N.},
     TITLE = {The {C}auchy problem for the {V}lasov-{D}irac-{B}enney
              equation and related issues in fluid mechanics and
              semi-classical limits},
   JOURNAL = {Kinet. Relat. Models},
  FJOURNAL = {Kinetic and Related Models},
    VOLUME = {6},
      YEAR = {2013},
    NUMBER = {4},
     PAGES = {893--917},
      ISSN = {1937-5093,1937-5077},
   MRCLASS = {35Q83 (35B35 35L45 76X05 82D10)},
  MRNUMBER = {3177634},
MRREVIEWER = {Jonathan\ Ben-Artzi},
       DOI = {10.3934/krm.2013.6.893},
       URL = {https://doi.org/10.3934/krm.2013.6.893},
}

@InCollection{Chemin-Leray,
 Author = {Chemin, J.-Y.},
 Title = {The incompressible {Navier}-{Stokes} system seventy years after {Jean} {Leray}},
 BookTitle = {Actes des journ\'ees math\'ematiques \`a la m\'emoire de Jean Leray, Nantes, France, juin 17--18, 2002},
 ISBN = {2-85629-160-0},
 Pages = {99--123},
 Year = {2004},
 Publisher = {Paris: Soci{\'e}t{\'e} Math{\'e}matique de France},
 Language = {French},
 zbMATH = {2159057},
 Zbl = {1075.35035}
}

@article{besse2011waterbag,
  title={On the waterbag continuum},
  author={Besse, N.},
  journal={Archive for rational mechanics and analysis},
  volume={199},
  number={2},
  pages={453--491},
  year={2011},
  publisher={Springer}
}

@article{ambrosio2014continuity,
  title="{Continuity equations and ODE flows with non-smooth velocity}",
  author={Ambrosio, L. and Crippa, G.},
  journal={Proceedings of the Royal Society of Edinburgh Section A: Mathematics},
  volume={144},
  number={6},
  pages={1191--1244},
  year={2014},
  publisher={Royal Society of Edinburgh Scotland Foundation}
}

@Article{GPI,
 Author = {Griffin-Pickering, M. and Iacobelli, M.},
 Title = {Global well-posedness for the {Vlasov}-{Poisson} system with massless electrons in the 3-dimensional torus},
 FJournal = {Communications in Partial Differential Equations},
 Journal = {Commun. Partial Differ. Equations},
 ISSN = {0360-5302},
 Volume = {46},
 Number = {10},
 Pages = {1892--1939},
 Year = {2021},
 Language = {English},
 DOI = {10.1080/03605302.2021.1913750},
 zbMATH = {7433741},
 Zbl = {1478.35205}
}

@article{griffin2020globalwholespace,
Author = {Griffin-Pickering, M. and Iacobelli, M.},
 Title = {Global strong solutions in {{\(\mathbb{R}^3\)}} for ionic {Vlasov}-{Poisson} systems},
 FJournal = {Kinetic and Related Models},
 Journal = {Kinet. Relat. Models},
 ISSN = {1937-5093},
 Volume = {14},
 Number = {4},
 Pages = {571--597},
 Year = {2021},
 Language = {English},
 DOI = {10.3934/krm.2021016},
 zbMATH = {7450818},
 Zbl = {1476.35274}
}

@article{bouchut1991global,
  title={Global weak solution of the {V}lasov--{P}oisson system for small electrons mass},
  author={Bouchut, F.},
  journal={Communications in partial differential equations},
  volume={16},
  number={8-9},
  pages={1337--1365},
  year={1991},
  publisher={Taylor \& F.}
}

@article{arsenev1975existence,
  title={Existence in the large of a weak solution to the {V}lasov system of equations},
  author={Arsenev, A.A.},
  journal={Zhurnal Vychislitelnoi Matematiki i Matematicheskoi Fiziki},
  volume={15},
  pages={136--147},
  year={1975}
}

@Article{BattRein-weak,
 Author = {Batt, J. and Rein, G.},
 Title = {A rigorous stability result for the {Vlasov}-{Poisson} system in three dimensions},
 FJournal = {Annali di Matematica Pura ed Applicata. Serie Quarta},
 Journal = {Ann. Mat. Pura Appl. (4)},
 ISSN = {0373-3114},
 Volume = {164},
 Pages = {133--154},
 Year = {1993},
 Language = {English},
 DOI = {10.1007/BF01759319},
 zbMATH = {495898},
 Zbl = {0791.49030}
}

@Article{BattRein-strong,
 Author = {Batt, J. and Rein, G.},
 Title = {Global classical solutions of the periodic {Vlasov}-{Poisson} system in three dimensions},
 FJournal = {Comptes Rendus de l'Acad{\'e}mie des Sciences. S{\'e}rie I},
 Journal = {C. R. Acad. Sci., Paris, S{\'e}r. I},
 ISSN = {0764-4442},
 Volume = {313},
 Number = {6},
 Pages = {411--416},
 Year = {1991},
 Language = {English},
 zbMATH = {21973},
 Zbl = {0741.35058}
}

@article{lions1991propagation,
  title={Propagation of moments and regularity for the 3-dimensional {V}lasov-{P}oisson system},
  author={Lions, P-L. and Perthame, B.},
 FJournal = {Inventiones Mathematicae},
 Journal = {Invent. Math.},
   volume={105},
  number={1},
  pages={415--430},
  year={1991},
  publisher={Springer}
}

@article{BardosDegond,
  title={Global existence for the {V}lasov-{P}oisson equation in 3 space variables with small initial data},
  author={Bardos, C. and Degond, P.},
 FJournal = {{Annales de l'Institut Henri Poincar\'e. Analyse Non Lin\'eaire}},
 Journal = {{Ann. Inst. Henri Poincar\'e, Anal. Non Lin\'eaire}},
   volume={2},
  number={2},
  pages={101--118},
  year={1985},
  organization={Elsevier}
}

@Article{Pfaffelmoser,
 Author = {Pfaffelmoser, K.},
 Title = {Global classical solutions of the {Vlasov}-{Poisson} system in three dimensions for general initial data},
 FJournal = {Journal of Differential Equations},
 Journal = {J. Differ. Equations},
 ISSN = {0022-0396},
 Volume = {95},
 Number = {2},
 Pages = {281--303},
 Year = {1992},
 Language = {English},
 DOI = {10.1016/0022-0396(92)90033-J},
 zbMATH = {33010},
 Zbl = {0810.35089}
}

@article{Schaeffer,
 author = {Schaeffer, J.},
 title = {Global existence of smooth solutions to the {Vlasov}-{Poisson} system in three dimensions},
 fjournal = {Communications in Partial Differential Equations},
 journal = {Commun. Partial Differ. Equations},
 issn = {0360-5302},
 volume = {16},
 number = {8-9},
 pages = {1313--1335},
 year = {1991},
 language = {English},
 doi = {10.1080/03605309108820801},
 zbMATH = {28022},
 Zbl = {0746.35050}
}

@book {dellacherie1978probabilities,
    AUTHOR = {Dellacherie, C. and Meyer, P.-A.},
     TITLE = {Probabilities and potential},
    SERIES = {North-Holland Mathematics Studies},
    VOLUME = {29},
 PUBLISHER = {North-Holland Publishing Co., Amsterdam-New York},
      YEAR = {1978},
     PAGES = {viii+189}
}

@article{benachour2003global,
  title={{Global existence for the {Vlasov}--{Darwin} system in R3 for small initial data}},
  author={Benachour, S. and Filbet, F. and Lauren{\c{c}}ot, P. and Sonnendr{\"u}cker, E.},
  journal={Mathematical methods in the applied sciences},
  volume={26},
  number={4},
  pages={297--319},
  year={2003},
  publisher={Wiley Online Library}
}

@article{pallard2006initial,
  title={The initial value problem for the relativistic {Vlasov}-{Darwin} system},
  author={Pallard, C.},
  journal={International Mathematics Research Notices},
  volume={2006},
  number={9},
  pages={57191--57191},
  year={2006},
  publisher={OUP}
}

@techreport{gardner1960similarity,
  title={Similarity in the asymptotic behavior of collision-free hydromagnetic waves and water waves},
  author={Gardner, C.S. and Morikawa, G.K.},
  year={1960},
  institution={New York Univ., New York. Inst. of Mathematical Sciences}
}

@article{alonso2024well,
  title={Well-posedness for an hyperbolic--hyperbolic--elliptic system describing cold plasmas},
  author={Alonso-Or{\'a}n, D. and Granero-Belinch{\'o}n, R.},
  journal={Applied Mathematics Letters},
  volume={147},
  pages={108863},
  year={2024},
  publisher={Elsevier}
}

@article{natalini2020mean,
title = {On the mean field limit for {Cucker}-{Smale} models},
  author={Natalini, R. and Paul, T.},
journal = {Discrete and Continuous Dynamical Systems - B},
volume = {27},
number = {5},
pages = {2873-2889},
year = {2022},
}

@article{natalini2023mean,
  title={The Mean-Field limit for hybrid models of collective motions with chemotaxis},
  author={Natalini, R. and Paul, T.},
  journal={SIAM Journal on Mathematical Analysis},
  volume={55},
  number={2},
  pages={900--928},
  year={2023},
  publisher={SIAM}
}

@article{HS-VSmeanfield2,
  title={Sedimentation of particles with very small inertia II: Derivation, {C}auchy problem and hydrodynamic limit of the {V}lasov-{S}tokes equation
},
  author={Höfer, R.M. and Schubert, R.},
  journal={arXiv preprint arXiv:2311.01891},
  year={2023}
}

@book{MB,
 author = {Majda, A. and Bertozzi, A.},
 title = {Vorticity and incompressible flow},
 fseries = {Cambridge Texts in Applied Mathematics},
 series = {Camb. Texts Appl. Math.},
 isbn = {0-521-63057-6; 0-521-63948-4},
 year = {2002},
 publisher = {Cambridge: Cambridge University Press},
 language = {English},
 doi = {10.1017/CBO9780511613203},
 zbMATH = {1644218},
 Zbl = {0983.76001}
}

@InCollection{BB2,
 Author = {Bardos, C. and Besse, N.},
 Title = {Hamiltonian structure, fluid representation and stability for the {Vlasov}-{Dirac}-benney equation},
 BookTitle = {Hamiltonian partial differential equations and applications.},
 ISBN = {978-1-4939-2949-8; 978-1-4939-2950-4},
 Pages = {1--30},
 Year = {2015},
 Publisher = {Toronto: The Fields Institute for Research in the Mathematical Sciences; New York, NY: Springer},
}

@article {BB3,
 Author = {Bardos, C. and Besse, N.},
 Title = {Semi-classical limit of an infinite dimensional system of nonlinear {Schr{\"o}dinger} equations},
 FJournal = {Bulletin of the Institute of Mathematics. Academia Sinica. New Series},
 Journal = {Bull. Inst. Math., Acad. Sin. (N.S.)},
 ISSN = {2304-7909},
 Volume = {11},
 Number = {1},
 Pages = {43--61},
 Year = {2016},
 Language = {English},
 URL = {web.math.sinica.edu.tw/bulletin/archives_articlecontent16.jsp?bid=MjAxNjEwMw==},
 zbMATH = {6562114},
 Zbl = {1336.35331}
}

@article {loeper2006uniqueness,
    AUTHOR = {Loeper, G.},
     TITLE = {Uniqueness of the solution to the {V}lasov-{P}oisson system
              with bounded density},
   JOURNAL = {J. Math. Pures Appl. (9)},
  FJOURNAL = {Journal de Math\'ematiques Pures et Appliqu\'ees. Neuvi\`eme
              S\'erie},
    VOLUME = {86},
      YEAR = {2006},
    NUMBER = {1},
     PAGES = {68--79},
      ISSN = {0021-7824},
   MRCLASS = {82D10 (35D05 35F20 35Q35)},
  MRNUMBER = {2246357},
MRREVIEWER = {Christophe\ Pallard},
       DOI = {10.1016/j.matpur.2006.01.005},
       URL = {https://doi.org/10.1016/j.matpur.2006.01.005},
}

@article {Miot2016uniqueness,
    AUTHOR = {Miot, E.},
     TITLE = {A uniqueness criterion for unbounded solutions to the
              {V}lasov-{P}oisson system},
   JOURNAL = {Comm. Math. Phys.},
  FJOURNAL = {Communications in Mathematical Physics},
    VOLUME = {346},
      YEAR = {2016},
    NUMBER = {2},
     PAGES = {469--482},
      ISSN = {0010-3616,1432-0916},
   MRCLASS = {35Q83 (35A02)},
  MRNUMBER = {3535893},
MRREVIEWER = {Silvia\ Caprino},
       DOI = {10.1007/s00220-016-2707-7},
       URL = {https://doi.org/10.1007/s00220-016-2707-7},
}

@incollection {holding2018uniqueness,
    AUTHOR = {Holding, T. and Miot, E.},
     TITLE = {Uniqueness and stability for the {V}lasov-{P}oisson system
              with spatial density in {O}rlicz spaces},
 BOOKTITLE = {Mathematical analysis in fluid mechanics---selected recent
              results},
    SERIES = {Contemp. Math.},
    VOLUME = {710},
     PAGES = {145--162},
 PUBLISHER = {Amer. Math. Soc., [Providence], RI},
      YEAR = {2018},
}

@Article{HKM3,
 Author = {Han-Kwan, D. and Miot, E. and Moussa, A. and Moyano, I.},
 Title = {Uniqueness of the solution to the 2D {Vlasov}-{Navier}-{Stokes} system},
 FJournal = {Revista Matem{\'a}tica Iberoamericana},
 Journal = {Rev. Mat. Iberoam.},
 ISSN = {0213-2230},
 Volume = {36},
 Number = {1},
 Pages = {37--60},
 Year = {2020},
}

@misc{HKM3-3D,
      title={{On Uniqueness for the three-dimensional {Vlasov}-{Navier}-{Stokes} system}}, 
      Author = {Han-Kwan, D. and Miot, E. and Moussa, A. and Moyano, I.},
      year={2025},
      Journal={2511.10099},
      archivePrefix={arXiv},
      url={https://arxiv.org/abs/2511.10099}, 
}

@article{danchin2024fujita,
 author = {Danchin, R.},
 title = {Fujita-{Kato} solutions and optimal time decay for the {Vlasov}-{Navier}-{Stokes} system in the whole space},
 fjournal = {Archive for Rational Mechanics and Analysis},
 journal = {Arch. Ration. Mech. Anal.},
 issn = {0003-9527},
 volume = {250},
 number = {2},
 pages = {45},
 note = {Id/No 13},
 year = {2026},
 language = {English},
 doi = {10.1007/s00205-026-02171-x},
 zbMATH = {8177170}
}

@book{BCD,
  title={Fourier analysis and nonlinear partial differential equations},
  author={Bahouri, H. and Chemin, J.-Y. and Danchin, R.},
  year={2011},
  publisher={Springer}
}

@Book{choquetbruhat,
 Author = {Choquet-Bruhat, Y.},
 Title = {General relativity and the {Einstein} equations.},
 FSeries = {Oxford Mathematical Monographs},
 Series = {Oxford Math. Monogr.},
 ISBN = {978-0-19-923072-3},
 Year = {2009},
 Publisher = {Oxford: Oxford University Press},
 Language = {English},
 zbMATH = {5382453},
 Zbl = {1157.83002}
}

@Article{Stokesweight,
 Author = {Specovius-Neugebauer, M.},
 Title = {Weak solutions of the {Stokes} problem in weighted {Sobolev} spaces},
 FJournal = {Acta Applicandae Mathematicae},
 Journal = {Acta Appl. Math.},
 ISSN = {0167-8019},
 Volume = {37},
 Number = {1-2},
 Pages = {195--203},
 Year = {1994},
 Language = {English},
 DOI = {10.1007/BF00995141},
 zbMATH = {709469},
 Zbl = {0814.35102}
}

@Article{EHKM,
 Author = {Ertzbischoff, L. and Han-Kwan, D. and Moussa, A.},
 Title = {Concentration versus absorption for the {Vlasov}-{Navier}-{Stokes} system on bounded domains},
 FJournal = {Nonlinearity},
 Journal = {Nonlinearity},
 ISSN = {0951-7715},
 Volume = {34},
 Number = {10},
 Pages = {6843--6900},
 Year = {2021},
 Language = {English},
 DOI = {10.1088/1361-6544/ac1558},
 zbMATH = {7396374},
 Zbl = {1486.35383}
}

@Article{HK-VNSR3,
 Author = {Han-Kwan, D.},
 Title = {Large-time behavior of small-data solutions to the {Vlasov}-{Navier}-{Stokes} system on the whole space},
 FJournal = {Probability and Mathematical Physics},
 Journal = {Probab. Math. Phys.},
 ISSN = {2690-0998},
 Volume = {3},
 Number = {1},
 Pages = {35--67},
 Year = {2022},
 Language = {English},
 DOI = {10.2140/pmp.2022.3.35},
 zbMATH = {7527322},
 Zbl = {1486.35384}
}

@article{choi2024revisit,
 author = {Choi, Y.-P. and Jung, J. and Kim, J.},
 title = {A revisit to the pressureless {Euler}-{Navier}-{Stokes} system in the whole space and its optimal temporal decay},
 fjournal = {Journal of Differential Equations},
 journal = {J. Differ. Equations},
 volume = {401},
 pages = {231--281},
 year = {2024},
}

@article{ChoiJungR3,
 author = {Choi, Y.-P. and Jung, J.},
 title = {On the {Cauchy} problem for the pressureless {Euler}-{Navier}-{Stokes} system in the whole space},
 fjournal = {Journal of Mathematical Fluid Mechanics},
 journal = {J. Math. Fluid Mech.},
 volume = {23},
 number = {4},
 pages = {16},
 note = {Id/No 99},
 year = {2021},
}

@article{HuangTangZouR3,
 author = {Huang, F. and Tang, H. and Zou, W.},
 title = {Global well-posedness and optimal time decay rates of solutions to the pressureless {Euler}-{Navier}-{Stokes} system},
 fjournal = {Communications in Mathematical Analysis and Applications},
 journal = {Commun. Math. Anal. Appl.},
 volume = {3},
 number = {4},
 pages = {582--623},
 year = {2024},
}

@article{ZZ,
 author = {Zhang, Z. and Zi, R.},
 title = {Convergence to equilibrium for the solution of the full compressible {Navier}-{Stokes} equations},
 fjournal = {Annales de l'Institut Henri Poincar{\'e}. Analyse Non Lin{\'e}aire},
 journal = {Ann. Inst. Henri Poincar{\'e}, Anal. Non Lin{\'e}aire},
 issn = {0294-1449},
 volume = {37},
 number = {2},
 pages = {457--488},
 year = {2020},
 language = {English},
 doi = {10.1016/j.anihpc.2019.09.001},
 zbMATH = {7180624},
 Zbl = {1435.35069}
}

@article{DW,
 author = {Danchin, R. and Wang, S.},
 title = {Exponential decay for inhomogeneous viscous flows on the torus},
 fjournal = {ZAMP. Zeitschrift f{\"u}r angewandte Mathematik und Physik},
 journal = {Z. Angew. Math. Phys.},
 issn = {0044-2275},
 volume = {75},
 number = {2},
 pages = {33},
 note = {Id/No 62},
 year = {2024},
 language = {English},
 doi = {10.1007/s00033-024-02198-8},
 zbMATH = {7839631},
 Zbl = {1540.35286}
}

@article{GWZ,
 author = {Guo, S. and Wu, G. and Zhang, Y.},
 title = {Global existence and large time behaviour for the pressureless {Euler}-{Navier}-{Stokes} system in {{\(\mathbb{R}^3\)}}},
 fjournal = {Proceedings of the Royal Society of Edinburgh. Section A. Mathematics},
 journal = {Proc. R. Soc. Edinb., Sect. A, Math.},
 volume = {154},
 number = {1},
 pages = {328--352},
 year = {2024},
}

@article{li2025global,
  title={Global existence and large-time behavior of solutions to the pressureless Euler/Navier-Stokes system},
  author={Li, H-L. and Tang, H-Z. and Zhang, Y.},
  journal={Science China Mathematics},
  pages={1--22},
  year={2025},
  publisher={Springer}
}

@article{LiShouZhang,
      title={Large-friction and incompressible limits for pressureless
Euler/isentropic Navier-Stokes flows}, 
      author={Li, H-L. and Shou, L-Y. and Zhang, Y.},
      year={2025},
      journal={arXiv preprint 	arXiv:2508.20730}, 
}

@article{bianchini2024hydrostatic,
  title={On the hydrostatic limit of stably stratified fluids with isopycnal diffusivity},
  author={Bianchini, R. and Duch{\^e}ne, V.},
  journal={Communications in Partial Differential Equations},
  volume={49},
  number={5-6},
  pages={543--608},
  year={2024},
  publisher={Taylor \& F.}
}

@article{fradin2024well,
  title={Well-posedness of the Euler equations in a stably stratified ocean in isopycnal coordinates},
  author={Fradin, T.},
  journal={Ann. Inst. H. Poincaré C Anal. Non Linéaire },
  year={2025},
  volume={published online first}
}

@article{bianchini-ertz2024review,
  title={Mathematical Insights into Hydrostatic Modeling of Stratified Fluids},
  author={Bianchini, R. and Ertzbischoff, L.},
  journal={Communications in Applied and Industrial Mathematics (CAIM)},
  volume={Accepted for publication},
  year={2024},
}

@phdthesis{hofer2020sedimentation,
  title={Sedimentation of particle suspensions in Stokes flows},
  author={Höfer, R.M.},
  year={2020},
  school={Universit{\"a}ts-und Landesbibliothek Bonn}
}

@book {silin,
    AUTHOR = {Silin, V.P.},
     TITLE = {Introduction to the Kinetic Theory of Gases [in Russian]},
    SERIES = {Applied Mathematical Sciences},
    VOLUME = {53},
 PUBLISHER = {Phys. Inst., Russ. Acad. Sci., Moscow},
      YEAR = {1971},
      }

@Book{Glassey-kinetic,
 Author = {Glassey, R.T.},
 Title = {The {Cauchy} problem in kinetic theory},
 Series = {Other Titles Appl. Math.},
 ISBN = {0-89871-367-6},
 Year = {1996},
 Publisher = {Philadelphia, PA: SIAM},
 Language = {English},
 zbMATH = {852535},
 Zbl = {0858.76001}
}

@Article{E1,
 Author = {Ertzbischoff, L.},
 Title = {Decay and absorption for the {Vlasov}-{Navier}-{Stokes} system with gravity in a half-space},
 FJournal = {Indiana University Mathematics Journal},
 Journal = {Indiana Univ. Math. J.},
 ISSN = {0022-2518},
 Volume = {73},
 Number = {1},
 Pages = {1--80},
 Year = {2024},
 Language = {English},
 DOI = {10.1512/iumj.2024.73.9538},
 zbMATH = {7834163},
 Zbl = {1540.35395}
}

@Article{GHKM,
 Author = {Glass, O. and Han-Kwan, D. and Moussa, A.},
 Title = {The {Vlasov}-{Navier}-{Stokes} system in a 2D pipe: existence and stability of regular equilibria},
 FJournal = {Archive for Rational Mechanics and Analysis},
 Journal = {Arch. Ration. Mech. Anal.},
 ISSN = {0003-9527},
 Volume = {230},
 Number = {2},
 Pages = {593--639},
 Year = {2018},
 Language = {English},
 DOI = {10.1007/s00205-018-1253-1},
 zbMATH = {6941735},
 Zbl = {1398.35242}
}

@article{desvillettes2010-model,
  title={Some aspects of the modeling at different scales of multiphase flows},
  author={Desvillettes, L.},
  journal={Computer methods in applied mechanics and engineering},
  volume={199},
  number={21-22},
  pages={1265--1267},
  year={2010},
  publisher={Elsevier}
}

@article{jeans1915theory,
  title={On the theory of star-streaming and the structure of the universe},
  author={Jeans, J.},
  journal={Monthly Notices of the Royal Astronomical Society, Vol. 76, p. 70-84},
  volume={76},
  pages={70--84},
  year={1915}
}

@article{vlasov1938vibrational,
  title={The vibrational properties of an electron gas},
  author={Vlasov, A.A.},
  journal={ J. Exp. Theor. Phys (In Russian)},
  volume={8},
  number={3},
  pages={291–318},
  year={1938},
}

@article{Landau,
  title={On the vibration of the electronic plasma},
  author={Landau, L.},
  journal={Zh. Eksp. Teor. Fiz.},
  volume={16},
  pages={574},
  year={1946},
}

@article{UO,
 author = {Ukai, S. and Okabe, T.},
 title = {On classical solutions in the large in time of two-dimensional {Vlasov}'s equation},
 fjournal = {Osaka Journal of Mathematics},
 journal = {Osaka J. Math.},
 issn = {0030-6126},
 volume = {15},
 pages = {245--261},
 year = {1978},
 language = {English},
 zbMATH = {3629435},
 Zbl = {0405.35002}
}

@article{Horst,
 author = {Horst, E.},
 title = {On the asymptotic growth of the solutions of the {Vlasov}--{Poisson} system},
 fjournal = {Mathematical Methods in the Applied Sciences},
 journal = {Math. Methods Appl. Sci.},
 issn = {0170-4214},
 volume = {16},
 number = {2},
 pages = {75--85},
 year = {1993},
 language = {English},
 doi = {10.1002/mma.1670160202},
 zbMATH = {166315},
 Zbl = {0782.35079}
}

@article{roberts1967nonlinear,
  title={Nonlinear evolution of a two-stream instability},
  author={Roberts, K.V. and Berk, H.L.},
  journal={Physical Review Letters},
  volume={19},
  number={6},
  pages={297},
  year={1967},
  publisher={APS}
}

@article{Jab,
  title={Large time concentrations for solutions to kinetic equations with energy dissipation},
  author={Jabin, P-E.},
 FJournal = {Communications in Partial Differential Equations},
 Journal = {Commun. Partial Differ. Equations},
   volume={25},
  number={3-4},
  pages={541--557},
  year={2000},
  publisher={Taylor \& F.}
}

@article{ChKw,
  title={Global well-posedness and large-time behavior for the inhomogeneous {V}lasov--{N}avier--{S}tokes equations},
  author={Choi, Y.-P. and Kwon, B.},
  journal={Nonlinearity},
  volume={28},
  number={9},
  pages={3309},
  year={2015},
  publisher={IOP Publishing}
}

@article{Arthur2,
 author = {Touati, A.},
 title = {The reverse {Burnett} conjecture for null dusts},
 fjournal = {Annals of PDE},
 journal = {Ann. PDE},
 issn = {2524-5317},
 volume = {11},
 number = {2},
 pages = {90},
 note = {Id/No 22},
 year = {2025},
 language = {English},
 doi = {10.1007/s40818-025-00213-3},
 zbMATH = {8084257}
}

@article{huneau2025,
      title={High-frequency backreaction for the {Einstein} equations under $\mathbb{U}(1)$ symmetry: from {Einstein-dust to Einstein-Vlasov}}, 
      author={Huneau, C. and Luk, J.},
      year={2025},
       journal = {https://arxiv.org/abs/2506.21779},
 arXiv = {arXiv:2506.21779}
}

@article{HuneauL2,
 Author = {Huneau, C. and Luk, J.},
 Title = {Burnett's conjecture in generalized wave coordinates},
 Year = {2024},
 journal = {Arxiv preprint https://arxiv.org/abs/2403.03470},
 arXiv = {arXiv:2403.03470}
}

@Article{HuneauL1,
 Author = {Huneau, C. and Luk, J.},
 Title = {Trilinear compensated compactness and {Burnett}'s conjecture in general relativity},
 FJournal = {Annales Scientifiques de l'{\'E}cole Normale Sup{\'e}rieure. Quatri{\`e}me S{\'e}rie},
 Journal = {Ann. Sci. {\'E}c. Norm. Sup{\'e}r. (4)},
 Volume = {57},
 Number = {2},
 Pages = {385--472},
 Year = {2024},
}

@article{HuneauL0,
 author = {Huneau, C. and Luk, J.},
 title = {High-frequency backreaction for the {Einstein} equations under polarized {{\(\mathbb U(1)\)}}-symmetry},
 fjournal = {Duke Mathematical Journal},
 journal = {Duke Math. J.},
 issn = {0012-7094},
 volume = {167},
 number = {18},
 pages = {3315--3402},
 year = {2018},
 language = {English},
 doi = {10.1215/00127094-2018-0035},
 zbMATH = {7009768},
 Zbl = {1412.35327}
}

@article{Friedmann,
  author  = {A. Friedmann},
  title   = {Über die Krümmung des Raumes},
  journal = {Zeitschrift für Physik},
  volume  = {10},
  pages   = {377--386},
  year    = {1922}
}

@article{Lemaitre,
  title={{Un Univers homog{\`e}ne de masse constante et de rayon croissant rendant compte de la vitesse radiale des n{\'e}buleuses extra-galactiques}},
  author={Lema{\^\i}tre, G.},
  journal={Annales de la Soci{\'e}t{\'e} Scientifique de Bruxelles, A47, p. 49-59},
  volume={47},
  pages={49--59},
  year={1927}
}

@book{Tolman,
  author    = {R. C. Tolman},
  title     = {Relativity, Thermodynamics, and Cosmology},
  publisher = {Oxford University Press},
  year      = {1934}
}

@article{Oppenheimer,
  title={On continued gravitational contraction},
  author={Oppenheimer, J. R. and Snyder, H.},
  journal={Physical Review},
  volume={56},
  number={5},
  pages={455},
  year={1939},
  publisher={APS}
}

@article{Godel,
  author  = {G\"odel, K.},
  title   = {An Example of a New Type of Cosmological Solutions of
             Einstein's Field Equations of Gravitation},
  journal = {Reviews of Modern Physics},
  volume  = {21},
  pages   = {447--450},
  year    = {1949},
  doi     = {10.1103/RevModPhys.21.447}
}

@Article{CDSM,
 Author = {Collot, C. and Danesi, E. and de Suzzoni, A-S. and Mal{\'e}z{\'e}, C.},
 Title = {Stability of homogeneous equilibria of the {Hartree}-{Fock} equation for its equivalent formulation for random fields},
 FJournal = {Probability and Mathematical Physics},
 Journal = {Probab. Math. Phys.},
 Volume = {6},
 Number = {1},
 Pages = {241--279},
 Year = {2025},
}

@Article{deSuzzoni,
 Author = {de Suzzoni, A.-S.},
 Title = {About systems of fermions with large number of particles: a probabilistic point of view},
 FJournal = {S{\'e}minaire Laurent Schwartz. EDP et Applications},
 Journal = {S{\'e}min. Laurent Schwartz, EDP Appl.},
 Volume = {2015-2016},
 Pages = {ex},
 Year = {2016},
}

@article{Pavlov,
 author = {Pavlov, M. V.},
 title = {Hamiltonian formalism of two-dimensional {Vlasov} kinetic equation},
 fjournal = {Proceedings of the Royal Society of London. Series A. Mathematical, Physical and Engineering Sciences},
 journal = {Proc. R. Soc. Lond., Ser. A, Math. Phys. Eng. Sci.},
 volume = {470},
 number = {2172},
 pages = {10},
 year = {2014},
}

@article{ChesnoPavlov,
 author = {Chesnokov, A. A. and Pavlov, M. V.},
 title = {Reductions of kinetic equations to finite component systems},
 fjournal = {Acta Applicandae Mathematicae},
 journal = {Acta Appl. Math.},
 volume = {122},
 number = {1},
 pages = {367--380},
 year = {2012},
}

@article{FeraPavlov,
 author = {Ferapontov, E. V. and Pavlov, M. V.},
 title = {Kinetic equation for soliton gas: integrable reductions},
 fjournal = {Journal of Nonlinear Science},
 journal = {J. Nonlinear Sci.},
 volume = {32},
 number = {2},
 pages = {22},
 note = {Id/No 26},
 year = {2022},
 language = {English},
}

@incollection{golse2016dynamics,
  title={On the dynamics of large particle systems in the mean field limit},
  author={Golse, F.},
  booktitle={Macroscopic and large scale phenomena: coarse graining, mean field limits and ergodicity},
  pages={1--144},
  year={2016},
  publisher={Springer}
}

@article{jabin2014review,
  title={{A review of the mean field limits for Vlasov equations}},
  author={Jabin, P.-E.},
  journal={Kinetic and Related models},
  volume={7},
  number={4},
  pages={661--711},
  year={2014},
  publisher={Kinetic and Related Models}
}

@article{MoussaSueur,
 author = {Moussa, A. and Sueur, F.},
 title = {On a {Vlasov}-{Euler} system for 2D sprays with gyroscopic effects},
 fjournal = {Asymptotic Analysis},
 journal = {Asymptotic Anal.},
 volume = {81},
 number = {1},
 pages = {53--91},
 year = {2013},
}

@article{LiShou,
 author = {Li, H.-L. and Shou, L.-Y.},
 title = {Global well-posedness of one-dimensional compressible {Navier}-{Stokes}-{Vlasov} system},
 fjournal = {Journal of Differential Equations},
 journal = {J. Differ. Equations},
 volume = {280},
 pages = {841--890},
 year = {2021},
}

@article{Choi-comp,
 author = {Choi, Y.-P.},
 title = {Large-time behavior for the {Vlasov}/compressible {Navier}-{Stokes} equations},
 fjournal = {Journal of Mathematical Physics},
 journal = {J. Math. Phys.},
 volume = {57},
 number = {7},
 pages = {071501, 13},
 year = {2016},
}

@article{Moschidis,
 author = {Moschidis, G.},
 title = {A proof of the instability of {AdS} for the {Einstein}-null dust system with an inner mirror},
 fjournal = {Analysis \& PDE},
 journal = {Anal. PDE},
 volume = {13},
 number = {6},
 pages = {1671--1754},
 year = {2020},
}

@article{Moschidis2,
 author = {Moschidis, G.},
 title = {A proof of the instability of {AdS} for the {Einstein}-massless {Vlasov} system},
 fjournal = {Inventiones Mathematicae},
 journal = {Invent. Math.},
 issn = {0020-9910},
 volume = {231},
 number = {2},
 pages = {467--672},
 year = {2023},
 language = {English},
 doi = {10.1007/s00222-022-01152-7},
 zbMATH = {7661328},
 Zbl = {1552.83015}
}

@article{danchin2025-TORUS2d,
 author = {Danchin, R. and Shou, L.-Y.},
 title = {Large-time asymptotics of periodic two-dimensional {Vlasov}-{Navier}-{Stokes} flows},
 fjournal = {Journal of the London Mathematical Society. Second Series},
 journal = {J. Lond. Math. Soc., II. Ser.},
 issn = {0024-6107},
 volume = {113},
 number = {1},
 pages = {43},
 note = {Id/No e70435},
 year = {2026},
 language = {English},
 doi = {10.1112/jlms.70435},
 zbMATH = {8153446}
}

@article{ELT,
 author = {Engelberg, S. and Liu, H. and Tadmor, E.},
 title = {Critical thresholds in {Euler}-{Poisson} equations},
 fjournal = {Indiana University Mathematics Journal},
 journal = {Indiana Univ. Math. J.},
 volume = {50},
 pages = {109--157},
 year = {2001},
}

@article {caffarelli1990interior,
    AUTHOR = {Caffarelli, L.A.},
     TITLE = {Interior {$W^{2,p}$} estimates for solutions of the
              {M}onge-{A}mp\`ere equation},
   JOURNAL = {Ann. of Math. (2)},
  FJOURNAL = {Annals of Mathematics. Second Series},
    VOLUME = {131},
      YEAR = {1990},
    NUMBER = {1},
     PAGES = {135--150},
}

@article{CCTT,
 author = {Carrillo, J.A. and Choi, Y.-P. and Tadmor, E. and Tan, C.},
 title = {Critical thresholds in 1D {Euler} equations with non-local forces},
 fjournal = {M\(^3\)AS. Mathematical Models \& Methods in Applied Sciences},
 journal = {Math. Models Methods Appl. Sci.},
 volume = {26},
 number = {1},
 pages = {185--206},
 year = {2016},
}

@incollection{LLS,
 author = {Lannes, D. and Linares, F. and Saut, J-C.},
 title = {The {Cauchy} problem for the {Euler}-{Poisson} system and derivation of the {Zakharov}-{Kuznetsov} equation},
 booktitle = {Studies in phase space analysis with applications to PDEs. In part selected papers based on the presentations at a meeting, Bertinoro, Italy, September 2011},
 pages = {181--213},
 year = {2013},
 publisher = {New York, NY: Birkh{\"a}user/Springer},
}

@article{BCK,
 author = {Bae, J. and Choi, J. and Kwon, B.},
 title = {Formation of singularities in plasma ion dynamics},
 fjournal = {Nonlinearity},
 journal = {Nonlinearity},
 volume = {37},
 number = {4},
 pages = {29},
 note = {Id/No 045011},
 year = {2024},
}

@article{HKapde,
 author = {Han-Kwan, D.},
 title = {On propagation of higher space regularity for nonlinear {Vlasov} equations},
 fjournal = {Analysis \& PDE},
 journal = {Anal. PDE},
 volume = {12},
 number = {1},
 pages = {189--244},
 year = {2019},
}

@book {Horm,
    AUTHOR = {H{\"o}rmander, L.},
     TITLE = {Linear partial differential operators},
 PUBLISHER = {Springer Verlag, Berlin-New York},
      YEAR = {1976},
     PAGES = {vii+285},
   MRCLASS = {35-XX (46FXX)},
  MRNUMBER = {0404822},
}

@article {Gerard,
    AUTHOR = {G{\'e}rard, P.},
     TITLE = {Moyennisation et r\'egularit\'e deux-microlocale},
   JOURNAL = {Ann. Sci. \'Ecole Norm. Sup. (4)},
  FJOURNAL = {Annales Scientifiques de l'\'Ecole Normale Sup\'erieure. Quatri\`eme
              S\'erie},
    VOLUME = {23},
      YEAR = {1990},
    NUMBER = {1},
     PAGES = {89--121},
}

@article{GIRR,
      title={From relativistic {V}lasov-{M}axwell to electron-{MHD} in the quasineutral regime}, 
      author={Gagnebin, A. and Iacobelli, M. and Rege, A. and Rossi, S.},
      year={2025},
      journal={arXiv preprint arXiv:2505.11428}, 
}

@article{Zhai,
 author = {Zhai, X. and Chen, Y. and Li, Y. and Zhao, Y.},
 title = {Optimal well-posedness for the pressureless {Euler}-{Navier}-{Stokes} system},
 fjournal = {Journal of Mathematical Physics},
 journal = {J. Math. Phys.},
 issn = {0022-2488},
 volume = {64},
 number = {5},
 pages = {13},
 note = {Id/No 051506},
 year = {2023},
 language = {English},
 doi = {10.1063/5.0136429},
 zbMATH = {7693291},
 Zbl = {1512.35454}
}

@article{lemarié2025,
      title={{The Pressureless Euler-Navier-Stokes System}}, 
      author={Lemari\'e, V.},
      year={2025},
     journal={arXiv preprint arXiv:2505.17577}, 
}

@article{MJ,
 author = {Hadzic, M. and Speck, J.},
 title = {The global future stability of the {FLRW} solutions to the dust-{Einstein} system with a positive cosmological constant},
 fjournal = {Journal of Hyperbolic Differential Equations},
 journal = {J. Hyperbolic Differ. Equ.},
 volume = {12},
 number = {1},
 pages = {87--188},
 year = {2015},
 language = {English},
}

@article{ChoiJung-comp,
 author = {Choi, Y.-P. and Jung, J.},
 title = {On regular solutions and singularity formation for {Vlasov}/{Navier}-{Stokes} equations with degenerate viscosities and vacuum},
 fjournal = {Kinetic and Related Models},
 journal = {Kinet. Relat. Models},
 volume = {15},
 number = {5},
 pages = {843--891},
 year = {2022},
}

@book{Carles-book,
 author = {Carles, R.},
 title = {Semi-classical analysis for nonlinear {Schr{\"o}dinger} equations},
 isbn = {978-981-279-312-6},
 year = {2008},
 publisher = {Hackensack, NJ: World Scientific},
 language = {English},
 zbMATH = {5243173},
 Zbl = {1153.35070}
}

@phdthesis{oro,
  title={Collective drop effects on vaporizing liquid sprays},
  author={O'Rourke, P.J.},
  year={1981},
  school={Los Alamos National Lab. $\&$ Princeton University}
}

@article{ReitzBook,
  title={Computer modeling of sprays},
  author={Reitz, R.D.},
  journal={Spray Technology Short Course, Pittsburgh, PA},
  year={1996}
}

@article{Jabin,
 author = {Jabin, P.-E.},
 title = {Various levels of models for aerosols.},
 fjournal = {M\(^3\)AS. Mathematical Models \& Methods in Applied Sciences},
 journal = {Math. Models Methods Appl. Sci.},
 issn = {0218-2025},
 volume = {12},
 number = {7},
 pages = {903--919},
 year = {2002},
 language = {English},
 doi = {10.1142/S0218202502001957},
 zbMATH = {1882902},
 Zbl = {1163.35460}
}

@article{BGLM,
 author = {Boudin, L. and Grandmont, C. and Lorz, A. and Moussa, A.},
 title = {Modelling and numerics for respiratory aerosols},
 fjournal = {Communications in Computational Physics},
 journal = {Commun. Comput. Phys.},
 issn = {1815-2406},
 volume = {18},
 number = {3},
 pages = {723--756},
 year = {2015},
 language = {English},
 doi = {10.4208/cicp.180714.200415a},
 zbMATH = {6799575},
 Zbl = {1373.76078}
}

@article{ERS,
 author = {E, W. and Rykov, Y.G. and Sinai, Y.G.},
 title = {Generalized variational principles, global weak solutions and behavior with random initial data for systems of conservation laws arising in adhesion particle dynamics},
 fjournal = {Communications in Mathematical Physics},
 journal = {Commun. Math. Phys.},
 volume = {177},
 number = {2},
 pages = {349--380},
 year = {1996},
}

@book{KrallT,
  title={{Principles of Plasma Physics}},
  author={Krall, N. A. and Trivelpiece, A. W.},
  year={1973},
  publisher={McGraw-Hill}
}

@book{Nicholson,
  title={{Introduction to Plasma Theory}},
  author={Nicholson, D. R.},
  year={1983},
  publisher={John Wiley $\&$ Sons}
}

@book{Williams,
  title={{Combustion Theory}},
  author={Williams, F.},
  year={1985},
  publisher={Benjamin Cummings}
}

@article{PierceH,
  title={A new type of high-frequency amplifier},
  author={Pierce, J. R. and Hebenstreit, W. B.},
  journal={The Bell system technical journal},
  volume={28},
  number={1},
  pages={33--51},
  year={1949},
  publisher={Nokia Bell Labs}
}

@article{BohmGross,
  title={{Theory of plasma oscillations. B. Excitation and damping of oscillations}},
  author={Bohm, D. and Gross, E. P.},
  journal={Physical Review},
  volume={75},
  number={12},
  pages={1864},
  year={1949},
  publisher={APS}
}

@article{Bun,
  title={{Instability, turbulence, and conductivity in current-carrying plasma}},
  author={Buneman, O.},
  journal={Physical Review Letters},
  volume={1},
  number={1},
  pages={8},
  year={1958},
  publisher={APS}
}

@book{chen2015introduction,
  title={Introduction to plasma physics and controlled fusion},
  author={Chen, F.},
  year={2015},
  publisher={Springer}
}

@book{evans2022partial,
  title={Partial differential equations},
  author={Evans, L.C.},
  volume={19},
  year={2022},
  publisher={American mathematical society}
}

@book{Santambrogio,
 author = {Santambrogio, F.},
 title = {Optimal transport for applied mathematicians. {Calculus} of variations, {PDEs}, and modeling},
 fseries = {Progress in Nonlinear Differential Equations and Their Applications},
 series = {Prog. Nonlinear Differ. Equ. Appl.},
 year = {2015},
 publisher = {Cham: Birkh{\"a}user/Springer},
}

@article{danchin2023-dissipativereview,
 author = {Danchin, R.},
 title = {Partially dissipative systems in the critical regularity setting, and strong relaxation limit},
 fjournal = {EMS Surveys in Mathematical Sciences},
 journal = {EMS Surv. Math. Sci.},
 issn = {2308-2151},
 volume = {9},
 number = {1},
 pages = {135--192},
 year = {2022},
}

@article{BF1,
 author = {Bessemoulin-Chatard, M. and Filbet, F.},
 title = {On the stability of conservative discontinuous {Galerkin}/{Hermite} spectral methods for the {Vlasov}-{Poisson} system},
 fjournal = {Journal of Computational Physics},
 journal = {J. Comput. Phys.},
 issn = {0021-9991},
 volume = {451},
 pages = {28},
 note = {Id/No 110881},
 year = {2022},
 language = {English},
 doi = {10.1016/j.jcp.2021.110881},
 zbMATH = {7517167},
 Zbl = {1560.65523}
}

@article{BF2,
 author = {Bessemoulin-Chatard, M. and Filbet, F.},
 title = {On the convergence of discontinuous {Galerkin}/{Hermite} spectral methods for the {Vlasov}-{Poisson} system},
 fjournal = {SIAM Journal on Numerical Analysis},
 journal = {SIAM J. Numer. Anal.},
 issn = {0036-1429},
 volume = {61},
 number = {4},
 pages = {1664--1688},
 year = {2023},
 language = {English},
 doi = {10.1137/22M1518232},
 zbMATH = {7713539},
 Zbl = {1519.65043}
}

@misc{BDFV,
 author = {Blaustein, A. and Dimarco, G. and Filbet, F. and Vignal, M.-H.},
 title = {A structure and asymptotic preserving scheme for the quasineutral limit of the {Vlasov}-{Poisson} system},
 year = {2025},
journal = {arXiv preprint {arXiv}:2504.04826},
 url = {https://arxiv.org/abs/2504.04826},
 arXiv = {arXiv:2504.04826}
}

@article{HC1,
 author = {He, C. and Chen, J.},
 title = {{Vlasov}--{Poisson} equation in weighted {Sobolev} space {{\({W}^{m, p}(w)\)}}},
 fjournal = {Cubo},
 journal = {Cubo},
 issn = {0716-7776},
 volume = {24},
 number = {2},
 pages = {211--226},
 year = {2022},
 language = {English},
 doi = {10.56754/0719-0646.2402.0211},
 zbMATH = {7578199},
 Zbl = {1493.35120}
}

@article{HC2,
 author = {He, C. and Chen, J.},
 title = {{Vlasov}--{Poisson} equation in {Besov} space},
 fjournal = {Taiwanese Journal of Mathematics},
 journal = {Taiwanese J. Math.},
 issn = {1027-5487},
 volume = {26},
 number = {5},
 pages = {1003--1028},
 year = {2022},
 language = {English},
 doi = {10.11650/tjm/220304},
 zbMATH = {7598975},
 Zbl = {1498.35541}
}

@article{HC3,
 author = {Chen, J. and He, C.},
 title = {{Vlasov}--{Poisson} equation in {{\({H}^{s, p}(w)\)}} space},
 fjournal = {Michigan Mathematical Journal},
 journal = {Mich. Math. J.},
 issn = {0026-2285},
 volume = {73},
 number = {3},
 pages = {557--569},
 year = {2023},
 language = {English},
 doi = {10.1307/mmj/20205973},
 zbMATH = {7720193},
 Zbl = {1518.35588}
}

@misc{JT,
 author = {Jeong, I.-J. and Tae, S.},
 title = {Low regularity {Sobolev} well-posedness for {Vlasov}--{Poisson}},
 year = {2025},
 journal = {arXiv preprint {arXiv}:2510.02112 [math.{AP}] (2025)},
 url = {https://arxiv.org/abs/2510.02112},
 arXiv = {arXiv:2510.02112}
}

@article{Tae,
  title={Well-posedness for {Vlasov}--{Poisson} on Low Regularity $H^{s,p}$ Spaces in All Dimensions},
  author={Tae, S.},
  journal={arXiv preprint arXiv:2606.15281},
  year={2026}
}

@article{Ngu,
  title={{Local Well-Posedness for {Vlasov}--{Poisson} with $L^{d+}$ Initial Density and Fractional Velocity Regularity}},
  author={Nguyen, Q.-H.},
  journal={arXiv preprint arXiv:2607.24400},
  year={2026}
}

@article{CISS,
 author = {Crippa, G. and Inversi, M. and Saffirio, C. and Stefani, G.},
 title = {Existence and stability of weak solutions of the {Vlasov}-{Poisson} system in localised {Yudovich} spaces},
 fjournal = {Nonlinearity},
 journal = {Nonlinearity},
 issn = {0951-7715},
 volume = {37},
 number = {9},
 pages = {26},
 note = {Id/No 095015},
 year = {2024},
 language = {English},
 doi = {10.1088/1361-6544/ad5bb3},
 zbMATH = {7892746},
 Zbl = {1545.35208}
}

@article{CKS,
 author = {Choi, Y-P. and  Koo, D. and Song, S.},
 title = {Global existence of {Lagrangian} solutions to the ionic {Vlasov}--{Poisson} system},
 year = {2025},
 journal = {arXiv preprint arXiv:2501.13872},
 url = {https://arxiv.org/abs/2501.13872},
 arXiv = {arXiv:2501.13872}
}

@article{pallard2014-VPtorus,
  title={Space moments of the {Vlasov}--{Poisson} system: propagation and regularity},
  author={Pallard, C.},
  journal={SIAM Journal on Mathematical Analysis},
  volume={46},
  number={3},
  pages={1754--1770},
  year={2014},
  publisher={SIAM}
}

@article{HKI-ions1d,
 author = {Han-Kwan, D. and Iacobelli, M.},
 title = {The quasineutral limit of the {Vlasov}-{Poisson} equation in {Wasserstein} metric},
 fjournal = {Communications in Mathematical Sciences},
 journal = {Commun. Math. Sci.},
 issn = {1539-6746},
 volume = {15},
 number = {2},
 pages = {481--509},
 year = {2017},
}

@article{Ger93,
 author = {G{\'e}rard, P.},
 title = {Remarks on the semiclassical analysis of the nonlinear {Schr{\"o}dinger} equation},
 fjournal = {S{\'e}minaire {\'E}quations aux D{\'e}riv{\'e}es Partielles},
 journal = {S{\'e}min. {\'E}qu. D{\'e}riv. Partielles, {\'E}c. Polytech., Cent. Math. Laurent Schwartz, Palaiseau},
 pages = {ex},
 year = {1993},
 language = {French},
 zbMATH = {1019008},
 Zbl = {0874.35111}
}

@article{AC,
 author = {Alazard, T. and Carles, R.},
 title = {Semi-classical limit of {Schr{\"o}dinger}--{Poisson} equations in space dimension {{\(n{{\geqslant}} 3\)}}},
 fjournal = {Journal of Differential Equations},
 journal = {J. Differ. Equations},
 issn = {0022-0396},
 volume = {233},
 number = {1},
 pages = {241--275},
 year = {2007},
 language = {English},
 doi = {10.1016/j.jde.2006.10.003},
 zbMATH = {5124485},
 Zbl = {1107.35018}
}

@article{Zhang03,
 author = {Zhang, P.},
 title = {Wigner measure and the semiclassical limit of {Schr{\"o}dinger}-{Poisson} equations},
 fjournal = {SIAM Journal on Mathematical Analysis},
 journal = {SIAM J. Math. Anal.},
 issn = {0036-1410},
 volume = {34},
 number = {3},
 pages = {700--718},
 year = {2003},
 language = {English},
 doi = {10.1137/S0036141001393407},
 zbMATH = {1987564},
 Zbl = {1032.35132}
}

@article{ZZM,
 author = {Zhang, P. and Zheng, Y. and Mauser, N.},
 title = {The limit from the {Schr{\"o}dinger}-{Poisson} to the {Vlasov}-{Poisson} equations with general data in one dimension},
 fjournal = {Communications on Pure and Applied Mathematics},
 journal = {Commun. Pure Appl. Math.},
 issn = {0010-3640},
 volume = {55},
 number = {5},
 pages = {582--632},
 year = {2002},
 language = {English},
 doi = {10.1002/cpa.3017},
 zbMATH = {1860573},
 Zbl = {1032.81011}
}

@article{CarlesM,
 author = {Carles, R. and Masaki, S.},
 title = {Semiclassical analysis for {Hartree} equation},
 fjournal = {Asymptotic Analysis},
 journal = {Asymptotic Anal.},
 issn = {0921-7134},
 volume = {58},
 number = {4},
 pages = {211--227},
 year = {2008},
 language = {English},
 zbMATH = {5521949},
 Zbl = {1161.35310}
}

@article{Mas,
 author = {Masaki, S.},
 title = {Local existence and {WKB} approximation of solutions to {Schr{\"o}dinger}-{Poisson} system in the two-dimensional whole space},
 fjournal = {Communications in Partial Differential Equations},
 journal = {Commun. Partial Differ. Equations},
 issn = {0360-5302},
 volume = {35},
 number = {10-12},
 pages = {2253--2278},
 year = {2010},
 language = {English},
 doi = {10.1080/03605301003717142},
 zbMATH = {5839295},
 Zbl = {1232.35155}
}

@article{Gr98,
 author = {Grenier, E.},
 title = {Semiclassical limit of the nonlinear {Schr{\"o}dinger} equation in small time},
 fjournal = {Proceedings of the American Mathematical Society},
 journal = {Proc. Am. Math. Soc.},
 issn = {0002-9939},
 volume = {126},
 number = {2},
 pages = {523--530},
 year = {1998},
 language = {English},
 doi = {10.1090/S0002-9939-98-04164-1},
 zbMATH = {1132286},
 Zbl = {0910.35115}
}

@article{LP,
 author = {Lions, P-L. and Paul, T.},
 title = {On {Wigner} measures},
 fjournal = {Revista Matem{\'a}tica Iberoamericana},
 journal = {Rev. Mat. Iberoam.},
 issn = {0213-2230},
 volume = {9},
 number = {3},
 pages = {553--618},
 year = {1993},
 language = {French},
 doi = {10.4171/RMI/143},
 url = {https://eudml.org/doc/39445},
 zbMATH = {482230},
 Zbl = {0801.35117}
}

@article{GP,
 author = {Golse, F. and Paul, T.},
 title = {The {Schr{\"o}dinger} equation in the mean-field and semiclassical regime},
 fjournal = {Archive for Rational Mechanics and Analysis},
 journal = {Arch. Ration. Mech. Anal.},
 issn = {0003-9527},
 volume = {223},
 number = {1},
 pages = {57--94},
 year = {2017},
 language = {English},
 doi = {10.1007/s00205-016-1031-x},
 zbMATH = {6688557},
 Zbl = {1359.35164}
}

@article{Lafleche,
 author = {Lafleche, L.},
 title = {Propagation of moments and semiclassical limit from {Hartree} to {Vlasov} equation},
 fjournal = {Journal of Statistical Physics},
 journal = {J. Stat. Phys.},
 issn = {0022-4715},
 volume = {177},
 number = {1},
 pages = {20--60},
 year = {2019},
 language = {English},
 doi = {10.1007/s10955-019-02356-7},
 zbMATH = {7115606},
 Zbl = {1426.82034}
}

@article{IacobelliLafleche,
 author = {Iacobelli, M. and Lafleche, L.},
 title = {Enhanced stability in quantum optimal transport pseudometrics: from {Hartree} to {Vlasov}-{Poisson}},
 fjournal = {Journal of Statistical Physics},
 journal = {J. Stat. Phys.},
 issn = {0022-4715},
 volume = {191},
 number = {12},
 pages = {18},
 note = {Id/No 157},
 year = {2024},
 language = {English},
 doi = {10.1007/s10955-024-03367-9},
 zbMATH = {7955012},
 Zbl = {1563.81051}
}

@article{LS,
 author = {Lafleche, L. and Saffirio, C.},
 title = {Strong semiclassical limits from {Hartree} and {Hartree}-{Fock} to {Vlasov}-{Poisson} equations},
 fjournal = {Analysis \& PDE},
 journal = {Anal. PDE},
 issn = {2157-5045},
 volume = {16},
 number = {4},
 pages = {891--926},
 year = {2023},
 language = {English},
 doi = {10.2140/apde.2023.16.891},
 zbMATH = {7713387},
 Zbl = {1515.35228}
}

@article{CLS,
 author = {Chong,  J. and Lafleche, L. and Saffirio, C.},
 title = {On the {{\(L^2\)}} rate of convergence in the limit from the {Hartree} to the {Vlasov}-{Poisson} equation},
 fjournal = {Journal de l'{\'E}cole Polytechnique -- Math{\'e}matiques},
 journal = {J. {\'E}c. Polytech., Math.},
 issn = {2429-7100},
 volume = {10},
 pages = {703--726},
 year = {2023},
 language = {English},
 doi = {10.5802/jep.230},
 zbMATH = {7680771},
 Zbl = {1522.81084}
}

@incollection{LS2,
 author = {Lafleche, L. and Saffirio, C.},
 title = {Uniqueness criteria for the {Vlasov}-{Poisson} system and applications to semiclassical analysis},
 booktitle = {From particle systems to partial differential equations. PSPDE X, Braga, Portugal, June 27 -- July 1, 2022},
 isbn = {978-3-031-65194-6; 978-3-031-65197-7; 978-3-031-65195-3},
 pages = {301--317},
 year = {2024},
 publisher = {Cham: Springer},
 language = {English},
 doi = {10.1007/978-3-031-65195-3_14},
 zbMATH = {8025156},
 Zbl = {1564.35262}
}

@article{Br00,
 author = {Brenier, Y.},
 title = {Convergence of the {Vlasov}-{Poisson} system to the incompressible {Euler} equations},
 fjournal = {Communications in Partial Differential Equations},
 journal = {Commun. Partial Differ. Equations},
 issn = {0360-5302},
 volume = {25},
 number = {3-4},
 pages = {737--754},
 year = {2000},
 language = {English},
 doi = {10.1080/03605300008821529},
 zbMATH = {1442880},
 Zbl = {0970.35110}
}

@article{HK11,
 author = {Han-Kwan, D.},
 title = {Quasineutral limit of the {Vlasov}-{Poisson} system with massless electrons},
 fjournal = {Communications in Partial Differential Equations},
 journal = {Commun. Partial Differ. Equations},
 issn = {0360-5302},
 volume = {36},
 number = {7-9},
 pages = {1385--1425},
 year = {2011},
 language = {English},
 doi = {10.1080/03605302.2011.555804},
 zbMATH = {5963960},
 Zbl = {1228.35251}
}

@article{HKR,
 author = {Han-Kwan, D. and Rousset, F.},
 title = {Quasineutral limit for {Vlasov}-{Poisson} with {Penrose} stable data},
 fjournal = {Annales Scientifiques de l'{\'E}cole Normale Sup{\'e}rieure. Quatri{\`e}me S{\'e}rie},
 journal = {Ann. Sci. {\'E}c. Norm. Sup{\'e}r. (4)},
 issn = {0012-9593},
 volume = {49},
 number = {6},
 pages = {1445--1495},
 year = {2016},
 language = {English},
 doi = {10.24033/asens.2313},
 url = {smf4.emath.fr/en/Publications/AnnalesENS/4_49/html/ens_ann-sc_49_1445-1495.php},
 zbMATH = {6680023},
 Zbl = {1361.35179}
}

@misc{GPI2,
 author = {Griffin-Pickering, M. and Iacobelli, M.},
 title = {Stability in {Quasineutral} {Plasmas} with {Thermalized} {Electrons}},
 year = {2023},
journal = {arXiv preprint {arXiv}:2307.07561 [math.{AP}] (2023)},
 url = {https://arxiv.org/abs/2307.07561},
 arXiv = {arXiv:2307.07561}
}

@article{Iacobelli,
 author = {Iacobelli, M.},
 title = {A new perspective on {Wasserstein} distances for kinetic problems},
 fjournal = {Archive for Rational Mechanics and Analysis},
 journal = {Arch. Ration. Mech. Anal.},
 issn = {0003-9527},
 volume = {244},
 number = {1},
 pages = {27--50},
 year = {2022},
 language = {English},
 doi = {10.1007/s00205-021-01705-9},
 zbMATH = {7505277},
 Zbl = {1507.35181}
}

@article{BreLoe,
 author = {Brenier, Y. and Loeper, G.},
 title = {A geometric approximation to the {Euler} equations: the {Vlasov}-{Monge}-{Amp{\`e}re} system},
 fjournal = {Geometric and Functional Analysis. GAFA},
 journal = {Geom. Funct. Anal.},
 issn = {1016-443X},
 volume = {14},
 number = {6},
 pages = {1182--1218},
 year = {2004},
 language = {English},
 doi = {10.1007/s00039-004-0488-1},
 zbMATH = {2162381},
 Zbl = {1075.35046}
}

@article{ABB,
 author = {Ambrosio, L. and Baradat, A. and Brenier, Y.},
 title = {Monge-Amp{\`e}re gravitation as a {{\({{\Gamma}} \)}}-limit of good rate functions},
 fjournal = {Analysis \& PDE},
 journal = {Anal. PDE},
 issn = {2157-5045},
 volume = {16},
 number = {9},
 pages = {2005--2040},
 year = {2023},
 language = {English},
 doi = {10.2140/apde.2023.16.2005},
 zbMATH = {7785255},
 Zbl = {1536.49017}
}

@article {shizuta1983classical,
    AUTHOR = {Shizuta, Y.},
     TITLE = {On the classical solutions of the {B}oltzmann equation},
   JOURNAL = {Comm. Pure Appl. Math.},
  FJOURNAL = {Communications on Pure and Applied Mathematics},
    VOLUME = {36},
      YEAR = {1983},
    NUMBER = {6},
     PAGES = {705--754},
      ISSN = {0010-3640,1097-0312},
   MRCLASS = {76P05 (35Q20 45K05 47D05 82A40)},
  MRNUMBER = {720591},
MRREVIEWER = {Reinhard\ Illner},
       DOI = {10.1002/cpa.3160360602},
       URL = {https://doi.org/10.1002/cpa.3160360602},
}

@article{CCM2,
 author = {Caprino, S. and Cavallaro, G. and Marchioro, C.},
 title = {The {Vlasov}-{Poisson} equation in {{\(\mathbb{R}^3\)}} with infinite charge and velocities},
 fjournal = {Journal of Hyperbolic Differential Equations},
 journal = {J. Hyperbolic Differ. Equ.},
 issn = {0219-8916},
 volume = {15},
 number = {3},
 pages = {407--442},
 year = {2018},
 language = {English},
 doi = {10.1142/S0219891618500157},
 zbMATH = {7059029},
 Zbl = {1416.82041}
}

@article{CCM1,
 author = {Caprino, S. and Cavallaro, G. and Marchioro, C.},
 title = {Time evolution of a {Vlasov}-{Poisson} plasma with infinite charge in {{\(\mathbb{R}^{3}\)}}},
 fjournal = {Communications in Partial Differential Equations},
 journal = {Commun. Partial Differ. Equations},
 issn = {0360-5302},
 volume = {40},
 number = {2},
 pages = {357--385},
 year = {2015},
 language = {English},
 doi = {10.1080/03605302.2014.944267},
 zbMATH = {6428304},
 Zbl = {1320.82062}
}

@article{J00,
 author = {Jabin, P.-E.},
 title = {The {Vlasov}-{Poisson} system with infinite mass and energy.},
 fjournal = {Journal of Statistical Physics},
 journal = {J. Stat. Phys.},
 issn = {0022-4715},
 volume = {103},
 number = {5-6},
 pages = {1107--1123},
 year = {2001},
 language = {English},
 doi = {10.1023/A:1010321308267},
 zbMATH = {1674997},
 Zbl = {1126.82327}
}

@article{Danchin-ENS2026,
      title={An elementary approach to the pressureless {E}uler-{N}avier-{S}tokes system}, 
      author={Danchin, R.},
      year={2026},
      eprint={2602.06821},
      journal={arXiv preprint arXiv:2602.06821}
}

@article{luk2020high,
  title={High-frequency limits and null dust shell solutions in general relativity},
  author={Luk, J. and Rodnianski, I.},
  journal={arXiv preprint arXiv:2009.08968},
  year={2020}
}

@article{BardosGolseNguyenSentis,
 author = {Bardos, C. and Golse, F. and Nguyen, T.T. and Sentis, R.},
 title = {The {Maxwell}-{Boltzmann} approximation for ion kinetic modeling},
 fjournal = {Physica D},
 journal = {Physica D},
 issn = {0167-2789},
 volume = {376-377},
 pages = {94--107},
 year = {2018},
}

@article{Tony,
  title={{Semi-classical limit of the massive Klein-Gordon-Maxwell system toward the relativistic Euler-Maxwell system via an adapted modulated energy method}},
  author={Salvi, T.},
  journal={arXiv preprint arXiv:2502.06622},
  year={2025}
}

@article {tadmor2022critical,
    AUTHOR = {Tadmor, E. and Tan, C.},
     TITLE = {Critical threshold for global regularity of the
              {E}uler-{M}onge-{A}mp\`ere system with radial symmetry},
   JOURNAL = {SIAM J. Math. Anal.},
  FJOURNAL = {SIAM Journal on Mathematical Analysis},
    VOLUME = {54},
      YEAR = {2022},
    NUMBER = {4},
     PAGES = {4277--4296},
      ISSN = {0036-1410,1095-7154},
   MRCLASS = {35Q35 (35B30 35K96 76N10)},
  MRNUMBER = {4451915},
       DOI = {10.1137/21M1437767},
       URL = {https://doi-org.ezproxy.math.cnrs.fr/10.1137/21M1437767},
}

@article {loeper2005quasi,
    AUTHOR = {Loeper, G.},
     TITLE = {Quasi-neutral limit of the {E}uler-{P}oisson and
              {E}uler-{M}onge-{A}mp\`ere systems},
   JOURNAL = {Comm. Partial Differential Equations},
  FJOURNAL = {Communications in Partial Differential Equations},
    VOLUME = {30},
      YEAR = {2005},
    NUMBER = {7-9},
     PAGES = {1141--1167},
      ISSN = {0360-5302,1532-4133},
   MRCLASS = {35Q35 (35B25 76W05)},
  MRNUMBER = {2180297},
MRREVIEWER = {Xiaoming\ Wang},
       DOI = {10.1080/03605300500257545},
       URL = {https://doi-org.ezproxy.math.cnrs.fr/10.1080/03605300500257545},
}

@article {ambrosio2008hamiltonian,
    AUTHOR = {Ambrosio, L. and Gangbo, W.},
     TITLE = {Hamiltonian {ODE}s in the {W}asserstein space of probability
              measures},
   JOURNAL = {Comm. Pure Appl. Math.},
  FJOURNAL = {Communications on Pure and Applied Mathematics},
    VOLUME = {61},
      YEAR = {2008},
    NUMBER = {1},
     PAGES = {18--53},
      ISSN = {0010-3640,1097-0312},
   MRCLASS = {37J05 (28A33 34F05 49J52 60B10)},
  MRNUMBER = {2361303},
MRREVIEWER = {C\'edric\ Villani},
       DOI = {10.1002/cpa.20188},
       URL = {https://doi-org.ezproxy.math.cnrs.fr/10.1002/cpa.20188},
}

@article {brenier2016double,
    AUTHOR = {Brenier, Y.},
     TITLE = {A double large deviation principle for {M}onge-{A}mp\`ere
              gravitation},
   JOURNAL = {Bull. Inst. Math. Acad. Sin. (N.S.)},
  FJOURNAL = {Bulletin of the Institute of Mathematics. Academia Sinica. New
              Series},
    VOLUME = {11},
      YEAR = {2016},
    NUMBER = {1},
     PAGES = {23--41},
      ISSN = {2304-7909,2304-7895},
   MRCLASS = {35Q82 (35Q83 49S05 60F10)},
  MRNUMBER = {3497744},
MRREVIEWER = {Giuseppe\ Genovese},
}

@article{levy2024monge,
  title = {Monge-Amp\`ere gravity: From the large deviation principle to cosmological simulations through optimal transport},
  author = {L\'evy, B. and Brenier, Y. and Mohayaee, R.},
  journal = {Phys. Rev. D},
  volume = {110},
  issue = {6},
  pages = {063550},
  numpages = {15},
  year = {2024},
  month = {Sep},
  publisher = {American Physical Society},
  doi = {10.1103/PhysRevD.110.063550},
  url = {https://link.aps.org/doi/10.1103/PhysRevD.110.063550}
}

@article{leonard2026monge,
  title="{Monge-Ampère gravitating fluids. Least action principles and particle systems}",
  author={Léonard, C. and Mohayaee, R.},
  journal={arXiv preprint arXiv:2503.01537},
  year={2026}
}

@article{Herau1,
 author = {H{\'e}rau, F.},
 title = {Hypocoercivity and exponential time decay for the linear inhomogeneous relaxation {Boltzmann} equation},
 fjournal = {Asymptotic Analysis},
 journal = {Asymptotic Anal.},
 issn = {0921-7134},
 volume = {46},
 number = {3-4},
 pages = {349--359},
 year = {2006},
}

@article{HerauNier,
 author = {H{\'e}rau, F. and Nier, F.},
 title = {Isotropic hypoelliptic and trend to equilibrium for the {Fokker}-{Planck} equation with a high-degree potential},
 fjournal = {Archive for Rational Mechanics and Analysis},
 journal = {Arch. Ration. Mech. Anal.},
 issn = {0003-9527},
 volume = {171},
 number = {2},
 pages = {151--218},
 year = {2004},
}

@article{DolbeaultMouhotSchmeiser,
 author = {Dolbeault, J. and Mouhot, C. and Schmeiser, C.},
 title = {Hypocoercivity for linear kinetic equations conserving mass},
 fjournal = {Transactions of the American Mathematical Society},
 journal = {Trans. Am. Math. Soc.},
 issn = {0002-9947},
 volume = {367},
 number = {6},
 pages = {3807--3828},
 year = {2015},
}

@article{MouhotNeumann,
 author = {Mouhot, C. and Neumann, L.},
 title = {Quantitative perturbative study of convergence to equilibrium for collisional kinetic models in the torus},
 fjournal = {Nonlinearity},
 journal = {Nonlinearity},
 issn = {0951-7715},
 volume = {19},
 number = {4},
 pages = {969--998},
 year = {2006},
}

@book{Villani,
 author = {Villani, C.},
 title = {Hypocoercivity},
 fseries = {Memoirs of the American Mathematical Society},
 series = {Mem. Am. Math. Soc.},
 issn = {0065-9266},
 volume = {950},
 isbn = {978-0-8218-4498-4; 978-1-4704-0564-9},
 year = {2009},
 publisher = {Providence, RI: American Mathematical Society (AMS)},
}

@article{Zhi,
 author = {Zhidkov, P. E.},
 title = {Existence of solutions to the {Cauchy} problem and stability of kink- solutions of the nonlinear {Schr{\"o}dinger} equation},
 fjournal = {Siberian Mathematical Journal},
 journal = {Sib. Math. J.},
 issn = {0037-4466},
 volume = {33},
 number = {2},
 pages = {73--79},
 year = {1992},
 language = {English},
 doi = {10.1007/BF00971094},
 zbMATH = {130774},
 Zbl = {0783.35074}
}

@article{Gallo,
 author = {Gallo, C.},
 title = {Schr{\"o}dinger group on {Zhidkov} spaces},
 fjournal = {Advances in Differential Equations},
 journal = {Adv. Differ. Equ.},
 issn = {1079-9389},
 volume = {9},
 number = {5-6},
 pages = {509--538},
 year = {2004},
 language = {English},
 zbMATH = {5054483},
 Zbl = {1103.35093}
}

@article{BF,
  title={Non linear electron plasma oscillation: the “water bag model”},
  author={Bertrand, P. and Feix, M.},
  journal={Physics Letters A},
  volume={28},
  number={1},
  pages={68--69},
  year={1968},
  publisher={Elsevier}
}

@article{NBwater,
  title={Multiple “water-bag” model and Landau damping},
  author={Navet, M. and Bertrand, P.},
  journal={Physics Letters A},
  volume={34},
  number={2},
  pages={117--118},
  year={1971},
  publisher={Elsevier}
}

@article{BDBF,
  title={Stability of inhomogeneous two-stream plasma with a water-bag model},
  author={Bertrand, P. and Doremus, J.P. and Baumann, G. and Feix, M.R.},
  journal={The Physics of Fluids},
  volume={15},
  number={7},
  pages={1275--1281},
  year={1972},
  publisher={AIP Publishing}
}

@article{EJK,
  title={Global dynamics in the {3D} pressureless {Euler}-{Poisson} system for ions},
  author={Ertzbischoff, L. and Jurja, C. and Widmayer, K.},
  journal={arXiv preprint arXiv:2606.19177},
  year={2026}
}

@article{song2026h-highMach,
  title={{High Mach number limit for the 3D Euler-Poisson equations of ion dynamics}},
  author={Song, Z.},
  journal={arXiv preprint arXiv:2606.17863},
  year={2026}
}

@article{LemarieMultiphaseVNS,
      title={Multiphase formulation of the {Vlasov}-{Navier}-{Stokes} equations}, 
      author={Lemari\'e, V.},
      year={2026},
      journal={arXiv preprint arXiv:2606.30022}, 
}

@article{Davidson,
      title={Kinetic waves and instabilities in a uniform plasma}, 
      author={Davidson, R. C.},
      year={1983},
      journal={Handbook of
Plasma Physics, Volume 1}, 
}

@article{bae2025emergence,
  title={{Emergence of Peaked Singularities in the Euler--Poisson System}},
  author={Bae, J. and Moon, S.-H. and Woo, K.},
 fjournal = {Journal of Nonlinear Science},
 journal = {J. Nonlinear Sci.},  volume={35},
  number={1},
  pages={25},
  year={2025},
}

@article{choi2025critical,
  title={{Critical thresholds in pressureless Euler--Poisson equations with background states}},
  author={Choi, Y.-P. and Kim, D.-H. and Koo, D. and Tadmor, E.},
  journal={Annales de l'Institut Henri Poincar{\'e} C},
  volume={43},
  number={1},
  pages={203--237},
  year={2025}
}

@article{BaeKimKwon-singularity,
 author = {Bae, J. and Kim, Y. and Kwon, B.},
 title = {Delta-shock for the pressureless {Euler}-{Poisson} system},
 fjournal = {SIAM Journal on Mathematical Analysis},
 journal = {SIAM J. Math. Anal.},
 issn = {0036-1410},
 volume = {57},
 number = {3},
 pages = {3255--3296},
 year = {2025},
}

@article{MouhotVillani,
 author = {Mouhot, C. and Villani, C.},
 title = {On {Landau} damping},
 fjournal = {Acta Mathematica},
 journal = {Acta Math.},
 issn = {0001-5962},
 volume = {207},
 number = {1},
 pages = {29--201},
 year = {2011},
 language = {English},
 doi = {10.1007/s11511-011-0068-9},
 zbMATH = {6012919},
 Zbl = {1239.82017}
}

@article{BMM,
 author = {Bedrossian, J. and Masmoudi, N. and Mouhot, C.},
 title = {Landau damping: paraproducts and {Gevrey} regularity},
 fjournal = {Annals of PDE},
 journal = {Ann. PDE},
 issn = {2524-5317},
 volume = {2},
 number = {1},
 pages = {71},
 note = {Id/No 4},
 year = {2016},
 language = {English},
 doi = {10.1007/s40818-016-0008-2},
 zbMATH = {6919581},
 Zbl = {1402.35058}
}

@article{GNR,
 author = {Grenier, E. and Nguyen, T. and Rodnianski, I.},
 title = {Landau damping for analytic and {Gevrey} data},
 fjournal = {Mathematical Research Letters},
 journal = {Math. Res. Lett.},
 issn = {1073-2780},
 volume = {28},
 number = {6},
 pages = {1679--1702},
 year = {2021},
 language = {English},
 doi = {10.4310/MRL.2021.v28.n6.a3},
 zbMATH = {7579208},
 Zbl = {1496.35080}
}

@misc{IPWW,
 author = {A. D. Ionescu and B. Pausader and X. Wang and K. Widmayer},
 title = {Nonlinear {Landau} damping and wave operators in sharp {Gevrey} spaces},
 year = {2024},
 journal = {arxiv preprint {arXiv}:2405.04473 (2024)},
 url = {https://arxiv.org/abs/2405.04473},
 arXiv = {arXiv:2405.04473}
}

@article{vidav1970,
title = {Spectra of perturbed semigroups with applications to transport theory},
journal = {Journal of Mathematical Analysis and Applications},
volume = {30},
number = {2},
pages = {264-279},
year = {1970},
issn = {0022-247X},
doi = {https://doi.org/10.1016/0022-247X(70)90160-5},
url = {https://www.sciencedirect.com/science/article/pii/0022247X70901605},
author = {Vidav, I.}
}

@article{GagnebinIacobelli,
 author = {Gagnebin, A. and Iacobelli, M.},
 title = {Landau damping on the torus for the {Vlasov}-{Poisson} system with massless electrons},
 fjournal = {Journal of Differential Equations},
 journal = {J. Differ. Equations},
 issn = {0022-0396},
 volume = {376},
 pages = {154--203},
 year = {2023},
 language = {English},
 doi = {10.1016/j.jde.2023.08.020},
 zbMATH = {7757708},
 Zbl = {1528.35202}
}

@article{HNX,
 author = {Huang, L. and Nguyen, Q.-H. and Xu, Y.},
 title = {Nonlinear {Landau} damping for the {2D} {Vlasov}-{Poisson} system with massless electrons around {Penrose}-stable equilibrium},
 fjournal = {SIAM Journal on Mathematical Analysis},
 journal = {SIAM J. Math. Anal.},
 issn = {0036-1410},
 volume = {57},
 number = {2},
 pages = {1939--1963},
 year = {2025},
 language = {English},
 doi = {10.1137/23M1595382},
 zbMATH = {8043976},
 Zbl = {1565.35326}
}

@article{BVR,
 author = {Bigorgne, L. and Velozo Ruiz, R.},
 title = {Late-time asymptotics of small data solutions for the {Vlasov}-{Poisson} system},
 fjournal = {Nonlinearity},
 journal = {Nonlinearity},
 issn = {0951-7715},
 volume = {39},
 number = {4},
 pages = {54},
 note = {Id/No 045007},
 year = {2026},
 language = {English},
 doi = {10.1088/1361-6544/ae55f5},
 zbMATH = {8194646}
}

@article{HRV,
 author = {Hwang, H. J. and Rendall, A. D. and Vel{\'a}zquez, J. J. L.},
 title = {Optimal gradient estimates and asymptotic behaviour for the {Vlasov}-{Poisson} system with small initial data},
 fjournal = {Archive for Rational Mechanics and Analysis},
 journal = {Arch. Ration. Mech. Anal.},
 issn = {0003-9527},
 volume = {200},
 number = {1},
 pages = {313--360},
 year = {2011},
 language = {English},
 doi = {10.1007/s00205-011-0405-3},
 url = {hdl.handle.net/11858/00-001M-0000-0013-5F9F-5},
 zbMATH = {5952976},
 Zbl = {1228.35252}
}

@article{Smu,
 author = {Smulevici, J.},
 title = {Small data solutions of the {Vlasov}-{Poisson} system and the vector field method},
 fjournal = {Annals of PDE},
 journal = {Ann. PDE},
 issn = {2524-5317},
 volume = {2},
 number = {2},
 pages = {55},
 note = {Id/No 11},
 year = {2016},
 language = {English},
 doi = {10.1007/s40818-016-0016-2},
 zbMATH = {6919588},
 Zbl = {1397.35033}
}

@article{CK-VP,
 author = {Choi, S.-H. and Kwon, S.},
 title = {Modified scattering for the {Vlasov}-{Poisson} system},
 fjournal = {Nonlinearity},
 journal = {Nonlinearity},
 issn = {0951-7715},
 volume = {29},
 number = {9},
 pages = {2755--2774},
 year = {2016},
 language = {English},
 doi = {10.1088/0951-7715/29/9/2755},
 zbMATH = {6629811},
 Zbl = {1351.82041}
}

@article{IPWW-scat,
 author = {Ionescu, A. D. and Pausader, B. and Wang, X. and Widmayer, K.},
 title = {On the asymptotic behavior of solutions to the {Vlasov}-{Poisson} system},
 fjournal = {IMRN. International Mathematics Research Notices},
 journal = {Int. Math. Res. Not.},
 issn = {1073-7928},
 volume = {2022},
 number = {12},
 pages = {8865--8889},
 year = {2022},
 language = {English},
 doi = {10.1093/imrn/rnab155},
 zbMATH = {7542580},
 Zbl = {1491.35082}
}

@article{IRW,
 author = {Iacobelli, M. and Rossi, S. and Widmayer, K.},
 title = {On the stability of vacuum in the screened {Vlasov}-{Poisson} equation},
 fjournal = {Journal of the London Mathematical Society. Second Series},
 journal = {J. Lond. Math. Soc., II. Ser.},
 issn = {0024-6107},
 volume = {113},
 number = {1},
 pages = {37},
 note = {Id/No e70426},
 year = {2026},
}

@article{Lin1,
 author = {Lin, Z.},
 title = {Instability of periodic {BGK} waves},
 fjournal = {Mathematical Research Letters},
 journal = {Math. Res. Lett.},
 issn = {1073-2780},
 volume = {8},
 number = {4},
 pages = {521--534},
 year = {2001},
 language = {English},
 doi = {10.4310/MRL.2001.v8.n4.a11},
 zbMATH = {1690765},
 Zbl = {0993.35084}
}

@article{Lin2,
 author = {Lin, Z.},
 title = {Nonlinear instability of periodic {BGK} waves for {Vlasov}-{Poisson} system},
 fjournal = {Communications on Pure and Applied Mathematics},
 journal = {Commun. Pure Appl. Math.},
 issn = {0010-3640},
 volume = {58},
 number = {4},
 pages = {505--528},
 year = {2005},
 language = {English},
 doi = {10.1002/cpa.20028},
 zbMATH = {2165077},
 Zbl = {1067.35012}
}

@article{GS-BGK,
 author = {Guo, Y. and Strauss, W. A.},
 title = {Instability of periodic {BGK} equilibria},
 fjournal = {Communications on Pure and Applied Mathematics},
 journal = {Commun. Pure Appl. Math.},
 issn = {0010-3640},
 volume = {48},
 number = {8},
 pages = {861--894},
 year = {1995},
 language = {English},
 doi = {10.1002/cpa.3160480803},
 zbMATH = {846110},
 Zbl = {0840.45012}
}

@article{GL,
 author = {Guo, Y. and Lin, Z.},
 title = {The existence of stable {BGK} waves},
 fjournal = {Communications in Mathematical Physics},
 journal = {Commun. Math. Phys.},
 issn = {0010-3616},
 volume = {352},
 number = {3},
 pages = {1121--1152},
 year = {2017},
 language = {English},
 doi = {10.1007/s00220-017-2873-2},
 zbMATH = {6721419},
 Zbl = {1373.82089}
}

@misc{BGHP,
 author = {Bian, D. and  Grenier, E. and  Huang, W. and  Pausader, B.},
 title = {Stability and instability of small {BGK} waves},
 year = {2026},
 journal = {arXiv preprint {arXiv}:2601.10030 [math.{AP}] (2026)},
 url = {https://arxiv.org/abs/2601.10030},
 arXiv = {arXiv:2601.10030}
}

@article{dziurzynski1987patches,
  title={Patches of electrons and electron sheets for the 1-D {Vlasov}-{Poisson} equation},
  author={Dziurzynski, R.S.},
  year={1987},
  journal={PhD Thesis - University of California, Berkeley}
}

@article{Roulley,
 author = {Roulley, E.},
 title = {Local and global bifurcation of electron-states},
 fjournal = {Discrete and Continuous Dynamical Systems},
 journal = {Discrete Contin. Dyn. Syst.},
 issn = {1078-0947},
 volume = {45},
 number = {8},
 pages = {2381--2419},
 year = {2025},
}

@article{DGR-deriv,
 author = {Desvillettes, L. and Golse, F. and Ricci, V.},
 title = {The mean-field limit for solid particles in a {Navier}-{Stokes} flow},
 fjournal = {Journal of Statistical Physics},
 journal = {J. Stat. Phys.},
 issn = {0022-4715},
 volume = {131},
 number = {5},
 pages = {941--967},
 year = {2008},
}

@article{HoferPHD,
  title={Sedimentation of particle suspensions in {S}tokes flows},
  author={H{\"o}fer, R.},
  year={2020},
  journal={Universit{\"a}ts-und Landesbibliothek Bonn}
}

@article{Dechicha,
 author = {Dechicha, D.},
 title = {Gevrey regularity and analyticity for the solutions of the {Vlasov}-{Navier}-{Stokes} system},
 fjournal = {SIAM Journal on Mathematical Analysis},
 journal = {SIAM J. Math. Anal.},
 issn = {0036-1410},
 volume = {56},
 number = {6},
 pages = {7903--7939},
 year = {2024},
}

@book{EHK-thick,
 author = {Ertzbischoff, L. and Han-Kwan, D.},
 title = {On well-posedness for thick spray equations},
 fseries = {Memoirs of the European Mathematical Society},
 series = {Mem. Eur. Math. Soc.},
 issn = {2747-9080},
 volume = {26},
 isbn = {978-3-98547-097-6; 978-3-98547-597-1},
 year = {2025},
 publisher = {Berlin: European Mathematical Society (EMS)},
 language = {English},
 doi = {10.4171/MEMS/26},
 zbMATH = {8144695}
}

@article{ACF,
 author = {Ambrosio, L. and Colombo, M. and Figalli, A.},
 title = {On the {Lagrangian} structure of transport equations: the {Vlasov}-{Poisson} system},
 fjournal = {Duke Mathematical Journal},
 journal = {Duke Math. J.},
 issn = {0012-7094},
 volume = {166},
 number = {18},
 pages = {3505--3568},
 year = {2017},
 language = {English},
}

@book{galdi2000introduction,
  title={An introduction to the {Navier}-{Stokes} initial-boundary value problem},
  author={Galdi, G.P.},
  booktitle={Fundamental directions in mathematical fluid mechanics},
  pages={1--70},
  year={2000},
  publisher={Springer}
}

@book{RobinsonRodrigoSadowskiNS3,
 author = {Robinson, James C. and Rodrigo, Jos{\'e} L. and Sadowski, Witold},
 title = {The three-dimensional {Navier}-{Stokes} equations. {Classical} theory},
 fseries = {Cambridge Studies in Advanced Mathematics},
 series = {Camb. Stud. Adv. Math.},
 volume = {157},
 isbn = {978-1-107-01966-9; 978-1-139-09514-3},
 year = {2016},
 publisher = {Cambridge: Cambridge University Press},
}
 \bibliographystyle{alpha}
  \end{document}